\documentclass[11pt]{amsart}

\usepackage{amstext}
\usepackage{amsfonts} 
\usepackage{amsthm,euscript}
\usepackage{amscd}
\usepackage{amssymb}
\usepackage{amsmath}
\usepackage{mathrsfs}
\usepackage{mathtools}
\usepackage[all]{xy}
\usepackage{paralist}
\usepackage{enumitem}
\usepackage{moreenum}
\usepackage{scalerel,stackengine}
\usepackage[colorlinks=true, pdfstartview=FitV, linkcolor=blue,
citecolor=blue,urlcolor=blue,breaklinks=true]{hyperref}
\usepackage{subfiles} 
\usepackage{mathdots}

\swapnumbers 
\newcommand{\lv}[1]{} 
\newcommand{\comments}[1]{}
\newcommand{\pcomments}[1]{}
\newcommand{\inparcom}[1]{}
\newcounter{question} 
\newcommand{\quest}[1]{} 
\newcommand{\new}{}
\newcommand{\enew}{}
\newcommand{\heute}{\tiny (\today)}
\newcommand{\arxiv}[1]{\href{http://arxiv.org/abs/#1}{\tt arXiv:\nolinkurl{#1}}}
\newtheorem{thm}[subsection]{Theorem}

\newtheorem{prop}[subsection]{Proposition}

\newtheorem{lem}[subsection]{Lemma}

\newtheorem{cor}[subsection]{Corollary}

\theoremstyle{definition}
\newtheorem{defn}[subsection]{Definition}

\newtheorem{remark}[subsection]{Remark}

\newtheorem{remarks}[subsection]{Remarks}

\newtheorem{example}[subsection]{Example}

\newtheorem{examples}[subsection]{Examples}

\newtheorem{leer}[subsection]{\bf}%

\numberwithin{equation}{subsection}

\newcommand\sm{\smallskip}
\newcommand\ms{\medskip}

\newcommand{\ppmatrix}[4]{%
  { \left(\begin{smallmatrix} #1 & #2 \\ #3 & #4 \end{smallmatrix}\right)}%
}

\def\limind{\mathop{\oalign{lim\cr\hidewidth$\longrightarrow$\hidewidth \cr}}}
\def\limproj{\mathop{\oalign{lim\cr\hidewidth$\longleftarrow$\hidewidth\cr}}}

\newcommand{\ol}{\overline} \newcommand{\ul}{\underline}

\newcommand{\simlgr}{\buildrel \sim \over \longrightarrow}

\newcommand{\simla}{\buildrel \sim \over \longleftarrow}

\newcommand{\wdh}{\widehat}
\newcommand\ch{^{\scriptscriptstyle\vee}}

\newcommand{\pure}{\mathrm{pure}}
\newcommand\pu{_{\pure}}
\newcommand{\fppf}{_\mathrm{fppf}}
\newcommand{\et}{_\mathrm{\acute{e}t}}
\newcommand{\rmet}{\mathrm{\acute{e}t}} 

\newcommand{\longto}{\longrightarrow}
\newcommand{\we}{\wedge}
    \newcommand{\bwe}{\bigwedge}

\newcommand{\ti}{^\times}
\newcommand{\lan}{\langle}\newcommand{\ran}{\rangle}
\def\co{\colon}
\def\ot{\otimes} 
\def\op{^{\rm op}}
\def\ideal{\triangleleft}
\newcommand{\me}{^{-1}}
\def\dar[#1]{\ar@<2pt>[#1]\ar@<-2pt>[#1]}
\def\wtl{\widetilde}
 
\newcommand{\alg}{_\mathrm{alg}}

\newcommand{\gq}{\mathbin{/\mkern-5mu/}} 

\stackMath
\newcommand\reallywidehat[1]{%
\savestack{\tmpbox}{\stretchto{%
  \scaleto{%
    \scalerel*[\widthof{\ensuremath{#1}}]{\kern.1pt\mathchar"0362\kern.1pt}%
    {\rule{0ex}{\textheight}}
  }{\textheight}%
}{2.4ex}}%
\stackon[-6.9pt]{#1}{\tmpbox}%
}
\newcommand{\cE}{\reallywidehat}

\newcommand{\rmA}{\mathrm{A}}

\newcommand{\Aut}{\operatorname{Aut}}
\newcommand{\ann}{\operatorname{Ann}}

\newcommand{\Br}{\operatorname{Br}} 
\newcommand{\rmB}{\mathrm{B}}

\newcommand{\rmC}{\mathrm{C}}
\newcommand{\Coker}{\operatorname{Coker}}

\newcommand{\Char}{\operatorname{char}} 

\newcommand{\can}{\operatorname{can}}
\newcommand{\Cent}{\mathrm{Cent}} 
\newcommand{\Cont}{\operatorname{Cont}} 

\newcommand{\rmD}{\mathrm{D}}
\newcommand{\diag}{\operatorname{diag}}

\newcommand{\dis}{{\operatorname{dis}}}
\newcommand{\hDisc}{{\operatorname{Disc}}}
\newcommand{\Disc}{\hDisc}
\newcommand{\Di}{{\operatorname{Dick}}}
\newcommand{\Dyn}{\mathrm{Dyn}} 
\newcommand{\Dqd}{\rmD(q)^{[d]}}
\newcommand{\Dqn}{\rmD(q)^{[n]}} 

\newcommand{\rmE}{\mathrm{E}}
\newcommand{\End}{\operatorname{End}}

\newcommand{\rmev}{{\operatorname{ev}}} 

\newcommand{\rmG}{\mathrm{G}}
\newcommand{\GO}{{\operatorname{GO}}}
\newcommand{\GL}{{\operatorname{GL}}}
\newcommand{\Gal}{\operatorname{Gal}}

\newcommand{\Gr}{{\operatorname{Grass}}}

\newcommand{\Hom}{\operatorname{Hom}}

\newcommand{\hyp}{\mathrm{hyp}}

\newcommand{\Id}{\operatorname{Id}}
\newcommand{\Ima}{\operatorname{Im}}

\newcommand{\Int}{{\operatorname{Int}}}
\newcommand{\sint}{{\operatorname{int}}}
\newcommand{\inc}{{\operatorname{inc}}}

\newcommand{\Jac}{{\operatorname{Jac}}}

\newcommand{\Ker}{\operatorname{Ker}}
\newcommand{\rmK}{\mathrm{K}}
\newcommand{\Kneb}{\operatorname{Kneb}}
\newcommand{\rmKE}{\operatorname{KE}}

\newcommand{\Lie}{\operatorname{Lie}}

\newcommand{\Mat}{{\operatorname{M}}}

\newcommand{\Mor}{\operatorname{Mor}}
\newcommand{\mult}{\operatorname{mult}} 
\newcommand{\Max}{\operatorname{Max}} 

\newcommand{\rmN}{\operatorname{N}}

\newcommand{\Nrd}{\operatorname{Nrd}} 

\newcommand{\Norm}{\operatorname{Norm}} 

\newcommand{\Orth}{{\operatorname{O}}}\newcommand{\orth}{\Orth}

\newcommand{\rmQ}{\operatorname{Q}}

\newcommand{\rmP}{\operatorname{P}}

\newcommand{\Prim}{\operatorname{Prim}} 
\newcommand{\pr}{\operatorname{pr}}

\newcommand{\Pic}{{\operatorname{Pic}}}

\newcommand{\Pc}{\operatorname{Pc}} 
\newcommand{\Pcrd}{\operatorname{Pcrd}} 
\newcommand{\Par}{{\operatorname{Par}}} 

\newcommand{\rank}{\operatorname{rank}}
\newcommand{\red}{\operatorname{red}}
\newcommand{\rad}{\operatorname{rad}}
\newcommand{\reg}{{{\operatorname{reg}}}} 
\newcommand{\res}{\operatorname{res}}

\newcommand{\Refl}{\operatorname{Refl}} 
\newcommand{\rmsm}{^{\rm sm}} 

\newcommand{\rmS}{\operatorname{S}}
\newcommand{\wS}{\widehat{\rmS}} 
\newcommand{\Spec}{\operatorname{Spec}}
\newcommand{\Sp}{\underline{\operatorname{Spec}}}
\newcommand{\Specmax}{\operatorname{Specmax}}
\newcommand{\SL}{{\operatorname{SL}}}
\newcommand{\SO}{{\operatorname{SO}}}

\newcommand{\Spin}{{\operatorname{Spin}}}
\newcommand{\Stab}{\operatorname{Stab}}

\newcommand{\Span}{\operatorname{Span}}
\newcommand{\Sym}{{\operatorname{S}}}

\newcommand{\SG}{{\operatorname{S\Gamma}}} 
\newcommand{\rmsn}{\operatorname{sn}} 
\newcommand{\spi}{\operatorname{s\pi}}
\newcommand{\SB}{\operatorname{SB}} 
\newcommand{\ssr}{^{\,\operatorname{ss,reg}}} 
\newcommand{\sw}{\operatorname{sw}} 
\newcommand{\uSpr}{\operatorname{Spr}} 
\newcommand{\Sur}{\operatorname{Sur}} 
\newcommand{\SN}{\operatorname{SN}} 

\newcommand{\Trd}{\operatorname{Trd}}
\newcommand{\Top}{\operatorname{\mathrm Top}}

\newcommand{\tr}{\operatorname{tr}}
\newcommand{\Tr}{\operatorname{Tr}}

\newcommand{\tot}{_{\,\operatorname{tot}}}

\newcommand{\W}{\operatorname{W}}
\newcommand{\Wq}{\operatorname{W}_q}

\newcommand{\hWq}{\widehat{\operatorname{W}}_q}

\DeclareMathOperator{\rmZ}{Z} 

\newcommand{\euN}{\EuScript{N}}

\newcommand{\euT}{\EuScript{T}}

    \newcommand{\Cliff}{{\mathcal C}\ell} 
    \newcommand{\Cli}{\Cliff}
\newcommand{\calD}{\operatorname{\mathcal D}}
    \newcommand{\Dis}{{\mathcal Dis}} 

\newcommand{\calO}{\operatorname{\mathcal O}}
\newcommand{\calP}{\operatorname{\mathcal P}}

\newcommand\sfP{{\sf P}}

\newcommand\sfR{{\sf R}}\newcommand\sfr{{\sf r}}

\newcommand{\scA}{\mathscr{A}}

\newcommand{\scC}{\mathscr{C}}

\newcommand{\scL}{\mathscr{L}}
\newcommand{\scM}{\mathscr{M}}
\newcommand{\scN}{\mathscr{N}}
\newcommand{\scO}{\mathscr{O}}
\newcommand{\scP}{\mathscr{P}}
\newcommand{\scQ}{\mathscr{Q}}

\newcommand{\ulcent}{\underline{\rm Cent}} 
\newcommand{\ulCent}{\ulcent}

\newcommand{\ulGr}{\underline{\mathrm{Grass}}} 

\newcommand{\uIsom}{\underline{\rm Isom}}

\newcommand{\ulL}{\underline{\rm L}}

\newcommand{\Of}{\mathrm{Of}} 
\newcommand{\ulO}{\underline{\rm O}}

\newcommand{\ulSO}{\underline{\rm SO}}
\newcommand{\ulSpr}{\underline{\mathrm{Spr}}} 

\newcommand{\ulU}{\underline{\rm U}} 

\newcommand{\ulV}{\underline{\rm V}}

\newcommand{\ulw}{\underline{\rm W}}
\newcommand{\ulW}{\ulw}

\newcommand{\uAut}{\mathbf {Aut}}

\newcommand{\uB}{\mathbf B}

\newcommand{\bfb}{\mathbf b}

\newcommand{\uCent}{\mathbf{Cent}}

\def\ue{\mathbf e} 

\newcommand{\uGL}{\mathbf{GL}}

\newcommand{\uGO}{\mathbf{GO}}
\newcommand{\uGSO}{\mathbf{GSO}}
\newcommand{\uGr}{\mathbf{Grass}} 

\newcommand{\uKneb}{\operatorname{Kneb}}

\newcommand{\uL}{\mathbf L}

\newcommand{\uO}{\mathbf {O}}

\newcommand{\uP}{\mathbf{P}}

\newcommand{\uPGL}{\mathbf {PGL}}

\newcommand{\uPGO}{\mathbf{PGO}}
\newcommand{\uPrim}{\mathbf{Prim}} 

\newcommand{\uQ}{\mathbf{Q}}
\newcommand{\bfq}{\mathbf{q}} 

\newcommand{\uS}{\mathbf S}
\newcommand{\bS}{\uS}

\newcommand{\uSL}{\mathbf {SL}}
\newcommand{\uSO}{\mathbf {SO}}%
\newcommand{\uSpin}{\mathbf {Spin}}

\newcommand{\uSG}{{\uS{\boldsymbol{\Gamma}}}} 
\newcommand{\uSB}{\mathbf{SB}} 

\newcommand{\uU}{\mathbf U}

\newcommand{\uV}{\mathbf V}

\newcommand{\uW}{\mathbf W}

\newcommand{\uX}{\mathbf X}

\newcommand{\uY}{\mathbf Y}\newcommand{\bfY}{\uY}

\newcommand{\uZ}{\mathbf Z}\newcommand{\bfZ}{\uZ}

\newcommand{\fra}{\mathfrak a} 

\newcommand{\frb}{\mathfrak b}

\newcommand{\m}{\mathfrak m} 
    \newcommand{\gm}{\m}

\newcommand{\gn}{\mathfrak n}
\newcommand{\frn}{\gn}

\newcommand{\p}{\mathfrak p}

\newcommand{\gp}{\mathfrak p}

\newcommand{\frq}{\mathfrak q}

\newcommand{\frR}{\mathfrak R}

\newcommand{\frS}{\mathfrak S}

\newcommand{\frt}{\mathfrak t}

\newcommand{\Aalg}{A\mathchar45\mathbf{alg}}

\newcommand{\Falg}{F\mathchar45\mathbf{alg}}
\newcommand{\kalg}{k\mathchar45\mathbf{alg}}

\newcommand{\Ralg}{R\mathchar45\mathbf{alg}}
\newcommand{\Rnalg}{R_0\mathchar45\mathbf{alg}}
\newcommand{\Rpalg}{R'\mathchar45\mathbf{alg}}
\newcommand{\Salg}{S\mathchar45\mathbf{alg}}

\newcommand{\ZZalg}{\ZZ\mathchar45\mathbf{alg}}

\newcommand{\hm}{\mathfrak{hm}}

\newcommand{\frmod}{\mathfrak{Mod}}

\newcommand{\qm}{\mathfrak{qm}}
    
    \newcommand{\qs}{\mathfrak{qs}}

\newcommand{\Sch}{\mathfrak{Sch}} 

\newcommand{\bbA}{{\mathbb A}}
\newcommand{\CC}{{\mathbb C}}
\newcommand{\FF}{{\mathbb F}}
\newcommand{\GG}{{\mathbb G}}
\newcommand{\HH}{{\mathbb H}}
\newcommand{\MM}{{\mathbb M}}
\newcommand{\NN}{{\mathbb N}}

\newcommand{\PP}{{\mathbb P}}

\newcommand{\RR}{{\mathbb R}}
\newcommand{\ZZ}{{\mathbb Z}}

\newcommand\al{\alpha}
\newcommand\be{\beta}
\newcommand\ga{\gamma} 
\newcommand\Ga{\Gamma}

\newcommand\de{\delta}\newcommand\del{\delta}
\newcommand\De{\Delta}
\newcommand\eps{\epsilon}
\newcommand\veps{\varepsilon} 
\newcommand\io{\iota}
\newcommand\ka{\kappa}
\newcommand\la{\lambda} \newcommand\La{\Lambda}
 \newcommand\vphi{\varphi}

\newcommand\si{\sigma}

\newcommand\ze{\zeta}

\newcommand{\bmu}{\boldsymbol{\mu}}

\begin{document}

\title[Quadratic Spaces over semilocal Rings]{Quadratic Spaces and  Orthogonal Groups over semilocal Rings}
\author[P. Gille]{Philippe Gille}
\address{
Institut Camille Jordan - Universit\'e Claude Bernard Lyon 1, 43 boulevard du 11 novembre 1918, 69622 Villeurbanne cedex - France and} 
\address{Institute of Mathematics "Simion Stoilow" of the Romanian Academy, 21 Calea Grivitei Street, 010702 Bucharest, Romania} 
\thanks{The first author was supported by the project ``Group schemes, root systems, and related representations'', funded by the European Union - NextGenerationEU through Romania’s National Recovery and Resilience Plan (PNRR), call no.~PNRR-III-C9-2023-I8, Project CF159/31.07.2023, and coordinated by the Ministry of Research, Innovation and Digitalization (MCID) of Romania.}
\email{gille@math.univ-lyon1.fr}

\author[E. Neher \heute]{Erhard Neher}
\address{ Department of Mathematics and Statistics, University of Ottawa,
Ottawa, Ontario, Canada, K1N 6N5}
\thanks{The second author was partially supported by NSERC (Canada) through a Discovery grant during the early stages of the project}
\email{neher@uottawa.ca}

\date{\today}
\maketitle
\noindent{\bf Abstract:}
We prove Springer's Odd Degree Theorem for quadratic forms over LG rings, and Scharlau's and Knebusch's norm principles for quadratic forms over semilocal rings. We present applications to the flat cohomology of spin groups and \'etale norm groups. 
\medskip 

\medskip
\setcounter{tocdepth}{1}

\tableofcontents

\addcontentsline{toc}{section}{Introduction}

\medskip

\[ \textbf{\large Introduction} \]

The algebraic theory of quadratic forms, initiated  
by Ernst Witt's foundational 1937 paper, has evolved into a vast and intricate domain intersecting algebra, geometry, and topology. Traditionally, 
this theory was developed over fields, where the structural behavior of quadratic spaces and their associated orthogonal groups admits a high degree of tractability. Over fields, the classification of forms, the structure of orthogonal groups, and the behavior of quadrics under field extensions follow well-established paradigms. However, the progression of mathematics necessitates extending these concepts to more general bases, replacing fields by arbitrary commutative rings. The transition from fields to rings introduces profound algebraic and geometric challenges. Properties that are trivial or vacuous over fields -—- such as the distinction between regular and nonsingular forms, or the non-freeness of projective modules --— become important challenges in the broader context.

These lecture notes provide a rigorous and detailed exposition of the theory of quadratic forms and orthogonal group schemes over general base rings, with a particular emphasis on semilocal rings and, more generally, rings satisfying the Local-Global (LG) condition. A ring $R$ is an LG ring if every polynomial which represents a unit locally at each maximal ideal also represents a unit globally. The LG condition, a property abstracting key features of semilocal rings, proves to be a powerful framework for extending classical theorems from base fields to base rings. This condition guarantees, for instance, that 
quadratic spaces have Witt cancellation, \ref{canqf}, and 
allow a decomposition into rank $2$ subspaces, \ref{nqf-LG}, thus facilitating inductive arguments that would otherwise fail outside the realm of fields.

The central narrative of these notes weaves together the algebraic structure of the orthogonal group schemes $\uO(q)$ and its special orthogonal subgroup $\uSO(q)$, the geometric action of these groups on the associated spheres and quadrics, and the deep constraints imposed by the norm principles of Knebusch and Scharlau under \'etale extensions. By systematically exploring these areas, we aim to bridge the gap between classical quadratic form theory over fields and the modern schematic approach, highlighting both the direct generalizations and the subtle divergences that arise over rings. We now describe the main advances of these notes in more detail.  \sm

{\em Orthogonal Groups and their Actions on Spheres and Quadrics}, \S\ref{sec:sphere}--\S\ref{sec:quadrics}. 
A fundamental problem in the study of quadratic spaces over rings concerns the structure and behavior of the orthogonal group scheme $\uO(q)$ and its subgroups. Over a field $F$, the generation of the orthogonal group $\orth(q) = \uO(q)(F)$ by reflections (the Cartan-Dieudonn\'e-Kneser Theorem~\ref{CDK}) and the transitivity of its action on spheres of given length, a consequence  of Witt's Extension Theorem,  are cornerstones of the theory. When the base ring is local or semilocal, these classical theorems require careful reformulation and the introduction of new geometric tools.

We begin by proving surjectivity of the canonical map $\uSO(q)(R) \to \uSO(q)(R/\fra)$ for a nonsingular $q$ and any ideal $\fra$ of the semilocal ring $R$ contained in the Jacobson radical of $R$, \ref{knex-c}. The proof relies on the smoothness of the group scheme $\uSO(q)$, demonstrating the immediate utility of the language of group schemes. An immediate consequence is the surjectivity of $\uO(q)(R) \to \uO(q)(R/\fra)$ for local rings $(R,\gm)$ when the quadratic form $q$ is regular or when $\fra= \gm$. This surjectivity is not merely a technical lemma; it is a vital lifting property that allows one to deduce global structural features of $\uO(q)$ and $\uSO(q)$ from their behavior over residue fields. 
However, this lifting fails when $q$ is not regular and $\mathfrak{a} \subsetneq \mathfrak{m}$, highlighting the subtle pathologies present over rings of characteristic 2 where the difference between regular and nonsingular forms becomes important.


Moving beyond the groups $\uO(q)$ and $\uSO(q)$, we deeply analyze their actions on spheres and quadrics. Let $\uS\rmsm_{q,a}$ denote the scheme of smooth points of the sphere $\uS_{q,a}$ with $R$--points $\{x \in M : q(x) = a\}$. Extending Witt's transitivity in the field case, we prove that for a quadratic space $(M, q)$ of rank at least 3 over a semilocal ring $R$, the group $\uSO(q)(R)$ acts transitively on $\uS\rmsm_{q,a}(R)$ for any unit $a \in R\ti$, \ref{prop_quadric2}. This transitivity is achieved by first verifying it over the residue fields (using Witt's theorem) and then lifting the isometries via the aforementioned surjectivity properties. This allows us to prove that, under mild conditions on the residue fields, $\uSO(q)(R)$ can be used to construct regular planes and split off hyperbolic direct summands, mimicking the classical orthogonal decomposition process, \ref{actsoLG} and \ref{plar}.

Furthermore, we explore the action of $\uO(q)$ and $\uSO(q)$ on higher-rank quadrics, i.e., the scheme $\uQ_\nu(q)$ parameterizing totally isotropic, complemented submodules of rank $\nu$. These higher-rank quadrics are shown to be smooth projective schemes exactly when the quadratic form is nonsingular, \ref{quadsm}. A central geometric result of these notes is the identification of the stabilizers of this action. We demonstrate that the stabilizer of an $R$-point in $\uQ_\nu(q)$ is a parabolic subgroup $P$ of the reductive group scheme $\uSO(q)$. Moreover, when the rank $\nu$ is not maximal, i.e., $\uQ_\nu(q)$ is not a Lagrangian, the quadric $\uQ_\nu(q)$ is isomorphic to the quotient space $\uSO(q)/P$, \ref{prop_hig}. This provides a complete translation between the geometric objects (quadrics) and the structural components of the algebraic group (parabolic subgroups).

The case of maximal totally isotropic submodules --- the Lagrangian quadric $\uL(q)$ --- requires a more nuanced treatment. Here, 
the $\uO(q)$ action is transitive, but no longer that of $\uSO(q)$. To understand the geometry of $L(q)$, we explicitly construct its Stein factorization, a result of Deligne for quadratic forms over base schemes. We detail Deligne's result showing that the structural morphism 
$\uL(q) \to \Spec(R)$ factors through the spectrum of the discriminant algebra $\Dis(q)$, a quadratic \'etale $R$-algebra. The morphism $\uL(q) \to \Spec\big(\Dis(q)\big)$ is smooth, projective, and has geometrically connected fibers, \ref{deligne}. This factorization beautifully captures the connected components of the Lagrangian quadric and explicitly links the geometry of maximal isotropic subspaces to the classical discriminant invariant of the quadratic form.
\sm 

{\em Springer's Theorem for LG Rings}, \S\ref{sec:springer}.  One of the most celebrated theorems in the classical theory of quadratic forms is Springer's Odd Degree Theorem (1952). It asserts that if a quadratic form over a field $F$ (of characteristic $\neq 2$) becomes isotropic over a finite field extension $K/F$ of odd degree, then the form must already be isotropic over the base field $F$. This result imposes severe constraints on the anisotropic behavior of quadratic forms under field extensions and is a fundamental tool in the algebraic theory of quadratic forms.

A major objective of these notes is to extend Springer's Theorem to the far broader context of LG rings: 
Let $(M, q)$ be a quadratic space over an LG ring $R$, and let $S$ be a finite, locally free $R$-algebra of constant odd rank. Under suitable conditions (for instance, if $S$ is generated by a single element, or if $S$ is \'etale and $R$ satisfies certain primitive criteria), if the extended quadratic space $(M, q)_S$ is isotropic over $S$, then the original space $(M, q)$ is already isotropic over $R$, \ref{spo}, \ref{spo-LGG}.

Our proof of this generalization is highly non-trivial and departs significantly from the classical field-theoretic arguments.
Instead, we employ a geometric approach centered on the construction and representability of what we term ``Springer functors'', \ref{desf}. For a one-generated odd-degree extension $S = R[x]$, the Springer functor parameterizes specific isotropic vectors in the extended space $(M, q)_S$ that satisfy certain conditions relative to the minimal polynomial of $x$. We prove that this functor is representable by a smooth, quasi-compact open subscheme of an affine space, \ref{respr}.

By establishing that this representing scheme has geometrically integral fibers and is non-empty over LG rings (provided the rank of the space is at least 3), we guarantee the existence of rational points. These rational points correspond to specific isotropic vectors that allow us to construct an intermediate algebra of strictly lower odd degree over which the form is also isotropic. The theorem then follows by an induction on the degree of the extension, ultimately reducing the problem to degree one, which is the base ring itself.

This geometric proof is powerful because it avoids the need for $2$ to be invertible and relies intrinsically on the structural properties of reductive group schemes and their homogeneous spaces. Furthermore, we demonstrate that this generalized Springer theorem implies several strong cancellation and restriction properties for quadratic forms and the Witt group of LG rings, \ref{spo}\eqref{spo-b}. We also provide an analogous result for Hermitian spaces over quadratic \'etale extensions, showing that isotropy over an odd-degree extension descends to the base ring, thereby illustrating the robust nature of this geometric method, \ref{sprher}.
\sm 

{\em  Norm Principles of Scharlau and Knebusch}, \S\ref{sec:scharlau} and \S\ref{sec:kneb}. Another major theme explored in these notes concerns the behavior of similarity factors and values of quadratic forms under finite \'etale extensions, formalized in the norm principles of Scharlau and Knebusch. These principles quantify the extent to which  data defined over an extension can be "pushed down" to the base ring via the norm map.

Scharlau's Norm Principle deals with the similarity of quadratic forms. Let $(M, q)$ be a quadratic space over a semilocal ring $R$, let $S$ be a finite \'etale $R$--algebra, and let $\rmG(q)$ denote the group of similarity factors of $q$ (units $u \in R\ti$ such that $uq \cong q$). The principle asserts that if $u \in S^\times$ is a similarity factor for the extended form $q_S$, then its norm $\rmN_{S/R}(u)$ is a similarity factor for the base form $q$, that is, $\rmN_{S/R}(\rmG(q_S)) \subset \rmG(q)$, \ref{thm_scharlau}. 

To establish Scharlau's principle, we first prove it for one-generated extensions and forms of odd rank, where $\rmG(q_S) = S\ti{}^2$.
The general case for semilocal rings is then achieved through a 
``Noetherian reduction''. We express the semilocal ring $R$ as a direct limit of semilocalizations of finitely generated $\mathbb{Z}$-algebras. By passing to an ind-\'etale extension where the residue fields are infinite, we can find a primitive element for the algebra and apply the one-generated case. The result then descends back to the original ring via limit arguments, demonstrating a quintessential application of algebraic geometry techniques to solve purely algebraic problems over rings.

Knebusch's Norm Principle is intimately related but focuses on the values represented by the quadratic form. Let $\rmD(q)$ be the set of values $q(x)$ for vectors $x \in M$ such that $q(x) \in R\ti$. We define 
$\rmD(q)^{[\rm ev]}$ (respectively $\rmD(q)^{[\rm od]}$) as the union of products of an even (respectively odd) number of values. Knebusch's principle states that for a semilocal $R$ and a finite \'etale extension $S/R$ of constant degree $d$, the norm of the values represented by $q_S$ is contained in $D(q)^{[\rm ev]}$ if $d$ is even, and in $D(q)^{[\rm od]}$ if $d$ is odd, \ref{thm-kneb}. 

The proof of Knebusch's principle is arguably the most technically demanding part of these notes. It relies on the construction of ``Knebusch functors''. which are intricate refinements of the Springer functors used earlier. These functors are represented by smooth quasi-compact schemes, \ref{lem_representability}, 
By proving that these Knebusch schemes are non-empty over semilocal rings (again using geometric properties like schematic denseness and Weil restrictions), \ref{lem_PR}, we can systematically decompose the norm of a represented value into a product of values represented over the base ring. 
\sm 

{\em Applications}, \S\ref{sec:consequences}--\S\ref{sec:normgroups}. 
These norm principles are not merely structural curiosities; they have profound cohomological consequences. For instance, we apply Knebusch's principle to evaluate the spinor norm map  $\rmsn \co \SG(q) \to R\ti$, $x \mapsto \ol x x$ where $\SG(q)$ is the special Clifford group, \eqref{mgs11}. 
We prove that over a semilocal ring $R$, the image of the spinor norm $\rmsn$ coincides  precisely with the group $\rmD(q)^{[\rm ev]}$, \ref{spinim}. This explicit identification is then utilized to study the flat cohomology of the spin group $\uSpin(q)$. Specifically, we prove that for a finite \'etale extension $S/R$ of odd degree, the base change map for the first fppf cohomology sets, 
$H^1\fppf(R, \uSpin(q)) \longto H^1\fppf (S, \uSpin(q_S))$, is injective, \ref{prop_sn}\eqref{prop_sn-b}. This provides a deep cohomological parallel to Springer's theorem 

An equally important application of Knebusch's principle is to \'etale norm groups. Slightly adjusting the Kato-Saito definition of norm groups over fields, we define the \'etale norm group of an $R$--scheme $X$ as the
subgroup $\rmN_X^{\rmet}(R)$ of $R^\times$  generated by the subgroups
$\rmN_{R'/R}\bigl( (R')^\times \bigr)$, where $R'$
varies over the finite \'etale extensions $R'$ of $R$ of positive rank for which $X(R') \not=\emptyset$. For a quadratic space $(M,q)$ over a semilocal ring $R$ with $\rank M \ge 2$ and $X=\uQ_1(q)$, the quadric associated with $q$, we show that $\rmN_X^{\et}(R) =  \rmD(q)^{[\rm ev]}$, \ref{prop_norm_quad}. Again for $R$ semilocal  and $X$ the Severi-Brauer scheme of an Azumaya $R$--algebra $A$ with reduced norm $\Nrd$ we prove that $\rmN_X^{\rmet}(R)= \Nrd(A\ti)$, \ref{prop_norm_SB}.  

These applications underscores the unifying power of the geometric and functorial methods developed throughout these notes.

\newpage

\comments{(2024-12) At first, Section \ref{sec:preliminaries-LG}  will simply be a collection of terms, notions, notation that will be used in the main body of the paper. This section should be written at the end, when we know what we actually need. It will be updated as we go along.}

\section{Preliminaries} \label{sec:preliminaries-LG}

\comments{
List of notation, definitions and facts we need:
\sm 
\begin{enumerate}
\item  $\NN=\{0, 1, 2, \ldots \}$, $\NN_+ = \{ 1, 2, \ldots\}$, $\sqcup$ is disjoint union of sets and schemes

\item $R$, $R\ti$; $X \subset Y$ means that $X$ is a subset of the set $Y$ and  
 $X \subsetneq Y$ means that $X$ is a proper subset of $Y$,
 
\item $M_S = M \ot_R S$ faithfully projective, norm (?), 
     $\rank M \ge 2$ means $\rank_\p M_\p \ge 2 $ for all $\p \in \Spec(R)$.  



\item $\FF_p = \ZZ/p\ZZ$  

\item 
finite projective = finitely generated projective = finite locally free; 

\item rank decomposition of finite projective modules

\item $\uW(M)$ is the scheme representing the $R$--functor $T \mapsto M_T$ of a locally free $R$--module of finite rank, \ref{ag}\eqref{ag-ex}. 
    
\item If $P\in R[X]$ is a polynomial over $R$ we write the quotient ring of $R[X]$ by the ideal $(P)$ generated by $P$ as $R[X]/P$, instead of $R[X]/(P)$. 

\item $A\in \Ralg$ {\em finite\/} means that $A$ is finitely generated as $R$--module, \cite[V, \S1.1, Def.~2]{BAC2}. 
    
\item $\Jac(A)$ is the Jacobson radical of a not necessarily commutative ring; for every ideal $\fra\ideal A$ lying in $\Jac(A)$ we have $\Jac(A/\fra) = \Jac(A) / \fra$, see (4.6) of Lam's book on rings, hence $(A/\fra) \big/ \Jac(A/\fra) \cong A/\Jac(A)$.\end{enumerate}
}

\subsection{Faithfully projective modules ({\cite[1.5]{GN-LG}}) }\label{fapmod} Recall that an $R$--module $M$ is {\em faithful\/} if the structure map $R\to \End_R(M)$, $r\mapsto r\Id_M$, is injective. It is a standard fact in commutative algebra,  see for example \cite[IX, Prop.~4.6, page 476]{Bas2},  
that the following conditions are equivalent for a finite projective $R$--module $P$:
\begin{enumerate}[label={\rm (\roman*)}]
\item $P$ is faithful;

\item every localization $P_\p$, $\p \in \Spec(R)$, is non-zero;


\item $P_{R/\m} \ne \{0\}$ for every maximal ideal $\m \ideal R$;


\item $P$ is faithfully flat;

\item there exists an $R$--module $Q$ such that $P\ot_R Q \cong R^n$ for some $n\in \NN_+$.
\end{enumerate}
In this case, $P$ is called {\em faithfully projective}.

\subsection{Direct products of base rings.} \label{dpb} Let 
$R = R_0 \times \cdots \times R_n $ be a direct product of rings. We put $\veps_i = (0, \ldots,0, 1_{R_i}, 0, \ldots, 0)$. Associating with an $R$--module $M$ the direct product 
\begin{equation}\label{dpb0}  
   (\veps_0 M, \ldots, \veps_n M) = (M_0,  \ldots, M_n)
\end{equation}   
and with an $R$--linear map $f \co M \to N$ the direct product $(\veps_1 f, \ldots, \veps_n f)$ of $R_i$--linear maps gives rise to an equivalence between the category $\frmod_R$ of $R$--modules and the direct product of the categories $\frmod_{R_i}$ of $R_i$--modules, which we take as an identification: 
\begin{equation} \label{dpb1} 
 \frmod_R = \frmod_{R_0} \times \cdots \times \frmod_{R_n}.
\end{equation}  
The standard algebraic constructions respect the decomposition \eqref{dpb1}, for example the symmetric algebras do so: 
\begin{equation} \label{dpb2}
 \Sym_R(M) = \Sym_{R_0}(M_0) \times \cdots \times \Sym_{R_n}(M_0).
\end{equation} 
Similarly, properties of $R$--modules $M$ are reflected by the corresponding properties of the $R_i$--modules $M_i$. For example, an $R$--module $M$ is finite projective (faithfully projective resp.) if and only if all $M_i$ are finite projective (faithfully projective resp.) $R_i$--modules. Or, $A$ is a unital associative $R$--algebra if and if the $A_i$ are unital associative $R$--algebras.  

Let $S=\Spec(R)$ and $S_i = \Spec(R_i)$. We know 
$ S= S_0 \sqcup \cdots \sqcup S_n, $ 
\ref{ag}\eqref{ag-du}. For a finite projective $R$--module $M$ the associated $S$--scheme $\uW(M)$, see \ref{ag}\eqref{ag-ex}, is the disjoint union of $S_i$--schemes, 
\begin{equation}\label{dpb3}
 \uW(M) = \uW(M_0) \sqcup \cdots \sqcup \uW(M_n),
\end{equation} 
which is immediate from $\uW(M) = \Spec\big( \Sym_R(M\ch)\big)$ and \eqref{dpb2}. 
  
A standard way to obtain ta situation as above occurs by letting $M=M_0\times \cdots \times M_n$ be the {\em rank decomposition\/} of a finite projective $R$--module $M$ for which $M_i$, $0\le i \le n$, is a finite projective $R_i$--module of constant rank $i$. The discussion above then describes the reduction of quadratic modules to quadratic modules of constant rank. We will refer to the process of passing from $M$ to the $M_i$ as {\em reduction to constant rank}.

\subsection{Unimodular vectors (\cite[9.13-9.17]{PRbook}, \cite[0.3]{Lo-genalg})}\label{unimod}
Let $M$ be an $R$--module and let $M^*$ be its dual space. An $x\in M$ is  called  {\em unimodular\/} if the following equivalent conditions are
fulfilled:
\begin{enumerate}[label={\rm (\roman*)}]
  \item\label{unimodi} $R\cdot x$ is a free $R$--module of rank $1$ and a
      direct summand of $M$,

\item\label{unimodii} there exists $\vphi \in M^*$ for which $\vphi(x) = 1$.
\end{enumerate}
If $M$ is finitely generated projective, then \ref{unimodi} and \ref{unimodii} are equivalent to \ref{unimodiii} and \ref{unimodiv} below: 
\begin{enumerate}[label={\rm (\roman*)}]\setcounter{enumi}{2}
\item\label{unimodiii} $x\ot1_{\ka(\p)} \ne 0$ for all $\p\in \Spec(R)$, where
  $\ka(\gp)$ is the residue field of the local ring $R_\gp$,  

\item \label{unimodiv} $x\ot 1_{R/\m} \ne 0$ for every maximal $\m \in \Spec(R)$.
\end{enumerate}
We denote by $M_u$ the set of unimodular vectors of $M$, and list some properties of unimodular vectors. \sm 

\begin{inparaenum}[(a)] \item \label{unimod-ab}       
{\em For $I \subset \Jac(R)$ an ideal, the map
\begin{equation}\label{unimod-2}
  M_u \twoheadrightarrow (M/IM)_u, \quad m \mapsto m\ot 1_{R/I}
\end{equation}
is surjective}. Indeed, this easily follows from surjectivity of $M \twoheadrightarrow M/IM$, bijectivity of the maximal spectra, given by $\Specmax(R) \to \Specmax(R/I)$, $\m \mapsto \m/I$, and condition \ref{unimodiv}. Of course, this is also a special case of the general criterion~\ref{opemaLG}, as explained in \ref{opema-rem}. \sm 

\item\label{unimod-dp} ({\em Direct products of base rings}) Let $R=R_0 \times \cdots \times R_n$ be a direct product of rings, and let $M=M_0 \times \cdots \times M_n $ be an $R$--module as in \ref{dpb}. Then $m=(m_0, \ldots, m_n)$ with $m_i \in M_i$ is unimodular in the $R$--module $M$ if and only if $m_i$ is unimodular in the $R_i$--module $M_i$ for $i=0, \ldots, n$. 
\sm 
    
\item ({\em Unimodular vectors in free modules}) Suppose $M$ is a free $R$--module of finite rank, say with basis $(e_1, \ldots, e_n)$. Then a vector $x=r_1 e_1 + \cdots + r_ne_n$ is unimodular if and only if the ideal $(r_1, \ldots r_n)$ generated by the coefficients $r_i$  equals $R$. 
\lv{
  Indeed, if this ideal is a proper ideal, there exists a maximal ideal $\m$ containing it. It follows that all coefficients $r_i \ot 1_{R/\m} = 0$, i.e., $x(\m) = 0$. Thus, $(r_1, \ldots, r_n) = R$ whenever $x$ is unimodular. Conversely, suppose that the ideal is $R$. We then have a relation $1_R = \sum_i s_i r_i$ for suitable $s_i \in R$. For every maximal ideal $\gm \in R$ we then $0 \ne 1_{R/\m} = \sum_i s_i(\m)r_i(\m)$ and hence $0 \ne x(\m) = \sum_i r_i(\m) e_i$ (else all $r_i(\m)=0$), contradiction.}
\sm 

\new
\item \label{unimod-d} ({\em Example liner forms}) Let $M$ be a reflexive, e.g., finite projective $R$--module. A linear form $\la \co M \to R$ is surjective if and only if $\la$ is a unimodular element of $M^*$. 
    
    Indeed, let $c_M \co M \to M^{**}$ be the canonical isomorphism. Then $\la$ is surjective if and only if there exists $m\in M$ such that $\la(m) = c_M(m)(\la) =1$. By \ref{unimodii}, the latter condition is equivalent to $\la$ being  unimodular in $M^*$.   
\enew
\end{inparaenum}

\begin{lem}[The scheme $\uW(M) \setminus \{0\} = \uW(M)_u$,  {\cite[0.12]{Lo-genalg}}] \label{spreq-a} Let $M$ be a finite locally free $R$--module. The $R$--functor $T \mapsto (M\ot_R T)_u$ is represented by the open subscheme $\uW(M) \setminus \{0\}$ of the affine $R$--scheme $\uW(M)$, often abbreviated as $\uW(M)_u$. It has the following properties. 
\sm 

\begin{enumerate}[label={\rm (\roman*)}]
  \item \label{spreq-ai} $\uW(M)_u$ is a quasi-compact, equivalently, a finitely presented open subscheme of $\uW(M)$. \sm 

\item \label{spreq-aii} Forming $\uW(\cdot)_u$ commutes with base change, i.e., for $T\in \Ralg$ we have
 $\uW(M)_u \times_R T \cong \uW(M\ot_R T)_u$. \sm 
 
 \item \label{spreq-aiii} If $M$ has constant rank $1$, then $\uW(M)_u$ is an affine scheme.  \sm      
 
\item \label{spreq-b} Recall that the projective space $\uP(M^*)$ represents the $R$--functor assigning to $T\in \Ralg$ the set of complemented projective submodules in $M_T$ of constant rank $1$, {\rm \ref{grap}\eqref{grape-dual}}. Assigning to a unimodular $x\in M_T$ the complemented line bundle $Tx$, gives rise to a morphism of $R$--schemes 
     \[ p \co \uW(M)_u \to \uP(M^*) \]
     which is a $\GG_m$--torsor for the Zariski topology.   
\sm 

\new
\item \label{spreg-v} $\uW(M^*)_u$ represents the $R$--functor whose $T$--points, $T\in \Ralg$, are $\{ \la \in M_T^* ; \text{ $\la$ is surjective}\}$. 
\enew 
\end{enumerate} 
\end{lem}

\begin{proof} \ref{spreq-ai} Let $\la_1, \ldots, \la_n$ be a spanning set of  the dual space $M\ch$ of $M$. By \ref{unimod}\ref{unimodii}, a vector $m\in M_T$, $T\in \Ralg$, is unimodular if and only if the ideal generated by $\la_1\ot 1_T, \ldots, \la_n\ot 1_T$ is all of $T$. Hence, the $R$--functor $T \mapsto (M_T)_u$ is represented by the union $\uW(M)_u$ of the finitely many principal open subschemes $\uW(M)_{\la_i}$, $i = 1, \ldots, n$. It follows that $\uW(M)_u$ is a quasi-compact scheme and, consequently, the open immersion $\uW(M)_u \to \uW(M)$ is quasi-compact, equivalently, is of finite presentation \cite[Tag 01TU]{St}. The composition $\uW(M)_u \to \uW(M) \to \Spec(R)$ is therefore also of finite presentation. \sm 

\ref{spreq-aii} Using the notation of \ref{spreq-ai} we have
\begin{align*}
  \uW(M)_u \times_R T &\cong \big( \uW(M)_{\la_1} \cup \cdots \cup \uW(M)_{\la_n}\big) \times_R T 
  \\& \cong ( \uW(M)_{\la_1} \times_R T ) \cup \cdots 
          ( \uW(M)_{\la_n} \times_R T )
  \\&\qquad  \qquad(\text{inside $\uW(M)\times_R T$, by \cite[Tag 01JS]{St}})
\\& \cong (\uW(M)\times_R T)_{\la_1 \ot 1_T} \cup \cdots \cup 
           (\uW(M)\times_R T)_{\la_n \ot 1_T}
\\& \qquad \qquad (\text{by \eqref{agbas1}})
\\& \cong \uW(M\ot_R T)_{\la_1 \ot 1_T} \cup \cdots \cup 
           \uW(M\ot_R T)_{\la_n \ot 1_T} 
\\& = \uW(M\ot_R T)_u.
    \end{align*}

\ref{spreq-aiii} 
This is  clear if $M$ is free of rank $1$, since then $\uW(M)_u\cong \GG_m$. In general, we can choose a standard Zariski cover $\{ U_i \to \Spec(R)\}_{i=1, \ldots, \ell}$,  cf.~\ref{Zarev}, and such that $M\ot_R R_i$ is free of rank $1$. Denoting by $f \co \uW(M)_u \to \Spec(R)$ the structure map, the inverse image $f\me\big( \Spec(R_i)\big)$ is the fibre product
\begin{equation*}  \vcenter{ 
\xymatrix@C=50pt{\ar @{} [dr] |{\small\qed} %
f\me\big( \Spec(R_i)\big)  \ar[r] \ar[d] & \Spec(R_i) \ar[d] 
 \\ \uW(M)_u \ar[r]^f & \Spec(R)} }  \quad.  \end{equation*}
Hence, by \ref{spreq-aii}, the scheme $f\me\big( \Spec(R_i)\big) = \uW(M)_u \times_R R_i \cong \uW(M\ot_R R_i)_u$ is affine. Thus, by \cite[Tag 01S8]{St}, the morphism $f$ is affine. Therefore $f\me\big( \Spec(R)\big) = \uW(M)_u$ is an affine scheme. \sm 

\ref{spreq-b} is straightforward by Zariski localization. \new \ref{spreg-v} follows from \ref{unimod}\eqref{unimod-d}.\enew  \end{proof}

\sm 

\textbf{Remark.} As a quasi-compact open subscheme of the affine scheme $\uW(M)$, the scheme $\uW(M)_u$ is quasi-affine \cite[Tag 01P6]{St}, but it is in general not affine. Indeed, let $k$ be a field. Then $\uW(k^2)_u = \bbA_k^2\setminus \{(0,0)\}$, which is well-known to be a non-affine scheme \cite[Tag 0IL]{St}.   

\comments{
Further results (from our Springer-paper): \item\label{isotrop-suff} (a) Let $R[X]$ the polynomial ring over $R$ in the variable $X$ and let $(M,q)$ be a quadratic module. For $v=v(X) \in M \ot_R R[X]$ define the affine $R$--scheme
    \[ Z_v = \{x\in \GG_{a,R} : v(x) = 0 \},  \]
    whose $T$--points, $T\in \Ralg$, is the set $Z_v(T) = \{ t\in T: v(t)=0\}$ where $v(t) \in M \ot_R T$ is obtained by substituting $t$ for $X$.  Then
    \begin{equation} \label{isotrop-suff1}  \text{\em $Z_v$ empty}\quad \implies \quad \text{$v$ unimodular.}
     \end{equation}
     Indeed, if $\p \in \Spec(R[X])$, then $v(\p) = v(X \ot 1_{\ka(\p)}) \ne 0$ since otherwise $X\ot 1_{\ka(\p)} \in Z_v(\ka(\p))$.
}

\subsection{Semilocal rings}\label{slr} Recall that  a unital commutative ring $R$  is {\em semilocal\/} if it has only a finite number of maximal ideals, equivalently, $R/\Jac(R)$ is a finite direct product of fields, see e.g. \cite[II, \S3.5]{BAC}.

We list several facts used later, but for which we could not find a convenient reference. Throughout, $R$ is a semilocal ring with maximal ideals $\m_i$, $i=1, \ldots, n$,  and residue fields $\ka_i = R /\m_i$. Thus $R/\Jac(R) = \ka_1 \times \cdots \times \ka_n$. \sm

\begin{inparaenum}[(a)]
   \item\label{slr-a} A finite direct product of semilocal rings is semilocal if and only if every factor is semilocal. \sm

 \item \label{slr-ai} For any ideal $I \ideal R$ the quotient ring $R/I$ is semilocal. 
In particular, if $I \subset \Jac(R)$, then
     \begin{equation}  \label{slr-ai1}
     (R/I)\big/ \Jac(R/I) \cong \textstyle \prod_i \ka_i
     \end{equation}
     since $\Specmax(R) \to \Specmax(R/I), \m \to \m/I$ is a bijection.
   \sm

\inparcom{(2025-07-15) added \eqref{slr-b}, it may be known, we need this}

\item\label{slr-b} If $A\in \Ralg$ is a finite $R$--algebra, then $A$ is semilocal, the canonical map $A\ti \to \prod_i (A\ot_R \ka_i)\ti$ is onto, and $\Jac(R) A \subset \Jac(A)$. 

{\em Proof.} That $A$ is semilocal follows for example from \cite[V, (1.1.1)]{K}
(it is proven in \cite[10.1.1]{Ford} for local rings).  Denoting by $\Jac(A)$ and $\Jac(R)$ the Jacobson radicals of $A$ and $R$ respectively, we have $\Jac(R) A \subset \Jac(A)$ by \cite[II, (4.2.4)]{K},
and hence $A\ti \to (A/\Jac(R) A)\ti$ is onto by \ref{nak}\eqref{nak-a}. 
But $A/\Jac(R) A \cong A\ot_R (R/\Jac(R)) \cong \prod_i A \ot_R \ka_i$.
\end{inparaenum}

\subsection{Unimodular rings}\label{unifap} In \S\ref{sec:quadratic-forms} we will encounter the following condition on a (commutative) ring $R$. We will say that a ring $R$ is {\em unimodular\/} if every faithfully projective $R$--module  contains a unimodular vector. \sm 

\begin{inparaenum}[(a)]  \item\label{unifap-a} ({\em Direct products}) If $R=R_1 \times R_2$ is a direct product of rings, then $R$ is unimodular if and only if $R_1$ and $R_2$ are unimodular. 
Indeed, this follows easily from the characterizations of faithfully projective modules and unimodular vectors in \ref{fapmod} and \ref{unimod}\eqref{unimod-dp} respectively. \sm 

\item\label{unifap-b} The following is an easy induction:
    \begin{equation} \label{unimod-b1} \begin{split}
    &\text{\em If $R$ is a unimodular ring and $M$ is finite projective $R$--module}  \\ &\qquad \text{\em of constant rank, then $M$  is free.}
   \end{split} \end{equation} 

\inparcom{(2025-08-15) The condition \ref{unimod-b1}  is not out of the world! The following is \cite[Cor~22.17]{PRbook}: \tt  Suppose every projective k-module of constant finite rank is free. Then the automorphism group of a composition algebra C over k acts transitively on the elementary idempotents of C.\sm 

\cite[Cor.~22.18]{PRbook} Suppose every projective k-module of constant finite rank is free. For a composition algebra C of constant rank > 1 over k, the following conditions are equivalent: (i) C is split.
(ii) The norm of C is isotropic. }

Moreover, using the rank decomposition of finite projective modules as well as \eqref{unifap-a} and  \eqref{unimod-b1} allows us to conclude: 
  \begin{equation} \label{unimod-b2} \begin{split}
    &\text{\em Two finite projective modules are isomorphic if and only if}\\
    &\quad \text{\em they have the same rank function $\Spec(R) \to \NN$.}
    \end{split}\end{equation}

\item The following are examples of unimodular rings: \end{inparaenum}
\enew
\begin{enumerate}[label={\rm (\roman*)}]
  \item LG rings, in view of \eqref{unimod-b0}, in particular $R=0$ is unimodular.  
  
  \item Rings over which every finite projective $R$--module is free, like
  PIDs (\cite[VII, \S3.1, Cor.~3]{BA5}), polynomial rings $k[X_1, \ldots, X_n]$ for  $k$ a field, and Laurent polynomial rings $[X_1^{\pm 1}, \ldots, X_n^{\pm 1}]$ (Serre's question). In particular, $\ZZ$ is unimodular, but not LG by \ref{revLG}\eqref{revLG-non}.
\end{enumerate}

\subsection{Review of LG rings {\cite{EG}, \cite{PRbook}, \cite{GN-LG}}} \label{revLG} Given $S\in \Ralg$, one says that a polynomial $g\in S[X_1, \ldots, X_n]$ represents a unit over $S$ if there exists $s_1, \ldots, s_n \in S$ such that $g(s_1, \ldots, s_n) \in S\ti$. A ring $R$ is an {\em LG ring\/} if for every $n\in \NN_+$ and every $f\in R[X_1, \ldots, X_n]$ the polynomial $f$ represents a unit over $R$ if and only if for every maximal $\gm \ideal R$ the induced polynomial $f_\gm \in R_\gm[X_1, \ldots, X_n]$ represents unit over $R_\gm$, equivalently, $f_{R/\gm}$ represents a unit over any residue field $R/\gm$ of $R$. 
Following are some facts that we will use. 
\sm 

\begin{inparaenum}[(a)] 

\item\label{revLG-a} (\cite[Thm.~2.10]{EG}, \cite[1.6]{GN-LG}, \cite[p.~457]{MW}) Suppose that $R$ is an LG ring. The following is crucial,  implying that $R$ is a unimodular ring:
\begin{equation}
  \label{unimod-b0} \text{\em A faithfully projective $R$--module contains a unimodular vector.}
\end{equation}

\item \label{revLG-ideals} ({\em Quotients}) Let $\fra$ be an ideal of an LG ring $R$. It is immediate from the definition that then $R/\fra$ is an LG ring too. Conversely, if $\fra \subset \Jac(R)$, the Jacobson radical of $R$, then $R$ is LG if and only if $R/\fra$ is LG. This follows from the characterization of invertibility and $\fra \subset \gm$ for every maximal ideal of $R$. 

\item \label{revLG-aa} ({\em Direct products}) $R=R_1 \times R_2$ is LG $\iff$ both $R_1$ and $R_2$ are LG, \cite[11.21]{PRbook}. 
\sm 

\item \label{revLG-int} ({\em Integral extensions}) If $R$ is an LG ring, then so is every integral extension of $R$ \cite[Cor.~2.3]{EG}. In particular, every finite $R$--algebra is an LG ring. \sm 

\item\label{revLG-ex} ({\em Examples}) Every  semilocal ring, for example a field, is an LG ring. But not every LG ring is semilocal. For example, $0$--dimensional rings and the ring of algebraic integers are LG rings, but are not semilocal. \sm  

\item \label{revLG-non} ({\em Non-examples}) The ring $\ZZ$ and the polynomial ring $k[X]$ for an integral domain $k$ are not LG rings (\cite[Exc.~11.42]{PRbook}).  \sm
\end{inparaenum}

\subsection{Rings satisfying the primitive criterion}\label{revLGG}
Recall that a polynomial in $R[X_1, \ldots, X_n]$ is {\em primitive\/} if its coefficients generate $R$ as ideal.
One says that a ring $R$ {\em satisfies the primitive criterion\/} \cite{EG, MW}
if the following equivalent conditions hold: 
\begin{enumerate}[label={\rm (\Roman*)}]
\item for every primitive polynomial $P\in R[X]$ there exists $r\in R$ such that $P(r) \in R\ti$;

\item for every primitive $Q\in R[X_1, \ldots, X_n]$ there exists $(r_1, \ldots, r_n)\in R^n$ such that $Q(r_1, \ldots, r_n) \in R\ti$;

\item \label{LG-ex-prii} $R$ is LG and all residue fields of $R$ are infinite.
\end{enumerate}

\sm
An example of a ring satisfying the primitive criterion, is the ring $S\me R[X]$ where $R$ is arbitrary and $S$ is the multiplicative subset of all primitive polynomials in the polynomial ring $R[X]$, \cite[1.13]{vdK}.

\begin{prop}[Characterization of LG rings {\cite[Prop.~1.4]{GN-LG}}] \label{prop_baire} Let $R$ be a LG-ring, let $M$ be a finite  locally free $R$-module, and let $U$ be an open quasi-compact subscheme of $\uW(M)$.
\sm

\begin{enumerate}[label={\rm (\alph*)}] \item \label{prop_baire-a}
$U(R) \not = \emptyset \iff U(R/\gm) \not = \emptyset$ for every maximal ideal $\gm \ideal R$.
\sm

\item  \label{prop_baire-b} If  $R$ satisfies the primitive condition as in  {\rm  \ref{revLGG}} and $U$ is $R$-dense \new (= universally schematically dense), \enew then $U(R) \not = \emptyset$.
\end{enumerate}\end{prop}

Property \ref{prop_baire-a} characterizes LG rings: if for every finite locally free $R$-module $M$ and every open quasi-compact 
$U \subset \uW(M)$ we have the equivalence \ref{prop_baire-a}, then $R$ is an LG ring.

Regarding \ref{prop_baire}\ref{prop_baire-b}: Since $\uW(M)$ is smooth, we can by \ref{tod}\eqref{tod-e} use the criterion \eqref{tod-d1} saying that $R$--dense is equivalent to universally schematically dense. 

\subsection{Characteristic polynomials, traces and norms}\label{trno} Let $M$ be a finite projective $R$--module. One can define the characteristic polynomial, trace and determinant of an endomorphism of $M$ in equivalent ways, either by embedding $M$ into  a finite free $R$--module \cite{Goldman} or by faithfully flat descent \cite[II, \S2.4]{KO}, \cite[5.3.3, 5.3.5]{Ford}. All of these generalize the case of a free $R$--module treated in \cite[III, \S8.11]{BA}. 

In the following, let $A$ be a unital associative $R$--algebra,  which is finite projective as $R$--module and let $a\in A$. The {\em characteristic polynomial $\Pc_{A/R}(a;X)$\/}, the {\em trace\/} $\Tr_{A/R}(a)$ and the {\em norm\/} $\rmN_{A/R}(a)$ of $a$ are defined as the characteristic polynomial, the trace and the determinant of the left multiplication $L_a\co A \to A$, $x \mapsto ax$ of $A$ by $a\in A$. They are related by
\begin{equation}\label{trno-0}
 \Pc_{A/R}(a;X) = \det (X\Id_A - L_a) = \rmN_{A[X]/R[X]}(X - a).
\end{equation} 
If $A$ has constant rank $n$ as $R$--module, then
\begin{equation}\label{trno-00} \begin{split}
  \rmN_{A/R}(a) &= (-1)^n \Pc_{A/R}(a; 0)\qquad \text{and} \\
\Pc_{A/R}(a;X) &= X^n - \Tr_{A/R}(a) X^{n-1} + \cdots + (-1)^n \rmN_{A/R}(a).
\end{split}\end{equation}
The following facts are easily established, see for example \cite[III, \S9.3, 9.4]{BA} in case $A$ is a free $R$--module. 
\sm

\begin{inparaenum}[(a)]
   \item \label{trno-d} ({\em Base change}) Characteristic polynomials, traces and norms respect base change in the obvious sense: $\Tr_{A\ot S/S}(a \ot 1_S) = \Tr_{A/R}(a) \ot 1_S$ and $\rmN_{A\ot S/S}(a \ot 1_S) = \rmN_{A/R}(a) \ot 1_S$ for any $S\in \Ralg$.
       \sm

   \item \label{trno-a} ({\em Multiplicativity}) The norm respects products: for $a_1, a_2 \in A$ we have $\rmN_{A/R}(a_1 a_2)= \rmN_{A/R}(a_1) \, \rmN_{A/R}(a_2)$. An element $a\in A$ is invertible if and only if $\rmN_{A/R}(a) \in R\ti$ \cite[III, \S9.4, Prop.~3]{BA}. \sm

  \item \label{trno-b} ({\em Transitivity}) Let $S\in \Ralg$ and let $A$ be a unital associative $S$--algebra. Assume that $S$ is finite projective as $R$--module and $A$ is finite projective as $S$--module. Then $A$ is finite projective as $R$--module 
      and
      \[ \rmN_{A/R}(a) = \rmN_{S/R}\big( \rmN_{A/S}(a)\big) \]
    holds for $a\in A$. This is for example proven in \cite[III, \S9.4, (26)]{BA} for free modules and follows in general by localization. \sm

  \item \label{trno-c} ({\em Direct products}) If $A=A_1 \times \cdots \times A_n$ is a direct product of $R$--algebras, then $\rmN_{A/R}(a_1, \ldots, a_n) = \prod_i \rmN_{A_i/R}(a_i)$. \sm

\item  \label{trno-e}  If $A$ is projective of constant rank $d$, then $\rmN_{A/R}(r) = r^d$ for every $r\in R$.
\sm

\item\label{trno-bou} Let $p = a_0 +  a_1X + \cdots + a_{d-1}X^{d-1} + X^d$ be a polynomial over $R$, put $S=R[X]/(p)$ and let $x$ be the class of $X$ in $S$. Then $p = \Pc_{S/R}(x;X)$, $\Tr_{S/R}(x) =  -a_{d-1}$ and $\rmN_{S/R}(x) = (-1)^d a_0$. In particular, by \eqref{trno-a}, %
 \begin{equation}\label{trno-bou1}   
    x \in S\ti \iff a_0 \in R\ti.
  \end{equation}%

 \item\label{trno-ds} By \eqref{trno-d}, the norm gives rise to a polynomial in the sense of \cite{Roby},
\begin{equation}\label{trno-ds0}  
   \rmN_{A/R} \co \uW(A) \to \GG_a, \quad \rmN_{A/R}(S) = \rmN_{A\ot S/S}.
\end{equation} 
    By \eqref{trno-a}, the inverse image of $\GG_m$ under $\rmN_{A/R}$ is the principal open subscheme
     \begin{equation}\label{trno-ds1}
            \uGL_1(A) :=  \uW(A)_{\rmN_{A/R}}
     \end{equation}
     representing the $R$--functor $S\mapsto (A\ot_R S)\ti$. It is an affine finitely presented $R$--group scheme \cite[II, \S1, 2.3]{DG} (where  $\uGL_1(A)$ is denoted $\bmu^A$),  
     and is universally schematically dense in $\uW(A)$ since $1_{A\ot_S S} \in (A\ot_R S)\ti$, cf.\ \ref{usd-d}. 
     \lv{
The constructions $A \mapsto \uW(A)$ and $\rmN_{A/R}$ commute with base change. Hence, for the proof of $R$--denseness it is enough to assume that $R$ is a field. We have then to show that the invertible elements of $A$ are Zariski-dense, see for example [Loos 2006, Generically algebraic.., 0.14]
}
If $M$ is a finite projective $R$--module, we abbreviate
\[ \uGL(M) = \uGL_1(\End_R(M)).\]
If $A$ is commutative, the $A$--group scheme $\GG_{m,A}$ should not be confused with the $R$--group scheme $\uGL_1(A)$. The two are related by the Weil restriction
\begin{equation} \label{trno-dsw} 
    \frR_{A/R}(\GG_{m,A}) = \uGL_1(A).
\end{equation} 

By \eqref{trno-a}, the restriction of the polynomial $\rmN_{A/R}$ to the $R$--group scheme $\uGL_1(A)$ gives rise to a homomorphism of $R$--group schemes, 
\begin{equation}\label{trno-ds2} 
  \rmN_{A/R} \co \uGL_1(A) \to \GG_m.
\end{equation}
That we use the same notation $\rmN_{A/R}$ for \eqref{trno-ds0} and \eqref{trno-ds2} will not lead to any confusion, since it will always be clear from the context what is meant. 
\end{inparaenum}

\subsection{Finite \'etale algebras} \label{fea} An $A \in \Ralg$ is an {\em \'etale $R$--algebra} if $A$ is a separable finitely presented $R$--algebra and flat as $A$--module. There are several equivalent conditions defining \'etale algebras, see for example \cite[IV$_4$, \S17, \S18]{EGA}, \cite{Ford}, \cite{Lenstra}, \cite{Ray-hensel}, \cite[Tag 00U0]{St} and \eqref{fea-a}--\eqref{fea-e} below. Below we list some of the facts we use or that may be helpful for understanding \'etale $R$--algebras.%

We call $A$ {\em finite \'etale} if $A$ is a finite and \'etale $R$--algebra, hence necessarily finite projective by \eqref{fea-gi}. In \cite[Tag 00U0]{St}, ``\'etale'' means finite \'etale. An {\em \'etale $R$--algebra $A$ of degree $d\in \NN_+$} is a finite \'etale $R$--algebra whose underlying $R$--module is projective of constant rank $d$. We say that $A$ is {\em faithfully \'etale\/} if $A$ is a faithfully projective and \'etale $R$--algebra, hence in particular finite \'etale.
\sm

\begin{inparaenum}[(a)] \item \label{fea-a} ({\em Base fields}) An \'etale $k$--algebra $E$ over a field $k$ is the same as an \'etale $k$--algebra in the sense of \cite[V, \S6]{BA5}, i.e., $E$ is a finite-dimensional $k$--algebra satisfying $E=K_1 \times \cdots \times K_n$ where $K_i/k$, $i=1, \ldots, n$, are finite separable field extensions, \cite[Tag 00U3]{St}. \sm

\item \label{trno-f} ({\em Finite base fields}) Let $F$ be a finite field and let $E$ be a finite \'etale $F$--algebra. Then $\rmN_{E/F}\co E \to F$ is surjective.
    Indeed, $E= L_1 \times \cdots \times L_n$ is a product of finite fields. Hence  $\rmN_{E/F}(\ell_1, \ldots, \ell_n) = \prod_i \rmN_{L_i/F}(\ell_i)$, $\ell_i \in L_i$ by \ref{trno}\eqref{trno-c}, and each $\rmN_{L_i/F}$ is surjective by \cite[V, \S12.2, Prop.~4]{BA5}.
\sm

\item \label{fea-gi} ({\em Finite \'etale})  A finite $R$--algebra $A$ is finitely presented as $R$--algebra if and only if $A$ is finitely presented as $R$--module \cite[IV$_1$, 1.4.7]{EGA}.  In particular, a finite \'etale $R$--algebra is finite projective as $R$--module.
    A finite projective $A\in \Ralg$ is (finite) \'etale if and only if $A$ is separable.
  \sm

\item\label{fea-fibc} ({\em Fiberwise Criterion}) Let $A\in \Ralg$ be a finitely presented $R$--algebra whose underlying $R$--module is flat. Then the following are equivalent.

    \begin{inparaenum}[(i)]
      \quad \item $A$ is an \'etale $R$--algebra;

      \quad \item for each $\p \in \Spec(R)$ the $\ka(\p)$--algebra $A \ot_R \ka(\p)$ is \'etale;

      \quad \item for each $R$--field $F$, the $F$--algebra $A\ot_R F$ is \'etale.
    \end{inparaenum} 
\sm 

\item \label{fea-e} ({\em Trace criterion}) Let $A\in \Ralg$ be finite projective and faithful. Then $A$ is \'etale if and only if the bilinear form
    $ A \times A \to R$, $(x,y) \mapsto \Tr_{A/S}(xy)$
     is regular in the sense of \ref{bfLG}\eqref{bfLG-ad}. In this case, $A$ is finite \'etale. This follows for example from \cite[V, \S8.3, Prop.~3]{BA5} and \eqref{fea-fibc}.     
\sm

\item \label{fea-bc} ({\em Base change}) If $A$ is an \'etale $R$--algebra, the $S$--algebra $A\ot_R S$ is \'etale for any $S\in \Ralg$; analogously for finite \'etale (\cite[9.2.5]{Ford}, \cite[II, Prop.~2]{Ray-hensel}, \cite[Tag 00U2]{St}).%
     \sm

\item\label{fea-ten} ({\em Tensor products}) If $A$ and $B$ are both \'etale $R$--algebras, then $A\ot_R B$ is an \'etale $R$--algebra (\cite[9.2.5]{Ford}, \cite[II, Prop.~3]{Ray-hensel}); analogously for finite \'etale algebras.
\sm

\item \label{fea-d} ({\em Direct products of base rings}) Let $R=R_1 \times \cdots \times R_n$ be a direct product of rings. Given $A\in \Ralg$, let $A=A_1 \times \cdots \times A_n$ be the corresponding decomposition of $A$ into the direct product of $R_i$--algebras. Then $A$ is an \'etale $R$--algebra if and only if every $A_i$ is an \'etale $R_i$--algebra; analogously for finite \'etale algebras.
\sm

\item\label{fea-trans} ({\em Transitivity}) If $A\in \Ralg$ is an \'etale $R$--algebra and $B\in \Aalg$ is an \'etale $A$--algebra, then the $R$--algebra $\frR_{A/R}(B)$ is \'etale. The analogous criterion holds for finite \'etale $R$--algebras,  \cite[Tag 00U2]{St}). \sm

\item \label{fea-g} ({\em Reduction mod $\Jac(R)$}) Let $A\in \Ralg$ be faithfully projective and let $I \subset \Jac(R)$ be an ideal of $R$. Then $A$ is an \'etale $R$--algebra if and only if $A/I$ is an \'etale $R/I$--algebra. This follows for example from \eqref{fea-e} and  \ref{nak}\eqref{nak-a}.
\sm

\item \label{fea-sp} ({\em Split \'etale}) Let $d\in \NN_+$. The $R$--algebra $R^d$ is \'etale, called the {\em split \'etale $R$--algebra of rank $d$}. The name is justified by the following fact:

   Suppose $A\in \Ralg$ is projective of constant rank $d$. Then $A$ is \'etale if and only if there exist a faithfully flat $R$--algebra $T$ such that $A \ot_R T \cong T^d$ as $T$--algebra, \cite[4.6.11]{Ford} or \cite[II, Prop.~4]{Ray-hensel}. One can even ``split'' $A$ by a Galois extension $T\in \Ralg$ with Galois group $\frS_d$, the symmetric group in $d$ letters. Then 
$A$ is obtained by Galois descent from the split case,
see e.g.\ \cite[Prop.~2.18 and its Remark]{Salt} or \cite[2.5.2.4]{CF}.
\end{inparaenum}

\comments{(2022-06-10) Added freeness criterion in \ref{nak}\eqref{nak-a}; deleted the following subsections of \ref{nak}:
\sm

(-) (\cite[II, (4.5.3)]{K}) 
  Let $M$ and $N$ be finite projective $A$--modules. Then any $\overline{A}$--linear map $\vphi \co \ol M \to \ol N$ lifts an $A$--linear map $f \co M \to N$ satisfying $\can_N \circ f = \vphi \circ \can_M$:
      \[ \xymatrix@C=40pt{ M \ar@{-->}[r]^f \ar[d]_{\can_M}& N \ar[d]^{\can_N}     \\
        \ol M \ar[r]^\vphi & \ol N} \]
\sm

(-) Let $M$ be a finite $A$--module, and let $N$ and $P$ be  finite projective $A$--modules. We further assume that $g \co M \to P$ and $\vphi \co \ol P \to \ol N$ are $A$-linear maps such that $\vphi \circ \ol g $ is an isomorphism of $\ol A$--modules. Let $f \co P \to N$ be the lift of $\vphi$ according to \eqref{nak-b}. Then $f \circ g$ is an isomorphism and $g(M)$ is complemented in $P$:
    \[ \xymatrix@C=40pt{
        M \ar[r]^g \ar[d]_{\can} & P \ar@{-->}[r]^f \ar[d]_{\can}
          & N \ar[d]^{\can}     \\
          \ol M \ar[r]^{\ol g} & \ol P \ar[r]^\vphi & \ol N }
    \]
}

\subsection{Around Nakayama's Lemma.}\label{nak}
Let $\Jac(R)$ be the Jacobson radical of $R$ (recall $\Jac(R)$ is the intersection of all maximal ideals of $R$),  and let $\fra \subset \Jac(R)$ be an ideal of $R$. We denote the reduction mod $\fra$ by a ``bar''. Thus $\ol R = R / \fra$, $\ol M = M/ M \fra = M \ot_R \ol R$ for any $R$--module $M$ and $\ol f = f \ot \Id_{\ol R} \co \ol M \to \ol N$ for an $R$-linear map $f \co M \to N$ of $R$-modules. We let $\can = \can_M \co M \to \ol M$ be the canonical map. \sm

\begin{inparaenum}[(a)]
  \item \label{nak-a} (\cite[II, (4.2.2), (4.2.3)]{K}) 
    Let $f\co M \to N$ be an $R$-linear map between $R$--modules with $N$ being
    finitely generated. Then $f$ is surjective $\iff \bar f$ is surjective.
    Moreover, if $M$ has finite type and $N$ is finitely generated projective, then  $f$ is bijective if and only if $\bar f$ is bijective. In particular, a finite projective $R$--module $M$ is free if and only if $\ol M$ is free.
 \sm

\item\label{nak-dec} (\cite[II, (4.4.1)]{K}) Let $M$ be a finite projective $R$--module and let $M_1$, $M_2$ be submodules of $M$ such that $\ol M = \ol{M_1} \oplus \ol{ M_2}$. Then $M = M_1 \oplus M_2$. \sm

\item  \label{nak-aa} If $A\in \Ralg$ is a finite $R$--algebra, then  $\Jac(R)A \subset \Jac(R)$.
 Hence, applying \eqref{nak-a} to the $R$--module $A$, the $R$--module map $f=L_a$ (= left multiplication by $a\in A$), and $\ol a = \can(a)$, we get
\begin{equation}\label{jac-rad-1}
 a\in A\ti \quad \iff \quad \ol a \in {\ol A}{}\ti .
\end{equation}

\item \label{nak-f} (\cite[II, \S3.2, Cor.~of~Prop.~6 and \S3.3, Thm.~1]{BAC})  Let $f\co M \to N$ be an $R$--linear map between finite projective $R$--modules $M$ and $N$. Then the following are equivalent:
\end{inparaenum}

\begin{enumerate}[label={\rm (\roman*)}]

  \item  $f$ is bijective,

  \item the induced $(R/\m)$--linear map $M/\m M \to N/\m N$ is bijective for all maximal ideals $\m$ of $R$ ,

   \item $f$ is surjective and $\rank M_\m = \rank N_\m$ holds for all maximal ideals $\m$ of $R$.
\end{enumerate}

\subsection{Invertible elements}\label{inel} 
Let $A$ be a unital associative $R$--algebra. Recall that $a\in A$ is invertible if and only if the left multiplication $L_a \co A \to A$ by $a$ is invertible, if and only if the left and right multiplication by $a$ are surjective. 
We put $A\ti = \{a\in A : \text{$a$ is invertible}\}$, and note that invertibility is preserved under base change, thanks to $L_a \ot 1_S = L_{a\ot 1_S}$ for every $S\in \Ralg$. We list some  criteria for invertibility, complementing the ones of \ref{trno}\eqref{trno-a} and \eqref{trno-bou1}.  
\sm 

An element $x\in A$ is invertible if and only if 
\begin{enumerate}[label={\rm (\roman*)}]
   \item \label{ineli}
    $x\ot 1_{R_\gm}\in A_{R_\gm}$ is invertible for all maximal ideals $\gm \ideal R$, \cite[II, \S3.3, Thm.~1]{BAC};
     
   \item \label{inelii} $x\ot 1_E\in A_E$ is invertible for some faithfully flat $E\in \Ralg$, \cite[I, \S3.1, Prop.~2]{BAC};
\end{enumerate}
If $A$ is a finite $R$--algebra, then $x\in A$ is invertible if and only if 
\begin{enumerate}[label={\rm (\roman*)}]\setcounter{enumi}{2}
 \item \label{ineliii} $x\ot 1_{R/\gm}\in A \ot_R (R/\gm) = A /\gm A$ is invertible for all maximal $\gm \ideal R$, \cite[II, \S3.3, Prop.~11]{BAC}; 
 
   \item \label{ineliv} $x \ot 1_{A/\fra A} \in A/\fra A = A/\fra A$ is invertible for some ideal $\fra \ideal R$ contained in $\Jac(R)$, (use $\Jac(R) A \subset \Jac(A)$ by \cite[II, (4.2.4)]{K} and 
       $x\in A\ti \iff \can(x) \in (A/\Jac(A))\ti$ by \cite[II, (4.2.2)]{K}).  
\end{enumerate}

 \newpage

\section{One-generated algebras}\label{one-generated}


\subsection{One- and unit-generated algebras} \label{plr}
Let $A\in \Ralg$. For $a\in A$ the submodule 
\[ R[a] = \Span_R \{a^n : n \in \NN \}
\]
is a unital subalgebra of $A$. It is the image of the evaluation homomorphism
\begin{equation}   \label{plr-0}
{\rm ev}_a \co R[X] \to A; \quad f(X) \to f(a),
\end{equation}
where $R[X]$ is the polynomial ring over $R$ in the variable $X$.
If $T\in \Ralg$, the base change ${\rm ev}_A \ot 1_T$ of the evaluation map and the evaluation map of $a\ot 1_T \in A_T$ are related by the commutative diagram
\begin{equation} \label{plr-00} \vcenter{
\xymatrix{
   R[X] \ot_R T \ar[rr]^{ {\rm ev}_a\ot 1_T} \ar[rd]_\cong
   && A_T \\ 
  & T[X] \ar[ur]_{{\rm ev}_{ a\ot 1_T} }   } } \quad . 
\end{equation} 
 
We call $A\in \Ralg$ {\em one-generated\/}, or {\em one-generated over $R$\/} in case $R$ is important to know,  if $A=R[x]$ for some $x\in A$,  referred to as a {\em primitive element}. We will say that $A$ is {\em unit-generated\/} if $A=R[x]$ for some $x\in A\ti$. The equivalences in \eqref{splica-3} show that a one-generated algebra need not be unit-generated.

One-generated algebras are called {\em simple\/}  or {\em monogeneous\/} in \cite{Fer, Lo-genalg}, and {\em monogenic} in \cite{BFP}.\sm

For an $R$--algebra $A$ we denote by $\Prim_R(A)$ the set of primitive elements of the $R$--algebra $A$. Thus, 
\begin{equation}\label{plr-triv} 
\text{\em $x\in \Prim_R(A) \iff \rmev_x$ is surjective.}
\end{equation}
We write $\Prim(A)$ instead of $\Prim_R(A)$ if $R$ is clear from the context, and put $\Prim_R(A)\ti = \Prim_R(A) \cap A\ti$. We use the following elementary facts where throughout $A\in \Ralg$. \sm

\begin{inparaenum}[(a)]
\item\label{plr-a} ({\em Direct products of base rings}) Let $R=R_0 \times \cdots \times R_n$ be a direct product of rings. Any $A\in \Ralg$ has the form $A=A_0 \times \cdots \times A_n$ where $A_i$ are $R_i$--algebras, cf.~\ref{dpb}. Since $R[(a_0, \ldots, a_n)] = R_0[a_0] \times \cdots \times R_n[a_n]$ we get
      \begin{equation}\label{plr-a1} \begin{split}
        \Prim_R(A) &= \Prim_{R_0}(A_0) \times \cdots \times \Prim_{R_n}(A_n), \quad \text{and} \\
        \Prim_R(A)\ti &= \Prim_{R_0}(A_0)\ti \times \cdots \times \Prim_{R_n}(A_n)\ti.
      \end{split}\end{equation}

 \item \label{plr-func} ({\em Functoriality}) Let $f\co A \to B$ be a {\em surjective\/} $R$--algebra homomorphism. Then $f\big(\Prim_R(A)) \subset \Prim_R(B)$. \sm %

  \item\label{plr-bb} ({\em Direct products of $R$--algebras}) Let $A=A_1 \times \cdots \times A_n$ be a direct product of $R$--algebras.  Applying \eqref{plr-func} to the projection $\pr_i \co A \to A_i$ homomorphism, yields 
      \begin{equation} \label{plr-b1}
        \Prim_R(A_1 \times \cdots \times A_n) \subset
          \Prim_R(A_1) \times \cdots \times \Prim_R(A_n).\end{equation}
      As \eqref{splica-3} shows, the inclusion \eqref{plr-b1} is in general not an equality.

The set $\Prim_R(A)$ has the following description as a subset of the right-hand side of \eqref{plr-b1}. Let $a=(a_1, \ldots, a_n)\in \Prim_R(A_1) \times \cdots \times \Prim_R(A_n)$, and let $\fra_i = \Ker({\rm ev}_{a_i})$, so that $A_i \cong R[X]/\fra_i$ for $i=1, \ldots, n$. By the Chinese Remainder Theorem, 
$a\in \Prim_R(A)$, i.e., ${\rm ev}_a$ is surjective, 
if and only if the ideals $\fra_i$, $i=1, \ldots, n$,  are relatively prime in the sense that $\fra_i + \fra_j = R[X]$ for all $i\ne j$.

A situation, where \eqref{plr-b1} is an equality, is described in Lemma~\ref{lem_generatornew}. \sm

\item\label{plr-bc} ({\em Base change, $R$--functors}) Let $S\in \Ralg$ and let $x\in \Prim_R(A)$, i.e., $\rmev_x \co R[X] \to A$ is surjective. Since then 
    $\rmev_x \ot \Id_S \co R[X] \ot_R S \to A \ot_R S$ is surjective, the commutative diagram \ref{plr-00} shows that $\rmev_{x\ot 1_S}$ is surjective, i.e., $x\ot 1_S \in \Prim_S(A_S)$.  \sm

For any not necessarily one- or unit-generated $R$--algebra $A$ we  define subfunctors $\ul \Prim_R(A)$ and $\ul \Prim_R(A)\ti$ of the $R$--functor\/ $\ulW(A)$ by assigning to $S\in \Ralg$ the subsets
 \[  \ul\Prim_{R}(A)\, (S)= \Prim_S(A_S) \quad \text{and} \quad
         \ul\Prim_R (A)\ti (S) = \Prim_R(A_S) \ti \]
   of $\ulW(A)(S) = A \ot_R S$ and to a homomorphism $f \co S \to T$ in $\Ralg$ the set map $\Prim_S(A_S) \to \Prim_T(A_T)$, $x \mapsto x \ot_S 1_T$. 
 We will show in \ref{prsch} that $\ul\Prim_R(A)$ is representable by an open finitely presented subscheme of $\uW(A)$ if $A$ is finite locally free. \sm

\item \label{plr-d}({\em Example})  Let $f\in R[X]$ be  a monic polynomial over $R$. Then $A=R[X]/(f)$ is a one-generated $R$--algebra with $x=X + (f)\in A$ as primitive element. The $R$--module $A$ is free of rank $\deg(f)$. In the finite locally free case, this example is actually the general case, see \ref{genolemLG}\ref{genolemLGv}. \sm 

\item\label{plr-f} ({\em \'Etale algebras over an infinite base field}) Let $R=k$ be an infinite field and suppose $E\in \kalg$ has only a finite number of subalgebras, e.g., $E$ is an  \'etale $k$--algebra as in  \ref{fea}\eqref{fea-a}. Then 
     $\Prim_k(E) \ne \emptyset$ by \cite[V, \S7.4, Prop.~7]{BA5}. 
     In the \'etale case, $E=K_1 \times \cdots \times K_n$ is a finite product of finite separable field extensions. Since then 
    $\Prim_k(K_i) = \Prim_k(K_i)\ti$, it follows 
    that $E$ is even unit-generated. See \eqref{splica-3} and Corollary~\ref{cor_primitive_etale} for unit-generation of finite \'etale $R$--algebras.  
\sm 

\item\label{plr-c} ({\em Reduction modulo the Jacobson radical})  We use the setting of \ref{nak}: $I \subset \Jac(R)$ is an ideal of $R$, and $A$ is a finite $R$--algebra. We put $\ol R = R/I$, $\ol A = A \ot_R \ol R = A/IA$, and $\ol a = \can(a) \in \ol A$ for $a\in A$. Then
\begin{equation} \label{plr-c11} \begin{split}
   a\in \rmP_R(A) \quad &\iff \quad \ol a \in \rmP_{\ol R}(\ol A) ,\\
   a\in \rmP_R(A)\ti \quad &\iff \quad \ol a \in \rmP_{\ol R}(\ol A)\ti.
\end{split} \end{equation}
In particular, $A$ is one- or unit--generated if and only if $\ol A$ is so.

For the proof of \eqref{plr-c11} consider the evaluation map ${\rm ev}_a \co R[X] \to A$ of \eqref{plr-0}. Then
\[ \xymatrix{
   R[X] \ot_R {\ol R} \ar[rr]^{\ol {{\rm ev}_a}} \ar[rd]_\cong
   && \ol A
   \\ & {\ol R}[X] \ar[ur]_{{\rm ev}_{\ol a} }} \]
is a commutative diagram. By \ref{nak}\eqref{nak-a} we then get $a\in \rmP_R(A) \iff {\rm ev}_a$ is surjective $\iff {\rm ev}_{\ol a}$  is surjective $\iff \ol a \in \rmP_{\ol R}(\ol A)$. The second equation in \eqref{plr-c11} follows from the first and \eqref{jac-rad-1}. 
\end{inparaenum}

\comments{(2026-01-28) New version of the lemma~\ref{lem_generatornew}. Before we assumed $R=k$ field. Simpler proof, due to use of Chinese Remainder Theorem in \ref{plr}\eqref{plr-bb}. The lemma has been used in the old version of the Knebusch preprint; I don't know yet about the new version.

(2026-04-11) New version of Lemma~\ref{lem_generatornew} to be used in the proof of Corollary~\ref{cor_primitive_etale}. 
}

The following Lemma~\ref{lem_generatornew} is inspired by \cite[Prop.~2.3]{FRS}. 

\begin{lem}\label{lem_generatornew} Let $(K_i)_{i\in I}$ be a finite family of $R$--fields and let ,  $I = \bigsqcup_{j=1}^t I_j$ be the unique partition of $I$ such that $K_i \cong K_m \iff $ there exists $j$ such that $i,m\in I_j$. Finally, put  
\[ A = \textstyle \prod_{i\in I} K_i, \qquad A_j = \prod_{i\in I_j} K_i 
\] 
so that $A= A_1 \times \cdots \times A_t$. Then 
\[  \Prim_R(A) =  \Prim_R(A_1) \times \cdots \times \Prim_R(A_t) = \Prim_R(A)\ti.\]
\end{lem}

\begin{proof} Regarding the first equality, by \eqref{plr-b1} we only need to prove that any $a=(a_1, \ldots, a_t) \in \Prim_R(A_1) \times \cdots \times \Prim_R(A_t)$ is a primitive element of $A$. Write $I_j = \{ i_{j1}, \ldots, i_{j r_j}\}$ and correspondingly $a_j = (a_{j1}, \ldots, a_{j r_j})$ where each $a_{j\ell}$, $1\le \ell \le r_j$, is a primitive element of $K_{j\ell}$, again by \eqref{plr-b1}. Let $\gm_{j\ell}$ be the kernel of the evaluation map $R[X]\to K_{jl}$ by $a_{j\ell}$. By the description of $\Prim_R(A)$ in \ref{plr}\eqref{plr-bb}, we know that $a\in \Prim_R(A)$ if and only if the family $(\gm_{j\ell})_{1\le j \le t, 1\le \ell \le r_j}$ of maximal ideals of $R[X]$ consists of relatively prime ideals, i.e., are pairwise distinct. But this is indeed the case: if $\gm_{j\ell} = \gm_{j'\ell'}$ then $K_{j\ell} = K_{j'\ell'}$ and therefore $j=j'$ be construction of $A_j$. The same argument, applied to the primitive element $a_j \in A_j$, shows that $\gm_{j\ell'} =\gm_{j\ell'}$ implies $\ell = \ell'$. The second equality follows $\Prim_R (A_i) = \Prim_R (A_i)\ti$. 
\end{proof}
\sm

In the following Lemma~\ref{genolemLG} we characterize primitive elements of a finite $A\in \Ralg$, although later we will only be interested in $A$ finite locally free. 

\begin{lem}[Characterization of primitive elements for $A$ finite] \label{genolemLG} Let $A\in \Ralg$ be a finite $R$--algebra. Then the following are equivalent for $x\in A$: \sm
\begin{enumerate}[label={\rm (\roman*)}]
   \item\label{genolemLGi} $x\in \Prim_R(A)$;  
   
  \item \label{fgenii} $x\ot 1_{R/\gm} \in \Prim_{R/\gm}( A _{R/\gm})$ for every maximal ideal $\gm \ideal R$;  
      
  \item\label{fgeniii} $x\ot 1_E\in \Prim_E(A_E)$ for some faithfully flat $E\in \Ralg$; 

\item \label{fgeniv} $\ol x \in \Prim_{R/I}(A/IA)$ in the setting of {\rm \ref{nak}}, i.e., $\ol x = x\ot 1_{A/IA} \in A/IA = A\ot_R (R/I)$ for some ideal $I \subset \Jac(R)$. 
\end{enumerate}
Moreover, if $A\in \Ralg$ is finite projective of constant rank $d\in \NN_+$ as $R$--module, then \ref{genolemLGi}--\ref{fgeniv} are equivalent to \ref{genolemLGv}--\ref{genolemLGiii}:

\begin{enumerate}[label={\rm (\roman*)}]\setcounter{enumi}{4}     
 \item \label{genolemLGv} the characteristic polynomial $\Pc_{A/R}(x;X)$ of $x\in A$, cf.~{\em \ref{trno-0}}, is annihilated by the algebra homomorphism $\rmev_x$, and $\rmev_x$ induces an isomorphism 
     \[ 
          R[X]/\big(\Pc_{A/R}(x; X)\big)\simlgr A, \quad X + \big(\Pc(x; X)\big) \mapsto x
     \]
     of $R$--algebras; 

\item \label{genolemLGii} $A$ is free of rank $d$ as $R$--module and $1_A = x^0, x, x^2, \ldots , x^{d-1}$ is a basis of the $R$--module $A$; 
       
  \item\label{genolemLGiii} $1_A \we x \we x^2 \cdots \we x^{d-1}$ is a basis of the invertible $R$--module  $\bigwedge^d_R A$, equivalently, is a unimodular vector. 

\end{enumerate}
Finally, if $A$ is free of rank $d$ as $R$--module, say with basis $e_0, \ldots, e_{d-1}$, define for $x\in A$ the matrix $M_x= (r_{ij})_{0\le i,j <d}\in \Mat_d(R)$ by $x^i = \sum_{i,j} r_{ij} e_j$. Then 
\begin{equation} \label{genolemLGb} 
    x \in \Prim_R(A) \iff \det M_x \in R\ti. 
\end{equation}
\end{lem}

\begin{proof} 
\ref{genolemLGi} $\iff$ \ref{fgenii}: By \cite[II, \S3.3, Prop.~11]{BAC}, the evaluation map $\rmev_x$ is surjective if and only if it induces a surjective linear  map $\ol{\rmev_x} \co R[X]/(\gm R[X]) \to A/\gm A$ for every maximal ideal $\gm \ideal R$. Under the canonical isomorphisms $R[X]/(\gm R[X]) \cong R[X]\ot_R (R/\gm)$ and $A/\gm A \cong A \ot_R (R/\gm)$, the induced map $\ol{\rmev_x}$ becomes $\rmev_{x (\gm)}$. 

\ref{genolemLGi} $\iff$ \ref{fgeniii} holds by \cite[I, \S3.2, Prop.~11]{BAC}, and \ref{genolemLGi} $\iff$ \ref{fgeniv} follows from \ref{nak}\eqref{nak-a} and the commutative diagram \eqref{plr-00} for $T=R/I$.

\ref{genolemLGi} $\implies$ \ref{genolemLGv}: The algebra homomorphism ${\rm ev}_x \co R[X] \to A$, $f(X) \mapsto f(x)$ is surjective. Hence, $A \cong R[X]/\Ker {\rm ev}_x =:A'$ is finite projective of constant rank $d$ as $R$--module. But then \cite[0.11]{Lo-genalg} shows \ref{genolemLGii}. (Indeed, let $z = X + \Ker {\rm ev}_x \in A'$, let $L_z$ be the left multiplication of $A'$ by $z$, and let $\Pc(z) = \Pc_{A/R}(z; X)$ be the characteristic polynomial of $z\in A'$, cf.~\eqref{trno-0}, which is a monic polynomial of degree $d$. Hence $R[X]/(\Pc(z))$ is free of rank $d$. By Cayley-Hamilton for $A$, \cite[IV, Cor.~2.3]{KO}, $\Pc(z)$ lies  in $\Ker {\rm ev}_x$, so that we get a surjective linear map $R[X]/(\Pc(z)) \to A$ between finite projective $R$--modules of rank $d$. Such a map is necessarily an isomorphism.)

The implications \ref{genolemLGv} $\implies $ \ref{genolemLGii} $\implies$ \ref{genolemLGi} are  obvious, and 
\ref{genolemLGii} $\iff$ \ref{genolemLGiii} is a special case of \cite[0.6]{Lo-genalg}. (Indeed, \ref{genolemLGii} $\implies$ \ref{genolemLGiii} is standard, see for example \cite[III, \S7.8, Thm.~1]{BA}. To see that \ref{genolemLGiii} $\implies$ \ref{genolemLGii}, let $\be \in (\bigwedge^d_R(A))\ch$, $\be(1_A \we x \we \cdots \we x^{d-1}) = 1$ and define linear forms $\al_i \in A\ch$ by 
\[ 
\al_i (a) = \be(1_A \we x \we \cdots \we x^{i-1} \we a \we x^{i+1} \we \cdots \we x^{d-1}) 
\] 
for $a\in A$. Then $\al_i (x^j) = \de_{ij}$ implies that $1_A, x, \ldots, x^{d-1}$ is a basis of the submodule $\Span_R \{ 1_A, x, \ldots, x^{d-1} \}\subset A$ and this submodule is complemented by $\bigcap \Ker \al_i$, hence it coincides with $A$.) 

The equivalence \eqref{genolemLGb} is clear from the characterization \ref{genolemLGii} of a primitive element. \end{proof}

\begin{remarks}
\label{genolemLGrem} \begin{inparaenum}[(a)] 
\item \label{genolemLGrem-a} ({\em $\Prim_R(A)\ti$}) 
One obtains a description of the elements in $\Prim_R(A)\ti$ by combing the characterizations of $\Prim_R(A)$ above with those of \ref{inel}. For example, in the setting of \ref{genolemLG}\ref{fgeniv} we get for a finite $A\in \Ralg$ that  
\begin{equation} \label{plr-c1} \begin{split}
   x\in \Prim_R(A)\ti \quad &\iff \quad \ol x \in \Prim_{\ol R}(\ol A)\ti.
\end{split} \end{equation}
Since $A \to \ol A$ is surjective, it follows that 
\begin{equation}\label{plr-c2}
\text{\em $A$ is one- or unit--generated if and only if $\ol A$ is so.}
\end{equation}

Also, in the setting of \ref{genolemLG}\ref{genolemLGv}, 
\begin{equation}\label{genolemLGa1}
 x \in \Prim_R(A)\ti \iff 
   \Pc_{A/R}(x;X) (0) \in R\ti. 
\end{equation}

\item \label{genolemLGrem-b} The equivalence \ref{genolemLGi} $\iff$ \ref{fgeniii} in \ref{genolemLG} says that $\ul\Prim_R(A)$ is a sheaf in the (fpqc)--topology. By a theorem of Grothendieck \cite[Tag 023Q]{St},
    this follows if $\ul \Prim_R(A)$ is representable by a scheme, e.g. if $A$ is finite locally free, as we show in \ref{prsch}.  
\end{inparaenum}
\end{remarks}


\begin{cor} \label{genolem-semi-LG}
Let $R$ be a semilocal ring with maximal ideals $\m_i$, $1\le i \le n$, and residue fields $\ka_i = R/\m_i$, and let $A\in \Ralg$ be a finite $R$--algebra. Then the following are equivalent:
\begin{enumerate}[label={\rm (\roman*)}]
 \item \label{genolem-semi-i} $A$ is one-generated (unit-generated respectively) over $R$,

 \item \label{genolem-iv} $A\ot_R R_{\m_i}$, $1\le i \le n$,  is one-generated
 (unit-generated respectively) over $R_{\m_i}$,

 \item\label{genolem-iii} $A \ot_R \ka_i$, $1\le i \le n$, is one-generated  (unit-generated respectively) over $\ka_i$.
\end{enumerate}
In this case, if $A \ot_R \ka_i = \ka_i[x_i]$ for $i=1, \ldots, n$, then $A=R[x]$ where $x$ is any lift of $(x_1, \ldots, x_n) \in \prod_i A\ot_R \ka_i$. 
\end{cor}

\begin{proof} The canonical algebra homomorphisms $R \to R_{\gm_i} \to \ka_i$ show that \ref{genolem-semi-i} $\implies$ \ref{genolem-iv} $\implies$ \ref{genolem-iii}. If \ref{genolem-iii} holds, say $A\ot_R \ka_i$ is generated by $x_i$ as $\ka_i$--algebra, then 
$(x_1, \ldots, x_n)$ is a primitive element of the $R/\Jac(R)$--algebra 
$A/\Jac(R)A = (A\ot_R \ka_1) \times \cdots \times (A \ot_R \ka_n)$, \eqref{plr-a1}. Any lift $x\in A$ of $(x_1, \ldots, x_n)$ is a primitive element of the $R$--algebra $A$ by \ref{genolemLG}\ref{fgeniv}. \end{proof}
\sm 

We will prove a related version of Corollary~\ref{genolem-semi-LG} for arbitrary $R$,  but $A$ finite locally free in \ref{loog}. 
\sm 

We will use the discriminant of $a\in A$ to describe primitive elements of finite \'etale $R$--algebras in \ref{lem_primitive}.

\subsection{Discriminant of an element}\label{disc} Let $A\in \Ralg$ be finite projective, and let $a\in A$.  \sm 

At first assume that $A$ has constant rank $d\in \NN$. Following \cite[IV; \S6.7]{BA5}, we define the {\em discriminant\/} of $a\in A$ as the discriminant of the characteristic polynomial $\Pc_{A/R}(a; X)$ 
\begin{equation}  \label{disc-1}
 \dis_{A/R}(a) =  \dis \big( \Pc_{A/R}(a;X)\big) =
  \det\big( (\Tr_{A/R}(a^{i+j}))_{0\le i,j< d}\big),
\end{equation}
where the second equation in \eqref{disc-1}  follows from formula (45) in {\em loc.\ cit.}.

By \ref{trno}\eqref{trno-d}, \begin{equation}
  \label{disc-2} \big(\dis_{A/R}(a)\big) \ot 1_S = \dis_{A\ot_R S/S}(a \ot 1_S)
\end{equation}
holds for any $S\in \Ralg$. Moreover, we have
\begin{equation}\label{disc-3}
 \dis_{A/R}(ra) = r^{d(d-1)} \dis_{A/R}(a) \quad (r\in R).
\end{equation}

If $A\in \Ralg$ does not have constant rank, we use the unique rank decomposition of $A$, \ref{dpb}. Thus,  $R=R_0 \times \cdots \times R_n$ and $A=A_0 \times \cdots \times A_n$ where $A_i$, $i=0, \ldots, n$, is a projective $R_i$--algebra of constant rank $i$, and define
\begin{equation}\label{disc-4}
\dis_{A/R}(a_1, \ldots, a_n)
   = (\dis_{A_1/R_1}(a_1), \ldots, \dis_{A_n/R_n}(a_n)\big)
\end{equation}

The formulas \eqref{disc-2} and \eqref{disc-3}, together with Roby's approach to polynomials \cite{Roby}, show that
\begin{equation}\label{disc-3a}
 \dis \co \uW(A) \to \GG_a, \quad a (\in A \ot_R S) \mapsto \dis_{A\ot S/S}(a)
\end{equation}
defines a polynomial on $\uW(A)$,  which is homogeneous of degree $d(d-1)$ in case $A$ has constant degree $d$. 
\sm 

\subsection{Example: the split \'etale case $A_0=R^d$, $d\in \NN_+$} \label{splica} Let $A_0$ be the split \'etale $R$--algebra of rank $d\in \NN_+$. Thus
$A_0= R\times \cdots \times R$, the direct algebra product of $d$ factors $R$, and let $a=(r_1, \ldots, r_d) \in A_0$. We have $\Pc_{A_0/R}(a; X) = \prod_{1\le i \le d} (X-r_i)$. Hence, by \eqref{disc-1} and \cite[IV, \S6.7, (46)]{BA5},
\begin{equation}  \label{splica-1}
 \dis_{A_0/R}(a) = \textstyle \prod_{i<j}\, (r_i - r_j)^2  = (-1)^{d(d-1)/2} \, \prod_{i\ne j} (r_i - r_j).
\end{equation}
One knows \cite[Ex.~1.2]{Fer}:
\begin{equation}\label{splica-2}
 \text{\em $a$ is a primitive element of $A_0$} \quad \iff \quad \dis_{A_0/R}(a) \in R\ti.
\end{equation}
Indeed, observe that the Vandermonde matrix
\[
 V = \begin{pmatrix}
1 & r_{1} & \cdots & r_{1}^{d-1} \\
1& r_{2} & \cdots & r_{2}^{d-1} \\
\vdots  & \vdots  & \ddots & \vdots  \\
1 & r_{d} & \cdots & r_{d}^{d-1}
\end{pmatrix}
\]
is the transition matrix from the standard basis of $R^d$ to the sequence $(1, a, a^2, \ldots , a^{d-1})$ in $A_0$, so that \eqref{splica-2} follows from \eqref{genolemLGb}.  
\sm

{\em Special case: $R=k$ is a field}. The formula \eqref{splica-1} and the equivalence \eqref{splica-2} imply for the $k$--algebra $k^d$:
 \begin{equation} \label{splica-3} \begin{split}
     \text{$k^d$ is one-generated} \quad &\iff \quad |k|\ge d , \\
      \text{$k^d$ is unit-generated} \quad &\iff \quad |k|> d, 
\end{split} \end{equation}
see Corollary~\ref{cor_primitive_etale} for the non-split case. 
\sm 

By \ref{fea}\eqref{fea-sp}, an \'etale $R$--algebra $A$ of constant rank $d\in \NN_+$ can be split by a  Galois extension $T\in \Ralg$ with Galois group $\frS_d$. The formula \eqref{splica-1} shows that $\dis_{A_0/R}$ is $\frS_d$--invariant and therefore descends to a function on $A$, which is $\dis_{A/R}$ by \eqref{disc-2}. 

\comments{(2025-12-02) In the old version we required in \ref{lem_primitive} that $A$ is faithfully projective. I think faithfulness is not necessary. The $0$ algebra is \'etale, generated by the primitive element $0$    }

\begin{lem} \label{lem_primitive} Let $A$ be a finite \'etale $R$--algebra. Then $a\in A$ is a primitive element of $A$ if and only if $\dis_{A/R}(a)\in R\ti$. \end{lem}

\begin{proof} In view of the definition \eqref{disc-4}, it is no harm to assume that $A$ is \'etale of rank $d\in \NN_+$.
The element $a\in A$ is primitive if and only if the linear map $R^d \to A$, $(r_0, \ldots, r_s) \mapsto \sum_{0 \le i < d} \, r_i a^i$ is an isomorphism of $R$--modules. This is the case if and only if it holds after a faithfully flat extension. Similarly, $\dis_{A/R}(a) \in R\ti$ if and only if $\dis_{A/R}(a)$ becomes invertible after some faithfully flat extension. Since $A$ becomes split after some faithfully flat extension, see \ref{fea}\eqref{fea-sp}, it is no harm to assume that $A$ is already split. But in this case the claim has been proven in \eqref{splica-2}. \end{proof}
\sm 

We can now show that the $R$--functor $\ul \Prim(A)$ is representable by an affine open subscheme of $\uW(A)$, which is even principal affine in case $A$ is finite \'etale or free as $R$--module, see the Remarks~\ref{prschrem} for some background. In the finite \'etale case the representing scheme is universally schematically dense, using the terminology of \ref{tod}\eqref{tod-b}. 

\begin{prop}[The scheme $\uPrim(A)$] \label{prsch} Let $A\in \Ralg$ be finite locally free. Then the $R$--functor $\ul\Prim(A)$, defined in {\rm \ref{plr}\eqref{plr-bc}}, is representable by a finitely presented open affine  subscheme $\uPrim (A)$ of $\uW(A)$. \sm 

\begin{inparaenum}[\rm (a)] \item \label{prscha} More precisely, decompose $R=R_0 \times \cdots \times R_n$ such that $A_d= A \ot_R R_d$, $0\le d \le n$, has constant rank $d$ and $A=A_0 \times \cdots \times A_n$, cf.~{\rm \ref{dpb}}. Then 
\begin{equation}  \label{prscha1}
 \uPrim (A) = \textstyle \bigsqcup_{\, 0 \le d \le n}\, \uPrim (A_d), 
\end{equation}
where $\uPrim(A_d)$ is the pull-back of the open subscheme $\uW\big(\bigwedge^d(A_d)\big)_u$ of $\uW(\bigwedge^d(A_d))$, see {\rm \ref{spreq-a}}, under the polynomial $a \mapsto a^0 \we a \we \cdots \we a^{d-1}$. \sm 

\item \label{prschb} Suppose $A$ is finite \'etale. Then $\uPrim(A)$ is isomorphic to the principal open subscheme of $\uW(A)$ determined by the discriminant, 
\begin{equation}  \label{prschb2}
   \uPrim(A) \cong \uW(A)_\dis\, ;  
\end{equation}
it is universally schematically dense in $\uW(A)$, in particular $\uPrim(A) \ne \emptyset$. \sm 

\item\label{prschc} Suppose $A$ is free of finite rank $d$. Let $\de$ be the polynomial of $\uW(A)$, given by $\det M_x$ for $x\in A_T$, $T\in \Ralg$ and $M_x$ as in {\rm \ref{genolemLG}}. Then $\uPrim(A)$ is isomorphic to the principal open subscheme of $\uW(A)$ determined by $\de$, 
    \begin{equation}  \label{prschc1}
   \uPrim(A) \cong \uW(A)_\de.  
\end{equation}
\end{inparaenum}\end{prop}

\begin{proof} Because of \eqref{plr-a1} the $R$--functor $\ul \Prim (A)$ is the direct product of the $R$--functors $\ul \Prim (A_d)$. Hence, by \eqref{dpb3}, it is sufficient to prove that the $R_d$--functor $\ul \Prim (A_d)$ is representable as specified in \eqref{prscha}. Without loss of generality we can therefore assume that $A$ has constant rank $d$ as $R$--module. In this case, we have the polynomial map $p_d \co \uW(A) \to \uW(\bigwedge^d A)$, given on $R$--points by assigning to $a\in T\in \Ralg$ the vector $1_A \we a \we a^2 \cdots \we a^{d-1}$. It then follows from Lemma~\ref{genolemLG}\ref{genolemLGiii} that $\ul \Prim(A)$ is representable by the fibre product
\begin{equation} \label{prscha2} \vcenter{ 
\xymatrix@C=50pt{\ar @{} [dr] |{\small\qed} %
\uPrim(A)  \ar[r] \ar[d] & \uW(A) \ar[d]^{p_d} 
 \\ \uW\big(\bigwedge^d(A)\big)_u \ar[r] & \uW\big(\bigwedge^d(A)\big)  } }  \quad.  \end{equation}
By \ref{spreq-a}, the bottom open immersion in \eqref{prscha2} is finitely presented. By base change, the same is therefore also true for the top horizontal map. 

Moreover, the bottom horizontal map is affine by \ref{spreq-a}\ref{spreq-aiii} and \cite[Tag 01SH]{St}. Hence, by base change for affine morphisms \cite[Tag 01SD]{St}, so is the morphism $\uPrim(A) \to \uW(A)$. In particular, $\uPrim(A)$ is an affine scheme.%
\sm

\eqref{prschb} That $\uW(A)_\dis$ represents the $R$--functor $\ul \Prim(A)$ follows from \ref{lem_primitive}. 
We apply the criterion~\ref{usd-d} to prove that $\uPrim(A)$ is universally 
schematically dense in $\uW(A)$: 
it suffices to show that $\ul \Prim(A)(K) =\Prim(A_K) \ne \emptyset$ for every algebraically closed field $K$ in $\Ralg$. But $A_K$ is split \'etale, so that $\Prim(A_K)\ne \emptyset$ follow from \eqref{splica-3}.  
\eqref{prschc} is a consequence of \eqref{genolemLGb}. 
\end{proof}

\begin{remarks}\label{prschrem} That the $R$--functor $\ul \Prim(A)$ is represented by a finitely presented open subscheme of $\uW(A)$ is a special case of \cite[Prop.~1.14]{Lo-genalg}, proven there for $A$ a unital Jordan algebra which is finite locally free as $R$--module. We have included a proof for the convenience of the reader. 

On the other hand, that $\ul \Prim(A)$ is represented by an affine $R$--scheme, is proven for Noetherian rings $R$ and $A$ in \cite[Thm.~1.4]{ABHS}. 
Moreover, also Theorem~\ref{prsch}\eqref{prschc} is shown in \cite[Thm.~3.5]{ABHS}.  
That paper's approach is based on the fact that $a\in A$ is primitive element if and only if the morphism $\bbA^1 \to \uW(A)$ corresponding the evaluation map ${\rm ev}_a$ of \eqref{plr-0} is a closed immersion. This allows the authors to consider primitivity of algebras over a non-affine base. 
\end{remarks}

\begin{cor}[The scheme $\uPrim(A)\ti$] \label{prischinv}
Let $A\in \Ralg$ be finite locally free. Then the $R$--functor $\ul\Prim(A)\ti$, defined in {\rm \ref{plr}\eqref{plr-bc}}, is represented by the finitely presented open subscheme 
\begin{equation}   \label{prischinv1}
 \uPrim(A)\ti = \uPrim(A) \cap \uGL_1(A)
\end{equation}
of $\uW(A)$.  If $A$ is finite \'etale, \new $\uPrim(A)\ti$ is a principal open and \enew universally schematically dense subscheme of $\uW(A)$. 
\end{cor}

\begin{proof}
The formula \eqref{prischinv1} follows from the definition of the $R$--functor $\ul \Prim(A)\ti = \ul \Prim(A) \cap \ul \GL_1(A)$. Thus $\uPrim(A)\ti$ is the fibre product
\[   \vcenter{ 
\xymatrix@C=50pt{\ar @{} [dr] |{\small\qed} %
\uPrim(A)\ti  \ar[r] \ar[d] & \uGL_1(A) \ar[d]^{\inc} 
 \\ \uPrim(A) \ar[r]^\inc & \uW(A) } } \quad.  
 \] 
Since ``finite presentation'' allows base change and composition, 
$\uPrim(A)\ti$ is a finitely presented open subscheme of $\uW(A)$. Finally, in the finite \'etale case, $\uPrim(A)\ti$ is the intersection of two universally schematically dense open subschemes of $\uW(A)$, \ref{prsch}\eqref{prschb} and \ref{trno}\eqref{trno-ds}. It is therefore also universally schematically dense by \cite[Tag 01RF]{St}. 
\new
It is a principal open subscheme of the affine scheme $\uW(A)$, since the intersection of two principal open subschemes of the affine scheme $\uW(A)$, \ref{trno}\eqref{trno-ds} and \eqref{prschb2}, is again principal open \cite[Tag 00E0(15)]{St}. 
\enew
\end{proof}
 
\comments{(2026-02-02) Regarding \ref{log}, we only use Zariski so far. It makes of course sense for other Grothendieck topologies, like \'etale, fppf, or fpqc. 
These are considered in \cite{ABHSII}. }
\pcomments{(2026-02-23) \ref{log}: Du point de vue du vocabulaire, Ferrand a raison,
« local » se rapporte \`a l’anneau de base. Maintenant que
les autres ont mis une d\'efinition non orthodoxe, c'est compliqu\'e
de changer les choses. Je pense que l'on reste comme on est.}

\subsection{Local one-generation}\label{log} 
We say that $A\in \Ralg$ is {\em Zariski--locally one-generated\/} if there exists a standard Zariski--cover $R'\in \Ralg$ such that $A\ot_R R'$ is one-generated as $R$--algebra. The notion of a {\em Zariski--locally unit-generated} $R$--algebra is the obvious one.  

These types of algebras have been considered in \cite{ABHSII} for $R$ replaced by a locally noetherian scheme and $A$ by a finite locally free morphism $S'\to S$ of constant degree. A weaker concept of a fpqc-one-generated $R$--algebra is used in \cite{Fer}: $A\ot_R R'$ is required to only be one-generated as $R'$--algebra.  

\sm 

The following Corollary~\ref{loog} of Proposition~\ref{prsch} is in the spirit of Corollary~\ref{genolem-semi-LG}. The equivalence \ref{loogi}--\ref{loogiii} is proven in \cite[Thm.~2.1]{ABHSII} in the setting of that paper. 

\comments{(2026-02-13) We use \ref{loog}\ref{loogc} in the proof of Corollary \ref{spo-LGG} (Springer section)}

\begin{cor}\label{loog} Let $A$ be a finite locally free algebra over an arbitrary base ring $R$. Then the condition
\begin{enumerate}[label={\rm (\roman*)}]
\item \label{loogv} $A$ is one-generated
\end{enumerate}
implies the equivalent conditions \ref{loogi}--\ref{loogiii} below:  
\begin{enumerate}[label={\rm (\roman*)}]\setcounter{enumi}{1}
 \item\label{loogi} $A$ is Zariski-locally one-generated,
 
 \item\label{loogin} $A_{R_\gp}$ is one-generated for every $\gp\in \Spec(R)$,  
 
 \item \label{loogii} $A_{R_\gm}$ is one-generated for every maximal ideal $\gm\ideal R$,  
     
\item \label{loogiii} $A_{R/\gm}$ is one-generated for every maximal ideal $\gm \ideal R$.    
\end{enumerate}
Moreover:  \sm 
 
\begin{enumerate}[label={\rm (\alph*)}]
\item\label{looga} If $|R/\gm| = \infty$ for all maximal ideals $\gm\ideal R$ and $\uPrim_R(A) \ne \emptyset$, then the conditions  \ref{loogi}--\ref{loogiii} are satisfied. Thus, $A$ is Zariski-locally one-generated. \sm 

\item \label{loogb} If $R$ is an LG ring, then the conditions \ref{loogv}--\ref{loogiii} are equivalent. \sm

\item \label{loogc} If $R$ is an LG ring with $|R/\gm|= \infty$ for all maximal $\gm \ideal R$, i.e., $R$ satisfies the primitive condition {\rm \ref{revLGG}},  and if $\uPrim(A) \ne \emptyset$, then the conditions \ref{loogv}--\ref{loogiii} are satisfied. Thus, $A$ is one-generated.  
\end{enumerate}
\end{cor}

\begin{proof} Since $\uPrim(A)$ is an open finitely presented (= quasi-compact) 
subscheme of $\uW(A)$, the equivalences \ref{loogi} $\iff$ \ref{loogin} $\iff$ \ref{loogii} follow from Corollary~\ref{fibn}.
To prove \ref{loogii} $\iff$ \ref{loogiii} we can assume that 
$R$ is a local ring with maximal ideal $\gm$, in which case the equivalence is a special case of Corollary~\ref{genolem-semi-LG}. 
\sm 

\ref{looga} Since $\uPrim_{R/\gm}(A_{R/\gm})$ is a non-empty open subvariety of an affine space, it is unirational so that condition~\ref{loogiii} holds by the criterion in \ref{zdens}\eqref{adc}. If $R$ is an LG ring, then \ref{loogv}--\ref{loogiii} are equivalent by Corollary~\ref{fibn}, which proves \ref{loogb}, and then \ref{loogc} follows by combining \ref{looga} and \ref{loogb}. 
\end{proof}
\sm 

\textbf{Remarks.} We note that $\uPrim(A) \ne \emptyset$ whenever $\uPrim_R(A)$ is universally schematically dense in $\uW(A)$, e.g., for $A$ a finite \'etale $R$--algebra. 

Replacing $\uPrim(A)$ by the scheme $\uPrim(A)\ti$ of Corollary~\ref{prischinv}, one sees that Lemma~\ref{loog} holds mutatis mutandis for ``unit-generated'' instead of ``one-generated''. We leave the details to the reader.

\new
\comments{Corollary~\ref{cor_primitive_etale} was a corollary in section 3 of the previous version 2026-04-03. It is used in the proof of Corollary~\ref{cor_gonflement}. Philippe suggested to leave it in the paper, I agree. It was previously stated for semilocal rings.}

We can now extend \eqref{splica-3} to the non-split case. 

\begin{cor} \label{cor_primitive_etale}
Let $R$ be an LG ring and let $E\in \Ralg$ be a finite \'etale $R$--algebra such that $\rank E < |k|$ for every residue field $k$ of $R$. Then $E$ is unit-generated. 
\end{cor}

\begin{proof} By Corollary~\ref{loog}\ref{loogb} and the remark above, it suffices to consider the case that $R=k$ is a field. By \ref{plr}\eqref{plr-f}, we can assume that $k$ is a finite field, say $k=\FF_q$. Also, we know from \ref{fea}\eqref{fea-a} that $E = K_1 \times \cdots \times K_r$,  where every $K_i$ is a finite separable extension of $k$. By Lemma~\ref{lem_generatornew} it is enough to consider the case that the $K_i$ are pairwise isomorphic, say $K_i \cong \FF_{q^n}$, thus $\rank E = rn$. In this case, the result follows from \ref{plr}\eqref{plr-bb} and the fact that the number of irreducible monic polynomials in $\FF_q[X]$ is bounded by $q/n$, or by specializing the formula \cite[Thm.~1.2]{FRS}. 
\end{proof}
\enew
 
In the remainder of this section we will consider special types of primitive elements. 

\comments{(2026-02-04) Changed $h^{u,m}$ to $h^{u,n}$ for better readability in  } 

\begin{cor}[The morphism $h^{u,n} \co \uW(E) \to \uW(E)$, $E$ finite \'etale] \label{corhumLG} Let $E$ be a finite \'etale $R$--algebra of positive rank. Also, let $u\in E\ti$ and let $n$ be a positive integer. 
\sm

\begin{inparaenum}[\rm (a)]
\item\label{corhumLGa} The morphism of schemes $h^{u,n} \co \uW(E) \to \uW(E)$, given on $T$--points, $T\in \Ralg$, by  
    \[ 
           h^{u,n}(T) \co E\ot_R T \to E\ot_R T, \quad e \mapsto (e \ot 1_T) a^n,
      \] 
    is finite locally free and faithfully flat. \sm 
    
\item \label{corhumLGb} Recall \eqref{prschb2}: $\uPrim(E) = \uW(E)_\dis$. The principal open subscheme 
    \[ 
          \uW(E)_{\dis \circ h^{u,n}} = (h^{u,n})\me \big(\uPrim (E)\big)
    \] 
    of $\uW(E)$ represents the $R$--functor, assigning to $T\in \Ralg$ the set 
    \[\{ v\in E_T : (u \ot 1_T)v^n \in \Prim(E_T)\}.
     \] It is universally schematically dense in $\uW(E)$.  
\sm 

\item \label{corhumLGc} The $R$--scheme
\[ (h^{u,n})\me \big(\uPrim (E)\ti \big) = (h^{u,n})\me \big(\uPrim (E)\cap \uGL_1(E) \big)
\]
is a \new principal open \enew subscheme of $\uW(E)$. It represents the $R$--functor which assigns to $T\in \Ralg$ the set 
\[
   \{ v\in (E_T)\ti : (u \ot 1_T)v^n \in \Prim(E_T)\ti\}.
 \] 
\end{inparaenum}
\end{cor}

\begin{proof}
Because $h^{u,n} = h^{u,1} \circ h^{1_E, n}$ and $h^{u,1}$ is an automorphism of $\uW(E)$, it suffices to consider the example $u=1_E$. \sm 

\eqref{corhumLGa} Since $h^{1_E, n}$ respects the rank decomposition of $E$, we can also assume that $E$ has constant rank $d\in \NN_+$. Furthermore, the property of a morphism of schemes to be faithfully flat or finite locally free allows flat descent \cite[Tags 02KV, 02L2, 02VO]{St}. 
\lv{
$\Spec(B) \to \Spec(C)$ is faithfully flat iff $C \to B$ makes $B$ a faithfully flat $C$--module \cite[I, \S2, 2.4]{DG}, which in turn can be checked after a faithfully flat extension by \cite[I, \S3.3, Prop.6]{BAC}.}
It is therefore enough to consider the case of a split $E$, i.e., $E \cong R^d$. In this case, the coordinate ring of $\uW(E)$ is the polynomial ring $B = R[X_1, \ldots, X_d]$ in $d$ variables and $h^{1_E, n}$ corresponds to the ring map $B \to B$, $X_i \mapsto X_i^n$, $i=1, \ldots , d$. It makes $B$ a non-zero free, hence faithfully projective $B$--module. 
\lv{
Idea for $n=1$: Every polynomial $f$ can be uniquely written as $f(X) = \sum_{i=0}^{n-1} p_i(X^n) X^i$ with $p_i(X) \in R[X]$. Write $n\in \NN$ in the form $n = \ell n + i$, $0 \le i < n$, and so the $n$th term becomes $r_n X^n = r_{\ell n + i} (X^n)^\ell X^i$, now put $p_i(X) = \sum_\ell r_{\ell n + i } X^\ell$. 

A scheme morphism $\Spec(B) \to \Spec(E)$ is faithfully flat iff $B$ is a faithfully flat E-algebra
}
This implies our claim. \sm

The first part of \eqref{corhumLGb} is obvious, the second follows from Corollary~\ref{usd-d} since $\dis\circ h^{1, n}$ is a non-zero polynomial over any $R$--field. 
\sm

\eqref{corhumLGc} By \ref{prischinv}, $\uPrim(E)\ti$ is a principal open subscheme of $\uW(E)$. Hence so is its inverse image $(h^{u,n})\me \big(\uPrim (E)\ti \big)$ by  \cite[I, (1.2.2)]{EGA-neu}. 
\end{proof}
\sm

If $A=E$ is finite \'etale, we can replace the scheme $\uPrim(A)$ in the proof of Corollary~\ref{loog}  by the scheme $\uW(E)_{\dis \circ h^{u,n}}$ of Corollary~\ref{corhumLG} and in this way obtain a result analogous to Lemma~\ref{loog} by applying again Corollary~\ref{fibn}. We leave the obvious formulation to the reader and only formulate the case of an LG base ring for later use.


\begin{lem}\label{hunco} Let $R$ be an LG ring and let $E$ be a finite \'etale $R$--algebra. We fix $u\in E\ti$ and $n\in \NN_+$. Then the following conditions \ref{huncoi}--\ref{huncoiii} are equivalent: 

\begin{enumerate}[label={\rm (\roman*)}]
 \item\label{huncoi}
  There exists $v\in E$ such that $uv^n \in \Prim_R(E)$. 
  
 \item For every maximal ideal $\gm \ideal R$ there exists $y_\gm \in E \ot_R R_\gm$ satisfying $(u\ot 1_{R_\gm}) \, y^n_\gm \in \Prim_{R_\gm} (E \ot_R R_\gm) $.
     
 \item \label{huncoiii} For every maximal ideal $\gm \ideal R$ there exists $z_\gm \in E \ot (R/\gm)$ satisfying $( u \ot 1_{R/\gm}) \, z_\gm^n \in \Prim_{R/\gm}\big(E \ot_R (R/\gm)\big)$ 
\end{enumerate} 
Also, the following conditions \ref{huncoI}--\ref{huncoIII} are equivalent:
\begin{enumerate}[label={\rm (\Roman*)}]
 \item\label{huncoI}
  There exists $v\in E\ti$ such that $uv^n \in \Prim_R(E)\ti$. 
  
 \item For every maximal ideal $\gm \ideal R$ there exists $y_\gm \in (E \ot_R R_\gm)\ti$ satisfying $(u\ot 1_{R_\gm}) \, y^n_\gm \in \Prim_{R_\gm} (E \ot_R R_\gm) \ti$.
     
 \item \label{huncoIII} For every maximal ideal $\gm \ideal R$ there exists $z_\gm \in \big( E \ot (R/\gm)\big){}\ti$ satisfying $( u \ot 1_{R/\gm}) \, z_\gm^n \in \Prim_{R/\gm}\big(E \ot_R (R/\gm)\big){}\ti$ 

\end{enumerate}
Furthermore, if $|R/\gm| = \infty$ for all maximal ideals of $R$, i.e., $R$ satisfies the primitive criterion {\rm \ref{revLGG}}, then all these conditions are fulfilled. 
\end{lem}

\begin{proof} Since $U=(h^{u,n})\me \big(\uPrim (E)\big)$ and $U'=(h^{u,n})\me \big(\uPrim (E)\ti\big)$ are quasi-compact open subschemes of $\uW(E)$, the equivalence of \ref{huncoi}--\ref{huncoiii} and of \ref{huncoI}--\ref{huncoIII} is just another application of Corollary~\ref{fibn}, and so is the last claim. 
\end{proof}

\comments{(2026-02-14) Proposition~\ref{ndnLG} is the old Proposition 2.5, label `ndn`, for $R$ a semilocal ring. I did not find where we used it in the old version, and it is not clear if we will use it later. For now, I will leave it here.} 

\begin{prop}\label{ndnLG} Let $R$ be an LG ring, let $R'\in \Ralg$ be finite \'etale over $R$, let $A'\in \Rpalg$ be finite \'etale over $R'$ {\em of constant rank}, and let $u'\in R'{}\ti$. Suppose that 
\begin{enumerate}[label={\rm (\roman*)}]
\item\label{ndnLGi} for every maximal ideal $\gm \ideal R$ with $k = R/\gm$ {\em finite\/}, the $k$--algebra $R'\ot_R k$ is unit-generated. 
\end{enumerate}  
Then there exists $a'\in A'$ such that $\rmN_{A'/R'}(a') u' \in \Prim_R(R')\ti$. 
\end{prop}

\begin{proof} Let $A = \frR_{R'/R}(A')$ be the Weil restriction of $A'$ (= the $R$--algebra obtained from $A'$ by restriction scalars to $R$). By \ref{fea}\eqref{fea-trans}, $A$ is a finite \'etale $R$--algebra. 
In part \eqref{ndnLGI} of the proof we will construct a quasi-compact open subscheme $U$ of the $R$--scheme $\uW(A)$ such that the proposition's claim means $U(R) \ne \emptyset$; in part \eqref{ndnLGII} we will then show that $U(R/\gm) \ne \emptyset $ for every maximal ideal $\gm$ of $R$, so that we can conclude with the characterization \ref{prop_baire}\ref{prop_baire-a} of LG rings. 
\sm 

\begin{inparaenum}[(I)] \item \label{ndnLGI} In this part of the proof it suffices to assume that $R'\in \Ralg$ and $A'\in \Rpalg$ are finite locally free as $R$-- and $R'$--modules respectively. Recall \ref{trno}\eqref{trno-ds} the norm morphism of $R'$--schemes, 
\[ 
  \rmN_{A'/R'} \co \uW_{R'}(A') \to \GG_{a,R'} = \uW_{R'}(R'). 
\]
Applying the Weil restriction functor $\frR = \frR_{R'/R}(\cdot)$, we get a morphism of $R$--schemes 
\[ 
\frR(\rmN_{A'/R'}) \co \frR\big( \uW_{R'}(A')\big) = \uW(A) \to \frR( \GG_{a, R'}) = \uW(R').
\]
The invertible element $u'\in R'$ induces an automorphism $f_{u'}$ of the $R$--scheme $\uW(R')$, given on $T\in \Ralg$ by  $R'\ot_R T \to R'\ot_R T$, $x \mapsto (u'\ot 1_T)\me x$.  By Corollary~\ref{prischinv}, the $R$--scheme $\uPrim(R')\ti$ representing invertible primitive elements is a quasi-compact open subscheme of $\uW(R')$. Hence so is $f_{u'} \big( \uPrim (R')\ti \big)$. We can now define the $R$--scheme $U$ as the fibre product 
\[   \vcenter{ \xymatrix@C=50pt{\ar @{} [dr] |{\small\qed} %
U \ar[r] \ar[d] & \uW(A) \ar[d]^{\frR(\rmN_{A'/R'})} 
 \\ 
f_{u'}\big( \uPrim(R')\ti \big) \ar[r] & \uW(R') } }   
 \]
By base change for open immersion and quasi-compact immersions, 
$U$ is a quasi-compact open subscheme of $\uW(A)$. Its $T$--points, $T\in \Ralg$ with $T'= R'\ot_R T$, are 
\begin{align*}
  U(T) & = \big\{ a'\in A'\ot_{R'}T': 
  \rmN_{A'\ot_{R'} T'/ T'} (a') \, (u'\ot 1_T) \in \Prim_T(T')\ti \big\}. 
\end{align*}
In particular, the claim of the proposition means $U(R) \ne \emptyset$. 
\sm 

\item\label{ndnLGII} Under the assumptions of the propositions we will show that $U(R/\gm) \ne \emptyset$ for every maximal  ideal $\gm \ideal R$. We define rings by the diagram below,  
    \[ \xymatrix{
       R \ar[r] \ar[d] & R' \ar[r] \ar[d] & A'\ar[d]  \\
       k=R/\gm \ar[r]  &E = R'\ot_R k \ar[r] &A'_E = A'\ot_{R'} E 
    } \] 
    and put $u = u'\ot 1_k \in E$. The description of $U(T)$ above for $T=k$ says that our claim is that  
\begin{equation}\label{ndnLG1}
  \text{there exists $a' \in (A'_E)\ti$ such that $\rmN_{A'_E/E}(a')  u \in \Prim_k (E)\ti$. }
\end{equation}
We distinguish the two cases $|k| = \infty$ and $|k|< \infty$. \sm 

$|k| = \infty :$ Since the $E$--algebra $A'_E$ is finite locally free of fixed rank, say of rank $n\in \NN_+$, we can apply Lemma~\ref{hunco} with $R$ replaced by $k$ and get that there exists $v \in E\ti$ such that $v^n u \in \Prim_k(E)\ti$. Since $v^n = \rmN_{A'_E/E}(v)$, we are done.  

$|k|< \infty$: The $k$--algebra $E$ is finite \'etale. Hence $E = L_1 \times \cdots \times L_m$ is a product of field extensions of $k$, necessarily of finite degree and therefore finite fields, \ref{fea}\eqref{fea-a}. The $E$--algebra $A'_E$ decomposes correspondingly, 
\[ A'_E = A'_1 \times \cdots \times A'_m,   
\]
where each $A'_i$ is a finite \'etale $L_i$--algebra, \ref{fea}\eqref{fea-d}. 
By \ref{trno}\eqref{trno-c}, $\rmN_{A'_E/E} = \prod_i \rmN_{A'_i/L_i}$ and by \ref{trno}\eqref{trno-f}, every $\rmN_{A'_i/L_i}$ is surjective. Hence $\rmN_{A'_E/E}$ is surjective too. By assumption \ref{ndnLGi}, the $k$--algebra $E$ is unit-generated, say by $z_E \in E\ti$. By \ref{trno}\eqref{trno-a} and surjectivity of the norm, there exists $a' \in A'_E$ such that $z_E u \me = \rmN_{A'_E/E}(a')$,  i.e., \eqref{ndnLG1} holds. 
\end{inparaenum}
\end{proof}
\sm

\textbf{Remark.} The assumption \ref{ndnLG}\ref{ndnLGi} is fulfilled if $R'\ot_R \ka = R'/\m R'$ is a field, hence a separable extension field of $\ka$ and therefore unit-generated. For example, this holds for the algebras $R_j$ in the setting of  Proposition~\ref{lem_gonflement}\eqref{semilocali-a}.

 \newpage

\section{Algebras over semilocal rings}\label{sec:semilocal}

The following result is folklore. We include a proof since we could not find a reference.

\comments{(2020-05) It is surprising that we could not find a reference for \ref{lem_cyclotomic}. Unfortunately, I do not know algebraic number theory. By Kummer theory, all that is needed in the case of a number field $F^\flat$ is that the group $F^\flat/(F^\flat )^n$ has an element of order $n$, see \cite[V, \S11.8, Ex.~3 and 4]{BA5}.}

 \begin{lem} \label{lem_cyclotomic}  Let $F$ be a field which is finitely generated over its prime field. Then for every $n\in \NN_+$ there exists a cyclic field extension of $F$ of degree $n$.
 \end{lem}

 \begin{proof} Let $P$ be the prime field of $F$ and let $F^\flat$ be the algebraic closure of $P$ in $F$. We will first show that there exists a cyclic extension $E^\flat/F^\flat$ of degree $n$.

 By \cite[V, \S14.7, cor.~1 of prop.~17]{BA5}, the field extension $F^\flat/P$ is finite, hence $F^\flat$ is a finite field or a number field. In the first case we take as $E^\flat/F^\flat$ the unique field  extension of degree $n$, well-known to exist and to be a cyclic extension. Let now $F^\flat$ be a number field.
 Let $p_1$ and $p_2$ be distinct prime numbers and let $E_1$ and $E_2$ be cyclic extensions of degree $p_1^{r_1}$ and $p_2^{r_2}$ respectively. The compositum of $E_1$ and $E_2$ of $F^\flat$ is a cyclic extension of degree $p_1^{r_1} p_2^{r_2}$ by \cite[V, \S10.8, Thm.~5]{BA5}. Hence, without loss of generality, it suffices to consider $n=p^r$ for some prime $p$. Let $F^\flat( \mu_{p^\infty} )$ be the extension of $F^\flat$, obtained by adjoining all roots of unity of $p$th--power order. It is known (\cite[XI.1]{NSW}) that $F^\flat( \mu_{p^\infty} )/F^\flat$ is a Galois extension whose Galois group is isomorphic to $\ZZ_p \times \Gamma$ where $\Gamma$ is a finite abelian group. The fixed points under the subgroup
$p^r \ZZ_p \times \Gamma$ is therefore a cyclic extension  $E^\flat/F^\flat$ of degree $p^r$.

Finally, we extend the cyclic extension $E^\flat/F^\flat$ to $F$: since
$F^\flat$ is algebraically closed in $F$, the tensor product $E^\flat \otimes_{F^\flat} F$ is a field \cite[V, \S17.5, Prop.~9]{BA5}. 
It is a cyclic extension of $F$ of degree $n=p^r$ by \cite[V, \S10.8, Thm.~5]{BA5}.
\end{proof}

\begin{lem}[{\cite[10.1.1]{Ford}} for local rings] \label{gonf-pre}
Let $R$ be a semilocal ring and let $S\in \Ralg$ be a finite $R$--algebra. Then $S$ is a semilocal ring and an integral extension of $R$.

Moreover, let $\m \ideal R$ be a maximal ideal of $R$ and let\/ $\frn_1, \ldots, \frn_a$ be the maximal ideals of $S$ lying over\/ $\m$. Put $\wtl S = S/\m S$ and denote by $\wtl \frn_i \ideal \wtl S$ the image of $\frn_i$ in $\wtl S$. Then the canonical map
\begin{equation}
  \label{gonf-pre1} \wtl S \; \simlgr \;  \wtl S_{\wtl \frn_1} \times \cdots \times \wtl S_{\wtl \frn_a}
\end{equation}
is an isomorphism of $S$--algebras. Each $\wtl S_{\wtl \frn_i}$, $1\le i \le a$, is a local Artinian ring. In particular, if $\wtl S \cong S \ot_R (R/\m)$ is a field, then $a=1$.
\end{lem}

\begin{proof} 
By \ref{slr}\eqref{slr-b}, $S$ is a semilocal ring. Also, $S/R$ is an integral extension, see e.g.\ \cite[V, \S1.1, Prop.~1]{BAC2}. It follows that for any prime ideal $\p \ideal R$ there exists a prime ideal $\frq\ideal S$ lying over $\p$, i.e., $\p = \frq \cap R$, and for such a pair $\p$ is maximal if and only if $\frq$ is maximal (\cite[V, \S2.1, Thm.~1 and Prop.~1]{BAC2}). Since $S$ is semilocal, there are only finite many ideals lying over $\m$.

Because $\wtl S$ is a finitely generated $(R/\m)$--vector space, it is an Artinian ring. 
Its maximal ideal are precisely the ideals $\wtl\frn_1, \ldots, \wtl \frn_a$. The isomorphism \eqref{gonf-pre1} is therefore part of the structure theorems for Artinian rings, see e.g.\ \cite[Thm.~3.2.11]{Ford}.
\end{proof}

\subsection{Semilocalization} \label{semilocali}
We remind  the reader of the semilocalization process.
 Let $A$ be a commutative (unital) ring and let
 $\gp_1, \dots \gp_c$ be prime ideals of $A$. The subset
 $S= \textstyle \bigcap_{i=1,\dots, c} (A \setminus \gp_i)=
 A \setminus \bigcup_{i=1,\dots,c} \gp_i$ is multiplicative
 and the localization $S\me A$ is a semilocal ring \cite[II, \S 3.5, Prop.~17]{BAC}.
 If $\frq_1, \ldots, \frq_b$ are the distinct maximal elements of the $\p_1, \ldots, \p_c$ with respect to inclusion, the maximal ideals of $S\me A$ are the $S\me \frq_i$, $1\le i \le b$, and these maximal ideals are pairwise distinct.  We call a ring $B$ a {\em semilocalization of $A$\/}, if there are prime ideal $\p_1, \dots, \p_s$ of $A$ such that $B\cong T\me A$ with $T=\bigcap_{i=j,\ldots, s} (A \setminus \p_j)$.

If $R$ is a semilocal ring with maximal ideals $\m_1, \ldots, \m_c$, then
 $R\ti = R\setminus S$ with $S=\bigcup_{i=1,\dots,c} \m_i$ by \ref{nak}\eqref{nak-a}. Hence the canonical map $R \simlgr S\me R$ is an isomorphism. In other words, $R$ coincides with the semilocalization at its maximal ideals.

\comments{(2026-04-01) Is there a result combining or generalizing the two Propositions~\ref{lem_gonflement} and Propositions~\ref{bfp-prop}?  

We use Prop.~\ref{lem_gonflement} in the proof of Prop.~\ref{prop_norm_SB} (the result that Philippe needs). We use Prop.~\ref{bfp-prop} in the proofs of \ref{cor-knex} and \ref{prop_norm_quad}.  
 }

\begin{prop} \label{lem_gonflement} Let $R$ be a semilocal ring and let
$\ell$ be a prime number. We denote by $\kappa_1, \dots, \kappa_c$
 the residue fields of the maximal ideals $\m_1, \ldots , \m_c$ of $R$.
\sm

\begin{inparaenum}[\rm (a)]
\item \label{semilocali-a} Then there exists a sequence
\[
R=R_0 \subset R_1 \subset R_2 \subset \cdots
\]
of semilocal rings such that for every $j \ge 0$ the following holds.
\end{inparaenum}

 \begin{enumerate}[label={\rm (\roman*)}]
  \item \label{lem_gonflement_i} $R_{j+1}$ is finite \'etale
  of degree $\ell$ over $R_j$, equivalently, $R_j$ is finite \'etale of  degree $\ell^j$ over $R$;

 \item \label{lem_gonflement_ii} $R_{j+1}=R_j[u_j]$  for some $u_j \in R_j^\times$;

 \item \label{lem_gonflement_iii} if $\ka_i$ is finite, $R_j \otimes_{R} \kappa_i$ is a field of degree $\ell^j$ over $\kappa_i$.
\end{enumerate}
\sm

\begin{inparaenum}[\rm (a)] \setcounter{enumi}{1}
\noindent \item\label{semilocali-b} Furthermore, if  $R$ is the semilocalization of a  finitely generated (commutative unital)  $\ZZ$--algebra, we can additionally require the following three conditions:
\end{inparaenum}
\begin{enumerate} [label={\rm (\roman*)}]\setcounter{enumi}{3}
  \item \label{lem_gonflement_iiia}
  For every $n\in \NN$ and every $i$, $1\le i \le c$, there exists exactly one maximal ideal of $R_n$ lying over $\m_i$, i.e., the rings $R_n$, $n\in \NN$, have exactly $c$ maximal ideals;

  \item \label{lem_gonflement_iv} the algebras
  $R_j \otimes_{R} \kappa_i$'s are fields  of dimension $\ell^j$ over $\kappa_i$ for all $i,j$.

  \item \label{lem_gonflement_v} $R_\infty = \limind R_i$ is a semilocal ring
  with $c$ maximal ideals, all of whose residue fields are infinite.
 \end{enumerate}
\end{prop}

\begin{proof} \eqref{semilocali-a} 
We construct the \'etale extensions $R_{j+1}$ satisfying \ref{lem_gonflement_i}--\ref{lem_gonflement_iii} inductively as $R_{j+1} = R[X]/(P_j)$ where $P_j$ is a suitable separable polynomial of degree $\ell^j$. Being finitely generated as $R$--modules, these extensions are semilocal rings by Lemma~\ref{gonf-pre}. Let $\bar R = R/\Jac(R) = \ka_1 \times \cdots \times \ka_c$. Since $R[X] \to \bar R[X] = \ka_1[X] \times \cdots \times \ka_c[X]$ is surjective, it is in view of Lemma~\ref{genolemLG} enough to construct appropriate separable polynomials $P_{ij}\in \ka_i[X]$ for $1\le i \le c$, $j\in \NN_+$, or, equivalently, appropriate \'etale extensions $K_{ij}$ of $\ka_i$. Let us first consider the case of a finite $\ka_i$.  In this case we let $K_{ij}$ be the (up to isomorphism) unique field extension of $\ka_i$ of degree $\ell^j$. It is a separable extension of $\ka_i$, hence a simple extension, i.e., one-generated by an obviously invertible element. Taking the $K_{ij}$ as subfields of a fixed algebraic closure of $\ka_i$, we indeed have $K_{ij} \subset K_{i, j+1}$. If $\ka_i$ is infinite, let $P_{ij}$ be a separable polynomial of degree $\ell^j$ with $P_{ij}(0) \ne 0$, for example a split polynomial. Such a polynomial exists since $\ka_i$ is infinite, and we can arrange them so that $P_{ij}$ divides $P_{i, j+1}$, i.e., $K_{ij} \subset K_{i,j+1}$. \sm

\eqref{semilocali-b} We now assume that $R$ is the semilocalization of a finitely generated $\ZZ$--algebra. We will modify the choice of the fields $K_{ij}$. The residue fields $\ka_1, \ldots, \ka_c$ are finitely generated over their prime fields. We let $K_{1i}$, $i=1, \ldots, c$, be a cyclic extension of $\ka_i$ of degree $\ell$, whose existence is guaranteed by Lemma~\ref{lem_cyclotomic}. As a separable extension of $\ka_i$, it is a simple extension. We can therefore construct $R_1$ as in the first part of the proof. It then follows from Lemma~\ref{gonf-pre} that there lies exactly one maximal ideal $\m_{i1}$ of $R_1$ over every maximal ideal $\m_i$ of $R$. Applying \cite[V, \S2.1, Cor.~3]{BAC2},
we get $\m_{i1} = \m_i R_1$. The residue fields of $R_1$ are finitely generated fields over their prime field, so that we can continue the process. This provides a tower satisfying properties \ref{lem_gonflement_iiia} and \ref{lem_gonflement_iv}.

The ring $R_\infty$ is an integral extension of $R$. Consequently, any maximal ideal $\m_\infty\ideal R_\infty$ lies over a unique maximal ideal $\m_i\ideal R$. If $\m_\infty'$ is a maximal ideal of $R_\infty$ lying over $\m_i$ and distinct from $\m_\infty$, there exists $n\in \NN$ such that $\m_\infty \cap R_n \ne \m'_\infty \cap R_n$. Hence $\m_\infty \cap R_n$ and $\m'_\infty\cap R_n$ are two distinct maximal ideals of $R_n$ lying over the same $\m_i$, contradicting \ref{lem_gonflement_iiia}. Thus $R_\infty$ has exactly $c$ maximal ideals, in particular $R_\infty$ is a semilocal ring.
Finally, the field  $R_\infty/ \m_i$ is the colimit of the $R_j/ \m_i$, so is
 infinite by taking into account \ref{lem_gonflement_iii}. We have therefore also established \ref{lem_gonflement_v}.
\end{proof}

\begin{cor} \label{cor_gonflement} Let $R$ be a semilocal ring and
let $S$ be an  \'etale extension of $R$ of degree $d$.
Then for each prime $\ell$ there exists a finite \'etale extension
$R'$ of $R$ such that $R'$ is a tower of unit-generated \'etale extensions of
degree $\ell$ and such that $S \otimes_R R'$ is a unit-generated \'etale extension of $R'$. \end{cor}

\begin{proof} We choose a tower $R= R_0 \subset R_1 \subset \cdots$ as in
Proposition~\ref{lem_gonflement}  and put $R'=R_d$. It then remains to prove that $S' =S \ot_R R'$ is unit-generated as $R'$--algebra. Since $R'$ is semilocal and $S'$ is an \'etale extension of $R'$ of degree $d$, we are in the setting of Corollary~\ref{cor_primitive_etale}. Thus, it suffices to show
\begin{equation} \label{cor_gonflement-1}
d < | R'/\m'  | \quad \text{for every maximal ideal $\m' \ideal R'$.}
\end{equation}
Because $R'/R$ is an integral extension, $\m = R\cap \m'$ is a maximal ideal of $R$ and $\ka' = R'/\m'$ is a field extension of $\ka = R/\m$. Hence \eqref{cor_gonflement-1} is clear if $\ka$ is an infinite field. But otherwise, by Proposition~\ref{lem_gonflement}\ref{lem_gonflement_iii} and \cite[V, \S2.1, Cor.~3]{BAC2}, $\m'$ is the only ideal lying over $\m$, hence $\m' = R'\m$ and so $\ka' = R'/\m' \cong R'\ot_R \ka$ is a field extension of $\ka$ of degree $\ell^d$. Putting $n = \ell^d$ and observing $n\ge 2^d > d$, we have $|\ka'| = |\ka|^n \ge 2^n >d$.
 \end{proof}

\comments{(2026-04-21) Moved the [BFP]--Proposition \ref{bfp-prop} to the end, so that it is clear that the results of this section do not depend on \ref{bfp-prop}}

\sm

For the next lemma we recall that $\Prim_R(A)$ denotes the set of primitive elements of an $R$--algebra $A$.

\begin{lem} \label{lem-btwo} Let $k$ be a finite field, let $L/k$ be a field extension of finite degree and let $a\in L\ti$. Then there exists $b \in L\ti$ such that $ab^2\in \Prim_k(L)$.
  \end{lem}

\begin{proof} Let $|k| = q$  and $L|= q^m$ for some $m\in \NN_+$. Since $\Prim_k(k) = k\ti$, we can assume $m>1$. We consider the map $f_a \co L\ti \to L\ti$, $b \mapsto ab^2$. Then $|\Ima(f_a)| = |\Ima  (f_1)| = |L\ti{}^2|$, where
  \[ | L\ti{}^2 | = \begin{cases}
     q^m -1 , & \Char(k) = 2, \\ (q^m-1)/2,  &\Char(k) \ne 2.
  \end{cases}
 \]
Let $\rmN = L\ti \setminus \Prim_k(L)$, the set of non-generators of $L$. Observe that $\ell\in L$ generates $L$ as $k$--algebra if and only if $\ell$ generates $L$ as a field. 
Hence, $\ell \in \rmN \iff \ell$ lies in a maximal subfield of $L$. To describe these, let
\[ m = p_1^{e_1} \cdots p_c^{e_c}, \qquad (\text{prime power decomposition})
\]
where the $p_i$ are pairwise distinct prime numbers and where $e_i \in \NN_+$. Let $m_i = m/p_i$.  For $1\le i \le c$ there exists a unique subextension $L_i/k$ of $L/k$ with $|L_i| = q^{m_i}$; the $L_i/k_i$ are precisely the maximal subextensions $L/k$. Hence $\rmN = \textstyle \bigcup_{i=1, \ldots c} \, L_i$.
Our claim is that $\Ima(f_a) \not\subset \rmN$. This holds whenever
\[ |\Ima (f_a)| > |\rmN| = |\textstyle \bigcup_{i=1, \ldots c} \, L_i |
\tag{*}\]
It is easily seen that (*) is satisfies if $\Char(k)= 2$. In this case, $\Ima(f_a) = L\ti{}^2$, while $\rmN \subsetneq L\ti{}^2$. In the following we suppose $\Char (k) \ne 2$, thus (*) becomes
\[  q^m -1 > 2|\rmN|= |\textstyle \bigcup_{i=1, \ldots c} \, L_i |
\]
A very rough estimate of $|\rmN|$ is $|\rmN| \le \sum_i |L_i| = \sum_i q^{m_i}$. Hence it suffices to show
\[ \tag{**} q^m - 1 > 2 \, \textstyle \sum_{i=1, \ldots, c} \, q^{m_i}. \]
{\em Reduction to $m=p_1 \cdots p_c$}: Let
\[ \wtl m = p_1^{e_1-1} \cdots p_c^{e_c-1} \ge 1.  \]
Then (**) is true iff it is true after division by $q^{\wtl m}\ge 3$, i.e.,
\[
  \frac{q^m}{q^{\wtl m}} - \frac{1}{q^{\wtl m}}
     > 2 \, \textstyle \sum_{i=1, \ldots, c} \, \frac{q^{m_i}}{q^{\wtl m}}.
\]
Since $ \frac{q^m}{q^{\wtl m}} - \frac{1}{q^{\wtl m}} > \frac{q^m}{q^{\wtl m}} -1$, it suffices to show (**) for $m=p_1 \cdots p_c$, i.e.,
\[
\tag{***} q^{p_1 \cdots p_c} - 1 > 2 \, \textstyle \sum_{i=1, \ldots, c} \,
               q^{p_1 \cdots \widehat{p_i} \cdots p_c} \]
where $\wdh{p_i}$ means that the factor $p_i$ has to be deleted. \sm

{\em Reduction to $c>1$}: If $c=1$, then (***) becomes $q(q^{p_1-1} -2) >1$, which is true because $q\ge 3$, $q^{p_1-1} -2 \ge q-2 \ge 1$, so $q(q^{p_1-1} -2) \ge 3 \cdot 1 >1$. \sm

We can now prove (***) in general: It is no harm to assume $2\le p_1 < p_2 < \cdots < p_c$, implying $p_2 \ge 3$ and $m_1> \cdots > m_c$. We will also use the obvious inequality
$   3^{p_2 \cdots p_c} > 2c$, easily established by induction. Then
\begin{align*}
  q^{m-m_1} &= q^{(p_1-1)(p_2 \cdots p_c)} \ge  q^{p_2 \cdots p_c} \ge
       3^{p_2 \cdots p_c} > 2c, \quad\text{hence} \\
 1&< q^{m_1}(q^{m-m_1} - 2c) = q^m - 2c q^{m_1}
 \\ &  < q^m - 2  \textstyle \sum_{i=1}^c \,
               q^{m_i} \end{align*}
where in the last inequality we used that $q^{m_1} > q^{m_i}$ for $i=2, \ldots, c$.
\end{proof}

\medskip

\begin{cor} \label{cor_pigeon} Let $k$ be a field, let $A$ be a multiplicity-free
\'etale $k$--algebra, and let $a \in A\ti$. Then there exists $b \in A\ti$
such that $ab^2\in \Prim_k(A)$.
\end{cor}

\begin{proof}
  We know $A=K_1 \times \cdots \times K_n$,  where $K_i/k$ are separable field extensions of finite degree. Since $A$ is multiplicity-free, Lemma~\ref{lem_generatornew} says $\Prim_k(A) =  \Prim_k(K_1) \times \cdots \times \Prim_k(K_n)$. Hence, it suffices to prove the claim for each factor $K_i/k$. For them the claim follows from Lemma~\ref{lem-btwo} if $k$ is finite, and from Lemma~\ref{hunco} 
  if $k$ is infinite. \end{proof}
\ms

Our next aim is to prove a version of Corollary~\ref{cor_pigeon} in case $A$ is not multiplicity-free, see Proposition~\ref{prop_tiroir}.

\begin{lem}\label{lem_multiplicity} Let $k$ be a finite field, and let $d \geq 1$ be  an integer. Then, for every integer $n \geq d^2+d+1$ and for every \'etale $k$--algebra $A$ of degree $\leq d$, there exists an \'etale $A$--algebra $B$ satisfying
the following:
\begin{enumerate}[label={\rm (\roman*)}]
 \item \label{lem_multiplicity_i} $B$ is a free $A$--module of rank $n$;
 \item \label{lem_multiplicity_ii} $B$ is multiplicity free as \'etale  $k$-algebra.
\end{enumerate}
\end{lem}

\begin{proof} We know $k=\FF_q$, $q=p^m$ for some prime $p$. Let $n$ and $A$ be as in the statement of the lemma. We can decompose $A= \FF_{q^{r_1}} \times \dots \times \FF_{q^{r_c}}$ with $r_1 \leq r_2 \leq \dots \leq r_c$. Note $r_1+ \dots + r_c = \dim_k A \leq d$ and $c\le d$ since all $r_i \ge 1$. We claim:

\begin{enumerate} \item the integers $r_1, r_1(n-1), 2 r_2, r_2(n-2), \dots, c r_c , r_c(n- c )$ are pairwise distinct.
\end{enumerate}
Indeed, by construction we have $ r_1  < 2 r_2  < \dots < c r_c \leq d^2$ and
$r_i (n-i) \geq (n-i) \geq n - d > d^2$ for $i=1,\dots,c$. In particular, $r_i (n-i) \not = r_j j$ for $i,j=1, \dots, c$. If  $r_i (n-i) = r_j (n - j)$ for $i<j$, then $(r_j-r_i) n = r_j j - r_i i >0$ so that $r_i < r_j$ and $n \leq
d^2$, which contradicts the assumption. The claim is established.

We define the \'etale $A$--algebra
\[
B= \prod\limits_{i=1}^c \big( \FF_{q^{r_i\, i }} \times  \FF_{q^{r_i (n-i)}} \Big),
\]
letting the $i$th factor $\FF_{q^{r_i}}$ of $A$ act diagonally on the $i$th factor of $B$.
Since  $i+n-i=n$ for $i=1, \dots, c$, $B$ is a free $A$-module of rank $n$.
Finally $B$ is a multiplicity free $k$-algebra according to the claim above.
\end{proof}

\comments{(2023-01) Deleted the assumption ``Let $e$ be a positive integer and let $q$ be a prime number.'' in \ref{prop_tiroir}. It is not used in the statement }

\begin{prop} \label{prop_tiroir}
Let $R$ be a semilocal ring, and let $S$ be a finite \'etale $R$--algebra which is locally free of rank $d$. Then, for each $a \in S\ti$ and for each $n \geq d^2+d+1$, there  exists a finite \'etale $S$-algebra $S'$ satisfying the following conditions:

\begin{enumerate}[label={\rm (\roman*)}]
  \item $S'$ is a locally free $S$--module of rank $n$;

  \item The $R$-algebra $S'$ admits a unit-generator  $a b^2$ for some
  $b \in (S')^\times$.
\end{enumerate}
\end{prop}

\begin{proof} We are given an integer $n \geq d^2+d+1$.
\sm

\noindent{\it First case: $R$ is a infinite field $k$.}
This case follows from \ref{hunco}. \sm

\noindent{\it Second  case: $R$ is a finite field $k$.}
 Lemma \ref{lem_multiplicity}  shows that there exists
an \'etale $S$--algebra $S'$ such that $S'$ is a free $S$--module of rank $n$ and
is  multiplicity free as $\FF_p$--algebra. Corollary \ref{cor_pigeon}
shows that $S'$  admits a unit-generator over $\FF_p$
  of the shape $ab^2$. A fortiori, this is a generator of the $S$--algebra $S'$.

\sm

\noindent{\it General  case.} Let $I=\Jac(R)$ be the Jacobson radical of the semilocal ring $R$ and decompose $R/I= \kappa_1 \times \kappa_2 \times \cdots \times \kappa_c$ where the $\kappa_i$'s are the residue fields of the semilocal ring $R$. We denote by $a_i$ the image of
$a$ in $S_i=S \otimes_R \kappa_i$ for $i=1,\dots,c$.
For $i=1,\dots, c$, there  exists
a finite \'etale $S_i$-algebra $S'_i$ satisfying the following conditions:

\begin{enumerate}[label={\rm (\roman*)}]
  \item $S'_i$ is a locally free $S_i$--module of rank $n$;

  \item The $\kappa_i$-algebra $S'_i$ admits a unit generator over $\kappa_i$ of the shape $a_i b_i^2$ with
  $u_i \in (S'_i)^\times$.
\end{enumerate}

In particular for $i=1,\dots ,c$ $a_i b_i^2$ is a unit-generator of $S'_i$ over $S_i$,
so there exists  an $S_i$-isomorphism $S_i[x]/P_i(x) \simlgr S'_i$,
$[x] \mapsto a_i b_i^2$ where $P_i$ is a monic $S_i$--polynomial of degree $q^n$
whose discriminant is invertible. We lift the $P_i$'s in a monic
polynomial $P \in S[x]$ of degree $q^n$. By construction we have $\dis(P) \in S^\times$
so that $S'=S[x]/P(x)$ is an \'etale $S$--algebra which is free of rank $n$.
Furthermore we have $S' \otimes_R \kappa_i \cong S' \otimes_S S_i \cong S'_i$ for $i=1,...,n$.

Since $(S')^\times \to (S'_1)^\times  \dots \times (S'_c)^\times$
is onto, we can pick a lift $b \in (S')^\times$ of $(b_1, \dots, b_c)$.
Since $b_i u_i^2$ generates the $\kappa_i$--algebra $ S' \otimes_R \kappa_i$
for $i=1, \dots, c$, the element $a b^2\in S'$ generates the $R$--algebra $S'$ 
by Corollary~\ref{genolem-semi-LG}. 
\end{proof}

\comments{(2026-04-09) PG suggests (my interpretation) trying to generalize \ref{bfp-prop} such that $\rank T$ in \ref{bfp-i} is not divisible by a given finite set of primes. This works except in case $k$ is a infinite field, and $S=F_1 \times \cdots \times F_c$. One needs to take $T=k^\ell$ since $k$ could be algebraically closed, hence $F_i\ot T= F_i^\ell$. But it is not clear that \ref{bfp-iii} holds. But this condition is used in the application of \ref{bfp-prop} in \ref{cor-knex} and \ref{prop_norm_quad}.  

(2026-04-21) I rewrote part of the proof of \ref{bfp-prop}. }

\begin{prop}[{\cite{BFP}}] \label{bfp-prop} Let $R$ be a semilocal ring, let $S$ be a faithful \'etale $R$--algebra and let $N\in \NN_+$. Then there exists an $R$--algebra $T$ such that
\begin{enumerate} [label={\rm (\roman*)}]
  \item \label{bfp-i} $T$ is finite \'etale, unit-generated and has constant odd rank,

  \item \label{bfp-ii} the $R$--algebra $S\ot_R T$ is unit-generated and

  \item  \label{bfp-iii} all residue fields of $S\ot_R T$ have cardinality $\ge N$.
\end{enumerate}\end{prop}
\ms

We recall, \ref{fea}\eqref{fea-ten}, that $S\ot_R T$ is a finite \'etale $R$--algebra. For ``one-generated'' instead of ``unit-generated'' and for $2\in R\ti$ the proposition is \cite[Prop.~7.3]{BFP}. The proof in the unit-generated case is a straightforward modification of the proof in \cite{BFP}, which in fact does not make use of $2\in R\ti$. It is included here for the convenience of the reader.

\begin{proof} After some preliminaries in \eqref{bfp-propI}, we prove the proposition in case $R=k$ is a finite field in \eqref{bfp-propII} and in case $R$ is a product of fields in \eqref{bfp-propIII}. Finally, in \eqref{bfp-propIV} we consider an arbitrary semilocal ring $R$. \sm

\begin{inparaenum}[(I)] \item \label{bfp-propI}
Let $k$ be a field, let $F/k$ and $L/k$ be finite field extensions of degree $n$ and $\ell$ respectively, which are relatively prime. Then the tensor product algebra $F\ot_k L$ is a field 
(\cite[Prop.~2.1]{Cohn}).

Suppose $k$ is a finite field of cardinality $q$, and let $g$ be an irreducible polynomial in $k[X]$ of degree $n$. Hence, $k[X]/(g)$ is a field of cardinality $q^n$, and $k[X]/(g) \cong k[X]/(h)$ for any other irreducible polynomial $h\in k[X]$ of degree $n$. Observe that this does not imply $(g) = (h)$. In fact, the number of monic irreducible polynomials $h\in k[X]$ of degree $n$ is $\frac{1}{n} q^n + {\mathrm O}(q^{n/2}) \approx \frac{1}{n} q^n$, see for example \cite[V, Exc.~22]{Lang}.    

\sm

\item\label{bfp-propII} Let again $R=k$ be a finite field, and let $S$ be a finite \'etale $k$--algebra. Thus $S=F_1 \times \cdots \times F_c$, where $F_1, \ldots, F_c$ are finite extension fields of $k$, say of degree $n_i = \dim_k F_i$, $1\le i \le c$. We choose an odd prime $\ell > \max\{n_1, \ldots, n_c, N\}$ such that there exist  monic irreducible polynomials $g_i\in k[X]$, $1\le i \le c$, of degree $n_i \ell$, which are pairwise distinct, and put $f = \prod_{i=1}^c g_i$. 
    That this is possible, follows from the formula in \eqref{bfp-propI}.
    We then get
    \[ k[X]/(f) \cong \textstyle \prod_{i=1}^c \, k[X]/(g_i) \]  
    by \ref{plr}\eqref{plr-bb}, because the maximal ideals $(g_1), \ldots, (g_c)$ are disitnct, hence  relatively prime. 
    
    Let $L$ be the finite field of cardinality $|k|^\ell$. Since $k[X]/(g_i)$ is a field of cardinality $n_i \ell$ and since $F_i \ot_k L$ is also such a field by \eqref{bfp-propI}, we get $k[X]/(g_i) \cong F_i \ot_k L$ as $k$--algebras and therefore 
    \[ S \ot_k L = \textstyle\prod_{i=1}^c F_i \ot_k L  \cong \prod_{i=1}^c k[X]/(g_i) \cong k[X]/(f).\]
     Since $k[X]/(f)$ is one-generated, so is $S\ot_k L$. Moreover, every primitive element of $\textstyle\prod_{i=1}^c F_i \ot_k L$ is invertible, because $F_i \ot_k L$ is a field. Finally, \ref{bfp-iii} holds with $T=L$ because every epimorphism of $S$ onto a field $K$ factors through (exactly) one of the $F_i \ot_k L$, hence is an extension field of $F_i \ot_k L$ whose cardinality is greater than $N$. This proves the proposition in case $R=k$ is a finite field.
    \sm

\item \label{bfp-propIII} Suppose $R=k_1 \times \cdots \times k_b$ is a finite product of fields $k_i$, $1\le i \le b$. Hence $S=S_1 \times \cdots \times S_b$, where each $S_i$ is a finite \'etale $k_i$--algebra.  Reordering the $k_i$ if necessary, we may  assume that $k_1, \ldots, k_a$ are finite fields and that $k_{a+1}, \ldots, k_b$ are infinite. First suppose $a\ge 1$. We apply \eqref{bfp-propII} to each of the finite fields $k_1, \ldots, k_a$ and thus get fields $L_1, \ldots, L_a$ of odd prime degree $\ell_i$ such that \ref{bfp-ii} and \ref{bfp-iii} hold for the $k_i$--algebra $S_i \ot_{k_i} L_i$. As pointed out in  \eqref{bfp-propII}, we can increase the degrees $\ell_i$ and in this way assume $\ell_1 = \cdots = \ell_a =: \ell$. We are done, if $a=b$.  In case $a< b$, we choose \'etale $k_i$--algebras $T_i$ of degree $\ell$ for $a < i \le b$, e.g., $T_i = k_i^\ell$, and have that $S_i \ot_{k_i} T_i$ is an \'etale $k_i$--algebra, which is unit-generated by \ref{plr}\eqref{plr-f}. By the same argument, also the $k_i$--algebras $T_i$ are unit-generated.
    The $R$--algebra $T= L_1 \times \cdots \times L_a \times T_{a+1} \times \cdots \times  T_b$  is \'etale of degree $\ell$ and unit--generated by \ref{plr}\eqref{plr-a}. Since     \[ S\ot_R T = (S_1 \ot_{k_1} T_1 )\times \cdots \times  (S_b \ot_{k_b} T_b), \]
    it follows again from \ref{plr}\eqref{plr-a} that $S\ot_R T$ is unit-generated. If $a=0$, i.e., every $k_i$ is infinite, then already $S$ is a unit-generated $R$--algebra. It satisfies \ref{bfp-iii} because every residue field of $S$ is an extension field of one of the $k_i$'s. So we are done by taking $T=R$. \sm

\item \label{bfp-propIV} Let now $R$ be an arbitrary semilocal ring. Thus $\ol R = R/\Jac(R)$ is a finite product of fields. By \eqref{bfp-propIII} with $\ol S = S \ot_R \ol R$, there exists a finite \'etale and unit-generated $\ol R$--algebra $\wdh T$ of constant odd degree such that \ref{bfp-ii} and \ref{bfp-iii} hold for the $\ol R$--algebra $\ol S \ot_{\ol R} \wdh T$. By \ref{genolemLG}, $\wdh T = \ol R [X]/(\wdh f)$ for some monic polynomial $\wdh f$ with $\wdh f(0) \in {\ol R}\ti$. Let $f\in R[X]$ be a monic lift of $\wdh f$ and put $T= R[X]/(f)$. Since $T \ot_R \ol R = \ol R[X]/(\wdh f) = \wdh T$, we know from \ref{fea}\eqref{fea-g}
    and \ref{plr}\eqref{plr-c} that $T$ is finite \'etale and unit-generated. It has the same constant degree as $\wdh T$. Since $(S\ot_R T ) \ot_R \ol R \cong \ol S \ot_{\ol R} \wdh T$, a second application of \ref{plr}\eqref{plr-c} shows that $S\ot_R T$ is unit-generated too. Finally, \ref{bfp-iii} holds because every residue field of $S\ot_R T$ factors through $\Jac(R)$ and hence is a residue field of $\ol S \ot_{\ol R} \wdh T$. \end{inparaenum} \end{proof}

 \newpage

\section{Regular and nonsingular quadratic forms}\label{sec:quadratic-forms}

We establish our terminology for quadratic and bilinear forms, review some known results and add several new ones. We will use the notation established in \S\ref{sec:preliminaries-LG}.

\comments{(2025-03-05) Deleted subsection 'detmo' on `Determinant and discriminant modules', used only in the characterization of regular bilinear forms as regular iff the determinant morphism $d_b \co \De(M) \ot_R \De(M) \to R$ is invertible, and then later in the section on Scharlau norms}

\subsection{Symmetric bilinear forms}\label{bfLG} 
A {\em symmetric bilinear $R$--module\/} is a pair $(M,b)$ consisting of a finitely generated projective $R$--module and a symmetric $R$--bilinear form $b\co M \times M \to R$. Since we will only consider bilinear forms that are symmetric, we often simply speak of {\em bilinear $R$--modules} or just {\em bilinear modules} if $R$ is clear from the context. When $M$ is clear from the context or unimportant, we will sometimes write $b$ for $(M,b)$.
    
Given two bilinear modules $(M_1, b_1)$ and $(M_2, b_2)$ we call an $R$--linear map $f \co M_1 \to M_2$ an {\em isometry\/} if $f$ is bijective and $b_2\big(f(m_1), \, f(m_1') \big) = b_1(m_1, m_1')$ holds for all $m_1, m_1' \in M_1$. We list some properties of quadratic modules. \sm   

\begin{inparaenum}[(a)] \item\label{bfLG-ad} ({\em Adjoints, radical, regularity}) 
Given a bilinear module $(M,b)$, its {\em adjoint\/} is the $R$--linear map
$\wdh b \co M \to M^*= \Hom_R(M,R)$, $m \mapsto b(m, \cdot)$. The {\em radical of $(M,b)$\/} is the kernel of $\wdh b$, i.e, $\rad(b) = \{ m \in M : b(m, M) = 0 \}$. 
We call $(M,b)$ {\em regular\/} if $\wdh b$ is an isomorphism.\footnote{ We warn the reader that the terminology of a regular bilinear form is not universally accepted, but it follows \cite{CF}, \cite{PRbook}, \cite{K} and \cite{Sc}, except that in the last two references ``regular'' and ``nonsingular'' are used interchangeably. A  regular bilinear form as defined here is called ``non singular'' in \cite{Ba} and ``nondegenerate'' in the book \cite{EKM}.}
\sm 

\item \label{bfLG-bc} ({\em Base change}) We associate with a bilinear $R$--module $(M,b)$ and $S\in \Ralg$ the bilinear $S$--module $(M,b)_S = (M_S, b_S)$, given by $b_S \co M_S \times M_S \to S$, $b_S(m_1 \ot s_1, \, m_2\ot s_2) = b(m_1, m_2) s_1 s_2$,  viewing $b(m_1, m_2)$ as element of $S$ canonically. Base change respects adjoints in the sense that
    \begin{equation}\label{bfLG-bc1} \vcenter{
      \xymatrix{M_S \ar[rr]^{(\wdh{b})_S} \ar[dr]_{\wdh{(b_S)}}
         && (M^*)_S \ar[dl]^{\nu_S}_\cong
        \\ & (M_S)^*
    }}\end{equation}
commutes, where $\nu$ is the canonical isomorphism. 
We will often identify $(M_S)^* \equiv (M^*)_S =: M_S^*$ and hence $\wdh{(b_S)}\equiv (\,\wdh{b}\, )_S = \wdh{b}_S$.  

\sm

\item \label{bfLG-rech} ({\em Characterizations of regularity}) The following are equivalent for a bilinear module $(M,b)$: \end{inparaenum}

\begin{enumerate}[label={\rm (\roman*)}]

\item\label{bfLG-rechi} $b$ is regular;

\item\label{bfLG-rechii} $b_S$ is regular for every $S\in \Ralg$;

\item\label{bfLG-rechiii} $b_K$ is regular for all algebraically closed fields
    $K\in \Ralg$;

\item\label{qfLG-rechiv} $b_{R_\m}$ is regular for all maximal\/ $\m \in
    \Spec(R)$;

\item\label{qfLG-rechv} $b_{R/\m}$ is regular for all maximal\/ $\m \in
    \Spec(R)$;

\item\label{bfLG-rechvi} there exists a Zariski cover $(f_1, \ldots, f_n)$ of $R$ such that every $b_{R_{f_i}}$, $i=1, \ldots, n$, is regular;

\item\label{bfLG-rechvii} $b_T$ is regular for some faithfully flat $T\in \Ralg$;

\item  \label{bfLG-rechCF} $\rad(b_S) = 0$ for all $S \in \Ralg$;

 \item \label{bfLG-rechix} $b_{R/I}$ is regular for some ideal $I \ideal \Jac(R)$.
\end{enumerate}
If $M$ is free, say with basis $(e_1, \ldots, e_n)$, then $b$ is regular if and only if 
\begin{enumerate}[label={\rm (\roman*)}]\setcounter{enumi}{9}

\item \label{bfLG-rechxi} $\det\big( b(e_i, e_j)\big) \in R^\times$. 
\end{enumerate}
\sm 

The proof that \ref{bfLG-rechii}--\ref{bfLG-rechvii} characterize regularity is standard, see for example \cite[C.2]{GN-LG}. The equivalence \ref{bfLG-rechi} $\iff$ \ref{bfLG-rechii} follows of course from the commutative diagram \eqref{bfLG-bc1}. The condition \ref{bfLG-rechCF} is the definition of regularity in \cite[2.6.0.26]{CF}; it is equivalent to regularity in our sense by \cite[2.6.0.28]{CF}.

\inparcom{(2025-03-05) Further equivalent conditions, likely not used here: 

{\tt ('qfba-rechdet')  the determinant morphism $d_b \co \De(M) \ot_R \De(M) \to R$ is invertible, i.e., $(\De(M), d_b)$ is a discriminant module in the sense of \cite[III, \S3]{K} } }

\begin{inparaenum}[(a)]\setcounter{enumi}{3}
\item \label{qfba-sum} ({\em Orthogonal sum}) Let $(M_i, b_i)$, $i=1,2$, be two bilinear $R$--modules. Their {\em orthogonal sum\/} is the bilinear $R$--module
          $(M_1, q_1)\perp (M_2, q_2) = (M_1 \oplus M_2, b_1 \perp b_2)$ with 
\[ (b_1 \perp b_2)(m_1 + m_2, m_1' + m_2') = b_1(m_1,m'_1) + b_2(m_2, m_2') \]
 for $m_i, m'_i \in M_i$.  We sometimes write $(M_1, b_1) \oplus (M_2, b_2)$ instead of $(M_1, b_1) \perp (M_2, b_2)$. The orthogonal sum $(M_1, b_1) \perp (M_2, b_2)$ is regular if and only if $(M_1, b_1)$ and $(M_2, b_2)$ are regular, see for example \cite[I, (3.6.2.2)]{K}. \sm 
 
\item\label{bfLG-odd} Let $(M,b)$ be a regular bilinear module. If $M$ has constant odd rank, then $2\in R\ti$ (\cite[I, (6.2.3)]{K}. \sm 

\item \label{meta} ({\em Metabolic spaces}) Let $(U,b)$ be a bilinear module.  The {\em metabolic space\/} associated with $(U,b)$ is the bilinear module $\MM(U,b) = (U \oplus U^*, b_{\MM(U,b)})$ whose bilinear form $b_{\MM(U,b)}$  is defined by
\begin{equation*} \label{meta-1}
 b_{\MM(U,b)} (u + \vphi, v + \psi) = b(u,v) + \vphi(v) + \psi(u).
 \end{equation*}
It is a regular bilinear form. We say a bilinear module $(M, b)$ is {\em metabolic\/} if there exists a bilinear module $(U, b_U)$ such that $(M, b)\cong \MM(U,b_U)$. A regular bilinear module $(M, b)$ is metabolic if and only if $M$ contains a totally isotropic complemented submodule $V$ with $V = V^\perp$, a so-called {\em Lagrangian}, \cite[I, Thm.~(4.6)]{Ba}, see \ref{Lag} for Lagrangians of quadratic modules.
\sm 

\inparcom{(2022-04-07) In Knebusch's 1976 Notes, these metabolic spaces are called ``split metaboliic". He defines ``metabolic'' as a regular quadratic $\calO_X$--module containing a Lagrangian and then shows that over an affine scheme the two notions coincide (Cor.~1 in I, \S3)}

\item \label{mx_lemb} ({\em Unimodularity}) Let $(M,b)$ be a bilinear module and assume $x,y\in M$ satisfy $b(x,y) = 1$. Then $x$ and $y$ are unimodular by \ref{unimodii} of \ref{unimod}, 
    \[ M = Rx \oplus (Ry)^\perp = Ry \oplus (Rx)^\perp,
     \] 
     and $(Rx)^\perp$ is finite projective with $\rank_\p (Rx)^\perp = \rank_\p (M) - 1$ for\/ $\p \in \Spec(R)$. 
    
\end{inparaenum} 

\begin{lem}\label{orthLG}
Let $(M,b)$ be a bilinear module and let $U \subset M$ be a submodule.
We put 
\begin{equation}\label{orthLGperp}  U^\perp = \{ m \in M : b(m, u) = 0 \text{ for all $u\in U$} \},
\end{equation} 
and abbreviate $b_U =b|_{U\times U}$ and $b_{U^\perp} = b|_{U^\perp \times U^\perp}$. \sm

\begin{inparaenum}[\rm (a)]\item\label{orthLG-a}
Suppose that $b_U$ is regular. Then
\begin{equation}\label{orthLG-0}
  (M,b) = (U, b_U) \; \perp \; (U^\perp, b_{U^\perp}).
\end{equation}
In particular, $U$ and  $U^\perp$ are finite projective, complemented submodules,  and
\begin{equation} \label{quadrepII-b1n}
 \rank_R M = \rank_R U + \rank_R U^\perp.
\end{equation}
Moreover, $(M,b)$ is regular if and only if $(U^\perp, b_{U^\perp})$ is regular, and in this case $U^{\perp\perp} = U$ holds. \sm

\item \label{orthLG-b} Suppose that $(M,b)$ is regular and that $U$ is complemented in $M$, hence finite projective.  Then $U^\perp$ is complemented, thus finite projective, the rank formula \eqref{quadrepII-b1n} holds and $U^{\perp\perp} = U$. 
\end{inparaenum} \end{lem}
\sm

Part~\eqref{orthLG-a} of \ref{orthLG} is proven in \cite[Lem.~(2.2)]{Bass-69}, and again in \cite[I, (3.2)]{Ba} and \cite[I, (3.6.2)]{K}. Part~\eqref{orthLG-b} is shown in \cite[Lem.~(2.1)]{Bass-69}; it is a special case of \cite[I, \S2, Prop.~1]{Knebusch-Queens}. We point out that \eqref{orthLG-0} need not be true in the setting of \eqref{orthLG-b}. 

\comments{More details for \ref{orthLG}\eqref{orthLG-b}: 

{\tt Moreover, identifying $(M/U)^* \equiv \{ \la \in M^* : \la(U) = 0 \}$, there exist unique isomorphisms
$ \al\co M/U^\perp \simlgr U^*$ and $\be \co U^\perp \simlgr (M/U)^*$
such that the diagrams
\[  \vcenter{\xymatrix{ M \ar[r]^{\hat b}_\cong  \ar[d]   & M^* \ar[d] \\
              M/{U^\perp} \ar[r]_\cong^\al & U^*}}
    \qquad\text{and} \qquad
\vcenter{ \xymatrix{ U^\perp \ar[r]^\be_\cong  \ar[d]  & (M/U)^* \ar[d] \\
              M  \ar[r]_\cong^{\hat b} & M^*}}
\]
commute. In these diagrams the vertical maps are canonical, surjective on the left and injective on the right.} }

\new
\comments{(2026-04-16) Added Lemma~\ref{sursch}, needed to prove that the smooth part of a sphere respects base change, which in turn is needed in the proof of representability of the Knebusch functor. }

\begin{lem}\label{sursch} Let $(M,b)$ be a bilinear module. We use the isomorphism $\nu_S$, $S\in \Ralg$,  of \eqref{bfLG-bc1} as identification, i.e., $M_S^* = (M^*)_S = (M_S)^*$, and thus get adjoint maps $\wdh b_S \co M_S \to M_S^*$ . \sm 

\begin{inparaenum}[\rm (a)] \item \label{sursch-a}
The family $(\wdh b_S)_{S\in \Ralg}$ defines a morphism $\ulW(M) \to \ulW(M^*)$ of $R$--functors, and hence a morphism $\wdh \bfb \co \uW(M) \to \uW(M^*)$ of the associated $R$--schemes. \sm 

\item \label{sursch-b} The inverse image of the subscheme $\uW(M^*)_u$ of $\uW(M)$, defined in {\rm \ref{spreq-a}}, under $\wdh \bfb$ is a quasi-compact open subscheme 
    \[ \Sur_b = \wdh \bfb \me \big(\uW(M^*)_u\big) 
    \] 
    of $\uW(M)$, representing the $R$--functor 
    \[  S \mapsto \{ m\in M_S : \wdh b_S(m) \co M_S \to M_S^* \text{ is surjective}\}.
    \] 

\item\label{sursch-c} The scheme $\Sur_b$ is stable under base change: for $R'\in \Ralg$ there exists a canonical isomorphism of $R'$--schemes
\begin{equation}\label{sursch-c1}
  \Sur_b \times_R R' \cong \Sur_{b\ot_R R'}\; .
\end{equation}     
\end{inparaenum}\end{lem}

\begin{proof} \eqref{sursch-a} is obvious. \eqref{sursch-b} The inverse image of a  quasi-compact open subscheme under a morphism between affine schemes is again quasi-compact open. By \ref{spreq-a}\ref{spreg-v}, $\Sur_b$ represents the functor described in \eqref{sursch-b}. For any $A'\in \Rpalg$ we have the canonical isomorphism $M_{R'}\ot_{R'} A' \cong M\ot_R A'$, hence \eqref{sursch-c} is clear at the level of $R'$--functors. \end{proof}
\enew

\comments{(2025-03-05) The subsection \ref{qf} is based on \cite[C.8]{GN-LG}.}

\subsection{Quadratic forms}\label{qf} 
\begin{inparaenum}[(a)]
 \item \label{qfLG-a} A {\em quadratic form (over $R$)\/} is a pair $(M,q)$ consisting of a finite projective $R$--module $M$ and a map $q \co M \to R$ satisfying $q(rm) = r^2 q(m)$ for all $r\in R$ and $m\in M$ and for which the {\em polar form $b_q$\/}, defined by $b_q(m, m') = q(m+m') - q(m) - q(m')$ for $m,m'\in M$, is a (symmetric) bilinear form on $M$. 
     
     We often abbreviate $(M,q)= q$ and refer to $(M,q)$ as a {\em quadratic module}.
     We call $(M,q)$ a {\em faithful\/} quadratic $R$--module if $M$ is a faithfully  projective $R$--module, \ref{fapmod}. We say that $(M,q)$ is a quadratic module of constant rank or of rank $n\in \NN$, if $M$ has that property. 

If $(M,q)$ is a quadratic module and $U\subset M$ is a submodule, then 
\[ U^\perp = \{m\in M : b_q(m,u) = 0 \text{ for all $u\in U$}\} \]
as in Lemma~\ref{orthLG}.

An {\em isometry\/} $f\co (M_1, q_1) \to (M_2, q_2)$ is an $R$--linear isomorphism $f\co M_1 \to M_2$ satisfying $q_2 \circ f = q_1$. If such an $f$ exists, we use the symbol $(M_1, q_1) \cong (M_2, q_2)$ to indicate this. 
If $(M_1, q_1) = (M_2, q_2) = (M,q)$, the isometries of $(M,q)$ form a group $\orth(M,q) = \orth(q)$, the {\em orthogonal group of $q$}.      \sm

\item ({\em Regularity}) \label{qf-regLG} By definition, a quadratic form $q$ is {\em regular\/}, if its polar form $b_q$ is regular in the sense of \ref{bfLG}\eqref{bfLG-ad}. Lemma~\ref{orthLG} holds {\em mutatis mutandis} for quadratic modules using \eqref{qf-perp} for the concept of orthogonal sums of quadratic forms.

  If $(M,q)$ is a regular quadratic module of constant odd rank, then $2\in R\ti$ by \cite[IV, (3.1)]{K}. 
\sm

\item ({\em Base change}) \label{qf-bc} Let $T\in \Ralg$ and let $(M,q)$ be a quadratic form. Analogous to \ref{bfLG}\eqref{bfLG-bc} we consider the $T$--module $M_T = M \ot_R T$. There exists a quadratic form $(M_T, q_T)$ over $T$  uniquely determined by the condition $q_T(m\ot t) = q(m) t^2$ for all $m\in M$ and $t\in T$ (\cite[Prop.~II.1]{Roby}, \cite[Thm.~1]{Sah}, or see \cite[11.5]{PRbook} for a recent reference). We often abbreviate $(M,q)_T = (M_T), q_T)$. The polar of $q_T$ is the base change of the polar $b_q$ of $q$, i.e., $b_{q_T} = (b_q)_T$. In particular, if $q$ is regular, then so is $q_T$. 
    
    If $f \co (M_1, q_1) \to (M_1, q_2)$ is an isometry of quadratic modules over $R$, the $T$--linear map $f_T = f \ot 1_T  \co (M_1, q_1)_T \to (M_2, q_2)_T$ is an isometry of the quadratic modules over $T$.   
   
    Let $\qm_R$ be the category whose objects are quadratic modules over $R$ and whose morphisms are isometries, and define $\qm_T$ accordingly. The assignments $(M,q) \mapsto (M,q)_T$ and $f \mapsto f_T$ define a {\em base change functor} $\qm_R \to \qm_T$.  
    \sm

\item\label{qf-rad} ({\em Radical}) The {\em radical\/} of a quadratic module $(M,q)$ is the submodule 
    \[  \rad(q)= \{m\in M : q(m) = 0 = b_q(m,M)\} \subset \rad(b_q) \]
     of $M$. Given $T\in \Ralg$, clearly $\rad(q) \ot_R T \subset \rad(q_T)$. This inclusion is in general not an equality, cf.~\eqref{qf-ns}. \sm

\item \label{qf-ns} ({\em Nonsingularity}) A quadratic form $(M,q)$ is called {\em nonsingular\/} if $\rad(q_F) = 0$ for all fields $F\in \Ralg$. In this case, we call $(M,q)$ a {\em quadratic space}.\footnote{As for regularity, \ref{bfLG}\eqref{bfLG-ad}, it is appropriate to point out that there is no universally accepted terminology for what we call a nonsingular quadratic form. Our terminology follows \cite{Sw}. A nonsingular quadratic form is called ``non-degenerate'' in \cite{Co1}, ``nondegenerate'' in \cite{EKM}, ``semiregular'' in \cite{CF,K} in case of odd rank, cf. \ref{qfnsp-vi} of \eqref{qfnsp}, ``separable'' in \cite{Lo0,Lo2,P-Fields}, 
     and ``ordinary'' in \cite[XII]{SGA7}.
     If $2\in R\ti$, then nonsingular = regular by \eqref{qf-ns2}. Hence our terminology coincides with that of Lam's book \cite{Lam-qf}.
    We point out that a ``quadratic space'' in the sense of \cite{Ba} or \cite{PRbook} is a quadratic module $(M,q)$ with a regular $q$.}

    If $q$ is regular, then $\rad(q) = 0$ and hence \ref{bfLG-rechii} of \ref{bfLG}\eqref{bfLG-rech} implies
\begin{equation}\label{qf-ns1}
   \text{\em $q$ regular} \implies \text{\em $q$ nonsingular.}
\end{equation}
The converse of \eqref{qf-ns1} is not true, see \eqref{qfba-one}. 
However, if $2\in R\ti$, then $\rad(q) = \{ m\in M : b_q(m,M) = 0 \}=\rad(b_q)$ for any quadratic form $q$ and so $\rad(q)$ is stable under base change. Hence,%
\begin{equation}\label{qf-ns2}
\text{\em  if  $2 \in R\ti$, then $q$ is nonsingular if and only if $q$ is regular.}
\end{equation}%

\item \label{qfba-one} ({\em $1$--dimensional forms}) Let $u\in R$. We define a bilinear  form $\lan u \ran_b$ and a quadratic form $\lan u \ran_q$ on $R$,
    \begin{equation}\label{qfba-one1} \begin{split}
    \lan u \ran_b & \co R \times R \to R, \quad (r_1, r_2) \mapsto u r_1 r_2 \\
    \lan u \ran_q & \co R \to R, \quad \qquad r \mapsto ur^2.
    \end{split}\end{equation}
Then 
\begin{equation}\label{qfba-one1}
  \text{$\lan u \ran_b$ is regular}\iff u\in R\ti 
     \iff \text{$\lan u \ran_q $ is nonsingular.} 
\end{equation}
Indeed, the first equivalence follows from \ref{bfLG-rechxi} of \ref{bfLG}\eqref{bfLG-rech}. Regarding the second, observe that $\rad\big(\lan u \ran_q)_S\big) = \{ s\in S : s^2 (u\ot 1_S) = 0 = 2 s (u\ot 1_S) \}$ for any $S\in \Ralg$ and this vanishes for every field $F$ if and only if $u \ot 1_F = 0$, if and only if $u\in R\ti$. 

The polar of the nonsingular quadratic form $\lan u \ran_q$ is $2 \lan u \ran_b$, whence {\em $\lan u \ran_q$ is regular if and only if $2\in R\ti$ and $u\in R\ti$.}  That $2\in R\ti$ if $\lan u \ran_q$ is regular, is a general fact for regular forms on spaces of odd rank, \ref{bfLG}\eqref{bfLG-odd}.

Any quadratic $R$--space $(M,q)$ with $M$ of constant rank $1$ has the form $(M,q) = \lan u \ran_q $ for some $u\in R\ti$ if $\Pic(R) = 0$, e.g., $R$ is unimodular or even LG, \ref{unifap} and \ref{revLG}. \sm 

\item ({\em $2$-dimensional forms}) 
Let $M$ be a free $R$--module of rank $2$ with basis $e_1, e_2$. For $a, b\in R$ we define the quadratic form
    \begin{equation}\label{qfba-two0}
         [a,b] \co M \to R, \quad r_1 e_1 + r_2 e_2 \mapsto a r_1^2 + r_1 r_2 + b r_2^2.
\end{equation}
If $a\in R\ti$, then $[a,b] = a n_E$ where $E$ is the free quadratic $R$--algebra $E=R[X]/(X^2 - a \me X + a\me b)$.  

By \ref{bfLG-rechxi} of \ref{bfLG}\eqref{bfLG-rech}, the quadratic form $[a,b]$ is regular (= nonsingular because of \eqref{quadfoe}) if and only
if $1-4ab \in R\ti$. 
Hence
\begin{equation}  \label{qfba-two1} \begin{split}
& \text{$[a,b]$ regular with $a\in R\ti$} \quad \implies \\
&\quad \text{$[a,b] = a n_E$ with $E$ free quadratic \'etale.}
\end{split} \end{equation}

\item ({\em Orthogonality}) \label{qf-perp} Given two quadratic forms $(M_1, q_1)$ and $(M_2, q_2)$, their {\em orthogonal sum\/} is the quadratic form $q_1 \perp q_2$, defined on $M= M_1 \oplus M_2$ by $(q_1 \perp q_2)(m_1, m_2) = q_1(m_1) + q_2(m_2)$ for $m_1 \in M_1$ and $m_2 \in M_2$. The polar form of $q_1 \perp q_2$ is the orthogonal sum $b_{q_1 \perp q_2} = b_{q_1} \perp b_{q_2}$ of \ref{bfLG}\eqref{qfba-sum}. As noted there, 
the quadratic form $q= q_1 \perp q_2$ is regular if and only if $q_1$ and $q_2$ are regular. Regarding nonsingularity,    one easily sees:
\begin{equation}\label{qf-perp0}
q_1 \perp q_2 \text{ nonsingular } \implies \text{$q_1$ and $q_2$ nonsingular. } 
\end{equation}
The converse of \eqref{qf-perp0} is not true. Indeed, assume $2R=0$. Then  the  quadratic form $\lan 1 \ran_q$ on $R$ is nonsingular by \eqref{qfba-one}, but their orthogonal sum $q= \lan 1 \ran_q \perp\lan 1 \ran_q$ is singular since $b_q = 0$ and $q(1,1) = 0$. However, 
\begin{equation} \label{qf-perp1} \begin{split}
 &\text{\em if $q_1$ is regular, then $q_1\perp q_2$ is nonsingular} \\
 & \iff \text{\em $q_2$ is nonsingular.}%
\end{split} \end{equation}
\lv{
Indeed, if $q$ is nonsingular, then $\rad(b_q) = \rad(b_{q_1}) \oplus \rad (b_{q_2})$ and $\rad(q_1) \oplus \rad(q_2) \subset \rad(q)$. Since $(M_S, q_S)$ is an orthogonal sum for any extension $S$ of $R$, nonsingularity of $q_1$ and $q_2$ follows. Conversely, if $(M_1, q_1)$ is regular,
then $\rad (b_q) = \rad (b_{q_2})$ and $\rad(q) = \rad(q_2)$
      
 Indeed, let $m \in \rad(q)$ and write $m=m_1 + m_2$ with $m_i in M_i$. Then
$0 = q(m, M_1) = q(m_1, M_1)$ implies $m_1 \in \rad(b_{q_1}) = 0$, whence
$m=m_2$, and $0 = q(m) = q(m_2)$ shows $m=m_2 \in \rad(q_2)$. The other
inclusion is equally obvious.
Since $q_1{}_S$ is again regular, we get $\rad(q_F) = \rad((q_2)_F)$ for any
field $F \in \Ralg$, whence the claim.}

\item \label{qf-redc} ({\em Direct products of base rings})
Let $R= R_0 \times \cdots \times R_n$ be a direct product of rings. 

Recall \ref{dpb}: Every $R$--module $M$ uniquely decomposes $M= M_0 \times \cdots \times M_n$ as a direct product of $R_i$-modules $M_i = R_i M$. Conversely, every family $(M_i)_{0\le i \le n}$ gives rise to an $R$--module $M=M_0 \times \cdots \times M_n$. The $R$--module $M$ is finite projective if and only if every $R_i$--module $M_i$ is finite projective.  

Let $(M_i, q_i)$, $i=0,\ldots, n$, be quadratic $R_i$--modules. Then
 \begin{equation}\label{qfba-red1} \begin{split}
     (M_0,q_0) \times \cdots \times (M_n, q_n) &= (M_0 \times \cdots \times M_n,  q_1 \times \cdots \times q_n) \\
    (q_0 \times \cdots \times q_n)(m_0, \ldots, m_n) &= (q_0(m_0), \ldots, q_n(m_n))
 \end{split} \end{equation}
is a quadratic $R$--module, referred to as the {\em direct product of the $(M_i, q_i)$}. The decomposition \eqref{qfba-red1} respects orthogonal groups:
  \begin{equation} \label{ortgr-bas-b1}
      \orth(q) = \orth(q_0) \times \cdots \times \orth(q_n).
   \end{equation}
Conversely, a quadratic $R$--module $(M,q)$ uniquely decomposes as a direct product
$(M,q) = (M_0, q_0) \times \cdots \times (M_n,q_n)$ of quadratic $R_i$--modules $(M_i, q_i)$. The quadratic $R$--module $(M,q)$ is regular (nonsingular respectively) if and only if every $(M_i, q_i)$ is a regular (nonsingular respectively) quadratic $R_i$--module.

A standard way to obtain the situation considered here occurs by letting $M=M_0\times \cdots \times M_n$ be the {\em rank decomposition\/} of a finite projective $R$--module $M$ for which $M_i$, $0\le i \le n$, is a finite projective $R_i$--module of constant rank $i$.
The discussion above then describes the reduction of quadratic modules to quadratic modules of constant rank. We will refer to this process as {\em reduction to constant rank}.
\sm

\item\label{qf-hyp} ({\em Hyperbolic spaces, hyperbolic planes, and hyperbolic pairs}) Let $U$ be a finite  projective $R$--module. The associated {\em hyperbolic space $\HH(U)$} is the quadratic module $(U^* \oplus U, \hyp)$ with quadratic form $\hyp_U(\vphi \oplus u) = \vphi(u)$, where $\vphi \in U^*$ and $u\in U$. The quadratic form $\hyp_U$ is regular, hence nonsingular by \eqref{qf-ns1}. In general, a {\em hyperbolic space\/}  is a quadratic module $(M,q)$ isometric to some $\HH(U)$. Any quadratic module isometric to $\HH(R)$ is called a {\em hyperbolic plane}. 
    
    Any isomorphism $g \co U \simlgr U'$ of $R$--modules extends to an isometry of the corresponding hyperbolic spaces,
\begin{equation}\label{qfba-hyp1}
    \HH(g) \co \HH(U) \simlgr \HH(U'), \quad \vphi \oplus u \mapsto ({^tg}{}\me)(\vphi) \oplus g(u).
 \end{equation}
Given two hyperbolic spaces $\HH(U)$ and $\HH(V)$, the canonical extension of a linear form on $U$ to a linear form on $U \oplus V$ shows that
\begin{equation}\label{qfba-hyp3}
  \HH(U) \perp \HH(V) \cong \HH(U \oplus V).
\end{equation}
    
Let $(M,q)$ be a quadratic module. A pair $(e,f)$ of elements of $M$ is called a {\em hyperbolic pair\/} if $q(e) = 0 = q(f)$ and $b_q(e,f) = 1$ holds. If $(e,f)$ is such a pair,  $\HH = Re + Rf$ is a free submodule of rank $2$ with basis $\{e,f\}$, the quadratic module $(\HH, q|_\HH)$ is isometric to $\HH(R)$, a hyperbolic plane, and 
    $(M,q) = (\HH, q|_\HH) \perp (\HH^\perp, q|_{\HH^\perp})$ by \ref{orthLG}\eqref{orthLG-a}, in particular $\rank M \ge 2$. See \ref{hps} for a generalization of hyperbolic pairs.
    
Let $U$ be a free $R$--module, say with an $R$--basis $e_1, \ldots, e_n$ and let $f_1, \ldots, f_n$ be the associated dual basis of $U^*$. Then $(i_1, f_1), \ldots, (e_n, f_n)$ are hyperbolic pairs such that $(e_i, f_i) \perp (e_j, f_j)$ for $i\ne j$.      
\sm

\item\label{quadfoc} ({\em Split quadratic forms}) Let $m\in \NN$. The quadratic  form $q_{0, 2m}$ is the hyperbolic form associated with the free $R$--module $R^m$. After identifying $R^m{}^* = R^m$, it is given on $R^{2m}$ by
\begin{equation}  q_{0, 2m} (r_{-m}, \ldots, r_{-1}, r_1, \ldots, r_m) = \textstyle \sum_{i=1}^m r_i r_{-i}, \label{quadfoc1}
\end{equation}
It is regular, hence also nonsingular by \eqref{qf-ns1}. The quadratic form
$q_{0, 2m+1} = \lan 1 \ran \perp q_{0, 2m}$ on $R^{2m+1}$, defined by
\begin{equation}
q_{0, 2m+1} (r_{-m}, \ldots, r_{-1}, r_0 ,r_1 \ldots, r_{m}) = r_0^2
+ \textstyle \sum_{i=1}^m r_i r_{-i},
\label{quadfoc2}
\end{equation}
is nonsingular, e.g.,  by the even rank case, by \eqref{qf-perp1} and by nonsingularity of $\lan 1 \ran$, see \eqref{qf-ns}. We will refer to $q_{0, n}$ for
$n$ even or odd as the {\em split quadratic forms\/}, see \eqref{qfnsp} for a justification for this terminology. \sm 

\item\label{qfnsp} ({\em Characterizations of nonsingularity} \cite[C.9]{GN-LG}) 
For a faithful quadratic module $(M,q)$ the following are equivalent:
\end{inparaenum} 
\begin{enumerate}[label={\rm (\roman*)}]
  \item \label{qfnsp-i} $q$ is nonsingular;

  \item \label{qfnsp-ii} $q_T$ is nonsingular for all $T\in \Ralg$;

  \item \label{qfnsp-iii} $q_{R/\gm}$ is nonsingular for all maximal ideals $\gm \ideal R$;

  \item\label{qfnsp-iv}  there exists a flat cover $(R_1, \ldots, R_n)$ such
    that each $M \ot_R R_i$ is a free $R_i$-module of finite rank $r_i$ and $q_{R_i}$  is the split quadratic form $q_{0,r_i}$ over $R_i$ defined in {\rm
    \ref{qf}\eqref{quadfoc}};

  \item \label{qfnsp-v} $q_S$ is nonsingular for some faithfully flat $S\in \Ralg$;
\end{enumerate}
If $M$ has constant rank $n\in \NN_+$, then \ref{qfnsp-i}--\ref{qfnsp-v} are equivalent to
\begin{enumerate}[label={\rm (\roman*)}]\setcounter{enumi}{5}
 \item\label{qfnsp-vi} $q$ is regular if $n$ is even, and $q$ is semiregular in the sense of {\rm \cite[IV, (3.1)]{K}} if $n$ is odd;
\end{enumerate}
If $R$ is a field, then $q$ is nonsingular if and only if
\begin{enumerate}[label={\rm (\roman*)}]\setcounter{enumi}{6}
  \item \label{qfnsp-vii} $q$ is {\em nondegenerate\/} in the sense of {\rm \cite[(7.17)]{EKM}}, i.e., one of the following two conditions hold:

   \begin{enumerate}[label={\rm (\alph*)}]
     \item $q$ is regular, or

     \item $\Char(R) = 2$, $\rad(q) = 0$, $\dim_R\{ m\in M : b_q(m,M) = 0 \} = 1$, and $\dim_R M$ is odd.
\end{enumerate}
\end{enumerate}

\begin{inparaenum}[(a)] \setcounter{enumi}{12} 
\item\label{quadfoe} {\em If $(M,q)$ is a quadratic module with $M$ of  constant even rank, then $q$ is regular if and only if $q$ is nonsingular}. 
    
Indeed, by \eqref{qf-ns1} we only need to show that under our assumption on $M$ a nonsingular form is regular. But this is \ref{qfnsp-vi} of \eqref{qfnsp}.  \sm

  \item\label{qf-detgr}  ({\em Orthogonal groups and determinants}) If $(M,q)$ is nonsingular and $g\in \orth(q)$, then $\det(g) \in \bmu_2(R) = \{r \in R : r^2 = 1_R \}$.   

Indeed, by \eqref{qf-redc} we can assume that $M$ has constant rank. In this case, the claim follows from \cite[IV, (5.1.1)]{K}, keeping in mind that $q$ is regular of even rank or semiregular of odd rank by \ref{qfnsp-vi} of \eqref{qfnsp}. 
\sm 

We always have a group homomorphism
\[ z_M \co \mu_2(R) \to \orth(q), \quad r \mapsto r\Id_M. \]
If $(M,q)$ has odd rank, $z_M$ is a section of the determinant homomorphism, 
\begin{equation}
  \label{qf-detgr1} \xymatrix@C=50pt{\orth(q) \ar@<0.5ex>[r]^{\; \det \; } &\mu_2(R) 
      \ar@<0.5ex>[l]^{\quad z_M \quad}. }
\end{equation}
\end{inparaenum}

\new
\comments{(2026-05-24) The following \ref{qfba-tens} is only used in the section on Scharlau's norm principle }

\subsection{The tensor product between a bilinear and a quadratic module}   \label{qfba-tens}  
Given an $R$--bilinear module $(M,b)$ and an $R$--quadratic form $(N,q)$, there exists a unique $R$--quadratic form $b\ot q \co M \ot_R N\to R$ satisfying
\begin{equation}\label{tenssq0}
\begin{split} (b \ot q)\, (m \ot n) &= b(m,m) \, q(n), \; \text{and}\\
  b_{b\ot q}(m \ot n, m' \ot n') &= b(m, m')\, b_q(n, n')
\end{split}
\end{equation}
for all $m, m'\in M$ and $n, n'\in N$, where $b_q$ is the polar form of $q$. It is called the {\em tensor product of $(M,b)$ and $(N,q)$}, \cite[Thm.~1]{Sah}.

{\em Sketch of proof.} Uniqueness is clear from \eqref{tenssq0}.
To show the existence of $b\ot q$, one first defines a quadratic form $q_F$ on the free $R$--module $F$ with basis $M \times N$ satisfying \eqref{tenssq0} with $m \ot n$ replaced by the basis element $(m,n)\in F$. One then verifies that the kernel of the canonical map $F \to M \ot_R N$ lies in the radical of $q_F$, defined in \ref{qf}\eqref{qf-rad}, which therefore descends to a quadratic form $b \ot q$ with the required properties.
\sm

A moment's thought will convince the reader that {\em the tensor product respects orthogonal sums}, i.e.,
\begin{equation}\label{qfba-tens-p} \begin{split}
  &(b_1 \perp b_2) \ot (q_1 \perp q_2) \cong
  \\& \quad (b_1\ot q_1) \perp (b_1 \ot q_2) \perp (b_2 \ot q_1) \perp (b_2 \ot q_2)
\end{split} \end{equation}
holds for bilinear modules $(M_i, b_i)$ and quadratic modules $(N_i, q_i)$, $i=1,2$. The analogous formula for the tensor product between orthogonal bilinear modules is of course also true. Moreover, {\em the tensor product is associative and commutative}: with obvious notation we have
\begin{equation} \label{qfba-tenso-2} \begin{split}
  (b_1 \ot b_2) \ot b_3 &\cong b_1\ot (b_2 \ot b_3), \quad
     (b_1 \ot b_2) \ot q \cong (b_1 \ot( b_2 \ot q),
 \\  b\ot q &\cong q \ot b
\end{split} \end{equation}
{\em The tensor product also respects isometries:} if $\vphi \co b \simlgr b'$ and $\psi \co q \to q'$ are isometries, then $\vphi \ot \psi \co b \ot q \simlgr b' \ot q'$ is an isometry.

If $(M,b)$ and $(N,q)$ are bilinear and quadratic modules respectively
   and $ub$ and $uq$ are the obvious scalar multiplications by $u$, then, using the definition \eqref{qfba-one},
\begin{equation} \label{qfba-tens-2}
 \lan u\ran_b \ot b \cong  ub, \quad \lan u \ran_b \ot q \cong uq
 \end{equation}
under the standard isomorphism $R \ot_R  M \simlgr M$. However,
$ \lan 1\ran_q \ot b = q_b$, where $q_b \co M \to R$, $m \mapsto b(m,m)$, is the quadratic form associated with $b$. Its polar is $b_{q_b} = 2b$.

The tensor product is compatible with base change: for $S\in \Ralg$ we have
\begin{equation}\label{qfba-tens-2}
(M,b)_S \ot_S (N,q)_S \simlgr \big((M,b) \ot_R (N,q)\big)_S
\end{equation}
with respect to $m\ot s_1 \ot n \ot s_2 \mapsto m\ot n \ot s_1 s_2$. We claim that
\begin{equation}\label{tenssq1} \begin{split}
  &\hbox{\em $(M,b)$ and $(N,q)$ are regular}
  \\ & \qquad \iff    \hbox{\em $(M\ot_R N, b\ot_R q)$ is regular,}
\\
 & \hbox{\em $(M_1,b_1)$ and $(M_2,b_2)$ are regular}
\\ &\qquad \iff
      \hbox{\em $(M_1\ot_R M_2, b_1\ot_R b_2)$ is regular.}
\end{split} \end{equation}
We prove the first equivalence. To this end, let $\mu \co M^* \ot_R N^* \to (M\ot_R N)^*$ be the canonical isomorphism
(\cite[II, \S4.4]{BA}). It fits into the commutative diagram
\[\xymatrix{ M^* \ot_R N^* \ar[rr]^\mu_\simeq && (M\ot_R N)^* \\
       & M\ot_R N \ar[ul]^{\wdh b \ot \wdh{b_q}} \ar[ur]_{\wdh{b\ot q}}
}\] where $\wdh{b\ot q}$ is the adjoint of the tensor product form $b\ot q$.
If both $(M,b)$ and $(N,q)$ are regular, the commutative diagram above
shows that $\wdh{b\ot q}$ is an isomorphism. For the proof of the converse we use
\ref{bfLG}\eqref{bfLG-rech}: 
regularity can be checked after localization. Thus we can assume that both $M$ and $N$ are free $R$--modules. But then invertibility of $\wdh{b\ot q}$, i.e., of $\wdh b \ot \wdh{b_q}$ by the diagram above, implies invertibility of $\wdh b$ and $\wdh{b_q}$ by  \cite[III, \S8.6 (33)]{BA}.

\comments{(2022-02-28) Added sketch of proof in \ref{qfba}\eqref{qfba-tens} and \eqref{tenssq1}, taken from the octonion book. Before we only had the trivial direction $\implies$ in \eqref{tenssq1}. }

\comments{(2016-06-11) We show here, in \ref{qfba}\eqref{qfba-tens}:
For $i=1,2$ let $f_i \co M_i \to N_i$ linear maps between finitely generated
projective $R$--modules. Then $f_1$ and $f_2$ are invertible if and only if
$f_1 \ot f_2$ is invertible.
This should be somewhere in Bourbaki. \sm

(2022-07-12) Indeed, this is \cite[VI, \S6.3, Lem.~1]{BA8}.}
\sm

\textbf{Remark.} We will not use the tensor product of two quadratic forms $(M,q_M)$ and $(N,q_N)$ constructed in \cite[\S8.3]{BA3} or \cite{MicRev} and required to satisfy $(q_M\ot q_N)(m\ot n) = 2 q_M(m) \, q_N(n)$.
\sm
\enew

\begin{lem} \label{quadco} Let $(M, q)$ be a faithful quadratic space over $R$.
 Then $(M,q)$ is primitive, i.e., $\Span_R\{ q(m): m \in M\} = R$. In other words, $q(M)$ generates $R$ as ideal. In particular, if $r_1, r_2 \in R$ then
  \begin{equation} \label{quadco-aa1}
   r_1q(m) = r_2 q(m) \text{ for all $m\in M$} \quad \implies \quad     r_1 = r_2.
   \end{equation}
\end{lem}

\begin{proof} Otherwise, there exists a maximal ideal $\m$ containing $q(M)$, whence $q_{R_\m}$ has values in $\m R_\m$ and therefore $q_{\ka(\m)} \co M_{\ka(\m)} \to \ka(\m)$ is the null form. Since $q_{\ka(\m)}$ is nonsingular by  \ref{qf}\eqref{qfnsp}, this forces $M_{\ka(\m)} =
\{0\}$ which in turn implies $M_{R_\m} = \{0\}$, contradicting our
assumption that $M$ be faithful.\end{proof}

\comments{(2025-07-03) Moved the subsection on unimodular vectors to the section on ``Preliminaries", since it does not need/use quadratic forms}
\sm

Primitivity is obvious in case $(M,q)$ contains an $m$ with $q(m) \in R\ti$. Such an $m$ always exists for LG rings:  

\begin{lem}[{\cite[C.13]{GN-LG}}]\label{LGqdi}  Let $(M,q)$ be a quadratic space over an LG-ring $R$. Then the following are equivalent: 

\begin{enumerate}[label={\rm (\roman*)}]

\item\label{LGqdi-i} $q(M) \cap R\ti \ne \emptyset$, i.e., 
the {\em extended sphere $\wS_q$},  defined by
    \begin{equation}   \label{quadco-aa2}
 \wS_q = \{ m\in M : q(m) \in R\ti \},
\end{equation}
is non-empty. \sm 

\item \label{LGqdi-iii} $M$ is faithfully projective.
\end{enumerate} 
\end{lem} 

\begin{lem}\label{mx_lem}
  Let $(M, q)$ be a faithful quadratic module. For $x\in M$ with $q(x)\in R\ti$ put
\begin{equation} \label{mx_lem-1}
 M_x = \{m\in M :  b_q(m,x) = 0 \}, \quad q_x = q|_{M_x}.
\end{equation}

\begin{inparaenum}[\rm (a)]
  \item \label{mx_lema} Then $x$ is unimodular and $q|_{R\cdot x}$ is a nonsingular quadratic
  form. \sm


\item \label{mx_lemc} Assume $q$ is regular. Then there exists $y\in M$ satisfying $b_q(x,y) = 1$, and the quadratic form $q_x = q|_{M_x}$ is nonsingular. \end{inparaenum}
\end{lem}

\begin{proof}
\eqref{mx_lema} We have $q_S(x \ot 1_S) \in S\ti$ for any $S\in \Ralg$, in
particular for $S=R/\m$, $\m \in \Spec(R)$ maximal. This shows $x\ot 1_{R/\m}
\ne 0$, and therefore $x$ is unimodular by definition in \ref{unimod}. It follows that $R x$ is free of rank $1$. Since $q(x) \in R\ti$, nonsingularity of $q|_{Rx}$ is a consequence of  \ref{qf}\eqref{qfba-one}. \sm 


\eqref{mx_lemc} Since $x$ is unimodular, there exists $\vphi \in M^*$ such
that $\vphi(x) = 1$. By regularity of $q$, this $\vphi$ has the form $\vphi =
b_q(y, -)$ for some $y\in M$. To prove nonsingularity of $q_x$, we will show that $\rad\big( (q_x)_F\big) = 0$ for all fields $F\in \Ralg$. To do so, it is no harm to assume that $R=F$ is a field since all our assumptions are stable under base change. We note
\begin{align*}
  \rad(q_x) &= \{ m \in M_x: q(m) = 0 = b_q(m , M_x) \} \\
      & \subset M_x^\perp = \{m \in M : b_q(m , M_x) = 0 \} = (Rx)^{\perp\perp} = Rx
\end{align*}
where the last equality follows from \ref{orthLG}\eqref{orthLG-b}. 
Thus $\rad(q_x) \subset R x$. But $0 = q(rx) = r^2 q(x)$ implies $r = 0$.
Therefore $\rad(q_x) = 0$.
\end{proof}
\sm 

\textbf{Remarks.} 
If in \ref{mx_lem}\eqref{mx_lemc} we only suppose that $q$ be nonsingular, it does not follow that $q_x$ is nonsingular, even if we assume that there exists $y\in M$ satisfying $b_q(x,y) = 1$. Indeed, suppose $2R=0$. Then $b_q(x,x) = 2q(x) = 0$, so $x\in M_x \cap M_x^\perp  = \rad(b_{q_x})$ and therefore $q_x$ is not regular. But in case $M$ has constant odd rank, equivalently, $M_x$ has constant even rank, nonsingularity of $q_x$ is equivalent with regularity by \ref{qf}\eqref{quadfoe}.

\begin{leer}[\bf Isotropic vectors and (an)isotropic modules] \label{isotrop} Let $(M,q)$ be a quadratic module over $R$. We call $m\in M$ {\em isotropic\/} or an {\em isotropic vector\/} if $m$ is unimodular and satisfies $q(m) = 0$. 

Following standard terminology, a submodule $N$ of $(M,q)$ is called a {\em totally isotropic submodule} if $q(N) = 0$. Thus, a totally isotropic submodule contains an isotropic vector if and only if it contains a unimodular vector, which may not be the case. We call a quadratic module $(M,q)$ {\em isotropic\/}, if $(M,q)$ contains a  faithful complemented totally isotropic submodule. A quadratic module is {\em anisotropic\/} if it is not isotropic. 

The definition of an isotropic quadratic module used here is more general than the one of \cite{GN-Sp}, where $(M,q)$ was called isotropic if it contains an isotropic vector. See \eqref{isotrop-unifap} below for a comparison of the two definitions, and \eqref{isotrop-Knebusch} and \eqref{isotrop-group} for a ``justification'' of the change of definition. We emphasize that the results of \cite{GN-Sp} remain true with the new terminology since they concern quadratic spaces over semilocal rings where the two definitions coincide. 

We will later show  that the $R$--functor,  assigning to $T\in \Ralg$ the set of isotropic vectors in $(M,q)_T$, is represented by a scheme $\uY_q$, \ref{sreq}\eqref{sreq-y}. Some properties of $\uY_q$ are established in \ref{spreq} and \ref{aqyle}. Below we list some useful immediate facts concerning isotropy. 
\sm 

\begin{inparaenum}[(a)] \item\label{isotrop-a} ({\em Base change})
If $m$ is an isotropic vector of $(M,q)$, then $m \ot 1_S$ is an isotropic vector of  $(M,q)_S$ for any $S\in \Ralg$.
\sm

\item\label{isotrop-b} ({\em Direct products}) Let $R=R_0 \times \cdots \times R_n$ be a direct product of rings $R_i$ and let $(M,q) = (M_0, q_0) \times \cdots \times (M_n,q_n)$ be a quadratic $R$--module where the $(M_i, q_i)$ are  quadratic $R_i$--modules, see \ref{qf}\eqref{qf-redc}. Then $m=(m_0, \ldots, m_n)\in M$ is an isotropic vector of $(M,q)$ if and only if every $m_i$ is an isotropic vector of $(M_i, q_i)$. 
    More generally, $(M,q)$ is an isotropic $R$--module if and only every $(M_i, q_i)$ is an isotropic $R_i$--module. 
   \  \sm

\item\label{isotrop-bc} (\cite[7.13]{EKM} for fields)
{\em Let $(M,q)$ be a quadratic module and assume $v\in M$ satisfies $q(v) = 0$. We consider the following conditions \eqref{isotrop-bci}--\eqref{isotrop-bciii}: 

\begin{inparaenum}[\rm (ci)]
\quad \item \label{isotrop-bci} the linear form $b_q(v, \cdot)$ is surjective,

\quad \item \label{isotrop-bcii} there exists $w\in M$ such that $(v,w)$ is a hyperbolic pair, see {\rm \ref{qf}\eqref{qf-hyp}}. 

\quad \item \label{isotrop-bciii} $v$ is unimodular, hence isotropic. 
\end{inparaenum}

\noindent Then 
\[ \eqref{isotrop-bci} \iff \eqref{isotrop-bcii} \implies \eqref{isotrop-bciii}.
\] 
If $q$ is nonsingular, then all three conditions are equivalent. 
}
\sm 

Indeed, if \eqref{isotrop-bci} holds, there exists $y\in M$ such that $b_q(v,y) = 1$. Then $(v, y - q(y) v)$ is a hyperbolic pair. The implication \eqref{isotrop-bcii} $\implies$ \eqref{isotrop-bci} holds by definition of a hyperbolic pair, and \eqref{isotrop-bcii} $\implies$ \eqref{isotrop-bciii} follows from \ref{unimod}\ref{unimodii}. 

Suppose that $q$ is nonsingular and \eqref{isotrop-bciii} holds, but not \eqref{isotrop-bci}. Then the ideal $b_q(v,M)\ideal R$ is proper, say contained in a maximal ideal $\m \ideal R$. Let $k=R/\m$. Then $b_{q_k}( v\ot 1_k, M_k) = b(v, M) \ot_R 1_k = 0$ and also $q_k (v \ot 1_k) = q(v) \ot 1_k = 0$. Hence $v\ot 1_k\in \rad(q_k) = \{0\}$ by nonsingularity of $q$. But $v\ot 1_k = 0$ contradicts  \ref{unimod}\ref{unimodiv}. Thus \eqref{isotrop-bcii} holds. \sm

\inparcom{(2025-03-24) \eqref{isotrop-bc} is generalized in \ref{quadrepII}. }

\item\label{isotrop-d} Specializing \eqref{isotrop-bc} and using \eqref{qf-perp1}, we obtain:    
    {\em If $(M,q)$ is a quadratic space and $v\in M$ is an isotropic vector, then $M$ contains $w\in M$ such that $(v,w)$ is a hyperbolic pair. Hence, $H = Rv + Rw$ is free of rank $2$ and isometric to a hyperbolic plane. Moreover, 
 \begin{equation}\label{isotrop-d0}   
   (M,q) = (H, q|_H) \perp (H^\perp, q|_{H^\perp})
\end{equation} 
with $(H^\perp, q|_{H^\perp})$ being nonsingular.} 
In particular, 
\begin{equation}  \label{isotrop-00}
 \rank M \ge 2,  
\end{equation}
and if $(M,q)$ is a quadratic space of constant rank $2$, then 
\begin{equation} \label{isotrop-d1}
     \text{$(M,q)$ contains an isotropic vector} \quad \iff \quad (M,q) \cong \HH(R).
     \end{equation}
 Thus, in this case $M$ is free of rank $2$. \sm 

\item\label{isotrop-unifap} If $(M,q)$ contains an isotropic vector, say $m\in M$, then $(M,q)$ is an isotropic module, since $Rm$ is a free, hence faithful, and  complemented totally isotropic submodule.     

Conversely, {\em a quadratic module $(M,q)$ over a unimodular ring $R$ is isotro\-pic if and only if $(M,q)$ contains an isotropic vector.}
Indeed, a faithful complemented submodule is faithfully projective and therefore contains a unimodular vector by definition of a unimodular ring in \ref{unifap}. 
\sm 

\item\label{isotrop-Knebusch} ({\em Equivalence with Knebusch's definition of an isotropic module}) Suppose $\Spec(R)$ is connected. Then {\em $(M,q)$ is isotropic if and only if $(M,q)$ contains a non-zero complemented totally isotropic submodule.} Indeed, a complemented submodule of $M$ is finite projective. In case $\Spec(R)$ is connected, it is therefore faithful if and  only if it is non-zero. 
    
    The canonical analogue of this characterization of isotropy for symmetric bilinear forms is Knebusch's definition of an isotropic symmetric bilinear module in \cite{Knebusch-habil} and \cite{Knebusch-Queens}.     
\sm 

\item\label{isotrop-group}  Let again $(M, q)$ be a quadratic space over $R$ and let $\uSO(q)$ be the associated reductive group scheme over $\Spec(R)$, see \ref{sogsc} for a review. In Proposition~\ref{prop_isotropic} we will prove: 
\begin{equation}\label{isotrop-group1} \begin{split}
&\text{\em A quadratic space $(M,q)$ with $\rank M \ge 3$ is isotropic 
} \\
&\quad \text{\em if and only if $\uSO(q)$ is an isotropic reductive group scheme.}  
\end{split}\end{equation}    
\end{inparaenum}
\end{leer}

\begin{lem}\label{quapreI}
  Let $(M,q)$ be a quadratic module and let $U$ and $W$ be  
  submodules of $M$. We will consider the $R$--linear map
\begin{equation}\label{quapreI-1}
  \be_{W,U} \co W \to U^*, \quad w \mapsto \big(u \mapsto b_q(w,u)\big)
\end{equation}
and its analogue $\be_{U,W} \co U \to W^*$. \sm

\begin{inparaenum}[\rm(a)]
 \item \label{quapreIa} Assume $U$ is finite projective and let $\be_{W,U}^* \co U^{**} \to W^*$ be the dual map of $\be_{W,U}$ and let $\can_U \co U \to U^{**}$ be the canonical isomorphism. Then $\be^*_{W,U} \circ \can_U = \be_{U,W}$. \sm

\item \label{quapreId} Let $S\in \Ralg$ and suppose that $U$ and $W$ are  complemented, which allows us to view $U \ot_R S $ and $W\ot_R S$ as complemented submodules of $M\ot_R S$. Denoting by $\be_{W,U} \ot 1_S \co W \ot_R S \to U^* \ot_R S$  the base change of $\be_{W,U}$ and by $\nu_U \co U^* \ot_R S \to (U\ot_R S)^*$ the canonical isomorphism, 
    we have  \[ \nu_U \circ (\be_{W,U} \ot 1_S) = \be_{W\ot_R S,\,  U\ot_R S}.\]
\end{inparaenum}
\end{lem}

\begin{proof} \eqref{quapreIa} and \eqref{quapreId} are straightforward from the definitions. \end{proof}
\lv{
Let $u\in U$, thus $\al = \can_U(u) = \big(\rho \mapsto \rho(u)\big)$ for $\rho \in U^*$. Also recall that $\be^*_{W,U}(\psi) = \psi \circ \be_{W,U}$ for $\psi \in U^{**}$. Hence $(\be^*_{W,U}\circ \can_U)(u) = \be_{W,U}^*(\al) = \al \circ \be_{W,U}$. We evaluate at $w\in W$ and get $\big(\be^*_{W,U}\circ \can_U\big)(u)(w) = \al \big( \be_{W,U}(w)\big) = \al \big( u \mapsto b_q(w,u)\big) = b_q(w,u) = \big(\be_{U,W}(u)\big)(w)$.

Proof of \eqref{quapreId}:
 \begin{align*}
  & \Big( \big(\nu_U \circ (\be_{W,U}\ot 1_S)\big)(w\ot s)\Big)(u\ot s_1) =
   \big(\nu_U (\be_{W,U}(w) \ot s) \big)(u \ot s_1)
 \\& \quad = \be_{W,U}(w)(u) \ot ss_1 = b_q(w\ot s, u\ot s_1)
    = \big(\be_{W\ot S, U\ot S}(w\ot s)\big)(u \ot s_1)
 \end{align*}
}
\sm

Proposition~\ref{quadrepII} improves \cite[Prop.~1.5]{GN-Sp}; it is proven in \cite[Lem.~(3.1)]{Bass-69} and again in \cite[I, Thm.~3.6]{Ba} for regular forms,  and in \cite[8.10]{EKM} for nonsingular forms over fields. 
A related result is \cite[I, \S3, Prop.~3a]{Knebusch-Queens}, which deals with regular symmetric bilinear modules. In that case, one can only conclude that $(U\oplus V, b_{U \oplus V})$ is metabolic. 

\begin{prop} \label{quadrepII} Let $(M,q)$ be a quadratic space and let $U \subset M$ be a totally isotropic and complemented submodule.  Then there exists a totally isotropic and complemented submodule $V\cong U^*$ satisfying
\begin{enumerate}[label= {\rm (\roman*)}]
  \item\label{quadrepII-ci} $(U \oplus V, q|_{U \oplus V}) \cong \HH(U)$,

  \item \label{quadrepII-cii} $M = U \oplus V \oplus (U \oplus V)^\perp$,

  \item \label{quadrepII-ciii} $U^\perp = U \oplus (U \oplus V)^\perp$, in particular $U^\perp$ is complemented by $V$.
\end{enumerate}
\end{prop}

\begin{proof}  The proof is divided in steps \eqref{quadrepII-cI}--\eqref{quadrepII-cIV}, using the notation of Lemma~\ref{quapreI}. First, in \eqref{quadrepII-cI}, we show that all claims follow, if we know that
 \[ \be_U = \be_{M,U} \co M \to U^*, \quad m \mapsto \big( u \mapsto b_q(m,u)\big)\]
is surjective. We then prove surjectivity of $\be_U$ in increasing generality in steps \eqref{quadrepII-cII}--\eqref{quadrepII-cIV}. In \eqref{quadrepII-cII} we treat the case that $R$ is a field or that $q$ is regular. Then, using \eqref{quadrepII-cII}, we can deal in \eqref{quadrepII-cIII} with the case of a local ring. Finally, in \eqref{quadrepII-cIV}, we handle the general case.
\sm

\begin{inparaenum}[(I)] \item \label{quadrepII-cI} As noted above, we will show in \eqref{quadrepII-cI} that it is sufficient to prove that $\be_U $ is surjective.
Indeed, if this holds, we get an exact sequence
\[ 0 \longto U^\perp \longto M \xrightarrow{\; \be_U\;} U^* \longto 0 \]
which is split-exact because $U^*$ is projective. Hence $U^\perp$ is complemented, say $M = U^\perp \oplus W$, and the restriction $\be_{W,U} = \be_U|_W \co W \to U^*$ is an isomorphism. Since $\be_{W,U}^*\circ \can_U = \be_{U,W}$ by Lemma~\ref{quapreI}\eqref{quapreIa}, we then know that also $\be_{U,W}$ is an isomorphism. 

We can then proceed as in the intermediate step (I) in the proof of \cite[Prop.~1.5]{GN-Sp}, which is itself inspired by the proof of \cite[I, Thm.~(3.6)]{Ba}. We repeat the short proof for the convenience of the reader.

By \cite[I, (1.7)]{Ba}, we can choose a not necessarily symmetric bilinear form $b_0$ satisfying $b_0(m,m) = q(m)$ for all $m\in M$. Since $\be_{U,W}$ is an isomorphism, for fixed $w\in W$ the linear form $w' \mapsto b_0(w,w')$ on $W$ is represented by a unique $u_w \in U$, i.e.,
$b_q(u_w, w') = b_0(w, w')$ holds for  all $w'\in W$. Because of uniqueness, the map $W \to U$, $w \mapsto u_w$, is $R$--linear; it is injective because $U \cap W = 0$. Then $V = \{ w - u_w : w \in W\}\cong W$ is a totally isotropic submodule: indeed, since $q(u_w)=0$ we have $q(w - u_w) = q(w) - b_q(w, u_w) = b_0(w,w) - b_q(u_w, w) = 0$.
Moreover, $U \cap V = 0$ and the map $\be_{U,V}$ is an isomorphism. Hence $(U\oplus V, q|_{U \oplus V}) \cong \HH(U)$, proving \ref{quadrepII-ci}.

Because $q|_{U \oplus V}$ is hyperbolic, thus regular, \ref{quadrepII-cii} is a consequence of \ref{orthLG}\eqref{orthLG-a}. Finally, $U^\perp \supset U \oplus (U \oplus V)^\perp$. Hence, writing an arbitrary $m\in M$  as $m = u + v + x$ with $u\in U$, $v\in V$ and $x\in (U \oplus V)^\perp$ by making use of  \ref{quadrepII-ci}, we get $m\in U^\perp \iff 0 = b_q(U, m) = b_q(U, v) \iff v=0$ because $\be_{U,V}$ is an isomorphism. Thus, also \ref{quadrepII-ciii} is satisfied.
\sm

\item \label{quadrepII-cII} We show in \eqref{quadrepII-cII} that the map $\be_U$ is surjective in case $q$ is regular or $R$ is a field.

 Indeed, since $U$ is complemented, the restriction map $\pi_U \co M^* \to U^*$, $\vphi \mapsto \vphi|_U$ is surjective.
   Let $\be_{M,M} \co M \to M^*$ be the adjoint map of $b_q$. Then $\be_U = \pi_U \circ \be_{M,M}$. Hence, if $q$ is regular, i.e., $\be_{M,M}$ is an isomorphism, then $\be_U$ is surjective.

Let now $R$ be a field. The kernel of the map
 \[ \be_{U,M} \co U \to M^*, \quad u \mapsto \big( m \mapsto b_q(m,u)\big)\]
is $U \cap M^\perp \subset \rad(q)$ since $U$ is totally isotropic. But if $R$ is a field, then $\rad(q) =\{0\}$ by nonsingularity of $q$. Thus $\be_{U,M}$ is injective. Hence its dual map $\be^*_{U,M} \co M^{**} \to U^*$ is surjective. From \ref{quapreI}\eqref{quapreIa}  we know that $\be_U = \be_{M,U} = \be_{U,M}^* \circ \can_U$ where $\can_M\co M \to M^{**}$ is the canonical isomorphism. Hence $\be_U$ is surjective too. \sm

\item \label{quadrepII-cIII} In this step we prove that the map $\be_U$ is surjective if $R$ is a local ring, say with maximal ideal $\m$ and residue field $k$.

Indeed, if $R$ is local, $U$ is free of finite rank, so that we can choose a basis $(z_1, \ldots, z_n)$ of $U$. We denote the reduction mod $\m M$ by $m \mapsto \ol m$. Then $\ol  U$ is a totally isotropic subspace of the quadratic space $(\ol M, \ol q)$ with basis $(\bar z_1, \ldots, \bar z_n)$. By \eqref{quadrepII-cII} we know that Proposition~\ref{quadrepII} is true over $k$. So we can choose a subspace $Z\subset \ol M$ such that $(\ol U \oplus Z, \ol q |_{\ol U \oplus Z}) \cong \HH(\ol U)$. Furthermore we can choose $z_{n+1}, \ldots, z_{2n}$ such that $(\ol z_{n+1}, \ldots, \ol z_{2n})$ is a basis of $Z$. Observe that $\det \big(b_q(\ol z_i, \ol z_j)\big) \ne 0$ by regularity of $\ol U \oplus Z$, hence $\det \big(b_q(z_i, z_j)\big) \in R\ti$.
This implies that $z_1, \ldots , z_{2n}$ are free in $M$: any relation $r_1 z_1 + \cdots + r_{2n} z_{2n}= 0$ with $r_i \in R$ gives rise to a homogeneous linear system for the unknowns $r_1, \ldots , r_{2n}$ with invertible coefficient matrix $\big( b_q(z_i, z_j)\big)$. Putting $W=\Span\{ z_{n+1}, \ldots, z_{2n}\}$, we get that $q|_{U \oplus W}$ is regular. By \eqref{quadrepII-cII} we then know that
\ref{quadrepII} holds for $U \oplus W$, in particular $q_{U\oplus W}$ is hyperbolic and $U \oplus W$ is complemented. Hence $\be_U$ is surjective. \sm
\sm

\item\label{quadrepII-cIV} Finally, to prove that $\be_U$ is surjective in general, it suffices to prove surjectivity for the localization of $\be_U$ in all maximal ideals. 
   By Lemma~\ref{quapreI}\eqref{quapreId} this follows from surjectivity of $\be_U$ in case $R$ is a local ring, which we have dealt with in \eqref{quadrepII-cIII}.\end{inparaenum}
\end{proof}

\subsection{Hyperbolic pairs of submodules}\label{hps} Let $(M, q)$ be a quadratic space. Generalizing \ref{qf}\eqref{qf-hyp}, we call a pair $(E,F)$ of submodules of $M$ a  {\em hyperbolic pair of submodules\/} if $q(E) = 0 = q(F)$ and $F \cong E^*$ via the adjoint map:
    \[ \be \co F \simlgr E^*, \quad f \mapsto b_q(f, -)|_E. \]
Observe that in this case the analogous map $\be_E \co E \simlgr F^*$ is  also an isomorphism since $\be_E = \be^* \circ \can$. Moreover, $E+F$ is a direct sum and  
\begin{equation}\label{hps1}
   (E \oplus F, q|_{E\oplus F}) \simlgr \HH(E), \quad e + f \mapsto \be(f) + e
\end{equation}
is an isometry. Finally, by \ref{orthLG}\eqref{orthLG-a} we have
\begin{equation}  \label{hps1a}
M = (E \oplus F) \perp (E \oplus F) ^\perp.
\end{equation}
Proposition~\ref{quadrepII} says that every complemented totally isotropic submodule $U$ of $(M,q)$ is part of a hyperbolic pair of submodules. \sm 

Let $(E,F)$ and $(E',F')$ be two hyperbolic pairs of submodules of a quadratic module $(M,q)$. {\em If $E$ and $E'$ are isomorphic as $R$--modules, then there exists an isometry
\begin{equation}\label{hps2}
 g \co (E \oplus F, q|_{E\oplus F}) \simlgr (E'\oplus F', q|_{E'\oplus F'}) 
\end{equation}
mapping $(E,F)$ to $(E', F')$.} Indeed, 
this follows by combining the isometries \eqref{hps1} and \eqref{qfba-hyp1}. 
If $R$ has Witt cancellation = extension of isometries, \ref{canqf}, the isometry $g$ extends to an isometry $\hat g$ of $(M,q)$. In this case, $\hat g$ induces an isometry of the orthogonal modules, 
\begin{equation}   \label{hps3}
 (E \oplus F)^\perp \simlgr (E' \oplus F')^\perp.
\end{equation} 
However, $g$ does in general not extend to an isometry, as the following example shows. 


\subsection{Example: The isometry \eqref{hps2} does not extend in general.} \label{hpsex}
Let $\HH$ be the real division quaternions. Its norm is the real  quadratic form $q_4 = \lan 1,1,1,1\ran$ of $4$ squares. Now consider $R=\RR[X,Y]$, the polynomial ring in 2 variables over $\RR$, and the Azumaya $R$--algebra $\HH[X,Y] = \HH \ot_\RR R$. Its norm is $q_{4, R}$, the base change of $q_4$ to $R$. 

Let $P$ be the non-free projective ideal of $\HH[X,Y]$, constructed in \cite[Prop.~1]{OS}. The endomorphism algebra $B = \End_{\HH[X,Y]}(P)$ is an Azumaya $R$--algebra of rank $4$. Let $q\co B \to R$ be its norm. It is known that $q$ is not extended from $\RR$ (\cite[Proof of Pro.~2.1]{KO77}). In particular, $q$ and $q_{4,R}$ are not isometric. 

However, by Karoubi's Theorem \cite[VII, Thm.~(4.1.3)]{K}, the quadratic forms $q$ and $q_{4,R}$ are stably isometric, i.e., there exist finite projective, hence free $R$--modules $E$ and $E'$ and an isometry $f \co q \perp \HH(E) \simlgr q_{4,R} \perp \HH(E')$. Since $E$ and $E'$ are isomorphic as $R$--modules, we can apply \eqref{hps2} and get an isometry $g\co f\me\big(\HH(E')\big)= \HH(f\me(E')) \simlgr \HH(E)$. Now suppose that $g$ extends to an isometry of $q \perp \HH(E)$. Then, by \eqref{hps3}, also the orthogonal spaces $\HH(E)^\perp$ and $f\me\big(\HH(E'))^\perp$ are isometric, i.e., $q$ and $q_{4,R}$ are isometric, contradiction. 
\ms 

Corollary \ref{carhyp}, 
and parts \eqref{carhyp-cor-a} and \eqref{carhyp-cor-aa} of \ref{carhyp-cor} are proven in \cite[Ch.~5, \S2]{BR} for regular quadratic modules. Corollaries \ref{carhyp} and \ref{carhyp-cor} generalize \cite[Cor.~1.6 and 1.7]{GN-Sp}. The case of a regular quadratic module is \cite[I, Thm.~4.6]{Ba}.  

\begin{cor}[Characterization of hyperbolic spaces] \label{carhyp}
Let $(M,q)$ be a quadratic $R$--space. Then the following are equivalent: \sm

\begin{enumerate}[label={\rm (\roman*)}]
  \item \label{carhyp-i} $(M,q)$ is hyperbolic;

  \item \label{carhyp-ii} $M$ admits a direct summand $L$ satisfying $q(L) = 0$
  and $2 \rank_\p L = \rank_\p M$ for all $\p \in \Spec(R)$;

  \item \label{carhyp-iia} $M$ admits a direct summand $L$ satisfying $q(L) = 0$
  and $2 \rank_\m L \ge \rank_\m M$ for all maximal $\m \in \Spec(R)$;

  \item \label{carhyp-iii} $M$  admits a direct summand $L$ satisfying $q(L) = 0$ and $L= L^\perp$.
\end{enumerate}
In this case $(M,q) \cong \HH(L)$.
\end{cor}

\begin{proof} The implications  \ref{carhyp-i} $\implies$ \ref{carhyp-ii} $\implies$ \ref{carhyp-iia} being obvious because $\rank_\p L = \rank_\p L^*$ for $\p \in \Spec(R)$, let us assume \ref{carhyp-iia} and prove \ref{carhyp-i}. By Proposition~\ref{quadrepII} there exists a submodule $V \subset M$ such that $(L \oplus V, q|_{L \oplus V})\cong \HH(L)$ is hyperbolic and $M = (L \oplus V) \oplus (L \oplus V)^\perp$. Since $V \cong L^*$ as $R$--modules, $\rank_\m (L \oplus V) =  \rank_\m M$ by \ref{carhyp-iia}, whence $L \oplus V = M$
and $L^\perp = L \oplus (L^\perp \cap V)$ with $L^\perp \cap V \cong  \{\vphi \in L^* : \vphi(L) = 0 \} 
= \{0\}$.
The proof of \ref{carhyp-iii} $\implies$ \ref{carhyp-i} follows the same pattern.  
\lv{
Again by Lemma~\ref{lift-qf}, there exists a submodule $V\subset M$ such that $L \oplus V\cong \HH(L)$ and $M = (L \oplus V)\oplus (L \oplus V)^\perp$. Since $(L\oplus V)^\perp \subset L^\perp = L$, we get $M = L \oplus V$.

Other proof: $(L\oplus V)^\perp$ is finite projective of constant rank $0$, so is itself zero. }
\end{proof}

\subsection{Lagrangians} \label{Lag}
A {\em Lagrangian\/} of a quadratic space $(M,q)$ is a complemented totally isotropic submodule $L$ with $L = L^\perp$. Corollary~\ref{carhyp} shows that the following three conditions \ref{Lagi}--\ref{Lagiii} are equivalent for a complemented submodule $L$ with $q(L) = 0$:
\begin{enumerate}[label={\rm (\roman*)}]
  \item\label{Lagi} $L$ is a Lagrangian, 
   \item\label{Lagii} $2 \rank_\p L = \rank_\p M$ for every $\p \in \Spec(R)$,
   \item \label{Lagiii} $2 \rank_\m L \ge \rank_\m M$ for every maximal $\p \in \Spec(R)$. 
\end{enumerate} 
We list some facts: \sm

\begin{inparaenum}[(a)] \item \label{Lag-a} By Proposition~\ref{quadrepII}, {\em any  Lagrangian $L$ is part of a hyperbolic pair $(L,L')$ of submodules of $(M,q)$\/} in the sense of \ref{hps}.\sm

\item \label{Lag-b} {\em Any $R$--module isomorphism of two Lagrangians $L$ and $L'$ of  $(M,q)$ extends to an isometry of $(M,q)$.} This follows from \eqref{Lag-a} and the isometry \eqref{hps2}.\sm

\item \label{Lag-c} {\em If $R$ is \new a unimodular \enew ring, then the orthogonal group $\orth(q)$ acts transitively on the set of Lagrangians of $(M,q)$. }

Indeed, by \ref{Lagii}, two Lagrangians $L$ and $L'$ of $(M,q)$ have the same rank function $\Spec(R) \to \NN$. Hence $L$ and $L'$ are isomorphic as $R$--modules by \eqref{unimod-b2}. Then \eqref{hps2} says that there exists an isometry of $(M,q)$ mapping $L$ onto $L'$.  \sm

\item \label{Lag-d} The group scheme $\uO(q)$ acts on the scheme of Lagrangians $\uL(q)$, defined in \ref{hrq}, in the obvious way. Using the terminology of \ref{transi}, we get: {\em The action of $\uO(q)$ on the Lagrangian quadric $\uL(q)$ is transitive on the big affine Zariski site of $\Spec(R)$.} 
    
    Indeed, since the schemes $\uO(q)$ and $\uL(q)$ are finitely presented, an application of the surjectivity criterion~\ref{surlem} reduces the proof to surjectivity in the fibres, where it follows for example from \eqref{Lag-c}.  
\sm 

\item It is in general not true that $\SO(q)$ acts transitively on the Lagrangians of a quadratic space, see \ref{eor}\eqref{eor-b} for a transitivity criterion.      
\end{inparaenum}

\begin{cor}[{\cite[I, (4.7.iii)]{Ba}}]\label{lem_metabolic} 
  For each regular quadratic module $(M,q)$ and bilinear module $(U,b)$ the $R$--quadratic space $\MM(U,b) \ot_R (M,q)$ is hyperbolic
 \begin{equation} \label{lem_metabolic1} \begin{split}
   \MM(U,b) \ot_R (M,q) &\cong \HH(U \ot_R M), \quad\text{ in particular,} \\
    \HH(U) \ot_R (M,q) &\cong \HH(U\ot_R M).
  \end{split}\end{equation}

\end{cor}

\begin{proof} The quadratic module \[ \MM(U,b) \ot (M,q) = \big( (U\ot_R M) \oplus (U^* \ot_R M), Q\big) = (V,Q)\]  is regular by \ref{qfba-tens}.  According to \eqref{tenssq0} we have
\begin{align*}
  Q(\vphi \ot m) &= b_{\MM(U)}(\vphi, \vphi)\, q(m) =0 \quad \text{and}\\
b_Q( \vphi \ot m, \vphi' \ot m') &= b_{\MM(U)}(\vphi, \vphi')\, b_q(m,m') = 0
\end{align*}
for $\vphi \in U^*$ and $m\in M$. It follows that $N=U^* \ot M$ satisfies the conditions of \ref{carhyp}\ref{carhyp-ii}. By loc.\ cit., we then get $(V,Q) \cong \HH(U^* \ot M) \cong \HH(U \ot M)$ because $U^* \ot M \cong U \ot M$ as $R$--modules.
\end{proof}

\subsection{Witt cancellation, extension of isometries} \label{canqf} 
We say that {\em a ring $R$ has Witt cancellation\/} if the following condition holds: 
{\em whenever $q_1$ and $q_2$ are nonsingular quadratic forms, and $q$ and $q'$ are {\em regular\/} quadratic forms, then}
\begin{equation}\label{canqf1} q_1 \perp q \cong q_2 \perp q' \;\text{ and }\; q\cong q'  \quad \implies \quad q_1 \cong q_2. \end{equation}

Witt cancellation is equivalent to the {\em extension property of isometries\/}:
\begin{equation} \label{canqf2}
  \begin{split}
  & \text{\em every isometry between regular submodules $M_1\simlgr M_2$}
   \\ &\text{\em of a quadratic space $(M,q)$ extends to an isometry of $M$.} 
  \end{split}
\end{equation}
 
Indeed, the proof of \eqref{canqf1} $\iff$ \eqref{canqf2} in \cite[(3.1), (3.2)]{Kneser} for fields works equally well for quadratic forms over rings. We present it here for the convenience of the reader. 

But let us first recall that a {\em regular submodule\/} of a quadratic module $(M,q)$ is a submodule $U\subset M$ such that $q|U$ is regular. By \ref{orthLG}\eqref{orthLG-a} and \ref{qf}\eqref{qf-regLG}, $U$ is complemented, finite projective and $M=U \perp U^\perp$. Moreover, by \eqref{qf-perp1}, the quadratic form $q|_{U^\perp}$ is nonsingular if and only if $q$ is nonsingular. 

Suppose \eqref{canqf1} and let $f \co (U, q_U) \simlgr (V, q_V)$ be an isometry between regular submodules of $(M,\wdh q)$. Since $\wdh q|_{U^\perp}$ and $\wdh q|_{V^\perp}$ are nonsingular quadratic forms, \eqref{canqf1} says that there exists an isometry $g \co U^\perp \simlgr V^\perp$. The isometry $f\perp g $ extends $f$. 

Conversely suppose \eqref{canqf2}, and let $f \co (M_1, q_1) \perp (U, q) \simlgr (M_2, q_2) \perp (U', q')$ and $g\co (U,q) \simlgr (U',q')$ be isometries. Then $h = g \circ (f|_U)\me $ is an isometry from $f(U)$ to $U'$. Let $\wdh h$ be an isometry of $M_2 \perp U$ extending $h$. Then 
\begin{align*} M_1 &\cong f(M_1) = f(U^\perp) = f(U)^\perp 
= \big((h \circ g)(U)\big)^\perp \\ & = \big( (\wdh h \circ g)(U)\big)^\perp = \big(\wdh h(U')\big)^\perp = \wdh h\big( (U')^\perp\big) = \wdh h(M_2) \cong M_2
\end{align*} 
shows $(M_1, q_1) \cong (M_2, q_2)$. \sm 

\textbf{Examples and non-examples.} \begin{inparaenum}[(a)]\item\label{canqf-a} The main example of rings with Witt cancellation are LG rings (\cite[5.5]{GN-LG}). \sm  

\item \label{canqf-b} Using \ref{qf}\eqref{qf-redc}, it is easily seen that a direct product ring $R=R_1 \times \cdots \times R_n$ has Witt cancellation if and only if all $R_1, \ldots , R_n $ have Witt cancellation.%
     \sm 
  
\item\label{canqf-iso} ({\em Witt cancellation and base change}) Let $\qs_R$ be the full subcategory of the category $\qm_R$ of \ref{qf}\eqref{qf-bc} whose objects are quadratic spaces over $R$ and whose morphisms are isometries. For $T\in \Ralg$ define $\qs_T$ accordingly.  
       
    Suppose that $T$ is faithfully flat, and that the restriction $\frb_{RT} \co \qs_R \to \qs_T$ of the base change functor of \ref{qf}\eqref{qf-bc}  is essentially surjective on objects and that isometries in $\qs_T$ descend to $\qs_R$. 
Since by \ref{bfLG}\eqref{bfLG-rech} the quadratic form $q_T$ is regular if and only if $q$ is regular and since $\frb_{RT}$ preserves orthogonal sums, the ring $R$ has Witt cancellation if and only if so does $T$. 
   
   An example of such a situation is $R=k$ a field of characteristic $\ne 2$ and $T=k[X]$, the polynomial ring over $k$. That $\frb_{RT}$ satisfies the conditions above, is Harder's Theorem. Thus, $k[X]$ has Witt cancellation (\cite[VII, (2.3.2)]{K}). 
   We point out that $k[X]$ is not an LG ring, \ref{revLG}\eqref{revLG-non}, while $k$ is. \sm 
   
\lv{Note on Harder's Theorem: We will show in \ref{wide} that every quadratic space has a Witt decomposition. That $\frb_{RT}$ is essentially surjective then boils down to proving that every hyperbolic and every anisotropic space is extended from $k$. The hyperbolic case is clear because finite projective modules over $k[X]$ are free. A proof of the anisotropic case is given in \cite[13.4.3]{Knebusch-habil} for symmetric bilinear forms. It applies here because of our concept of an anisotropic module coincides with the one used in \cite{Knebusch-habil} by \ref{isotrop}\eqref{isotrop-Knebusch}. It remains to show that for quadratic spaces $M$ and $N$ over $k$, every $k[X]$--linear isometry $M\ot_k k[X] \to N \ot_k k[X]$ gives rise to an isometry $M \to N$. This follows from passing from $k[X]$ to $k[X]/(X) \cong k$.} \sm

\item \label{canqf-non} ({\em A non-example}) The example \ref{hpsex} together with the remark in \ref{hps} shows that $\RR[X,Y]$ does not have Witt cancellation. \sm 

\item \label{canqf-c} ({\em Another non-example}) Witt cancellation is in general not true if $q$ and $q'$ in \eqref{canqf1} are not regular, as the following examples shows, lifted from \cite[(8.7)]{EKM}. 
    
Let $R$ be a ring with $2R=0$ and let $a\in R\ti$. Then, for any $b\in R$, 
\begin{equation}  \label{canqf-c1}
 [a,b]\perp \lan a \ran_q \quad \cong \quad \HH \perp \lan a \ran_q. 
\end{equation}  
Indeed, let $M$ be the $R$--module with basis $(x,y,z)$ and let $q$ be the quadratic form given by $q(r_1 x + r_2 y) = a r_1^2 + r_1 r_2 + b r_2^2$, $q(z) = a$ and $(Rx \oplus Ry) \perp Rz$, i.e., the quadratic module of the left-hand side of \eqref{canqf-c1}. One easily verifies that $\HH= R(x+z) \oplus R\big(y-b(x+z)\big)$ is a hyperbolic plane for which $\HH^\perp = Rz$. Thus, the quadratic module has the description of the right-hand side of \eqref{canqf-c1}. Since $[a,b] \not \cong \HH$ in general, one cannot cancel $\lan q \ran_q$ in \eqref{canqf-c1}.    
\end{inparaenum}
 
\comments{(2025-07-06) Below are parts of your comments regarding rings that have Witt cancellation, but are not LG rings, put here to avoid losing them:
\sm 

\pcomments{(2) Je ne sais plus o\`u cela en est avec l'\'etudiant de Seidon qui devait \'etudier le cas de l'anneau A des fonctions continues \`a valeur r\'eelle d'un vari\'et\'e diff\'erentiable X (compacte ou même paracompacte).
D'apr\`es Swan-Serre, on a en effet une \'equivalence de cat\'egories entre la cat\'egorie des A-modules projectifs de type fini et celle des fibr\'es vectoriels.
Cela donne lieu \`a une \'equivalence de cat\'egories entre la cat\'egorie des formes A-quadratiques r\'eguli\`eres et celle des fibr\'es quadratriques.
On cherche donc un espace topologique X telle que la cat\'egorie des fibr\'es quadratiques soit tr\`es simple, je ne suis pas sur que cela existe.}}
\ms 

We will establish some results for rings with Witt cancellation, but not necessarily LG rings, in \ref{carhyp-cor}\eqref{carhyp-cor-b} and \ref{equi}. 


\begin{cor}\label{carhyp-cor} Let $(M,q)$ and $(M',q')$ be quadratic modules and {\em assume that $(M,q)$ is  regular}. \sm

\begin{inparaenum}[\rm (a)] \item {\rm \cite[I, (4.7.i)]{Ba}}
\label{carhyp-cor-a} If $(M,q)$ and $(M', q')$ are isometric, the quadratic module $(M , q) \perp (M', -q')$ is hyperbolic: $(M,q) \perp (M',-q') \cong \HH(M)$.  \sm

\item \label{carhyp-cor-aa} $(M',q')$ is regular if and only if $(M',q') \perp (M', -q') \cong \HH(M')$. \sm

\item\label{carhyp-cor-b} If $R$ has Witt cancellation, for example if $R$ is an  LG ring, and $(M',q')$ is nonsingular, then%
\[ (M,q) \perp (M', -q') \cong \HH(M) \quad 
      \iff \quad (M,q) \cong (M', q'). \]%
\end{inparaenum}
\end{cor}

\begin{proof} \eqref{carhyp-cor-a} Let $f \co (M,q) \to (M',q')$ be an isometry. The quadratic form $q\perp (-q')$ is regular by \ref{bfLG}\eqref{qfba-sum} since $q$ and $q'$ are regular. The diagonal submodule $U = \{\big(m,f(m)\big) : m \in M\}\subset M \oplus M'$ is complemented by $\{ (m,0): m \in M\}$, has $q(U) = 0$  and thus satisfies \ref{carhyp}\ref{carhyp-ii}. Hence $(M,q) \perp (M',-q') \cong \HH(M)$.

\eqref{carhyp-cor-aa} If $(M',q') \perp (M', -q') \cong \HH(M')$, then $(M',q')$ is an orthogonal summand of the regular quadratic module $\HH(M')$ and therefore regular by \ref{bfLG}\eqref{qfba-sum}. The converse is \eqref{carhyp-cor-a}.%

\eqref{carhyp-cor-b} We have $(M,q) \perp (M,-q) \cong \HH(M) $ by \eqref{carhyp-cor-a}. Assuming  that $(M,q) \perp (M', -q') \cong \HH(M)$,  we therefore get $(M,q) \perp (M,-q) \cong (M,q) \perp (M',-q')$. Hence $(M,-q) \cong (M',-q')$ by Witt cancellation \ref{canqf}. 
The other direction follows from \eqref{carhyp-cor-a}.
\end{proof}
\sm

\textbf{Remark.} Corollary~\ref{carhyp-cor} is not true for nonsingular quadratic forms, even over fields. For example, let $(M,q) = (F, \lan u\ran_q) = (M', q')$ with $F$ a field of characteristic $2$ and $u\in F\ti$. Then $(M,q)$ is nonsingular by \eqref{qfba-one1}, but $0 \ne (1_F, 1_F) \in \rad(q \perp (-q))$, so that $q\perp (-q)$ is singular, hence in particular not hyperbolic. 
\ms

Taking a maximal faithful complemented totally isotropic submodule $U$ in \ref{quadrepII} we obtain the existence part of following Corollary~\ref{wide}, well-known for quadratic spaces over fields \cite[Thm.~8.5]{EKM} and for regular symmetric bilinear modules \cite[I, \S3, Prop.~3]{Knebusch-Queens}.    

\comments{(2025-07-19) In a zoom meeting yesterday, we agreed to assume in \ref{wide} that $R$ is connected, in order to avoid some counter-intuitive examples. }

\begin{cor}[Witt Decomposition] \label{wide} Every quadratic $R$--space $(M,q)$ \new over a {\em connected\/} $R$ \enew has an orthogonal decomposition, referred to as\/ {\em Witt decomposition}, 
\begin{equation}\label{wide1} 
 (M,q) = (M_h, q|_{M_h}) \perp (M_a, q|_{M_a})
\end{equation} 
such that $(M_a, q|_{M_a})$ is anisotropic and that $(M_h, q|_{M_h}\cong \HH(U)$ is hyperbolic with $U=0$, i.e., $(M,q)$ anisotropic, or $U$ is a maximal faithful complemented totally isotropic submodule of $M$.%
\sm 

Suppose in addition that $R$ is unimodular and has Witt cancellation, for example assume that $R$ is LG or that $R=k[X]$ for $k$ a field. Then the quadratic submodules $(M_h, q|_{M_h})$ and $(M_a, q|_{M_a})$ are unique, up to isometry. 
\end{cor}

\begin{proof} For the uniqueness proof let $(M,q) = \HH(U) \perp (M_a, q_a)$ and $(M,q) = \HH(U') \perp (M_a', q'_a)$ be two Witt decompositions. We can assume that $U$ or $U'$ are non-zero, say $U \ne 0$. Then $U$ is faithfully projective, hence contains a unimodular vector $x$, so that $U = Rx \oplus U_1$ for some complementary submodule $U_1$. From \eqref{qfba-hyp3} we get $\HH(U) = \HH(Rx) \perp \HH(U_1)$. Since  $\HH(Rx)$ is a faithful complemented totally isotropic submodule of $(M,q)$, necessarily $U' \ne 0$, so that the same argument yields $\HH(U') = \HH(Rx') \perp \HH(U'_1)$ for some unimodular $x'\in U'$. Since $R$ has Witt cancellation, we can cancel $\HH(Rx)$ and $\HH(Rx')$ and get Witt decompositions $\HH(U_1) \perp (M_a, q_a) = \HH(U'_1) \perp (M'_a, q'_a)$. We can then conclude by induction. 
\end{proof}
\enew
\sm 

If $R$ is a field, existence and uniqueness of quadratic submodules $(M_h, q|_{M_h}$ and $(M_a, q|_{M_a})$ of \eqref{wide1} is proven in \cite[Thm.~8.5]{EKM} with essentially the same proof. For a general $R$, we only have a weak ``uniqueness'', namely when the maximal complemented totally isotropic submodules of $(M,q)$ are isomorphic, see \ref{hps} and \ref{Lag}.
\sm 

\textbf{Examples:} \begin{inparaenum}[(a)]\item A nonsingular $1$--dimensional form $\lan a \ran_q$ is anisotropic. A $2$--dimensional quadratic form $[a,b]$ is either ansiotropic or hyperbolic. Hence, in these two cases the Witt ``decomposition'' is unique for any $R$. 

\item The Witt decomposition of the quadratic form of the Example~\ref{canqf}\eqref{canqf-c} is given in \eqref{canqf-c1}. 

\item \label{wide-c} One has of course a decomposition \eqref{wide1} of a quadratic space  over an arbitrary $R$. We have chosen to stick to the connected case since in general such a decomposition has rather counter-intuitive properties. A case in point is the direct product $R=R_1\times R_2$ of fields $R_1$ and $R_2$ with $(M_1, q_1)$ hyperbolic over $R_1$ and $(M_2, q_2)$ anisotropic over $R_2$. 
\end{inparaenum}

\begin{lem}[Nonsingular quadratic forms over LG rings]\label{nqf-LG} Let $(M,q)$ be a quadratic space over an LG ring $R$ which has constant rank $n\in \NN_+$.   

Then there exist free rank-$2$-submodules $M_1$, \ldots, $M_{[n/2]}$ of $M$  with $q_i = q|_{M_i}$ regular, and, if $n$ is odd, a free rank-$1$-submodule $M_0$ with $q_0 =q |_{M_0}$ nonsingular such that 
      \[ (M,q) = \begin{cases}
        (M_1, q_1) \perp \cdots \perp (M_m, q_m), & \text{$n=2m$ even} \\ (M_1, q_1) \perp \cdots \perp (M_m, q_m) \perp
                 (M_0, q_0), & \text{$n=2m+1$ odd} 
      \end{cases}  \;   . 
      \]
By {\rm \ref{qf}\eqref{qfba-one}}, $(M_0, q_0) \cong \lan u \ran_q$ for some $u\in R\ti$.  Every rank-2-submodule $(M_i, q_i)$, $1\le i \le m$, in the decomposition above is isometric to a quadratic form $[a_i, b_i]$, defined in \eqref{qfba-two0}, for suitable $a_i, b_i \in R$ with $1-4a_i b_i \in R\ti$. \end{lem}

\begin{proof} 
That $(M,q)$ decomposes as claimed, is proven in \cite[Thm.~7.3(a)]{EG} for regular forms. The proof for nonsingular forms is the same, in view of the field case treated in \cite[Prop.~7.29, Cor.~7.32]{EKM}. 

It remains to identify the regular rank-2-submodules $(M_i, q_i)$, $1 \le i \le m$. For simpler notation, let $(M_i, q_i) = (M,q)$, and let $f$ and $g$ be the polynomials on $\uW(M\oplus M)$ whose $R$--points are $f(u_1, u_2) = \det \big( b_q(u_i, u_j)\big)$ and $g(u_1 , u_2) = b_q(u_1, u_2)$. Our claim is that there exists $(u_1, u_2) \in M\oplus M$ such that $(u_1, u_2)$ is a basis of $M$, i.e., $f(u_1, u_2) \in R\ti$, with $g(u_1, u_2) \in R\ti$, since then $(u_1, b_q(u_1, u_2)\me u_2)$ is a basis of $M$ as required in the lemma. Equivalently, we claim that the quasi-compact open subscheme $U=D(f) \cap D(g) = D(fg) 
\subset \uW(M\oplus M)$ has an $R$--point. As $R$ is LG, it suffices to prove this over fields in $\Ralg$. But this is \cite[Rem.~7.4]{EKM}. 
\end{proof} 

\lv{
Let $R$ be an LG ring. We can assume that $n \ge 2$. Let $v_1, \ldots, v_n$ be a basis of $M$. Let $X_1, \ldots, X_n, Y_1, \ldots, Y_n$ be indeterminates, define $w_1 = X_1 v_1 + \cdots + X_n v_n$ and $w_2 = Y_1 v_1 + \cdots Y_n v_n \in M \ot_R R[X_1, \ldots, X_n, Y_1, \ldots, Y_n]$ and  $f(X_1, \ldots, X_n, Y_1, \ldots, Y_n) = \det \big( b_q(w_i, w_j)\big)_{1\le i,j\le 2}$. By the field case, the polynomial $f$ represents a unit over every field in $\Ralg$. Hence, it represents a unit over $R$. This means that there exist $u,v\in M$ such that $M_1 = Ru + Rv$ is free of rank $2$ and $q_1 = q|_{M_1}$ is regular. By \eqref{quadrepIIbn} we can therefore split off $(M_1, q_1)$ and continue with its orthogonal complement, which is nonsingular by \eqref{qfba-sing-vii}. 

Previous: For a semilocal $R$, part \eqref{nsfLG-b} is proven in \cite[IV, (2.2.2) and (3.1.7)]{K} in general and for regular forms in \cite[I, (3.4)]{Ba}.
}

\comments{(2021-03-02) Taking into account \eqref{qfba-two1} and reduction to constant rank \ref{qf}\eqref{qf-redc}, the following corollary is an immediate consequence of \eqref{nsf1}. It is no longer  needed: \\

\textbf{Corollary} 
Let $(M,q)$ be a quadratic space over a semilocal ring $R$ with $\rank_\p M \ge 2$ for all $\p \in \Spec(R)$. Then there exists $a \in R\ti$, a quadratic \'etale $R$--algebra $E$ with norm $n_E$ and a nonsingular quadratic form $q'$ such that $q\cong a n_E \perp q'$. \\

(2026-02-27) However, we use in the proof of \ref{knevid}\ref{knevid-b} that a quadratic space $(M,q)$ over an LG ring $R$ with $\rank M \ge 2$ contains a submodule $N$ which is free of rank $2$ and for which $q|_N$ is regular. This follows from the rank decomposition of $(<q)$ and then using \ref{nqf-LG}. 
}

\comments{(2021-11-17) Bass also  defines a Witt ring as the cokernel of $K_0$ of the category of regular quadratic modules under the canonical map from $K_0$ of finite projective modules}

\subsection{Witt-Grothendieck group $\hWq(R)$ and the Witt group $\Wq(R)$}\label{wgg}
We denote by $\hWq(R)$ the Witt-Grothendieck group of $R$ \cite[I, \S2]{Ba}. The elements of $\hWq(R)$ are formal differences $[q_1] - [q_2]$ of isometry classes $[q_i]$, $i=1,2$, of {\em regular\/} quadratic modules $q_i$. Two elements $[q_1] - [q_2]$ and $[q_1'] - [q_2']$ are equal in $\hWq(R)$ if and only if there exists a regular quadratic module $p$ such that $q_1 \perp q_2' \perp p \cong q_1' \perp q_2 \perp p$. The group operation of $\hWq(R)$ is induced by the orthogonal sum $[q_1] + [q_2] = [q_1 \perp q_2]$ (this is why we only use regular quadratic modules, see \ref{qf}\eqref{qf-perp}). The zero of $\hWq(R)$ is represented by the zero quadratic form.
Given $S\in \Ralg$, base change induces a group homomorphism
\begin{equation}\label{wgg1}
 \wdh r_{S/R} \co  \hWq(R) \longto \hWq(S), \qquad [q] \mapsto [q_S],
 \end{equation}
called the {\em restriction homomorphism}. \sm

Following \cite[I, \S4]{Ba}, the Witt group $\Wq(R)$ is the quotient group of $\hWq(R)$ by the subgroup of elements of the form $[\HH(P_1)] - [ \HH(P_2)]$, $P_1$ and $P_2$ finite projective $R$--modules. It follows that the elements of $\Wq(R)$ are represented by isometry classes $[q]$ of regular quadratic forms $q$. Equality $[q_1] = [q_2]$ in $\Wq(R)$ means that there exist finite projective $R$--modules $P_1$ and $P_2$ such that
\begin{equation}\label{wgg11}
q_1 \perp \HH(P_1) \cong q_2 \perp \HH(P_2).
\end{equation}
The equivalence relation \eqref{wgg11} is called {\em Witt equivalence}. As before we have group homomorphisms
\begin{equation}\label{wgg2}
 r_{S/R} \co \Wq(R) \longto \Wq(S), \qquad [q] \mapsto [q_S]
 \end{equation}
for any $S\in \Ralg$. \ms

The following Proposition~\ref{equi} is inspired by \cite[Prop.~1.2]{Collio}, which proves it for semilocal rings $R$ with $2\in R\ti$. It puts the Corollaries 2.7, 2.8, 2.9 and 2.12 of \cite{GN-Sp} into a general framework. We will say that {\em a quadratic module $(M,q)$ contains a quadratic module $(M_1, q_1)$\/} if there exists a complemented submodule $N\subset M$ such that
$(M_1,q_1) \cong (N,q|_N)$. We recall that if $(M_1, q_1)$ is regular, e.g., a hyperbolic space, then $(M, q) = (M_1, q_1) \perp (M_1, q_1)^\perp$ by \ref{orthLG}\eqref{orthLG-a} (and \ref{qf}\eqref{qf-regLG}). 

\comments{(2023)
Note that over fields, by Springer's Theorem, condition \ref{equi-A} holds for ANY quadratic form, since a quadratic form is either anisotropic or isotropic. Then we also have condition \ref{equi-Bhyp}, \ref{equi-Breg}, \ldots \ref{equi-C}. This could be addressed in \ref{equisa}.}

\comments{(2024)
MISSING are examples that the implications cannot be reversed.
Example: \eqref{wgg2} is in general not injective (\cite[VIII, (2.1.1)]{K} for $R$ a noetherian domain of dimension one, and again in \cite[VIII, (2.2)]{K} for $R$ a regular domain).
\sm

In \cite[I, (10.3.3)]{K} it is shown that \ref{equi-E} $\rightarrow$ \ref{equi-Dreg} for $\veps$--hermitian forms over a division ring over a field of characteristic $\ne 2$. The proof works whenever $R$ and $S$ are fields.}

\comments{(2025-07-07) I completely agree with your comment on \ref{equi}: 
\pcomments{(2025-06-26) C'est assez indigeste.}

\comments{(2025-10-09) I use \ref{equi} in the section on Springer's Theorem \ref{sec:springer}. So, I am now for keeping \ref{equi} and the remarks \ref{equisa}.} }

\begin{prop}\label{equi} We fix a pair $(R,S)$ consisting of a ring $R$ and 
$S\in \Ralg$,  and consider the following conditions  
on quadratic $R$--spaces:

\begin{enumerate}[label={\rm (\Alph*)$_\hyp$}]
 \item \label{equi-Bhyp} Given a quadratic $R$--space $(M,q)$ and a finite projective $R$--module $N$ such that $(M,q)_S$ contains $\HH(N)_S$, then $(M,q)$ contains $\HH(N)$. \sm
\end{enumerate}

\begin{enumerate}[label={\rm (\Alph*)$_\reg$}]
 \item\label{equi-Breg} Let $(M,q)$ be a quadratic $R$--space and let $(M_1,q_1)$ be  a\/ {\em regular} quadratic module such that $(M,q)_S$ contains $(M_1, q_1)_S$. Then $(M,q)$ contains $(M_1, q_1)$. \sm
\end{enumerate}

\begin{enumerate}[label={\rm (\Alph*)$_\hyp$}]\setcounter{enumi}{1} 
\item\label{equi-Dhyp} If $(M,q)$ is a quadratic $R$--space and $N$ is a finite projective $R$--module such that $(M,q)_S \cong \HH(N)_S$, then $(M,q)\cong \HH(N)$.  \sm 
\end{enumerate}

\begin{enumerate}[label={\rm (\Alph*)$_\reg$}]\setcounter{enumi}{1}
\item \label{equi-Dreg} If $(M,q)$ is an $R$--space and $(M',q')$ a regular quadratic $R$--module such that $(M,q)_S \cong (M', q')_S$, then $(M,q) \cong (M', q')$. \sm
\end{enumerate}

\begin{enumerate}[label={\rm (\Alph*)$_{{\rm rk} 1}$}]\setcounter{enumi}{1}

 \item \label{equi-D1} If $(M,q)$ and $(M',q')$ are quadratic $R$--spaces of rank $1$ such that $(M,q)_S \cong (M', q')_S$, then $(M,q) \cong (M', q')$.
     \sm
\end{enumerate}

\begin{enumerate}[label={\rm (\Alph*)}]\setcounter{enumi}{1}
\item \label{equi-D} If $(M,q)$ and $(M',q')$ are quadratic $R$--spaces such that $(M,q)_S \cong (M', q')_S$, then $(M,q) \cong (M', q')$. \sm

\item\label{equi-E} The homomorphisms $\hWq(R) \to \hWq(S)$ and
$\Wq(R) \to \Wq(S)$ of  \eqref{wgg1} and \eqref{wgg2} are injective. \sm
\end{enumerate}

\begin{enumerate}[label={\rm (\Alph*)}]\setcounter{enumi}{3}
 \item \label{equi-A} If $(M,q)$ is a quadratic $R$--space such that $(M,q)_S$ contains an isotropic vector, then so does $(M,q)$. \sm
\end{enumerate}

\begin{enumerate}[label={\rm (\Alph*)}]\setcounter{enumi}{4}
  \item \label{equi-C} If the\/ {\em regular\/} quadratic module $(N,q)$ and $u\in R\ti$ are such that $u \ot 1_S$ is represented by $(N,q)_S$, then $u$ is represented by $(N,q)$. \sm
\end{enumerate}
These conditions are related in the diagram below with the following explanation of arrows: 

\quad $\Rightarrow$ holds whenever $R$ has Witt cancellation, e.g., if $R$ is an LG ring, 

\quad $\rightarrow$ holds if $R$ and $S$ are LG rings, 

\quad $\rightsquigarrow$ holds if $R$ is semilocal, and 

\quad $\dashrightarrow$ holds if $2\in R\ti$:
\[ \xymatrix@C=20pt{
   \text{\ref{equi-C}}\ar@<0.5ex>@{-->}[r]
   & \text{\ref{equi-A}} \ar@<0.5ex>[r] \ar@{~>}@<0.5ex>[l]
   & \text{\ref{equi-Bhyp}} \ar@2[d]\ar@2{<->}[r]  \ar@2@<0.5ex>[l]
  &\text{\ref{equi-Breg}} \ar@2[d] 
\\
& & \text{\ref{equi-Dhyp}} \ar@2{<->}[r]  & \text{\ref{equi-Dreg}} \ar[r] \ar@{-->}@<0.5ex>[d]
   &\text{\ref{equi-E}}
\\
 && \text{\ref{equi-D1}}+\text{\ref{equi-Breg}}
   \ar[r] & \text{\ref{equi-D}} \ar@2@<0.5ex>[u]
}\]
\end{prop}

\begin{proof}

\ref{equi-Bhyp} $\implies$ \ref{equi-A}: If $(M,q)_S$ contains an isotropic vector, it contains a hyperbolic plane $\HH$ by \ref{isotrop}\eqref{isotrop-bc}. Hence so does $(M,q)$ by \ref{equi-Bhyp}, in particular \ref{equi-A} holds.%
\sm

\ref{equi-Bhyp} $\implies$ \ref{equi-Breg}: Let ($M,q)$ and $(M_1, q_1)$ be as in the assumption of \ref{equi-Breg}. Since $(M_1, q_1)_S$ is regular, there exists a nonsingular $S$--space $(\wtl M_2, \wtl q_2)$ such that $q_S \cong q_{1,S} \perp \wtl q_2$. Also recall that $q_{1,S} \perp (- q_{1,S}) \cong \HH(M_1)_S$ by \ref{carhyp-cor}\eqref{carhyp-cor-a}. Hence
 \[ \big (q \perp (-q_1)\big)_S \cong q_S \perp (-q_1)_S \cong q_{1,S} \perp (-q_1)_S \perp \wtl q \cong \HH(M_1)_S \perp \wtl q.
 \]
Thus we can apply \ref{equi-Bhyp} to the nonsingular form $q \perp (-q_1)$ and $N = \HH(M_1)$, and can conclude that $q \perp (-q_1)$ contains $\HH(M_1) \cong q_1 \perp (-q_1)$, i.e., there exists a quadratic $R$--space $q_2$ such that $q \perp (-q_1) \cong q_1 \perp (-q_1) \perp q_2$. Cancelling the regular form $-q_1$ using \ref{canqf}, we get $q \cong q_1 \perp q_2$, finishing the proof of \ref{equi-Bhyp} $\implies$ \ref{equi-Breg}, and establishing \ref{equi-Bhyp} $\iff$ \ref{equi-Breg} since
\ref{equi-Bhyp} is a special case of \ref{equi-Breg}, as a hyperbolic space is regular. \sm

The implications \ref{equi-Bhyp} $\implies$ \ref{equi-Dhyp} and \ref{equi-Breg} $\implies$ \ref{equi-Dreg} are clear, and \ref{equi-Dhyp} $\iff$ \ref{equi-Dreg} is a special case of \ref{equi-Bhyp} $\iff$ \ref{equi-Breg}. Obviously \ref{equi-Dreg} $\implies $ \ref{equi-D}.\sm

\ref{equi-A} $\rightarrow$ \ref{equi-Bhyp}: We assume that $S$ is LG, so that we can use Witt cancellation \ref{canqf} for forms over $S$. 

Let $(M,q)$ and $N$ be as in the assumption of \ref{equi-Bhyp}. Applying the rank decomposition to $N$ and correspondingly to $(M,q)$, cf.\ \ref{qf}\eqref{qf-redc}, we see that it suffices to establish \ref{equi-Bhyp} for $N$ of constant rank, hence for free $R$--modules of rank $r>0$ by \eqref{unimod-b1} using that $R$ is LG. For doing so, we use induction on $r$. 
We observe that $\HH(N)_S$ contains an isotropic vector and that therefore so does $(M,q)_S$. By \ref{equi-A}, we can choose an isotropic vector $e\in M$. Thus, by \eqref{isotrop-d0}, $(M,q) = \HH(Re) \perp (M_1, q_1)$ with a nonsingular $(M_1, q_1)$. Also, since $N$ is free, we can choose a decomposition $N = Re_1 \oplus N_1$ with $R_1$ and $N_1$ being free of rank $1$ and $r-1$ respectively. Correspondingly $\HH(N) \cong \HH(Re_1) \oplus \HH(N_1)$. This finishes the proof for $r=1$. If $r>1$, we get
  \begin{align*}
    \HH(Se) \perp (M_1, q_1)_S & \cong (M,q)_S \cong \HH(N)_S \perp \HH(N)_S^\perp \\ & \cong \HH(S e_1) \perp \HH(N_1)_S \perp \HH(N)_S^\perp.
  \end{align*}
By \ref{canqf} we can cancel $\HH(Se) \cong \HH(Se_1)$, 
  thus $(M_1, q_1)_S \cong \HH(N_1)_S \perp \HH(N)_S^\perp$. Now \ref{equi-Bhyp} follows by induction.
\sm

\ref{equi-Dreg} $\rightarrow$ \ref{equi-E}: To prove injectivity of $\hWq(R) \to \hWq(S)$, let $q_i$, $i=1,2$ be two regular quadratic modules such that $[q_{1,S}] - [q_{2,S}] =[0] \in \hWq(S)$. This means that there exists a regular quadratic $S$--module $p$ such that $q_{1,S} \perp p \cong q_{2,S}  \perp p$. Applying again Witt cancellation \ref{canqf}, we obtain $q_{1,S} \cong q_{2,S}$, and then $q_1 \cong q_2$ by \ref{equi-Dreg}. In particular $[q_1] - [q_2] = [0] \in \hWq(R)$.
\sm

For proving injectivity of $\Wq(R) \to \Wq(S)$, we can assume that $q$ is a regular quadratic $R$--form for which there exist finite projective $S$--modules $\wtl P$ and $\wtl Q$ such that $q_S \perp \HH(\wtl P) \cong \HH(\wtl Q)$. Decomposing $S$ into a suitable finite product of rings permits to  assume that $\wtl P$ and $\wtl Q$ have constant rank, hence are free $S$--modules. Since $\rank_S \wtl P \le \rank_S \wtl Q$, there exists free $R$--modules $P$, $P'$ and $Q$ such that $P_S \cong \wtl P$ and $(P \oplus P')_S \cong \wtl Q$. Hence $q_S \perp \HH(P)_S \cong q_S \perp \HH(P_S) \cong q_S \perp \HH(\wtl P) \cong \HH(\wtl Q) \cong \HH(P)_S \perp \HH(P')_S$. Canceling $\HH(P)_S$, using that $S$ is assumed to be LG, we get $q_S \cong \HH(P')_S$ and then $q \cong \HH(P')$.
\sm

\ref{equi-D1} $+$ \ref{equi-Breg} $\rightarrow$ \ref{equi-D}:
Let $(M,q)$ and $(M',q')$ be quadratic $R$--spaces such that $(M,q)_S \cong (M',q')_S$. We can assume that both $M$ and $M'$ have constant rank, necessarily equal and denoted by $r$. If $r$ is even, then $q$ and $q'$ are regular by \ref{qf}\eqref{quadfoe}, and the claim $q \cong q'$ is a special case of \ref{equi-Breg}. 
We can therefore assume that $r$ is odd. By \ref{nqf-LG} there exist $u$, $u' \in R\ti$ and regular quadratic forms $(M_0, q_0)$ and $(M'_0, q'_0)$ such that
\begin{equation} \label{prop_odd1}
(M,q) \cong  (R, \lan u \ran) \perp (M_0, q_0) \text{ and }
   (M',q') \cong  (R, \lan u' \ran) \perp (M'_0, q'_0).
\end{equation}
If $r=1$, i.e., $M_0=0=M_0'$, the assumptions $\lan u \ran_S \cong \lan u'\ran_S$ and \ref{equi-D1} yield
$(R, \lan u \ran) \cong (R, \lan u'\ran)$. Let now $r\ge 3$.
By assumption, the quadratic form $q_{0,S}$ is contained in $q'_S$. Hence, by \ref{equi-Breg}, $q_0$ is contained in $q'$. We can therefore replace $(M'_0, q_0')$ in the decomposition \eqref{prop_odd1} by $(M_0, q_0)$. Thus, $\lan u\ran_S \perp q_{0,S} \cong \lan u'\ran_S \perp q_{0,S}$.  Since $S$ is assumed to be LG, we can now apply Witt cancellation~\ref{canqf} and get $\lan u \ran _S \cong \lan u'\ran_S$. The case $r=1$ then shows $\lan u \ran \cong \lan u'\ran$. But then $q \cong q'$ by \eqref{prop_odd1}. This finishes the proof of \ref{equi-Breg} $+$ \ref{equi-D1} $\rightarrow$ \ref{equi-D}.
\sm

\ref{equi-A} $\rightsquigarrow$ \ref{equi-C}: Let $(N,q)$ be a regular quadratic module and let $u\in R\ti$ such that $u_S = u\ot 1_S$ is represented by $q_S$, say $u_S = q_S(x)$ for some $x\in N_S$. Then $(N', q') = (N,q) \perp (R, \lan -u \ran)$ is a quadratic space by \eqref{qfba-one1} and \eqref{qf-perp1} such that $(x,1_S)$ is an isotropic vector of $(N', q')_S$. By \ref{equi-A}, the quadratic space $(N', q')$ contains an  isotropic vector, i.e., there exists $y\in N$ and $r\in R$ such that $y'= (y,r) \in N \oplus R$ is unimodular and satisfies $q(y)= u r^2$. If $r\in  R\ti$, then $y_0 = r\me y\in N$ satisfies $q(y_0) = u$. It therefore suffices to find $y_1 \in N$ such that $q(y_1) = r^2 u$ with $r\in R\ti$. That this is possible, is shown in \cite[Lem.~2.13]{GN-Sp}. We repeat the proof for the convenience of the reader.

The goal is to construct a $v' \in N'$ with $q'(v') \in R\ti$ such that the vector $\si_{v'}(y') $ has an invertible $R$--component. Here $\si_{v'}$ is the reflection in $v'$, reviewed in \ref{refle}. Writing such a $v'$ as $v'=(v,s)$ with $v\in N$ and $s\in R$, we have \begin{align*}
 \si_{v'}(y') &= y' - b_q(y',v')q'(v')\me v' \\
     &= y' - ( b_q(y,v) - 2urs) (q(v) - u s^2)\me (v,s) \\
    &= (\star, r- ( b_q(y,v) - 2urs) (q(v) - u s^2)\me s) =: (\cdot, t)
\end{align*}
where the last entry $t$ is the $R$--component of $\si_{v'}(y')$.
We know $t\in R\ti \iff t\ot 1_{R/\m} \ne 0$ for every maximal ideal $\m \ideal R$.  Since $N \to N/\Jac(R)$ is surjective, we can therefore choose $v$ and $s$ ``componentwise'', i.e., for each maximal ideal $\m \ideal R$. Thus, let us fix a maximal ideal $\m$ of $R$ and put $\ka = R/\m$. For $p\in N$ we write  $p_\ka = p \ot 1_\ka$ and analogously for elements of $R$.

Case $r_{\ka} \ne 0$: we put $s_\ka= 0$ and $v_\ka = y_\ka$. Then $q'_\ka(v'_\ka) = q_\ka(y_\ka) - u_\ka s_\ka^2 = q(y_\ka) = u_\ka r_\ka^2 \ne 0$ and $t_\ka = r_\ka \ne 0$.

Case $r_\ka = 0$: we put $s_\ka = 1_\ka$. Since then $t_\ka = b_{q_\ka}(y_\ka,v_\ka)(q_\ka(v_\ka) - u_\ka)\me$ the conditions on $v_\ka$ are
\begin{equation}\label{equi-b1}
  q(v_\ka) \ne a_\ka, \qquad b_q(y_\ka, v_\ka) \ne 0.
\end{equation}
Observe that $y_k$ is an isotropic vector of $(M_\ka, q_\ka)$, which therefore embeds into a hyperbolic plane $(y_\ka= e, f)$. Indeed, $(M,q)$ is a quadratic space, hence so is $(M,q)_\ka$, so that the claim follows from \ref{isotrop}\eqref{isotrop-bc}.  For $v_\ka = f$ we get $q(f) = 0 \ne u_\ka$ (since $u\in R\ti$) and $b(y_\ka , f) = 1 \ne 0$. Thus \eqref{equi-b1} holds, and this finishes the proof of \ref{equi-A} $\rightsquigarrow$ \ref{equi-C}. \sm

\ref{equi-C} $\dashrightarrow$ \ref{equi-A}: We assume $2\in R\ti$, and
let $(M, q)$ be a quadratic space for which $(M,q)_S$ contains an  isotropic vector, hence, by \ref{isotrop}\eqref{isotrop-bc}, contains a hyperbolic plane $P= Se \oplus Sf$ with $q(e) = 0 = q(f)$ and $b_q(e,f) = 1$. Then $x'=\frac{1}{2}(e+f)$ and $y' = \frac{1}{2}(e-f)$ are a basis of $P$ satisfying $q(x) = 1$, $b_q(x,y) = 0$ and $q(y) = -1$. Applying \ref{equi-C}, we know that there exists $x\in M$ with $q(x) = 1$. By \ref{mx_lem}, $(M_x = (Rx)^\perp, q_x = q|_{M_x})$ is a quadratic space. Since $y'\in (M_x, q_x)_S$, another application of \ref{equi-C} shows that there exists 
$y\in (Rx)^\perp$ such that $q(y) = -1$. But then $x+y\in M$ is an isotropic vector.
\sm

Under the assumption $2\in R\ti$ we know that regular and nonsingular quadratic forms are the same. Therefore 
\ref{equi-Dreg} = \ref{equi-D}. \end{proof}


\subsection{Instances when some of the conditions of Proposition~\ref{equi} are satisfied}\label{equisa} Let $R$ be arbitrary and let $S\in \Ralg$; further conditions on $R$ and/or on $S$ will be specified below. The conditions \ref{equi-Bhyp},...,\ref{equi-C} of \ref{equi} make sense for any pair $(R, S)$. We give below a likely incomplete account of what is known about their validity; see \cite[2.3]{Collio} for a discussion covering results up to 1979.
\sm

\begin{inparaenum}[(a)]\item\label{equisa-x} By \ref{isotrop}\eqref{isotrop-bc}, condition \ref{equi-A} of \ref{equi} is equivalent to the condition
\end{inparaenum}
\begin{enumerate}[label={\rm (\Alph*$'$)}]\setcounter{enumi}{3}
 \item \label{equi-A'} If $(M,q)$ is a quadratic $R$--space such that $(M,q)_S$ contains a hyperbolic plane, then already $(M,q)$ contains a hyperbolic plane.
\end{enumerate}
This is the formulation of condition \ref{equi-A} used in several papers, e.g. in \cite{PS}.%
\sm 

\begin{inparaenum}[(a)]\setcounter{enumi}{1}
\item (Condition~\ref{equi-D1})  \label{equisa-a} 
  Suppose $\Pic(R) = 0$, e.g., assume $R$ is unimodular or even LG. Then condition \ref{equi-D1} of \ref{equi} is equivalent to the following condition: Given $u$, $u'\in R\ti$ and $v\in S\ti$ such that $u'\ot 1_S = (u\ot 1_S)v^2$, then there exists $y\in R\ti$ satisfying $u' = uy^2$.
\lv{
Indeed, since $M$ and $M'$ are free of rank $1$, we can assume that $(M,q) = \lan u \ran$ and $(M', q') = \lan u' \ran$ for some $u$, $u'\in R\ti$. An isomorphism over $S$ is given by $v\in S\ti$ such that $u'\ot 1_S = v^2 u\ot 1_S$. Similarly, and isomorphism over $R$ is given by $y\in R\ti$ satisfying $u' = y^2 u$.
}

For example, {\em if $S$ is finite projective of constant {\em odd\/} rank, then \ref{equi-D1} holds}. Indeed, applying the norm $\rmN_{S/R}$ of the extension $S/R$ to the equation
$u'\ot 1_S = (u\ot 1_S)v^2$, we get $u'{}^d = u^d \rmN_{S/R}(v)^2$ for $d=\rank S$. Since $d$  is odd, $u'= u y^2$ for some $y\in R\ti$ follows.
\sm

\item\label{equisa-d} {\em  If $R$ is an integrally closed domain, i.e., a normal ring,  with $\Pic(R) = 0$ and field of factions $S$, then \ref{equi-D1} holds.} Indeed, as can be seen from  the argument in \eqref{equisa-a}, proving \ref{equi-D1} boils down to showing that a quadratic equation that has a solution in S already has a solution in $R$. \sm

\new
\item\label{equisa-cn} (The LG case) {\em If $R$ is an LG ring and $S$ is a finite locally free $R$--algebra of constant odd degree for which 
\begin{equation}\label{equisa-cn1} \begin{split}      
    &\text{ there exists a finite locally free $T\in \Ralg$ of constant odd} \\ &\text{ degree such that $S\ot_R T$ is one-generated as $R$--algebra,}
\end{split}\end{equation}    
    then  \ref{equi-Breg}, \ref{equi-D}, \ref{equi-E} and \ref{equi-A} hold.}

    Indeed, since $R$ is unimodular and therefore has $\Pic(R) = 0$, the hypothesis \ref{equi-D1} holds by \eqref{equisa-a}. We will prove the validity of hypothesis \ref{equi-A} in Corollary~\ref{spo-LGG}, which follows immediately from 
    Theorem~\ref{spo}, commonly referred to as Springer's Odd Degree Theorem. Since $S$ is also LG by \ref{revLG}\eqref{revLG-int}, the two conditions \ref{equi-D1} and \ref{equi-A} imply all the others. 
    
   We will also show in Corollary~\ref{spo-LGG} that the assumption \eqref{equisa-cn1} is satisfied if $R$ satisfies the primitive criterion and $\uPrim(S) \ne \emptyset$.  \sm 
\enew
      
\item \label{equisa-b} (The semilocal case) If $R$ is a semilocal ring and $S$ is finite \'etale of constant odd degree, then \ref{equi-A}, and hence  \ref{equi-Breg}, \ref{equi-D}, \ref{equi-E}, and \ref{equi-C} hold.

     Indeed, \ref{equi-A} follows from \ref{spo-LGG} and \ref{equi-D1} from \eqref{equisa-a}. \sm 
     
\item \label{equisa-c} (The field case) If $R$ is a field and $S/R$ is a field extension, which is either purely transcendental or finite of odd degree, then \ref{equi-A} holds  \cite[7.15, 18.5]{EKM}. It is also clear from  \eqref{equisa-a} that \ref{equi-D1} is satisfied. Thus, under the conditions of \eqref{equisa-c}, also \ref{equi-Breg}, \ref{equi-D}, \ref{equi-E} and \ref{equi-C} hold.    \sm

\item\label{equisa-oja82} If $R$ be a regular local ring of dimension $2$ with $2\in R\ti$ and let $S$ be its fraction field, then \ref{equi-Dhyp} holds by \cite[Thm.~1]{Oj82}, \ref{equi-Dreg} by \cite[Cor.~2 of Thm.1]{Oj82}, \ref{equi-C}
     by \cite[Thm.~2]{Oj82}, and $\Wq(R) \to \Wq(S)$ is injective, i.e., part of \ref{equi-E}, by \cite[Cor. 1 of Thm.~1]{Oj82}. \sm

\item\label{equisa-z} The condition \ref{equi-D} means of course that the base change map \[ \frb_{\uO, S} \co H^1\fppf(R, \uO(q)) \longto H^1\fppf(S, \uO(q_S))\]
of \eqref{cohCT-0} is injective.  We will revisit this condition in \ref{cohCT}--\ref{prop_odd}. 
\sm

\item \label{esquisa-h} Part of condition \ref{equi-E} of \ref{equi} is the statement that 
\begin{equation}\label{esquisa-h1}
 \text{the canonical map $\Wq(R) \to \Wq(S)$ is injective.}
\end{equation}
Here is an incomplete list of pairs $(R,S)$, for which \eqref{esquisa-h1} is true; in all references $S$ is the field of fractions of $R$:
\end{inparaenum} \begin{enumerate}[label={\rm (\Roman*)}]
  \item $R$ is a regular local ring with $2\in R\ti$ and which is essentially of finite type over a field $K$ \cite{O};
      
  \item $R=A_f$ where $A$ is a regular local ring containing a field of characteristic $\ne 2$, $f$ is a local parameter of $A$, \cite[Thm.~B]{OP};      
  
  \item $R$ a complete regular noetherian local domain with $2\in R\ti$,  \cite[VIII, (2.2.3)]{K}; 
       
  \item  $R$ is a regular local ring with $2\in R\ti$, and $R$ is henselian, e.g., complete, or $S$ is Pythagorean \cite[Cor.~2.3 and Cor.~2.6]{CRW}; 
     
  \item $R$ is a regular domain of dimension $2$ or $3$, \cite[Thm.~17 and Thm.~24]{Ojan82};

   \item $R$ is a regular local ring of Krull dimension $\le 4$ with $2\in R\ti$ \cite[Cor.~10.4]{Balmer-Walter}. 
\end{enumerate}

In Springer's Odd Degree Theorem\ref{spo}\eqref{spo-b} we will \eqref{esquisa-h1} for any LG ring $R$ and $S\in \Ralg$ one-generated and finite locally free of odd degree. 


\begin{lem}  \label{hypdes} Let $R$ be a semilocal domain and let $K$ be its fraction field. Suppose the following condition \ref{hypdes-i} on quadratic $R$-spaces is satisfied:

\begin{enumerate}[label={\rm (\Alph*)$_\hyp$}]\setcounter{enumi}{3}

\item  \label{hypdes-i} If $(N,q)$ is a quadratic $R$--space and if $U$ is a finite projective $R$--module such that $(N,q)_K \cong \HH(U)_K$, then $(N,q) \cong \HH(U)$.
\end{enumerate}
Then \ref{equi-D} holds for\/ {\em regular} quadratic modules, and hence \ref{equi-D} and there also \ref{equi-E} hold in full generality if $2\in R\ti$.

\end{lem}

\begin{proof} We need to show: if $(M,q)$ is a regular quadratic module and $(M',q')$ is a quadratic space such that $(M',q')_K \cong (M,q)_K$, then $(M',q') \cong (M,q)$. Thus let $(M',q')$ be such a quadratic space. By Corollary~\ref{carhyp-cor}\eqref{carhyp-cor-a} we have $(M,q) \perp (M,-q) \cong \HH(M)$ because $(M,q)$ is regular. Hence
\begin{align*}
  \HH(M)_K &\cong \big( (M,q) \perp (M,-q)\big)_K \cong (M,q)_K \perp (M,-q)_K
  \\&  \cong (M,q)_K \perp (M', -q')_K \cong \big( (M,q) \perp (M', - q')\big)_K.
\end{align*}
Applying \ref{hypdes-i}, we get $(M,q) \perp (M',-q') \cong \HH(M)$ and then $(M,q) \cong (M', q')$ by Corollary~\ref{carhyp-cor}\eqref{carhyp-cor-b}. \end{proof}
\ms

\subsection{Remarks.} \label{hypdesrem} The proof of Lemma~\ref{hypdes} goes back to \cite[Prop.~1.2, (F) $\implies$ (D)]{Collio}, stated there for a semilocal normal $R$ with $2\in R\ti$. The corollary has been retaken in \cite[Thm.~3]{Fedorov}, where $R$ is assumed to be a regular local ring with $2\in R\ti$.
\sm

A similar local-global approach as in \ref{hypdes} is used in \cite{PS} for $R$ a noetherian domain of dimension $1$ with $2\in R\ti$ and a quadratic $R$--space $(M,q)$.
\begin{enumerate}[label={\rm (\roman*)}]
\item (\cite[Thm.~2.1]{PS}) If in addition $R$ is semilocal and if $(M,q)\ot_R R_\p$ is isotropic for all $\p \in \Spec(R)$, then $(M,q)$ is isotropic.

\item  (\cite[Thm.~3.1]{PS}) If the singular set $\mathrm{Sing}(R)$ is finite and non-empty and if $(M,q)\ot_R R_\p$ is isotropic for all $\p \in \mathrm{Sing}(R)$, then $(M,q)$ is isotropic.
\end{enumerate}

\pcomments{(2025-06-26) A mon avis, la bonne hypoth\`ese pour 
le lemma en bas est que $(R,I)$ est un couple hens\'elien (= henselian pair) ce qui couvre tous les cas. Pour un sch\'ema lisse $X/R$ on a que $X(R) \to X(R/I)$ est surjective.
Ainsi  $q$ est isotrope (au sens de Knebusch) sur $R$ si et seulement si
$q $  l'est sur $R/I$. Pour les d\'ecomopositions, un argument du m\^eme genre doit marcher. Je peux r\'ediger cela.}


\comments{(2026-03-28) Are henselian pairs (see above) still in the plan? See your comment above from 2025-06-26.}
\pcomments{(2026-04-09) Je crois que 4.23 servait \`a se ramener au cas d'un anneau int\`egre, cela n'intervient plus et donc tu peux l'enlever.}

\comments{(2026-04-22) Here are the two lemmas to keep in memory what they say. 
\\

\textbf{Lemma}[Lifting orthogonal decompositions, {\cite[I, Cor.~(3.4)]{Ba}}] 
{|em Let $I\ideal R$ be an ideal contained in the Jacobson radical of $R$, and put $\ol R = R/I$. Also, let $(M,q)$ be a quadratic $R$--module for which $(M,q)_{\ol R} = (\ol M, \ol q)$ has an orthogonal decomposition $\ol M = N_1  \perp N_2$ with $N_1$ being free and $\ol q|_{N_1}$ being regular. Then $M$ has an orthogonal decomposition $M= M_1 \perp M_2$ such that $M_1$ is free, $q|_{M_1}$ is regular, and the $M_i$, $i=1,2$, map onto $N_i$ under the canonical map $M \to \ol M$.} \\

\textbf{Lemma}[$R$ complete, {\cite[Lem.~1.1]{PS}} if $2\in R\ti$] 
{\em Let $R$ be a ring which is complete with respect to an ideal $I\ideal R$. Then
condition \ref{equi-A} holds for $S=R/I$: if a nonsingular $q$ is isotropic over $S$, then it is already isotropic over $R$. In particular: \sm

\begin{inparaenum}[\rm (a)]
\item\label{completAi} {\rm (\cite[V, Lem.~1.4]{Ba} for $R$ semilocal and $q$ regular)} Let $R$ be LG and complete with respect to $\Jac(R)$ and put $S=R/\Jac(R)$. Then \ref{equi-A} and all other conditions \ref{equi-Bhyp},  $\ldots$, \ref{equi-A} of {\rm \ref{equi}} hold, except possibly \ref{equi-C}. However, if $R$ is semilocal and complete with respect to $\rad(R)$, then \ref{equi-C} holds and in this case \eqref{wgg2} is an isomorphism: 
      \begin{equation} \label{completAi1}
       r_{S/R} \co \Wq(R) \simlgr \Wq(S).     
      \end{equation}  

\item\label{completAii} {\rm ([PS; Lem.~1.2] for $2\in R\ti$)} 
If $q$ is a nonsingular quadratic form over an artinian ring $R$ and $q_{R/\m}$ is isotropic for every maximal ideal $\m\ideal R$, then $q$ is isotropic. Moreover, all conditions \ref{equi-Bhyp},  $\ldots$, \ref{equi-C} of {\rm \ref{equi}} hold.
\end{inparaenum}
}}

\comments{(2025-04-09) There is more stuff on metabolic forms, and tensor products of quadratic forms with bilinear forms here in 'lv', 2 pages. This is likely not needed for the Knebusch project.}

\lv{
\subsection{Metabolic spaces} \label{meta} Let $(U,b)$ be a bilinear module.  The {\em metabolic space\/} associated with $(U,b)$ is the bilinear module $\MM(U,b) = (U \oplus U^*, b_{\MM(U,b)})$ whose bilinear form $b_{\MM(U,b)}$  is defined by
\begin{equation*} \label{meta-1}
 b_{\MM(U,b)} (u + \vphi, v + \psi) = b(u,v) + \vphi(v) + \psi(u).
 \end{equation*}
It is a regular bilinear form. We say a bilinear module $(M, b)$ is {\em metabolic\/} if there exists a bilinear module $(U, b_U)$ such that $(M, b)\cong \MM(U,b_U)$. A regular bilinear module $(M, b)$ is metabolic if and only if $M$ contains a totally isotropic complemented submodule $V$ with $V = V^\perp$, a so-called {\em Lagrangian}, \cite[I, Thm.~(4.6)]{Ba}, see \ref{Lag} for Lagrangians of quadratic modules.

\comments{(2022-04-07) In Knebusch's 1976 Notes, these metabolic spaces are called ``split metaboliic". He defines ``metabolic'' as a regular quadratic $\calO_X$--module containing a Lagrangian and then shows that over an affine scheme the two notions coincide (Cor.~1 in I, \S3)}

\subsection{Tensor products of bilinear and quadratic forms} \label{qfba-tens} 
Let $(M_i, b_i)$, $i=1,2$, be  bilinear modules. There exists a unique symmetric bilinear form $b_1 \ot b_2$ on $M_1 \ot_R M_2$ satisfying $(b_1 \ot b_2)(m_1 \ot m_2, m_1'\ot m_2') = b_1(m_1, m'_1)\, b_2(m_2, m_2')$ for $m_i, m_i' \in M_i$.

Given an $R$--bilinear module $(M,b)$ and an $R$--quadratic form $(N,q)$, there exists a unique $R$--quadratic form $b\ot q \co M \ot_R N\to R$ satisfying
\begin{equation}\label{tenssq0}
\begin{split} (b \ot q)\, (m \ot n) &= b(m,m) \, q(n), \; \text{and}\\
  b_{b\ot q}(m \ot n, m' \ot n') &= b(m, m')\, b_q(n, n')
\end{split}
\end{equation}
for all $m, m'\in M$ and $n, n'\in N$, where $b_q$ is the polar form of $q$. It is called the {\em tensor product of $(M,b)$ and $(N,q)$}, \cite[Thm.~1]{Sah}.
The polar form of $b \ot q$ is the tensor product of the symmetric bilinear forms $b \ot b_q$.

\sm

A moment's thought will convince the reader that {\em the tensor product respects orthogonal sums}, i.e.,
\begin{equation}\label{qfba-tens-p} \begin{split}
  &(b_1 \perp b_2) \ot (q_1 \perp q_2) \cong
  \\& \quad (b_1\ot q_1) \perp (b_1 \ot q_2) \perp (b_2 \ot q_1) \perp (b_2 \ot q_2)
\end{split} \end{equation}
holds for bilinear modules $(M_i, b_i)$ and quadratic modules $(N_i, q_i)$, $i=1,2$. The analogous formula for the tensor product between orthogonal bilinear modules is of course also true. Moreover, {\em the tensor product is associative and commutative}: with obvious notation we have
\begin{equation} \label{qfba-tenso-2} \begin{split}
  (b_1 \ot b_2) \ot b_3 &\cong b_1\ot (b_2 \ot b_3), \quad
     (b_1 \ot b_2) \ot q \cong (b_1 \ot( b_2 \ot q),
 \\  b\ot q &\cong q \ot b
\end{split} \end{equation}
{\em The tensor product also respects isometries:} if $\vphi \co b \simlgr b'$ and $\psi \co q \to q'$ are isometries, then $\vphi \ot \psi \co b \ot q \simlgr b' \ot q'$ is an isometry.

If $(M,b)$ and $(N,q)$ are bilinear and quadratic modules respectively
   and $ub$ and $uq$ are the obvious scalar multiplications by $u$, then, using the definition \ref{qf}\eqref{qfba-one},
\begin{equation} \label{qfba-tens-2}
 \lan u\ran_b \ot b \cong  ub, \quad \lan u \ran_b \ot q \cong uq
 \end{equation}
under the standard isomorphism $R \ot_R  M \simlgr M$. However,
$ \lan 1\ran_q \ot b = q_b$, where $q_b \co M \to R$, $m \mapsto b(m,m)$, is the quadratic form associated with $b$. Its polar is $b_{q_b} = 2b$.

The tensor product is compatible with base change: for $S\in \Ralg$ we have
\begin{equation}\label{qfba-tens-2}
(M,b)_S \ot_S (N,q)_S \simlgr \big((M,b) \ot_R (N,q)\big)_S
\end{equation}
with respect to $m\ot s_1 \ot n \ot s_2 \mapsto m\ot n \ot s_1 s_2$. We claim that
\begin{equation}\label{tenssq1} \begin{split}
  &\hbox{\em $(M,b)$ and $(N,q)$ are regular}
  \\ & \qquad \iff    \hbox{\em $(M\ot_R N, b\ot_R q)$ is regular,}
\\
 & \hbox{\em $(M_1,b_1)$ and $(M_2,b_2)$ are regular}
\\ &\qquad \iff
      \hbox{\em $(M_1\ot_R M_2, b_1\ot_R b_2)$ is regular.}
\end{split} \end{equation}
We prove the first equivalence. To this end, let $\mu \co M^* \ot_R N^* \to (M\ot_R N)^*$ be the canonical isomorphism
(\cite[II, \S4.4]{BA}). It fits into the commutative diagram
\[\xymatrix{ M^* \ot_R N^* \ar[rr]^\mu_\simeq && (M\ot_R N)^* \\
       & M\ot_R N \ar[ul]^{\wdh b \ot \wdh{b_q}} \ar[ur]_{\wdh{b\ot q}}
}\] where $\wdh{b\ot q}$ is the adjoint of the tensor product form $b\ot q$.
If both $(M,b)$ and $(N,q)$ are regular, the commutative diagram above
shows that $\wdh{b\ot q}$ is an isomorphism. For the proof of the converse we use
\ref{qfLG-rechiv}  of \ref{bfLG}\eqref{bfLG-rech}: regularity can be checked after localization. Thus we can assume that both $M$ and $N$ are free $R$--modules. But then invertibility of $\wdh{b\ot q}$, i.e., of $\wdh b \ot \wdh{b_q}$ by the diagram above, implies invertibility of $\wdh b$ and $\wdh{b_q}$ by  \cite[III, \S8.6 (33)]{BA}.

\begin{cor}[{\cite[I, (4.7.iii)]{Ba}}]\label{lem_metabolic} 
  For each regular quadratic module $(M,q)$ and bilinear module $(U,b)$ the $R$--quadratic space $\MM(U,b) \ot_R (M,q)$ is hyperbolic, 
 \begin{equation} \label{lem_metabolic1} \begin{split}
   \MM(U,b) \ot_R (M,q) &\cong \HH(U \ot_R M), \quad\text{ in particular,} \\
    \HH(U) \ot_R (M,q) &\cong \HH(U\ot_R M).
  \end{split}\end{equation}

\end{cor}

\begin{proof} The quadratic module \[ \MM(U,b) \ot (M,q) = \big( (U\ot_R M) \oplus (U^* \ot_R M), Q\big) = (V,Q)\]  is regular by \eqref{tenssq1}.  According to \eqref{tenssq0} we have
\begin{align*}
  Q(\vphi \ot m) &= b_{\MM(U)}(\vphi, \vphi)\, q(m) =0 \quad \text{and}\\
b_Q( \vphi \ot m, \vphi' \ot m') &= b_{\MM(U)}(\vphi, \vphi')\, b_q(m,m') = 0
\end{align*}
for $\vphi \in U^*$ and $m\in M$. It follows that $N=U^* \ot M$ satisfies the conditions of \ref{carhyp}\ref{carhyp-ii}. By loc.\ cit., we then get $(V,Q) \cong \HH(U^* \ot M) \cong \HH(U \ot M)$ because $U^* \ot M \cong U \ot M$ as $R$--modules.
\end{proof}
}

\newpage

\section{Orthogonal transformations}\label{sec:orthgroup-LG}


\subsection{Reflections}\label{refle} Let $(M, q)$ be a quadratic
module over $R$ and let $b_q$ be the polar form of $q$. We recall that $\orth(M,q) = \orth(q)$ denotes the orthogonal group of $q$, \ref{qf}\eqref{qfLG-a}, that $\wS_q=\{ x\in M : q(x) \in R\ti \}$ is the extended sphere of $(M,q)$, \eqref{quadco-aa2}. For $x\in \wS_q$ the {\em reflection in $x$\/} is the orthogonal transformation $\rho_x\in \GL(M)$, given by
\begin{equation*}  \label{def-orth-1} 
    \rho_x (m) = m -  
        b_q(x,m)q(x)\me \, x  \quad (m\in M).
\end{equation*}
Some of the properties satisfied by reflections are listed below: 
\sm

\begin{enumerate}[label={\rm (\roman*)}]

\item \label{refle-ii} $\rho_x \in \orth(q)$, $\rho_x^2 = \Id_M$,  \sm
 \lv{
 since
\[
q\big(\rho_u(m)\big) = q(m) - b_q\big( m, \frac{b_q(m,u)}{q(u)}u\big) -
    \frac{b_q(m,a)^2}{q(u)^2} q(u) =q(m).
\]}

\item\label{reflei} $\rho_x(x) = -x$, and \\
$\rho_x m = m \iff m\in M_x  = \{ m\in M: b_q(m,x) = 0\}$, hence \sm

\item\label{refle-iia} $\rho_x = \Id_M \iff x\in \rad(b_q)$, \sm

\item\label{refle-det} $\det(\rho_x) = -1$ (\cite[2.6.1.18]{CF}),
 \lv{
  If $2\in R\ti$, then $M = Rx \oplus M_x$ and so $\det(\rho_x) = -1$. If $2R=0$, then $x\in M_x$. In the field case, either $M=M_x$ (and so $\det(\rho_x) = 1 = -1$), or $M_x$ has codimension $1$. Since $M_x$ is complemented, we can write $M_x = R y \oplus M_x$ and then $\rho_x(y) = y - rx$ for suitable $r\in R$. Again $\det(\rho_x) = 1$ follows. }
\item \label{refle-iii} $g \,\rho_x \, g\me = \rho_{g(x)}$ for $g\in \orth(q)$,
\sm

\item $\rho_{ux} = \rho_x$ for $u\in R\ti$.
\end{enumerate}
We put $\sfR(q) = \{ \rho_x : x\in \wS_q\}$ and $\sfR^+(q) = \{ \rho_x \rho_y : x,y\in \wS_q\}$, and denote by
\begin{equation*}
    \label{refle-0} \Refl(q) \quad \text{and} \quad \Refl^+(q) \end{equation*}
the subgroups of $\orth(q)$ generated by $\sfR(q)$ and $\sfR^+(q)$ respectively. Their elements are products of reflections and products of an even number of reflections respectively, the later are sometimes referred to as {\em rotations}. Some facts:  \sm

\begin{inparaenum}[(a)] \item \label{refle-a} By \ref{refle-iii}, the groups $\Refl(q)$ and $\Refl^+(q)$ are normal subgroups of $\orth(q)$ satisfying $\Refl^+(q) \ideal \Refl(q)$ with $[\Refl(q): \Refl^+(q)] \le 2$ since $\rho_y = \rho_x (\rho_x \rho_y)$
for $x,y\in \wS_q$. But note, by \ref{refle-iia},  
\begin{equation}\label{refle-a1}
  \Refl^+(q) = \Refl(q), \quad \text{if $\wS_q \cap \rad (b_q) \ne \emptyset$.}
\end{equation}%
This can only happen if $2=0$ in $R$ since $2q(x) = b(x,x) = 0$ for $x\in \wS_q \cap \rad(b_q)$. 
\lv{
Indeed, $\rho_e = \Id_M$ for any $e\in \wS_q \cap \rad(b_q)$. Therefore $\sfR(q) \subset \sfR^+(q)$ and then $\Refl(q) = \lan \sfR(q)\ran \subset \lan \sfR^+(q)\ran = \Refl^+(q) \subset \Refl(q)$.
}
\sm 

\item \label{refle-one} ({\em Example rank $1$}) Let $(M,q)$ be a quadratic space with $M$ of constant rank $1$. Since $\GL_R(M) = R\ti \Id_M$, it follows that
\begin{equation} \label{refle-one1}
  \mu_2(R) \simlgr \orth(q), \quad x \mapsto x \Id_R
\end{equation}
is an isomorphism of groups. An $x\in \wS_q$ is unimodular by \ref{mx_lem}\eqref{mx_lema}, which forces $M$ to be free of rank $1$ with basis $\{x\}$, and then $\rho_x = -\Id_M$. Thus, assuming $\wS_q \ne \emptyset$, we get $\sfR(q) = \{-\Id\}$, $\Refl^+(q) = \{\Id\}$, and $\Refl(q) = \{\pm \Id_R\}\subsetneq \orth(q)$ in general, while $\orth(q) = \Refl(q)$ if, for example,  $R$ is an integral domain. 
\sm

\item \label{refle-c} ({\em Base change}) Let $S\in \Ralg$ and let $x\in \wS_q$. Then $x\ot 1_S \in S\ti$ and $\rho_{x\ot 1_S} = \rho_x \ot \Id_S$. Hence we have canonical group homomorphisms
\begin{equation}\label{refle2}
   \Refl(q) \to \Refl(q_S) \quad \text{and} \quad\Refl^+(q) \to \Refl^+(q_S),
\end{equation}
induced by $g \mapsto g\ot \Id_S$.

In particular, let $S=\bar R = R/\fra$ where $\fra \subset \Jac(R)$ is an ideal of $R$, and denote reduction mod $\fra$ by a bar, thus $(M,q)_{\bar R} = (\bar M, \bar q)$ and $M \to \bar M$, $m \mapsto \bar m$ is the canonical map. Given $x\in \wS_{\bar q}$, there exists $m\in M$ such that $\bar m = x$. We have $\overline{q(m)} = \bar q(\bar m) = \bar q(x) \in R\ti$. Hence, by Nakayama, 
$q(m) \in R\ti$. Thus the maps $\wS_q \to \wS_{\bar q}$, $m \mapsto \bar m$, and $\sfR(q) \to \sfR(\bar q)$, $\rho_m \mapsto \overline{\rho_m} = \rho_{\bar m}$,  are surjective. Consequently, the maps \eqref{refle2} are surjective too:
\begin{equation}\label{refle4}
   \Refl(q) \twoheadrightarrow \Refl(\bar q) \quad \text{and} \quad
    \Refl^+(q) \twoheadrightarrow \Refl^+(\bar q).
\end{equation}

\item Let $R= R_1 \times \cdots \times R_n$, hence $(M,q) = (M_1, q_1) \times \cdots \times (M_n,q_n)$ as in \ref{qf}\eqref{qf-redc}, where the $(M_i, q_i)$ are quadratic $R_i$--modules which are regular (nonsingular respectively) if $(M, q)$ is so. Then
\begin{equation*} 
\wS_q = \wS_{q_1} \times \cdots \times \wS_{q_n} \quad\text{and} \quad
      \sfR(q) = \sfR(q_1) \times \cdots \times \sfR(q_n).
\end{equation*}
Indeed, $R\ti = R_1\ti \times \cdots \times R_n\ti$ and therefore $x=(x_1, \ldots, x_n) \in M$ has $q(x) = (q_1(x_1), \ldots, q_n(x_n))\in R\ti \iff q_i(x_i) \in R_i\ti$ for $i=1, \ldots, n$. For such an $x$ we have $\rho_x = \rho_{x_1} \times \cdots \times \rho_{x_n}$. It follows that
\[
\Refl(q)  \subset \Refl(q_1) \times \cdots \times \Refl(q_n),
\]
an inclusion which is in general not an equality.

For example, let $R=R_1 \times R_2$ with $2\in R\ti_i$ and let $(M_i, q_i) = (R_i, \lan 1_{R_i} \ran)$. Then $\sfR(q) = \{-\Id_M\}$ by \eqref{refle-one} and thus $\Refl(q) = \{ \pm \Id_M\}$, while the group $\Refl(q_1) \times \Refl(q_2) = \{ \pm\Id_{R_1}\} \times \{ \pm \Id_{R_2}\}$ has order $4$. Or, from a more general perspective, let $(M_i, q_i)$ be faithful quadratic $R_i$--modules over $R$--fields $R_i$ with $2\in R\ti_i$. The determinant induces a surjective map $\Refl(q_1) \times \Refl(q_i) \to \ZZ/2\ZZ \times \ZZ/2\ZZ$, while the image of $\Refl(q)$ is the diagonal $\ZZ/2\ZZ$.
\sm

However, we always have
\begin{equation} \label{refle5}
 \Refl^+(q) = \Refl^+(q_1) \times \cdots \times \Refl^+(q_n).
\end{equation}
The inclusion $\Refl^+(q) \subset \Refl^+(q_1) \times \cdots \times \Refl^+(q_n)$ follows from $\sfR^+(q) = \sfR^+(q_1) \times \cdots \times \sfR^+(q_n)$. For the other inclusion, it suffices by symmetry to show $\Refl^+(q_1) \times \{\Id_{M_2}\} \times \cdots \times \{\Id_{M_n}\} \subset \Refl^+(q)$. We can assume that $\wS_q \ne \emptyset$. For $1\le i \le n$  let $x_i \in \wS_{q_i}$ and let $y_1 \in \wS_{q_1}$.
Then
\[ (\rho_{x_1} \rho_{y_1}, \Id_{M_2}, \ldots, \Id_{M_n})
   = (\rho_{x_1}, \rho_{x_2}, \ldots, \rho_{x_n}) \,
      (\rho_{y_1}, \rho_{x_2}, \ldots, \rho_{x_n})
\]
lies in $\Refl^+(q)$, which implies our claim.
\end{inparaenum}

\begin{thm}[Cartan-Dieudonn\'e-Kneser] \label{CDK}
Let $(M,q)$ be a quadratic space over a field $F$. Then $\orth(q) = \Refl(q)$, unless $F=\FF_2$ is the field of two elements and $(M,q)$ is hyperbolic of dimension $4$.
\end{thm}

\textbf{Remarks.} \begin{inparaenum}[(a)] \item For regular quadratic forms, Theorem~\ref{CDK} is proven in \cite[II, \S4, \S10]{Dieu} and for nonsingular forms in \cite{kneser-Sem}. An exposition is given  in \cite[(3.5)]{Kneser} (see \cite[page 19]{Kneser} for the proof in characteristic $2$). 
In the exceptional case of \ref{CDK} it is easily verified that $\Refl(q) \subsetneq \orth(q)$, see the Example  \ref{refplusfield-exa}, in fact $[\orth(q) : \Refl(q)] = 2$ holds in that case (\cite[page 14]{Kneser}). We describe $\Refl^+(q)$ over fields in \ref{replusLG}. \sm 

\item\label{refplusfield-exb} One may wonder if $\orth(q) = \Refl(q)$ holds in greater generality. Let $R$ be a local ring with maximal ideal $\gm$, and let $(M, q)$ be a regular quadratic $R$--module. Then indeed $\orth(q) = \Refl(q)$ holds, unless $R/\m = \FF_2$ and $\rank M \le 4$, \cite[4.6]{Kneser}. But this is not true for an arbitrary nonsingular $q$. Indeed, if $R=\FF_2[x]/(x^2)$ and $q= q_{0, 2m+1}$ is the split quadratic form of odd rank, then $\orth(q)\ne \Refl(q)$, as shown in \cite[2.6.1.22]{CF}. 
    
\end{inparaenum}

\comments{
(2021-09-30) In Faulkner's Projective Geometry book he proves for a regular $(V,q)$ over a field that $\orth(V,q)$ is generated by reflections and Eichler transformations (Thm.~9.8) and that an Eichler transformation is a product of two reflections unless $V= H_1 \perp H_2$ with $H_i$ hyperbolic over the field $\FF_2$ of $2$ elements (Lemma~9.9). 

(2025-05-01) Faulkner's Eichler transformations are called Siegel transformations in \cite[p.~60]{Ba} and \cite[1.4]{Knebusch}. The term Eichler transformation seems to be traditional. Background: Eichler introduced these maps on page 13 of his book with a note to the ``Anmerkungen'' at the end of this book where he says that these types of transformations have first been used by Siegel in his paper {\em \"Uber die analytische Theorie der quadratischen
Formen II}, Annals of Maths, Princeton \textbf{36} (1935) 230--263. }
\ms

To define the special orthogonal group $\SO(q)$, we use the notion of the discriminant algebra of a nonsingular quadratic form $q$, which is a subalgebra of the Clifford algebra of $q$. We review these concepts in \ref{clialg} and \ref{qfdi}.
 
\subsection{Clifford algebras}\label{clialg}
Let $(M,q)$ be a quadratic $R$--module and $\Cli(M,q)= \Cli(q)$ be its Clifford algebra, see e.g.\ \cite[IV, \S1]{K}. It is a $(\ZZ/2\ZZ)$--graded algebra,
\[ \Cli(q) = \Cli_0(q) \oplus \Cli_1(q); \]
its even part $\Cli_0(q)$ is referred to as the even Clifford algebra. We can and will identify $M$ with a submodule of $\Cli_1(q)$. \sm 

An isometry $g \co (M_1, q_1) \to (M_2, q_2)$ of quadratic $R$--modules gives rise to an isomorphism $\Cli(g) \co \Cli(M_1, q_1) \simlgr \Cli(M_2, q_2)$ of Clifford algebras which respects the $(\ZZ/2\ZZ)$--gradings. In particular, denoting by $\Aut_R
\big(\Cli(q), \Cli_1(q)\big)$ the group of automorphisms $\ga$ of $\Cli(q)$ satisfying $\ga\big(\Cli_1(q)\big)= \Cli_1(q)$ we have a group homomorphism
\[ \orth(M,q) \to \Aut_R\big(\Cli(q), \Cli_1(q)\big), \quad g \mapsto \Cli(g).
\]%

The Clifford algebra respects base change: for $S\in\Ralg$ the canonical map $M_S \to \Cli(q)_S$ extends to an  isomorphism $\Cli(q_S) \simlgr \Cli(q)_S$ of $(\ZZ/2\ZZ)$--graded $S$--algebras.

Let $R=R_0 \times \cdots \times R_n$ be a direct product of rings. As we have seen in \ref{qf}\eqref{qf-redc} any quadratic $R$--module $(M,q)$ is uniquely a direct product \[ (M,q) = (M_0, q_0) \times \cdots \times (M_n,q_n) \]
of quadratic $R_i$--modules $(M_i, q_i)$, and conversely.
The Clifford algebra of $q$ respect this decomposition: 
\begin{equation}  \label{clialg1}
 \Cli(q) = \Cli(q_0) \times \cdots \times \Cli(q_n),  
\end{equation}
i.e., the $(\ZZ/2\ZZ)$--graded $R$--algebra $\Cli(q)$ can be identified with the direct product of the $(\ZZ/2\ZZ)$--graded $R_i$--algebras $\Cli(q_i)$. In particular, the above says that we can use the rank decomposition of a quadratic module.

\subsection{Discriminant algebras} \label{qfdi}  
Let $(M,q)$ be a faithful quadratic $R$--space\footnote{We will not venture into defining discriminant algebras for arbitrary quadratic $R$--modules, see \cite[IV, (4.8.4)]{K} which shows that $\Cli(q)^{\Cli(q_0)}$ does not commute with scalar extensions.}, and let $\Cli(M,q)= \Cli(q)$ be its Clifford algebra, \ref{clialg}. The {\em discriminant algebra $\Dis(q)$ of $(M,q)$\/} is the subalgebra of $\Cli(q)$ centralizing $\Cli_0(q)$:
\[ \calD:= \Dis(q) = \Cli(q)^{\Cli_0(q)}. \]

\begin{inparaenum}[(a)] \item ({\em Basic properties}) A discriminant algebra $\calD$ is a quadratic $R$--algebra in the sense of \cite[I, (1.3.6) and III, \S4]{K}, i.e., its underlying $R$--module is projective of rank $2$. In particular, it is commutative and carries a so-called {\em standard involution\/} $\si_{\calD}$.   Discriminant algebras respect base change and direct products of base rings.

In more detail, let $R=R_1 \times \cdots \times R_n$ be a direct product of rings, and let
\begin{equation} \label{discralg0}
    (M,q) = (M_1, q_1) \times \cdots \times (M_n,q_n)
\end{equation}
be the corresponding decomposition. The $(M_i, q_i)$ are faithful quadratic $R_i$--spaces, which are all regular if $(M,q)$ is so. The Clifford and even Clifford algebra decompose correspondingly,
and this gives rise to the decomposition
\begin{equation}
        \label{discralg-a1} \Dis(q_1 \times \cdots \times q_n) = \Dis(q_1) \times \cdots \times \Dis(q_n).
\end{equation}
We can use \eqref{discralg-a1} and the rank decomposition of $M$ to infer properties of $\Dis(q)$ from the known case of $M$ having constant rank.
\sm

\item ({\em The group homomorphism $\Dis$}) By the universal property of the Clifford algebra $\Cli(q)$, every $g\in \orth(q)$ induces an automorphism $\Cli(g)$ of the algebra $\Cli(q)$ stabilizing $\Cli_0(q)$ and $\Cli_1(q)$, and hence an automorphism of the $\ZZ/2\ZZ$--graded algebra $\calD = \calD_0 \oplus \calD_1$. Thus, we get a homomorphism of groups,
\begin{equation}\label{sog-dishom}
 \Dis \co \orth(q) \longto \Aut (\calD, \calD_1), \quad g \mapsto
    \Cli(g)|_{\calD} =:\Dis(g).
\end{equation}

\item \label{discralg-even} ({\em Even rank}) Assume $M$ has constant even rank. Thus $q$ is regular by  \ref{qf}\eqref{quadfoe}. In this case,  $\calD$ is the centre of $\Cli_0(q)$, and a quadratic \'etale  $R$--algebra \cite[IV, (2.2.3)]{K}. Hence, by \cite[III, (4.1.2)]{K}, the automorphism group scheme of $\calD$ can be identified with  the $R$--group scheme of locally constant functions with values in $\ZZ/2\ZZ = \{0,1\}$,  
    \begin{equation}\label{evdi1}
 \uAut(\calD) \simla (\ZZ/2\ZZ)_R. 
\end{equation}
\inparcom{(2020-12-11) We cannot conclude:  ``{\sf In this case $\si_{\Dis(q)}  \ne \Id_{\Dis(q)}$}'' because this is not true for the zero ring. A rank-$2$ module over the zero ring is zero.   Of course  $\si_{\Dis(q)}  \ne \Id_{\Dis(q)}$  when $R\ne 0$. }
\sm

\item \label{discralg-d} ({\em Odd rank}) Assume $M$ has constant odd rank. Then the   discriminant algebra $\calD $ is the centre of $\Cli(q)$. So $\calD$ inherits the $(\ZZ/2\ZZ)$--grading of $\Cli(q)$, 
  \[ \calD = \calD_0 \oplus \calD_1, \qquad \calD_j = \calD \cap\, \Cli_j(q), \; j = 0,1.   
    \]
  Moreover,   \end{inparaenum} 
  \begin{enumerate}[label={\rm (\roman*)}]
 \item $\calD_0= R \, 1_{\Cli(q)}$ is free of rank $1$.

 \item $(\calD_1, \theta)$ is a discriminant module in the sense of  \cite[III, \S3]{K}, 
     where $\theta\co \calD_1\ot_R \calD_1 \to \calD_0 = R$ is the restriction of the multiplication of $\Cli(q)$; it is called the {\em discriminant module of $q$}.
     By \cite[III, (3.2.1)]{K}, the canonical map
     \begin{equation}\label{discrald-d0}
     \mu_2(R) \simlgr \Aut(\calD, \calD_1), \quad x \mapsto \Id_{\calD_0} \oplus \, x\Id_{\calD_1} 
     \end{equation}
    is an isomorphism of groups. 
    


 \item Let $g\in \orth(M,q)$. Then the automorphism $\Cli(g)$ of $\Cli(q)$ stabilizes $\calD$ and the homogeneous parts $\calD_j$, $j=0,1$. It acts on $\calD_1$ by the determinant of $g$. Hence 
     \begin{equation} \label{discralg-d1} 
        \orth(q) \to \Aut(\calD, \calD_1) , \quad g \mapsto \Cli(g)|_{\calD} 
     \end{equation} 
    is a homomorphism of groups, which is surjective by \eqref{discrald-d0} and \eqref{qf-detgr1}.
    
 \item By \cite[IV, (4.3.1)]{K}, the standard involution of the quadratic $R$--algebra $\si_{\Dis(q)}$ is the grading automorphism,
\begin{equation}\label{discralg-b1}
 \si_{\Dis(q)}(d_i ) = (-1)^id_i, \quad d_i \in \Dis_i(q), i=0,1.
\end{equation}

\item \label{discralg-dv}  If $2\in R\ti$, then $\Dis(q)$ is an \'etale (= separable) $R$--algebra, \cite[IV, (3.2.5)]{K}.

\end{enumerate} 

\begin{inparaenum}[(a)] \setcounter{enumi}{4}
  \item {\em If $(M,q)$ is regular, then $\Dis(q)$ is an \'etale $R$--algebra}. Indeed, after applying the rank decomposition, we can assume that $M$ has constant rank $n$. If $n$ is even, the claim follows from \eqref{discralg-even}. If $n$ is odd, then $2\in R\ti$ by \ref{bfLG}\eqref{bfLG-odd}, and the claim follows from  \ref{discralg-dv} above. 
\end{inparaenum}

\subsection{Special orthogonal groups.} \label{sog} Let $(M,q)$ be  a faithful quadratic $R$--space. By the universal property of the Clifford algebra $\Cli(q)$, every $g\in \orth(q)$ induces an automorphism $\Cli(g)$ of $\Cli(q)$ stabilizing $\Cli_0(q)$ and hence also $\Dis(q)$. It is immediate that we get a homomorphism of groups,
\begin{equation}\label{sog-dishom}
 \Dis \co \orth(q) \longto \Aut \big(\Dis(q)\big), \quad g \mapsto
    \Cli(g)|_{\Dis(q)} =:\Dis(g).
\end{equation}
The {\em special orthogonal group $\SO(q)$\/} is defined as
\begin{equation}
  \label{sog1} \SO(q) = \Ker(\Dis) = \{ g\in \orth(q): \Dis(g) = \Id_{\Dis(q)}\}.
\end{equation}

\begin{inparaenum}[(a)]
\item\label{sog-ba} ({\em Functoriality}) Let $S\in \Ralg$. Both $\orth(q)$ and $\Dis(q)$ respect base change. It follows that so does $\SO(q)$. \sm

\item ({\em Rank decomposition}) \label{sog-ra} As before, we first consider the case where $R$ is a direct product of rings, $R=R_1 \times \cdots \times R_n$. The decompositions \eqref{discralg0} and \eqref{discralg-a1} are decompositions into quadratic $R_i$--spaces and $R_i$--algebras respectively. This implies that $\SO(q)$ decomposes correspondingly:
\begin{equation}  \label{sog2}
\SO(q_1 \times \cdots \times q_n) = \SO(q_1) \times \cdots \times \SO(q_n).
      \end{equation}
Thus, we can often reduce proofs to the case of $M$ having constant rank.
\sm

\item ({\em $\SO$ and determinants}) \label{sog-b} We always have
    \[ \SO(q) \subset \{ g\in \orth(q) : \det(g) = 1\}, \]
    with equality if $M$ has odd rank or if $2\in R\ti$.
    This is for example proven in \cite[IV, (5.1.1)]{K}, but note the misprint in (3) of loc.\ cit., where ``$\subset$'' should be ``$=$".   \sm

\item \label{sog-c} ({\em Knebusch's\/ $\SO(q)$ definition}) Let $(M,q)$ be a regular quadratic module with a faithful $M$, and let $g\in \orth(q)$. Then $g\in \SO(q) \iff g\ot_R \Id_{R/m} \in \SO(q_{R/\m})$ for all maximal $\m \in \Spec(R)$.

    This characterization will be proven later in Proposition~\ref{somax}.  It is the definition of $\SO(q)$ in Knebusch's paper \cite{Knebusch}.
\sm

\item \label{sog-d} ({\em Dickson homomorphism}) We denote by
\[ \ZZ/2\ZZ(R) = \{ \eps \in R : \eps = \eps^2 \}
\]
the group of idempotents with the operation $\eps \star \eps' = \eps + \eps' - 2 \eps \eps'$. Every $\eps \in \ZZ/2\ZZ(R)$ induces a decomposition $\calD :=\Dis(q) = \eps \calD\times (1-\eps) \calD$ of the $R$--algebra $\calD$. By \cite[III, (4.1.2)]{K} the map
\begin{equation} \label{sog-d00}
\psi \co \ZZ/2\ZZ(R) \to \Aut(\calD), \quad
 \eps \mapsto \psi(\eps)
    = \begin{cases} \rho_{\calD} & \text{on $\eps \calD$} \\
                    \Id_{\calD} &\text{on $(1-\eps) \calD$}
     \end{cases}
\end{equation}
 is an injective  homomorphism of groups. Note
\begin{equation} \label{sog-d0}
\psi(0_R) = \Id_{\calD} \quad \text{and} \quad \psi(1_R) = \rho_{\calD}.
\end{equation}
{\em Assume $q$ is regular, whence $\calD$ is \'etale.} Then, by loc.\ cit., $\psi$ is an isomorphism and we can define the {\em Dickson map\/} $\Di$ as $\psi\me \circ \Dis$:
\begin{equation} \label{sog-d1} \begin{split}
 \xymatrix@C=20pt{\orth(q) \ar[rr]^{\Di} \ar[dr]_{\Dis} && \ZZ/2\ZZ(R) \ar[dl]^\psi_\cong \\ & \Aut(\Dis(q))}
\end{split}\quad .  \end{equation}
We then have
\begin{equation} \label{sog-d2}
 \SO(q) = \Ker(\Di), \qquad(\text{$q$ regular}).
\end{equation}
Moreover, {\em if $(M,q)$ is regular and contains a hyperbolic plane as direct summand, then
\begin{equation}  \label{sog-d3}
1 \longto \SO(q) \longto \orth(q) \xrightarrow{\Di} \ZZ/2\ZZ(R) \longto 1
\end{equation}
is split exact} (\cite[IV, (5.2.2)]{K}).
\sm

\item\label{sog-e} ({\em $\orth(M,q)$ with $M$ of odd rank}) Let $(M,q)$ be a faithful quadratic space. Independent of the rank, the map
    \[ z_M \co \mu_2(R) \longto \orth(q), \quad r \mapsto r\Id_M
    \]
    is an injective group homomorphism with central image. If $M$ has constant odd rank, then $z_M$ is a section of $\det \co \orth(q) \to \mu_2(R)$, and
    \begin{equation} \label{sog-e1}
        \orth(q) \cong \mu_2(R) \times \SO(q) \qquad (\text{$M$ odd rank}).
    \end{equation}
In particular, if $M$ has constant rank $1$, then $\orth(q) \cong \mu_2(R)$, cf.~\ref{refle}\eqref{refle-one}, and hence 
\begin{equation}\label{sog-e2}
  \SO(q) = \{\Id_M\}= \Refl^+(q) \qquad(\rank M = 1).
\end{equation}

\item\label{sog-hyp} ($(M,q) = \HH$) Let $(M,q) = \HH$ be a hyperbolic plane with hyperbolic pair $(e,f)$. Every $u\in R\ti$ induces an orthogonal transformation $\rho_{e-f} \rho_{e - uf} \in \Refl^+(q) \subset \SO(q)$,  
acting as  $e \mapsto ue$, $f\mapsto u\me f$. One knows 
\begin{equation}  \label{sog-hyp1}
 \SO(\HH) = \{ \rho_{e-f}\,\rho_{e-uf} : u \in R\ti\} \cong R\ti
\end{equation}
(\cite[III, (2.2)]{Ba}, \cite[V, (2.6.3)]{K}). 
\end{inparaenum}

\pcomments{(2020-12-11 for later use) The $\ZZ/2\ZZ$-grading defines an action  of $\bmu_2$ on $\Spec(\Dis(q))$, see \cite[VIII.1]{SGA3}. In other words,
$\Spec(\Dis(q)) \to  \Spec(R)$ is a $\bmu_2$--torsor and one needs to check that it is the same as the $\bmu_2$--torsor provided by the map $\det\co \uO(q_0) \to \bmu_2$. }

\begin{lem}[Realizing the standard involution of $\Dis(q)$]\label{refso}\label{refso} 
Let $(M,q)$ be a faithful quadratic space,
let $x\in \wS_q$ and let $\rho_x$ be the associated reflection. Then the automorphism $\Dis(\rho_x) \in \Aut(\Dis(q))$ is the standard involution of\/ $\Dis(q)$:
\begin{equation}\label{refso1}
 \Dis(\rho_x) = \si_{\Dis(q)}.
 \end{equation}
In particular,
\begin{equation}\label{refso2} \begin{split}
 \Refl^+(q) &\; \subset \; \SO(q), \quad\text{and} \\
   \text{$q$ regular} \quad & \implies \quad \Di(\rho_x) = 1_R \in \ZZ/2\ZZ(R).
\end{split}   \end{equation}
\end{lem}

\comments{(2026-04-29) The proof of \eqref{refso1} is taken from \cite[C.12]{GN-LG}. We use the formula \eqref{reso0} again later. I felt it does not make sense to refer the reader for  the proof of \eqref{refso1} to \cite[C.12]{GN-LG}. }

\begin{proof}
 The element $x\in M \subset \Cli_1(q)$ is invertible in $\Cli(q)$ with inverse $x\me = q(x)\me x$. We will use the well-known formula relating $\rho_x(m)$, $m\in M$, with the inner automorphism of $\Cli(q)$ induced by $x$:
\begin{equation}\label{reso0}
 \rho_x(m) = - xmx\me
\end{equation}
which follows from $ xmx\me = (xm)(xq(x)\me) = (-mx + b_q(m,x))(xq(x)\me)
= - m + b_q(m,x) q(x)\me x = - \si_x(m)$.
It implies
\begin{equation*}\label{reso00}
 \Dis(\rho_x)(c_j) = (-1)^j xc_jx\me, \qquad (c_i \in \Cli_j(q), \, j = 0,1).
\end{equation*}
For the proof of \eqref{refso1}, we can without loss of generality assume that $M$ has constant rank.

Suppose $M$ has constant even rank. By \cite[IV, (4.3.1.4)]{K}, $\si_{\calD}(d) x = x d$ holds for $d\in \calD$. Since $\calD \subset \Cli_0(q)$ we get $\Dis(\rho_x)(d) = xdx\me = \si_{\calD}(d)$.

Suppose $M$ has constant odd rank. Since then $\calD = \rmZ(\Cli(q))$, we obtain for $d_j \in \calD_j$, $j=0,1$,  that $\Dis(\rho_x)(d_j) = (-1)^j xd_j x\me = (-1)^j d_j x x\me = (-1)^j d_j = \si_{\calD}(d_j)$ by \eqref{discralg-b1}.
\sm 

The formula \eqref{refso1} 
implies the first part of \eqref{refso2}. If $q$ is regular, then $\Di(\rho_x) = \psi\me(\si_{\Dis(q)}) = 1$ by \eqref{sog-d0}, proving the implication \eqref{refso2}.
\end{proof}
\sm

The next three lemmata are applications of the formulas in Lemma~\ref{refso}.
In \ref{Lag} we have seen that $\orth(q)$ acts transitively on the set of Lagrangians of the same rank in the LG case. It is well-known, see for example \cite[Prop.~3.7]{Co3}, that over fields the action of $\SO(q)$ on this set has two orbits. It is also shown in loc.\ cit.\ that $\SO(q)$ acts transitively on the set of totally isotropic subspaces that are not Lagrangians. As a first application of \eqref{refso2}, we will extend this transitivity result to the LG case in Lemma~\ref{traq}. 
\comments{(2025-04-08) Previously we had the following statement here: 

{\tt While the proof of \ref{traq} is ad-hoc, we revisit transitivity in Corollary~\ref{traco} from a broader perspective: it is a special case of Demazure's Conjugacy Theorem~\ref{thm_conj_demazure} for $G=\uSO(q)$.} 

The last sentence to be removed if appendix H on parabolic subgroups and Severi-Brauer schemes is not included in the paper. The application of Demazure's Theorem requires identifying totally isotropic complemented submodules as points of a parabolic subgroup of $\uO(q)$ or $\uSO(q)$.}

While the proof of \ref{traq} is ad-hoc, we revisit transitivity in Corollary~\ref{traco} from a broader perspective: it is a special case of Demazure's Conjugacy Theorem~\ref{thm_conj_demazure} for $G=\uSO(q)$.

\begin{lem} \label{traq} Let $R$ be an LG ring, and let $N$ and $N_1$ be two  totally isotropic direct summands of a faithful quadratic $R$--space $(M,q)$ for which $2 \rank_\p N= 2 \rank_\p N_1 < \rank_\p M$ holds for all $\p \in \Spec(R)$.
Then there exists $f\in \SO(q)$ satisfying $f(N) = N_1$.
\end{lem}

\begin{proof} By Proposition~\ref{quadrepII} there exist totally isotropic submodules $N'$ and $N_1'$ of $M$ such that $(N \oplus N', q|_{N \oplus N'}) \simlgr \HH(N)$ and $(N_1 \oplus N_1', q|_{N_1 \oplus N_1'}) \simlgr \HH(N_1)$ via $b_q$. By \eqref{unimod-b2}, there exists an $R$--module isomorphism $g \co N \simlgr N_1$. Since $(N,N')$ and $(N_1, N_1')$ are hyperbolic pairs, it follows from \eqref{hps2} that the $R$--module isomorphism $g$ extends to an isometry $\wtl g $ of hyperbolic spaces
\[  \wtl g \co (N \oplus N', q|_{N \oplus N'}) \simlgr (N_1 \oplus N_1', q|_{N_1 \oplus N_1'}).\]
Let $P$ and $P_1$ be the orthogonal complement of $N \oplus N'$ and $N_1 \oplus N_1'$ respectively. Thus,
\[
 (N \oplus N') \perp P = M = (N_1 \oplus N_1') \perp P_1.
 \]
Note that $P$ and $P_1$ are faithful projective $R$--modules, that $q|_P$ is nonsingular by \eqref{qf-perp0} and that $q_{N\oplus N'} \cong q_{N_1 \oplus N_1'}$ is regular. Hence, by Witt cancellation \ref{canqf}, there exists an isometry $h \co (P, q|_P) \simlgr (P_1 , q|_{P_1})$. Then $f=\wtl g \oplus h$ is an isometry of $(M,q)$ with $f(N)=N_1$.

To see that we can modify $f\in \orth(q)$ to get an $\wtl f \in \SO(q)$ with $\wtl f(N) = N_1$, it is now harm to assume that $M$ and $N$ have constant rank, cf.~\ref{qf}\eqref{qf-redc} and \ref{sog}\eqref{sog-ra}. We will treat the cases of odd and even rank of $M$ separately.

Assume that $M$ has constant odd rank, and let $r= \det(f)\me$. Then $r\in \bmu_2(R)$ and $r\Id_M \in \orth(q)$ by \ref{sog}\eqref{sog-e}. Hence $\wtl f = rf\in \SO(q)$ by \ref{sog}\eqref{sog-b}, and, obviously, $\wtl f(N) = N_1$.

Finally, let $M$ have constant even rank. In this case, $q$ is regular and $\Di(f) \in (\ZZ/2\ZZ)(R)$ by \eqref{sog-d1}. If $\Di(f) = 0$, i.e., $\Dis(f) = \Id_{\Dis(q)}$, then $f\in \SO(q)$ by \eqref{sog1}. Otherwise, after possibly further decomposing $R$, we can assume $\Di(f) = 1$. Since $q|_P$ is nonsingular (even regular), there exists $x\in P$ with $q(x) \in R\ti$. The reflection $\rho_x$ of $M$ has $\Di(\rho_x) = 1$ by \eqref{refso2}. Hence $\wtl f = f \rho_x \in \SO(q)$ and, clearly, $\wtl f(N) = N_1$. \end{proof}

\begin{lem}\label{lem_transitivity_SO} Let $R$ be an \new LG ring \enew and let
$(M,q)$ is a faithful quadratic $R$--space. Furthermore, let $x,y\in \wS_q$.

 \begin{enumerate}[label={\rm (\alph*)}]
  \item  \label{lem_transitivity_SO2} The following are equivalent:
 \begin{enumerate}[label={\rm (\roman*)}]
 \item \label{lem_transitivity_SO2_i}  $y \in \SO(q)\cdot ( \mu_2(R) \cdot x)$;
 \item \label{lem_transitivity_SO2_ii} $y \in \orth(q)\cdot x$.
 \end{enumerate}

\item \label{lem_transitivity_SOb}   Suppose $\rank_R M \ge 2$. Then 
 \[ y \in \SO(q)\cdot x \quad \iff \quad y \in \orth(q)\cdot x.
\]
under any of the following conditions:
\begin{enumerate}[label={\rm (\Roman*)}]
 \item \label{lem_transitivity_SO1} $q$ is regular, or  
 
 \item \label{lem_transitivity_SO3} $R$ is a field. 
\end{enumerate}
\end{enumerate}
%
\end{lem}
\sm

We note that the equivalence in \ref{lem_transitivity_SO1} or \ref{lem_transitivity_SO3} is not true if $\rank M = 1$, since then $\SO(q) = \{\Id_M\}$ while $\orth(q) \ne \{\Id_M\}$ in general, see \ref{sog}\eqref{sog-e}.

\begin{proof} For \ref{lem_transitivity_SO2} and \ref{lem_transitivity_SOb},  after 
applying the rank decomposition \ref{qf}\eqref{qf-redc} and using \ref{revLG}\eqref{revLG-aa}, we can assume that $M$ has constant rank $r$.

\ref{lem_transitivity_SO2} Since $(\mu_2(R) \Id_M)) \cdot \SO(q) \subset \orth(q)$, we only need to prove \ref{lem_transitivity_SO2_ii} $\implies$ \ref{lem_transitivity_SO2_i}. If $r$ is odd, then $\orth(q) = (\mu_2(R) \Id_M) \cdot \SO(q)$ by \eqref{sog-e1}. If $r$ is even, then $r\ge 2$ and $q$ is regular by \ref{qf}\eqref{quadfoe}. It is therefore enough to prove $\Longleftarrow$ in \ref{lem_transitivity_SOb}\ref{lem_transitivity_SO1}.

\ref{lem_transitivity_SOb}\ref{lem_transitivity_SO1} We write $y = g\cdot x$ for some $g\in \orth(q)$ and assume $g\notin \SO(q)$. Since $R$ is connected and $g\notin \SO(q)$, we have $\Di(g) = 1$. Moreover, by Lemma~\ref{mx_lem},  $M_x = \{m\in M :  b_q(m,x) = 0 \}$ is projective of rank $r-1\ge 1$ and $q_x = q|_{M_x}$ is nonsingular. Hence, by  \ref{LGqdi}, there exists $z \in  M_x$ such that $q(z) \in  R^\times$. The reflection $\rho_z$ satisfies $\rho_z(x) = x$, so that $(g\circ \rho_z)(x) = y$. By \eqref{refso2}, $\Di(g \circ \rho_x) = \Di(g) \star \Di(\rho_z) = 1 \star 1 = 0\in \ZZ/2\ZZ(R)$, proving $g\circ \rho_z \in \SO(q)$ by applying \eqref{sog-d2}.

\ref{lem_transitivity_SOb}\ref{lem_transitivity_SO3} If $R$ has characteristic $\ne 2$, the quadratic form is regular, and so the equivalence follows from \ref{lem_transitivity_SO1}. On the other hand, if $R$ has characteristic $2$, then $\mu_2(R) = \{1\}$ and the equivalence follows from \ref{lem_transitivity_SO2}. \end{proof}

\comments{(2025-04-14) The previous Lemma \ref{replusLG} only considered fields, which was enough for the use in these notes. But its proof had a gap. }

\begin{lem}\label{replusLG} Let $(M,q)$ be a faithful quadratic $R$--module. \sm 

\begin{inparaenum}[\rm (a)] \item\label{replusLGa} Then $\Refl^+(q)$ is a normal subgroup of $\SO(q) \cap \Refl(q)$ of index $\le 2$. \sm 

\item\label{replusLGb} We have 
\begin{equation}   \label{replusLGb1} 
            \SO(q) \cap \Refl(q) = \Refl^+(q) 
\end{equation}
 in any one of the following cases:
\end{inparaenum}
\begin{enumerate}[label={\rm (\roman*)}]
  \item\label{replusLGbi} $M$ has constant even rank;
  
  \item\label{replusLGbii} $M$ has constant odd rank and $2$ is not a zero-divisor in $R$; 
      
  \item \label{replusLGbiii} $R=F$ is a field. 
\end{enumerate}
\sm 

\noindent More precisely, let $R=F$ be a field. Then
\begin{enumerate}[label={\rm (\Roman*)}]
  \item \label{replusLGb4} $\SO(q) = \Refl^+(q)$, 
   
   \noindent unless $F=\FF_2$ and $(M,q)$ is a hyperbolic space of dimension $4$, 
   \sm 
   
  \item \label{replusLGb3} $\SO(q) = \orth(q) = \Refl(q) = \Refl^+(q)$, 
  
  \noindent if $F=\FF_2$ and $\dim_F M $ is odd.
\end{enumerate}
\end{lem}

\begin{proof} \eqref{replusLGa} follows from $[\Refl(q) : \Refl^+(q)] \le 2$ and \eqref{refso2}. \sm 

\eqref{replusLGb} Putting $\calD = \Dis(q)$, we first show 
\begin{equation}  \label{replusLGa1} \begin{split}
 \SO(q) \cap &\Refl(q) = \Refl^+(q) \\ & \cup \big(\{ \rho_x : \Dis(\rho_x) = \Id_{\calD}\}\big) \cdot \Refl^+(q).  
\end{split}\end{equation}
For the proof of \eqref{replusLGa1}, let $g\in \SO(q) \cap \Refl(q)$, say $g=\rho_{x_1} \cdots \rho_{x_n}$. We can assume that $n$ is odd. By definition of $\SO(q)$ and \eqref{refso1}, we have $\Id_{\calD} = \Dis(g) = \Dis(\rho_{x_1})\cdots \Dis(\rho_{x_n}) = \si_{\calD}^n = \si_{\calD}$, which settles the case $n=1$. For $n>1$, observe that $g=\rho_{x_1} h$ with $h=\rho_{x_2} \cdots \rho_{x_n}\in \Refl^+(q)$, so that $\Dis(\rho_{x_1}) = \Id_{\calD}$ follows because $\Dis(h) = \Id_{\calD}$. For the proof of the other direction, we know that $\Refl^+(q) \subset \SO(q) \cap \Refl(q)$ by \eqref{refso2}, and $\rho_{x}\in \SO(q)$ if (and only if) $\Dis(\rho_x) = \Id_{\calD}$. In particular, \eqref{replusLGa1} implies that \eqref{replusLGb1} holds whenever $\si_{\calD} \ne \Id_{\calD}$.

If $M$ has constant even rank, then $\Dis(q)$ is a quadratic \'etale (= quadratic Galois) $R$--algebra, and therefore $\si_{\calD} \ne \Id_{\calD}$. By \eqref{discralg-b1}, the same argument works if $M$ has odd rank and $2$ is not a zero-divisor in $R$. 

Finally, let $F$ be a field. By \ref{replusLGbi} and \ref{replusLGbii} we can assume that $F$ is a field of characteristic $2$ and $\dim_F M$ is odd. Then $\rad(b_q)$ is $1$--dimensional and $q|_{\rad(b_q)} \ne 0$ by \ref{qfnsp-vii} of \ref{qf}\eqref{qfnsp}.   Hence $\wS_q \cap \rad(b_q) \ne \emptyset$, and so $\Refl^+(q) = \Refl(q)$ by \eqref{refle-a1}. Applying the Cartan-Dieudonn\'e-Kneser Theorem~\ref{CDK} we therefore have $\orth(q) = \Refl(q) = \Refl^+(q)$. On the other hand, $\Aut(\calD, \calD_1) = \mu_2(R)$ by \cite[III, (3.2.1)]{K}. Since $\mu_2(R) = \{1\}$, we get $\SO(q) = \Ker(\Dis) = \orth(q)$. Hence altogether, we get \ref{replusLGb3} and so also \eqref{replusLGb1}. 

The more precise formula \ref{replusLGb4} follows the Cartan-Dieudonn\'e-Kneser Theorem~\ref{CDK}. 
\end{proof}

\subsection{Example: $F=\FF_2$, $(M,q)$ hyperbolic, $\dim_F M = 4$} \label{refplusfield-exa} 
We consider the exceptional case in \ref{replusLG}\ref{replusLGb4}: $(M,q)= \HH_1 \perp \HH_2$ is the orthogonal sum of two hyperbolic planes $\HH_i$, $i=1,2$, over $F=\FF_2$ with bases formed by hyperbolic pairs $(e_i, f_i)$. Put $x_i = e_i + f_i$. 

One easily sees that $\wS_q = \{x_1, x_2\}$ and that 
\[ \rho_{x_1} = \sw_1 \times \Id_{\HH_2} , \quad \rho_{x_2} = \Id_{\HH_1} \times \sw_2\]
where $\sw_i$ exchanges $e_i$ and $f_i$. Thus, 
\begin{equation}\label{refplusfield-exa1} \begin{split} 
\Refl(q) &= \{ \Id, \, \sw_1 \times \Id_{\HH_2}, \, \Id_{\HH_1} \times \sw_2, \,   \sw_1 \times \sw_2\}
\\ 
\Refl^+(q) &= \{ \Id, \, \sw_1 \times \sw_2\},
\end{split}\end{equation}
in particular, $\Refl(q)$ is a Klein $4$-group. Since $\Dis(q)$ is quadratic \'etale, we know $\Dis(\rho_{x_i}) = \si_{\Dis(q)} \ne \Id_{\Dis(q)}$, and therefore get \eqref{replusLGb1}: $\SO(q) \cap \Refl(q) = \Refl^+(q)$. %
\sm 

That $\Dis(\rho_{x_i}) \ne \Id_{\Dis(q)}$ can also be seen as follows. By \cite[(6.22)]{Kneser}, the discriminant algebra $\calD = \Dis(q)$ is free of rank $2$, spanned by $1_{\Cli(q)}$ and  $t = e_1f_1 e_2 f_2 + f_2 e_2 f_1 e_1$ where the product on the right hand side is the multiplication in $\Cli(q)$. One easily verifies that $e_if_i$, $f_ie_i$ and $t$ are non-zero idempotents. It follows (\cite[I, (1.3.6)]{K} that $\si_{\calD} = 1 - t$, while 
$\Dis(\rho_{x_1})(t) = f_1 e_1 e_2 f_2 + e_1 f_2 f_2 e_1 \ne 1-t$. 

We can also see that
\begin{equation} \label{refplusfield-exa2}
 \Refl(q) \subsetneq \orth(q), \quad \SO(q) \setminus \Refl(q) \ne \emptyset. 
\end{equation} 
Indeed, the  map defined by $(e_1, f_1, e_2, f_2) \mapsto (f_2, e_2, f_1, e_1)$ is an orthogonal transformation of $(M,q)$, even more, it lies in $\SO(q)$ because it fixes $t$. Clearly, $g\not\in \Refl(q)$. Of course, 
\eqref{refplusfield-exa2} is obvious from \eqref{refplusfield-exa1} and the formula for the order of $\orth(q)$, see e.g.  \cite[(13.3)]{Kneser}, 
which yields $|\orth(q)| = 72$ and implies $|\SO(q)| = 36$ by \eqref{sog-d3}. 

\comments{(2025-05-28) The remaining part of this section in the previous version has been moved to three places: the cohomological results to \S\ref{sec:consequences}, the (problematic) results in the semilocal case to the new section\S\ref{sec:trans-semi}, and the results on isomorphisms or orthogonal group schemes to appendix~\ref{sec:iso-orth-group}. In this way, they can easily be removed if we decide not to publish some of the material. }

\subsection{Orthogonal group scheme $\uO(q)$.}\label{orthsc}
Let $(M,q)$ be a quadratic module over $R$. The $R$--group functor
$\underline{\orth(q)}$, assigning to $S\in \Ralg$ the group $  \underline{\orth(q)}(S) = \orth(q_S)$ is represented by
an affine finitely presented $R$--group scheme denoted $\uO(q)$. It is the special case $S=\Spec(R)$ of the $S$--group scheme $\uO(q)$ associated in \cite[Def.~4.1.0.2]{CF} or in \cite[page 364]{Co1} with a quadratic module $(M,q)$ over a scheme $S$. Below we list some facts that we will use later. \sm

\begin{inparaenum}[(a)]
\item ({\em Centre}) \label{orthsc-zen} 
Suppose $M$ is faithfully projective.
Then the map
$s \mapsto s\Id_{M\ot S}$ for $s\in \bmu_{2,R}(S)$, $S\in \Ralg$, induces a monomorphism
\begin{equation}\label{orthsc-zen1}
z_M \co \bmu_{2,R} \longto \uO(q),
\end{equation}
of $R$--group schemes. If $(M,q)$ is faithfully quadratic space, $z$ is an isomorphism onto the (schematic) centre of $\uO(q)$ (the latter claim because $\GG_m$ is the centre of $\uGL(M)$), or see \cite[Cor.~C.3.9]{Co1}.
By \eqref{refle-one1}, it is an isomorphism if $M$ has constant rank $1$.
\sm

\item \label{orthsc-a} ({\em Direct products}) Suppose $R=R_1 \times \cdots \times R_n$ is a direct product of rings. A quadratic module $(M,q)$ uniquely decomposes into the direct product $(M,q) = (M_1, q_1) \times \cdots \times (M_n,q_n)$ of quadratic $R_i$--modules $(M_i, q_i)$, \ref{qf}\eqref{qf-redc}. Moreover, by \eqref{ortgr-bas-b1}, orthogonal groups respect this decomposition. Hence
\begin{equation} \label{orthsc0} \begin{split}
  &\uO(q)(R_1 \times \cdots \times R_n) =
    \uO(q)(R_1) \times \cdots \times \uO(q)(R_n)
 \\ & \quad = \orth(q_{R_1}) \times \cdots \times \orth(q_{R_n})
   = \uO(q_1)(R_1) \times \cdots \times \uO(q_n)(R_n),
   \end{split} \end{equation}
where we consider $\uO(q_i)$ as an $R_i$--group scheme.
In particular, we will apply this reduction in case $R= R_1 \times \cdots
\times R_n$ corresponds to the rank decomposition of $M$.

Let $p_i \co R \to R_i$ be the canonical projection, view $R_i$ as an $R$--algebra and let $p_{i*}(M_i, q_i)$ be the quadratic $R$--module obtained from $(M_i, q_i)$ by restricting scalars to $R$ via $p_i$. If $q$ is nonsingular or regular, each $p_{i*}(M_i, q_i)$ is nonsingular or regular respectively. Letting $p_{i*}\big( \uO(M_i, q_i)\big)$ be the base ring restriction of the $R_i$--group scheme $\uO(M_i, q_i)$ to an $R$--group scheme, 
the decomposition \eqref{orthsc0} leads to an isomorphism of $R$--group schemes
\begin{equation}\label{orthsc1}
 \uO(M,q) \cong p_{1*}\big(\uO(M_1, q_1)\big) \times \cdots \times
             p_{n*}\big( \uO(M_n, q_n)\big).
\end{equation}

\item \label{orthsc-det} ({\em Determinants}) For arbitrary $q$ we know from \ref{qf}\eqref{qf-detgr} that $\det(g) \in \mu_2(R)$ for any $g\in \orth(q)$. Assigning $\det(g)$ to $g\in \orth(q_T)$, $T\in \Ralg$, induces a homomorphism of $R$--group schemes
\begin{equation}\label{orthsc-det0}   
    \det \co \uO(q) \longto \bmu_{2,R} 
\end{equation}
where $\bmu_{2,R}$ is the $R$--group scheme with $\bmu_2(T) = \{ u\in T : u^2 = 1_T \}$. Following \cite[p.~364]{Co1} we define the $R$--group scheme 
\begin{equation}\label{orthsc-det1} 
    \uSO'(q) := \Ker (\det),  
\end{equation} 
and call $\uSO'(q)$ the {\em naive special orthogonal group\/}, not to be confused with the special orthogonal group of \ref{sogsc}. \sm 

\item \label{orthsc-d} ({\em Smoothness})    
    It is instructive to recall the following special case of \cite[Thm.~C.1.5]{Co1}). Let    $(M, q)$ be a quadratic space of constant positive rank $n$. Then $\uO(q)$
    is smooth if and only if either $n$ is even or $n$ is odd and $2\in R\ti$. 
    The naive special orthogonal group $\uSO'(q)$ is smooth if either $n$ is odd or $n$ is even and $2\in R\ti$. 
\end{inparaenum}

%

%

\subsection{The Dickson homomorphism}\label{dickhom}
Let $(M,q)$ be a faithful quadratic $R$--space and let $\Dis(q)$ be its discriminant algebra (\ref{qfdi}). 
We denote by $\uAut(\Dis(q))$ the automorphism group scheme representing the $R$--group functor $S \mapsto \Aut(\Dis(q)_S)= \Aut(\Dis(q_S))$. The group homomorphism $\Dis \co \orth(q) \to \Aut(\Dis(q))$ of \eqref{sog-dishom}
induces a homomorphism of $R$--group schemes
\begin{equation}  \label{dickhom1}
\Dis \co \uO(q) \longto \uAut\big(\Dis(q)\big)
\end{equation}
Let $(\ZZ/2\ZZ)_R$ be the constant $R$--group scheme associated with the abstract group $\ZZ/2\ZZ$. Its $S$--points, $S\in \Ralg$, can be identified with the idempotents of $S$, cf.\ \ref{sog}\eqref{sog-d}.
The map $\psi$ of \eqref{sog-d00} 
gives rise to a homomorphism of $R$--group schemes
\begin{equation}\label{dickhom0}
   \psi \co (\ZZ/2\ZZ)_R \to \uAut\big( \Dis(q) \big),
\end{equation}
which is an isomorphism if $q$ is regular. In this case (and only in this case) the {\em Dickson homomorphism $\Di$} is defined as $\psi\me \circ \Dis$, thus rendering the diagram \eqref{dickhom11} commutative
\begin{equation} \label{dickhom11} \vcenter{ \xymatrix{
     & \uAut\big( \Dis(q) \big) \\
  \uO(q) \ar[rr]^{\Di} \ar[ur]^{\Dis} && (\ZZ/2\ZZ)_R \ar[ul]_\psi^\cong
} } .
\end{equation}
Finally, we have a canonical $R$--group scheme homomorphism
\begin{equation}\label{dickhomc}
\chi \co (\ZZ/2\ZZ)_R \longto \bmu_{2,R}
\end{equation}
sending the idempotent $\eps \in S$ to the element $1-2 \eps\in \bmu_2(S)$. The map $\chi$ is an isomorphism if $2\in R\ti$. 
These maps fit into the commutative diagram
\begin{equation} \label{dickhom2}\vcenter{
\xymatrix@C=40pt{ & \uAut\big( \Dis(q)\big)\\
 \uO(q) \ar@{-->}[rr]^{\Di} \ar[ur]^{\Dis} \ar[dr]_\det && (\ZZ/2\ZZ)_R \ar[ul]_\psi
   \ar[dl]^{\chi}
 \\ &\bmu_{2,R}
}
}\end{equation}
in which, we recall, the dashed arrow $\Di$ is only defined if $q$ is regular, and 
$\det$ is the homomorphism \eqref{orthsc-det0}. Commutativity of the lower triangle is for example proven in \cite[Cor.~C.3.2]{Co1} and \cite[IV, (5.1.2)]{K}, or \cite[III, 5.2.9]{DG} for split forms. \sm 

Let $(M,q) = (M_1, q_1) \perp (M_2, q_2)$ be the orthogonal sum of two regular quadratic modules. We then know from \ref{bfLG}\eqref{qfba-sum} that $q$ is regular too. Moreover,  the diagram below commutes 
\begin{equation}\label{dickm-dir1} \vcenter{
 \xymatrix@C=50pt{
\uO(q_1) \times_R \uO(q_2) \ar[d]_{\inc} \ar[r]^{{\Di \times \Di }}&
(\ZZ/2\ZZ)_R \oplus (\ZZ/2\ZZ)_R \ar[d]^{\sum} \\
\uO(q)  \ar[r]^{\Di }& (\ZZ/2\ZZ)_R
}
}\end{equation} 
where $\inc$ is the natural closed immersion and $\sum$ is the sum map. Indeed, this can be proven over fields $F$. Moreover, since Dickson maps are stable
under base change, we may assume $|F|>2$ which allows us to apply the
Cartan-Dieudonn\'e-Kneser Theorem \ref{CDK}. Thus, it suffices to check
commutativity for $g_1=(s_1, \Id)$ and $g_2=(\Id, s_2)$ where $s_i$ are reflections
of $M_i$, $i=1,2$. Since $g_1$ and $g_2$ are reflections too, commutativity follows from \eqref{refso2}.

\subsection{Special orthogonal group scheme $\uSO(q)$} \label{sogsc}
Let $(M,q)$ be a faithful quadratic $R$--space and let $\Dis  \co \uO(q) \to \uAut\big(\Dis(q)\big)$ be the $R$--group homomorphism \eqref{dickhom1}. We define the $R$--group scheme $\uSO(q)$ as its kernel,
\begin{equation} \label{sogsc0} \uSO(q) = \Ker( \Dis). \end{equation}
Thus, for $S\in \Ralg$ we have
\[ \uSO(q)(S) = \{ g\in \orth(q_S) : \Dis(g) = \Id_{\Dis(q_S)}\} = \SO(q_S)\]
where we used \eqref{sog1} for the second equality. 
 
\begin{inparaenum}[(a)]
  \item \label{sogsc-dp} ({\em Direct products}) Let $R= R_1 \times \cdots \times R_n$ be a direct product of rings, and let $(M,q) = (M_1, q_1) \times \cdots \times (M_n, q_n)$
be the corresponding decomposition into a direct product of $R_i$--modules.
Analogously to \eqref{orthsc0} and keeping in mind \eqref{sog2}, we obtain
\begin{equation}\label{sogsc1} \begin{split}
  &\uSO(q)(R_1 \times \cdots \times R_n) =
    \uSO(q)(R_1) \times \cdots \times \uSO(q)(R_n)
 \\ & \quad = \SO(q_{R_1}) \times \cdots \times \SO(q_{R_n})
 \\ & \quad    = \uSO(q_1)(R_1) \times \cdots \times \uSO(q_n)(R_n),
   \end{split} \end{equation}
where we view $\uSO(q_i)$ as an $R_i$--group scheme. Moreover,  using the notation of \ref{orthsc}\eqref{orthsc-a}, we have the decomposition
\begin{equation}\label{sogsc11}
 \uSO(M,q) \cong p_{1*}\big(\uSO(M_1, q_1)\big) \times \cdots \times
             p_{n*}\big( \uSO(M_n, q_n)\big)
\end{equation}
into a direct product of $R$--group schemes, analogous to \eqref{orthsc1}.

Applying \eqref{sogsc11} to the rank decomposition of $(M,q)$ and using the identification of the structure of $\uSO(M,q)$ in \eqref{sogsc-odd} and \eqref{sogsc-even} below, we find that $\uSO(q)$ is a direct product of reductive (semisimple if $\rank M \ge 3$) $R$--group schemes, hence itself reductive.
\sm

\item \label{sogsc-det} ({\em $\uSO$ and determinants}) \label{sogsc-det} By \ref{sog}\eqref{sog-b},  {\em there exists a canonical monomorphism $\uSO(q) \to \Ker(\det)=\uSO'(q)$ of $R$--group schemes, which is an isomorphism 
\begin{equation} \label{sogsc22}
  \uSO(q) \simlgr \uSO'(q) \end{equation}
if $M$ has odd rank or if $2\in R\ti$.} \sm

\item ({\em $q$ regular}) If $q$ is regular, the map $\psi$ of \eqref{dickhom0}  is an isomorphism. Hence, from \eqref{dickhom11} we get the traditional definition of $\uSO(q)$ as kernel of the Dickson homomorphism,
      \begin{equation} \label{sogsc2}
           \uSO(q) = \Ker (\Di) .
       \end{equation}
Recall from \eqref{qf-ns2} 
that any nonsingular $q$ is regular if $2\in R\ti$. In this case, $\chi$ is obviously an isomorphism which explains \eqref{sogsc2}.

It follows from \eqref{sog-d3} that in case $(M,q)$ is regular of rank $\ge 2$, the sequence of $R$--group schemes
\begin{equation}  \label{sogsc2a}
1 \longto \uSO(q) \longto  \uO(q) \xrightarrow{\Di} \ZZ/2\ZZ \longto 1
\end{equation}
is exact in the flat topology \cite[IV, (5.2.2)]{K}, 
and even in the \'etale topology \cite[4.3.0.28]{CF}. 
\sm

\item ({\em Odd rank})\label{sogsc-odd}  Let $(M,q)$ be a quadratic space of odd rank.
Then the morphism $z$ of \eqref{orthsc-zen1} is a
section of $\det$. Hence, by \eqref{sogsc22}, see also \eqref{sog-e1},
\begin{equation}
  \label{orthsc-2} \uO(q) \cong \bmu_2 \times_R \uSO(q).
\end{equation}
If $M$ has constant rank $1$, then $\uSO(q) = \{\star \}$, and if $M$ has constant odd rank $2n+1\ge 3$, then $\uSO(q)$ is an adjoint semisimple $R$--group scheme of type ${\rm B}_n$ (${\rm B}_1 = {\rm A}_1$ for $n=1$), in particular it is a smooth affine group scheme (\cite[Prop.~C.3.10]{Co1}).
\sm

\item({\em Even rank}) \label{sogsc-even}  Let $(M,q)$ be a quadratic space of positive even rank. Then $q$ is regular by \ref{qf}\eqref{quadfoe} and the following hold.
\end{inparaenum}
\begin{enumerate}[label={\rm (\roman*)}]
 \item  \label{def-sosII} $\uSO(q)$ is the identity  component (\ref{neutrevLG}) of the smooth      $R$--group scheme $\uO(q)$; it is an open and closed subgroup scheme of $\uO(q)$ (\cite[Remark after Thm.~C.2.11]{Co1}). \sm

 \item \label{sogsc-even-ii} If $M$ has constant rank $2$,  then $\uSO(q)$ is a rank one torus.  Indeed, by descent, we are reduced to the hyperbolic case $q= xy$  where $\uSO(q)\cong \GG_m$, cf.~\eqref{sog-hyp1}.  
\sm

\item \label{sogsc-even-iii}  If $M$ has constant rank $2n\ge 4$, then $\uSO(q)$ is a semisimple $R$--group scheme of type ${\rm D}_n$ 
     (\cite[Prop.~C.3.10]{Co1}).%
 \sm

 \item \label{sogsc-even-iv} The centre of $\uSO(q)$ is $\bmu_2 \cong z_M(\bmu_2)\subset \uSO(q)$, see
     \eqref{orthsc-zen1}. Moreover,  $z_M(\bmu_2) = \GG_m \cap \uSO(q)= \GG_m
     \cap \uO(q)$.
\end{enumerate}
We point out that our definition of $\uSO(q)$ is does not coincide with the definition of $\uSO(q)$ in \cite[4.3.0.22]{CF}. Rather, our group $\uSO(q)$ is denoted $\uO_q^+$ in \cite{CF}, defined in \cite[4.3.0.27]{CF}. The notation $\uO^+(M,q)$ for our group $\uSO(q)$ is also employed in \cite{KMRT}, see \cite[page~358]{KMRT}.
\sm

\begin{inparaenum}[(a)]\setcounter{enumi}{5} \item \label{sogsc-cent} Let again $(M,q)$ be an arbitrary faithful quadratic space. We can now identify the image of the group monomorphism $z_M \co \bmu_{2,R} \to \uO(q)$ of \eqref{orthsc-zen1}: it represents the centralizer functor $\ulCent_{\ulO(q)}(\ulSO(q))$ in the sense of \cite[II, \S1, 3.4]{DG},
\begin{equation}   \label{sosc-cent0}
 z_m \co \bmu_{2,R}  \simlgr \uCent_{\uO(q)}(\uSO(q)).
\end{equation}
Indeed, by \cite[C.3.9]{Co1}, the functorial centre of $\uO(q)$ is represented by $\bmu_{2,R}$,  so that the sequence of monomorphisms $\bmu_{2,R} \to  \uCent_{\uO(q)}(\uSO(q)) \to \uCent(\uO(q))$
provides isomorphisms
\begin{equation} \label{sosc-cent1}
 \bmu_{2,R} \simlgr \uCent_{\uO(q)}(\uSO(q)) \simlgr \uCent(\uO(q)).
\end{equation}%
\lv{
OLD VERSION: the isomorphism \eqref{sosc-cent1} is implicit in the proof of \cite[C.3.9]{Co1}. Let us nevertheless give some hints. We can assume that $M$ has constant rank. In case it is odd, the isomorphism follows from \eqref{orthsc-2} and the fact that in this case $\uSO(q)$ is an adjoint group and therefore centre-less. In the even rank case, first observe that the image of $z_M$ is clearly a subfunctor of $\ulCent_{\ulO(q)}(\ulSO(q))$. To prove equality we can thus assume that $q$ and therefore also $\uSO(q)$ is split. One then verifies that the standard torus in $\uSO(q)$ is self-centralizing, which easily implies our claim.
}

\item\label{sogsc-ss} Let $(M,q)$ be a faithful quadratic $R$--space. Then 
$\uSO(q)$ is a reductive (even semisimple if $\rank_R M \ge 3$) $R$--group scheme with a central torus of rank $\le 1$. This follows from the discussion in \eqref{sogsc-odd} and \eqref{sogsc-even}. \sm

\item \label{sogsc-aut} By \eqref{sogsc-ss}, the automorphism group functor of $\uSO(q)$ is represented by a smooth affine $R$--group scheme $\uAut\big( \uSO(q)\big)$     
    (\cite[7.1.9]{Co1}, \cite[XXIV, Cor.~1.9]{SGA3}).
    Also, conjugation by elements of $\uO(q)$ gives rise to an automorphism of $\uSO(q)$, hence to a homomorphism
    $ \Int \co \uO(q) \to \uAut\big(\uSO(q))$
    whose kernel is $\bmu_{2,R}$ by \eqref{sosc-cent1}. If $M$ has constant rank $\ne 2$, the ensuing sequence
    \begin{equation}      \label{sogsc-aut1}
    1\longto \bmu_{2,R} \xrightarrow{\; z_M\; } \uO(q) \xrightarrow{\; \Int \;}
      \uAut\big(\uSO(q)\big) \longto  1
    \end{equation}
   is exact in the flat topology. Indeed, after identifying $\uO(q)/\bmu_{2,R} \equiv \mathbf{PGO(q)}$, this is \cite[C.3.13]{Co1}, a result of Dieudonn\'e if $R$ is a field of characteristic $\ne 2$. For $\rank_R M = 2$, see \cite[V; (2.6)]{K} or \cite[C.3.15]{Co1} and the remarks preceding it.
\end{inparaenum}

\ms

We can now show that Knebusch's definition of $\SO(q)$ in \ref{sog}\eqref{sog-c} coincides with our definition.

\begin{prop} \label{somax} Let $(M,q)$ be a {\em regular} faithful quadratic module.
Then $\uSO(q)$ is the identity component of $\uO(q)$, in
particular it is an
  open subscheme of\/ $\uO(q)$ and satisfies for all $A\in \Ralg$ that 
\begin{equation}\label{somax1} \begin{split}
    \uSO(q)(A) &= \{ g\in \uO(q)(A) : \;   g_{A/\m}\in \uSO(q)(A/\m)
    \\ &\hbox{ for all maximal  $\m\in \Spec(A)$}\}.
\end{split}\end{equation}

\end{prop}

\begin{proof} We apply the rank decomposition for $\uO(q)$ and $\uSO(q)$, see
\eqref{sogsc11}, and in this way reduce the claim to be proven to the constant rank case (recall from \ref{qf}\eqref{qf-redc} that the quadratic forms $q_i$ inherit regularity
from $q$). In the odd rank case, $2\in R\ti$ by \ref{qf}\eqref{qf-regLG}, whence 
the decomposition \eqref{orthsc-2} implies that $\uSO(q)$ is the identity component
of $\uO(q)$. In the even rank case, this follows from
\ref{sogsc}\ref{def-sosII}. Since the identity component of $\uO(q)$ is representable, namely by $\uSO(q)$, it is an open subscheme of $\uO(q)$ by \ref{neutrevLG}. The equation \eqref{somax1} therefore is a special case of \eqref{somax1g}. \end{proof}

\textbf{Example.} Proposition~\ref{somax} is not true, if $q$ is not regular. For example, let $R$ be the local ring $R=\ZZ/4\ZZ$ with maximal ideal $\m = 2\ZZ/4\ZZ$, and let $q $ be the nonsingular $1$-dimensional form $q= \lan 1_R \ran$. Then $- \Id_R \in \orth(q) \setminus \SO(q)$ because $-\Id_R \ne \Id_R$ and $\SO(q) = \{ \Id_R\}$ by \ref{sog}\eqref{sog-e}, but $(-\Id)_{R/\m} = \Id_{R/\m} \in \SO(q_{R/\m})$.


\begin{lem}\label{ogn} Let $(M,q)$ be a quadratic $R$--space, and let $N\subset M$ be a submodule satisfying $M = N \oplus N^\perp$. We note that both $q_N=q|_N$ and $q_{N^\perp} = q|_N^\perp$ are nonsingular by \eqref{qf-perp0}, and define the closed subgroup scheme $\uO(q)_N$ by requiring that the $T$--points, $T\in \Ralg$, are given by 
\[ \uO(q)_N (T) = \{ g\in \uO(q)(T\emph{}) : g(N_T) = N_T \text{ and } g|_{N_T^\perp} = \Id \}.   
\]   
We put $\uSO'(q)_N = \uSO'(q) \cap \uO(q)_N$ and $\uSO(q)_N = \uSO(q) \cap \uO(q)_N$. 
\sm 

\begin{inparaenum}[\rm (a)] \item\label{onga} The natural extension maps 
\begin{equation}\label{onga1} 
  \uO(q_N) \simlgr \uO(q)_N  \quad \text{and} \quad 
\uSO'(q_N) \simlgr \uSO'(q)_N,
\end{equation} 
given on $T$--points, $T\in \Ralg$, by $g\mapsto g \perp \Id_{N^\perp}$, are isomorphisms of $R$--group schemes. Their inverses are given by  restriction. \sm

\item\label{ongc} Let $N$ be a faithful $R$--module. 
Then the group scheme homomorphism $\uSO(q_N) \to \uSO'(q)$, $g \mapsto g \perp \Id_{N^\perp}$, factors through the monomorphism $\uSO(q) \to \uSO'(q)$ of {\rm \ref{sogsc}\eqref{sogsc-det}} and thus gives rise to a well--defined monomorphism 
of $R$--group schemes,
\begin{equation} \label{ongc1} 
\uSO(q_N) \hookrightarrow \uSO(q), \quad g \mapsto g \perp \Id_{N^\perp}.
\end{equation}%
\sm 

\item \label{ongb} The isomorphisms \eqref{onga1} restrict to an isomorphism of $R$--group schemes, 
\begin{equation} \label{ongb1}
 \uSO(q_N) \simlgr \uSO(q)_N.
\end{equation}    
\end{inparaenum}
\end{lem}

\begin{proof} The obvious proof of \eqref{onga} is left to the reader. \sm 

\eqref{ongc} The given map $\uSO(q_N) \to \uSO'(q)$ can be factored as follows
\[\xymatrix{
   \uSO(q_N) \ar[d]_{\io_N} \ar@{-->}[rr] && \uSO(q) \ar[d]^{\io} \\
   \uSO'(q_N) \ar[r]^\al_\cong & \uSO'(q)_N \ar[r]^\be & \uSO'(q)
}\]
where $\io_N$ and $\io$ are the monomorphisms of \ref{sogsc}\eqref{sogsc-det}, where $\al$ is the isomorphism of \eqref{onga} and where $\be$ is the obvious closed immersion. Since all these maps are in particular monomorphisms, the claim follows whenever $\io$ is an isomorphism. By \eqref{sogsc22}, this is the case if $2\in R\ti$ or if $\rank M$ is odd. To proceed, we use the rank decomposition of the quadratic spaces $(N,q_N)$ and $(N^\perp, q_{N^\perp})$. Without loss of generality, we can therefore assume that both $N$ and $N^\perp$ have constant positive rank and discuss the possible parities of these ranks. 

- $(\rank N, \rank N^\perp) = (\text{even}, \text{odd})$ or $(\text{odd}, \text{even})$:  Then $M$ has odd rank, and we are done by \eqref{sogsc22}.

- $(\rank N, \rank N^\perp) = (\text{odd}, \text{odd})$: Then $M$ has even rank, 
the quadratic form is regular by \ref{qf}\eqref{quadfoe}, which, by \ref{bfLG}\eqref{qfba-sum}, implies that both $q_N$ and $q_{N^\perp}$ are regular too. But then $2\in R\ti$ by \ref{bfLG}\eqref{bfLG-odd}, and we are again done. 

- $(\rank N, \rank N^\perp) = (\text{even}, \text{even})$:  In this case, we link the exact sequence \eqref{sogsc2a} for $q_N$ and $q$ with the restriction of the commutative diagram \eqref{dickm-dir1} to $\uSO(q_N) \times \{1\}$. Thus we get the diagram 
\[\xymatrix{
   1 \ar[r] & \uSO(q_N) \ar@{-->}[d] \ar[r] & \uO(q_N) \ar[r] \ar[d]^f & (\ZZ/2\ZZ)_R \ar[r] \ar@{=}[d]    &1 
  \\
 1 \ar[r] & \uSO(q) \ar[r] & \uO(q) \ar[r] & (\ZZ/2\ZZ)_R \ar[r] &1 
   }\]
in which the right square commutes and $f$ is a closed immersion onto $\uO(q)_N$ by \eqref{onga}. It follows that the restriction of $f$ maps $\uSO(q_N)$ to $\uSO(q) \cap \uO(q)_N = \uSO(q)_N$, and we are again done. 
\sm 

\eqref{ongb} follows from \eqref{ongc} and the definition of $\uSO(q)_N$. \end{proof}

\ms

We can now characterize isotropic quadratic spaces in terms of isotropy of the associated special orthogonal group. We refer the reader to \ref{rgs} for some background on isotropic or reducible reductive groups, and to \ref{quadint}, \ref{hrq} for the concept of higher rank quadrics.       

\begin{prop} \label{prop_isotropic} Let $(M,q)$ be a faithful quadratic $R$--space. Then the following are equivalent:

\begin{enumerate}[label={\rm (\roman*)}]

\item \label{prop_isotropic_i}  $(M,q)$ is isotropic.  

\item \label{prop_isotropic_ii}  There exists 
a locally constant function $\nu : \Spec(R) \to \ZZ$
with $1 \leq  \nu  \leq \frac{\mathrm{rank}(M)}{2}$
and such that the higher rank quadric $\uQ_\nu(q)$ admits an $R$--point. 

\item  \label{prop_isotropic_iii}  $\uSO(q)$ is an isotropic reductive $R$--group scheme.

\end{enumerate}
Furthermore, if $M$ is of rank $\geq 3$, this is also  equivalent to the following
\begin{enumerate} [label={\rm (\roman*)}]\setcounter{enumi}{3}
 
\item  \label{prop_isotropic_iv}  $\uSO(q)$ is reducible, i.e., admits an 
everywhere proper parabolic $R$--subgroup scheme.
\end{enumerate}
\end{prop}

\begin{proof} \ref{prop_isotropic_i} $\iff$ \ref{prop_isotropic_ii}: 
By definition in \ref{isotrop}, the quadratic space $(M,q)$ is isotropic if and only if $(M, q)$ admits a totally isotropic, faithful direct summand $U$.  The rank function $\nu$ of $U$ satisfies $1 \leq \nu \leq \frac{\mathrm{rank}(M)}{2}$ in view of Proposition \ref{quadrepII}, i.e., $U$ defines an $R$--point of the higher rank quadric $\uQ_\nu$.  Thus, \ref{prop_isotropic_i} $\implies$ \ref{prop_isotropic_ii}. The converse is immediate from the definition of $\uQ_\nu$. The implication \ref{prop_isotropic_i} $\implies$ \ref{prop_isotropic_iii} follows from Lemma~\ref{gmu}\eqref{gmu-a} below, taking into account  Proposition~{\rm \ref{quadrepII}} and \ref{hps}: given a totally isotropic and complemented submodule $U\subset M$, there exists a submodule $V$ such that $(U,V)$ form a hyperbolic pair of submodules and $M= (U \oplus V) \perp (U \oplus V)^\perp $ with $U\oplus V \cong \HH(U)$. \sm 

 
\ref{prop_isotropic_iii} $\implies$ \ref{prop_isotropic_i}: 
We first reduce to the case of $(M,q)$ being a quadratic space of constant rank. To do so, we use the standard rank decomposition \ref{qf}\eqref{qf-redc} of $(M,q)$ and the corresponding decomposition \eqref{sogsc1} of $\uSO(q)$. 
By assumption, there exists  a monomorphism $\GG_{m,R} \hookrightarrow   \uSO(q) \subset \uGL(M)$. It gives rise to monomorphisms $\GG_{m,R_i} \hookrightarrow \uSO(q_{R_i})$ for every $i$ (notation of \eqref{sogsc1}). 

According to \cite[I, 4.7]{SGA3}, we have a decomposition in eigenspaces
$M= \bigoplus_{n \in \ZZ} \, M_n$ 
where $\GG_{m,R}$ acts by isometries on each $M_n$ through the character 
$\chi_n(t)=t^n$. Let $m_n \in M_n$, $n\ne 0$. Then  $q(m_n)=q(t \, . \, m_n)  = q(t^n  m_n)=  t^{2n} \, q(m_n) \in R[t,t^{-1}]$, so that $q(m_n)=0$. Hence every $M_n$, $n\ne 0$, is a totally isotropic, complemented submodule. 
Since the rank of a finite projective $R$--module is locally constant, for each $n\in \ZZ$, $n\ne 0$, the set 
\begin{align*} U_n &= \{ s\in \Spec(R): \rank M_{n, \ka(s)} \ge 1 \} 
 \\   &= \textstyle \bigcup_{i=1, \ldots, d} \, \{s\in \Spec(R) : \rank M_{n, \ka(s) } = i \}
\end{align*} 
is a clopen subset of $\Spec(R)$. The action of $\GG_{m,R}$ being faithful, the $U_n$'s cover $\Spec(R)$. But $\Spec(R)$ is quasi-compact, so that finitely many clopen subsets $U_{n_1} = \Spec(R_1), \ldots, U_{n_c} = \Spec(R_c)$ cover $\Spec(R)$. Without loss of generality we can then assume that $U_{n_1} = \Spec(R)$. Thus,  $M_{n_1}$ is a totally isotropic, complemented submodule of positive rank. \sm 

We now assume $\rank M \ge 3$. By \ref{sogsc}\eqref{sogsc-ss}, the $R$-group scheme $\uSO(q)$ is then semisimple. Hence, by \ref{rgs}\eqref{rgs-aa}, the group scheme $\uSO(q)$ is isotropic if and only if it is reducible. This proves  \ref{prop_isotropic_iii} $\iff$ \ref{prop_isotropic_iv}. 
\end{proof}

\comments{(2025-09-14) We only use \ref{gmu}\eqref{gmu-a} in this section. We use \ref{gmu}\eqref{gmu-d} in the Springer section. I felt it is best to put it here, together with \ref{gmu}\eqref{gmu-b} and \ref{gmu}\eqref{gmu-c} since the "technic" (matrices) is used here. For now, we do not need \eqref{gmu-d} for general $U$ and $V$, so I did not do it.}

\begin{lem}[The cocharacter of a hyperbolic pair]\label{gmu} Let $(M,q)$ be a quadratic space, and let $(U,V)$ be a hyperbolic pair of submodules. Thus, by {\rm \ref{hps}}, we have $M= (U \oplus V) \perp M'$, $M'=(U \oplus V)^\perp $ and $U\oplus V \cong \HH(U)$. \sm 

\begin{inparaenum}[\rm (a)] \item \label{gmu-a}
Then the map, given on $T$--points, $T\in \Ralg$, as  
\begin{equation} \label{gmu1} 
    t \cdot (u,v, m) =  (t u, \, t\me v, m' ),
\end{equation} 
where $u \in U_T$, $v\in V_T$, $m' \in M'_T$, defines a monomorphism of $R$--group schemes
\[ \la_{U,V} \co \GG_{m,R} \to \uSO(q), 
\]
which factors through the canonical monomorphism $\uSO( U \oplus V) \to \uSO(q)$. \sm 

\item \label{gmu-b} The parabolic subgroup $P = \rmP_{\uSO(q)}(\la_{U,V})$ of $\uSO(q)$, \eqref{pare-1}, coincides with the stabilizer $\Stab_{\uSO(q)}(U)$ of $U$ in $\uSO(q)$. \sm

\item\label{gmu-levi} The action of $\uGL(U)$ on $\uW(M)$ by $g \cdot (u,v, m) = (gu, {^t g}\me v, m) $ gives rise to a closed immersion 
    \begin{equation} \label{gmu-levi1}
     \uGL(U) \times_R \uSO(q|_{M'} )\hookrightarrow \uSO(q) 
    \end{equation}
    whose image is the centralizer of $\la_{U,V}$ in $\uSO(q)$ and hence a Levi subgroup of the parabolic subgroup $P$ of {\rm \eqref{gmu-b}}. \sm

\item \label{gmu-c} The stabilizer $P' = \Stab_{\uSO(q)}(V)$ is a parabolic subgroup opposite to $P$. If $R$ is an LG ring and $Q$ is a parabolic subgroup of $\uSO(q)$ opposite to $P$, then there exists $V'\subset M$ such that $(U,V')$ is a hyperbolic pair of submodules and $Q=\Stab_{\uSO(q)}(V')$. 
    \sm      

\item \label{gmu-d} Suppose $U=Ru$ and $V= Rv$ for a hyperbolic pair $(u,v)$ of vectors in $M$. Then the unipotent radical $\rad^u(P)$ of $P$ is isomorphic to $\uW(M')$. The $T$--points $\big(\rad^u(P)\big)(T)$, $T\in \Ralg$, are precisely those endomorphisms of $M_T$ for which there exists $g'\in M'_T$ such  that 
\begin{equation}\label{gmu-d1}   \begin{split}
    &v_T \mapsto v_T - q(g') u_T + g', \\ 
    &u_T  \mapsto u_T, \quad    m' \mapsto m' - b_q(g',m')\, u_T
\end{split}\quad .\end{equation} 
\end{inparaenum}
\end{lem}


\begin{proof} We start by proving part of \eqref{gmu-levi1}: the given action of $\uGL(U)$ induces a monomorphism $\al \co \uGL(U) \to \uSO(q)$. 

It is immediate that $\uGL(U)$ acts by orthogonal transformations: for $T$--points $(u,v)$ of $\uW(U \oplus V)$ we have $q( gu + ({^t g}\me) v) = b_q(gu, ({^t g}{}\me) v) = b_q(g\me g u,v) = q(u + v)$. We thus get a monomorphism $\al \co \uGL(U) \to \uO(q)$, which factors through $\uO\big( \HH(U)\big) \to \uO(q)$. We claim that $\al$ factors through the homomorphism $\uSO( \HH(U)\big) \to \uSO(q)$, which is a monomorphism by Lemma~\ref{ogn}\eqref{ongc}. It suffices to consider the case that $M = \HH(U)$ is hyperbolic, in particular regular. 
Then Proposition~\ref{somax} further reduces the proof to showing that the map $\al (g)$, $g\in \uGL(M)(k)$, lies in $\SO(q_k)$ for $k$ a field in $\Ralg$. But by the same result, $\uSO(q)$ is the identity component of $\uO(q)$. Since $\uGL(M)$ is connected, our claim follows. 

\lv{
OLD proof: 
By definition of $\SO(q)$ in \ref{sog}, this means that $\al(g)$ fixes the discriminant algebra $\Dis(q_k)$. Let $(e_1, \ldots, e_n)$ be a basis of $U$ and let $(f_1, \ldots, f_n)$ be the basis of $V$,  which is the image of the dual basis of $e_1, \ldots, e_n$ under the isomorphism $V \simlgr  U^*$, $v \mapsto b_q(-, v)$. Then $1_k$ and $z = e_1 \cdots  e_n \cdot f_n \cdots f_n$ form a basis of $\Dis(q_k)$ by \cite[(7.4)]{Kneser}. One can show that $z$ is fixed under the action of $\al(g)$. So $\al(g)\in \SO(q_k)$.}
\sm

\eqref{gmu-a} The map $\la_{U,V}$ is the composition of the canonical monomorphism $\GG_{m,R} \to \uGL()$ and the monomorphism $\al \co \uGL(U) \to \uSO(q)$ above. 
\sm 

\comments{(2025-09-14) The proof of \eqref{gmu-b} is inspired by \cite[C.7]{GN-LG} (whose published version has changed the matrix below to an array.}

\eqref{gmu-b} We write $g\in \uGL(M)(T)$ as a matrix $g=(g_{ij})_{1\le i,j \le 3}$ of endomorphisms with respect to the family $(U_T, V_T, M'_T)$ of submodules of $M_T$. Then 
\[
  \la(t) \, (g_{ij}) \, \la(t)\me = \begin{pmatrix}
     g_{11} & t^2 g_{12} & t g_{13} \\ t^{-2} g_{21} & g_{22} & t\me g_{23} 
     \\ t\me g_{31} & t g_{32} & g_{33} 
  \end{pmatrix} 
\]
for $t\in T$. Thus, the $T$--points of $\rmP_{\uGL(M)}(\la)$ are those $g$ satisfying $g_{21} = g_{31} = g_{23} = 0$. The condition that $g$ be invertible is equivalent to $g_{11}$, $g_{22}$ and $g_{33}$ be invertible. On the other side, matrices of the form 
\begin{equation}\label{gmu-b1}\begin{pmatrix}
     g_{11} &  g_{12} &  g_{13} \\ 0 & g_{22} & 0 \\ 0  &  g_{32} & g_{33} 
  \end{pmatrix} 
\end{equation}
are exactly the matrices in $\uGL(M)(T)$ stabilizing $U_T$ and $U_T^\perp = U_T \oplus M'_T$. Therefore 
\[ \Stab_{\uGL(M)}(U) = \Stab_{\uGL(M)}(U) \cap \Stab_{\uGL(M)}(U^\perp). \] 
Since an orthogonal transformation stabilizes $U$ if and only if it stabilizes $U$ and $U^\perp$, it follows that $P_{\uO(q)}(\la) = \Stab_{\uO(q)}(U)$ which implies \eqref{gmu-b}. 
\sm

\eqref{gmu-levi} It is a general fact that the centralizer of a cocharacter in a reductive group scheme is a Levi subgroup of the parabolic subgroup given by the cocharacter, \ref{pare}\eqref{pare-b}. Obviously, the $T$--points of the centralizer  $\Cent_{\uGL(M)}(\la_{U,V})$ consists of the matrices \eqref{gmu-b1} with $g_{ij} = 0$ for $i\ne j$. Such a matrix is orthogonal if and only if $g_{33}$ is orthogonal and $g_{22} = {^t g_{11}}{}\me$. Hence with $q' = q|_{M'}$ we get 
\begin{align*}
\Cent_{\uO(q)}(\la_{U,V})  &=  \al\big(\uGL(U)\big) \times_R \uO(q'),  \\
 \Cent_{\uSO'(q)}(\la_{U,V})  &=  \al\big(\uGL(U)\big) \times_R \uSO'(q'),\\
 \Cent_{\uSO(q)}(\la_{U,V})  &=  \al\big(\uGL(U)\big) \times_R \big( \uSO'(q') \cap \uSO(q)\big)  \\ 
   & \cong \uGL(U)  \times_R \uSO(q') ,
\end{align*}  
where we identify $\uO(q') = \uO(q)_{M'}$, $\uSO'(q') = \uSO'(q)_{M'}$ by \eqref{onga1} and use  $\uSO(q') \cong \uSO(q)_{M'}$ in view of \eqref{ongb1}.  
\sm 

The first part of \eqref{gmu-c} follows from \eqref{gmu-a} and \ref{pare}\eqref{pare-b} since $P' = \rmP_{\uSO(q)}(\la\me)$ is a parabolic subgroup of $\uSO(q)$ opposite to $P$. If $R$ is an LG ring, any two parabolic subgroups opposite to $P$ are conjugate by an element of $\uSO(q)(R)$. This implies the second part of \eqref{gmu-c}.  \sm 

\eqref{gmu-d}  By general theory, $P \cap P'$ is a Levi subgroup of $P$ and $P'$. Its $T$--points are given by diagonal matrices, i.e., $g_{12} = g_{13} = g_{32} = 0$ in \eqref{gmu-b1}. The $T$--points of $\rad^u(P)$ are those matrices in $\uSO(q)(T)$ which are the identity on $U$, $V$ and $M'$. The conditions that such a matrix lies in $\uO(q)(T)$ are  
\begin{enumerate}[label={\rm (\roman*)}] 
\item \label{gmu-di} $b_q(v, g_{13} m') + b_q(g_{32} v, m') = 0$ for all $v\in V_T$ and $m'\in M'_T$, and 

\item\label{gmu-dii} $q(g_{32}v) + b_q(g_{12} v, v) = 0$ for all $v\in V_T$.  
\end{enumerate}
Specializing $U = Ru$ and $V=Rv$ for a hyperbolic pair $(u,v)$ and writing 
\begin{align*}
  & g_{12} \co R v \to Ru, \quad v \mapsto ru \text{ for some $r\in R$}, \\
  & g_{32} \co Rv \to M', \quad v \mapsto g' \text{ for some $g'\in M'$}, \\
  & g_{13} \co M' \to Ru, \quad m' \mapsto \vphi(m') u \text{ for some $\vphi \in M^{\prime\, *}$},
\end{align*}
we obtain that \ref{gmu-di} and \ref{gmu-dii} are equivalent to 
$\vphi(m') + b_q(g',m') = 0$ for all $m' \in M'_T$ and $q(g') + r = 0$, which implies \eqref{gmu-d}. \end{proof}
\sm 

\textbf{Remarks.} By \ref{quadrepII}, every totally isotropic and complemented submodule $U\subset M$ is part of a hyperbolic pair of submodules. We will give another proof of \ref{gmu}\eqref{gmu-b} in \ref{prop_hig}\ref{prop_higa}. We have already established the special case $U=Ru$ of \ref{gmu}\eqref{gmu-b} in \cite[Prop.~A.5]{GN-Sp}. That $\rad^u(P) \cong \uW(E)$ in \ref{gmu}\eqref{gmu-d}  for some finite projective $R$--module $E$ is a general fact, \ref{pare}\eqref{pare-urad}.

\comments{(2026-04-29) Subsections \ref{mgs}, \ref{spinimm} and \ref{abc} were previously part of section~\ref{sec:consequences}. We use \ref{abs2} in the new proof of Theorem~\ref{knex-c} that $\SO(q) \to \SO(q_{R/\Jac(R)})$ is surjective. 
}

\subsection{Group schemes $\uSpin(q)$ and $\mathbf{S\Ga}(q)$.}\label{mgs} Let $(M,q)$ be a faithful quadratic $R$--space. Besides the group schemes associated with $(M,q)$ in this section, we can associate two more group schemes with $(M,q)$, namely $\uSpin(q)$ and $\uSG(q)$. We summarize their construction below; details can be found in \cite{Bass-74}, \cite[IV, \S6]{K} and \cite[4.5]{CF}. \sm

Let $\Cli(q) = \Cli=  \Cli_0 \oplus \Cli_1$ be the Clifford algebra of $q$. It is a $(\ZZ/2\ZZ)$--graded $R$--algebra. We can and will identify $M$ with a submodule of $\Cli_1$. Any $x\in \Cli\ti_0 $ gives rise to the inner automorphism $i_x$ of $\Cli$, defined by $i_x(c) = x c x\me$ for $c\in \Cli$. It respects the grading of $\Cli$. The {\em special Clifford group\/} is defined as
\[ \SG(q) = \{ x\in \Cli_0\ti: i_x(M) \subset M \}
           = \{ x\in \Cli_0\ti: i_x(M) = M\}. \]
For $x\in \SG(q)$ we put 
\[ \spi (x) = i_x|_M. 
\]
Then $\spi(x) \in \orth(q)$, since $q\big( \spi(x)(m)\big) = \big( \spi(x)(m)\big)^2 = x m^2 x\me = x q(m) x\me \allowbreak = q(m)$. In fact, $\spi(x) \in \SO(q) = \Ker (\Dis)$ for $\Dis\co \orth(q) \to \Aut\big(\Dis(q)\big)$, cf.\ \eqref{sog1}, because
the discriminant algebra $\Dis(q)$ is the centralizer algebra of $\Cli_0$ in $\Cli$.
We have $\spi(x) = \Id_M \iff x\in \rmZ(\Cli)$, the centre of $\Cli$. By the structure theory of the Clifford algebra $\Cli$, we know  $\rmZ(\Cli) \cap \Cli_0 = R \cdot 1_{\Cli}$. Hence
\begin{equation} \label{mgs1}
  1 \longto R\ti \xrightarrow{\; \inc\; } \SG(q) \xrightarrow{\; \spi\; } \SO(q)
\end{equation}
is an exact sequence of groups. \sm

Let $c\mapsto \ol c$ be the {\em standard involution\/} of $\Cli$, uniquely determined by $\ol m = - m$ for $m\in M$. Any $x\in \SG(q)$ satisfies $xmx\me = - \ol {xmx\me} = \ol x\me m \ol x$, whence $(\ol x x) m (\ol x x)\me = m$, so that $\ol x x \in \rmZ(\Cli) \cap \Cli_0 = R \cdot 1_{\Cli}$. It is then immediate that
\begin{equation} \label{mgs11}
 \rmsn \co \SG(q) \to R\ti, \quad x \mapsto \ol x x 
\end{equation}
is a group homomorphism, called the {\em spinor norm.} We put
\[ \Spin(q) = \Ker(\rmsn),  \]
called the {\em spin group\/},  and thus get two more exact sequences,
\begin{equation} \label{mgs2}
  1 \longto \Spin(q) \xrightarrow{\; \inc\; } \SG(q) \xrightarrow{\; \rmsn \; }
 R\ti
\end{equation}
 and, by restricting \eqref{mgs1},
\begin{equation} \label{mgs2}
  1 \longto \bmu_2(R) \xrightarrow{\; \inc\; } \Spin(q) \xrightarrow{\spi|_\Spin}
  \SO(q).
\end{equation}

It is immediate that the constructions above are stable under base change, thus giving rise to group functors $\underline{\SG}(q)$ and $\underline{\Spin}(q)$. They are represented by affine finitely presented $R$--group schemes $\uSG(q)$ and $\uSpin(q)$. In fact, if $\rank M \ge 3$, then $\uSpin(q)$ is a semisimple simply connected $R$--group scheme \cite[6.4.0.34, 8.2.0.59]{CF}, and $\uSG(q)$ is a reductive group scheme.
\comments{ (2026-06-25) Deleted: "for example by \ref{redss}\eqref{redss-i}"

(2026-04-29)
The reference \ref{redss}\eqref{redss-i} above refers to the following fact: \sm 

{\tt 
Let $1 \to G' \to G \to G'' \to 1$ be an exact sequence of $S$--groups in the flat topology. If $G'$ and $G''$ are reductive (semisimple respectively), then so is $G$. } \sm 

Previously, this was in a review of reductive and semisimple groups.}

Similarly, the maps $\spi$ and $\rmsn$ respect base change and therefore induce homomorphisms of group schemes. In this way we arrive at the diagram
\begin{equation} \label{mgs3} \begin{split}
  \xymatrix@C=40pt{ & 1 \ar[d] & 1\ar[d] \\
     1 \ar[r] & \bmu_2 \ar[r]^{\inc\;\;  } \ar[d] & \uSpin(q) \ar[r]^{\spi|_\Spin} \ar[d]^{\inc}& \uSO(q) \ar[r] \ar@{=}[d]& 1 \\
   1 \ar[r]& \GG_m \ar[r]^{\inc} \ar[d]^{\times 2} & \uSG(q) \ar[r]^{\spi} \ar[d]^{\rmsn}& \uSO(q) \ar[r]& 1 \\
   & \GG_m \ar@{=}[r]\ar[d]  & \GG_m \ar[d] \\
   & 1 & 1
  }
\end{split}
\end{equation}
It is immediate from the constructions above that all three squares commute. Furthermore, it is shown in \cite[4.6.0.8, 4.6.0.9]{CF} that all rows and columns are exact in the flat topology.
\comments{(2021-02-13) The proof of exactness in the even rank case in  \cite[4.6.0.9]{CF} is sketchy. It reduces to surjectivity of $\uSG(q) \to \uSO(q)$ and then the latter is killed by a reference to \cite[IV, (6.2.3)]{K}. So it is perhaps worth to have a clean direct proof. }

We present a slightly different proof here. First, it is standard that the left column is exact in the flat topology. Hence, by commutativity of the diagram, so is the middle column. It now suffices to show that the first row is exact in the flat topology, since this will imply exactness of the second row. Moreover, by \eqref{mgs1}, we are left with proving that $\spi \co \uSpin(q) \to \uSO(q)$ is an epimorphism in the flat topology. For simple notation we show this for the $R$--points. Thus, let $g\in \uSO(q)(R)= \SO(q)$.

Let us first suppose $g\in \Refl^+(q)$, say $g=\rho_{m_1} \cdots \rho_{m_n}$. Then $x=m_1 \cdots m_n \in \SG(a)$ with $\spi(x) = g$ by Lemma~\ref{spinimm} below. Passing to an appropriate flat cover we may assume that all $q(m_i)$ are squares, say $q(m_i) = u_i^{-1}$. Then $y=(u_1 m_1) \cdots (u_n m_n) \in \Spin(q)$ because $\ol{u_i m_i}u_i m_i = u_i^2 q(m_i)$, and $\spi(y) = g$ because $\rho_{u_i m_i} = \rho_{m_i}$.

It remains to prove that, up to passing to a flat cover, we can always assume $g\in \Refl^+(q)$. But this follows from Corollary~\ref{cor-knex}.

\comments{(2026-04-29) Need to check if Corollary~\ref{cor-knex} is true, or if the reference here to that corollary can be avoided. See \cite[Prop.~4.3.0.25]{CF}, which shows: \sm 

{\tt $R$ local ring, $q= \lan 1 \ran \perp q'$ with $q'$ regular, then every $g\in \SO(q)$ is a product of reflections.} }
\ms

The following well-known Lemma~\ref{spinimm}
relates the image of the spinor norm  $\rmsn$ with the values $\rmD(q) = R\ti \cap q(M)$ of the quadratic form $q$, reviewed in \ref{dqd-ele}. Recall $\rmD(q)^{[d]}  = \rmD(q) \cdots \rmD(q)$ ($d$ factors) and $  \rmD(q)^{[\rm ev]} = \textstyle \bigcup_{0 < d \text{ even}} \Dqd$, a subgroup of $R\ti$ containing $R\ti{}^2$.


\comments{(2026-04-29) Lemma~\ref{spinimm} was previously part (a) of Lemma~\ref{spinim}. } 
\begin{lem} \label{spinimm} Let $(M,q)$ be a faithful quadratic $R$--space.
Let $m_1, \ldots, m_n \in M$ with $q(m_1), \ldots, q(m_n)  \in R\ti$ and $n$ even. Then $x= m_1 \cdots m_n \in \SG(q)$ and
 \begin{equation} \label{spinim1} \begin{split}
    \spi(m_1 \cdots m_n) &= \rho_{m_1} \cdots \rho_{m_n},  \\
    \rmsn(m_1\cdots  m_n) &= q(m_1) \cdots q(m_n) \in \rmD(q)^{[\rm ev]}.
  \end{split} \end{equation}
\end{lem}

\begin{proof} Viewed as element of $M \subset \Cli_1 \subset \Cli$ we have $m_1^2 = q(m_1) \in R\ti$. Hence $m_1$ is invertible in $\Cli$ with inverse $q(m_1)\me m_1$. It follows that $m_1 m_2 \in \Cli\ti_0$. We have seen in \eqref{reso0} that $i_{m_1}(m) = - \rho_{m_1}(m)$ holds for $m\in M$, implying $m_1 m_2 \in \SG(q)$ with $\spi(m_1 m_2) = \rho_{m_1} \rho_{m_2}$. Its spinor norm is $\rmsn(m_1 m_2) = \ol m_2 \ol m_1 m_1 m_2 = m_2 m_1^2 m_2 = q(m_1) q(m_2)$.
\end{proof}

\subsection{The spinor norm $\SN \co \uSO(q) \to H^1(R, \bmu_2)$} \label{abc} 
\begin{inparaenum}[(a)] \item \label{abc-a}
Since $\bmu_2$ is central in $\uSG(q)$, we can apply Lemma \ref{lem_snake}\ref{lem_snake4} to the diagram \eqref{mgs3} and get an exact diagram of pointed sets
\begin{equation}\label{diag_sn} \vcenter{
 \xymatrix@C=20pt{
 & \uSG(q)(R) \ar[d]^{\rmsn} \ar[r]^\spi & \uSO(q)(R)   \ar[d]^{\SN}
 \ar[r] & \Pic(R) \ar@{=}[d]  \\
R\ti \ar[r]^{\times 2} &R^\times \ar[r]^{\delta} \ar[d]& H^1\fppf(R, \bmu_2)  \ar[d] \ar[r] & \Pic(R) \\
& H^1\fppf( R, \uSpin(q)) \ar[r]^\sim & H^1\fppf( R, \uSpin(q))
}}
\end{equation}
where $\SN$ and $\del$ are the characteristic maps associated with the exact sequences $1 \to \bmu_2 \to \uSpin(q) \to \uSO(q) \to 1$ and $1\to \bmu_2 \to \GG_m \to \GG_m \to 1$ of \eqref{mgs3}. The upper left square in \eqref{diag_sn} anti-commutes and the upper right and lower square commute. The relations in \eqref{diag_sn} explain why the characteristic map $\SN$ is also called the spinor norm, \cite[p.~232]{K}. It should not be confused with the spinor norm $\rmsn$ of \eqref{mgs11}. Both $\delta$ and $\SN$ are group homomorphisms.

We identify $H^1\fppf( R, \bmu_2)$ with the abelian group $\Disc(R)$ of discriminant modules (\cite[III, \S3]{K}). Then $\de(u) = [ (R, \lan u \ran_q)]$ for $u\in R\ti$. Also we know from \eqref{spinim1}  and anti-commutativity of the left upper square in \eqref{diag_sn} that 
\begin{equation}\label{abc-a1}
\SN(\rho_{m_1} \cdots \rho_{m_n}) = \de\big( q(m_1) \cdots q(m_n)\big){}\me
\end{equation}
for $m_1, \ldots, m_n \in M$ with $q(m_i)\in R\ti$ and $n$ even. 
\sm

\item \label{abc-b} {\em We now assume $\Pic(R) = \{1\}$}. Hence (for example by considering the middle row of \eqref{diag_sn}) we have
\begin{equation} \label{abs2}
   H^1\fppf(R, \bmu_2) \cong R\ti / R\ti{}^2
 \end{equation}
and $\de$ can be identified with the homomorphism $u \mapsto u R\ti{}^2$. Also, $\Pic(R) = \{1\}$ implies that $\spi$ is surjective. Therefore \begin{equation} \label{abc3}
\Ima(\SN) = \de\big( \Ima(\rmsn)\big)= \Ima(\rmsn)/R\ti{}^2 \supset \rmD(q)^{[\rm ev]}/ R\ti{}^2,
\end{equation}
where the last inclusion follows from \eqref{abc-a1}. 
\sm

\item\label{abc-c}  {\em Let $R$ be a semilocal ring.} Hence $\Pic(R) = \{1\}$, so that \eqref{abc-b} applies. By \ref{spinim}\eqref{spinim-cii}  we know that  $\Ima(\rmsn) = \rmD(q)^{\rm [ev]}$, hence \eqref{abc3} can be improved to
\begin{equation}
  \label{abc-c1}    \Ima(\SN) = \rmD(q)^{[\rm ev]}/R\ti{}^2.
\end{equation}
\end{inparaenum}
%
%

\subsection{}\label{oddi-prep} As a preparation for the following Lemma~\ref{oddi},  we recall 
that the discriminant algebra $\Dis(q_o)$ of a quadratic $R$--space $(M_o, q_o)$ with $M_o$ of constant odd rank  is $(\ZZ/2\ZZ)$--graded,
\[ \Dis(q_o) = \Dis_0(q_o) \oplus \Dis_1(q_o), \quad \Dis_0(q_o)= R, \]
with $\Dis_1(q_o)$ being a discriminant module with respect to the multiplication $h_0 \co \Dis_1(q_o) \times \Dis_1(q_o) \to \Dis_0(q_o) = R$ of $\Dis(q_o)$, which is itself induced from the multiplication of the Clifford algebra $\Cli(q_o)$. If $(M_o, q_o)=(L,q_L)$ is a discriminant module, then $\Cli(q_L) = \Dis(q_L)$ with $\Cli_1(q_L) = \Dis_1(q_L) =  L$ as $R$--module and $h_L(\ell, \ell) = q_L(\ell)$.

We also recall that 
the even Clifford algebra of an orthogonal sum is $(\ZZ/2\ZZ)$--graded. In particular, for any quadratic $R$--module $(M,q)$ and discriminant module $(L, q_L)$ we have
\begin{equation}\label{oddi-prep1}  \Cli_0\big(q \perp (-q_L)\big) = \big(\Cli_0(q) \wdh \ot_R \, R \big)
    \oplus \big( \Cli_1(q) \wdh \ot_R \, L \big).
\end{equation}

\begin{lem}[Discriminant modules of tensor product forms {\cite[IV, (7.3.2), (7.3.3)]{K}}] \label{oddi} We use the setting of\/ {\rm~\ref{oddi-prep}:} $(M,q)$ is a quadratic $R$--module and $(L,q_L)$ is a discriminant module. Then the following hold. \sm

The map $M\ot_R L \to \Cli_1(q) \wdh \ot_R L$, $m\ot \ell \mapsto m \wdh \ot \ell$, induces an isomorphism
\begin{equation}  \label{oddi1}
 \Cli(q\ot q_L) \simlgr \Cli_0\big(q \perp (-q_L)\big)
\end{equation}
of $\ZZ/2\ZZ$--graded $R$--algebras.

If $(M_o,q_o)$ is a quadratic $R$--space of constant odd rank, the discriminant modules $\Dis_1(q_o\ot q_L)$ and $\Dis_1(q_o)$ of the quadratic forms $q_o\ot q_L$ and $q_o$ are related by the isometry
\begin{equation}  \label{oddi2}
 \Dis_1(q_o\ot q_L) \simlgr \Dis_1(q_o) \ot_R (L, q_L)
\end{equation}
of discriminant modules, obtained by restriction of \eqref{oddi1}.
\end{lem}
\lv{
\begin{proof}
  Since $(m\wdh \ot \ell)^2 = q(m) q_L(\ell)$ in $\Cli_0\big(q\perp (-q_L)\big)$, the universal property of Clifford algebras yields a unique homomorphism
  \[ \vphi \co \Cli(q \ot q_L) \to \Cli_0\big(q \perp (-q_L)\big)
  \]
  of $\ZZ/2\ZZ$--graded $R$--algebras, extending $m\ot \ell \mapsto m \wdh \ot \ell$. Because its image contains a generating set of $\Cli_0\big(q \perp (-q_L)\big)$, the homomorphism $\vphi$ is surjective. It is in fact an isomorphism because the rank functions on both sides coincide.

Both $M_o \ot L$ and $M_o$ have the same odd rank. Also, by definition of discriminant algebra and \ref{discralg}\eqref{discralg-d},
\[ \Dis(q_o\ot q_L) = \Cli(q_o\ot q_L)^{\Cli_0(q_o\ot q_L)}
= R \oplus \big(\Cli_1(q_o\ot q_L)^{\Cli_0(q_o\ot q_L)}\big). \]
Since $\vphi$ is $\ZZ/2\ZZ$--graded, it maps $\Dis_1(q_o \ot q_L)$ onto $(\Cli_1(q_o)\wdh \ot L) ^{\Cli_0(q_o)} = \Dis_1(q_o) \wdh \ot L$. Let $h_{q_o\ot q_L}$ and $h_{q_o}$ be the symmetric bilinear forms of the discriminant modules $\Dis_1(q_o\ot q_L)$ and $\Dis_1(q_o)$. Then, by definition of $\wdh \ot$, the product of elements of $\Dis_1(q_o) \wdh \ot L$ is
\[ (d_1 \ot \ell_1) \cdot (d_2\ot \ell_2) = - \big( (d_1 d_2)(-h_{q_L}(\ell_1, \ell_2))\big) = h_{q_o}(d_1, d_2) \, h_{q_L}(\ell_1, \ell_2) \]
which is the product of the discriminant algebra $\Dis(q_o\ot q_L)$,  implying \eqref{oddi2}. \end{proof}}

\subsection{Similitudes} \label{simi-defi} Given two quadratic $R$--modules $(M, q)$ and $(M', q')$, a {\em similitude\/} $(M, q) \to (M', q')$ is pair $(f,\mu)$ consisting of a bijective $R$--linear map $f \co M \to M'$ and $\mu \in R\ti$  such that
\[   q'\big(f(m)\big) = \mu\, q(m)
\]
holds for all $m\in M$, i.e., $f \co (M,\mu q) \simlgr (M', q')$ is an isometry.
We denote by $\GO(q, q')$ the set of similitudes $(M,q) \to (M',q')$. If a similitude $(M, q) \to (M',q')$ exists, we 
say that $(M,q)$ and $(M',q')$ are {\em similar\/} and call $\mu=
\mult(f)$ the {\em multiplier of $f$}. For example, an isometry is a similitude with multiplier $1$.

Below we list some easily established facts, where $f\co (M, q) \to (M',q')$ is a similitude with multiplier $\mu$. \sm

\begin{inparaenum}[(a)] \item \label{simi-defi-x} Let $(R, \lan \mu \ran_q)$ be the discriminant module associated with $\mu \in R\ti$.
Since the canonical isomorphism $M \ot_R R \simlgr M$ is an isometry
$\nu_u \co (M,q) \ot_R (R \lan u \ran_q) \simlgr (M,q)$, a similitude $(f,\mu)$ is essentially the same as an isometry $(M,q) \ot (R, \lan \mu \ran_q) \simlgr
(M', q')$:
\begin{equation} \label{simi-defi-1}\vcenter{
 \xymatrix{(M,q) \ot_R (R, \lan \mu \ran_q) \ar[dr]_{\wtl f}\ar[rr]^{\nu_\mu}_\cong
   && (M,\mu q)\ar[dl]^f  \\ & (M',q')}
}\end{equation}

\item\label{semi-defi-a} If  $\ka \in R\ti$, then $\ka f\co ( M, q) \to (M', q')$ is a similitude with multiplier $\ka^2\mu$. In particular,
    {\em $x\Id_M$ is a similitude of $(M,q)$ with multiplier $x^2$}. We have a central homomorphism
    \begin{equation}
      \label{semi-defi-a3} z'_M \co R\ti \to \GO(q), \quad x \mapsto x\Id_M.
    \end{equation}

\item\label{simi-difi-y} Let $g \co (M',q') \to (M'', q'')$  be a similitude with multiplier $\mu(g)$. Then $g \circ f$ is a similitude with  multiplier $\mu(g \circ f) = \mu(g) \mu(f)$. In particular, the inverse of $f$ is a similitude with multiplier $\mu(f\me) = \mu(f)\me$. It follows that the similitudes of the quadratic module $(M, q)$ form a group under composition, i.e., $(f,\mu) \cdot (g,\mu') = (g\circ f, \mu'\mu)$, denoted
    \[ \GO(q) =\GO(q,q), \]
and called the {\em orthogonal similitude group}.
The map
\begin{equation}  \label{semi-defi2}
 \mult \co G(q) \to R\ti, \quad (f,\mu) \mapsto \mu
\end{equation}
is a group homomorphism whose image we denote by
\begin{equation} \label{semi-defi22} \rmG(q) = \{ \mu \in R\ti : q \cong \mu q \}
\end{equation}
and call the group of {\em similarity factors of $q$}. Thus, by \eqref{semi-defi-a},
\begin{equation}\label{simi-difi-y1}
R\ti{}^2 \subset \rmG(q) \subset R\ti.
\end{equation}
Both extreme cases occur naturally: it is immediate from \ref{carhyp} that $G(q) = R\ti$ for any hyperbolic form, and we have
\begin{equation}\label{simi-difi-y3}
\rmG(q_o) = R\ti{}^2 \end{equation}
{\em for a quadratic $R$--space $(M_o,q_o)$ of constant odd rank.}

For the proof of \eqref{simi-difi-y3} note that any similitude of $(M_o, q_o)$ with multiplier $\mu$, i.e., by \eqref{simi-defi-x}, an isometry $q_o \ot \lan \mu \ran_q \simlgr q_o$, gives rise to an isometry of the associated discriminant modules $\Dis_1(q_o \ot \lan \mu \ran_q) \simlgr \Dis_1(q_o)$. But $\Dis_1(q_o \ot \lan \mu \ran_q) \cong \Dis_1(q_o) \ot (R, \lan \mu\ran_q)$ by \eqref{oddi2}. Thus, $\Dis_1(q_o) \cong \Dis_1(q_o) \ot (R, \lan  \mu\ran_q)$, which
implies $(R, \lan \mu\ran_q) \cong (R,\lan 1_R \ran_q)$ and then $\mu\in R\ti{}^2$.

We will say more about $\rmG(q)$ in 
the following subsections of \S\ref{sec:scharlau}, see \eqref{simi-defisch-b11} for an extension of the exact sequence
\[ 1 \longto \, \orth(q) \, \xrightarrow{\inc} \, \GO(q) \, \xrightarrow{\mult} \, \rmG(q)\,  \longto 1 \]
to group schemes.
\sm

\item\label{simi-def-aa} Assume $M$ is faithfully projective and $(M,q)$ is primitive, e.g. nonsingular. Then the multiplier of $f$ is unique. This follows for example from the criterion \eqref{quadco-aa1}. In this case, we will write $f$ instead of $(f,\mu)$.
    \sm

\item Since isometries preserve regularity and nonsingularity, the interpretation of similitudes in terms of isometries together with \eqref{tenssq1} 
    shows
 \begin{equation} \label{simi-difi-z3} \begin{split}
     \text{\em$(M,q)$ is regular} \quad &\iff \quad
    \text{\em $(M',q')$ is regular.}
   \\
    \text{\em $(M,q)$ is nonsingular} \quad &\iff \quad
      \text{\em $(M',q')$ is nonsingular.}
  \end{split}\end{equation}%
\end{inparaenum}

\subsection{The group schemes $\uGO(q)$, $\uPGO(q)$ and $\uGSO(q)$} \label{simi-defisch}
We compile some known facts regarding the similitude group scheme $\uGO(q)$ and related group schemes, following \cite[C.3]{Co1} and \cite[\S4.4]{CF}. But contrary to \cite{Co1} we will not consider line bundle-valued quadratic forms, and contrary to \cite{CF} we will not consider group schemes associated with quadratic pairs. On the other hand, we will allow quadratic forms of arbitrary rank and not restrict to regular quadratic forms as in \cite{CF}. Also, the group scheme $\uPGO^+$ of \cite{CF} is denoted $\uGSO(q)$ here. \sm

Let $(M,q)$ be a faithful quadratic $R$--space. The $R$--group functor $\underline{\GO}(q)$, assigning to $S\in \Ralg$ the group $\underline{\GO}(q_S)$ is represented by an affine $R$--group scheme
\[ \uGO(q).\] It is smooth by \cite[C.3.12]{Co1}, or see
\ref{simi-defisch-a-ii} and \ref{simi-defisch-b-iii} below.
The multiplier homomorphism of \ref{simi-defi}\eqref{simi-difi-y} extends to a homomorphism
\[  \mult \co \uGO(q) \to  \GG_m \]
of $R$--group schemes. We identify $\uO(q)$ with its kernel.

The central monomorphism $z_M\co \bmu_{2,R} \to \uO(q)$ of \eqref{orthsc-zen1} extends to a central monomorphism
\begin{equation}\label{simi-defisch1}
 z'_M \co \GG_{m,R} \to \uGO(q) \end{equation}
of $R$--group schemes, representing the homomorphism of group functors $\GG_m(S) = S\ti \to \GO(q_S)$, $x \mapsto x\Id_{M \ot S}$, $S\in \Ralg$. Since  
$x\Id_{M\ot S}$ is a similitude with multiplier $x^2$, it follows that
 the homomorphism $\mult$ is surjective in the flat topology:
\begin{equation}  \label{simi-defisch-b11}
1\longto \uO(q) \xrightarrow{\; \inc\; } \uGO(q) \xrightarrow{\; \mult \;}
      \GG_{m,R} \longto  1.
\end{equation}
The flat quotient of $\uGO(q)$ by the central $\GG_m \cong z'_M(\GG_m)$ of \eqref{simi-defisch1}  is denoted
\[ \uPGO(q) = \uGO(q)/ \GG_m\]
and called the {\em projective similitude group.} It is smooth affine by \cite[C.3.12]{Co1}. Since the central $\GG_m \subset \uPGO(q)$ intersects $\uO(q)$ in the central $\bmu_{2,R}$, we get the commutative diagram of group homomorphisms
\begin{equation}\label{simi-defisch-b2}\vcenter{
  \xymatrix@C=40pt{\uO(q) \ar[r]^{\inc} \ar[d]& \uGO(q) \ar[d] \\
       \uO(q)/\bmu_{2,R} \ar[r]^\cong & \uPGO(q)}
}\end{equation}
where the vertical maps are the canonical quotient homomorphisms and where the bottom horizontal map is a monomorphism by construction. That it is also an epimorphism, hence an isomorphism, follows from surjectivity of $\mult$ in \eqref{simi-defisch-b1}. Finally, following \cite[p.~386]{Co1}, we denote by
\[ \uGSO(q) \]
the subgroup sheaf of $\uGO(q)$ generated by $\uSO(q)$ and the central $\GG_m\subset \uGO(q)$, see \eqref{simi-defisch-a} and \eqref{simi-defisch-b} below for a discussion of this group scheme.
\lv{
(2023-04-01, PG)
Au del\`a du cas des corps, il n'y a pas de notion de sous-sch\'ema en groupes engendr\'e par des sous-groupes. Ici heureusement cela marche. On consid\`ere le morphisme produit $\uSO(q) \times \GG_m \to \uGO(q)$. Si $q$ est de rang impair, c'est un monomorphisme et donc une immersion ferm\'ee (\cite[XVI.1.5]{SGA3}) et cela d\'efinit $\uGSO(q)$. Si $q$ est de rang pair, le noyau est le $\bmu_2$ diagonal et on note alors $\uGSO(q)$ le quotient de $\uSO(q) \times \GG_m$ par $\bmu_2$. C'est l\'egitime car $\bmu_2$ est diagonalisable. Tu peux faire tout en meme temps en quotientant par le groupe diagolanisable  $\uSO(q) \cap \GG_m$. Je n'avais pas vu que tu discutes cela plus loin.}

Since all the group schemes defined above respect the orthogonal rank decomposition, 
it is no harm to assume that $M$ has constant rank  in order to investigate these group schemes more closely.
\ms

\begin{inparaenum}[(a)] \item\label{simi-defisch-a}
{\em $(M,q)=(M_o,q_o)$ is a  quadratic space with $M_0$ of constant odd rank $\ge 3$}: Since $\uSO(q_o) = \Ker(\det)$ by \eqref{sogsc22}, we obtain $\GG_m \cap \uSO(q_o) = \{1\}$ and then the first isomorphism of \eqref{simi-defisch-a1},
\begin{equation}  \label{simi-defisch-a1}
 \uGSO(q_o) \cong \GG_{m,R} \times \uSO(q_o) \cong \uGO(q_o),
\end{equation}
the second being a consequence of \eqref{orthsc-2} and the exact sequence \eqref{simi-defisch-b1}. Moreover, \eqref{simi-defisch-a1} implies
\begin{equation}  \label{simi-defisch-a2}
 \uPGO(q_o) \cong \uSO(q_o).
\end{equation}
Consequence of \eqref{simi-defisch-a1} and \eqref{simi-defisch-a2}:
\end{inparaenum}\begin{enumerate}[label=\rm (\roman*)]
  \item\label{simi-defisch-a-i} $\uGO(q_o)$ is a reductive $R$--group scheme, while $\uPGO(q_o)$ is a semisimple $R$--group scheme of adjoint type.

  \item \label{simi-defisch-a-ii} In particular, both $\uGO(q_o)$ and $\uPGO(q_o)$ are smooth affine and have connected geometric fibres.
\end{enumerate}
Here is a diagrammatic summary:
\[\xymatrix@C=35pt{
 1 \ar[r] & \bmu_2 \ar[r]\ar@<2pt>[d]^{z_M} & \GG_m \ar[d]^{z'_M}\ar[r]^{\times^2} & \GG_m \ar[r]\ar@{=}[d] & 1
  \\
  1 \ar[r] & \uO(q_o)\ar@<2pt>[u]^\det \ar[r]\ar@<2pt>[d] \ar[r] &\uGSO(q_o)= \uGO(q_o) \ar[r]^>>>>>>{\mult}\ar[d] & \GG_m \ar[r] &1
 \\
   & \uSO(q_o)\ar@<2pt>[u]^{\inc} \ar[r]^\cong & \uPGO(q_o)}\]
\sm

\begin{inparaenum}[(a)]\setcounter{enumi}{1}
\item \label{simi-defisch-b}  {\em $(M,q)=(M_e, q_e)$ is a  quadratic space with $M$ of constant even rank $\ge 2$}: In this case $\GG_m \cap \uSO(q_e) = \GG_m \cap \uO(q_e) \cong \bmu_{2,R}$, so that
\begin{equation}
  \label{simi-defisch-b1} \uGSO(q_e) \cong \big( \GG_m \times \uSO(q_e)\big)/\bmu_{2,R}.
\end{equation}
As $\uGSO(q_e)$ is a closed and normal subgroup of $\uGO(q_e)$ this implies
\begin{equation}
  \label{simi-defisch-b2}
\begin{split}  \uGO(q_e) / \uGSO(q_e) & \cong \big(\big( \GG_m \times \uO(q_e)\big)/\bmu_2\big)\big/
  \big(\big( \GG_m \times \uSO(q_e)\big) / \bmu_2\big)
  \\ &\cong \uO(q_e) / \uSO(q_e) \cong (\ZZ/2\ZZ)_R,
\end{split} \end{equation}
where the last isomorphism follows from the exact sequence \eqref{sogsc2a}:
\[ 1 \longto \uSO(q_e) \longto  \uO(q_e) \xrightarrow{\Di} \ZZ/2\ZZ \longto 1.
\]
We thus obtain  a quotient homomorphism
\begin{equation}\label{simi-defisch-dg}  \mathrm{GDick} \co \uGO(q_e) \longto (\ZZ/2\ZZ)_R \end{equation}
extending the Dickson homomorphism $\Di$ whose kernel is $\uGSO(q_e)$. Thus
\[ \uGSO(q_e) = \uPGO^+(q_e)\]
using the notation of \cite[4.4.0.35]{CF}.
Consequences:
\end{inparaenum}
\begin{enumerate}[label=\rm (\roman*)]\setcounter{enumi}{2}
\item \label{simi-defisch-b-iii} $\uGO(q_e)$ and $\uPGO(q_e)$ are smooth affine $R$--group schemes whose geometric fibres have two connected components.

\item $\uGSO(q_e)$ is a reductive $R$--group scheme.
\end{enumerate}
Diagrammatic summary:
\[ \xymatrix@C=35pt{
      & 1 \ar[d] & 1 \ar[d]
  \\
  1 \ar[r] & \uSO(q_e)\ar[d]_\inc \ar[r]^{\inc} & \uGSO(q_e) \ar[d]_\inc \ar[r]^{\mult} & \GG_m \ar@{=}[d]\ar[r] & 1
  \\
  1 \ar[r] & \uO(q_e) \ar[r] \ar[d]_{\Di}& \uGO(q_e)\ar[d]_{\mathrm{GDick}} \ar[r]^\mult & \GG_m \ar[r] & 1
 \\
  & (\ZZ/2\ZZ)_R \ar@{=}[r] \ar[d]& (\ZZ/2\ZZ)_R \ar[d]
\\ & 1 & 1
}\]

\newpage

\section{Spheres and their smooth loci}\label{sec:sphere}

\comments{Section summary needs to be written}

We start with a general Lemma \ref{old*} to be used in the proof of Lemma~\ref{q-Faser}, describing the value schemes of quadratic forms. We then specialize to invertible values in Lemma~\ref{ussl}. We will consider another specialization of Lemma~\ref{q-Faser} later in section \S\ref{sec:springer} on Springer's Odd Degree Theorem by taking zero values. 


\begin{lem}\label{old*} Let $X,Y,Z_1, \ldots, Z_n$ be indeterminants over an integral doamin PID $A$ and let $0 \ne b \in B=A[Z_1, \ldots, Z_n]\subset A[Z_1, \ldots, Z_n,X,Y] = B[X,Y]$. Then $B[X,Y]/(XY+b)$ is an integral domain. Moreover, if $A$ is a Pr\"ufer domain, e.g. a Dedekind domain, then $B[X,Y]/(XY + b)$ is a flat $A$--module under the canonical $A$--action.   
\end{lem}

\begin{proof} It is straightforward to verify that $XY + b$ 
is an irreducible element of the polynomial ring $B[X,Y]$. Hence $B[X,Y]/(XY+b)$ is an integral domain. 
Since $A$ embeds into $B[X,Y]/(XY+b)$, this is a torsion-free and therefore flat $A$--module \cite[VII, \S2, Exc.~12]{BAC2}. 
\end{proof}
\sm 

For the notion of the content of a polynomial, like a quadratic form, see \ref{cop}; the concept of a universally schematically subscheme is recalled in \ref{tod}.   

\comments{(2026-01-22) We need $\Cont(q)$ to ensure that $q_F \ne 0$ for any $R$--field $F$. This is used in the proof of \eqref{q-Faser0}. }

\begin{lem}\label{q-Faser} Let $(M,q)$ be a faithful quadratic module over a ring $R$ with $\Cont(q) = R$, see {\rm \ref{cop}\eqref{cop-c}}, and let $r\in R$. Then the $R$--functor $\ulV_{\, q,r}$, given on $T$--points, $T\in \Ralg$, by  
\[ 
   \ulV_{\, q,r}(T) = \{ m \in M_T : q_T(m) = r \ot 1_T\},
\]
is represented by an affine, finitely presented $R$--scheme $\uV_{q,r}$. \sm

\begin{inparaenum}[\rm (a)]\item\label{q-Fasera} 
The smooth locus $\uV_{q,r}\rmsm$ of $\uV_{q,r}$ satisfies 
 \begin{equation}\label{q-Faser0} \begin{split}
  \uV\rmsm_{q,r}(T)& = \big\{  x  \in \uV_{q,r}(T) : \,  b_q(x, \cdot) \co  M_T \to T \hbox{\enskip  is surjective}  \big \} \\
   &= \big\{  x  \in \uV_{q,r}(T) : \,  b_q(x, y) = 1_T \hbox{\enskip for some $y\in M_T$} \big \} 
\end{split} \end{equation}
for each $R$--algebra $T$. Hence, $\uV\rmsm_{q,r}(T)$ consists of unimodular vectors. 
\sm 

\new
\item \label{q-Faserb} The smooth locus $\uV_{q,r}\rmsm$ is the intersection of $\uV_{q,r}$ and the open subscheme $\Sur_{b_q}$ of $\uW(M)$.  It is a quasi-compact open subscheme of $\uV_{q,r}$, hence also finitely presented and quasi-affine.  \sm 
    
\item\label{q-Faserbb} Both $\uV_{q,r}$ and $\uV\rmsm_{q,r}$ are stable under base change: for $R'\in \Ralg$ there exists a canonical isomorphism of $R'$--schemes
\begin{equation}\label{q-Faserb1} 
 \uV_{q,r}\times_R \, R' \cong \uV_{q_{R'}, r \ot 1_{R'}}\, , \qquad
  \uV_{q,r}\rmsm\times_R \, R' \cong \uV\rmsm_{q_{R'}, r \ot 1_{R'}}. 
 \end{equation}  
\enew

\item \label{q-Faserc} Moreover, {\em if $q$ is nonsingular\/}, the following additional properties hold.   
\end{inparaenum}

\begin{enumerate}[label={\rm (\roman*)}]
  \item \label{q-Faserii} $\uV_{q,r}$ is a flat $R$--scheme. Furthermore, if $q=q_{0,n}$ is the split quadratic form of rank $n\in \NN_+$, the structure morphism $\uV_{q_{0,n},r} \to \Spec(R)$ is  essentially free, {\rm \ref{ag}\eqref{ag-essen}}. 
      \sm 
      
 
  \item \label{q-Faseri} If $\rank M\ge 3$, then $\uV_{q,r}$ has geometrically integral fibres, and $\uV_{q,r}\rmsm$ is universally schematically dense in $\uV_{q,r}$.  
\end{enumerate}
\end{lem}%


\begin{proof} We abbreviate $\ulV_{\, r} = \ulV_{\, q,r}$ and $\uV_r = \uV_{q,r}$. Clearly $\ulV_{\, r}$ is represented by the affine and finitely presented $R$--scheme $\Spec\big( R[\uW(M)]/(q-r)\big)$, viewing $q-r$ as a polynomial on $\uW(M)$. \sm 

\eqref{q-Fasera} To prove \eqref{q-Faser0}, we first assume that $T=k$ is a field. Then $M_T$ is free, say of rank $n$, and $k[\uV_r] = k[X_1, \ldots, X_n]/(q-r)$.   
According to the Jacobi criterion \cite[\S 2.2, Prop. 7]{BLR},
a point $v \in \uV_r(k)$ belongs to $\uV_{r}\rmsm (k)$ if and only if
the partial derivatives $\partial q/ \partial X_i$'s do not vanish simultaneously at $v$, i.e., $b_q(v, \cdot): M_k \to k$ is nonzero, equivalently, $b_q(v, \cdot)$ is surjective.

Let now $T\in \Ralg$ be arbitrary. Since $\uV_{r}\rmsm$ is open in $\uV_{r}$, a point $v\in \uV_{r}(T)$ lies in $\uV_{r}\rmsm(T)$ if and only if $v_{T/\m}$ belongs to $\uV_{r}\rmsm(T/\m)$ for every maximal ideal $\m\ideal T$, due to Lemma \ref{opemaLG}. 
It follows that $v \in \uV_{r}\rmsm(T)$ if and only if $b_q(v_{T/\m}, \cdot) \co M_T/\m M_T \to T/\m$ is surjective. This last condition is equivalent to the surjectivity of $b_q(x,\cdot)$. The second equation of \eqref{q-Faser0} is obvious, and it implies unimodularity by \ref{unimod}\ref{unimodii}. 
\sm 

\new
\eqref{q-Faserb} The first part follows from \eqref{q-Faser0} and \ref{sursch}\eqref{sursch-b}. The open immersion $\uV_{q,r}\rmsm \to \uV_{q,r}$ is quasi-compact by base change from the quasi-compact open immersion $\Sur_{b_q} \to \uW(M)$. In particular, $\uV\rmsm_{q,r}$ is a quasi-compact scheme. The remaining assertions then follow. 
\sm 

\eqref{q-Faserbb} The first equation in \eqref{q-Faserb1} is clear. The second then follows, using the analogous equation for $\Sur_b$ proven in \eqref{sursch-c1}. \sm  
\enew

\eqref{q-Faserc} For the remaining proof, it will be useful to first identify the coordinate ring $R[\uV_r]$ in case $q$ is the split quadratic form $q_{0,n}$ of rank $n$. Depending on $n$, this ring has the following structure: \sm 

\begin{tabular}{l l}
$n=1,$ & $R[X]\big/(X^2-r)$,   \sm \\ 
 $n=2,$ & $R[X,Y]\big/(XY-r)$, \sm \\
$n \ge 3$, & $B_n[X,Y]/(XY + b_n),$
\end{tabular} \sm

\noindent where $B_n$ is a polynomial ring over $R$ and $0 \ne b_n \in B_n$ are given as follows: 

\begin{tabular}{l l l}
 $n=3$, & $B_3 = R[X_0]$, &   $b_3 = X_0^2 - r$, \sm \\  
  $n=2m \ge 4,$ & $B_n =  R[X_2, Y_2, \ldots X_m, Y_m],$
        & $b_n = \sum_{i=2}^m  X_i Y_i - r$,
 \sm  \\
 $n=2m+1 \ge 5,$ & $B_n =  B_{2m}[X_0],$
        & $b_n = X_0^2 + b_{2m}$.
\end{tabular} \sm 

\ref{q-Faserii} For the proof of flatness, i.e., flatness of the $R$--module $R[\uV_r]$, we note that the schemes $\uV_r$ and flatness respect the rank decomposition of $(M,q)$. We can therefore suppose that $M$ has constant rank $r$ for proving flatness. Also, since flatness is stable under base change and faithfully flat descent, \cite[Tags 01U9, 02JZ]{St}, 
we  can assume that $q$ is the split quadratic form $q_{0,n}$ over $R$. It is then enough to show that $R[\uV_{q,r}]$ for $q=q_{0,n}$ is a free $R$--module, which then also proves the second claim, i.e., that the structure morphism of $\uV_{q_{0,n}, r} \to \Spec(R)$ is essentially free in the sense of \ref{ag}\eqref{ag-essen}.
Indeed, using the description of this ring given above, its $R$--module structure is isomorphic to: \sm

\begin{tabular}{l l}
$n=1,$ & $R[X]\big/(X^2-r) \cong R \oplus R$,  
\sm \\ 
$n=2,$ & $R[X_1, Y_1]\big/ (X_1Y_1-r) \cong R[X_1] \oplus Y_1 R[Y_1],$
\sm \\ 
$n=3,$ & $R[X_0, X_1, Y_1]\big/(-r + X_0^2 + X_1 Y_1) $
\sm \\
& \quad $\cong
     R[X_1, Y_1] \oplus X_0 R[X_1, Y_1], $
\sm \\
 $n = 2m \ge 4,$ & $R[X_1, Y_1, \ldots, X_m, Y_m] \big /
    (-r + X_1 Y_1 + \sum_{i=2}^m X_i Y_i)$ 
\sm \\  
& \quad $\cong R[X_1, X_2, Y_2, \ldots, X_m, Y_m]$
\sm \\
& \qquad $ \oplus Y_1 R[Y_1, X_2, Y_2, \ldots, X_m , Y_m],$
\sm \\
$n=2m+1 \ge 5, $ & $R[X_0, X_1, Y_1, \ldots, X_m , Y_m] \big/ (-r + X_0^2 + \sum_{i=1}^m X_i Y_i )$ 
\sm \\
& \quad $\cong  R[X_1, Y_1, \ldots, X_m, Y_m]$
\sm \\
& \qquad $\oplus X_0 R[X_1, Y_1, \ldots, X_m, Y_m].$
\end{tabular} 
\sm

\noindent We note that in case $R$ is a Pr\"ufer domain, flatness for the cases $n\ge 3$ also follows from Lemma~\ref{old*}. \sm 


\ref{q-Faseri} Let $K$ be an algebraically closed extension of $R$. By \ref{ag}\eqref{ag-c}, it suffices to show that $\Spec(R[\uV_r])_K = \Spec(K[\uV_r])$ is an integral scheme, equivalently, that $K[\uV_r]$ is an integral domain. Since $q$ is necessarily split, $K[\uV_r] = B_n[X,Y]\big/ (XY + b_n)$ is as explained above with $R=K$. Thus, the claim follows from Lemma~\ref{old*}. 

Since $\uV_r$ is a flat, finitely presented $R$--scheme with geometrically integral fibres, we can apply the characterization \ref{todle}\eqref{todle-b} to prove universally schematically denseness of $\uV_{r}\rmsm$. Thus, we need to show that $\uV_{r}\rmsm(K)= \uV_{q_r, r\ot 1_K}\rmsm (K) \ne \emptyset$ for any algebraically closed $R$--field $K$. Because $(M_K, q_K)$ contains an isotropic vector, the claim follows from the Example~\ref{lem_smooth_locus_hyp} below. \end{proof}

\subsection{Example $\uV_{q,r}\rmsm(R) \ne \emptyset$} 
\label{lem_smooth_locus_hyp} 
{\em If $(M,q)$ contains a hyperbolic plane, e.g., $(M,q)$ is a quadratic space containing an isotropic vector, then $\uV_{q,r}\rmsm (R) \ge |R\ti|\ge 1$ for every $r\in R$.} \sm 

Indeed, let $(e,f)$ be a hyperbolic pair in $(M,q)$. Then 
\[ 
   R\ti \to \uV_{q,r}\rmsm (R), \quad c \mapsto v_c = rce + c\me f 
\]
is an injective map to $\uV_{q,r}(R)$. Its image lies in the smooth part since $b_q(v_c,c e) = 1$. 

If $(M,q)$ is a quadratic space containing an isotropic vector, 
the quadratic space $(M,q)$ contains a hyperbolic pair $(e,f)$ by \ref{isotrop}\eqref{isotrop-d}. 

\subsection{The extended sphere $\wdh\uS_q$, the sphere $\uS_{q,a}$, and  its smooth locus $\uS^{\rm sm}_{q,a}$} \label{ussl}
We now define the main objects of study in this section. Let $(M,q)$ be a quadratic $R$--module. \sm 

\begin{inparaenum}[(a)] \item ({\em The extended sphere}) 
Viewing $q$ as a polynomial on $\uW(M)$, the {\em extended sphere\/} is the principal open subscheme 
\begin{equation}\label{ussl1}  
      \wdh \uS_q = D(q)
\end{equation}  
of the affine $R$--scheme $\uW(M)$, sometimes also denoted $\uW(M)_q$ or $\uW(M) \setminus \{q=0\}$. Its functor of points is $T \mapsto \wS_{q_T}= \{ m \in M_T: q_T(m) \in T\ti\}$, the extended sphere of the quadratic $T$--module $(M_T, q_T)$, see \eqref{quadco-aa2}. \sm 

\item ({\em Spheres})  For $a \in R^\times$, the {\em $a$--sphere} is  $\uS_{q,a} = \uV_{q,a}$, as defined in \ref{q-Faser}. We call $\uS_q := \uS_{q,1}$ the {\em unit sphere $\uS_q$} of $(M,q)$. If $a$ is not important, we will often refer to $\uS_{q,a}$ as a {\em sphere}.
    \sm 
    
\item ({\em The smooth sphere}) We denote by $\uS^{\rm sm}_{q,a}$ the smooth locus of the morphism $\bS_{q,a} \to \Spec(R)$ \cite[Tag 01V5]{St}, that is, the set of points where $\bS_{q, a}$ is smooth. It is an open subscheme of $\uS_{q,a}$.
\end{inparaenum}
\ms

The following lemma presents some geometric properties of the schemes introduced in \ref{ussl}. For easier reference we re-state those already established in \ref{q-Faser} for $\uS_{q,a} = \uV_{q,a}$. 

\comments{(2026-02-26) Deleted the previous part (a) of \ref{lem_smooth_locus-LG}: "$\wdh \uS_q$ is a smooth affine $R$--scheme." since this is obvious from \eqref{ussl1}. 

Old proof: As a principal open subscheme of the smooth affine scheme $\uW(M)$ the $R$--scheme is quasi-compact, separated and smooth, hence also of finite presentation. (A principal open subscheme of an affine scheme is qc \cite[Prop. 2.5]{GW}.) }

\begin{lem}[Geometric properties of spheres] \label{lem_smooth_locus-LG}
Let $(M,q)$ be a faithful quadratic space and let $a  \in R^\times$.
Then the following hold. \sm 

\begin{enumerate}[label={\rm (\alph*)}] 

\item \label{lem_smooth_locus-LGb}  $\uS_{q,a}$ is a flat affine $R$--scheme of finite presentation.
 \sm 
   
\item \label{lem_smooth_locus-LGc} The smooth locus of $\uS_{q,a}$ satisfies
 \begin{align*}
  \uS^{\rm sm}_{q,a}(T)&= \big\{  x  \in \uS_{q,a}(T) : \, b_q(x, \cdot) \co  M_T \to T \hbox{\enskip  is surjective}  \big \} \\
   &= \big\{  x  \in \uS_{q,a}(T) : \,  b_q(x, y) = 1_T \hbox{\enskip for some $y\in M_T$} \big \} 
 \end{align*} 
for each $R$--algebra $T$. \sm

\item \label{lem_smooth_locus-LGd} $(M,q)$ is regular if and only if\/ $\uS_{q,a}$ is smooth over $R$. \sm

\item \label{lem_smooth_locus-LGe} $\uS_{q,a}\rmsm$ is a quasi-compact open subscheme of $\uS_{q,a}$, hence a quasi-affine scheme of finite presentation.
    \sm

\item \label{lem_smooth_locus-LGf} Suppose $M$ has rank $\ge 2$.
\sm
  \begin{enumerate}[label={\rm (\roman*)}]

    \item \label{lem_smooth_locus-LGfi} The $R$--scheme $\uS_{q,a}$ has geometrically integral fibres. \sm

   \item \label{lem_smooth_locus-LGfiii} If $M$ has constant rank $r\ge 2$, the geometric fibres of $\uS_{q,a}$ have dimension $r-1$ and $\uS_{q,a}\rmsm$ is smooth  of relative dimension $r-1$.%
       \sm 

 \item \label{lem_smooth_locus-LGfii} The $R$--scheme $\uS^{\rm sm}_{q,a}$ is universally schematically dense in $\uS_{q,a}$. \sm

\end{enumerate} \end{enumerate}
\end{lem}

\begin{proof}
\ref{lem_smooth_locus-LGb}, \ref{lem_smooth_locus-LGc} and \ref{lem_smooth_locus-LGe}  are special cases of results in Lemma~\ref{q-Faser}. 
\sm

\ref{lem_smooth_locus-LGd} Regularity of a quadratic form can be checked on the fields $R/\m$ for $\m$ a maximal ideal of $R$, and, according to the fiberwise smoothness criterion \cite[IV$_4$, 17.8.2]{EGA}, so can be smoothness of the finitely presented flat $R$-scheme $\uS_{q,a}$. We are thus reduced to the field case and even to an algebraically closed field $k$ and $a=1$. In this case, \ref{lem_smooth_locus-LGc} states that
\begin{equation} \label{lem_smooth_locus31} \begin{split}
\uS^{\rm sm}_q(k) &= \{  x  \in \uS_q(k) : b_q(x,\cdot): M \to k \hbox{\enskip is surjective}\} \\ & =
\uS_q(k) \setminus \rad(b_q).
\end{split} \end{equation}
If $q$ is regular, then $\rad(b_q) =0$, so all closed points of $\uS_q$ are smooth
and hence $\uS_q$ is smooth. Conversely, if $\uS_q$ is smooth, we have $\uS_q(k) \cap \rad(b_q) = \emptyset$, i.e., $q\big( \rad(b_q)\big) \ne 1$, and then $q( \rad(b_q) )=0$, so that $\rad(b_q)= \rad(q)$. Since $q$ is nonsingular, $\rad(q)=0$, whence $\rad(b_q)=0$. Thus $q$ is regular. \sm

\ref{lem_smooth_locus-LGf} If $M$ has rank $\ge 3$, parts \ref{lem_smooth_locus-LGfi} and \ref{lem_smooth_locus-LGfii} follow from Lemma~\ref{q-Faser}\ref{q-Faseri}. In case $\rank M = 2$, the coordinate ring $K[\uS_{q,a}]$, $K$ an algebraically closed $R$--field, is $K[X,Y]/(XY-a) \cong K[X,X\me]$, cf.~the proof of Lemma~\ref{q-Faser}, proving \ref{lem_smooth_locus-LGfi}. The claim \ref{lem_smooth_locus-LGfii} is trivial in this case, as $q$ is regular and therefore $\uS_{q,a}\rmsm = \uS_{q,a}$ by \ref{lem_smooth_locus-LGd}. Finally, \ref{lem_smooth_locus-LGfiii} follows from \ref{lem_smooth_locus-LGfi}. \end{proof}

\subsection{Examples (smooth spheres $\uS_{q,a}\rmsm$)} \label{lem_smooth_locus_exam} \begin{inparaenum}[(a)] \item\label{lem_smooth_locus_one}  ({\em Rank $1$}) Let $(M,q) = (R, \lan 1 \ran_q)$ and let $a\in R\ti$. Then $R[\uS_{q,a}] = R[X]/(X^2 -a)$, which is in general not an integral domain, e.g., if there exists $b\in R$ with $b^2 =a$. This explains the rank assumption in Lemma~\ref{lem_smooth_locus-LG}\ref{lem_smooth_locus-LGf}, parts \ref{lem_smooth_locus-LGfi} and \ref{lem_smooth_locus-LGfiii}. 
Also, if $2R=0$, then $\uS_{q,a}^{\rm sm}  = \emptyset$, while $\uS_{q,a}\ne \emptyset$.
This justifies the rank assumption in \ref{lem_smooth_locus-LGfii}  of   Lemma~\ref{lem_smooth_locus-LG}\ref{lem_smooth_locus-LGf}.
\sm

\item\label{lem_smooth_loci} Let $(M,q)$ be a faithful quadratic module and let $M'\subset M$ a complemented submodule. Put $q' = q|{M'}$. Then 
\[ \uS_{q', a}(T) \subset \uS_{q,a}(T) \quad \hbox{and}\quad
      \uS_{q', a}\rmsm(T) \subset \uS_{q,a}(T)\rmsm
\]
for any $T\in \Ralg$. The second inclusion follows from the second equation in  \ref{lem_smooth_locus-LG}\ref{lem_smooth_locus-LGc}. \sm 

\item \label{lem_smooth_locus_hypn} 
{\em If $(M,q)$ contains a hyperbolic plane, e.g., $(M,q)$ is a quadratic space containing an isotropic vector, then $|\uS_{q,a}\rmsm (R)|\ge |R\ti| \ge 1$ for all $a\in R\ti$.} Indeed, this is a special case of Example~\ref{lem_smooth_locus_hyp}. \sm

\item \label{lem_smooth_locus_isof} {\em Let $(M,q)$ be a quadratic space over a field $k$  with $\dim_k M \ge 3$. If $\uS_{q,a}\rmsm(k) \ne \emptyset$, then $|\uS_{q,a}\rmsm(k) | \ge 3$ for any $a\in k\ti$.} \sm 
    
    {\em Proof.} Suppose first that $(M,q)$ is isotropic. By \ref{isotrop}\eqref{isotrop-d}, $(M,q)$ contains a hyperbolic pair $(e,f)$. Its orthogonal space $(ke + kf)^\perp$ is non-zero. Therefore, by \eqref{qf-perp1} and \ref{LGqdi}, there exists $g\in (ke + kf)^\perp$ with $q(g) = b \in k\ti$. By \eqref{lem_smooth_loci} it then suffices to consider $M = ke \oplus kf \oplus kg$, so that  $q(xe + yf + zg) = xy + bz$. If $a=b$, then 
    $(a,1,0), (0,1,1), (1,0,1)$ lie in $\uS_{q,a}\rmsm(k)$, and if $a\ne b$, then $\uS_{q,a}\rmsm(k)$ contains the three points 
    $(a,1,0), (1, a-b, 1),  (\frac{b}{a}, \frac{a}{b}(a-b), 1)$.
    
    To prove our claim in general, note that by \cite[12.3)]{Kneser}, any quadratic space $(M,q)$ over a finite field is isotropic, so that the above applies. If $k$ is infinite, we will see in \ref{smoLG}\eqref{smoLG-e} that $\uS_{a,q}\rmsm (k)$ is infinite. (We won't use \eqref{lem_smooth_locus_isof} in the meantime.)   
    \sm

\item \label{lem_smooth_locus_exam-c} {\em If $(M,q)$ is a quadratic space over an LG ring $R$ with $\rank M \ge 2$, then $\uS_{q,a}\rmsm \ne \emptyset$ for some $a\in R\ti$.}

    To see this, it is no harm to assume that $M$ has constant rank. The structure of such forms, described in \ref{nqf-LG}, shows that $(M,q) = (M', q') \perp (M'', q'')$ with $M'$ of constant rank $2$. Because $q'$ is regular, $\uS_{q',a}(R) \ne \emptyset$ for some $a\in R\ti$ by \ref{LGqdi}, and $\uS_{q',a}(R) = \uS_{q',a}\rmsm$ by 
    \ref{lem_smooth_locus-LG}\ref{lem_smooth_locus-LGd}.
       Since $  \uS_{q',a}\rmsm \subset   \uS_{q,a}\rmsm$ by 
    \ref{lem_smooth_locus-LG}\ref{lem_smooth_locus-LGc}, we are done.
\end{inparaenum} \ms

In the following we will use the obvious terminology that a {\em nonsingular submodule\/} or a {\em regular submodule\/} of a quadratic $R$--space $(M,q)$ is a  submodule $V$ of $M$ for which $q|_V$ is nonsingular or regular respectively. A {\em regular plane\/} is a regular submodule of constant rank $2$.  
Our next aim is Proposition~\ref{embp} showing that one can embed vectors in quadratic spaces over semilocal rings into quadratic submodules of rank $3$.  We first consider the field case.   

\begin{lem}\label{embk} Let $(M,q)$ be a quadratic space over a field $k$ of dimension $\dim_k M \ge 3$ and let $m\in M$. Then there exists a nonsingular subspace $W \subset M$ with $\dim_k W = 3$ and $m\in W$.  
\end{lem}

\begin{proof} For later use we will consider the following more precise statement \ref{embk1} below which, as we will see, holds in all cases except case~\eqref{embk-e}:
\begin{enumerate}[label=($\star$)]
  \item\label{embk1} {\em There exist a regular plane $P$ in $(M,q)$ such that $m \in P$.} 
\end{enumerate}
If \ref{embk1} is true, we are done. Indeed, we have the orthogonal decomposition $M = P \perp P^\perp$ where $P^\perp$ is nonsingular by \eqref{qf-perp1} and of positive dimension by assumption. Therefore, by ~\ref{LGqdi}, $P^\perp$ contains $m'$ such that $q(m') \in k\ti$. Since $k m'$ is a nonsingular submodule by \ref{mx_lem}\eqref{mx_lema}, the submodule $V = P \oplus km'$ is a submodule we are looking for. 
\sm 

\begin{inparaenum}[(a)]
  \item\label{embk-a} $m=0$: Then \ref{embk1} holds since $(M,q)$ contains a regular plane by \cite[7.29 and 7.32]{EKM}.  
\sm 

\item $m\ne 0$, $q(m) = 0$: Since then $m$ is isotropic, it is part of a hyperbolic pair $(m, n)$ by \ref{isotrop}\eqref{isotrop-d}, so that $P =km \oplus kn$ satisfies \ref{embk1}. \sm
    
\item\label{embk-c} $m \in \uS_{q,a}(k)$, $\Char(k) \ne 2$: Then regular = nonsingular by \eqref{qf-ns2}, and hence $M = km \perp M'$ with $M'$ a nonsingular submodule of positive dimension by \ref{mx_lem}\eqref{mx_lemc}. As recalled above, there exists $m' \in M'$ with $q(m') \ne 0$. Then $P=km \oplus km'$ satisfies \ref{embk1}.  \sm 
    
\item \label{embk-d} $m\in \uS_{a,q}\rmsm(k)$, $\Char(k) = 2$: Then there exists $m'\in M$ such that $b_q(m,m') = 1$ and $P=km \oplus km'$ is regular. \sm 
    
\item\label{embk-e} $m\in \uS_{a,q}(k) \setminus \uS_{q,a}\rmsm(k)$, $\Char(k) = 2$: Then $m\in \rad b_q$ and $M= km \perp M'$ with $M'$ a regular submodule,  \cite[7.32]{EKM}. As $M' \ne 0$, it contains a regular plane $P'$, so that $W = km \perp P'$ fulfills the requirements of the lemma. 
\end{inparaenum}  
\end{proof}
\sm

We can extend Lemma~\ref{embk} to semilocal rings.   

\begin{prop}\label{embp} Let $R$ be a semilocal ring, let $(M,q)$ be a quadratic $R$--space with $\rank M \ge 3$, and let $m\in M$.  Then there exists a complemented nonsingular submodule $W$ of $M$ with $\rank W = 3$ and $m\in W$. Moreover, $\uS\rmsm_{q|_W, a}(R) \ne \emptyset$ for some $a\in R\ti$. \end{prop}

\begin{proof} Since the claim is stable with respect to the standard rank decomposition of $(M,q)$, see \ref{qf}\eqref{qf-redc}, we can assume that $M$ has constant rank.%
\sm 

Let $\ol R = R /\Jac(R)$, let $\Max$ be the set of maximal ideals of $R$, and let $M \to \ol M = M\ot_R \ol R = \prod_{\gm \in \Max} M_{R/\gm}$, $m \to \ol m = (m_{R/\gm})_{\gm \in \Max}$ be the canonical map. By Lemma~\ref{embk}, for every $\gm \in \Max$, we can choose a $3$--dimensional nonsingular submodule $W[\gm]$ of the $R/\gm$--space $(M,q)_{R/\gm}$ containing $m_{R/\gm}$. Furthermore, we can choose a submodule $V[\gm]$ of the $R/\gm$--module $M_{R/\gm}$ complementing $W[\gm]$. Then $W' = \prod_{\gm \in \Max} W[\gm]$ is a nonsingular free submodule of $(M,q)_{\ol R}$ of $\rank 3$ containing $m_{\ol R}$. It is complemented by the free $\ol R$--submodule $V' = \prod_{\gm \in \Max} V[\gm]$ of $\ol M$. We can lift $W'$ to a nonsingular submodule $W$ of $(M,q)$ which has constant rank $3$,  contains $m$ and maps onto $W'$ under the canonical map $M \to \ol M$. Similarly, we can lift $V'$ to a submodule $V$ of $M$ mapping onto $V'\in \ol M$. Since $\ol M = W' \oplus V'$, it follows from \ref{nak}\eqref{nak-dec} that $M = W \oplus V$. 
That $\uS\rmsm_{q|_W, a}(R) \ne \emptyset$ for some $a\in R\ti$, follows from Example~\ref{lem_smooth_locus_exam}\eqref{lem_smooth_locus_exam-c}. 
\end{proof}
\ms

We have seen in \ref{canqf} that an isometry between regular submodules of a quadratic space $(M,q)$ over an LG ring $R$ extends to an isometry of $(M,q)$. If $R$ is a local ring, isometries extend from more general submodules due to the following result of Kneser, see \cite[Th.~8.3]{EKM} for the field case.

\begin{thm}[{\cite[(4,4), (4.5) with $H=M$]{Kneser}}]\label{knesth}
Let $(M,q)$ be a quadratic module over a local ring $(R,\m)$, and let $M_1$ and $M_2$ be two free submodules of $M$ satisfying $b_{M_1} (M) = M_1^*$ and $b_{M_2}(M) = M_2^*$, where for a submodule $N$ of $M$ we put  $b_N \co M \to N^*, m \mapsto b_q(m, \cdot)|_N$. 

Then any isometry $t \co M_1 \to M_2$ extends to an isometry of $(M,q)$. Moreover,  
\begin{enumerate}[label={\rm (\roman*)}]
  \item \label{knesthi} if $R/\m \ne \FF_2$ and $q(M) \ne 0$, or 
  \item \label{knesthii} if $R/\m = \FF_2$, and $q(M^\perp) \ne 0$,
\end{enumerate}
then there exists an extension of $t$ in $\Refl(q)\subset \orth(q)$.  
\end{thm}
\sm 

The assumptions $b_{M_i}(M) = M_i^*$ in \ref{knesth} are fulfilled if $M_1$ and hence also $M_2$ are regular submodules of $(M,q)$, because then $b_{M_i}(M_i) = M_i^*$ and any linear form of $M_i$ extends to a linear form of $M$ since $M_i$ is complemented by \ref{orthLG}\eqref{orthLG-a}. For regular submodules, Theorem~\ref{knesth} could be derived from \cite[Satz~(0.3)]{Knebusch}. However, the weaker assumptions $b_{M_i}(M) = M_i^*$ are crucial for the application of Theorem~\ref{knesth} in Lemma~\ref{lem_sphere_local} where we deal with nonsingular, not necessarily regular submodules.

\begin{cor}\label{lem_sphere_local} Let $R$ be a local ring, let $(M,q)$ be a faithful quadratic space over $R$, and let $a\in R\ti$. \sm 

\begin{enumerate}[label={\rm (\alph*)}] 
\item \label{lem_sphere_local-a}
Then the group $\orth(q)$ acts transitively on the sphere $\uS\rmsm_{q,a}(R)$. \sm 

\item \label{lem_sphere_local-b} Let $R=k$ be a field, and suppose $\dim_k(M) \ge 2$. Then already the group $\SO(q)$ acts transitively on $\uS_{q,a}\rmsm(k)$.    
\end{enumerate}
\end{cor}

\begin{proof} Let $u_1$ and $u_2 \in \uS_{q,a}\rmsm(R)$. By \ref{mx_lem}\eqref{mx_lema}, the submodules $M_1=Ru_1$ and $M_2 = R u_2$ are free of rank $1$. Moreover, since the maps $b_q(u_i, \cdot) \co M\to R$, $i=1,2$, are surjective by \ref{lem_smooth_locus-LG}\ref{lem_smooth_locus-LGc},  the conditions $b_{M_i}(M) = M_i^*$ of Kneser's Theorem \ref{knesth} are fulfilled. Hence, the isometry $M_1 \to M_2$, given by $u_1 \mapsto u_2$, extends to an isometry of $(M,q)$, proving \ref{lem_sphere_local-a}. Part \ref{lem_sphere_local-b} then follows from Lemma~\ref{lem_transitivity_SO}\ref{lem_transitivity_SO3}.
\end{proof}

\lv{
\begin{lem}\label{lem_sphere_field1} Let $k$ be a field, let $(V,q)$ be
a quadratic $k$--space of dimension $\ge 2$, and let $a \in k^\times$. Then the group $\orth(q)$ acts transitively on $\uS\rmsm_{q,a}(k)$. If $\dim_k V \ge 2$, then already $\SO(q)(k)$ acts transitively on $\uS\rmsm_{q,a}(k)$.  
\end{lem}

\begin{proof} Let $x, y \in \uS^{\rm sm}_{q,a}(k)$.
Then $kx$ and $ky$ are isometric subspaces of $(V,q)$
such that $kx \cap \rad(b_q)=0$ and $ky \cap \rad(b_q)=0$. The Witt Extension Theorem \cite[Th.~8.3]{EKM} then says that there exists $g \in \orth(q)(k)$ such that $g(x)=y$. The claim therefore follows from Lemma~\ref{lem_transitivity_SO}\ref{lem_transitivity_SO3}. \end{proof}
}
 
In \ref{sotrans}\eqref{sotrans-a} we will generalise \ref{lem_sphere_local}\ref{lem_sphere_local-b} to local rings. Its proof uses the case of fields, see the proof of \ref{smoLG}\eqref{smoLG-f}. But first a consequence of \ref{lem_sphere_local}\ref{lem_sphere_local-a}, using the general concept of transitivity on the big affine Zariski site of $\Spec(R)$, reviewed in \ref{transi}.

\begin{cor} \label{cor_sphere_homog-orth} Let $R$ be an arbitrary base ring, let $(M,q)$ be a faithful quadratic $R$--space, and let $a\in R\ti$. Then the action of $\uO(q)$ on $\bS^{\rm sm}_{q,a}$ is transitive on the big affine Zariski site of\/ $\Spec(R)$.
\end{cor}

\begin{proof} Since $\uO(q)$ and $\uS_{q,a}\rmsm$ are finitely presented $R$--schemes, Lemma~\ref{surlem} says that it suffices to prove transitivity on local rings. But this follows from  Corollary~\ref{lem_sphere_local}\ref{lem_sphere_local-a}. 

Indeed, we are given an $R$--ring $A$ and two points $y,z \in \bS^{\rm sm}_{q,a}(A)$. We want to show that there exists a Zariski cover $A'$ of $A$ such that $z_{A'} \in \uO(q)(A').y_{A'}$. In other words, we  claim that the orbit map $f_y: \uO(q)(A) \to  \bS^{\rm sm}_{q,a}(A)$, $g \mapsto g.y$, is an epimorphism of Zariski sheaves. Since $\uO(q)$ and $\bS^{\rm sm}_{q,a}$ are of finite presentation, we can apply the surjectivity criterion \ref{surlem}. It is then enough to check that $f_y(A_\p): \uO(q)(A_\gp) \to  \bS^{\rm sm}_{q,a}(A_\gp)$ is onto for each prime ideal $\gp$ of $A$. Put $B= A_\gp$. Since $q_B$ is nonsingular, and since $\uO(q)(B)=\orth(q_B)$ and $\bS\rmsm_{q,a}(B) = \bS\rmsm_{q_B, a \ot 1_B}(B)$,  this is a consequence of Corollary~\ref{lem_sphere_local}\ref{lem_sphere_local-a} applied to the local ring $B$. We conclude that there does indeed exist a Zariski cover $A'$ of $A$ such that $z_{A'} \in f_y( \uO(q)(A')) = \uO(q)(A').y_{A'}$.
\end{proof}

\begin{examples}[\bf More examples of transitivity of $\orth(q)$ on spheres]\label{exatr} 
Transitivity of $\orth(q)$ on spheres holds beyond the case of local rings. We present some easy examples. \sm

\begin{inparaenum}[\rm (a)]  \item \label{exatr-a} 
{\em $(M,q)$ is a quadratic space over an LG ring $R$ with $2\in R\ti$}: As the proof of \ref{lem_sphere_local}\ref{lem_sphere_local-a} shows, transitivity can be viewed as a question whether isometries between nonsingular (= regular) rank $1$ submodules extend to isometries of $(M,q)$. In the given setting, transitivity of $\orth(q)$ on spheres follows from the  extension property of isometries in \ref{canqf}. \sm 

\item {\em $(M,q)$ is a quadratic space free rank $1$ over an arbitrary ring $R$}: Then $(M,q) = \lan z \ran_q$ for some $z\in R\ti$ by \eqref{qfba-one1}. Thus, for $u,u'\in \uS_{q,a}(R)$ we have $a = zu^2 = z u'^2$,  and $(u'/u) \in \mu_2(R)$. Therefore $(u'/u)\Id \in \orth(q)$ sends $u$ to $u'$.     
\end{inparaenum}
\sm 

One of the goals in the remaining part of this section is to prove transitivity of $\SO(q)$ acting on smooth spheres over local rings, see \ref{sotrans}\eqref{sotrans-a}. 
\end{examples}


\comments{(2026-02-26) Improved \ref{smoLG}\eqref{smoLG-a}: $\wdh \uV_{q,u}$ is not only a principal open subscheme of $\wdh \uS_q$, but also of $\uW(M)$, in particular it is smooth affine. Reformulated the proof of \ref{smoLG}\eqref{smoLG-f}: Added that $\wdh \uV_{q,u}$ is a quasi-compact open subscheme of the affine space $\uW(M)$, which is needed in order to apply the characterization of LG rings.} 

\begin{prop}[The open subscheme $\uU_{q,u}$ of $\uS_{q,a}\rmsm$]\label{smoLG}
Let $(M,q)$ be a faithful quadratic $R$--space, let $a\in R\ti$ and suppose there exists $u \in \uS_{q,a}\rmsm(R)$. 

\begin{inparaenum}[\rm (a)] 

\item\label{smoLG-a} The map $\wdh{\vphi} \co \wdh \uS_q \to \uS_{q,a}\rmsm$, given on ${\wdh \uS_q}(T)$, $T\in \Ralg$, by 
\[ v \mapsto \wdh \vphi(v) = - \rho_v(u) = -u + \frac{b_q(u,v)}{q(v)}\, v
\]  
is a well-defined morphism of $R$--schemes which is smooth on the 
principal open subscheme 
\[ 
                 \wdh \uV_{q,u} := (\wdh \uS_q)_{b_q(u, \cdot)} = \uW(M)_q \cap \uW(M)_{b_q(u, \cdot)}
\] 
of the affine space $\uW(M)$.  Clearly, 
\[ \wdh \uV_{q,u}(T) = \{ m\in \wdh \uS_q(T): b_{q_T}(u \ot 1_T, m) \in T\ti\}
\] 
for $T\in \Ralg$. The $R$--scheme $\wdh \uV_{q,u}$ is smooth affine.  
\sm 

\item\label{smoLG-b} The action of $\GG_m$ by scalar multiplication stabilizes $\wdh  \uV_{q,u}$. The flat quotient $(\wdh \uV_{q,u})/\GG_m$ is a scheme, which we can identify 
    with the open affine subscheme 
    \[ 
    \uV_{q,u} := (\wdh \uV_{q,u})/ \GG_m 
    \equiv \PP(M\ch) \setminus \big((\{q=0\}\cup \{ b_q(u,\cdot)=0\})/\GG_m\big)
    \] 
    of\/ $\uW(M)_u/\GG_m\simlgr \PP(M\ch)$. 
    \sm 
    
\item\label{smoLG-c} The map $\wdh \vphi$ of {\rm \eqref{smoLG-a}} descends to an open immersion 
    \begin{equation}\label{smolG-c1} 
     \vphi \co \uV_{q,u} \to \uS_{q,a}\rmsm. 
     \end{equation} 
     We denote by $\uU_{q,u}$ its image, which is also the image of $\wdh \uV_{q,u}$ under $\wdh \vphi$. The $R$--scheme $\uU_{q,u}$ is separated, smooth and finitely presented. It is an open neighbourhood of $u$. \sm 
     
\item\label{smoLG-d} Assume $\rank M \ge 2$. Then $\uU_{q,u}$ has geometrically integral fibers and is universally schematically dense, equivalently  $R$--dense, in $\uS\rmsm_{q,a}$.
    \sm 
   
\item\label{smoLG-e} Let $R=k$ be an infinite field and suppose $\dim_k M \ge 2$. Then $\wdh \uV_{q,u}(k)$, $\uV_{q,k}$ and $\uU_{q,u}(k)$ are  infinite, and $\uU_{q,u}(k)$ and therefore also $\uS\rmsm_{q,a}(k)$ are Zariski dense in $\uS_{q,a}$. \sm      

\item\label{smoLG-f} If $R$ is an LG ring and $\rank_R M \ge 3$, then $\wdh \uV_{q,u}(R)$, $\uV_{q,u}(R)$ and $\uU_{q,u}(R)$ are all non-empty.
\end{inparaenum} 
\end{prop}


\begin{proof} Since $q$ and $u$ are fixed, we can abbreviate $\wdh\uV = \wdh \uV_{q,u}$, $\uV = \uV_{q,u}$ and $\uU = \uU_{q,u}$. \sm 

\eqref{smoLG-a} Recall \eqref{ussl1}: $\wdh\uS_q =\uW(M)_q$. Hence $\wdh \uV$ is the intersection of two principal open subschemes of the affine scheme $\uW(M)$, and is therefore itself a principal open subscheme of $\uW(M)$. 
In particular, it is a smooth affine scheme. It remains to show that $\wdh \vphi$ is a well-defined smooth morphism of schemes.  

Let $v\in \wdh \uS_q(R)$. Since $q\big(-\rho_v(u)\big) = q(u) = a$ and since $b_q(-\rho_v(u), \cdot) = - b_q\big(u, \rho_v(\cdot)\big)$ is a  surjective linear form, Lemma~\ref{lem_smooth_locus-LG}\ref{lem_smooth_locus-LGc} says that the element $-\rho_v(u)$ lies in $\uS_{q,a}\rmsm(R)$. The same argument works for any $T\in \Ralg$, implying that $\wdh \vphi$ is a well-defined morphism of $R$--schemes.

We will first prove smoothness of $\wdh \vphi$ over fields in $\Ralg$. Without loss of generality, we can therefore temporarily assume that $R$ is a field. Let us fix $v\in \wdh  \uV(R)$, thus $q(v) \in R\ti$ and $b_q(u,v) \in R\ti$, put $S=\Spec(R)$, and identify the tangent spaces $T_{\wdh \uV/S}(v) \equiv M$ and 
\[ 
  T_{\uS_{q,a}\rmsm/S}\big(\wdh \vphi(v)\big) \equiv \{ m\in M : b_q(\wdh \vphi(v), m) = 0 \}. 
\]   
Computing with the dual numbers $R[\epsilon]$, the tangent map of $\wdh \vphi$ at $\wtl v$ sends $h\in M$ to  
\begin{eqnarray} \nonumber
&\wdh \vphi(v+ \epsilon h)- \wdh \vphi(v) =
 \Bigl(  \frac{b_q(u,v)}{q(v)} h   +  \frac{b_q(u,h)}{q(v)}v  -
 \frac{b_q(u,v)b_q(v,h)}{q(v)^2}v \Bigr) \epsilon \\  \nonumber
&=
q(v)^{-2}  \, \,   \Bigl(  q(v) \,  b_q(u,v)h +  \big(  q(v) b_q(u,h)- b_q(u,v) b_q(v,h)  \big)v \Bigr) \epsilon.
\end{eqnarray}
The kernel of the linear map $h \mapsto \wdh \vphi(v + \eps h) - \wdh \vphi(v)$ is easily seen to be $Rv$. It is also surjective, being a linear map from an $r$--dimensional vector space to an $(r-1)$--dimensional vector space with a $1$--dimensional kernel. Finally, since $\wdh  \uV$ is smooth, it follows from (d) $\implies$ (b) of \cite[IV$_4$, (17.11.1)]{EGA} that $\wdh \vphi$ is smooth at $v\in V$.

Let now $R$ again be arbitrary. Since $\wdh \uV$ and $\uS_{q,a}\rmsm$ are finitely presented by \eqref{smoLG-a} and \ref{lem_smooth_locus-LG}\ref{lem_smooth_locus-LGe}, and $\wdh \uV$ is flat, the fibrewise smoothness criterion \cite[IV$_4$, 17.8.2]{EGA} together with the field case, proves that $\wdh \vphi$ is smooth in general.    
\sm

The first part of \eqref{smoLG-b} is clear. Also, that $\uW(M)_u/\GG_m  \simlgr \PP(M\ch)$ is well-known, see for example \cite[I, (5.5.3)]{Jan} for a free $M$. Since $\GG_m$ acts freely on the affine scheme $\wdh \uV$, the flat quotient $\wdh \uV/\GG_m$ is represented by an affine scheme $\uV$ and $\wdh  \uV \to \uV$ is a $\GG_m$--torsor by \cite[VIII, Thm.~5.1]{SGA3} (or see \cite[Tag 07S7]{St}). That $\uV$ can be realized as described in \eqref{smoLG-b}, follows from \cite[I, 5.7(3)]{Jan}.


\sm 

In \eqref{smoLG-c} observe that $\wdh \vphi(v) = \wdh \vphi(v')$ for $v,v'\in \wdh  V(R)$ is equivalent to $v' = xv$ for some $x\in R\ti = \GG_m(R)$, so that $\wdh \vphi$ descends to a smooth monomorphism $\vphi \co \uV \to \uS_{q,a}\rmsm$, which is the same as an open immersion \cite[IV$_4$,17.9.1]{EGA}.   \sm 

\eqref{smoLG-d} The first part follows from part \ref{lem_smooth_locus-LGfi} of \ref{lem_smooth_locus-LG}\ref{lem_smooth_locus-LGf}. 
For the proof of $R$-denseness, which by \ref{tod}\eqref{tod-e} is equivalent to universally schematically denseness, it suffices to consider $\uU$ over fields in $\Ralg$. Without loss of generality, let us assume that $R=k$ is a field. Then $\uU$ is an open subscheme of the irreducible $k$--variety $\uS_{q,a}\rmsm$ 
and is therefore Zariski dense as soon as it is non-empty. 
But non-emptiness follows from \eqref{smoLG-e}.
\sm 

\eqref{smoLG-e} By \eqref{smoLG-a}, $\wdh \uV$ is a $k$--variety. It is rational, being an open subscheme of an affine space. Hence $\wdh \uV(k)$ is infinite by \ref{zdens}\eqref{adc}. The same argument works for $\uV(k)$ and $\uU(k)$.  \sm 

\eqref{smoLG-f} We first prove that $\wdh \uV(k)\ne \emptyset$ in case $R=k$ is a field. Because of \eqref{smoLG-e} we can assume that $k$ is finite. Since $\rank_k M \ge 3$, the quadratic from $q$ is isotropic. 
By \eqref{isotrop-d0}, we then know that $(M,q) = H \perp (M',q')$ where $H$ is a hyperbolic plane and where $q'=q|_{M'}$ is nonsingular. By 
\ref{lem_sphere_local}\ref{lem_sphere_local-b}, 
the group $\SO(q)$ acts transitively on $\uS_{q,a}\rmsm(k)$. We can therefore assume that $u=(a,1,0)$. Let $v=(0,1,v')$ for some $v'$ such that $q'(v') \not = 0$; such a $v'$ exists by Lemma~\ref{LGqdi}. Since $q(v) = q(v')$ and since $b_q(u,v) = a$, we have $v\in \wdh \uV(k)$, hence 
$\rho_v(-u)= ( -a,-1 + q(v')\me, q(v')\me v')$ belongs to $\uU(k)$. Thus also $\uU_{q,u}(k)$ is non-empty. 

Let now $R$ be an LG ring. Since $\wdh \uV$ is a principal open, in particular quasi-compact subscheme of the affine space $\uW(M)$, we can apply the criterion 
\ref{prop_baire}\ref{prop_baire-a} and the field case above to conclude that $\wdh \uV(R) \ne \emptyset$. This implies $\uV(R) \ne \emptyset$ and hence also $\uU(R) \ne \emptyset$. \end{proof}

\begin{example}[$\wdh \uV_{q,a}(R) = \emptyset$] \label{smoLG-ex} Let $\HH=R e \oplus Rf$ be the hyperbolic plane over $R=\ZZ/2\ZZ$. In the setting of \ref{smoLG} we then have $a=1$, $u=e + f \in \uS\rmsm_{q,a}(R) = \uS_{q,a}(R) = \{u\}$, but $\wdh \uV_{q,a}(R) = \emptyset$. 
This justifies the assumption $\rank M \ge 3$ in \ref{smoLG}\eqref{smoLG-f}. \end{example}

\begin{lem}\label{knhi} 
Let $(M,q)$ be a quadratic $R$--module, $R$ arbitrary, and let $u, u', w\in M$ satisfying
 \begin{enumerate}[label=\rm (\roman*)] 
 \item\label{knihii} $q(u) = q(u')$, $b_q(u,w) \in R\ti$, $q(w) \in R\ti$, and 
 
 \item\label{knhiii} $v:=\rho_w(u) - u' \in \wS_q$. 
\end{enumerate}  
Then 
\begin{equation}  \label{knihi1}
 \la(u,u',w) := \rho_{\rho_w(u) - u'} \,\rho_w \in \Refl^+(q) 
\end{equation}
sends $u$ to $u'$, and fixes any $m\in M$ satisfying $b_q(m, u-u') = 0 = b_q(m,w)$. If \ref{knihii} holds, then the assumption \ref{knhiii} is fulfilled under any one of the following conditions: 

\begin{enumerate}[label=\rm (\greek*)] 
  \item \label{knhia} $u \equiv u' \mod \fra M$ where $\fra\ideal R$ is an ideal with $\fra \subset \Jac(R)$. 
      
  \item \label{knhib} $q(u-u') = 0$, but $b_q(u',w) \in R\ti$.     
      
      \end{enumerate}
\end{lem}

\begin{proof} This is \cite[(4.1)]{Knebusch} assuming \ref{knhia}.  It is \cite[Lem.~8.2]{EKM} assuming \ref{knhib}. 
\end{proof}
 
\comments{(2025-05) Not clear if we use \ref{knhi}\ref{knhib}. }
\comments{(2025-05) In the previous version we had Prop.~\ref{sotrans}\eqref{sotrans-a} for a semilocal ring, see \ref{prop_quadric2}. But it uses \cite{Knebusch}. 
We then proved \ref{sotrans}\eqref{sotrans-b} as a corollary to transitivity of $\SO(q)$ over semilocal rings. The new version avoids this.}

\begin{prop}[Transitivity of $\uSO(q)$]\label{sotrans}
Let $(M,q)$ be a quadratic space over a ring $R$ satisfying $\rank_R M \ge 3$, and let $a\in R\ti$. \sm

\begin{inparaenum}[\rm (a)] \item\label{sotrans-a} If $R$ is a local ring, then $\SO(q)=\uSO(q)(R)$ acts transitively on the sphere $\uS_{q,a}\rmsm(R)$. \sm  

\item \label{sotrans-b} For any $R$, the action of $\uSO(q)$ on the smooth sphere  $\uS_{q,a}\rmsm$ is transitive on the big affine Zariski site of $\Spec(R)$. 
\end{inparaenum}
\end{prop}

\begin{proof} \eqref{sotrans-a} Let $\m$ be the maximal ideal of $R$. If $R/\m$ has characteristic $\ne 2$ or if $\rank_R M$ is even, then $q$ is regular by \eqref{qf-ns2} 
and \ref{qf}\eqref{quadfoe} respectively. Transitivity then follows from 
transitivity of $\orth(q)$ in \ref{lem_sphere_local}\ref{lem_sphere_local-a}  and Lemma~\ref{lem_transitivity_SO}\ref{lem_transitivity_SO1}. 

It remains to consider the case that $R/\m$ has characteristic $2$ and that $\rank_R M $ is odd. Let $u$ and $u'\in \uS_{q,a}\rmsm(R)$. Again by 
\ref{lem_sphere_local}\ref{lem_sphere_local-a} there exists $f\in \orth(q)$ such that $f(u') = u$. By \eqref{sog-e1} there exists $x\in \mu_2(R)$ and $f' \in \SO(q)$ such that $f=x\,f'$. Hence it suffices to show that there exists $g\in \SO(q)$ satisfying $g(u) = xu$. Because $\mu_2(R/\m) = \{1\}$, we get $x = 1 + y$ for some $y\in \m$, and therefore $xu = u + yu \equiv u \mod \m M$. By \ref{smoLG}\eqref{smoLG-f} 
we can choose $w\in M$ satisfying $b_q(u,w) \in R\ti$ and $q(w)\in R\ti$. We can now use Lemma~\ref{knhi} with case \ref{knhia} and get that $\la(u,xu,w)\in \SO(q)$ maps $u$ to $xu$. \sm 

Part~\eqref{sotrans-b} can be proven in the same way as Corollary~\ref{cor_sphere_homog-orth}. \end{proof}
\sm

In Theorem~\ref{prop_quadric2} we will generalize \ref{sotrans}\eqref{sotrans-a} to semilocal base rings.

\newpage

\section{Transitivity of $\uSO(q)$ on spheres in the semilocal case}\label{sec:trans-semi}

\comments{(2025-05-28) In this new section I combined all results based on Knebusch's paper \cite{Knebusch} which previously where part of section 5 on spheres. Having everything together in one section makes it easier to correct and work on this section. I think, Knebusch's paper has a gap (Satz (4.2) does not seem to be correct: his choice of elements $g_i$ in the proof does not work in general, for example in the Example~\ref{smoLG-ex}). 

For now (2025-09-06) I have eliminated references to this section \ref{sec:trans-semi} as much as possible.
}

\comments{(2025-05) Perhaps \ref{replusLG} and \ref{refplusfield-exa} are better suited here}

\begin{thm} \label{knex-c} 
 Let $(M,q)$ be a faithful quadratic space over a\/ {\em semilocal} ring $R$ and let $\fra \ideal R$ be an ideal of $R$ contained in the Jacobson radical  $\Jac(R)$ of $R$. Then the canonical map $\SO(q) \to \SO(q_{R/\fra})$ is surjective.
\end{thm}

For a regular $q$, this is \cite[Satz~(0.4)]{Knebusch},  proven in \S3 of that paper. The proof is elaborated in \cite[III, Thm.~(3.8)]{Ba} for $\fra = \Jac(R)$. The proof there, although somewhat involved, could be extended to the case of a nonsingular $q$. We leave this to the interested reader. Rather, we prefer to give a more conceptual proof, using the group schemes $\uSO(q)$ and $\uSpin(q)$ which we have reviewed in \S\ref{sec:orthgroup-LG}. 

\comments{(2025-05) In Baeza's notes, Theorem~\ref{knex-c} is stated as \cite[III, Thm.~(3.8)]{Ba} and proven in \cite[III, (3.17)]{Ba}, using the intermediate results (3.9)--(3.16). I checked that all intermediate results hold for nonsingular $q$. They do not require the problematic results on extending isometries, which we used in our previous (wrong) proof of \ref{knex-c}. 

(2026-05-01) I did not see a proof of Theorem~\ref{knex-c} for arbitrary $\fra \subset \Jac(R)$ in \cite{Ba}. \sm 

(2026-05-01) There is no explicit proof of \ref{knex-c} in \cite{Knebusch}. He only says the proof of \ref{knex-c} follows similarly as his proof of \cite[Liftungssatz (0.3)]{Knebusch} using \cite[Lem. (3.1)]{Knebusch}, which says that for $R$ semilocal and $(M,q)$ an orthogonal sum of hyperbolic planes, $\SO(q)$ is generated by Siegel (= Eichler) transformations and products of certain reflections. }

\begin{proof} It is interest to point out that our proof works not only in the semilocal case, but also for an LG ring $R$ for which \eqref{knex-c2} holds. 
\sm 

Using the rank decomposition \eqref{sog2} of $\SO(q)$, we can assume that $M$ has constant rank. Moreover, since $\SO(q)$ is trivial for $\rank M = 1$ by \eqref{sog-e2}, we can further suppose that  $\rank M \ge 2$. Then $\uSO(q)$ is a reductive $R$--group scheme.  

Next, we will reduce to the case of an isotropic $q$. To do so, let $\HH=R u \oplus Rv$ be the hyperbolic plane over $R$, \ref{qf}\eqref{qf-hyp}, and let $(M',q') = (M, q) \perp \HH$. Then $(M',q')$ is a quadratic space by \eqref{qf-perp1}. We put $\ol R = R/\fra$ and claim: 
\begin{equation}\label{knex-c1} \begin{split}
  &\text{$\uSO(q)(R) \to \uSO(q)(R/\fra)$ is onto} 
  \\ &\iff
  \text{$\uSO(q')(R) \to \uSO(q')(R/\fra)$ is onto.} 
\end{split}\end{equation}
To prove \eqref{knex-c1}, let $\la$ be the cocharacter constructed in Lemma~\ref{gmu} for the hyperbolic pair $(U,V) = (Ru, Rv)$,  let $P' = \rmP_{\uSO(q')}(\la)$ be the associated parabolic subgroup and let $L' = \Cent_{\uSO(q')}(\la)$, a Levi subgroup of $P'$. Then $L' = \GG_m \times_R \uSO(q)$ by \eqref{gmu-levi1}, and the group homomorphism induced by $R\to \ol R$ is the direct product homomorphism
\[ L'(R) = R\ti \times \uSO(q)(R) \to {\ol R}\ti \times \uSO(q)(\ol R) = L'(\ol R).  
\]
Since $R \to {\ol R}\ti$ is onto by \ref{inel}\ref{ineliv}, the homomorphism $\uSO(q)(R) \to \uSO(q)(\ol R)$ is onto if and only if the homomorphism $L'(R)\to L'(\ol R)$ is onto, which by Lemma~\ref{lem_onto}\eqref{lem_onto-a} is in turn equivalent to $\uSO(q')(R) \to \uSO(q')(\ol R)$ is onto. 


Having thus established \eqref{knex-c1}, we will in the following assume that $q$ is isotropic. In particular, $\rmD(q) = q(M) \cap R\ti = R\ti$ and then $\rmD^{\rm ev}(q) = R\ti$, which by  \eqref{abc3} and \eqref{abs2} implies that the spinor norm $\SN \co \uSO(q)(R) \to H^1(R, \bmu_2) = R\ti/ R\ti{}^2$ is surjective. Moreover, we can also suppose that $\rank M \ge 4$ so that $\uSO(q)$ and $\uSpin(q)$ are semisimple group schemes. The exact sequence $1 \to \bmu_2 \to \uSpin(q) \to \uSO(q) \to 1$ then induces a commutative diagram of exact sequences of groups 
\begin{equation*} \label{diag_star} \vcenter{\xymatrix{
   \uSpin(q)(R) \ar[r] \ar[d] & \uSO(q)(R) \ar[r]^\SN \ar[d]  & R^\times/ (R^\times)^2 \ar[r] \ar[d] &1 
   \\
   \uSpin(q)(\ol R) \ar[r]  & \uSO(q)(\ol R) \ar[r]^\SN & (\ol R)^\times/ ((\ol R)^\times)^2 \ar[r]  &1.
}}\quad .
\end{equation*} 
Since $R\ti \to {\ol R}\ti$ is surjective, 
the right vertical arrow is surjective. A standard diagram chase then shows that 
the middle vertical arrow is surjective, i.e., our claim is true, as soon as  the left vertical arrow is surjective. Thus, we are left with proving
\begin{equation} \label{knex-c2}
  \Spin(q) \to \Spin(q_{\ol R}) \text{ is surjective.}
\end{equation}
To do so, we use Lemma~\ref{lem_onto}\eqref{lem_onto-b} with the family $(P, P^-)$ where $P$ is a maximal parabolic subgroup of $\uSpin(q)$ and $P^-$ its opposite, e.g., the one given by the cocharacter $\la$ above. Over any field $k$ one knows (\cite[Prop.~6.2]{BT}) that the subgroup $\big\lan P(k) \cup P^-(k) \big\ran$ of $G(k)$ is the subgroup $G^+$ used in the Kneser-Tits problem. However, by \cite[1.1.5 and 1.1.6]{Tits} (or by \cite[Thm.~6.1]{Gille-Bourbaki}) one knows that $G^+ =G(k)$ for $G= \uSpin(q)$. 
\end{proof}

\begin{cor}\label{knexaa} In the setting of\/ {\rm \ref{knex-c}} suppose that $\Refl^+(q) = \SO(q)$. Then also $\Refl^+(q_{R/\fra}) = \SO(q_{R/\fra})$. 
\end{cor}

\begin{proof}
This is an easy diagram chase in 
\[ \vcenter{\xymatrix{ \Refl^+(q) \ar[r] \ar[d] & \Refl^+(q_{R/\fra}) \ar[d]
\\  \SO(q) \ar[r] & \SO(q_{R/\fra})
}  }. \qedhere\]
\end{proof}

\begin{examples} \label{knex-cex} \begin{inparaenum}[(a)] \item The interest in the generality of Theorem~\ref{knex-c} lies in the following application. Let again $R$ be a semilocal ring,  let $S\in \Ralg$ be a finite $R$--algebra  and let $(M,q)$ be a faithful quadratic $S$--space. Then $S$ is semilocal and $\Jac(R)S \subset \Jac(S)$ by \cite[II, (4.2.4)]{K}. Hence \ref{knex-c} with $R$  replaced by $S$ and $\fra=\Jac(R)S$ says that $\SO(q) \twoheadrightarrow \SO(q_{S/\Jac(R)S})$ is onto.\sm 

\item \label{knex-cex-b} 
For special quadratic modules $(M, q)$, surjectivity of $\SO(q) \to \SO(q_{R/\fra})$ holds for more general rings. For example, let $(M,q) = \HH$ be the hyperbolic plane over an arbitrary ring $R$. Then surjectivity holds in view of $\SO(q) \cong R\ti$ by \ref{sog-hyp1}, and surjectivity of $R\ti \to (R/\fra)\ti$ by \ref{inel}\ref{ineliv}.
\end{inparaenum}
\end{examples}

\ms

The analogue of Theorem~\ref{knex-c} for orthogonal groups instead of special orthogonal groups is not true, as pointed out by Knebusch in a remark before \cite[Satz~(0.5)]{Knebusch}: the problem lies in passing from the local to the semilocal case, where $|\orth(q)/\SO(q)|> 2$ may happen. However, in the local case the following holds.

\begin{thm}\label{Knesur} Let $(M,q)$ be a quadratic space over a\/ {\em local} ring $R$ with maximal ideal $\m$ and let $\fra \ideal R$ be an ideal contained in $\m$. Assume that $q$ is regular or that $\fra = \m$. Then the canonical map $\orth(q) \to \orth(q_{R/\fra})$ is surjective.
\end{thm}

\begin{proof}
In case $q$ is regular, this is the special case $E=F$ of \cite[(0.3)]{Knebusch}.
If $q$ is not regular, then $M$ has constant odd rank by  \ref{qf}\eqref{quadfoe} and $2\notin R\ti$ by \eqref{qf-ns2}\emph{}, so that $R/\m$ has characteristic $2$. Hence
$\orth(q_{R/\m}) = \SO(q_{R/\m})$ by \eqref{sog-e1},  
and surjectivity follows from  surjectivity of $\SO(q) \to \SO(q_{R/\fra})$ (Theorem~\ref{knex-c}) and the commutative diagram 
\[\vcenter{ \xymatrix@C=50pt{ \SO(q) \ar[r] \ar[d] & \SO(q_{R/\fra}) \ar@{=}[d] \\
             \orth(q) \ar[r] & \orth (q_{R/\fra}) }
} . \qedhere\]
\end{proof}

Theorem~\ref{Knesur} is not true in case $q$ is not regular and $\fra\subsetneq \m$. Indeed, as in the proof above, surjectivity of $\orth(q) \to \orth(q_{R/\fra})$ reduces to surjectivity of $\mu_2(R) \to \mu_2(R/\fra)$. Here is an example where the latter map is not surjective: let $k$ be a field of characteristic $2$, let $k[X]$ be the polynomial ring over $k$ and let $R=k[X]/(X^3)$. Put $t = X + (X^3)\in R$, and let $\fra = t^2R$. Then $1+t\in \mu_2(R/\fra)$
is not in the image of $\mu_2(R) \to \mu_2(R/\fra)$.
\lv{
Details:
 since any candidate in $R\ti$ has the form $x=a + b t + c t^2$ with $a,b,c\in k$, hence $(a + bt + c t^2)^2 = a^2 + b^2 t^2$, so that $a^2=1$, $b^2 = 0$ and therefore $a=1$, $b=0$ holds if $x^2 = 1$. Thus $x = 1 + ct^2$, which maps to $1$ in $R\fra$. }


\comments{(2025-05) It would be nice to have Thm.~\ref{prop_quadric2} for LG rings, in particular since \ref{smoLG}\eqref{smoLG-f} holds for LG rings. One could for example expect a theory for a smaller $\wdh \uV_{q,a}$ satisfying \ref{knhi}\ref{knhib}. I have not thought about this.}

Using Theorem~\ref{knex-c} we can extend Proposition~\ref{sotrans}\eqref{sotrans-a} from local to semilocal base rings. 

\begin{thm}[Transitivity in the semilocal case]  \label{prop_quadric2} Let $R$ be a semilocal ring, let $(M,q)$  be a quadratic $R$--space with $\rank M \ge 3$, and let $a\in R\ti$. Then the group $\uSO(q)(R)$ acts transitively on $\uS_{q,a}\rmsm (R)$.
\end{thm}

\begin{proof} 
Let $\Max$ be the set of maximal ideals of $R$, and let $\m \in \Max$. We can apply Lemma~\ref{lem_sphere_local}\ref{lem_sphere_local-b} to the quadratic space $(M,q)_{R/\m}$ over the field $R/\m$, and get $g_{[\m]} \in \SO(q_{R/\m})$ such that   $g_{[\m]}\cdot u_{R/\m} = u'_{R/\m}$. The family $(g_{[\m]})_{\m \in \Max}$ yields $\ol g \in \uSO(q)(R/\Jac(R))$ satisfying $\ol g \cdot u_{R/\Jac(R)} = u'_{R/\Jac(R)}$. By \ref{knex-c} we can lift $\ol g$ to a $g\in \SO(q)$. Hence $g \cdot u \equiv u' \mod \Jac(R)M$. Then, by \ref{smoLG}\eqref{smoLG-f}, there exists $w\in M$ satisfying the condition \ref{knhi}\ref{knihii} for $g \cdot u$ and $u'$. Hence $\la(g\cdot u, u',w))(g \cdot  u) = u'$, and we are done because $\Refl^+(q) \subset \SO(q)$ by \eqref{refso2}. \end{proof}
\enew 
\lv{
In view of the standard rank decomposition of $(M,q)$, $\uO(q)$ and $\uSO(q)$, see \ref{qf}\eqref{qf-redc}, \eqref{ortgr-bas-b1} and \eqref{sogsc1}, we can assume that $M$ has constant rank $r \ge 1$, in particular $M$ is free. Let $u, u' \in \uS_{q,a}\rmsm(R)$. \sm

In the following we will assume $r\ge 2$ and prove transitivity under $\SO(q)$. We first deal with the case that $q$ is regular. Both $u$ and $u'$ are unimodular by Lemma~\ref{mx_lem}\eqref{mx_lema}. We thus have a well-defined isometry $Ru \to Ru'$, $u \mapsto u'$. We claim that it can be extended to an isometry of $M$. Indeed, $q|_{M_u}$, $M_u = (Ru)^\perp$, is nonsingular by Lemma~\ref{mx_lem}\eqref{mx_lemc}. Hence the ideal of $R$ generated by $q(M_u)$ is $R$ by \ref{quadco}. We can therefore apply Knebusch's Extension Theorem (\cite[Satz~5.2]{Knebusch}, \cite[III, Thm.~4.1]{Ba}), proving $u'\in \uO(q)(R)\cdot u = \orth(q) \cdot u$. Finally, Lemma~\ref{lem_transitivity_SO}\ref{lem_transitivity_SO1} shows $u'\in \uSO(q)(R)\cdot u$.}
\sm   

The following Proposition~\ref{prop_quadric1} is another application of Theorem~\ref{knex-c}.  


\begin{prop} \label{prop_quadric1} Let $(M,q)$ be a quadratic space over the semilocal ring $R$ with $\rank M \ge 3$, let $a\in R\ti$, and let $\fra \subset \Jac(R)$ be an ideal of $R$. Then the map  $\uS_{q,a}^{\rm sm}(R) \to \uS_{q,a}^{\rm sm}(R/\fra)$ is onto.  
\end{prop}

\begin{proof} By \ref{smoLG}\eqref{smoLG-f}, there exists $v \in \uU_{q, u}(R)$. We want to lift a given element $x\in \uS_{q,a}^{\rm sm}(R/\fra)$ to $\uS_{q,a}^{\rm sm}(R)$. Recall from \eqref{slr-ai1} that
\[ R/\fra \twoheadrightarrow (R/\fra)\big/ \Jac(R/\fra) \cong R/\Jac(R) =
    \textstyle \prod_{\m \in \Max} R/\m .
\]
Pushing $x$ down to $x_{R/\Jac(R)}= (x_{R/\m})_{\m \in \Max}$, the argument used in the proof of \ref{prop_quadric2} 
shows that there exists $(g[\m])_{\m \in \Max} \in \uSO(q)(R/\Jac(R))$ satisfying $g[\m] \cdot x_{R/\m} = v_{R/\m}$ for all $\m \in \Max$. Let $g\in \uSO(q)(R)$ be the lift of the family $(g[\m])_{\m \in \Max}$, which exists by \ref{knex-c}. Since then $(g \cdot x)_{R/\m} = v_{R/\m}$ and since $v\in \uS_{q,a}(R)$, we can without loss of generality assume that $x_{R/\m} = v_{R/\m}$ holds for all $\m \in \Max$. Now Lemma~\ref{opemaLG}, applied to the open subscheme $\uU_{q,u}$ of $\uS_{q,a}^{\rm sm}$ and the ring $R/\fra$, shows that $x\in \uU_{q,a}(R/\fra)$. Hence $x=\vphi([\ol m])$ for some unimodular $\ol m \in M\ot_R (R/\fra) = M/\fra M$ with $q_{R/\fra}(\ol m)\in (R/\fra)\ti$. By \eqref{unimod-2} we can lift $\ol m$ to an unimodular $m\in M$ satisfying $q(m) \in R\ti$ because $q(m)_{R/\fra} = q_{R/\fra}(\ol m)$. It defines a point $[m] \in \big(\PP(M^\vee)\setminus \{q=0\}\big)(R)$. Since $\vphi([m])_\m = v_\m$ for all $m\in \Max$, a second application of Lemma~\ref{opemaLG} shows that $\vphi([m]) \in \uU_{q,a}^{\rm sm}(R)$. We have thus lifted $x= \vphi([\ol m])$ to $\uU_{q,a}^{\rm sm}(R) \subset \uS_{q,a}^{\rm sm}(R)$.
\end{proof}

%
%

The following Lemma~\ref{acsoLG} is yet another application of Theorem~\ref{knex-c}. 

\begin{lem}\label{acsoLG}
  Let $R$ be a semilocal ring, let $(M,q)$ be a quadratic $R$--space of constant rank $\ge 3$, let $u\in \uS_{q,a}(R)$ and let $v\in \uS_{q,z}\rmsm(R)$. Then there exists $g\in \uSO(q)(R)$ such that $Ru + R(g\cdot v)$ is a complemented submodule of $M$ which is free of rank $2$ with basis $(u, g \cdot v)$.  
\end{lem}

\begin{proof} We first prove the lemma in case $R=k$ is a field. We can assume that $v\in ku$. Since $|ku \cap \uS_{q,z}\rmsm(k)| \le 2$ while $|\uS_{q,z}\rmsm(k)|\ge 3$ by Example~\ref{lem_smooth_locus_exam}\eqref{lem_smooth_locus_isof}, there exists $v'\in \uS\rmsm_{q,z}(k) \setminus ku$. Then transitivity of $\uSO(q)(k)$ on $\uS_{q,z}\rmsm(k)$, i.e., Proposition~\ref{sotrans}\eqref{sotrans-a} or even Theorem~\ref{prop_quadric2}, finishes the proof in the field case. 

Let now $R$ be arbitrary semilocal, let $\ol R = R/\Jac(R)$, and let $\Max$ be the finite set of maximal ideals of $R$. For $\gm \in \Max$ let $u_{R/\gm} \in M_{R/\gm} = M / \gm M$ be the canonical image of $u$ and use the analogous notation for $v$. By the first part of the proof there exists a family 
\[ 
(g[\m])_{\m \in \Max} \in \textstyle\prod_{\m \in \Max}\, \uSO(q)(R/\m) =    \uSO(q)(R/\Jac(R))
\]
such that 
\[
   (u', v') = \big(\textstyle  \prod_{\gm \in \Max} u_{R/\gm}, \, \prod_{\gm \in \Max} g[\gm]\cdot v_{R/\gm} \big) \in \ol M = M/\Jac(R) M 
\] 
is a basis of the  rank-$2$-submodule of the $\ol R$--module $\ol M$ spanned by $u'$ and $v'$. It is complemented by the free $\ol R$--submodule $N' = \prod_{\gm \in \Max} N[\gm]$ where $N[\gm]$ is a complement of the $(R/\gm)$--submodule spanned by $u_{R/\gm}$ and $g[\gm]\cdot v_{R/\gm}$. By \eqref{sog2} and Theorem~\ref{knex-c}, we can lift the family  $(g[\gm])_{\gm \in \Max}$ to some $g\in \uSO(q)(R)$. Under the canonical map $M \to \ol M$ the pair $(u, g\cdot v)$ maps onto $(u',v')$. It is therefore $R$--free. We lift the free $\ol R$--module $N'$ to an $R$--submodule $N$ of $M$. Since $N$ maps onto $N'$, it follows from \ref{nak}\eqref{nak-dec} that $N$ is a complement of the $R$--submodule spanned by $u$ and $g \cdot v$. 
\end{proof}
\sm 

The plane $Ru \oplus R(g\cdot v)$ in Lemma~\ref{acsoLG} need not be regular. This is however the case in some special situations, like in the following Lemma~\ref{actsoLG} and Proposition~\ref{plar}. 

\comments{(2025-06-06) I formulated \ref{actsoLG} in this way, hoping that transitivity in \ref{prop_quadric2} can be shown in a more general setting}

\begin{lem}\label{actsoLG} Let $(M,q)$ be an\/ {\em isotropic} quadratic space over a ring $R$ for which $\SO(q)$ acts transitively on smooth spheres $\uS_{q,z}\rmsm(R)$ for any $z\in R\ti$, e.g., assume $R$ is semilocal and $\rank_R M \ge 3$, see {\rm \ref{prop_quadric2}}. Furthermore, suppose there exists $u \in \uS_{q,a}\rmsm$ for some $a \in R\ti$. Then, for $v\in M$, there exists $g\in \SO(q)$ such that 
\[ Ru \oplus R ( g \cdot v) \]
is a hyperbolic plane, if \ref{actsoLGi} or \ref{actsoLGii} below holds:  
\begin{enumerate}[label={\rm (\roman*)}]
 \item \label{actsoLGi} $v$ is isotropic; 
 
 \item \label{actsoLGii} $R$ is an LG ring with $|R/\m|\ge 4$ for every maximal ideal $\m \ideal R$, and $v\in \uS_{q, q(v)}\rmsm$, $q(v) \in R\ti$. 
\end{enumerate} 
\end{lem}

\begin{proof} \ref{actsoLGi} By \ref{isotrop}\eqref{isotrop-d} there exists $w\in M$ such that $(v,w)$ is a hyperbolic pair. Let $H = Ru \oplus Rw$. Then $u'=v + aw \in \uS_{q|H,a} = \uS_{q|H, a}\rmsm(R) \subset \uS_{q,a}\rmsm(R)$ by \ref{lem_smooth_locus-LG}\ref{lem_smooth_locus-LGd} and \ref{lem_smooth_locus-LG}\ref{lem_smooth_locus-LGc}. Choose $g\in \SO(q)$ such that $g\cdot u = u'$. Then $Rv \oplus Rw = Rv \oplus Ru' = R (g \cdot u) \oplus Rv = g \big( Ru \oplus R (g \cdot v)\big)$. \sm     

\ref{actsoLGii} Since $(M,q)$ is isotropic, there exists a hyperbolic pair $(e,f) \in M$. Then $u' = ae + f \in \uS_{q,a}\rmsm(R)$ and by transitivity $gu' = u$ for some $g\in \SO(q)$, so that $u \in R(g \cdot e) \oplus R(g \cdot f)$, Thus, without loss of generality, we can assume that $u=u' \in H$, the hyperbolic plane spanned by a hyperbolic pair $(e,f)$. Observe that $v' = e + q(v) f \in \uS_{q, q(v)}\rmsm(R)$. Hence, again by transitivity, we can assume that $v=v'$. By \eqref{sog-hyp1}, for any $x\in R\ti$, the vector $v_x = xe + x\me q(v) f \in \SO(q) \cdot v$, and $Ru \oplus R v_x = H$ if and only if $x^2 - q(v)a \in R\ti$.   
Let $f$ be the polynomial $f(X,Y) = (X^2 - q(v)a) Y \in R[X,Y]$. Since $|R/\m| \ge 4$ for every maximal ideeal $\m \ideal R$, the polynomial $f$ has a solution over any $R/\m$. Because $R$ is LG, $f$ has a solution in $R$, finishing the proof. \end{proof}

\begin{prop}\label{plar} Let $R$ be a semilocal ring for which all residue fields $R/\gm$, $\gm\ideal R$ maximal, are infinite, thus a semilocal ring satisfying the primitive condition. Let $(M,q)$ be a quadratic $R$--space with $\rank M \ge 3$, and suppose that there exist $u\in \uS_{q,a}\rmsm(R)$ and $v\in \uS_{q,z}\rmsm(R)$. \sm 

\begin{inparaenum}[\rm (a)] \item\label{plar-a}
Then there exists $g\in \uSO(q)(R)$ such that  
\[ P= Ru \oplus R(g\cdot v)
\]
is a regular plane. \sm 

\item\label{plar-b} There exists $w\in M$ such that $W=P \oplus Rw$ is a complemented nonsingular submodule of constant rank $3$ with basis $\{ u, g\cdot v, w\}$. 
\end{inparaenum}
\end{prop}

\begin{proof} \eqref{plar-a} Without loss of generality, we can assume that $M$ has constant rank. In parts \eqref{plar-I}, \eqref{plar-II} and \eqref{plar-III} below we will prove the proposition in case that $R=k$ is an infinite field. \sm 
  
\begin{inparaenum}[(I)] \item\label{plar-I} We claim that {\em there exists $v'\in \uS_{q,z}\rmsm(k)$ for which $b(u,v') \ne 0$}. Assume otherwise, i.e.,  $\uS_{q,z}\rmsm(k) \subset M_u = (ku)^\perp$. Hence, $\uS_{q,z}\rmsm\subset \uW(M_u)$ by Zariski-density of $\uS\rmsm_{q,z}(k)$ in $\uS_{q,z}$ and then in $\uW(M_u)$. We know $\dim M_u = \dim_k M - 1 = \dim_k \uS_{q,z}$ by \ref{bfLG}\eqref{mx_lemb} and \ref{lem_smooth_locus-LG}\ref{lem_smooth_locus-LGf} respectively. As $\uS_{q,z}$ is a closed subscheme of the irreducible scheme $\uW(M_u)$, we get the contradiction $\uS_{q,z} = \uW(M_u)$.%
\sm   

\item \label{plar-II} We continue with the assumption that $R=k$ is a field. Using transitivity of $\uSO(q)(k)$ on $\uS_{q,z}\rmsm(k)$, we can by \eqref{plar-I} assume $b(u,v) \ne 0$ in the following. Our next claim is that {\em $m \mapsto b_q(m,u)$ is not constant on $\uS_{q,z}\rmsm(k)$.} To prove this, we consider the morphism of schemes 
\[ \psi \co \uS_{q,z}\rmsm \to \GG_a, \quad m \mapsto b_q(m,u)
\]    
which we compose with the morphism
\[ \wdh{\vphi} \co \wdh \uS_q \to \uS_{q,z}\rmsm, \quad x \mapsto \wdh \vphi(m) = - \rho_x(v) = -v + b_q(v,x)\, q(x)\me \, x
\]  
of \ref{smoLG}\eqref{smoLG-a}. Assume that $\psi$ is constant on $\uS_{q,z}\rmsm(k)$. Then so is $\psi \circ \wdh \vphi$. Since $\psi \circ \wdh \vphi(v) = b_q(u,v)$ we get, using again Zariski-density of $\uS_{q,z}\rmsm(k)$ in $\uW(M)$, that 
$b_q(u,v)  = \psi \circ \wdh \vphi (x)  
      = b_q(u,v) -  q(x)\me b_q(u,x) b_q(v,x)$,
equivalently, 
\begin{equation}\label{plar1} 
 2 q(x) b_q(u,v) = b_q(u,x)\, b_q(v,x) 
 \end{equation}  
holds for all $x\in M$.
If $\Char(k) \ne 2$, then \eqref{plar1} says that $q$ is a product of two non-zero linear forms, contradicting the assumptions $q$ nonsingular and $\dim_k M \ge 3$. On the other side, if $\Char(k) = 2$, then \eqref{plar1} says that $M$ is the union of two hyperplanes, 
again a contradiction. 
Thus, $\psi \circ \wdh \vphi$ is not constant on $k$. So $\psi$ is not constant either. 
\sm 

\item\label{plar-III}
We now consider the restriction of $b_q$ to $Ru + Rv$. The determinant of the matrix representing $b$ in the spanning set $\{u,v\}$ is 
\[ \det \begin{pmatrix} 2a & b(u,v) \\ b(u,v) & 2z \end{pmatrix} = 4az - b(u,v)^2. 
\]
Since $\psi$ is not constant, we can, after possibly replacing $v$ by some $v''\in \uS_{q,z}\rmsm(k)$ and using again transitivity of $\uSO(q)(k)$, without loss of generality assume that $4az - b(u,v)^2 \ne 0$. This implies that $\{u,v\}$ is a basis of $P[k]=ku \oplus kv$ and that  $q|_{P[k]}$ is regular. \sm 

\item\label{plar-IV} It remains to deal with a general semilocal $R$. We proceed as in the proof of Theorem~\ref{prop_quadric2} and Lemma~\ref{acsoLG}. For every $\gm \in \Max$, the set of maximal ideals of $R$, we get $g[\gm]\in \uSO(q)(R/\gm)$ such that $u_{R/\gm}$ and $g[\gm] \cdot v_{R/\gm}$ are a basis of a regular plane in $M_{R/\gm}$. Let $\ol R = R/\Jac(R)$. Then $g' = (g[\gm]) = \prod_{\gm \in \Max} \uSO(q)(R/\gm) = \uSO(q)(\ol R)$ has the property that $u_{\ol R}$ and $g' \cdot v_{\ol R}$ is the basis of a regular plane in $\ol M = M_{\ol R}$. By Theorem~\ref{knex-c} we can lift $g'$ to an element $g\in \uSO(q)(R)$. Since the canonical images of $u$ and $g \cdot v$ in $\ol M$ are a basis of a regular plane, the same holds for $u$ and $g \cdot v$, thus finishing the proof of \eqref{plar-a}. 
\end{inparaenum}  
\sm 

\eqref{plar-b} The orthogonal complement $P^\perp$ is a faithfully projective $R$--module for which $q|_{P^\perp}$ is nonsingular. By Lemma~\ref{LGqdi}, there exists $w\in P^\perp$ with $q(w) \in R\ti$. Such a $w$ is unimodular in $P^\perp$ by Lemma~\ref{mx_lem}, hence $Rw \subset P^\perp$ is free of rank $1$ and complemented by \ref{unimod}. This implies \eqref{plar-b}.  
\end{proof}

\sm 

The final topic of this section is the question under which conditions do we have $\Refl^(q) = \SO(q)$ where, we recall, $\Refl^(q)$ is the subgroup of $\orth(q)$ generated by an even number of reflections, \ref{refle}, and where $\SO(q) = \uSO(q)(R)$ which always contains $\Refl^(q)$ by \eqref{refso2}. We have dealt with the field case in Lemma~\ref{replusLG} and the Example~\ref{refplusfield-exa}. This will now be extended to semilocal rings in Theorem~\ref{Bazth} and Corollary~\ref{cor-knex}, using the following results of Knebusch. 

\begin{thm}\label{Bazth}
Let $R$ be a semilocal ring, let $(M,q)$ be a quadratic $R$--module and let
$\ga \co (N, q|_N) \to (M,q)$ be an isometric embedding with $\ga(n) - n \in \Jac(R)M$ for all $n\in N$.   Then there exists $g\in \Refl^+(q)$ such that $g|_N = \ga$ under any one of the following two conditions: \sm

\begin{enumerate}[label={\rm (\roman*)}]
  \item \label{Bazthi} {\rm (\cite[Satz~(0.1c)]{Knebusch})} $q|_N$ is regular and $q(N^\perp)$ generates $R$ as ideal, or \sm 
      
  \item\label{Bazthii} {\rm (\cite[Satz~(4.2)]{Knebusch}}, see also {\rm \cite[III, Thm.~3.19]{Ba})} $q$ is regular and $N$ is complemented.     
\end{enumerate}
\end{thm}
\sm 

We point out that the condition $\lan q(N^\perp) \ran = R$ in \ref{Bazthii} above is fulfilled in case $q$ is nonsingular and $N^\perp$ is faithfully projective, see \ref{LGqdi}.  

For a regular $q$, part \ref{knex-a} of the following theorem is proven in \cite[Satz (4.3)]{Knebusch}.

\begin{thm}\label{knex} Let $(M,q)$ be a quadratic space over a\/ {\em semilocal} ring $R$. \sm

\begin{enumerate}[label={\rm (\alph*)}]
\item\label{knex-a} Then $\Refl^+(q) = \SO(q) \iff \Refl^+(q_{R/\m}) = \SO(q_{R/\m})$ holds for all maximal ideals $\m \ideal R$.  \sm

\item\label{knex-b}  We have $\Refl^+(q) = \SO(q)$ whenever

 \begin{enumerate}[label={\rm (\roman*)}]
    \item\label{knex-bi} $M$ has constant rank $\ne 4$, or

    \item $|R/\m| \ge 3$ for all maximal ideals $\m \ideal R$.
\end{enumerate} 
\end{enumerate}
\end{thm}

\begin{proof}
  We will prove this theorem in steps \eqref{Knex-I}--\eqref{Knex-V} below, but not follow the order in which the claims are stated. \sm

\begin{inparaenum}[(I)]
  \item\label{Knex-I} We suppose $R=R_1 \times \cdots \times R_n$ is a direct product of rings and let $q=q_1 \times \cdots \times q_n$ be the corresponding decomposition of $q$, cf.\ \ref{qf}\eqref{qf-redc}. We know that $q$ is nonsingular (regular) if and only if every $q_i$, $i=1, \ldots, n$, is  nonsingular (regular respectively). It follows from \eqref{refle5} and \eqref{sog2} that
      \[ \Refl^+(q)= \SO(q)  \iff \Refl^+(q_i) = \SO(q_i) \text{ for $i=1, \ldots, n$}. \]
In particular, we can apply the above to the rank decomposition of $(M,q)$ and in this way reduce the proof of $\Refl^+(q) = \SO(q)$ to the case of $M$ having constant rank. \sm

\item\label{Knex-II} Proof of $\Longleftarrow$ in \ref{knex-a}: Applying the rank decomposition we can assume that $M$ has constant rank. We put $\ol R = R/\Jac(R)$ and $(M,q)_{\bar R} = (\ol M, \ol q)$. Then $\Refl^+(\ol q) = \SO(\ol q)$ by \eqref{Knex-I}. Recall $\Refl^+(q) \subset \SO(q)$ by \eqref{refso2}. For the proof of the other inclusion, let $g\in \SO(q)$. Then $\ol g = g \ot \Id_{\ol R} \in \SO(\ol q)$ by functoriality \ref{sog}\eqref{sog-ba}. Therefore $\ol g \in \Refl^+(\ol q)$. Since reflections can be lifted to $(M,q)$, there exists $f\in \Refl^+(q)$ such that $\ol g = \ol f$, i.e., $(f\me g)(m) \equiv m \mod \Jac(R) M$.

    If the rank of $M$ is odd, we use \ref{nqf-LG} and write $(M,q) = (Re, \lan u \ran_q) \perp (N,q|_N)$ with $N$ free of even rank, $q|_N$ regular and $u\in R\ti$. By Theorem~\ref{Bazth}\ref{Bazthi} there exists $h\in \Refl^+(q)$ such that $h|_{M'} = (f\me g)|_{M'}$. It follows that $h\me f\me g\in \SO(q)$ stabilizes $(M')^\perp = Re$, thus $(h\me f\me g)(e) = xe$ for some $x\in R\ti$. But $x= \det (h\me f\me g) = 1$ by \ref{sog}\eqref{sog-b}. Hence $g=fh \in \Refl^+(q)$.

    The proof in case $M$ has even rank can be done along the same lines, applying \ref{Bazth}\ref{Bazthii} with $M=N$. \sm 


\item\label{Knex-III} Proof of \ref{knex-b}: For every maximal ideal $\m$ of $R$ we get $\Refl^+(q_{R/\m}) = \SO(q_{R/\m})$ from Lemma~\ref{replusLG}. Therefore \ref{knex-b} follows from \eqref{Knex-II}.
\sm

\item\label{Knex-V} Proof of $\Longrightarrow$ in \ref{knex-a}: 
This follows from \eqref{Knex-I} and \ref{knexaa}. 
\end{inparaenum}
\end{proof}
\sm

\comments{(2022-09-29) First proves in \cite[Thm.~2]{Fir21}: Let $R$ be a commutative semilocal ring with $2\in R\ti$, let $A$ be an Azumaya $R$--algebra and let $\si \co A \to A$ be an orthogonal involution. Then the natural map
\[ \SO(A, \si) \to \textstyle \prod_{\m \in {\rm Max}(R)} \SO(A(\m), \si(\m))
\]
is surjective, where $\orth(A, \si) = \{ u \in A: u \si(u) = 1_A\}$ and
$\SO(A, \si) = \Ker (\Nrd)$ for $\Nrd \co A \to R$ is the reduced norm of $A$.

For $A=\End_R(M)$ and $\si$ the adjoint involution, this follows from our result \ref{knex}\ref{knex-c}; we are more general since we do not assume $2\in R\ti$. On the other hand, First is more general since he works with arbitrary Azumaya $R$--algebras}

\begin{cor} \label{cor-knex} Let $(M,q)$ be a faithful quadratic space over the semilocal ring $R$. Then there exists an \'etale extension $T\in \Ralg$ of constant odd degree with the following properties:
  \begin{enumerate}[label={\rm (\roman*)}]
    \item \label{cor-kneb-i} all residue fields of $T$ have cardinality $\ge 3$, and
    \item \label{cor-kneb-ii} $\Refl^+(q_T) = \SO(q_T)$.
  \end{enumerate}
\end{cor}

\begin{proof}
We first use Proposition~\ref{bfp-prop} with $S=R$ to get the existence of an \'etale $T\in \Ralg$ of constant odd rank such that \ref{cor-kneb-i} holds. We can then apply Theorem~\ref{knex}\ref{knex-b} to get \ref{cor-kneb-ii}. \end{proof}

\comments{There are (at least) two papers that prove Witt's Extension Theorem for certain noncommutative semilocal rings. Perhaps they could be used here: 

First, Uriya A., {\em  Witt's extension theorem for quadratic spaces over semiperfect rings} J. Pure Appl. Algebra 219 (2015), no. 12, 5673--5696. and \sm 

H. Reiter, {\em Witt's theorem for noncommutative semilocal rings}, J. Algebra 35 (1975) 483--499. }

\newpage

\section{The $\uSO(q)$ action on higher rank quadrics}\label{sec:quadrics}

In this section we study the $R$--scheme $\uQ_\nu(q)$ of higher rank quadrics whose points are totally isotropic, complemented submodules of a quadratic space $(M,q)$. We show that such a quadric is always smooth projective, \ref{quadint} and \ref{quadsm}, and prove in \ref{prop_hig} that the stabilizer of an $R$--point is a parabolic subgroup $P$ of the reductive $R$--group $\uSO(q)$. We identify the type of $P$ in \ref{Colem} and \ref{prohigqu}. We also show in \ref{prop_hig} that $\uSO(q)/P \cong \uQ_\nu(q)$ in case the rank function $\nu$ of the points of $\uQ_\nu(q)$ is not maximal. The latter case, referred to as Lagrangians $\uL(q)$, requires a special treatment which we present in the second part of this section. Specifically, we identify the Stein factorization of the structural morphism $\uL(q) \to \Spec(R)$. While this is known, even over an arbitrary base (Deligne \cite[XII, Prop.~2.8]{SGA7}), we present a detailed version of Deligne's result in the affine setting in Proposition~\ref{deligne}. This will allow us to describe the $\uSO(q)$--action on $\uL(q)$ in \ref{eor} and \ref{prohla}. \ms

\begin{lem}\label{quadint} Let $(M,q)$ be a faithful quadratic module and let  $\nu \co \Spec(R) \to \NN$ be a locally constant function. Recall $\nu_{R'} = \nu \circ (\Spec(R') \to \Spec(R))$ for $R'\in \Ralg$. Following {\rm \ref{grap}\eqref{grape-dual}}, we identify $\big({\ul \Gr}_\nu(M\ch)\big)(R')$  with the set of complemented, hence finite projective $R'$--submodules $N \subset M_{R'}$ of rank $\nu_{R'}$. We claim that the $R$--subfunctor  ${\underline{\rmQ}}_\nu(q)$ of ${\ul \Gr}_\nu(M\ch)$ whose $R'$--points are defined as
\[
\underline{\rmQ}_{\nu}(q)(R') = \{  N \in {\ul \Gr}_\nu(M\ch)(R') :
q_{R'}(N)=0 \}
\]
is a closed subfunctor of $\ulGr_\nu(q)$, representable by a closed $R$--subscheme $\uQ_\nu(q)$ of the Grassmannian scheme $\uGr_\nu(M)$. Furthermore, $\uQ_\nu(q)$ is projective and of finite presentation over $R$. \end{lem}

\begin{proof} In view of \cite[I, \S2, 4.1]{DG}, 
representability of $\ul \rmQ = {\ul \rmQ}_\nu(q)$ follows, once we have shown that ${\ul \rmQ}$ is a closed subfunctor of ${\ul \rmG} = {\ul \Gr}_\nu(M\ch)$. To see this, 
let $R'\in \Ralg$, denote by $\Sp(R')$ the $R$--functor associated with the affine scheme $\Spec(R')$ and let $f \co \Sp(R') \to {\ul \rmG}$ be a morphism of $R$--functors. We need to show that $f\me({\ul \rmQ})$ is a closed subfunctor of $\Sp(R')$. By Yoneda, $f$ is determined by a unique $N\in {\ul \rmG}(R')$ in the following way: to $A \in \Ralg$ corresponds the map
\[
 f(A)\co \Sp(R')(A) = \Hom_{\Ralg}(R', A) \to {\ul \rmG}(A), \quad \al \mapsto
 N \ot_{R'} A\]
where $A$ is viewed as $R'$--module via $\al$. Hence
\begin{align*} f\me({\ul \rmQ})(A) &=
  \{ \al \in \Hom_{\Ralg}(R', A): q_A(N \ot_{R'} A) = 0 \}
\\ &= \{\al \in \Hom_{\Ralg}(R', A): \al\big(q_{R'}(N)\big)A = 0 \}.
\end{align*}
Thus, $f\me ( {\ul \rmQ})$ is the closed subfunctor of $\Sp(R')$ associated with the ideal $\lan q_{R'}(N)\ran \ideal R'$.

As a closed subscheme of the projective $R$--scheme $\uGr_\nu(M\ch)$, the scheme $\uQ_\nu(q)$ is projective, hence proper
and therefore separated and quasi-compact. It is then of finite presentation as soon as it is locally of finite presentation.  But this is clear since the functor ${\ul \rmQ}$ commutes with direct limits of rings 
\ref{ag}\eqref{ag-lp}. 
\end{proof}

\subsection{Higher rank quadrics} \label{hrq} To avoid degeneracies in rank $1$, we assume that $\rank_R M \ge 2$ as locally constant function, so $\rank_R M$ is  not necessarily constant.

We will refer to the scheme $\uQ_\nu(q)$ of \ref{quadint} as the {\em quadric of rank $\nu$} or as a {\em higher rank quadric\/} if the rank is unimportant. If $\nu$ is the constant function assigning the value $d$ to all $\p \in \Spec(R)$ we write $\uQ_d(Q)$ for $\uQ_\nu(q)$. Two special cases will have their own name and notation: 

\begin{enumerate}[label={\rm (\roman*)}]
  \item If $\nu$ is the constant function $1$, we abbreviate $\uQ(q)=\uQ_1(q)$ and call it the {\em quadric associated with $q$}. 
      
  \item If $\rank_R M = 2\nu$ we abbreviate $\ulL(q) = \underline{Q}_\nu(q)$ and $\uL(q) = \uQ_\nu(q)$, and call $\ulL(q)$ the {\em Lagrangian functor} and $\uL(q)$ the {\em Lagrangian quadric}, cf.\ \ref{Lag}.
\end{enumerate}

The group scheme $\uO(q)$ acts on the higher rank quadrics $\uQ_\nu(q)$ in the obvious way. Using the terminology of \ref{transi}, we note:
{\em The action of $\uO(q)$ on the Lagrangian quadric $\uL(q)$ is transitive on the big affine Zariski site of $\Spec(R)$.} Indeed, given $R'\in \Ralg$ and $L$, $L' \in \uL(q)(R')$, there exists a Zariski cover $R''$ of $R'$ such that $L_{R''}$ and $L'_{R''}$ are free, hence are isomorphic $R''$--modules. By \ref{Lag}\eqref{Lag-b}, there exists an isometry of $(M,q)_{R''}$ mapping $L_{R''}$ onto $L'_{R''}$.

%

\begin{thm}\label{quadsm}  Let $(M,q)$ be a quadratic $R$--module with $\rank_R M \ge 2$. Then
\begin{equation}\label{quadsm1}
\text{$q$ is nonsingular}\quad \iff \quad \text{$\uQ(q)$ is a smooth $R$--scheme}. \end{equation}
In this case, every higher rank quadric $\uQ_\nu(q)$ is a smooth projective $R$--scheme.
\end{thm}

\begin{proof} Let $q$ be nonsingular. We aim to prove that then any higher rank quadric $\uQ_\nu(q)$ is smooth.

Since smoothness is a local property for the flat topology, we can in view of \ref{qf}\eqref{qfnsp} assume that $(M,q)$ is split of constant rank $r>0$ and by a possible further refinement that $\nu$ is constant, say $\nu = n \in \NN_+$. Note $2n \le r$ by Proposition~\ref{quadrepII}. Since we already know that $\uQ_\nu(q)$ is finitely presented \eqref{quadint}, smoothness is equivalent to the lifting property 
(\cite[IV$_4$, (17.5.1)]{EGA}, \cite[XI, 1.5]{SGA3}): 
for any local $A\in \Ralg$ and any nilpotent ideal $\fra \ideal A$ the canonical map $\uQ_n(q)(A) \to \uQ_n(q)(A/\fra)$ is surjective. We can even assume that $\fra^2 = 0$. It is no harm to replace $A$ by $R$ and put $\bar R = R/\fra$.

We fix $N \in \uQ_n(q)(\bar R)$ and need to show that $N$ lifts to some $\wtl N \in \uQ_n(q)(R)$. By \ref{quadrepII} again, there exists a totally isotropic submodule $N'$ of $(M,q)_{\bar R}$ such that $N \cap N' = 0$ and $N' \simlgr N^*$ via $b_{\bar q}$, thus $N\oplus N' \cong \HH(N)$, and $(M,q)_{\bar R} = (N \oplus N')\perp (N \oplus N')^\perp$.

We first deal with the Lagrangian case $2n=r$.  Then $(M,q)_{\bar R} = N \oplus N'$ and $N$ is a Lagrangian submodule of $(M,q)_{\bar R}$. Since $q$ is split hyperbolic, there exists a Lagrangian $N_0$ of $(M,q)$. It descends to a Lagrangian $\bar N_0$ of $(M,q)_{\bar R}$. Since $R$ and hence $\bar R$ is local, both $\bar N_0$ and $N$ are free of rank $n$, and are therefore isomorphic as $\bar R$--modules. By \ref{Lag}\eqref{Lag-b}, any $\bar R$--isomorphism $\bar N_0 \simlgr N$ can be extended to an isometry $g$ of $(M,q)_{\bar R}$.
\lv{
Then $(M,q)_{\bar R} = N \oplus N'$ and $N$ is a Lagrangian submodule of $(M,q)_{\bar R}$. Since $q$ is split, there exists a hyperbolic pair $(N_0, N_0')$ of submodules of $M$ analogous to $(N,N')$, i.e., $q(N_0) =0 =  q(N_0')$, $N_0' \cong N_0^*$ via $b_q$ and $M = N_0 \oplus N'_0$. The pair $(N_0, N_0')$ descends to a hyperbolic pair $(N_0\ot_R \bar R, N'_0 \ot_R \bar R)$ of $(M,q)_{\bar R}$. As $R$ is local, both $N_0 \ot_R \bar R$ and $N$ are free of rank $n$, and are therefore isomorphic as $\bar R$--modules. By \eqref{qfba-hyp3}, any isomorphism $N_0 \ot_R \bar R \simlgr N$ can be extended to an isometry $g$ of $(M_q)_{\bar R}$ mapping $N'_0 \ot_R \bar R$ to $N'$.}
Since $M$ has constant even rank, $\uO(q)$ is a smooth $R$--group scheme, see \ref{orthsc}\eqref{orthsc-d}. Therefore, $\uO(q) \to \uO(q_{R/\fra})$ is surjective and $g$ can be lifted to an isometry $g_0$ of $(M,q)$.  



We know that $\wtl N = g_0(N_0)\in \uQ_n(q)(R)$, and by construction $\wtl N$ is a lift of $N$. This settles the case $2n=r$.

Let now $n$ be arbitrary. Since $P=N \oplus N'$ is regular and free of rank $2n$, the decomposition $M\ot_R \bar R = P \perp P^\perp$ lifts to an orthogonal decomposition $M = \wtl P \perp \wtl P ^\perp$ with $\wtl P$ regular and free of rank $2n$, \cite[I, Cor.~3.4]{Ba}. Putting $\wtl q = q|_{\wtl P}$ we have a closed embedding $\uQ_n(\wtl q) \to \uQ_n(q)$ of $R$--schemes with $N$ belonging to $\uQ_n(\wtl q)(\bar R)$. By the previous paragraph we know that $\uQ_n(\wtl q)$ is smooth. Therefore $\uQ(\wtl q)(R) \to \uQ(\wtl q)(\bar R)$ is surjective, in particular  $N$ belongs to the image of $\uQ_n(q)(R) \to \uQ_n(\bar R)$.

Summarizing, we have shown that if $q$ is nonsingular, then every higher rank quadric $\uQ_\nu(q)$ is smooth, in particular this holds for $\uQ(q)$. We have already seen in \ref{quadint} that all higher rank quadrics $\uQ_\nu(q)$ are projective $R$--schemes.

The direction from right to left in \eqref{quadsm1} follows from \cite[XII]{SGA7},  where smoothness of $\uQ(q)$ is taken as definition of a so-called ordinary quadratic form. It is then shown in Prop.~1.2 of loc.\ cit. that after an \'etale extension ordinary quadratic forms become exactly the split quadratic forms of \ref{qf}\eqref{quadfoc}, i.e., they are nonsingular in our sense.  Alternatively, one can argue as follows. Let $F\in \Ralg$ be a field. Since smoothness is stable under base change we know that $\uQ(q)_F  = \uQ(q_F)$ is smooth. But then $q_F$ is nonsingular by \cite[22.1]{EKM}. By definition of nonsingularity, this implies that $q$ is nonsingular.
\end{proof}
\ms

Our next goal is to investigate the actions of $\uO(q)$ and $\uSO(q)$ on 
higher rank quadrics. There are two cases that are already present over fields, quadrics not of full rank and Lagrangians. In both cases it will be important that stabilizers are parabolic in $\uSO(q)$, which we establish in \ref{prop_hig}. In \ref{prop_hig}--\ref{traco} we will assume $\rank_R M \ge 3$. The reason for this restriction is that if $M$ has constant rank $2$, then 
$\uSO(q)$ is a torus by \ref{sogsc-even-ii} of \ref{sogsc}\eqref{sogsc-even}, rendering the theory of parabolic subgroups of $\uSO(q)$ uninteresting.

\comments{(2025-09-05) Above (blue), I replaced the previous statement {\tt $\uSO(q)$ is a group of multiplicative type} by ``torus''.}
Moreover, if $(M,q)$ has rank $2$ and contains an isotropic complemented submodule $N$ of rank $1$, i.e., a Lagrangian, then $(M,q) \cong \HH(N)$ by Corollary~\ref{carhyp} and we have already seen in \ref{Lag}\eqref{Lag-c} that $\uO(q)$ acts transitively on Lagrangians in the unimodular case, see \ref{prohla} for the $\uSO(q)$ action. 
\sm 

Proposition~\ref{prop_hig} extends \cite[Prop.~A.5(a)(i),(ii)]{GN-Sp}, dealing with the case of quadrics. The proof given here is different from the one of \cite{GN-Sp}. A different proof of \ref{prop_hig}\ref{prop_higa} is given in \ref{gmu}\eqref{gmu-b}. 

\begin{prop}[Stabilizers]\label{prop_hig} Let $(M,q)$ be a quadratic space satisfying $\rank_R M \ge 3$, let $\nu$ be a locally constant function with $1\le \nu(\gp)$ for every $\gp \in \Spec(R)$,  and let $N \in \uQ_\nu(q)(R)$. Then the following hold. \sm

\begin{enumerate}[label={\rm (\alph*)}]
  \item \label{prop_higa} The stabilizer $P=\Stab_{\uSO(q)}(N)$ of $N$ is a parabolic subgroup of $\uSO(q)$.\sm

\item \label{prop_higb} Assume\/ $2\nu < \rank_R M$. Then the orbit map induces an isomorphism $\uSO(q)/P \simlgr \uQ_\nu(q)$ of $R$--schemes.\sm

\item \label{prop_higc} Suppose\/ $2\nu = \rank_R M$. Then $P = \Norm_{\uO(q)}(P)= \Stab_{\uO(q)}(N)$ and the orbit map induces an isomorphism $\uO(q)/P \simlgr \uL(q)$ of $R$--schemes.
\end{enumerate}
\end{prop}

\pcomments{(2021-12-20, PG) Later ask someone, Fasel for example, whether there is no easier argument for \ref{prop_hig}\ref{prop_higa}.\sm

(EN, 2022-01-06) The proof using the dynamic method is still here in 'lv'}

\begin{proof} \ref{prop_higa} We have proven this in \ref{gmu}\eqref{gmu-b}, using the dynamic method \ref{pare}\eqref{pare-b}. We offer a different proof here. 
\lv{
Indeed, by Proposition~\ref{quadrepII}\eqref{quadrepII-c}
there exists a submodule $F\subset M$ such that $q_{E\oplus F} \cong \HH(E)$. Let $C=(E\oplus F)^\perp$, thus \begin{equation} \label{prop_qquadric1-i1}
M = E \oplus F \oplus C, \qquad E^\perp = E \oplus C.
\end{equation}
The stabilizer $P$ of $E$ is also the stabilizer of the flag $(E, E^\perp)$. It coincides with the parabolic subgroup $\rmP_G(\la)$ of \eqref{pare-1} where $\la(t)$ for $t\in R'{}\ti$ is defined by $\la(t)|_{E\ot R'} = t \Id$, $\la(t)|_{F \ot R'} = t\me \Id$ and $\la(t)|_C = \Id$. The morphism $\la$ is indeed a one-parameter subgroup of $G=\uSO(q)$ since $\la(t)$ is obviously orthogonal and since $G$ is the identity component of $\uO(q)$ in the even rank case and is the kernel of the determinant morphism in the odd rank case, cf.\ \ref{sogsc}.
} 

Recall \ref{pare}\eqref{pare-aa} 
that we need to show that $P$ is a smooth subgroup scheme for which all geometric fibres are parabolic subgroups of the appropriate fibres of $\uSO(q)$. Since these conditions are local for the flat topology, we can apply 
\ref{qf}\eqref{qfnsp} 
and assume that $q$ is a split quadratic form and that $\nu$ is constant, say $\nu = n\in \NN_+$. Since split forms arise by base change from $\ZZ$ to $R$, we can further suppose that $R=\ZZ$.

Let $P_{univ}$ over $\uQ_n(q)$ be the stabilizer of the universal point of the $\ZZ$--scheme $\uQ_n(q)$, i.e., in the setting of \ref{ag}\eqref{ag-i} 
we have $S=\Spec(\ZZ)$, $G= \uSO(q)$, $X=\uQ_n(q)$ and $H= P_{univ}$. Then $\uQ_n(q)$ is a reduced $\ZZ$--scheme, cf.\ the remark at the end of \ref{ag}\eqref{ag-red}. 
Moreover, by \ref{ag}\eqref{ag-ifp}, 
the $\uQ_n(q)$--scheme $P_{univ}$ is locally of finite presentation since both $\uSO(q)\to S$ and $\uQ_n(q)\to S$ are so. The geometric fibres of $P_{univ}$ are stabilizers of the corresponding fibres of $N$ and are therefore parabolic subgroups, hence smooth connected. They have the same dimension by \cite[\S3]{Co3}, see also \ref{Colem} below. Thus, the conditions \ref{ag-fI} and \ref{ag-fII} of \ref{ag}\eqref{ag-f} are fulfilled, allowing us to conclude that $P_{univ}$ is smooth over $\uQ_n(q)$ and then that $P$ is smooth over $R$. It follows that $P$ is a parabolic subgroup of $\uSO(q)$. \sm

\ref{prop_higb} According to \cite[XXI, 5.8.5]{SGA3}, the fppf quotient
$\uSO(q)/P$ is representable by a smooth projective $R$-scheme.
The orbit map gives rise to a monomorphism $f:\uSO(q)/P \to \uQ_\nu(q)$
between proper $R$-schemes (recall that projective schemes are in particular proper 
\cite[Tag 01WC]{St}. 
It follows that $f$ is a proper monomorphism, 
hence a closed immersion \cite[Tag 04XV]{St}. 
For each point $s$ of $\Spec(R)$, we know that $f_s:(\uSO(q)/P)_s \to (\uQ_\nu(q))_s$ is surjective by \cite[Cor. 3.4]{Co3}. Since $\uQ_\nu(q)_s$ is geometrically integral
(as quotient of the geometrically integral $\kappa(s)$--scheme $\uSO(q)_s$), it follows that $f_s$ is an isomorphism. The fibrewise isomorphism criterion \ref{ag}\eqref{ag-d} enables us to conclude that $f$ is an isomorphism. \sm

\ref{prop_higc} Since $M$ has even rank, $\uSO(q)$ is the identity component of $\uO(q)$, part~\ref{def-sosII} of \ref{sogsc}\eqref{sogsc-even}. Hence, according to \cite[Prop.~3.4.3]{G2}, the normalizer $ \Norm_{\uO(q)}(P)=: X$ is representable  by  a smooth $R$--scheme and $\uO(q)/X$ is representable. Moreover, from \cite[Prop.~4.5.5(3)]{G2} we know that $X/P$ is representable and that the composition $X \to \uO(q)  \xrightarrow{\Di} \ZZ/2\ZZ$ induces a closed and open immersion  $X/P \to \ZZ/2 \ZZ$.  Over an $R$--field $k$, we have $P(k)= X(k)$ since otherwise $\uSO(q)(k)$ would act transitively on $\uL(k)$ contradicting \cite[Th.~3.9]{Co3}. It follows that $P_k= X_k$ and then that $P=X$, in particular $\uO(q)/P$ is representable.

Furthermore, $\Stab_{\uO(q)}(N)$ is an $R$--subfunctor of $X$ and is therefore  representable by $P$. It then follows that the orbit map induces an $\uO(q)$--equivariant monomorphism $\uO(q)/P \to \uQ_\nu(q)$ between $R$--schemes. It is an epimorphism by \ref{hrq}, and therefore an isomorphism.
\end{proof}

\subsection{Dynkin diagrams}\label{Dyn} To describe the type of the parabolic subgroups of Proposition~\ref{prop_hig}\ref{prop_higa}, we use the Dynkin diagrams of type $\rmB_n$, $n\in \NN_+$, and $\rmD_n$, $n\ge 2$, with the standard enumeration, except in type $\rmD_3=\rmA_3$. They are depicted below.
\begin{align*} 
\rmB_1 = \rmA_1:\quad
&\begin{picture}(100,20)
\put(30,0){\circle*{3}}
\put(25,5){$\alpha_1$}
\end{picture}
\\ \rmB_n, n \ge 2: \; \; 
&\begin{picture}(200,20)  
\put(30,00){\line(1,0){20}}
\put(60,00){\dots}
\put(80,00){\line(1,0){40}}
\put(130,00){\dots}
\put(150,00){\line(1,0){20}}
\put(170,1.1){\line(1,0){20}}
\put(170,-1.2){\line(1,0){20}}
\put(176,-2.5){$>$}
\put(30,0){\circle*{3}} 
\put(50,0){\circle*{3}}
\put(80,0){\circle*{3}}
\put(100,0){\circle*{3}}
\put(120,0){\circle*{3}}
\put(150,0){\circle*{3}}
\put(170,0){\circle*{3}}
\put(190,0){\circle*{3}}
\put(25,5){$\alpha_1$}
\put(45,5){$\alpha_2$}
\put(185,5){$\alpha_n$}
\end{picture}
\\
\rmD_2 = \rmA_1 \times \rmA_1 : \; \; 
&\begin{picture}(100,20)
\put(30,0){\circle*{3}} 
\put(50,0){\circle*{3}}
\put(25,5){$\alpha_1$}
\put(45,5){$\alpha_2$}
\end{picture}
\\
\rmD_3=\rmA_3: \;\;
&\begin{picture}(100,20)  
\put(50,0){\circle*{3}} 
\put(70,10){\circle*{3}} 
\put(70,-10){\circle*{3}} 
\put(50,00){\line(2,1){20}} 
\put(50,00){\line(2,-1){20}} 
\put(75,8){$\alpha_2$}
\put(75,-13){$\alpha_3$}
\put(35,0){$\alpha_1$}
\end{picture}%
\\
\rmD_n, n\ge 4:
&\begin{picture}(200,20)  
\put(30,00){\line(1,0){20}}
\put(60,00){\dots}
\put(80,00){\line(1,0){40}}
\put(130,00){\dots}
\put(150,00){\line(1,0){20}}
\put(170,00){\line(2,1){20}} 
\put(170,00){\line(2,-1){20}} 
\put(30,0){\circle*{3}} 
\put(50,0){\circle*{3}}
\put(80,0){\circle*{3}}
\put(100,0){\circle*{3}}
\put(120,0){\circle*{3}}
\put(150,0){\circle*{3}}
\put(170,0){\circle*{3}}
\put(190,10){\circle*{3}}
\put(190,-10){\circle*{3}}
\put(25,5){$\alpha_1$} 
\put(45,5){$\alpha_2$}
\put(195,8){$\alpha_{n-1}$}
\put(195,-13){$\alpha_n$}
\end{picture}%
\end{align*}
\ms

\noindent
The following Proposition~\ref{Colem} can be extracted from \cite[3.8--3.13]{Co3}:

\begin{prop}\label{Colem}
  Let $k$ be a field, let $(M,q)$ be a split quadratic $k$--space of dimension $r\ge 3$, and put $n= [\frac{r}{2}]$, thus $\uSO(q)$ is a split $k$--group scheme of type $\De=\rmB_n$, $n\ge 1$ or $\De=\rmD_n$, $n\ge 2$, depending on $r$ being odd or even. Furthermore, let $d\in \NN$ satisfying $1\le d\le n$, let $N\in \uQ_d(q)(k)$ and let $P=\Stab_{\uSO(q)}(N)$.

  Then the type $t(P)$ of the parabolic subgroup $P$ of $\uSO(q)$ is the subset $\De\setminus I(d)$ where $I(d)\subset \De$ is given in the following table:
\begin{equation*}
 \begin{tabular}{|c||c|c|c|c|}
  \hline
  $\De$ & $\rmB_n$ &  $\rmD_n$ & $\rmD_n$ &  $\rmD_n$ \\
  \hline $d$ & $1\le d \le n$   & $1\le d \le n-2$ & $d=n-1$ & $d=n$ \\
  \hline $I(d)$ & $\{\al_d\}$ & $\{\al_d\}$
        & $\{\al_{n-1}, \al_n\}$ & $\{\al_{n-1}\}$
          or $\{\al_n\}$ \\
  \hline {\rm Case} &  {\rm (I)} & {\rm (II)} & {\rm (III)} & {\rm (IV)} \\
  \hline
  \end{tabular}
\end{equation*}
Conversely, any parabolic subgroup of $\uSO(q)$ of type $\De\setminus I(d)$ and $I(d)$ as above is the stabilizer of some $N' \in \uQ_d(q)(k)$.
\end{prop}

Some special cases are worth pointing out: In case $\De=\rmB_1$, so $d=1$, the type of $P$ is empty, meaning that $P$ is a Borel subgroup, see \cite[page~13]{Co3}.
For $\De=\rmD_n$ and $d=n$, the scheme $\uQ_n(q) = \uL(q)$ has two connected components and the type of the stabilizer depends on the connected component containing $N$.
The example $\De=\rmD_2$ is discussed in \cite[Ex.~3.8]{Co3}: for $d=1$ the parabolic subgroup $P$ is a Borel subgroup, while for $d=2$ we get two non-isomorphic maximal parabolic subgroups. \ms

We return to the setting of quadratic spaces over an arbitrary $R$ and consider the relation between the higher rank quadrics and the scheme $\Par\big(\uSO(q)\big)$ of parabolic subgroups of $\uSO(q)$, cf.\ \ref{pare}\eqref{pare-d}.

\begin{cor}\label{prohigqu} Let $(M,q)$ be a quadratic $R$--space of constant rank $r\ge 3$, and let $1\le d \le n= [\frac{r}{2}]$. We abbreviate $G=\uSO(q)$. \sm

\begin{inparaenum}[\rm (a)]
  \item \label{prohigqu-a} If $r$ is odd or if $r$ is even but $2d< r$,
  the higher rank quadric $\uQ_d(q)$ is $\uO(q)$--isomorphic to the scheme $\Par(G)_t$ of parabolic subgroups of $G$ of constant type $t=\De \setminus I(d)$, where $I(d)$ is as in the cases {\rm (I)--(III)} of {\rm \ref{Colem}}.
   \sm

\item\label{prohigqu-c} The $\uO(q)$--equivariant map
$\uQ_d(q) \to \Par(G)$, defined by assigning to each $R$--algebra $A$ and each $N\in \uQ_d(q)(A)$ the stabilizer $\Stab_{G_{A}}(N)$, is an open and closed immersion.
\end{inparaenum}
\end{cor}

\begin{proof}
  \eqref{prohigqu-a} The statement is local for the fppf topology, so that we can assume $\uQ_d(q)(R) \ne \emptyset$, say $N\in \uQ_d(q)(R)$. The type of the parabolic subgroup $P = \Stab_G(N)$ can then be calculated by passing to an algebraic closure of $\ka(\gp)$, $\gp\in \Spec(R)$. The claim therefore follows from Propositions~\ref{prop_hig} and \ref{Colem} together with the isomorphism \eqref{pare-d2}. \eqref{prohigqu-c} is a consequence of \eqref{prohigqu-a}.
 \end{proof}
\sm

We will deal with the case $r$ even and $n=d=\frac{r}{2}$ later in Corollary~\ref{prohla}.
\sm

The case $d=1$ of Corollaries~\ref{prohigqu} and \ref{traco} is proven in parts (a)(iii) and (b) of \cite[Prop.~A.5]{GN-Sp}. Corollary~\ref{traco} deals with transitivity of the action of $\uSO(q)$ on higher rank quadratics that are not Lagrangians, complementing the Lagrangian case considered in \ref{Lag}\eqref{Lag-b} and \ref{hrq}. With the exception of part \ref{traco}\eqref{traco-bii}, Corollary~\ref{traco} is essentially the same as Lemma~\ref{traq}. Our reason for stating it here again is to emphasize the broader perspective: it is the special case of Demazure's Theorem~\ref{thm_conj_demazure} with $G=\uSO(q)$, once one knows Corollary~\ref{prohigqu}\eqref{prohigqu-a} saying that
$X=\Par(G)_t \cong \uQ_d(q)$.

\begin{cor}\label{traco}
  Let $(M,q)$ be a quadratic $R$--space of constant rank $r\ge 3$, let $d\in \NN$ satisfy $1\le d < [\frac{r}{2}]$, and abbreviate $\uQ_d=\uQ_d(q)$. Then the following hold. \sm

\begin{inparaenum}[\rm (a)]
\item \label{traco-a} The action of $\uSO(q)$ on the higher rank quadric $\uQ_d$ is transitive on the  small affine Zariski site of $\Spec(R)$.
  \sm

\item \label{traco-b} \label{traco-bi} Let $R$ be an LG ring. Then 
$\uSO(q)(R)$ acts transitively on $\uQ_d(R)$.
\sm

\item \label{traco-bii} If $R$ is semilocal and $R'$ is a
finite $R$-algebra such that $\uQ_d(R') \not = \emptyset$, the map
\[ \uQ_d(R') \longto \textstyle \prod_{\gm \in \Spec(R) \text { maximal}}\,  \uQ_d(R'/ \gm R') \]
is onto.
\end{inparaenum}
\end{cor}

In view of the restriction $d< [\frac{r}{2}]$ in \ref{traco}, it is natural to ask about the $\uSO(q)$--orbits in the excluded case $2d = \rank M$. They are described in \ref{prohla} as a corollary of the Stein factorization of the projective morphism $\uL(q) \to \Spec(R)$, presented in Proposition~\ref{deligne}. Its proof requires some preparation, notably complemented left ideals in separable algebras, \ref{sep-li}, and elementary idempotents in composition algebras, \ref{eica}, in particular in quadratic \'etale algebras, \ref{elid} and \ref{scelid}. 

\comments{(2025-09-04) Lemma~\ref{sep-li} is taken from the old file `parabolic-SB-LG', proposed to be eliminated; it is used in \ref{elid} for quadratic \'etale and again later for quaternions. }

\begin{lem}[Left ideals in separable algebras]\label{sep-li}
Let $B$ be a separable $R$--algebra, and let $L \subset B$ be a left ideal such that $B/L$ is projective as $R$--module. Then there exists an idempotent $c\in B$ such that $L = Bc$ and $B=Bc \oplus B(1_B-c)$. Hence, $L$ and $B/L \cong B (1_B -c)$ are projective as $B$--modules.
\end{lem}

\begin{proof} The exact sequence $0 \to L \to B \to B/L \to 0$ of $B$--modules is split-exact as exact sequence of $R$--modules. Hence, by \cite[4.4.1]{Ford}, it is split-exact as sequence of $B$--modules. Thus, there exists a left ideal $L' \subset B$ such that $L \oplus L' = B$. Decomposing $1_B = c + c'$ with $c\in L$, $c'\in L'$, one finds that $c$ is an idempotent with $L = Bc$ and $L' = B(1_B - c)$ \cite[1.1.20]{Ford}.
\end{proof}

\comments{(2025-09) The new stuff in blue is used to identify in \ref{scelid} that the image of the Stein factorization in \eqref{Dem2} is $\Spec(\Dis(q))$ -- something that was missing in the old version. Of course, Deligne did not say anything about this -- it must have been obvious for him.}

\subsection{Elementary idempotents in composition algebras}\label{eica}
Let $C$ be a composition algebra over $R$, see for example \cite[19.5]{PRbook}. Examples are quadratic \'etale algebras and quaternion algebras that we will consider in \ref{eica} and \ref{spq}. We denote by $1_C$ the identity element of $C$, and by $\Tr_{C/R}$ and $\rmN_{C/R}$ its  trace and norm respectively, cf.~\ref{trno}. The following are equivalent for an element $e\in C$: 
\begin{enumerate}[label={\rm (\roman*)}]
   
 \item\label{eicai} $\Tr_{C/R}(e) = 1$ and $\rmN_{C/R}(e) = 0$, 
 
  \item \label{eicaii} $e$ is an idempotent satisfying $ 0 \ne (e\ot_R 1_{R'}) \ne 1_{R'}$ for all $R'\in \Ralg$, 
  
 \item \label{eicaiii}  $e$ and $e'=1_C - e$ are unimodular idempotents, 
 
 \item \label{eicaiv} $e$ and $e' = 1_C - e$ are idempotents whose Peirce spaces $eCe$ and $e'Ce'$ are free of rank $1$ as $R$--modules with bases $e$ and $e'$,

 \item \label{eicav} the map $\vphi_e \co R \times R \to C$, $(r_1, r_2) \mapsto r_1 e + r_2e'$ is an injective $R$--algebra homomorphism. 
     
\end{enumerate}
Indeed, the equivalences \ref{eicai}--\ref{eicaiii} are \cite[Exc.~16.23]{PRbook}, the implication \ref{eicaiii} $\implies$ \ref{eicaiv} is \cite[Exc.~19.35]{PRbook}, the
implication \ref{eicaiv} $\implies$ \ref{eicav} follows from $e e' = 0$, and the 
implication \ref{eicav} $\implies$ \ref{eicaii} is clear. 

An element $e\in C$ satisfying the equivalent conditions above is called an {\em elementary idempotent}. If $e$ is an elementary idempotent, they so is $e'$. 

\begin{lem}[Elementary idempotents in quadratic \'etale algebras]\label{elid} Let $D$ be a quadratic \'etale $R$--algebra, see for example\/ \cite[19.19]{PRbook}. We use the notation of {\rm \ref{eica}}, and observe that $D=De \oplus De'$ for any idempotent $e$ with $e'=1_D - e$. Then the map  
\begin{equation}\label{eled3} \begin{split} 
\{ e \in D : \text{ $e$ elementary idempotent}\} \longto &\Hom_{\Ralg}(D,R), \\
     e \mapsto &(r_1e + r_2 e' \mapsto r_1) 
\end{split}\end{equation} 
is a bijection, and $D \cong Re \times Re'$ for any elementary idempotent $e$.   
\end{lem}

\begin{proof} The map is well-defined and injective by \ref{eica}. To prove surjectivity, let $\al \in \Hom_{\Ralg}(D, R)$. Since $\al$ is surjective, $\Ker(\al)$ satisfies the assumptions of Lemma~\ref{sep-li}. Thus, there exists an idempotent $e\in D$ such that $D=De \oplus De'$ and $\Ker(\al) = De'$ where $e' = 1_D - e$. 
The idempotent $e$ is elementary, for example by \ref{eica}\ref{eicaii}. The last claim follows \ref{eica}\ref{eicav}. \end{proof}

\subsection{The scheme of elementary idempotents of a quadratic \'etale algebra}\label{scelid} Let $D$ be a quadratic \'etale $R$--algebra. We define an $R$--functor $\underline{\rm Splid}(D)$ whose $R'$--points, $R'\in \Ralg$, are
\begin{equation}\label{scelid1} \begin{split}
  \big(\underline{\rm Splid}(D)\big)(R') &=
     \{ e \in D\ot_R R' : \text{ $e$ elementary idempotent}\}
\\  &= \Hom_{\Rpalg}(D\ot_R R', R')  = \Hom_{\Ralg}(D, R') .
\end{split}\end{equation}  
The second equality above is the identification \ref{elid}, and the third equality is standard. In particular, the third equality shows that  
$\underline{\rm Splid}(D)$ is the functor of points associated with the $R$--scheme $\Spec(D)$, which is a quadratic \'etale cover of $\Spec(R)$. 

Lemma~\ref{elid} says that elementary idempotents of $D$ give rise to a splitting of $D$, and, conversely, every splitting of $D$ is given by an elementary idempotent. Quadratic \'etale $R$--algebras are the same as composition algebras of rank $2$ \cite[19.19]{PRbook}. Splitting data for arbitrary composition algebras $C$ of constant rank are introduced in \cite[25.2]{PRbook}; for rank $2$, they are precisely the elementary idempotents, as in\eqref{scelid1}. It is proven in \cite[26.3, 26.5]{PRbook} that splitting data of arbitrary constant rank composition algebras are represented by a smooth affine scheme $\mathbf{Splid}(C)$. Indeed, for $C=D$ quadratic \'etale, $\mathbf{Splid}(C)= \Spec(D)$ is even quadratic \'etale by the above identification. 

\comments{(2025-08-30) I did not find a proof in \cite{PRbook} that $\mathbf{Splid}(C)$ is quadratic \'etale for $C=D$. Once you have read \ref{elid} and \ref{scelid}, I suggest asking Skip and Holger about this. Maybe, it is hidden somewhere in their book. 

(2026-03-27) Discussed with PG; no action needed}
\sm 

We now come back to Lagrangian submodules of quadratic spaces. We treat the case of hyperbolic spaces in \ref{dgan}--\ref{eor} and, finally, the general case in \ref{deligne}. 
The following Lemma~\ref{dgan}, dealing with open subfunctors of the Grassmannian functor $\ulGr_\nu(M)$, is well-known in the constant rank case.  

\begin{lem}\label{dgan}
Let $M$ be an $R$--module and let $N\subset M$ be a submodule. $M=N \oplus N'$.
\sm

\begin{inparaenum}[\rm (a)] \item\label{dgan-a} Suppose that $N$ is a complemented submodule, say $M = N \oplus N'$. For $\vphi \in \Hom_R (N',N)$ let $\Ga_\vphi = \{ \vphi(n') + n' : n'\in N'\}$ be the graph of $\vphi$.
Then $\Ga_\vphi$ is a submodule of $M$ complementing $N$, and the map
\begin{align*} \Ga \co \Hom_R(N', N) &\simlgr U_N := \{ P \subset M: M= N \oplus P \}, 
\\ \vphi &\; \mapsto \quad \Ga_\vphi\; ,
\end{align*}
is a bijection. \sm

\item\label{dgan-b} Suppose that $N$ is a finite projective submodule of $M$ with rank function $\nu \in \NN(R)$. Note that such a submodule always exists if\/ $\ulGr_\nu(M)(R) \allowbreak \ne \emptyset$. We define a subfunctor $\ulU_N$ of $\ulGr_\nu(M)$ by putting 
    \[ \ulU_N(T) = \{P\in \ulGr_\nu(M)(T) : \text{ $N_T \to M_T \to M_T/P$ is surjective}\} 
    \]
    for $T\in \Ralg$, where $N_T \to M_T$ is the base change of the inclusion $N \subset M$ and $M_T \to P$ is the quotient map.  Then $\ulU_N$ is an open subfunctor of $\ulGr_\nu(M)$. 
    
    If, in addition, $N$ is complemented by a finite projective submodule $N'$,  then the bijection $\Ga$ of {\rm \eqref{dgan-a}} induces an isomorphism of $R$--functors, 
    \begin{equation}\label{dgan-bHom}
       \ulW\big(\Hom(N', N)\big) \simlgr \ulU_N. 
    \end{equation} 
  Thus $\ulU_N$ is represented by the smooth affine $R$--scheme $\uW\big( \Hom_R(N', N)\big)$.
\end{inparaenum}
\end{lem}
\lv{
\begin{proof}
It is obvious that $\Ga_\vphi$ is a submodule of $M$. To see well-definedness of the map $\Ga$, i.e., $M = \Ga_\vphi \oplus N$, observe that $\Ga_\vphi \cap N = 0$ because $n' \oplus \vphi(n') \in N$ implies $n' \in N \cap N'= 0$ since $\vphi(n') \in N$,  so that $n' \oplus \vphi(n') = 0$. Also, we can write any $m\in M$, say $m = n + n'$ with $n\in N$ and $n' \in N'$, in the form $n + n' =  (n' + \vphi(n')) + (n - \vphi(n'))\in \Ga_\vphi \oplus N$, proving that $M = \Ga_\vphi \oplus N$ and hence that $\Ga$ is well-defined.

The map $\Ga$ is surjective: Given a submodule $N''$ with $N'' \oplus N = M$, any $n''\in N''$ can be uniquely written as $n'' = n + n'$ with $n\in N$ and $n'\in N'$. This defines a map $\vphi \co N' \to N$, $n'\mapsto n$, which is linear by uniqueness of $n$ and $n'$. From $n'' = n'+ \vphi(n')\in \Ga_\vphi$ we get $N'' \subset \Ga_\vphi$. Because $M = \Ga_\vphi \oplus N = N'' \oplus N$, this implies $N''= \Ga_\vphi$.

The map $\Ga$ is injective: Given $\vphi$, $\psi \in \Hom_R(N',N)$ such that $\Ga_\vphi = \Ga_\psi$, for any $n'\in N'$ there exists $n_1' \in N'$ such that $n' + \vphi(n') = n'_1 + \psi(n'_1)$, so that $n'=n'_1\in N'$ and $\vphi(n') = \psi(n'_1)$ and then $\vphi = \psi$ follows.
\end{proof}
}
\begin{proof}
We leave the easy proof of \eqref{dgan-a} to the reader. The proof of \eqref{dgan-b} is a straightforward adjustment of the proof of \cite[I, (9.7.4.6)]{EGA-neu} (see also \cite[I, \S1, 3.9]{DG}), given there for $\Spec(R)$ replaced by an arbitrary scheme but assuming $N$ of constant rank. We present the proof of \eqref{dgan-b} for the sake of completeness.
\sm


To prove that $\ulU_N$ is an open subfunctor, we fix $T\in \Ralg$ and $P \in \ulGr_\nu(M)(T)$. We then have to show that there exists an ideal $\frt \ideal T$ such that for every $T'\in \Ralg$ and for every $R$--algebra homomorphism $\psi \co T \to T'$ we have the equivalence
\begin{equation} \label{dgan-b1}
T' \psi(\frt) = T' \iff P_{T'} \in \ulU_{N}(T').
\end{equation}
Let $\al \co N_T \to M_T/P$ be the composition $N_T \to M_T \to M_T/P$. By definition, $\al$ is surjective. But since $N_T$ and $P$ are both finite projective of the same rank $\nu_T$, the map $\al$ is in fact an isomorphism. Also, $P$ is a complemented submodule of $M_T$. Hence, we can and will view $P_{T'}$ as a submodule of $M_{T'}$. It follows that $N_{T'} \to M_{T'} \to M_{T'} / P_{T'}$ is the base change of the isomorphism $\al$. Thus the right hand side of \eqref{dgan-b1} $\iff \al_{T'} = \al \ot 1_{T'}$ is surjective $\iff \Coker(\al_{T'}) = \Coker (\al) \ot_{T} T' = 0$. Since $M_T/P$ and therefore also $\Coker(\al)$ has finite type, it follows from  
\cite[II, \S4.4, Prop.~17 and 19]{BAC} that \eqref{dgan-b1} holds for $\frt$ the annihilator of $\Coker(\al)$. 

For any $T\in \Ralg$, we have $M_T = N_T \oplus N'_T$. Also, $P\in \ulU_N(T) \iff \al_T$ is an isomorphism $\iff M_T = P \oplus N_T$. By \eqref{dgan-a} this is in turn equivalent to $P = \Ga_\vphi$ for some $\vphi \in \Hom_T(N'_T, N_T)$. But $\Hom_R(N', N) \ot_R T \simlgr \Hom_T(N'_T, N_T)$ since $N'$ is finite projective.
\end{proof}
\sm

\textbf{Remark.} If $M$ is free of finite rank,  the family $(\ulU_N)$ indexed by all complemented projective submodules $N\subset M$ of constant rank $n\in \NN$ covers $\ulGr_n(M)$, and thus yields a proof that 
$\ulGr_n(M)$ is represented by a smooth scheme.
\sm

Recall \ref{qf}\eqref{qf-hyp} and \ref{Lag}: A Lagrangian of a hyperbolic $R$--space $(M,q)$ is a direct summand $L$ of $M$ such that $2\rank_R L = \rank_R M$ and $q(L)=0$. Any submodule $L'$ satisfying $M = L \oplus L'$ and $q(L')= 0$ is then a Lagrangian too, we have $\HH(L) \cong (M,q) \cong \HH(L')$ as quadratic spaces, $L' \cong L^*$ as $R$--modules and $(L, L')$ is a hyperbolic pair of submodules. We call a hyperbolic pair $(N,N')$ of submodules of $(M,q)$ a {\em Lagrangian pair\/} if $2 \rank_R L = \rank_R M$, in which case $(N', N)$ is also a Lagrangian pair.

\begin{lem}[{\cite[XII, 1.6]{SGA7}}]\label{lav} Let $(L,L')$ be a Lagrangian pair of submodules of the hyperbolic $R$--space $(M,q)$.
\sm

\begin{inparaenum}[\rm (a)] \item  \label{lav-a}
  Let $b_q$ be the polar form of $q$ and let $\scL_2(L'; R)$ be the $R$--module of bilinear forms $L'\times L' \to R$. Then
  \[ \La : \Hom_R(L', L) \simlgr \scL_2(L';R), \quad
    \vphi \mapsto \big((\ell'_1, \ell'_2) \mapsto b_q(\ell_1', \vphi(\ell_2')) \big)
    \]
  is an isomorphism of $R$--modules.
  \sm

\item\label{lav-aa} Let $\vphi \in \Hom_R(L', L)$. We associate with $\vphi$ the endomorphism  $\wdh \vphi$ of $M$, defined by $\wdh\vphi(\ell + \ell') = (\ell + \vphi(\ell') ) + \ell'$ for $\ell \in L$ and $\ell' \in L'$. Then the following are equivalent:
\end{inparaenum}

\begin{enumerate}[label={\rm (\roman*)}]
  \item \label{lav-aa-i} $\La(\vphi)$ is alternating,

  \item \label{lav-aa-ii} $\Ga_\vphi$ is totally isotropic,

  \item \label{lav-aa-iii} $\wdh \vphi$ is orthogonal. 
\end{enumerate}
In this case, $(L, \Ga_\vphi)$ is a Lagrangian pair and $\wdh \vphi \in \SO(q)$.
We call $\vphi$ {\em alternating\/} if \ref{lav-aa-i}--\ref{lav-aa-iii} hold. \sm

\begin{inparaenum}[\rm (a)]\setcounter{enumi}{2}
  \item\label{lav-b} Let $\scA_2(L'; R)$ be the submodule of alternating bilinear forms on $L'$ and let $V_L$ be the set of hyperbolic pairs of $(M,q)$ of the form $(L,L'')$. Also recall the isomorphism $\Ga \co \Hom_R(L', L) \simlgr U_L$ of Lemma~{\rm \ref{dgan}}. Then
      \[ \Ga \circ \La\me \co \scA_2(L'; R) \; \simlgr \; V_L \]
      is a bijection.
\sm

\item\label{lav-c} We define a subfunctor $\ulV_L$ of the Lagrangian functor $\ulL(q)$ of\/ {\rm \ref{hrq}} by
\[ \ulV_L(T) = \{ P\subset  M_T : M_T = L_T \oplus P, \; q(P) = 0 \} \]
for $T\in \Ralg$. Then $\ulV_L$ is an open subfunctor of $\ulL(q)$, represented by an open subscheme $\uV_L$ of the Lagrangian quadratic $\uL(q)$ which is isomorphic to the smooth affine scheme $\uW(\scA_2(L; R)\big)$.
\end{inparaenum}
\end{lem}

\begin{proof} \eqref{lav-a} follows by replacing $L$ in $\Hom_R(L', L)$ by $L'{}^*$, using the isomorphism $L \simlgr L'{}^*$, $\ell \mapsto b_q(\ell, - )$, and then applying the canonical isomorphism
  \[ \Hom_R (L', L'{}^*) \simlgr \scL_2(L'; R), \quad
     \psi \mapsto \big( (\ell_1', \ell_2') \mapsto \psi(\ell_1')(\ell_2) \big).
  \]

\eqref{lav-aa} Let $\vphi\in \Hom_R(L', L)$, and let $\ell \in L$, $\ell'\in L'$. Then
\begin{equation} \label{lav-aa1}
 q\big( (\ell + \vphi(\ell'))  + \ell'\big) = b_q(\ell + \vphi(\ell'), \ell'),
\quad q(\ell + \ell') = b_q(\ell, \ell').
\end{equation}
Hence \ref{lav-aa-i} $\iff b_q(\vphi(\ell'), \ell') = 0$ for all $\ell'\in L' \iff$ \ref{lav-aa-ii} by the first equation in \eqref{lav-aa1} with $\ell = 0$, and \ref{lav-aa-i} $\iff$ \ref{lav-aa-iii} by comparing the two equations in \eqref{lav-aa1}.

We postpone the proof that $\wdh \vphi \in \SO(q)$, but note that clearly $\det(\wdh \vphi) =1$, proving $\wdh \vphi \in \SO(q)$ in case $2\in R\ti$ or $M$ has constant odd rank by \ref{sog}\eqref{sog-b}.
\sm

\eqref{lav-b} follows from \eqref{lav-aa} and the definitions.
\sm

\eqref{lav-c} The functor $\ulU_L$ of Lemma~\ref{dgan}\eqref{dgan-b} is an open subfunctor of $\ulGr_\nu(M)$ for $\nu = \rank_R L$
Hence, $\ulV_L$, being the intersection of $\ulU_L$ and $\ulL(q)$, is an open subfunctor of $\ulL(q)$.
The claimed representability of $\ulV_L$ follows from \eqref{lav-b}, analogous to the argument in \ref{dgan}\eqref{dgan-b} keeping in mind that $\scA_2(L';R)$ is a finite projective module since $L'$ is so.
\sm

We can now prove that $\wdh \vphi \in \SO(q)$ for an alternating $\vphi \in \Hom_R(L',L)$. We have a homomorphism $f\co \uW(\scA_2(L; R)\big) \to \uO(q)$ of $R$--group schemes which on $S$--points, $S\in \Ralg$, is given by $\La\me(\vphi) \mapsto \wdh\vphi$. Since $\uW(\scA_2(L; R)\big)$ is locally of finite presentation and has geometrically connected fibres, Lem\-ma~\ref{lem_neutral0LG}\eqref{lem_neutral1} says that $f$ factors through the identity component of $\uO(q)$, which is $\uSO(q)$ by \ref{somax}.
\end{proof}
\sm

In the next Lemma~\ref{clic} we consider the $(\ZZ/2\ZZ)$--graded Clifford algebra $\Cli(M,q)$ of a hyperbolic space $(M,q)$ with a Lagrangian pair $(L,L')$. The  endomorphism algebra $\End_R(\bwe L)$ of the exterior algebra $\bwe L$ of $L$ also has a $(\ZZ/2\ZZ)$--grading, induced by the decomposition $\bwe L = \bwe_0 L \oplus \bwe_1 L$ into even and odd elements. We recall, see e.g. \cite[IV, (2.1.1)]{K}, that both algebras are isomorphic as graded algebra under the isomorphism
\begin{equation} \label{lav-2}
 \Psi' \co \Cli(M,q) \simlgr \End_R(\bwe L),
\end{equation}
uniquely determined by
\begin{equation} \label{lav-3}
\begin{split}
\Psi'(\ell)(\ell_1\we \cdots \we \ell_r) &= \ell \we \ell_1 \we \cdots \we \ell_r, \\ \Psi'(\ell')( \ell_1\we \cdots \we \ell_r) &=
  \textstyle \sum_{i=1}^r \, b_q(\ell', \ell_i) \, \ell_1 \we \cdots \we \wdh {\ell_i} \we \cdots \we \ell_r
\end{split}\end{equation}
for $\ell$, $\ell_1, \ldots, \ell_r \in L$ and $\ell'\in L'$. The notation $\Psi'$, rather than $\Psi$, should indicate that $\Psi'$ depends on $L'$.

\begin{lem}\label{clic} Let $(L, L')$ and $(L, L'')$ be two hyperbolic pairs of the faithful hyperbolic space $(M,q)$, and let
\[   \Cli(M,q) \xrightarrow{\; \Psi'\;} \End_R(\bwe L) \xleftarrow{\; \Psi''\,}
\Cli(M,q)\]
be the isomorphisms \eqref{lav-2} defined with respect to the pairs $(L,L')$ and $(L,L'')$ respectively. By {\rm \ref{dgan}\eqref{dgan-a}} and {\rm \ref{lav}\eqref{lav-aa}} we know $L''= \Ga_\vphi$ for a unique alternating $\vphi \in \Hom_R(L', L)$. Then
\begin{equation}\label{clic1}
 \Psi' = \Psi'' \, \circ \, \Cli(\wdh \vphi),
\end{equation}
where $\Cli(\wdh \vphi)$ is the automorphism of $\Cli(M, q)$ associated with the orthogonal transformation $\wdh \vphi$ of $(M,q)$. In particular, for the discriminant algebra $\Dis(q)$ of $(M,q)$, {\rm \ref{qfdi}}, we have
\begin{equation} \label{clic2}
\Psi'|_{\Dis(q)} \, = \, \Psi''|_{\Dis(q)}.
\end{equation}
\end{lem}

\begin{proof}
For the proof of \eqref{clic1} it suffices to show that both sides coincide when applied to $\ell \in L$ and $\ell'\in L'$.
By \eqref{lav-3} this is clear for $\ell$ because $\Cli(\wdh\vphi)(\ell) = \wdh\vphi(\ell) = \ell$ and both $\Psi'(\ell)$ and $\Psi''(\ell)$ are just left multiplication by $\ell$. For $\ell'$ the right-hand side of \eqref{clic1} becomes
\begin{align*}
 &\big( \big(\Psi''\circ \Cli(\wdh \vphi)\big)(\ell')\big)( \ell_1\we \cdots \we \ell_r) =
     \Psi'' (\vphi(\ell') + \ell')( \ell_1\we \cdots \we \ell_r)    \\
 &\quad = \textstyle \sum_{i=1}^r \, b_q(\vphi(\ell') + \ell' , \ell_i) \, (\ell_1 \we \cdots \we \wdh {\ell_i} \we \cdots \we \ell_r)
  = \Psi'(\ell')( \ell_1\we \cdots \we \ell_r).
\end{align*}
The equation \eqref{clic2} follows from \eqref{clic1} and $\Cli(\wdh \vphi)|_{\Dis(q)} = \Id_{\Dis(q)}$ since $\wdh \vphi \in \SO(q)$ by \ref{lav}\eqref{lav-aa} and $\SO(q) = \{ g\in \orth(q): \Dis(g) = \Id_{\Dis(q)}\}$ by \eqref{sog1}.
\end{proof}

\subsection{The elementary idempotent of $\Dis(q)$ associated with a Lagrangian \cite[XII, 1.6, 1.7]{SGA7}}\label{eidLa} 
Let again $(M,q)$ be a faithful hyperbolic space with a Lagrangian $L$. We can embed $L$ in a Lagrangian pair $(L, L')$. The algebra isomorphism $\Psi'$ of \eqref{lav-2} maps $\Dis(q)$, the centre of $\Cli_0(M,q)$, onto the centre of the even part of $\End_R(\bwe L)$. Hence, there exists a unique idempotent $e=e(L) \in \Dis(q)$ such that $\Psi'(e)$ acts as the identity on $\bwe_0(L)$ and as zero on $\bwe_1(L)$. 
While $\Psi'$ depends on $(L,L')$, the formula~\eqref{clic2} shows that $e(L)$ only depends on $L$, thus justifying the notation.

\begin{lem}\label{eor} Let $(M,q)$ be a faithful hyperbolic $R$--space, let $L$ be a Lagrangian of $(M,q)$, and let $e(L)$ be the associated elementary idempotent of\/ {\rm \ref{eidLa}}.\sm

\begin{inparaenum}[\rm (a)] \item \label{eor-a}
  Let $g\in \orth(q)$. Then
\begin{equation}  \label{eor-a1}
   e\big( g(L) \big) = \Dis(g) \big( e(L) \big),
\end{equation}
where $\Dis(g)$ is the automorphism of $\Dis(q)$ induced by $g$, see \eqref{sog-dishom}. \sm

\item\label{eor-b} Let $R$ be a unimodular ring. Then the following are equivalent for a Lagrangian $L_1$ of $(M,q)$:
\end{inparaenum}
\begin{enumerate}[label={\rm (\roman*)}] 

\item\label{eor-bi} $e(L) = e(L_1)$,

\item\label{eor-bii} there exists $g\in \SO(q)$ such that $g(L) = L_1$.
 \end{enumerate}
\sm

\noindent If $R=k$ is a field, then \ref{eor-bi} and \ref{eor-bii} are equivalent to
\sm

\begin{enumerate}[label={\rm (\roman*)}] \setcounter{enumi}{2}
 \item \label{eor-biii} $\dim_k\, L/(L\cap L_1)$ is even.
\end{enumerate}
\end{lem}

\begin{proof} \eqref{eor-a} Let $L'$ be a Lagrangian such that $(L, L')$ is a hyperbolic pair of $(M,q)$. Then $\big(g(L), g(L')\big)$ is also a hyperbolic pair, and as in \cite[IV, (2.1.1)]{K} we get a commutative diagram
\[ \xymatrix@C=50pt{
   \Cli(M,q) \ar[d]_{\Cli(g)}^\cong  \ar[r]^{\Psi'}_\cong & \End_R(\bwe L)
   \ar[d] ^{\bwe(g|_L)}_\cong
   \\
   \Cli(M,q) \ar[r]^{\Psi_1'}_\cong & \End_R(\bwe g(L))
}\]
where $\Psi'$ and $\Psi_1'$ are the isomorphisms \eqref{lav-2} with respect to $(L, L')$ and $(g(L)_, g(L'))$ respectively, and where $\bwe(g|_L)$ is the extension of $g|_L$ to an isomorphism $\bwe(L) \simlgr \bwe\big(g(L)\big)$. The formula \eqref{eor-a1} is then immediate from the definitions. \sm

\eqref{eor-b} \ref{eor-bi} $\implies$ \ref{eor-bii}: By \ref{Lag}\eqref{Lag-c} there exists $g\in \orth(q)$ satisfying $g(L) = L_1$. Hence $e(L) = e(L_1) = e(g(L)) = \Dis(g)\big(e(L)\big)$ by \eqref{eor-a1}. Thus, the automorphism $\Dis(g)$ of $\Dis(q)$ fixes the elementary idempotent $e(L)\in \Dis(q)$. By \ref{elid}, $\Dis(g)$ fixes a basis of $\Dis(q)$, i.e., $\Dis(q) = \Id_{\Dis(q)}$ and therefore $g\in \SO(q)$. The implication \ref{eor-bii} $\implies$ \ref{eor-bi} is immediate from \eqref{eor-a1} (and holds for arbitrary $R$).

In the field case, the equivalence \ref{eor-bi} $\iff$ \ref{eor-biii} is shown in \cite[XII, Prop.~1.12]{SGA7} (and \ref{eor-bii} $\iff$ \ref{eor-biii} in \cite[Ex.~3.5]{Co3}).
\end{proof}

\subsection{The morphism $\ue\co \uL(q) \to \Spec(\Dis(q))$ \cite[XII, 1.7]{SGA7}} \label{Dem} Let again $(M,q)$ be a faithful hyperbolic space with a Lagrangian $L$. In \ref{eidLa} we have associated an elementary idempotent $e(L) \in \Dis(q)$ with $L$. 

Replacing $q$ by $q_{R'}$, $R'\in \Ralg$, we get a map
\[ \ulL(q)(R') = \ulL(q_{R'}) \xrightarrow{\;\; \ul{e}(R') \;\; } \Dis(q_{R'}) = \Dis(q) \ot_R R',\]
which is functorial, thus giving rise to a morphism $\ul{e}\co \ulL(q) \to \ulW\big( \Dis(q) \big)$ of $R$--functors and then to a morphism of the corresponding $R$--schemes
\begin{equation} \label{Dem1}
   \ue \co \uL(q) \to \Spec\big(\Dis(q)\big),
\end{equation}
where $\Spec\big(\Dis(q)\big)$ is the spectrum of the $R$--algebra $\Dis(q)$.
We then get a commutative triangle of scheme morphisms
\begin{equation}  \label{Dem2} \vcenter{
\xymatrix{\uL(q) \ar[rr]^\ue \ar[dr]_f && \Spec\big( \Dis(q)\big) \ar[dl]^p  \\ & \Spec(R)
}}\end{equation}
where $f$ and $p$ are the canonical structure maps.

In the following result we use the concept of the Stein factorization of a proper morphism of schemes,  which the reader can find in \cite[III, (4.3.3)]{EGA} in the noetherian case and in \cite[Tag 03H2]{St} in general. Any projective morphism is proper by \cite[II, (5.5.3)]{EGA} (or see \cite[Tag 01WC]{St}). Thus, the structure morphism $f \co \uL(q) \to \Spec(R)$ has a Stein factorization, since $f$ is projective by \ref{quadint}.


%

\begin{thm}[Deligne for general base {\cite[XII, Prop.~2.8]{SGA7}}]\label{deligne} Let $(M,q)$ be a quadratic $R$--space of positive even rank.  \sm

\begin{inparaenum}[\rm (a)] \item \label{deligne-a} There exists an $\uO(q)$--equivariant morphism of $R$--schemes
\[ \ue \co \uL(q) \longto \Spec\big(\Dis(q)\big),  \]
extending the definition  \eqref{Dem1} in the hyperbolic case.
The morphism $\ue$ is projective, smooth, and surjective. It has geometrically connected fibres. \sm

\item\label{deligne-b} Let $f \co \uL(q) \to \Spec(R)$ and $p \co \Spec\big( \Dis(q)\big)\to \Spec(R)$ be the structure morphism of the $R$--schemes $\uL(q)$ and $\Spec\big( \Dis(q)\big)$ respectively. Then $f= p \circ e$ is the Stein factorization of the proper morphism $f$.
    \end{inparaenum}
\end{thm}

\begin{proof}   \eqref{deligne-a}
We define $\ue$ by descent from the hyperbolic case; here are the details. 

We can assume that $(M,q)$ has rank $2n$ and consider the $\uO(q_0)$--torsor $E=\uIsom(q_0,q)$ where $q_0$ is the split hyperbolic form of rank $2n$. Then $\uL(q)$ is the twist of $\uL(q_0)$ by $E$ and $\Dis(q)$ is the twist of $\Dis(q_0)$ by $E$ as well. The morphism $\uL(q_0) \to \Spec\big(\Dis(q_0)\big)$ is $\uO(q_0)$-equivariant by \eqref{eor-a1}.

The morphism $f \co \uL(q) \to \Spec(R)$ is projective and smooth by \ref{quadint}  and \ref{quadsm} respectively, while the morphism $p \co \Spec\big(\Dis(q)\big) \to \Spec(R)$ is finite \'etale because $\Dis(q)$ is a quadratic \'etale $R$--algebra. Then smoothness and projectivity of $e$ follows from $f= p \circ e$ and cancellation, \cite[IV$_4$, (17.3.5)]{EGA} and \cite[II, (5.5.5)]{EGA}, see also \cite[Tag 0C4Q]{St}.

We can also establish surjectivity by descent. 
Thus, without loss of generality, we can assume that $q$ is hyperbolic, in which case $\Dis(q)$ is split. To prove surjectivity in this case, it is further no harm to assume that $R$ is a field \cite[I, 3.6.2]{EGA-neu}. 
In this case $\Spec\big( \Dis(q)\big)$ has two points, while $\uL(q)$ has two connected components, namely two $\SO(q)$--orbits, see e.g.\ \cite[Prop.~3.7]{Co3}, so that $\ue$ is surjective by \ref{eor}\eqref{eor-b}.

To see that $\ue$ has geometrically connected fibres, we can assume that we are over an algebraically closed field, so that $q$ is hyperbolic. The claim then follows from \ref{eor}, \ref{eor-bi} $\Leftrightarrow$ \ref{eor-bii}, keeping in mind that $\uSO(q)$ is the identity component of $\uO(q)$ by \ref{sogsc}\eqref{sogsc-even}. (It also follows from \cite[Tag 03H2]{St} once we have established \eqref{deligne-b}.) \sm

\eqref{deligne-b} By definition in \cite[Tag 03H2]{St}, see also the proof of \cite[Tag 03GY]{St}, the Stein factorization of $f$ has the form
\[
   \uL(q) \xrightarrow{\; f'\; } S' \xrightarrow{\; \pi'\; } S=\Spec(R)
\]
where $\pi'$ is the normalization of $\Spec(R)$ in $\uL(q)$, as defined in \cite[035H]{St}. The scheme $S'$ is the spectrum of a quasi-coherent $\scO_S$--algebra $\scA$. But since $S$ is affine, we know $\scA$ is the quasi-coherent $\scO_S$--algebra associated with an $R$--algebra $A$.
Since $f$ is smooth and projective, the assumptions of \cite[Tag 0BUN]{St} are fulfilled.
\lv{
Because $f$ is smooth, it is flat and lfp, and since $f$ is proper, it is qcqs. But lfp + qcqs is fp = finitely presented by definition in \cite[I, \S3, 1.6]{DG}.

Smooth allows base change. So to check that the fibres are geometrically reduced, we can assume that we are over a field. But a smooth scheme over a field is geometrically reduced by \cite[Tag 056T]{St}.}
Hence, applying loc.\ cit., we know that $A$ is finite \'etale.

By the universal property of normalizations \cite[Tag 035I]{St}, there exists a unique morphism $h \co \Spec(A) \to \Spec (\Dis(q))$ such that the diagram below is commutative:
\[ \vcenter{\xymatrix@C=40pt@R=20pt{
   & S'= \Spec(A) \ar[dr]^\pi \ar[dd]^h\\ \uL(q) \ar[ur]^{f'} \ar[dr]_{\ue} && S= \Spec(R) \\ &  \Spec(\Dis(q)) \ar[ur]_p
}} \quad . \]
We can apply the Fiberwise Isomorphism Criterion \ref{ag}\eqref{ag-d}, to prove that $h$ is an isomorphism. Since the Stein factorization commutes with base change \cite[Tag 03GY]{St}, we can assume that $R=k$ is an algebraically closed field. Then $\Dis(q) = k \times k$ and surjectivity of $\ue$ implies that $\uL(q)= \uL_1 \sqcup \uL_2$ is also a disjoint union of two $k$--schemes, both connected and projective by \eqref{deligne-a}. But in this case, it is standard that $\Spec\big(\Dis(q)\big) = \Spec(k) \sqcup \Spec(k)$ is the Stein factorization of $f$.
\lv{
we know that $f'$ is qcqs, surjective and has geometrically connected fibres, $\pi$ is integral, $f'_*\scO_{\uL(q)} = \scO_{S'}$. Some properties of $h$:
\begin{enumerate}

  \item $h$ is an affine morphism, i.e., $A$ is an integral $\Dis(q)$--algebra,

  \item $h$ is surjective, since $e$ is surjective. What does this mean for the homomorphism $\Dis(q) \to A$?

  \item the fibres of $h$ are geometrically connected, since this holds for $e$ and $f'$.

  \item We know $\Dis(q) = R \times R$, thus $D=D_1 \sqcup D_2$ is a disjoint union of schemes, each isomorphic to $S=\Spec(R)$. By surjectivity of $e$ it follows that $\uL(q)= \uL_1 \sqcup \uL_2$ is also a disjoint union of two $R$--schemes. We can see this differently: the $\uL_i$, $i=1,2$, are the orbits of the $\uSO(q)$--action. Since $h$, rather $S'$, is the normalization of $D$ in $\uL(q)$, it follows from \cite[03GO]{St} that also $S'=S'_1 \sqcup S'_2$ is the disjoint union od two schemes. This is ``obvious'', using that $h$ is induced by a ring morphism $\Dis(q) \to A$.

 \item It suffices to show that $A= \Dis(q)$. From the set-up and using (4), it is enough to show: on a connected smooth projective scheme $X$ over $S$ with structure morphism $f \co X \to S$ the globally defined regular functions are "constant", i.e., $f_* \scO_X = \scO_S$.
\end{enumerate}}
\end{proof}
\sm

We can now treat the missing case of Corollary~\ref{prohigqu}, devoted to the scheme $\Par\big(\uSO(q)\big)$ of parabolic subgroups of $\uSO(q)$.  We use the Dynkin scheme $\Dyn\big(\uSO(q)\big)$, the type morphism $\mathbf{t} \co \Par(G) \to \Of\big(\Dyn(G)\big)$ and its fibres $\Par(G)_t$, reviewed in \ref{pare}\eqref{pare-d}.

\begin{cor}\label{prohla} Let $(M,q)$ be a quadratic $R$--space of constant even rank $2n \ge 4$. We put $G=\uSO(q)$.  \sm

\begin{inparaenum}[\rm (a)] \item \label{prohla-a} Let $(M,q)$ be hyperbolic, thus $G$ has type $\rmD_n$. Then the Lagrangian quadric $\uL(q)$ is $\uO(q)$--isomorphic to the disjoint union of  $\Par(G)_{t_{n-1}}$ and $\Par(G)_{t_n}$ where $t_{n-1} = \rmD_n \setminus \{\al_{n-1}\}$ and $t_n = \rmD_n \setminus \{\al_n\}$. \sm

\item\label{prohla-b} In general, we can identify $\Dyn(G)$ with the disjoint union of $n-2$ copies of $\Spec(R)$ and of $\Spec\big(\Dis(q)\big)$ such that the $R$--scheme $\uL(q)$ becomes $\uO(q)$--isomorphic to the $R$--scheme of parabolic subgroups of $G$ of type $\Dyn(G)\setminus \Spec\big(\Dis(q)\big)$.
\end{inparaenum}
\end{cor}

\begin{proof}
  \eqref{prohla-a} We can argue as in the proof of Corollary~\ref{prohigqu}\eqref{prohigqu-a}: In this case, $\uL(q)(R) \ne \emptyset$, $\Dis(q)$ is split and $\uL(q) = \uL_1 \cup \uL_2$ is a disjoint union of two $R$--schemes with an $\uSO(q)$--action, as we have seen in the proof of \ref{deligne}\eqref{deligne-b}. The identification of the $\uL_i$ as $\Par(G)_{t_j}$ for appropriate $t_j$ follows from Proposition~\ref{prop_hig} and case (IV) of Proposition~\ref{Colem}. \sm

\eqref{prohla-b} Let $q_0$ be the split quadratic form of rank $2n$ and let $G_0=\uSO(q_0)$. Recall that $q_0$ is hyperbolic and that $G_0 = \uSO(q_0)$ is split of type $\rmD_n$. Thus, its Dynkin scheme $\Dyn(G_0)$ is the constant scheme $(\rmD_n)_R$, i.e., the disjoint union of $n$ copies of $\Spec(R)$ indexed by the roots of the Dynkin diagram $\rmD_n$. The Dynkin scheme $\Dyn(G)$ is the twist of the Dynkin scheme $\Dyn(G_0)$ by the torsor $\uIsom(G_0, G)$. But since $G_0$ acts trivially on $\Dyn(G_0)$, the Dynkin scheme $\Dyn(G)$ can be identified with the twist of $\Dyn(G_0)$ by the torsor $\uIsom\big(\Dis(q_0), \Dis(q)\big)$, using the morphism $\ue$ of \eqref{Dem1} and \eqref{eor-a1}. Thus, in view of \eqref{prohla-a}, it follows that $\Dyn(G)$ is isomorphic to the disjoint union of $n-2$ copies of $\Spec(R)$ together with a copy of $\Spec\big(\Dis(q)\big)$.
\end{proof}

\subsection{Example: Quaternions} \label{exqu} Let $Q$ be a quaternion algebra. Thus $Q$ is a unital associative $R$--algebra whose underlying $R$--module has constant rank $4$ and for which there exists a nonsingular (= regular) quadratic form $n_Q$ permitting composition in the sense that $n_Q(ab) = n_Q(a)\, n_Q(b)$ holds for all $a,b\in Q$.

It is known that the identity element $1_Q$ satisfies $n_Q(1_Q) = 1_R$, and is therefore unimodular; for example, this follows from $n_Q(a) = n_Q(a) \, n_Q(1_Q)$ and \eqref{quadco-aa1}. It is also known that quaternion $R$--algebras are precisely the Azumaya $R$--algebras of constant rank $4$, in particular, the quadratic form $n_Q$ is the reduced norm of $A$ in the sense of 
\cite[IV; \S2]{KO} or \cite[4.3]{Salt}. 
We say that $Q$ is split, if it is so as Azumaya algebra, i.e., if it is isomorphic to $\Mat_2(R)$.  The automorphism group scheme $\uAut(Q)$ of an arbitrary $Q$ is a semisimple adjoint $R$--group scheme of type $\rmA_1$ \cite[(3.5.0.82)]{CF}. The Severi-Brauer scheme $\uSB(Q)$ is isomorphic to the scheme of Borel subgroups of $\uAut(A)$.

Let $\tr_Q = b_{n_Q}(1_Q, -)$, which is the reduced trace of the Azumaya algebra $Q$ (the notation $\tr_Q$ instead of $\Trd_{Q/R}$ is traditional). We put
\[ Q\pu = \{ a\in Q : \tr_Q(a)  = 0 \}= (R\, 1_Q)^\perp \quad \text{and} \quad
 n\pu = n_Q|_{Q\pu}, \]
which is a nonsingular quadratic form by \ref{mx_lem}\eqref{mx_lemc}. We also know from loc.\ cit.\ that $Q\pu$ has constant rank $3$. Hence,
the group $\uSO(n\pu)$ is a semisimple adjoint $R$--group scheme of type $\rmA_1$, and therefore isomorphic to $\uAut(Q)$, see Remark~\ref{cimrem}\eqref{cimc-a} and Lemma~\ref{cimc} for a ``concrete'' isomorphism.

\comments{(2025-09-02) The two references above refer to the section on "Cohomoligical consequences" of the old version; not clear at present if this section stays.}

The group scheme $\uAut(Q)$ acts on the quadric $\uQ(n\pu)$ by restriction. The image of this restriction is $\uSO(n\pu)$. This follows for example from Lemma~\ref{cimc}\ref{cimc-b} 

\comments{(2025-09-02) the reference above is a Lemma in the section on "Cohomological Consequences"; it does not use cohomology, and could be put here}

and the fact the even Clifford algebra $\Cli_0\big( Q\pu, n\pu\big)$ is $R$--isomorphic to $Q$, \cite[V, (3.2.4)]{K}.

By 
\ref{prohigqu}\eqref{prohigqu-a}, the quadric $\uQ(n\pu)$ is isomorphic to the scheme of Borel subgroups of $\uSO(n\pu)$. Hence, we have an isomorphism of $R$--schemes
\begin{equation}\label{exqu1}
   \uQ(n\pu) \simlgr \uSB(Q), \end{equation}
which is well-known, see for example \cite[XII, 2.4]{SGA7}. In Corollary~\ref{sqpc}  we construct a ``concrete'' isomorphism \ref{exqu1}. Its proof requires Lemma~\ref{spq} characterizing split quaternions. This lemma complements \cite[22.9]{PRbook} which says that a quaternion algebra $Q$ is reduced (but not necessarily split) if and only if its norm $n_Q$ contains an isotropic vector. 

\comments{(2025-09-03)Lemma~\ref{spq} does not (seem to) contradict \cite[22.9]{PRbook} or \cite[22.19]{PRbook}, since $n_Q$ isotropic does not imply that $n\pu$ is isotropic, only the other implication is true. }

\begin{lem}[Split quaternions]\label{spq} Let $Q$ be a quaternion $R$--algebra. We use the notation of {\rm \ref{exqu}}.  Then
\begin{equation}\label{spq1}  \text{$Q$ is split} \quad \iff \quad \text{$n\pu$ contains an isotropic vector.} \end{equation}
More precisely, given an isotropic vector $u$ of\/ $(Q\pu, n\pu)$, there exists an isotropic $v\in (Q\pu, n\pu)$ such that $c=vu$ is an elementary idempotent of $Q$ satisfying \begin{enumerate}[label={\rm (\roman*)}]
  \item \label{spqi} $Qu=Qc =Q_{11}\oplus Q_{01}$, where the $Q_{ij}$ are the Peirce spaces with respect to the two complementary orthogonal idempotents $(c_1, c_0) = (c, 1_Q - c)= (vu, uv)$;

   \item\label{spqii} the $Q_{ij}$ are free $R$--modules of rank $1$ and are given by  $Q_{11}= Rc$, $Q_{10} = Rv$, $Q_{01} = Ru$, $Q_{00} = R(1_Q-c)$.
\end{enumerate} \end{lem}

\begin{proof} If $Q=\Mat_2(R)$ is split, then $\tr_Q$ and $n_Q$ are the trace and determinant respectively, so that ``$\implies$'' in \eqref{spq1} is obvious. The converse will follow from the remaining statements of the lemma,  which we are going to prove now. \sm

The algebra $Q$  is a composition $R$--algebra in the sense of 
\cite[19.5]{PRbook}. 
 Thus for $a,b\in Q$ it satisfies
\begin{align} \label{sqp2}
  \begin{split}
    & a^2 - \tr_Q(a) a + n_Q(a)1_Q = 0, \\
    & aba = b_{n_Q}(a, \ol b) a - n_Q(a) \ol b,
 \end{split}
\end{align}
where $\ol b = \tr_Q(b) - b$. Since $\tr_Q = b_{n_Q}(1_Q, - )$, the submodule $Q\pu$ is complemented by \ref{mx_lem}. 
Hence $u$ is also an isotropic vector of $(Q, n_Q)$. By regularity of $n_Q$ and \ref{quadrepII}, there exists $v\in Q$ satisfying $1= b_{n_Q}(v, \ol u) = \tr_Q(vu)$. Because $n_Q(vu) = n_Q(v) n_Q(u)= 0$, the element $c=vu$ is an elementary idempotent of $Q$, as defined in \ref{eica}. It then follows from the first equation in \eqref{sqp2} that $u^2=0$, hence $cu=vu^2=0$, and $uc=uvc=u$. Thus $u\in Q_{01}$ and $Qu=Qc$ since $Qu=Quc\subset Qc = Qvu \subset Qu$. This in turn implies $Qu=Q_{11} \oplus Q_{01}$.

We decompose $v$ into its Peirce components,  $v= v_{11} + v_{10} + v_{01} + v_{00}$. Since then $c=vu = v_{10} u$, it is no harm to replace $v$ by $v_{10}$. Because $1=\tr_Q(c) = b_{n_Q}(v,\ol u)$ the element $v$ is an isotropic vector in $Q_{10} \subset Q\pu$. Also observe
\[ uv + vu = \tr_Q(u) v + \tr_Q(v)u  - b_{n_Q}(u,v) 1_Q = b_{n_Q}(v, \ol u) 1_Q = 1_Q \]
proving $vu = 1_Q - c$. It remains to show \ref{spqii}. To this end, we have $cxc=b_{n_Q}(c,\ol x) c \in Rc$ by \eqref{sqp2},  so that $Q_{11} = Rc$ and then, by symmetry $Q_{00} = R(1_Q - c)$. Finally, $vxv= b_{n_Q}(v,\ol x) v \in Rv$ proving $Q_{10} = Rv$ and, analogously, $Q_{01} = Ru$. That the Peirce spaces are free,  follows from unimodularity of $c,u,v$ and $1_Q - c$. \end{proof}
\ms

We can now give an elementary proof of the isomorphism \eqref{exqu1}.
\begin{cor}\label{sqpc} We use the notation of\/ {\rm \ref{spq}}: $Q$ is a quaternion $R$--algebra, $Q\pu = \{ x\in Q: \tr_Q(x) = 0 \}$, $n\pu = n_Q|_{Q\pu}$. Recall that $(Q\pu, n\pu)$ is a quadratic space. The natural transformation $\underline{\Psi}\co \underline{\rmQ}(n\pu)\to \underline{\SB}(Q)$ of $R$--functors, defined on $R'$--points by
\[ \underline{\Psi}(R') \co \underline{\rmQ}(n\pu) (R') \to \underline{\rm SB}(Q)(R'), \quad L \mapsto Q\cdot L, \]
is an isomorphism, and hence induces an isomorphism of $R$--schemes,
\begin{equation} \label{sqpc0}
          \Psi \co \uQ(n\pu) \simlgr \uSB(Q),
\end{equation}
which is equivariant with respect to the action of $\uAut(Q)$.
\end{cor}

\begin{proof} 
(I) {\em $\underline{\Psi}(R')$ is well-defined:} For simpler notation, let $R'=R$.
We need to prove that for $L\in \underline{\rmQ}(n\pu)(R)$ the $R$--module $M=Q/QL$ is projective of rank $2$.
Clearly, $M$ is a finitely  generated $R$--module. By \cite[II, \S5.3, Thm.~2]{BAC}, it then suffices to establish that for any $\p \in \Spec(R)$ the $R_\p$--module
$M_\p = Q_\p/(Q_\p L_\p)$ 
is free of rank $2$. To do so, we can replace $R_\p$ by $R$ and thus assume that $R$ is a local ring. Then $L$ is free and therefore spanned by an isotropic vector $u\in Q\pu$. In particular, we are in the setting of Lemma~\ref{spq}: $QL=Qu=Qc$, $c$ an idempotent, is free of rank $2$ and so is $Q/QL\cong Q(1_Q - c)$. \sm

\lv{
(II) {\em $\underline{\Psi}(R')$ is injective.} This will follow from
\begin{equation}
  \label{sqpc1} Q\pu \cap Q\cdot L = L
\end{equation}
Clearly, $L \subset Q\pu \cap Q\cdot L$. This inclusion is an equality if it is so after localization in every maximal ideal $\m$ of $R$. But over $R_\m$ we are again in the setting of Lemma~\ref{spq}, where \eqref{sqpc1} becomes the obvious $Q\pu \cap (Rc \oplus Ru) = Ru$. \sm

(III) {\em $\underline{\Psi}(R')$ is surjective.} Let $I\in \underline{\SB}''(Q)(R)$, i.e., the $Q$--module $Q/I$ is projective of rank $2$. By \ref{sep}\eqref{sep-li}, there exists an idempotent $c\in Q$ such that $I = Qc$ and $Q/I \cong Q(1_Q-c)$ as $Q$--modules. It follows that $I$ is projective of rank $2$ as $R$--module and then that $I=Q_{11} \oplus Q_{21}$ (Peirce decomposition with respect to $c$) with $Q_{21}$ projective of rank $1$ as $R$--module satisfying $n_Q(Q_{12}) = 0 = \tr_Q(Q_{21})$. Thus $Q_{21} \in \underline{\rmQ}(n\pu)(R)$. We have $\big(\underline{\Psi}(R)\big)(Q_{21}) = Q \cdot Q_{21} \subset I$. That this inclusion is an equality, follows as in the previous step by localization.  \sm

(IV) By Yoneda, the isomorphism $\underline{\Psi}$ induces an isomorphism $\Psi \co  \uQ(n\pu) \simlgr \uSB(Q)$ of $R$--schemes. Equivariance of $\Psi$ follows from the fact that an automorphism of $Q$ leaves $\tr_Q$ and $n_Q$ invariant, cf.\ \eqref{sqp2}. }

(II) Using Lemma~\ref{spq} and Lemma~\ref{sep-li}, it is now straightforward to show that $\underline{\Psi}(R')$ is a bijection for all $R'\in \Ralg$ and hence induces an isomorphism $\Psi$ as claimed in \eqref{sqpc0}.

However, it is quicker and more instructive to proceed as follows. The natural transformation $\underline{\Psi}$ gives rise to an $\uAut(Q)$--equivariant morphism $\Psi \co \uQ(n\pu) \simlgr \uSB(Q)$ of $R$--schemes. Both schemes are $\uAut(Q)$--homoge\-neous with respect to the flat topology.
A fortiori, $\Psi$ is a morphism between $\uGL_1(Q)$--homogeneous spaces. To show that $\Psi$ is an isomorphism, we can localize with respect to the flat topology and then suppose that $Q=\Mat_2(R)$, so that $\tr_Q$ and $n_Q$ are the usual trace and determinant. Let us consider the isotropic line $L_0 = \ppmatrix R 0 0 0 $ and the corresponding left ideal $I_0 = Q \cdot L_0  = \ppmatrix R 0 R 0$.
The morphism $\Psi$ becomes
\[
     \uGL_1(Q)\big/ \Stab_{\uGL_1(Q)}(L_0)\longto \uGL_1(Q)\big / \Stab_{\uGL_1(Q)}(I_0)
\]
where $\Stab_{\uGL_1(Q)}(L_0)$ and $\Stab_{\uGL_1(Q)}(I_0)$ are the stabilizers of $L_0$ and $I_0$ respectively under the corresponding actions. An easy computation shows $\big(\Stab_{\uGL_1(Q)}(L_0)\big)(R') = \big(\Stab_{\uGL_1(Q)}(I_0)\big)(R') = \ppmatrix {(R')\ti} 0 R' {(R')\ti} $ for each $R$--algebra $R'$, proving that $\Psi$ is indeed an $\uAut(Q)$--equivariant isomorphism. \end{proof}

\newpage

\section{Springer's Odd Degree Theorem for LG rings}\label{sec:springer}

In this section we are concerned with the following question: given a quadratic space $(M,q)$ over a base ring $R$
and some $S\in \Ralg$ such that the base change $(M,q)_S$ is an isotropic quadratic space over $S$, under which conditions is $(M,q)$ already isotropic? 

In Theorem~\ref{spo} 
we will give a positive answer to this question for $R$ an LG ring and for certain $S\in \Ralg$ which force $S$ to be an LG ring too. In this context it is appropriate to recall \ref{isotrop}\eqref{isotrop-unifap}: if $R$ is a unimodular ring, e.g., an LG ring or even a semilocal ring, a quadratic $R$--module $(M,q)$ is isotropic if and only if $(M,q)$ contains an isotropic vector. In the proof of Theorem~\ref{spo} we will therefore consider isotropic vectors, rather than dealing with arbitrary totally isotropic and complemented submodules. We will also employ the terminology that ``$q$ is $S$--isotropic''  to mean that the quadratic $S$--module $(M,q)_S$ is isotropic.

We refer to our Theorem~\ref{spo} as ``Springer's Odd Degree Theorem'' because it was first proven by Springer in \cite{springer-qf} for $R$ and $S$ being fields of characteristic $\ne 2$. This theorem has turned out to be a fundamental result in the theory of quadratic forms for which  several generalizations exist, some of them reviewed in the introduction of our paper \cite{GN-Sp}. In that paper we have proven Theorem~\ref{spo} for $R$ and $S$ being semilocal. The proof given here is different from the one of \cite{GN-Sp}, except that we will use \cite[Lem.~2.2]{GN-Sp}, dealing with the case of a quadratic space of constant rank $2$. For convenience, we restate that lemma below as Lemma~\ref{sptwo}. 

\comments{(2026-01-30) Added that $S$ is finite locally free. Before,we just said that {\tt $S$ is  a finite one-generated $R$--algebra of odd degree}. Here, ``degree'' only makes sense for finite projective =finite locally free. But this may not be obvious, in particular since Ferrand \cite{Fer} studies finite $R$--algebras that are fppf-locally one-generated, but not necessarily free} 

\begin{thm}[Springer's Odd Degree Theorem]\label{spo} Let $(M, q)$ be a  quadratic space over an LG ring $R$, and let $S$ be a one-generated and finite locally free $R$--algebra of odd degree. \sm 

\begin{inparaenum}[\rm (a)] \item \label{spo-a} If $q$ is $S$-isotropic, then $q$ is $R$-isotropic. Equivalently, if the reductive $S$--group scheme $\uSO(q)_S$ is isotropic, then already the $R$--group scheme $\uSO(q)$ is isotropic. \sm 

\item \label{spo-b} Moreover, the conditions \ref{equi-Breg}, \ref{equi-D} and \ref{equi-E} of {\rm \ref{equi}} hold: \end{inparaenum}

\begin{enumerate}
 \item[\ref{equi-Breg}] Let $(M,q)$ be a quadratic $R$--space and let $(M_1,q_1)$ be  a\/ {\em regular} quadratic module such that $(M,q)_S$ contains $(M_1, q_1)_S$. Then $(M,q)$ contains $(M_1, q_1)$. \sm  

\item[\ref{equi-D}] If $(M,q)$ and $(M',q')$ are quadratic $R$--spaces such that $(M,q)_S \cong (M', q')_S$, then $(M,q) \cong (M', q')$. \sm

\item[\ref{equi-E}] The restriction homomorphisms $\hWq(R) \to \hWq(S)$ and
$\Wq(R) \to \Wq(S)$ of  \eqref{wgg1} and \eqref{wgg2} are injective. 
\end{enumerate}\end{thm}


The proof of Theorem~\ref{spo} will be given in \ref{spopr}, using the concept of a Springer functor, \ref{desf}, which we prove in \ref{respr} to be representable. The basic idea of Springer functors can be traced back to \cite[Prop.~1.1]{panin-rehmann}. In that paper the authors prove \ref{spo} for certain noetherian local domains $R$ containing $2\in R\ti$. The approach was refined in \cite{PP} to allow certain semilocal rings $R$, but still with $2\in R\ti$. We have proven \ref{spo} for arbitrary semilocal rings $R$ in \cite[Thm.~2.1]{GN-Sp}; the proof here is different from the one in \cite{GN-Sp}.  

We immediately present two corollaries of Theorem~\ref{spo}. The first eliminates the condition that $S$ be one-generated in two circumstances.  We recall that a ring satisfies the primitive criterion if it is LG and all its residue fields are infinite, \ref{revLGG}.

\comments{(2026-02-19) The new formulation of Cor.~\ref{spo-LGG} makes it clearer what "drives the game", namely \ref{spo-LGG}\ref{spo-LGGi}. I immediately proved the corollary since its proof is very short. \sm 

Related thoughts: An easy argument in \cite[1.2]{Fer} shows that every \'etale $R$--algebra $E$ is ``locally simple'' in his sense, i.e., there exists a faithfully flat and finitely presented $R$--algebra $R'$ such that $E'=E \ot_R R'$ is one-generated as $R'$--algebra. If one can show that $R'$ is LG, then an application of \ref{spo} shows that $q_{R'}$ is isotropic, which one could call ``fpqc-isotropic''. }
\pcomments{(202603-01) Corolary~\ref{spo-LGG} est parfait.}

\begin{cor}\label{spo-LGG} Let $(M,q)$ be a quadratic space over an LG ring $R$, and let $S\in \Ralg$ is finite locally free of constant odd rank. 
Furthermore, assume that 
\begin{enumerate}[label={\rm (\roman*)}]
  \item \label{spo-LGGi} there exists a finite locally free $T\in \Ralg$ of constant odd rank such that $S\ot_R T$ is a one-generated $R$--algebra.    
\end{enumerate}
Then $q$ is $S$--isotropic if and only if $q$ is $R$--isotropic.

Furthermore, assume in addition that $\uPrim(S) \ne \emptyset$, e.g., assume $S$ is finite \'etale. Then the condition {\rm \ref{spo-LGGi}} is fulfilled in the following two circumstances: \sm 

\begin{inparaenum}[\rm (a)]
 \item\label{spo-LGGa} $R$ satisfies the primitive criterion {\rm \ref{revLGG}}. 
 
 \item \label{spo-LGGb} $R$ is semilocal and $S$ is finite \'etale.   
\end{inparaenum}
\end{cor}

\begin{proof} Assuming \ref{spo-LGGi}, the $R$--algebra $S\ot_R T$ is finite locally free of constant odd rank and $q$ is $(S\ot_R T)$--isotropic. 
Hence Theorem~\ref{spo} applies and proves the first claim. 

For \eqref{spo-LGGa} and \eqref{spo-LGGb} recall \ref{prsch}\eqref{prschb} that $\uPrim_R(S) \ne \emptyset$ whenever $S$ is finite \'etale. In the situation \eqref{spo-LGGa} we know by  \ref{loog}\ref{loogc} that $S$ is already one-generated, so that $T=R$ does the job. In the situation \eqref{spo-LGGb}, the existence of $T\in \Ralg$ as in \ref{spo-LGGi} is proven in \cite[Prop.~7.3]{BFP}, see \ref{bfp-prop} for an account. 
\end{proof}
\sm 

The second corollary of Theorem~\ref{spo} concerns hermitian spaces. As for quadratic forms, we will say that a hermitian space over an LG ring is isotropic if it contains an isotropic vector. 

\begin{cor}\label{sprher} Let $R$ be an LG ring, let $R'/R$ be a quadratic \'etale extension, let $(M', h')$ be a hermitian space over $R'/R$, and let $S$ be a finite one-generated $R$--algebra of odd degree. \sm 

\begin{inparaenum}[\rm (a)]
\item\label{sprher-a} If the hermitian space $(M', h')_{R'\ot_R S}$ over the quadratic \'etale extension $(R'\ot_R S)/S$ is isotropic, then already $(M', h')$ is isotropic.\sm 
    
\item\label{sprher-b}  Let $(M_1', h_1')$ be another hermitian space. 
Then 
\[ (M', h')_{R'\ot_R S} \cong (M'_1, h'_1)_{R'\ot_R S} \quad \implies \quad 
    (M', h') \cong (M'_1, h'_1).
\]
\end{inparaenum}
\end{cor}

We will prove 
Corollary~\ref{sprher} in \ref{herpr} after a review of hermitian spaces in \ref{herev} and \ref{herle}. Corollary~\ref{sprher} was first proven in \cite{BFL} in the setting of $\veps$--hermitian spaces over fields.
\comments{(2025-10-13) Need to say more here about other papers proving \ref{sprher}. } 
\ms 

Springer's Odd Degree Theorem for quadratic spaces of constant rank $2$ holds for more general base rings than LG rings. We have shown this in \cite[Lemma~2.2]{GN-Sp}, written in the context of semilocal rings (the semilocality assumption is perhaps not clear from the statement of \cite[2.2]{GN-Sp}). In fact, the proof of \cite[2.2]{GN-Sp} works in the more general setting of rings $R$ with $\Pic(R) = 0$, thus in particular LG rings. We state this as Lemma~\ref{sptwo}: 


\begin{lem}[Rank $2$] \label{sptwo} Let $R$ be a base ring with $\Pic(R) = 0$ and let $S$ be a finite locally free $R$--algebra of constant odd degree.  Furthermore, let $(M,q)$ be a quadratic $R$--space of constant rank $2$. If $q$ is $S$--isotropic, then $q$ is $R$--isotropic. More precisely, $(M,q) \cong \HH(R)$ is a hyperbolic plane.  
\end{lem}

%
%

\subsection{Some schemes related to a quadratic $R$--module $(M,q)$}\label{sreq}
For the theory leading up to the proof of Theorem~\ref{spo} it is not always necessary that $R$ be an LG ring. Hence, until further notice, we let $R$ be an arbitrary base ring and $(M,q)$ a quadratic $R$--module. 

We start by defining some of the schemes we will be using.
As usual, $\uW(M)$ is the affine scheme associated with the finite locally free $R$--module $M$, \ref{ag}\eqref{ag-ex}. Its functor of points is $T\mapsto M_T = M \ot_R T$, $T\in \Ralg$. 
Also recall (Lemma~\ref{spreq-a}) the scheme $\uW(M) \setminus \{0\} =: \uW(M)_u$ representing the $R$--functor which assigns to $T\in \Ralg$ the set $(M\ot_R T)_u$ of unimodular vectors of the $T$--module $M \ot_R T$. 
\sm

\begin{inparaenum}[(a)] \item \label{sreq-b} ({\em The scheme $\bfZ_q := \uV_{q,0}$ of zero vectors}) Recall from Lemma~\ref{q-Faser} that $\uZ_q$ represents the $R$--functor $T\mapsto \{v\in M_T: q_T(v) = 0 \}$, $T\in \Ralg$. \sm 

\item\label{sreq-y} ({\em The scheme $\bfY_q$ of isotropic vectors}) 
 We define
\[  \uY_q = \uW(M)_u \cap \uZ_q = \uW(M)_u \times_{\uW(M)} \uZ_q \] 
and call $\uY_q$ the {\em affine quadric associated with $q$}: 
\begin{equation} \label{aqy-c1} \vcenter{ 
\xymatrix@C=50pt{\ar @{} [dr] |{\small\qed} %
\uY_q \ar[r] \ar[d] & \uZ_q \ar[d] 
 \\ \uW(M)_u \ar[r] & \uW(M)  } }  .  \end{equation}
By definition, 
the $T$--points of $\uY_q$ are
\begin{equation} \label{aqy1} \begin{split}
 \uY_q(T) &= \{ v\in M_T: q(v) = 0, \text{ $v$ unimodular} \}
   \\ &= \{v \in M_T: \text{ $v$ isotropic} \} .
\end{split}\end{equation} 

\item\label{sreq-x} ({\em The projective quadratic $\uX_q$}) Let $\uX_q$ be the projective quadric associated with $q$, denoted $\uQ(q)$ in \S\ref{sec:quadrics}, see \ref{hrq}. Recall that its $R$--functor of points is 
    \[ 
      T \mapsto  \uX_q(T) = \{ L \in \PP(M\ch)(T): q_T(L) = 0 \}, 
        \]
and that $\PP(M\ch)(T)$ can be identified with the complemented invertible submodules of $M$. 
\end{inparaenum}
\ms

We now describe some geometric properties of the schemes introduced above, first in Lemma~\ref{spreq} for quadratic modules, and then in Lemma~\ref{aqyle} for quadratic spaces. Part of the results are special cases of Lemma~\ref{q-Faser}, which we however repeat for the convenience of the reader. 


\begin{lem}[Some properties of the schemes of {\ref{sreq}}, $q$ arbitrary]\label{spreq} Let $(M,q)$ be a faithful quadratic $R$--module satisfying $\Cont(q) = R$.  \sm 
  
\begin{inparaenum}[\rm (a)]
\item\label{spreq-z} The $R$--scheme $\uZ_q$ is an affine scheme of finite presentation. The open subscheme $\uZ_q\rmsm$ of smooth points of $\uZ_q$, 
    {\rm \cite[Tag 01V5]{St}},  
    is given by 
    \begin{equation} \label{spreq-z1}  
    \uZ_q \rmsm (T) = \{v\in M_T: q_T(v) = 0, \, b_{q_T}(v, \cdot) \text{ is surjective}\}
\end{equation}
for $T\in \Ralg$. In particular, $\uZ_q\rmsm $ is an open subscheme of $\uY_q$. 
\sm

\item \label{spreq-y}  The $R$--scheme $\uY_q$ is quasi-affine and of finite presentation. Moreover, 
\begin{equation} \label{spreq-y1}
\text{$\uY_q$ is smooth}\quad \iff \quad \text{$q$ is nonsingular.}
\end{equation}   

\item\label{spreq-x} The $R$--scheme $\uX_q$ is projective and of finite presentation. If $\rank_R M  \ge 2$, then $\uX_q$ is a smooth $R$--scheme if and only if $q$ is nonsingular. \sm 
    
\item \label{spreq-t} With respect to the canonical action of\/ $\GG_m$ on $\uY_q$, the morphism
\begin{equation}  \label{spreq-t1}
\vphi \co \uY_q \to \uX_q, \quad v (\in M_T) \mapsto Tv =: \underline{v}
\end{equation}
makes $\uY_q$ a $\GG_m$--torsor over $\uX_q$ for the Zariski topology.     
\end{inparaenum}
\end{lem}


\begin{proof} \eqref{spreq-z} That $\uZ_q$ is an affine, finitely presented $R$--scheme whose smooth locus is given by \eqref{spreq-z1}, is a special case of Lemma~\ref{q-Faser}. Any $v\in \uZ\rmsm_q(T)$ is isotropic by \ref{isotrop}\eqref{isotrop-bc}, implying that the open subscheme $\uZ_q\rmsm$ of $\uZ_q$ is in fact an open subscheme of $\uY_q$. \sm 

\eqref{spreq-y} The open immersion $\bfY_q \hookrightarrow \bfZ_q$ is quasi-affine and of finite presentation by base change from $\uW(M)_u \to \uW(M)$, \cite[Tags 01SO, 01TS]{St}. 
The composition $\bfY_q \to \bfZ_q \to \Spec(R)$ is therefore also quasi-affine and of finite presentation, \cite[Tags 01SN, 01TR]{St}.

To prove \eqref{spreq-y1}, let us first assume that $\uY_q$ is smooth. Then $\uZ_q\rmsm = \uY_q$ by \eqref{spreq-z}. Moreover, any base change, for example $\uY_q \times_R F$ for an $R$-field $F$, is smooth too. By definition of nonsingularity in \ref{qf}\eqref{qf-ns}, it is then no harm to assume that $R=F$ is a field, in which case we need to prove that $\rad(q) = 0$. Suppose there exists a non-zero $x\in \rad(q)$. Then $x$ is an isotropic vector, hence $x\in \uY(q_F) = \uZ_q\rmsm(F)$. But then $b_q(x, \cdot) = 0$ and $b_q(x, \cdot )$ surjective forces $M=0$, contradicting faithfulness of $M$.     

Conversely, assume $q$ is nonsingular. We know $\uZ_q\rmsm(T) \subset \uY_q(T)$ for any $T\in \Ralg$ by \eqref{spreq-z}. The other inclusion follows from \ref{isotrop}\eqref{isotrop-bc} saying that for a nonsingular $q$ a vector $v\in M_T$ with $q_T(v) = 0$ is unimodular (= isotropic) if and only if $b_{q_T}(v, \cdot)$ is surjective. 
\sm 

\eqref{spreq-x} These properties are special cases of \ref{quadint} and \ref{quadsm}. \sm 

\eqref{spreq-t} Omitted. 
\end{proof}

\begin{lem}[More properties of the schemes of {\ref{sreq}}, $q$ nonsingular]\label{aqyle} Let $(M,q)$ be a quadratic $R$--space. 
 \sm

\begin{inparaenum}[\rm (a)] \item\label{aqyle-a} The scheme $\uY_q$ is smooth and coincides with the scheme $\uZ_q\rmsm$ of smooth points of $\uZ_q$. 
\sm 

\item \label{aqyle-b} We have 
\begin{align}  \label{aqyle-b1} 
 \uY_q(R) \ne \emptyset &\iff (M,q) \text{ contains a hyperbolic plane.} \\
\label{aqyle-b2} 
     \uY_q \ne \emptyset &\iff \rank M \ge 2. 
\end{align}%

\item \label{aqyle-c} Suppose $\rank M \ge 2$. \end{inparaenum} 

  \begin{enumerate}[label={\rm (\roman*)}]

    \item \label{aqyle-ci}  The $R$--schemes $\uZ_q$ and $\uY_q$ have geometrically integral fibres. \sm

   \item \label{aqyle-ciii}   If $M$ has constant rank $r\ge 2$, the geometric fibres of $\uZ_{q}$ and $\uY_q$ have dimension $r-1$, and $\uY_{q}$ is smooth of relative dimension $r-1$.

 \item \label{aqyle-cii}   The scheme $\uY_q$ is universally schematically dense in $\uZ_{q}$. \sm

\end{enumerate}
\end{lem}

\begin{proof} \eqref{aqyle-a} follows from \eqref{spreq-z} and \eqref{spreq-y} of \ref{spreq}. 
\sm 

\eqref{aqyle-b} By \ref{isotrop}\eqref{isotrop-bc}, $(M,q)$ contains an isotropic  vector if and only if $(M,q)$ contains a hyperbolic plane, which shows \eqref{aqyle-b1}. To see \eqref{aqyle-b2}, assume $\uY_q \ne \emptyset$. Then $(M,q)_T$ is isotropic for some $T\in \Ralg$, and therefore $\rank M_T\ge 2$ by \ref{isotrop}\eqref{isotrop-bc}, which implies $\rank M \ge 2$.  Conversely, if $\rank M \ge 2$, then $(M,q)_T$ is isotropic for some $T\in \Ralg$, e.g., for some algebraically closed field $T$. \sm
 
\eqref{aqyle-c}  is a special case of Lemma~\ref{q-Faser}\ref{q-Faseri}. 
\end{proof}

\subsection{Example: Rank $1$}\label{aqyex} 
Let $(M,q) = (R, \lan 1 \ran_q)$. Then the ring of global functions of $\uZ_q$ is $R[\uZ_{q}] = R[X]/(X^2)$, which is not an integral domain. This explains the rank assumption in Lemma~\ref{aqyle}\eqref{aqyle-c}. Also, if $2R=0$, then $\uY_{q}\rmsm  = \emptyset$, while $\uZ_{q}\ne \emptyset$ in general, see Lemma~\ref{aqyle}\eqref{aqyle-b}.

\subsection{Setting in case $\uY_q(S) \ne \emptyset$ for some $S\in \Ralg$}\label{sunn} We assume that $(M, q)$ is a quadratic space for which there exist $S\in \Ralg$ and $v\in \uY_q(S)$. 
By \ref{isotrop}\eqref{isotrop-bc} we know that $v$ is part of a hyperbolic pair $(u,v)$. Moreover, by \ref{orthLG}\eqref{orthLG-a} and \eqref{qf-perp1}  
\begin{equation} \label{sunn1} 
    M = (Su \oplus Sv) \perp M', \quad \text{$q' = q|_{M'}$ is nonsingular}.
\end{equation}
We fix these data in the following. \sm 

The definitions and results of \ref{sreq}, \ref{spreq} and \ref{aqyle} apply to the quadratic space $(M,q)_S = (M_S, q_S)$ over $S$. For notational purposes we will write $\uY_{q,S} = \uY_{q_S}$ and $\uX_{q,S} = \uX_{q_S}$. 
\sm 

By \ref{gmu}\eqref{gmu-b} and \ref{gmu}\eqref{gmu-c} the stabilizers 
   \[ \Stab_{\uSO(q_S)} (Su) \quad \text{and}\quad 
       \Stab_{\uSO(q_S)} (Sv) 
   \]
are opposite parabolic subgroups of $\uSO(q_S)$, and by \ref{gmu}\eqref{gmu-d} the unipotent radical of $\Stab_{\uSO(q_S)}(Su)$ is a vector group,
\begin{equation}\label{sunn2}
 \rad^u\big(\Stab_{\uSO(q_S)}(Su) \big) \cong \uW(M').
\end{equation}  
Its $A$--points, $A\in \Salg$, are uniquely determined by vectors $n'\in M'_A$ as follows :
\begin{equation} \label{sunn3}\begin{split}
 u_A &\mapsto u_A, \qquad   m'  \mapsto m' -b_q(m', n') u_A  \\
  v_A &\mapsto -q_A(n') u_A + v_A + n'. 
\end{split}\end{equation}
By \cite[XXVI, 4.3.2(b), (vi$^\prime$)]{SGA3}, the orbit map
\[ i_{S v}\co \rad^u\big(\Stab_{\uSO(q_S)}(Su)\big) \to \uX_{q,S}, \quad
   g \mapsto g. (Sv)
\]   
is an open immersion; let $V$ be its image. Thus, $V$ is an open subscheme of $\uX_{q,S}$ isomorphic to $\uW(M')$. We let $U$ be the preimage  of $V$ under the $\GG_{m,S}\,$--torsor $\vphi \co \uY_{q,S} \to \uX_{q,S}$ and
$j \co U \to \uY_{q,S}$ the canonical open immersion. 

We observe that the restriction of $\vphi$  makes $U \to V$ a $\GG_{m,S}$--torsor. It is a trivial torsor, 
\begin{equation}  \label{sun1} U \cong \GG_{m,S} \times_S V \cong \GG_{m,S} \times_S \uW(M')
\end{equation}
as $S$--schemes, since the orbit map 
\[ 
     \rad^u\big(\Stab_{\uSO(q_S)}(Su)\big) \to \uY_{q,S}, \quad g\mapsto g.v
\] 
is a lift of the morphism $i_{Sv}$ and induces a splitting of $U \to V$. In particular, \eqref{sun1} shows that $U$ is (isomorphic to) a principal open subscheme of $\uW(S \oplus M)$. By \eqref{sunn3}, for $A\in \Salg$ the $A$--points of $U$ are determined by $(x,n') \in A\ti \times M'_A$ as follows
\begin{equation} \label{sunn5}
    -x\me q_A(n')u_A + xv_A + n'. 
\end{equation}  

\comments{(2026-03-27) There was a misprint in the previous (pre-2026-03-27) version of the formula \eqref{sunn5}; it read $-xq_A(n')u_A + xv_A + n'$ which is not isotropic in general. Corrected here and again in  \eqref{sunres0}.}

To summarize:
\[\xymatrix@C=50pt{
  \GG_{m,S} \times_S V \ar[d] \ar[r]^\cong &U \ar[r]^j \ar[d] & \uY_{q,S}   
      \ar[d]_\vphi
 \\
\uW(M')  \ar@/_1pc/[rr]|{\; i_{Sv}\; } \ar[r]^>>>>>>>>>\cong & V \ar[r]& \uX_{q,S} 
   }\]
\ms

We will pass from the $S$--schemes introduced above to $R$--schemes using Weil restriction, see \ref{weilres} for a review. 

\begin{lem}\label{sunres} In the setting of {\rm \ref{sunn}}, let $S$ be a finite projective $R$--algebra. Then the Weil restrictions $\frR_{S/R}(\cdot)$ of $U$ and $\uY_{q,S}$ are representable by finitely presented $R$--schemes. In particular, 
\begin{equation}\label{sunres1} 
    \frR_{S/R}(U) \cong \frR_{S/R}(\GG_m) \times \uW\big( \frR_{S/R}(M')\big) 
\end{equation}
is a smooth affine $R$--scheme whose $T$--points, $T\in \Ralg$, are 
\begin{equation}  \label{sunres0} \begin{split}
   U(S\ot_R T) 
   = &\{ - x\me q_{S\ot T}(n')u_{S\ot T} + x v_{S\ot T} + n' : \\
   &\qquad     x \in (S\ot T)\ti, \, n'\in M'_{S\ot T} \}
\end{split}\end{equation}
where $(u_{S\ot T}, \, v_{S\ot T})$ is the image of the $(u, v)$ under the base change $S \to S\ot T$. Moreover, the Weil restriction $\frR_{S/R}(j) \co \frR_{S/R}(U) \to \frR_{S/R}(\uY_{q,S})$
is a quasi-compact open immersion. \end{lem}

\begin{proof} The $S$--schemes $U$ and $\uY_{q,S}$ are affine, respectively quasi-affine of finite presentation by \eqref{sun1} and \ref{spreq}\eqref{spreq-y}. Hence their Weil restrictions exist by \ref{weilres}\eqref{weilres-a} and are finitely presented as $R$--schemes by \ref{weilres}\eqref{weilres-b0}. The formula \eqref{sunres1} follows from \eqref{sun1} and \eqref{weilres-3}, and implies smoothness of  $ \frR_{S/R}(U)$. The formula \eqref{sunres0} is a consequence of \eqref{sunres1} and \eqref{sunn3}.  The last claim is a special case of \ref{weilres}\eqref{weilres-d}. 
\end{proof}
\sm

Our proof of Springer's Theorem uses $R$--schemes representing $R$--functors which we will introduce in \ref{desf} in a setting again used again in the definition of Knebusch functors. 
We present this next. 

\subsection{A preliminary setting}\label{aps}  We consider data 
\[ (S,a, q) \]
consisting of 
\begin{enumerate}[label={\rm (\roman*)}] 
 \item \label{aps-i} an $S\in \Ralg$ which is free of rank $d\in \NN_+$ with basis 
 $1_S, a , \ldots, a^{d-1}$, thus one-generated, and 
 
\item \label{aps_ii} a quadratic $R$--module $(M,q)$.  
\end{enumerate}
Regarding \ref{aps-i},  we recall Lemma~\ref{genolemLG}: any $S\in \Ralg$ which is one-generated and finite locally free of constant rank $d$ has the form \ref{aps-i} for some $a\in S$. Let $P = \Pc_{S/R}(a; X)\in R[X]$ be the characteristic polynomial of $a$, a monic polynomial of degree $d$. We can and will identify 
\[ S = R[X]/(P) , \qquad a = X + (P).\] 
Given $T\in \Ralg$, we let $P_T$ be the canonical image of $P$ in $T[X]$. Then $S\ot_R T \cong T[X]/P_T$ is free of rank $d$ as $T$--module with $T$--basis \[ a_T^0 = 1_{S\ot T}, a_T = a \ot 1_T, \ldots , a_T^{d-1}= a^{d-1} \ot 1_T.
\]

We will consider $M \ot_R S$ as $S$--module as well as $R$--module by restriction of scalars. To be more precise, this is the finite projective $R$--module $\frR_{S/R}(M\ot_R S)$ in the latter case. The $S$--module $M\ot_R S$ gives rise to an $S$--functor and $S$--scheme, both   here denoted as $\uW_S(M\ot_R S)$ and  defined in \ref{ag}\eqref{ag-ex}. Analogously, the $R$--module $\frR_{S/R}(M\ot_R S)$ gives rise to an $R$--functor and $R$--scheme, which we here denote as $\uW_R\big(\frR_{S/R}(M\ot_R S)\big)$, or $\uW_R(M\ot_R S)$ to simplify. 
Observe
\begin{equation}\label{aps-0}
 \frR_{S/R}\big( \uW_S(M\ot_R S)\big) = \uW_R\big(\frR_{S/R}(M\ot_R S)\big)
   = \uW_R(M\ot_R S) 
\end{equation} 
by \eqref{weilres-3}. For $T\in \Ralg$ we have    
\[ \uW_R(M\ot_R S) (T)  = M \ot_R S \ot_R T .
 \]
Any $w\in M\ot_R S \ot_R T$ can be uniquely written in the form 
\begin{equation}\label{aps-1}
w = \textstyle \sum_{i=0}^{d-1} m_i a_T^i, \quad m_i \in M_T = M\ot_R T.
\end{equation}  
By functoriality of the representation \eqref{aps-1}, we get a morphism of $R$--schemes, 
\[ \pi_{d-1} \co \uW_R(M\ot_R S) \to \uW_R(M)
\]
given on $T$--points by 
\[ \textstyle \sum_{i=0}^{d-1} m_i a_T^i \quad \mapsto \quad
      m_{d-1}.
\]

The quadratic form $q\co M \to R$ gives rise to another morphism of $R$--schemes, 
\[ 
  \bfq \co \uW_R(M) \to \uW_R(R) = \GG_{a,R}
\]
which assigns to $T\in \Ralg$ the base change $q_T \co M_T \to T$ of $q$. 
The composition \begin{equation}\label{aps-2}
     \uW_R(M\ot_R S) \xrightarrow{\pi_{d-1}} \uW_R(M) \xrightarrow{\;\; \bfq \; \;} \GG_{a,R}.
\end{equation} 
is a morphism of affine spaces. 

The Springer and Knebusch functors are certain subfunctors of the $R$--functor associated with the principal open subscheme  
\[ \uW_R(M\ot_R S)_{\bfq \circ \pi_{d-1}} \]
of $\uW_R(M\ot_R S)$. Its $T$--points are the $w=\sum_i m_i a_T^i\in M\ot_R S\ot_R T$, written in the form \eqref{aps-1} with $q(m_{d-1}) \in T\ti$. 

\comments{(2026-03) I have stopped here with the generalities, although one could continue, for example: \sm 

As in Poonen's 2008-paper (beginning of section 1), the set of one-generated, finite projective $R$--algebras of constant rank $d$ are represented by an $R$--scheme $\scA_d$. 

In both the Springer and Knebusch functors we are given a subscheme $U$ of $\uW_R(M\ot_R S)_{\bfq \circ \pi_{d-1}}$ and a morphism $f\co U \to \scA_{d'}$ where $d'=d-2$ in the Springer case and $d' = d-1$ for Knebusch. The condition (ii) for both functors can be interpreted as a condition on the image of $f$. I believe one can prove that these generalized Springer/Knebusch functors are represented by $R$--schemes. }

\comments{(2026-03-16) I changed the previous $\ulSpr_{\, u,v}$ to $\ulSpr_{\, (u,v)}$ to better indicate a hyperbolic pair, and I changed the previous $c$, $S'$ and $a'$ in \ref{desfi} and \ref{desfii} of \ref{desf} to $c_w$, $S_w$ and $a_w$ to indicate the dependance on $w=w(X)$. Changes are not put in blue.}

\subsection{The Springer functor $\ulSpr_{\, (u,v)}$}\label{desf} We continue with the setting of \ref{sunn} and \ref{aps}, and define the Springer functor
as the $R$--functor $\ulSpr_{\, (u, v)}(q,S)$ as the subfunctor of $\frR_{S/R}(U)$ which on $T\in \Ralg$ consists of elements    
\[ w\in \frR_{S/R}(U)(T) = U(S\ot_R T) \subset M \ot_R S \ot_R T   
\]
satisfying \ref{desfi} and \ref{desfii} below: 
\begin{enumerate}[label={\rm (\roman*)}] 
  \item \label{desfi} If $w(X)\in M \ot_R T[X]= M_T \ot_T T[X]$ with $m_i \in M_T$ is the unique lift of $w$ of degree $\le d-1$ such that $w(a_T) = w$  using \eqref{aps-1},
       then 
      \begin{equation}\label{desfi1} 
           q_{T[X]} \big( w(X)\big) = c_w \, P_T(X)\, Q_w(X) \in T[X]
      \end{equation}
       where $c_w\in T\ti$ and $Q_w(X) \in T[X]$ is a monic polynomial of degree $d-2$.  
      
  \item \label{desfii} Putting $S_w=T[X]/Q_w(X)$ and $a_w = X + (Q_w) \in S_w$,     the vector  $w(a_w)$ is an isotropic vector of the $S_w$--module $M_T \ot_T S_w$  
    (note that $S_w$ is a one-generated $T$--algebra which is free of rank $d-2$).    
\end{enumerate}

Let us add some remarks to elucidate this definition. Since $w(X)$ maps onto the isotropic $w$, the monic polynomial $P_T$ divides  $q_{T[X]}\big(w(X)\big)$. Hence,  the decomposition \eqref{desfi1} exists with $c_w\in T$ and $Q_w(X)\in T[X]$ being well-defined and unique (Euclidean algorithm). Thus
\begin{equation}\label{desf4} \begin{split} 
 \ref{desfi} \quad &\iff \quad c_w\in T\ti \text{ and } \deg Q_w = d-2
  \\ & \iff \quad q(m_{d-1}) \in T\ti
\end{split} \end{equation}
where the second equivalence follows from the representation \eqref{aps-1} of $w$. 
Assuming only \ref{desfi}, we always have $q_{S_w}(a_w)= 0$. Hence
\[  \ref{desfii} \quad \iff w(a_w) \text{ is unimodular in the $S_w$--module $M_T \ot_T S_w$. }
\]
Since unimodularity is stable under base change, we can use
\begin{equation*}\label{desf5}
  \ref{desfii} \quad \Longleftarrow \quad \text{$w(X)$ is unimodular  in the $T[X]$--module $M_T \ot_T T[X]$.}
    \end{equation*}

The subscript $(u,v)$ in $\ulSpr_{\, (u,v)}(q,S)$ should remind the reader of the dependance on the hyperbolic pair $(u,v)$ via $U$.    

\begin{lem}[Example $\ulSpr_{\, (u,v)}(q,S)(T) \ne \emptyset$] \label{sprex} 
Suppose $T\in \Ralg$ satisfies the following conditions
\begin{enumerate}[label={\rm (\roman*)}] 
 \item \label{sprexi} the quadratic space $(M,q)_T$ contains a hyperbolic pair $(e,f)$, 
 
 \item \label{sprexii} there exists $\wtl m \in (Te \oplus Tf)^\perp \subset M_T$ with $q_T(\wtl m) \in T\ti$. 
     
 \item \label{sprexiii}   $S\ot_R T$ has Witt cancellation.        
\end{enumerate}
Then $\ulSpr_{\, (u,v)}(q,S)(T) \ne \emptyset$. \end{lem}
\sm

 Before we present the proof of \ref{sprex}, we note that \ref{sprexi}--\ref{sprexiii} hold whenever \ref{sprexI}--\ref{sprexIII} below are satisfied: 
\begin{enumerate}[label={\rm (\Roman*)}] 
 \item \label{sprexI} $\rank M \ge 3$, 
 \item \label{sprexII} $(M,q)_T$ is isotropic, and 
 \item \label{sprexIII} $T$ is an LG ring. 
\end{enumerate}
Indeed, by \ref{sprexI} and \ref{sprexII},  there exists a hyperbolic pair $(e,f)\subset M_T$ with $(T e \oplus T f)^\perp$ faithfully projective. The existence of $\wtl m $ then follows from \ref{LGqdi}. The finite extension $S\ot_R T$ of the LG ring $T$ is an LG ring by \ref{revLG}\eqref{revLG-int}. It has Witt cancellation by \ref{canqf}\eqref{canqf-a}. 

\begin{proof} 
We have hyperbolic pairs $(u,v)_{S\ot T}$ and $(e,f)_{S\ot T}$, where the first is obtained by base change $S \to S\ot T$ and the second by base change $T \to T[X] \to S\ot T$ using the surjective $T$--algebra homomorphism $M_T \ot_T T[X] \to M_T \ot_T (S\ot_R T)$. 
By \ref{sprexiii} and \eqref{hps2} there exists $g\in \orth(q_{S\ot T})$ mapping 
$(S\ot T \cdot u_{S\ot T}, \, S\ot T \cdot v_{S\ot T})$ onto $(S\ot T \cdot e_{S\ot T}, S\ot T \cdot f_{S\ot T})$. By \eqref{sunres0}, $g$ maps $U (S\ot_R T)$ onto the corresponding set with respect to the hyperbolic pair $(e,f)_{S\ot T}$. A fortiori, $g$ maps $\ulSpr_{\, (u,v)}(q,S)(S\ot_R T)$ onto the set $\ulSpr_{\, (e,f)}(q_T, S\ot_R T)(S\ot_R T)$ for $(e,f)$, observing that the $T$--algebra $S\ot_R T$ satisfies the conditions of \ref{desf}. It is therefore enough to show that the latter is non-empty. We can even assume that $(u,v)_{S\ot T} = (e,f)_{S\ot T}$. 

We now use polynomials $Q(X)$ and $R(X)\in T[X]$, uniquely determined by the condition that $X^{2d-2} = P_T(X) Q(X) + R(X)$ with $\deg Q(X) = d-2$ and $\deg R(X) \le d-1$.  Since $\wtl m $ is unimodular (\ref{mx_lem}), hence 
\[ M_T = (Te \oplus Tf) \perp (T \wtl m  \oplus N) \] 
for some $T$--module $N$, we can define 
\[ w(X) = -q_T(\wtl m)R(X) e + f + X^{d-1} \wtl m \in M_T[X]. \]
Then $w(a_T) \in U(S\ot_R T)$ by \eqref{sunres0}, and $w(X)$ is the unique lift of $w(a_T)$ in the setting of \ref{desfi}. The condition \eqref{desfi1}  holds because 
\[ q_{T[X]}\big(w(X)\big) = - q_T(\wtl m) R(X) + X^{2d-2} q_T(\wtl m) = 
    q_T(\wtl m) P_T(X) Q(X). \] 
Moreover, also \ref{desfii} holds because $w(X) \in M_{T[X]}$ is unimodular and therefore so is its image in $T[X]/Q(X)$. 
\end{proof}

\begin{prop}[Representability of $\ulSpr$] \label{respr} The $R$--functor  $\ulSpr_{\, (u,v)}(q,S)$ of {\rm \ref{desf}} is representable by an $R$--scheme 
$\uSpr_{(u,v)}(q,S)$ with the following properties: 
\begin{enumerate}[label={\rm (\roman*)}] 
\item \label{respri} $\uSpr_{(u,v)}(q,S)$ is a quasi-compact open subscheme of $\frR_{S/R}(U)$, in particular it is smooth. \sm 

\item \label{resprii} The $R$--module $N = S\oplus \frR_{S/R}(M')$ is finite locally free;  $\uSpr_{(u,v)}(q,S)$ is isomorphic to a quasi-compact open subscheme of $\uW(N)$.  \sm 

\item \label{respriii} If $\rank M \ge 3$, then $\uSpr_{(u,v)}(q,S)$ has  geometrically integral 
    fibres and the embedding $\uSpr_{(u,v)}(q,S) \hookrightarrow  \frR_{S/R}(U)$ is universally schematically dense. \sm 
    
\item\label{respriiv} If $R$ is an LG ring and $\rank M \ge 3$, then $\big(\uSpr_{(u,v)}(q,S)\big)(R) \ne \emptyset$.     
\end{enumerate}
\end{prop}
\sm 

In view of Lemma~\ref{sptwo}, the restriction $\rank M \ge 3$ in \ref{respriii} and \ref{respriiv} above is not serious. 

\begin{proof} Let us abbreviate $\ulSpr_{\, (u,v)}(q,S)=\ulSpr$.  
The subfunctor $\ul{\rm Sp}$ of $\frR_{S/R}(U)$ defined by condition~\ref{desf}\ref{desfi}, equivalently by \eqref{desf4} only is represented by the principal open subscheme $\frR_{S/R}(U)_f$ of the affine scheme $\frR_{S/R}(U)$ where $f$ is the restriction of the polynomial $\bfq \circ \pi_{d-1}$ of \eqref{aps-2}. It is therefore itself an affine scheme, say $\frR_{S/R}(U)_f = \Spec(C)$ for some $C\in \Ralg$.  
\lv{
We first prove representability of the subfunctor $\ul{\rm Sp}$ of $\frR_{S/R}(U)$ defined by condition~\ref{desf}\ref{desfi}, equivalently by \eqref{desf4}. Since every $w\in \frR_{S/R}(U)(T), T\in \Ralg$, can be uniquely written in the form \eqref{desf1n}, we get a well-defined map $f_T \co \frR_{S/R}(U)(T) \to T$, $w \mapsto q_T(m_{d-1})$, which is functorial in $T\in \Ralg$ and hence defines a regular function $f \co \frR_{S/R}(U) \to \bbA_R^1$. It follows that $\ul{\rm Sp}$ is represented by the principal open subscheme $\frR_{S/R}(U)_f$ of the affine scheme $\frR_{S/R}(U)$ and is therefore itself a affine scheme, say $\frR_{S/R}(U)_f = \Spec(C)$ for some $C\in \Ralg$.  
}
\sm 

We can apply the above to the universal point $w_\sharp\in (\frR_{S/R}(U)_f)(C)$. This provides an element $c_\sharp \in C\ti$ and a monic polynomial $Q_\sharp \in C[X]$ of degree $d-2$ such that $q(w_\sharp(X)) = c_\sharp P_C(X)Q_\sharp (X)$.  We put $C_\sharp = C[X]/Q_\sharp$, 
and denote by $\uW(M_{C_\sharp})$ the $C_\sharp$--scheme associated with the $C_\sharp$--module $M_{C_\sharp} = M \ot_R C_\sharp$. Since $C_\sharp$ is a $C$--algebra, we can apply Weil restriction to $\uW(M_{C_\sharp})$ and get a $C$--scheme $\frR_{C_\sharp /  C}\big( \uW(M_{C_\sharp})\big)$. By Yoneda,
\[ 
\Mor\big( \Spec(C), \, \frR_{C_\sharp/C}\big( \uW(M_{C_\sharp})\big)  =  \big( \frR_{C_\sharp/C}\big( \uW(M_{C_\sharp})\big)\big)\, (C) = M_{C_\sharp}. \] 
Thus, we have a morphism of $C$--schemes 
\[ 
s \co \Spec(C) = \frR_{S/R}(U)_f \to \frR_{C_\sharp/C}\big( \uW(M_{C_\sharp})\big)
 \]
corresponding to $w_\sharp \big(X + (Q_\sharp)\big) \in M_{C_\sharp}$. It is a section of the structure morphism $\frR_{C_\sharp/C}\big( \uW(M_{C_\sharp})\big) \to \Spec(C)$.  
We denote by $Z$ the fiber product of $C$--schemes
\begin{equation}\label{respr00} \vcenter{
\xymatrix{
 Z  \ar@{^{(}->}[r]^{i'} \ar[d] & \Spec(C) = \frR_{S/R}(U)_f \ar[d]^s   \\ 
 \frR_{C_\sharp/C}\big( \uW(M_{C_\sharp}) \setminus \{0 \} \big)   \ar@{^{(}->}[r]^i &  \frR_{C_\sharp/C}\big( \uW_{C_\sharp}(M_{C_\sharp})  \big)
} } \quad .
\end{equation} 
Thus, by definition, $Z$ is a $C$-scheme and hence, after composing with the structure morphism $\Spec(C) \to \Spec(R)$, also an $R$-scheme. We claim that the $R$--scheme $Z$ represents $\underline{\rm Spr}$, i.e., $Z(T) = \underline{\rm Spr}(T)$ holds for all $T\in \Ralg$.%

\comments{(2026-04) Before you only proved $Z(T) = \underline{\rm Spr}(T)$ for $T=R$. The changes from $R$ to an arbitrary $T$ are minimal. Therefore I worked out the general case.}
\sm 

The functor of points of the open subset $\frR_{S/R}(U)_f \subset \frR_{S/R}(\uW_S(M_S))$ is described in \ref{desf}. 
On the other side,  $Z(T)= \Mor_R(\Spec(T), Z)$ can be identified with the set of $R$--scheme morphisms $w \co \Spec(T) \to \Spec(C)$ such that the diagram \eqref{respr00} can be completed to a commutative diagram by the dotted arrow:
\[\xymatrix{
\Spec(T) \ar@{.>}[rd]\ar@/^/[drr]^w   \\
 & Z  \ar@{^{(}->}[r]\ar[d] & \Spec(C) = \frR_{S/R}(U)_f \ar[d]^s   
 \\ & \frR_{C_\sharp/C}\big( \uW(M_{C_\sharp}) \setminus \{0 \} \big)   \ar@{^{(}->}[r] & 
 \frR_{C_\sharp/C}\big( \uW(M_{C_\sharp})  \big).
}\]
We do not need to understand the map  $s\circ w: \Spec(T) \to \frR_{C_\sharp/C}\big( \uW(M_{C_\sharp})  \big)$
but only the condition it imposes.

The morphism $w \co \Spec(T) \to \Spec(C)$ is given by a unique $R$--algebra homomorphism $w^\flat \co C \to T$, allowing us to define the $R$--algebra
\[
C_\sharp^w := C_\sharp \ot_C^{w^\flat} T = (C[X]/Q_\sharp)\ot_C^{w^\flat} T,  
\] 
where the $C$--action on $T$ on the right-hand side is given by $w^\flat$. 

Next we aim to pull-back the cartesian square of $C$-schemes \eqref{respr00} 
by the morphism $w \co \Spec(T) \to \Spec(C)$. For the lower right corner we get, using 
\eqref{weilres-ex2},
\[ 
  \frR_{C_\sharp/C}\big( \uW(M_{C_\sharp})\big) \times_C T \cong 
   \frR_{ (C_\sharp \ot_C^{w^\flat} T)/T } \big( \uW(M_{C_\sharp} \ot_C^{w^\flat}T)\big).   
\]
The morphism $w$ corresponds to a point in $M_{S\ot T}$ satisfying \ref{desf}\ref{desfi} such that, using the notation of \ref{desf}\ref{desfii}, the pullback of $s \circ w$  is given by $w(a_w)$ and $C_\sharp \ot_C ^{w^\flat} T \cong S_w$ (isomorphism of $R$--algebras). Hence, 
\begin{equation}\label{respr01a}
 \frR_{ (C_\sharp \ot_C^{w^\flat} T)/T } \big( \uW(M_{C_\sharp} \ot_C^{w^\flat}T)\big)
 \cong \frR_{S_w/T}\big( \uW(M_{S_w} )\big). 
\end{equation} 
Analogously, by \eqref{weilres-bc1} and \ref{spreq-a}\ref{spreq-aii},  the pull-back of the left lower corner of \eqref{respr00} is 
\begin{equation}\label{respr01b} 
 \frR_{C_\sharp/C}\big( \uW(M_{C_\sharp}) \setminus \{0 \} \big) \times_C T 
  \cong  \frR_{S_w/T}\big( \uW(M_{S_w}) \setminus \{0 \} \big). 
\end{equation}    
Altogether we arrive at the cartesian square \eqref{respr02} of $R$--schemes
\begin{equation}\label{respr02} \vcenter{
\xymatrix{ 
  Z_w =Z\times_C T \; \ar@{^{(}->}[r]\ar[d] & \Spec(T)  
      \ar[d]^{s_w }  \ar@/_1pc /@{.>}[l]
 \\  
  \frR_{S_w/T}\big( \uW(M_{S_w}) \setminus \{0 \} \big)   \ar@{^{(}->}[r] & 
 \frR_{S_w/T}\big( \uW(M_{S_w})  \big) 
}} \quad .  
\end{equation} 
If the morphism $w\co \Spec(T) \to \Spec(C)$ lies in $Z(T)$, the dotted arrow above exists and shows that  the open immersion $Z_w \to \Spec(T)$ is an isomorphism or, equivalently, that  $w(a_w)$ is unimodular. Therefore, we have $Z(T) \subset \underline{\rm Spr}(T)$. The other inclusion follows by retracing the steps above. We have now shown that $\underline{\rm Spr}$ is represented by the $R$--scheme $Z$ defined in \eqref{respr00}. 
\sm 

\ref{respri} We now check the quasi-compactness of the open immersion $i'$ of \eqref{respr00}, which by \cite[Tag 01TU]{St} is equivalent to $i'$ being of finite presentation. Since finite presentation allows base change \cite[Tag 01TS]{St}, it is enough to show that the open immersion $i$ is finitely presented. But since $\uW(M_{C_\sharp})\setminus \{0\}$ and  
$\uW(M_{C_\sharp})$ are finitely presented $C_\sharp$--schemes by \ref{spreq-a}, this is again a special case of \ref{weilres}\eqref{weilres-d}. Thus, we established that $i'$ is a quasi-compact open immersion and this in turn implies that $\uSpr \to \frR_{S/R}(U)$ is a quasi-compact open immersion. Since $\frR_{S/R}(U)$ is smooth by \ref{sunres}, so is $\uSpr$. By the same reason, its geometric fibres are integral or empty. 
\sm 

\ref{resprii} We have the following diagram of $R$--schemes
\[ \uSpr \hookrightarrow  \frR_{S/R}(U) \simlgr \uGL_{1,R}(S) \times \uW(\frR_{S/R}(M')) \hookrightarrow \uW(N)
\]
where the first arrow is a quasi-compact open immersion by  \ref{respri}, the second arrow is the isomorphism \eqref{sunres1}, and the last arrow is the canonical quasi-compact open immersion. Since quasi-compact morphisms allow composition \cite[Tag 01KG]{St}, we can view $\uSpr$ as a quasi-compact open subscheme of $\uW(N)$. The $R$--module $\frR_{S/R}(M')$ is finite locally free = finite projective by transitivity of ``finite projective'', see for example \cite[1.1.8]{Ford}. Hence $N$ is finite locally free too, and so $\uW(N)$ is an affine $R$--space.  
\sm 

\ref{respriii} The first part follows from \ref{respri} and the formula \eqref{sunres1}. For establishing the second part, it is in view of \ref{respri} and the criterion in \ref{todle}\eqref{todle-b}
enough to show that $\uSpr(k) \ne \emptyset$ for each algebraically closed field $k \in \Ralg$. Since $k$ is an LG ring, the conditions \ref{sprexI}--\ref{sprexIII} of \ref{sprex} are fulfilled for $T=k$, so that $\uSpr(k) \ne \emptyset$ by \ref{sprex}. \sm 

\ref{respriiv} Since by \ref{resprii} the scheme $\uSpr$ is isomorphic to a quasi-compact open subscheme of $\uW(N)$, the fundamental property \ref{prop_baire}\ref{prop_baire-a} of LG rings says that it suffices to show 
\begin{equation}\label{respriiv1}
  \uSpr (R/\gm) \ne \emptyset \text{ for every maximal ideal $\gm\ideal R$.}
\end{equation}
If $|R/\gm|< \infty$, then  $(M,q)_{R/\gm}$ is split and $R/\gm$ is an LG ring. Hence, in this case, \eqref{respriiv1} follows from \ref{sprex}. On the other side, if $|R/\gm| = \infty$, then $\uSpr(R/\gm)$ is even infinite by \ref{zdens}\eqref{adc}, using that $\uSpr\ne \emptyset$, for example by \ref{respriii}. 
\end{proof}

Finally, we can prove Springer's Odd Degree Theorem~\ref{spo}. 

\subsection{Proof of Theorem~\ref{spo}} \label{spopr} \eqref{spo-a} We first reduce  to the case $\rank M \ge 3$, which will then allow us to apply Proposition~\ref{respr}\ref{respriiv} and finish by induction on the rank of $S$. \sm 

{\em Reduction to $M$ of rank $r\ge 3$}: Let $R= R_0 \times \cdots \times R_n$ and $(M,q) = (M_0, q_0) \perp \cdots \perp (M_n,q_n)$ be the rank decomposition of $(M,q)$ as in \ref{qf}\eqref{qf-redc}. Thus, $M_i$ is a projective $R_i$--module of rank $i$ and each $q_i \co M_i \to R_i$ is a nonsingular quadratic form. The $R$--algebra $S$ decomposes correspondingly, $S=S_0 \times \cdots \times S_n$  where each $S_i$ is a finite $R_i$--algebra of degree $d=\deg S$. We have
  \[ M \ot_R S \cong (M_0 \ot_{R_0} S_0) \times \cdots \times (M_n \ot_{R_n} S_n) \]
with each $M_i \ot_{R_i} S_i$ being projective of rank $i$ as $S_i$--module. Since $q_{M\ot_R S}$ is isotropic, so is every $q_{M_i \ot_{R_i} S_i}$ by \ref{isotrop}\eqref{isotrop-b}. Hence, by \ref{isotrop}\eqref{isotrop-d}, we know $M_i \ot_{R_i} S_i = 0$ for $i=0,1$. Therefore, by faithfully flatness of the $R_i$--modules $S_i$, we get $M_0=0=M_1$. Without loss of generality we can therefore assume $R_0=0 = R_1$.

In the decomposition $R= R_2 \times \cdots \times R_n$, each $R_i$ is an LG ring, \ref{revLG}\eqref{revLG-aa}. The case $r=2$ has been dealt with in Lemma~\ref{sptwo}. We can therefore assume that $M$ has rank $r\ge 3$.
\sm 

{\em The induction argument}: By \ref{respr}\ref{respriiv}, we know that 
$\big(\uSpr_{(u,v)}(q,S)\big)(R) \ne \emptyset$. In the notation of \ref{desf} this says that $q$ is $S_w$--isotropic for $S_w=R[X]/Q_w$ and some monic  polynomial $Q_w$ of degree $d-2$. Since $\big(\uSpr_{(u,v)}(q,S_w)\big)(R) \ne \emptyset$ by another application of \ref{respr}\ref{respriiv}, we can continue by induction until the degree of $S$ is one, i.e., $(M,q)$ is isotropic.  

The interpretation of the quadratic form result in terms of the group scheme $\uSO(q)$ is a consequence of the characterization \ref{prop_isotropic}.  \sm 

\eqref{spo-b} 
We have already noted in \ref{equisa}\eqref{equisa-cn} that Springer's Odd Degree Theorem, i.e., condition \ref{equi-A} of \ref{equi}, implies the conditions  \ref{equi-Breg}, \ref{equi-D} and \ref{equi-E} of \ref{equi}. 
\qed

\subsection{Review of hermitian modules and spaces}\label{herev} We review hermitian spaces over rings following \cite[I]{K}, see also \cite[C.1]{GN-LG} (but note the change of notation in order to accommodate the notation of Corollary~\ref{sprher}).  
\sm 

Throughout $R$ is an  arbitrary base ring and $R'\in \Ralg$ is a quadratic \'etale extension whose standard involution and trace we denote by $\si$ and $\tr = \Id + \si$ respectively. 

A hermitian module over $R'/R$ is a pair $(M',h')$ consisting of a finite projective $R'$--module $M'$ and a hermitian form $h' \co M'\times M' \to R'$, which means that $h'$ is a biadditive map and satisfies 
\begin{align*} h'(m_1 r_1',\, m_2 r_2') & 
        = \si(r')\, h'(m_1, m_2) \, r_2' \quad \text{and} \\
      h'(m_1, m_2) & = \si \big( h'(m_2, m_1) \big) 
  \end{align*}   
for all $m_1, m_2 \in M'$ and $r_1', r_2' \in R'$. An isotropic vector of $(M',h')$ is a unimodular $m\in M'$ satisfying $h'(m, m) = 0$. 

We twist the canonical $R'$--action on the dual module $\Hom_{R'}(M', R')$ by putting $\vphi \star r' = \vphi \, \si(r')$ for $\vphi \in \Hom_{R'}(M', R')$ and denote the resulting $R'$--module by $M^{\prime \, *}$. Any hermitian module $(M',h')$ over $R'$ gives rise to an $R'$--linear map $\wdh{h'} \co M' \to M^{\prime\, *}$, $m' \mapsto h'(m' , \cdot)$. We say that $h'$ regular if $\wdh{h'}$ is a bijection, in which case we call $(M',h')$ a hermitian space. 

A hermitian module $(M',h')$ over $R'/R$ gives rise to the $R$--quadratic form  
\[ \tr_{R'/R}(h') \co \frR_{R'/R}(M') \to R, \quad m' \mapsto h'(m',m'), \]
called the trace form of $(M',h')$. The assignment $h' \mapsto \tr_{R'/R}(h')$ naturally extends to a trace functor $\tr_{R'/R}\co \hm_{R'/R} \to \qm_R$ from the category $\hm_{R'/R}$ of hermitian modules over $R'/R$ to the category $\qm_R$ of quadratic forms over $R$. 
It preserves orthogonal sums and sends hermitian spaces to quadratic spaces. The latter claim is a special case of \cite[I, (7.2.4)]{K} since by \cite[I, (7.3.4)]{K} the trace $\tr \co R' \to R$ is an involution trace. 

Let $S\in \Ralg$. Then $S'=R'\ot_R S$ is a quadratic \'etale $S$--algebra. To any hermitian module $(M',h')$ over $R'$ we can associate the hermitian module 
$(M',h')_{S'} = (M'\ot_{R'} S', h'_{S'})$ over $S'$ where $M'\ot_{R'}S' = \frR_{R'/R}(M') \ot_R S$ and $h'_{S'}$ is given by 
$h'_{S'}(m_1' \ot s_1, \, m_2'\ot s_2) = h'(m_1', m_2') \ot s_1s_2$  for $m'_1, m'_2 \in M'$ and $s_1, s_2 \in S$.  If $h'$ is regular, then $h'_{S'}$ is regular too. Moreover, the trace functor commutes with base change:
\begin{equation} \label{herev1}
    \tr_{S'/S}\big((M',h')_{S'}\big) = \big(\tr_{R'/R}(M',h')\big)_S  
 \end{equation}
where on the left-hand side we use the base change of hermitian forms, see for example \cite[C.1(f)]{GN-LG}, and on the right-hand that of quadratic forms, \ref{qf}\eqref{qf-bc}.

\begin{lem}\label{herle} Let $R'/R$ be a quadratic \'etale extension, and let $(M',h')$ be a hermitian module over $R'/R$. 
\sm 

\begin{inparaenum}[\rm (a)] \item\label{herle-a}  {\em (Unimodularity, isotropy)} Let $m'\in M'$. Then $m'$ is unimodular in $M'$ if and only if $m'$ is unimodular in the $R$--module $\frR_{R'/R}(M')$. In particular, $m'$ is an isotropic vector of $(M',h')$ if and only if $m'$ is an isotropic vector of the quadratic form $\tr(M',h')$. \sm 

\item \label{herle-b} {\em (Orthogonal complement)} Let $U'\subset M'$ such that $h'|_{U'\times U'}$ is regular. Then $M' = U'\perp U^{\prime\, \perp}$ and the orthogonal complement $U^{\prime\, \perp}$ of $U'$ is 
    \begin{equation}\label{herle-b1}\begin{split}
     U^{\prime\, \perp} &= \{m'\in M' : h'(m', U') = 0 \} 
     \\ & = \{ m\in \frR_{R'/R}(M') : b_{\,\tr (h')} (m, \, \frR_{R'/R}(U')) = 0 \} 
\end{split}\end{equation}
where $b_{\, \tr(h')}$ is the polar form of the trace form $\tr(h')$. \sm 

\item\label{herle-c} Assume that $(M',h')$ is diagonalizable and that $(M_1, h'_1)$ is another hermitian space over $R'/R$ such that the quadratic modules $\tr(M', h')$ and $\tr(M_1', h_1')$ are isometric. Then  $(M',h') \cong (M_1, h'_1)$ as hermitian spaces. 
\end{inparaenum}\end{lem} 
 
\begin{proof} \eqref{herle-a} It suffices to prove the claim regarding unimodularity. First, suppose that $m'$ is a unimodular vector of $M'$. Thus, by \ref{unimod}\ref{unimodi}, $M' = R' m' \oplus M''$ for some $R'$--submodule $M''$ of $M'$ and the canonical map $R' \simlgr R'm'$ is an isomorphism of $R$--modules. We know, e.g., by \cite[2.4.6]{Ford}, that $1_{R'}\in R'$ is unimodular, so $R' = R 1_{R'} \oplus R''$ for some $R$--line bundle $R''$ and $R\simlgr R1_{R'}$. Writing $M' = R m' \oplus \big(R'' \oplus \frR_{R'/R}(M'')\big)$ shows that $m'$ is unimodular in $\frR_{R'/R}(M')$.  

Conversely, assume that $m'\in \frR_{R'/R}(M')$ is unimodular, and let $\gm' \ideal R'$ be a maximal ideal of $R'$. We know, e.g., by \cite[2.4.2]{Ford}, that $\gm = \gm' \cap R$ is a maximal ideal of $R$. Hence $R'/\gm'$ is an extension field of the field $R/\gm$. As $m'\ot 1_{R/\gm} \ne 0$ by \ref{unimod}\ref{unimodiv}, also $m'\ot 1_{R'/\gm'} \ne 0$. Since this holds for every maximal ideal of $R'$, the vector $m'$ is unimodular in $M'$, again by \ref{unimod}\ref{unimodiv}.  

\eqref{herle-b} That $M' = U' \perp U^{\prime\, \perp}$ is standard, see for example \cite[I, (3.6.2)]{K}, where the assumption that $U$ be finite projective is not necessary. Since $\tr(h')$ is regular, we also have the decomposition $\frR_{R'/R}(M') = \frR_{R'/R}(U') \perp V$ with $V= \{ m\in \frR_{R'/R}(M') : b_{\,\tr (h')} (m, \, \frR_{R'/R}(U')) = 0 \}$, see \ref{orthLG}\eqref{orthLG-a}. It then suffices to observe $U^{\prime \, \perp} \subset V$ because $b_{\tr(h')} = (1+ \si)\circ h'$.  \sm 

\eqref{herle-c} The hermitian module $(M',h')$ is diagonalizable if and only if there exists a decomposition $M' = U'_1 \perp \cdots \perp U'_n$ where each $U'_i$ is free of rank $1$ with $R'$--basis $\{u'_i\}$ satisfying $h'(u'_i, u'_i) \in R\ti$. If $g \co \tr(M',h') \simlgr \tr(M'_1, h'_1)$ is an isometry, then $\frR_{R'/R}(M_1') = g(U'_1) \perp \cdots \perp g(U'_n)$ by \eqref{herle-b}, and each  $g(U'_i)$ contains $g(u'_i)$ satisfying  $h'_1\big( g(u'_i), \, g(u'_i)\big) = h'(u'_i, u'_i)\in R\ti$. The claim then easily follows. \end{proof}

\subsection{Proof of Corollary~\ref{sprher}}\label{herpr} \eqref{sprher-a} Let us first note that $\tr_{R'/R}(M',h')$ is a regular (= nonsingular) quadratic form over $R$. Since by assumption $(M',h')_{R'\ot S}$ is isotropic, Lemma~\ref{herle}\eqref{herle-a} says that $\tr_{S'/S}((M',h')_{S'})$ is an isotro\-pic $S$--quadratic form, i.e., by \eqref{herev1}, the $S$-quadratic form $(\tr_{R'/R}(M',h'))_S$ is isotropic. By \ref{spo}, the quadratic form $\tr_{R'/R}(M',h')$ is isotropic. Then another application of Lemma~\ref{herle} finishes the proof. 
\sm 

\eqref{sprher-b} We can assume that $M'$ has constant rank. If the hermitian spaces $(M',h')_{R'\ot S}$ and $(M'_1, h'_1)_{R'\ot S}$ are isometric, then the $S$--quadratic forms $\big(\tr_{R'/R}(M',h')\big)_S$ and $\big(\tr_{R'/R}(M_1',h_1')\big)_S$ are isometric by \eqref{herev1}. But then, by \ref{equi-D} of Theorem~\ref{spo}, the $R$--quadratic forms $\tr(h')$ and $\tr(h'_1)$ are isometric. As $\tr(h')$ is diagonalizable by \cite[C.4]{GN-LG}, Lemma~\eqref{herle-c} applies and yields $(M',h') \cong (M'_1, h_1')$. 
\qed

\lv{

\section{Abfall} 
\comments{(2024-01) Lemma~\ref{aqyom} is the analogue of Lemma~\ref{lem_sphere_field}\ref{lem_sphere_field2} for fields and of Lemma~\ref{lem_sphere} for arbitrary base rings.}

\quest{(2024-12-12) Lemma~\ref{aqyom} does not seem to be correct: $\Psi$ is not injective as $R$--functor since $m$ and $r m$ for $r\in R\ti$ have the same image under $\Psi$. \ms

(2024-01) It is not clear to me why we take $-\si_m(y)$ and not $\si_m(y)$. We did not use the lemma}

\begin{lem}\label{aqyom} 
Let $(M,q)$ be a quadratic $R$--space and let $y\in \uY_q(R)$, so that $(M,q)$ contains  a hyperbolic plan by \eqref{aqyle-b1}. Then the $R$--functor defined by 
\[ \{ m \in M_S: q(m) \in S\ti\} \to \uY_q(S), \quad m \mapsto - \si_m(y) \]  
for $S\in \Ralg$ induces an open immersion 
\begin{equation}  \label{aqyom1}
\Psi \co  \uW(M) \setminus \big(\uZ_q \cup \Ker b_q(y, \cdot)\big) \; \hookrightarrow \; \uY_q\, .
 \end{equation}
\end{lem}

\begin{proof} Since the reflection $\si_m$ is an orthogonal transformation and since reflections are stable under base change, we get a well-defined morphism of schemes
$\wtl \Psi \co \uW(M) \setminus \uZ_q \to \uY_q$, whose restriction to the open subscheme 
\[ 
\uU := \uW(M) \setminus \big(\uZ_q \cup \Ker b_q(y, \cdot)\big)
\] 
of $\uW(M) \setminus \uZ_q$ is the map $\Psi$ of \eqref{aqyom1}. The purpose of restricting $\wtl \Psi$ is to make $\Psi$ a monomorphism. \sm

We will first prove that $\Psi$ is an open immersion in case $R=k$ is a field. To do so, we consider the tangent map 
$T(\Psi)_x \co T_x( \uU) \to T_{\Psi(x)} (\uY_q)$ 
of $\Psi$ at $x\in \uU(R)$, namely
\[ 
T(\Psi)_x \co M \to \Ker\big( b_q(\Psi(x), \cdot ) \big) 
\] 
because $T_x(\uU) = M$ and $T_{\Psi(x)}(\uY_q) = T_{\Psi(x)} (\uZ_q) = \Ker b_q(\Psi(x), \cdot)$.  Observe that $T_{\Psi(x)}(\uY_q)$ is $(r-1)$-dimensional where $r= \dim_k M$.  Computing with the dual numbers $k[\veps]$, we have $q(x+\veps h)\me = q(x)\me - \frac{b_q(x,h)}{q(x)^2}\veps$ and therefore  
\begin{align*} 
&\Psi(x+ \veps h)- \Psi(x) =
 \Bigl(  \frac{b_q(y,x)}{q(x)} h   +  \frac{b_q(y,h)}{q(x)}x  -
 \frac{b_q(y,x)b_q(x,h)}{q(x)^2}x \Bigr) \veps 
\\ &= q(x)^{-2}  \, \,   \Bigl(  q(x) \,  b_q(y,x)h +  \big(  q(x) b_q(y,h)- b_q(y,x) b_q(x,h)  \big)x \Bigr) \veps 
\\ &= q(x)^{-1}  \, \Bigl( b_q(y,x)h -  
    \big( b_q (\Psi(x), h)\big)x \Bigr) \veps.
\end{align*}
The formula above shows that $kx \subset \Ker (T(\Psi)_x)$ with equality whenever $q(x) \ne 0 \ne b_q(\Psi(x), \cdot )$. In the latter case, $T(\Psi)_x$ is surjective. Observe that this condition is stable when passing to an algebraic closure of $k$. Since $\uU$ is smooth at $x$, we can apply \cite[IV$_4$, 17.11.1(d)]{EGA} and conclude from \cite[IV$_4$, 17.11.1(a)]{EGA} that $\Psi$ is smooth at $x$. It then follows that $\Psi$ is a smooth monomorphism, hence an open immersion in view of \cite[IV$_4$, 17.9.1]{EGA}. \sm 

We can now prove \eqref{aqyom1} for an arbitrary $R$. As an open subscheme of $\uW(M)$, the scheme $\uU$ is flat and of finite presentation. Hence, the field case together with the Fibrewise Criterion \ref{ag}\eqref{ag-d} shows that $\Psi$ is an open immersion. \end{proof}

\quest{What is the right generalization of Springer's Theorem? The main part of Springer's proof deals with simple extensions (and not with arbitrary odd degree extension). Hence, it may be of interest to have an account of Springer's proof of Springer's Theorem in that case. This can be deleted later.}

\comments{Let $k$ be a field, so that $R=k[X,Y]$ is factorial. Then $M= (X,Y) \ideal R$ is a finitely generated torsion-free, hence finite projective $R$--module by the Structure Theorem of finite type modules over Dedekind domains. But $M$ is not free. This justifies the assumption ``$M$ free'' in Lemma~\ref{spriaux} below. }

\begin{lem}[Lifting unimodular vectors]  \label{spriaux} Let $R$ be a UFD, let $M$ be a free $R$--module of finite rank, and let $P\in R[X]$ be monic. Put $S= R[X]/(P)$. Then any unimodular $v_S \in M_S$ lifts to a unimodular $v(X) \in M_{R[X]}$. 

Moreover, assume that $(M,q)$ is a quadratic module and that $v_S$ is isotropic. Then there exists a polynomial $\wtl Q \in R[X]$ satisfying
  \[ q_{R[X]}\big(v(X)\big) = P(X) \wtl Q(X), \quad \deg\wtl Q \le \deg P - 2.\]
\end{lem}

\begin{proof} We have $S= \bigoplus_{i=0}^{d-1} Ra^i$ with $a= X + (P)$ and $d = \deg P$. Let $e_1, \ldots, e_n$ be a basis of $M$ and let $v_S\in M_S$ be unimodular. Write $v_S = \sum_i s_i e_i \in M_S$ with $s_i \in S$ and then write $s_i = \sum_{j=0}^{d-1} r_{ij} a^j$ with unique $r_{ij} \in R$. Define $g_i(X) = \sum_{j=0}^{d-1} r_{ij} X^j \in R[X]$ and $v(X) = \sum_{i=1}^n g_i(X) e_i$. Then clearly $v(X)$ lifts $v_S$. 

Since $v_S$ is unimodular, the ideal $(s_1, \ldots , s_n) = S$. 
The ideal $I = \big( g_1(X), \ldots, g_n(X)\big) \ideal R[X]$ is generated by some $d\in R[X]$ because $R[X]$ is a UFD. We define $g'_i \in R[X]$ by $g_i(X) = g'_i(X) d(X)$. Note
\[ I|_{X=a} = (s_1, \ldots, s_n) = S = d(a)\big( g_1'(a), \ldots g'_n(a)\big)\]
which implies $d(a) \in R\ti$ and $(g'_1(a), \ldots, g'_n(a)\big) = S$. 

Put $v'(X) = \sum_{i=1}^n g'_i(X) e_i$. If $\deg d(X) \ge 1$, we can lower $\sum_{i=1}^n \deg g_i(X)$ by replacing $g_i(X)$ by $g'_i(X)$; since then $v_S = d(a) v'_S$, it is enough to lift $v'_S$. Hence, without loss of generality,  we can assume $\deg d(X) = 0$, i.e., $d(X) \in R\ti$, so that $d(X) = d(a) \in R\ti$. We then have $I = \big( g_1(X), \ldots, g_n(X) \big) = R[X]$ because the gcd of the $g_i$ is a unit, proving that $v(X)$ is unimodular.


Suppose now that $v_S$ is isotropic. Since $P$ is monic, we can apply division with remainder and get polynomials $\wtl Q(X)$ and $U(X)$ such that 
\[ q_{R[X]}\big( v(S) \big) = P(X) \wtl Q(X) + U(X), \quad \deg U(X) < \deg P(X). 
\] 
Substituting $a$ for $X$, leads to $0 = q_S(v_S) = U(a)$. But $U(X)$ can be written in the form $U(X) = \sum_{i=0}^{d-1} r_i X^i$. Hence $0 = \sum_{i=1}^{d-1} r_i a^i$, so that all $r_i = 0$. Equivalently, $U(X) = 0$. Since $\deg g_i < d$, we can write $v(X)$ in the form $v(X) = \sum_{i=0}^{d-1} m_i X^i$. Then $q_{R[X]}\big( v(X)\big)$ has degree $\le 2(d-1)$. Since $q_{R[X]}\big(v(X)\big) = P(X) \wtl Q(X)$ has degree $d + \deg \wtl Q$, it follows that $\deg \wtl Q \le d-2$. \end{proof}
\sm

In order to use Lemma~\ref{spriaux}, i.e., replace $S=R[X]/(P)$ by $T=R[X]/(\wtl Q)$ we need to know that $T$ is a free $R$--module. This is guaranteed if $\wtl Q(X) = cQ(X)$ for some $c\in R\ti$ and $Q(X) \in R[X]$ monic. Lemma~\ref{spror} gives a sufficient (but not very interesting) condition for this. 

\begin{lem}[Springer's original proof]\label{sprior} Assume
\begin{enumerate}[label={\rm (\roman*)}]
  \item $R$ is a UFD, 
  
  \item $M$ is a free $R$--module, 
  
  \item \label{sprior-iii} $(M,q)$ is a quadratic module with an {\em anisotropic} quadratic form $q\co M \to R$ in the sense that $0 \ne m \in M \implies q(m) \in R\ti$.
\end{enumerate}
Then for any one-generated free (= projective) $S\in \Ralg$ of {\em odd degree\/}, the quadratic module $(M,q)_S$ is not isotropic.
 \end{lem}

\begin{proof}
Assume to the contrary that $(M,q)_S$ is isotropic for some $S= R[X]/(P)$ f odd degree $d$. It suffices to show that thee exists $T=R[X]/(Q)$ with $Q\in R[X]$ monic of degree $< d$ such that $(M,q)_T$ is isotropic. Then we can either proceed by induction on $d$ or get a contradiction by choosing $P$ of minimal degree such that $(M,q)_S$ is isotropic. 

We apply Lemma~\ref{spriaux} to an isotropic $v_S \in M_S$, get an unimodular left $v(X)\in M_{R[X]}$ and an factorization $q_{R[X]}\big(v(X)\big) = P(X) \wtl Q(X)$ with $\deg \wtl Q \le \deg P - 2$. We write $v(X)$ in the form $v(X) = \sum_{i=0}^k v_i X^i$ with $v_i \in M$, $v_k \ne 0$ and $k \le d-1$. Thus $q(v_k) \in R\ti$ by anisotropy of $q$. Let $c$ be the leading coefficient of $\wtl Q$ and let $e = \deg \wtl Q\le d-2$. We then have
\begin{align*}
  q_{R[X]}\big( v(X)\big) &= q(v_k) X^{2k} + (\text{lower terms})
    \\ & = cX^{d+e} + (\text{lower terms}).
\end{align*}
Hence by comparison of the leading terms, $c = q(v_k) \in R\ti$ and $2k = d+e$, which proves $\wtl Q(X) = c Q(X)$ with $Q(X) \in R[X]$ monic of odd degree $e= 2k-d\le 2(d- 1) - d = d-2$. \end{proof}
}
\newpage

\section{Scharlau's Norm Principle}\label{sec:scharlau}

The theme of this section is similarity of quadratic forms and Scharlau's norm principle, proven in Theorem~\ref{thm_scharlau} for finite \'etale extensions and nonsingular quadratic forms.  \sm

We start with some essentially known facts on the Scharlau transfer of quadratic forms.

\subsection{Frobenius extensions}\label{frobex} We call $(S, t)$ a {\em Frobenius extension\/} if
\begin{enumerate}[label=(FE\roman*)]   \item \label{frobex-i}
$ S$ is an $R$--ring, whose underlying $R$--module is finite projective, and
 \item \label{frobex-ii}
$t\in S^* = \Hom_R(S, R)$ has the property that the symmetric bilinear form
 $ \bar t \co S \times S \to R, (s_1, s_2) \mapsto t(s_1 s_2) $ 
 is regular.
\end{enumerate}
Some elementary properties:
\sm

 \begin{inparaenum}[(a)]
\item\label{frobex-c} Suppose $S\in \Ralg$ satisfies \ref{frobex-i} and let $t\in S^*$. By the standard characterization of regular bilinear forms, the following are equivalent. \end{inparaenum}

      \begin{enumerate}[label=(\roman*)]
         \item $(S,t)$ is a Frobenius $R$--extension,

         \item $(S_\m, t_\m)$ is a Frobenius $R_\m$--extension for all maximal ideals $\m\ideal R$,

        \item $(S/\m S, t_{S/\m S})$ is a Frobenius $R/\m$--extension for all maximal ideals $\m\ideal R$.
     \end{enumerate}

 \begin{inparaenum}[(a)]\setcounter{enumi}{1}
 \item \label{frobex-d}
(\cite[I, (2.10)]{Ba}) Let $(S, t)$ be a Frobenius extension of $R$. Recall that for an $S$--module $X$ we denote by $_R X$ the $R$--module obtained from $X$ by restricting scalars to $R$. Then for any $S$--module $M$ the map
    \begin{equation}\label{frobex-d1}
     t_* \co {_R (\Hom_S(M,S))} \simlgr \Hom_R( {_R M}, R), \quad f \mapsto
     t\circ f,
     \end{equation}
 is an isomorphism of $R$--modules.

Indeed, let $g\in \Hom_R({_R M}, R)$. For fixed $m\in M$, the map $s \mapsto g(sm)$ is a linear $R$--form on $S$. Hence, there exists a unique $g'(m) \in S$ satisfying $g(sm) = \bar t( g'(m), s)$ for all $s\in S$. It is immediate that $g'\in \Hom_S(M,S)$. We have $t_*(g') = g$ because $g(m) = \bar t( g'(m), 1_S) = (t \circ g')(m)$. Thus $t_*$ is surjective. It is also injective: if $t\circ f = 0$, then for all $m\in M$, $s\in S$, we get $0 = (t \circ f)(sm) = t\big( s f(m)\big)$, hence $f(m) = 0$ and then $f=0$.
\end{inparaenum}
\ms

\textbf{Examples.} \begin{inparaenum}[(I)] \item Let $S\in \Ralg$ be finite projective and let $t\in S^*$. The adjoint of the bilinear form $\tau:= \bar t$ has kernel $\Ker(\wdh \tau) = \{ s\in S: t(sS) = 0 \}$, which is an ideal contained in $\Ker (t)$. We have $\Ker(\wdh\tau) = \Ker(t)$ as soon as $\Ker(t)$ is an ideal of $S$. Thus, in general not every $t\in S^*$ makes $(S,t)$ a Frobenius extension. 

If $R=F$ is a field, then $(S,t)$ is a Frobenius extension if and only if $\Ker(\wdh \tau) = 0$. In particular, if $K/F$ is a finite field extension, then $(K, t)$ is a Frobenius extension if and only if $t\ne 0$.
\lv{
Indeed, it suffices to show that $(K,t)$ is a Frobenius extension whenever $t\ne 0$.
To see this, we show that the adjoint of the bilinear form $\bar t$ is injective. Assume, for contradiction, that $t(xK) = 0$ for some $0 \ne x$. Since $t\ne 0$, there exists $y\in K$ such that $t(y) \ne 0$. Then $0 = t(x(x\me y)) = t(y) \ne 0$, contradiction.
}
\sm

\item \label{frobex-I} Let $E$ be a finite \'etale $R$--algebra and let $\tr\co E \to R$ be its trace map. Then $(E, \tr)$ is a Frobenius extension. Indeed, if $R$ is a field, this holds by \cite[V, \S8.2, prop.~1]{BA5}. The general case then follows from \eqref{frobex-c}. \sm

\item \label{frobex-even} Let $S=R[X]/(p)$ where $p$ is a monic polynomial of degree $d\ge 1$. Let $u\in R\ti$ and define the linear form $t_u$ by
    \[ t_u \co S \to R, \quad t_u(x^i) = \begin{cases} 0, & 0\le i \le d-1,
                                         \\    u, & i =d-1 .
    \end{cases}
    \]
Then $(S, t_u)$ is a Frobenius extension of $R$.

\inparcom{(2022-12) For $d=2$ this extension is used in Lam's book \cite[VII, 3.4]{Lam-qf} to calculate the Scharlau transfer $t_{u, *}(\lan x \ran_b)$ of any $x\in S\ti$, which in turn is used in \cite[VII, 4.5]{Lam-qf} to prove an ``iff'' version of Scharlau's norm principle: let $K/F$ be a quadratic field extension (recall characteristic $\ne 2$), let $q$ be any $F$--form,  and let $x\in K$. Then $\rmN_{K/F}(x) \in \rmG(q) \iff x\cdot q \cong q'_K$ for some $F$--form $q'$. Scharlau's norm principle is $\Leftarrow$ by taking $q'=q$.}

\item See Lemma~\ref{lem_transfer}\eqref{lem_transfer-a} for another example of a Frobenius extension.
\end{inparaenum}

\comments{(2021-01-21) Not needed:\\

\begin{inparaenum}[(a)]
  \item 
  ({\em Base change}) Frobenius extensions are stable under arbitrary base change: if $(S,t)$ is a Frobenius extension of $R$, then for any $U \in \Ralg$ the pair $(S \ot_R U, t\ot \Id_U)$ is a Frobenius extension of $U$.
\sm

\end{inparaenum} }

\comments{(2021-03-21) In \cite{OP} the authors call $t$ of a Frobenius extension $(S,t)$ an {\em Euler trace\/} since Euler has considered this in a special case (this seems to be the case considered in Ch.III, \S6, Lemma 2 of Serre's Corps Locaux). They prove several properties, like base change, preservation of orthogonal sums, hyperbolic. They refer to \cite[I, \S7]{K}.

They use a different trace map in \ref{lem_transfer}, defined by projecting onto the highest power and show that for a one-generated extension of odd degree the composition of the canonical maps
\[ \W(R) \longto \W(S) \longto \W(R)\]
is the identity, i.e., $\W(R) \longto \W(S)$ is a section of the canonical maps on Witt groups induced by the Scharlau trace.}

\subsection{Scharlau transfer} \label{Schtra} Let $(S,t)$ be a Frobenius extension.
Given a bilinear module $(M,b)$ over $S$, we denote by
 $t_*(M,b)= ({_RM},B)$ the bilinear module over $R$  defined by
 $B(m_1,m_2)=t\bigl( b(m_1,m_2) \bigr)$. This is well-defined since ${_R M}$ is a
  finite projective $R$--module (\cite[1.1.8]{Ford}).
  Similarly, the Scharlau transfer of a quadratic module $(M,q)$ over $S$
  with respect to $t$ is $t_*(M,q)= ({_R M},Q)$ where $Q(m)= t(q(m))$ for each $m \in M$. Its  associated bilinear form  is $t_*(M,b_q)$. We will use the following facts. \sm

\begin{inparaenum}[(a)]
\item\label{Schtra-z} The Scharlau transfer preserves isometry and orthogonal sums. 
\sm

\item \label{Schtra-b} The transfer respects metabolic spaces,  $t_*\big(\MM(U,b)\big)\allowbreak  = \MM({_R U}, t_*(b)\big)$, and hyperbolic quadratic forms: $t_*\big(\HH(N)\big) = \HH({_R N})$. 
    \sm

Proof of \eqref{Schtra-b} for metabolic spaces:  We put $\frR_{S/R}(\cdot) = {_R (\cdot)}$ and then have $_R(U \oplus U^*)= {_RU} \oplus {_R(} U^*)$ and, by \ref{bfLG}\eqref{meta-1},
\begin{align*} &t_*\big( b_{\MM(U)} (u + \vphi, v + \psi)\big) = t_*(b)(u,v) + (t \circ  \vphi)(v) + (t\circ \psi)(u).
\end{align*}
Since $_R(U^*) \simlgr ({_R U})^*$, $\vphi \mapsto t\circ \vphi$ is an isomorphism by \eqref{frobex-d1}, and the claim follows. The claim for hyperbolic quadratic forms can be proven in the same way. \sm

\item \label{Schtra-a} {\rm (\cite[I, (2.9)]{Ba})} Let $(M,b)$ be a bilinear module over $S$. Then {\em $b$ is regular if and only if $t_*(b)$ is regular}. Similarly, if $(M,q)$ is a quadratic module over $S$, then {\em $q$ is regular if and only if $t_*(q)$ is regular}.

Proof of \eqref{Schtra-a}: Let $\wdh b \co M \to M^*$, $m \mapsto b(m, \cdot)$ be the adjoint of $b$. It is $S$--linear. The claim then follows from the fact that $t_* \circ (_R{\wdh b})$ is the adjoint of $t_*(b)$, where $t_*$ is the $R$--linear isomorphism of \eqref{frobex-d1}. The claim regarding quadratic forms can be shown in the same way. \sm

\item\label{Schatra-w} ({\em Witt groups\/})
In view of \eqref{Schtra-z}, \eqref{Schtra-b} and \eqref{Schtra-a}, the transfer map induces group homomorphisms
\begin{equation}
  \label{Schatra-w1} t_* \co  \hWq(S) \to  \hWq(R) \quad \text{and} \quad
 t_*\co \Wq(S) \to \Wq(R)
\end{equation}
of the Witt-Grothendieck groups and Witt groups, defined in \ref{wgg}.
\sm

\item \label{prop_scharlau_1} ({\em Frobenius reciprocity)}  Let $(M, b)$ be a bilinear module over $S$ and let $(N, q)$ be a quadratic module over $R$.   Then the Scharlau transfer satisfies
    \begin{equation} \label{Schtra-fr1}
         t_*\bigl( ( M, b) \otimes_S (N, q)_S    \bigr) \cong t_*( M, b) \otimes_R (N, q),
    \end{equation}
    see \ref{qfba-tens} for the tensor product between a symmetric and a quadratic form. Frobenius reciprocity is stated in \cite[I, (2.12)]{Ba} for regular bilinear and quadratic modules, but the proof there works in general.
    It is based on the natural isomorphism $_R\big(M \ot_S (S \ot_R N)\big) \simlgr (_RM) \ot_R N$ of $R$--modules.  \end{inparaenum}

\ms

We remind the reader of the notation \eqref{qfba-one1} regarding bilinear forms on $R$: for $x\in R$, we let $\lan x \ran_b  \co R \times R \to R$ be  the bilinear form $(r_1, r_2) \mapsto x r_1 r_2$. We will also use the $2$--dimensional bilinear form
\begin{equation}\label{Schatrano}
\lan x, y\ran_b = \lan x \ran_b \perp \lan y \ran_b\, ,
\end{equation}
which is regular if and only if both $x,y\in R\ti$.

\comments{(09-14) Lemma~\ref{lem_transfer} is also stated in \cite[Lem.~6.2]{Ho2}.  }
\begin{lem}\label{lem_transfer} Let $p= a_0 + a_1X + \cdots + a_{d-1}X^{d-1} + X^d$ be a monic polynomial over $R$ of degree $d\ge 2$. Put $S= R[X]/(p)$ and $x = X + (p) \in S$, so that $S=R[x]$. We denote by $\rmN_{S/R}$ the norm of the $R$--algebra $S$, {\rm \ref{trno}}, and let $t \co S \to R$ be the linear form defined by
\begin{equation}\label{lem_transfer0}
t(1_S) = 1_R \quad\text{and}\quad t(x^i) = 0, i=1, \ldots, d-1.
\end{equation}

\begin{inparaenum}[\rm (a)]

  \item \label{lem_transfer-a} Then $(S, t)$ is a Frobenius extension if and only if $a_0\in R\ti$  if and only if\/  $\rmN_{S/R}(x) \in R\ti$ if and only if $x\in S\ti$.
       \sm

  \item\label{lem_transfer-b} If\/ $a_0 \in R\ti$, the transfer of the $S$--bilinear forms $\lan 1_S\ran_b$ and $\lan x \ran_b$ satisfy
      \begin{align}  \label{lem_transfer-1}
       t_*(\lan  1_S\ran_b)  & \cong
       \begin{cases}
        \lan 1_R \ran_b \perp \MM(U_1), & \text{if $d$ is odd} \\
             \lan 1_R,  - \rmN_{S/R}(x)\ran_b \perp \MM(U_2), & \text{if $d$ is even} \\
          \end{cases}
    \\ \label{lem_transfer-2}
     t_*(\lan x \ran_b) &\cong \begin{cases}
     \lan \rmN_{S/R}(x)\ran_b  \perp \MM(U_3), & \text{if $d$ is odd} \\
     \MM(U_4), & \text{if $d$ is even,}
     \end{cases}
   \end{align}
   where  in all cases the metabolic space $\MM(U_i)$ is an orthogonal sum of free metabolic spaces $\MM(U_{ij})$ of rank $2$. \sm

\item \label{schatraco-a} Let $(M,q)$ be a regular faithful quadratic module over $R$. Then, using the notation of \eqref{lem_transfer-b} and the tensor product from {\rm \ref{qfba-tens}},
\begin{align*} t_*(\lan 1_S\ran_b) \ot q
   &\cong \begin{cases} q \perp \HH(U_1\ot_R M), & \text{if $d$ is odd,} \\
                   (\lan 1, - \rmN_{S/R}(x)\ran_b \ot_R q ) \perp \HH(U_2\ot_R M)  & \text{if $d$ is even},     \end{cases}
 \\
   t_*(\lan x \ran_b) \ot_R q &\cong
    \begin{cases} \lan \rmN_{S/R}(x) \ran_b \ot_R q \perp \HH(U_3\ot_R M),
                     & \text{if $d$ is odd,} \\
                   \HH(U_4 \ot M) &\text{if $d$ is even}.     \end{cases}
\end{align*}
\end{inparaenum}
\end{lem}


\begin{proof} \eqref{lem_transfer-a} is proven in \cite[V, \S3]{Ba}. The matrix representing $\bar t$ in the basis $1, x, \ldots, x^{d-1}$ of the $R$--module $S$ is
\[
  \big( {\bar t}(x^i, x^j)\big) =
    \begin{pmatrix}
      1  &  0  & \cdots  & 0 & 0  \\
      0  &  0  & \cdots  & 0 & -a_0 \\
      \vdots  &  \vdots  & \iddots  & -a_0 & 0 \\
      0  &  0  & \iddots  & 0 & 0 \\
      0  & -a_0  & \cdots  & 0 & 0
    \end{pmatrix}
\]
which has determinant $(-a_0)^{d-1}$. Regularity of $\bar t$ is equivalent to the matrix being invertible, \ref{bfLG}\eqref{bfLG-rech}. The second equivalence follows from \ref{trno}\eqref{trno-bou}. \sm

\eqref{lem_transfer-b} Over fields, the formulas \eqref{lem_transfer-1} and \eqref{lem_transfer-2} are established in \cite[20.9 and 20.12]{EKM}. The proofs over rings are the same, see \cite[V, (3.3)]{Ba} for \eqref{lem_transfer-1}. Let us sketch the proof of \eqref{lem_transfer-2}, setting $b= t_*(\lan x \ran)$.

Assume $d=2e+1$ is odd and put $K=\Ker(\tr) = Rx \oplus \cdots \oplus Rx^{d-1}$. Since $b(x^i, x^j) = t(x^{i+j+1})$ for $i,j\in \ZZ$, we find $b(x^e, x^e) = t(x^d) = -a_0$ and $b(x^e, x^{e-d}K)=0$. Thus $_R S = Rx^e \perp x^{e-d}K$. The submodule $V=R1_S \oplus \cdots \oplus Rx^{e-1}$ is a complemented Lagrangian in $x^{e-d}K$ with respect to $b$. Hence $x^{e-d}K= \MM(U)$, $U \cong V^*$ is metabolic by \ref{bfLG}\eqref{meta}. This reference also proves the last part.
If $d=2e$ is even, then $R1_S \oplus \cdots \oplus Rx^{e-1}$ is a Lagrangian of $({_R S}, b)$, proving the second formula in \eqref{lem_transfer-2}.
\sm

\eqref{schatraco-a} The formulas follow from the corresponding ones in \eqref{lem_transfer-b} taking into account \eqref{qfba-tens-p} and Lemma~\ref{lem_metabolic}.
\lv{
We give the proof for $d$ odd; the even rank case follows in the same way. We have
\begin{align*}
  t_*(q_S) & \cong t_*( \lan 1_S \ran_b \ot_S q_S) & (\text{by \ref{qfba-tens-2}}) \\
   &\cong t_*(\lan 1_S \ran_b) \ot_R q &(\text{by \ref{Schtra}\eqref{prop_scharlau_1}})\\
  &\cong \big(\lan 1_R\ran_b \perp \MM(U)\big) \ot_R q
      &(\text{by \eqref{lem_transfer-1}})\\
   &\cong (\lan 1_R\ran_b \ot_R q) \perp (\MM(U)\ot q)
   &(\text{by \eqref{qfba-tens-p}})
    \\
   & \cong q \perp \HH(U \ot_R M) &(\text{by \ref{lem_metabolic}})
\end{align*}
Proof for $d$ even:
\begin{align*}
  t_*(q_S) & \cong t_*( \lan 1_S \ran_b \ot_S q_S) & (\text{by \ref{qfba-tens-2}}) \\
   &\cong t_*(\lan 1_S \ran_b) \ot_R q &(\text{by \ref{Schtra}\eqref{prop_scharlau_1}})\\
  &\cong \big(\lan 1_R, - \rmN_{S/R}(x)\ran_b \perp \MM(U)\big) \ot_R q
      &(\text{by \eqref{lem_transfer-1}})\\
   &\cong (\lan 1_R\ran_b \ot_R q) \perp (\MM(U_2)\ot q)
    &(\text{by \eqref{qfba-tens-p}})
   \\
   & \cong (\lan 1, - \rmN_{S/R}(x)\ran_b \ot q ) \perp \HH(U_2 \ot_R M) &(\text{by \ref{lem_metabolic}})
\end{align*}
}
\end{proof}

%


\begin{cor}[{\cite[3.8]{OP}} if $2\in R\ti$] \label{schatraco} We use the notation of {\rm \ref{lem_transfer}},  and
denote by $r_{S/R} \co \Wq(R) \to \Wq(S)$, $[q] \mapsto [q_S]$, the restriction homomorphism of \eqref{wgg2}.
Then \[ t_* \circ r_{S/R} = \Id_{\Wq(R)}, \]
in particular $r_{S/R}\co \Wq(R) \to \Wq(S)$ is injective.
\end{cor}
\begin{proof} For a regular quadratic $R$--form $q$ we have $t_*(q_S) = t_*(\lan 1_S \ran_b) \ot q$ by Frobenius reciprocity \ref{Schtra}\eqref{prop_scharlau_1}. Hence, the first formula  in \ref{lem_transfer}\eqref{schatraco-a} shows that $q$ and $t_*(q_S)$ are Witt-equivalent. \end{proof}
\ms

\textbf{Remarks.} (1) Injectivity of $r_{S/R}\co \Wq(R) \to \Wq(S)$ in \ref{schatraco} is part of condition \ref{equi-E} of \ref{equi}, which is implied 
by condition \ref{equi-D} of \ref{equi}. We have seen in \ref{prop_odd}\eqref{prop_odd-a} that \ref{equi-D} and hence \ref{equi-E} hold,  whenever $R$ is semilocal and $S\in \Ralg$ is one-generated and of odd degree. Thus Corollary~\ref{schatraco} establishes a part of \ref{equi-E} for a not necessarily semilocal $R$,  at the expense of requiring $S$ to be unit-generated.
\sm

(2) That $r_{S/K} \co \Wq(R) \to \Wq(S)$ is a split monomorphism and hence injectivity of $r_{S/R}$ holds in greater generality, namely whenever $S\in \Ralg$ can be written as a tower $R=R_0 \to R_1 \to \cdots R_{n-1} \to R_n$ where each $R_i\to R_{i+1}$ is a unit-generated algebra of odd rank. For example, this is so if $R$ is a field and $S/R$ is a field extension of odd degree \cite[]{Lam-qf}.

\subsection{Similarity factors}\label{sf} Let $(M,q)$ and $(M,q')$ be quadratic modules over $R$. Keeping in mind \eqref{qfba-tens-2}, we denote by
\[ \rmG(q,q') = \{ u \in R\ti: q\cong u q'  \} = \{u \in R\ti: q \cong \lan u \ran_b \ot q'\}  \]
the set of {\em similarity factors of $q$ and $q'$}.
We mention some obvious facts: \sm

\begin{inparaenum}[(a)] \item \label{sf-a} $\rmG(q,q') = \rmG(q',q)$, and if $\rmG(q,q') \ne \emptyset$, say $u\in \rmG(q,q')$, then $\rmG(q,q') = u \rmG(q')$ is a coset of the group $\rmG(q')\subset R\ti$. Hence
\begin{equation}
  \label{sf-10} R\ti{}^2 \rmG(q,q') \subset \rmG(q,q').
\end{equation}
Thus, $u\in R\ti$ is an element of $\rmG(q,q')$ if and only if there exists $v\in R\ti$ such that $uv^2\in \rmG(q,q')$.
\sm

\item \label{sf-b} Suppose $R=R_1 \times \cdots \times R_n$ is a direct product of rings.
Then $(M,q)$ uniquely decomposes as a direct product $(M,q)= (M_1, q_1) \perp \cdots \perp (M_n,q_n)$, where each $(M_i, q_i)$ is a quadratic module over $R_i$ with $(M,q)$ being regular (or nonsingular) if and only if every $(M_i, q_i)$, $i=1, \ldots, n$, is regular (nonsingular respectively), \ref{qf}\eqref{qf-redc}. 
The same decomposition of $M$ gives rise to the direct product $(M,q')= (M_1, q'_1) \perp \cdots \perp (M_n,q'_n)$, and then to the decomposition
 \begin{equation}
  \label{sf-2} \rmG(q,q') = \rmG(q_1,q_1') \times \cdots \times \rmG(q_n,q'_n).
\end{equation}
We will use this to reduce the study of $\rmG(q,q')$ of arbitrary forms to that of $\rmG(q,q')$ of forms defined on finite projective $R$--modules of constant rank. \sm

\item \label{sf-d} Let $S\in \Ralg$ whose underlying $R$--module is projective of finite type. We will be interested in the connection between $\rmG(q,q')$ and $\rmG(q_S, q'_S)$, in particular in the implication
\begin{equation}\label{sf3}
 u\in \rmG(q_S,q'_S) \quad \implies \quad \rmN_{S/R}(u)\in \rmG(q,q')
\end{equation}
for specific elements $u$, and for all similarity factors, i.e.,  \begin{equation}\label{sf4}
 \rmN_{S/R}\big( \rmG(q_S, q'_S))\subset  \rmG(q,q').
\end{equation}
Some useful reductions:

\begin{inparaenum}[(i)] \item \label{sf-di} Let $u\in S\ti$.
Since $\rmN_{S/R}(S\ti{}^2) \subset R\ti{}^2$,
it follows from \eqref{sf-10} that $u\in \rmG(q_S, q'_S)$ satisfies
\eqref{sf3} if and only if there exists $v\in S\ti$ such that
\eqref{sf3} holds for $uv^2$.

\item \label{sf-dii} Let $S=S_1 \times \cdots \times S_n$ be a direct product of $R$--algebras. Hence each $S_i$ is a finite projective $R$--module and the decomposition    \eqref{sf-2} holds for $\rmG(q_S, q'_S)$. Applying \ref{trno}\eqref{trno-c} to $u=(u_1, \ldots, u_n)$ with $u_i \in \rmG(q_{S,i}, q'_{S,i})$, $i=1, \ldots, n$, we obtain that $\rmN_{S/R}(u) \in \rmG(q,q') \iff$ every $\rmN_{S_i/R}(u_i) \in \rmG(q,q')$.
\end{inparaenum} \end{inparaenum}
\sm

The next lemma establishes \eqref{sf3} in a special case. It is a first step towards towards the proof of Proposition~\ref{prop_oddd}.

\begin{lem} \label{T} Let $R$ be a semilocal ring, and let $S\in \Ralg$ be finite projective and unit-generated by $x\in S\ti$, hence\/ {\rm (\ref{genolemLG})} we know $S=R[X]/(P)$ where $P$ is a monic polynomial, say of degree $d\in \NN_+$. Further assume that $(M,q)$ and $(M',q')$ faithful regular quadratic $R$--modules. Then for any $v\in S\ti$ the implication
\begin{equation}\label{T1}
    v^2x \in \rmG(q_S,q_S') \quad \implies \quad \rmN_{S/R}(v^2x) \in \rmG(q,q')
\end{equation}
holds under any one of the following conditions,
\begin{enumerate}[label={\rm (\alph*)}]
 \item\label{Ta} $d$ is odd, or
 \item \label{Tb} $q=q'$.
\end{enumerate}
\end{lem}

\begin{proof} By \ref{sf}\eqref{sf-di} we know $v^2x\in \rmG(q_S, q'_S) \iff x\in \rmG(q_S, q'_S)$ and $\rmN_{S/R}(x) \in \rmG(q,q') \iff \rmN_{S/R}(v)^2\rmN_{S/R}(x) = \rmN_{S/R}(v^2x) \in G(q,q')$. It therefore suffices to prove \eqref{T1} for $v=1$.

We abbreviate $\lan \,\cdot \, \ran_b = \lan \,\cdot \,\ran$, and
use the Frobenius algebra structure $(S,t)$ of Lemma~\ref{lem_transfer}. Since $q_S \cong  \lan 1_S \ran \ot q_S \cong \lan x \ran \ot q'_S$, Frobenius reciprocity~\ref{Schtra}\eqref{prop_scharlau_1} shows
\begin{equation} \label{Ta1}
 t_*(\lan 1_S \ran) \ot q \cong t_*(\lan x \ran) \ot q'.
\end{equation}
At this point the proofs of \ref{Ta} and \ref{Tb} diverge. But observe that it is
enough to prove \ref{Ta} and then \ref{Tb} with $d$ even.
\sm

\ref{Ta} By the formulas in \ref{lem_transfer}\eqref{schatraco-a}
there exists free $R$--submodules $U\subset S$ and $U'\subset S$ such that
\[
q \; \perp \; \HH(U\ot_R M) \cong \lan \rmN_{S/R}(x) \ran_b \ot q' \; \perp \; \HH(U'\ot_R M').\]
Since $U\cong U'$ and $M \cong M'$ as $R$--modules, the hyperbolic spaces $\HH(U\ot_R M)$ and $\HH(U'\ot_R M'$ are isometric, so that our claim follows from Witt cancellation \ref{canqf}\eqref{canqf-a} .
\sm

\ref{Tb} with $d$ even: By the formulas in \ref{lem_transfer}\eqref{schatraco-a}
there exist free $R$--submodules $U_1$, $U_2$ and $L$ of $S$  with $L$ of rank $1$ such that
\begin{align*}
    &q \, \perp \, \big(\lan -\rmN_{S/R}(x) \ran \ot q\big) \, \perp \,
      \big(\HH(U_1 \ot M)\big)
    \\ &\quad  \cong \; \big(\HH(L \ot_R M \big)\, \perp \, \big(\HH(U_2 \ot_R M)\big).
\end{align*}
 Since $U_1$ and $U_2$ are free $R$--modules of the same rank, the hyperbolic spaces $\HH(U_1 \ot_R M)$ and $\HH(U_2 \ot_R M) $ are isometric.
We can therefore apply Witt cancellation \ref{canqf}\eqref{canqf-a}, conclude
$ q  \, \perp \, \big(\lan -\rmN_{S/R}(x) \ran \ot q \big) \, \cong \, \HH(L\ot_R M)$, and then get
\begin{align*}
& q  \perp  \big(\lan -\rmN_{S/R}(x) \ran \ot q \big)  \perp  \big(\lan \rmN_{S/R}(x) \ran \ot q \big)  \\ &\qquad \cong   \HH(L \ot_R M)  \perp \big(\lan \rmN_{S/R}(x) \ran \ot q\big) \end{align*}
By Corollary~\ref{carhyp-cor},  $\big(\lan -\rmN_{S/R}(x) \ran \ot q\big)  \perp  \big(\lan \rmN_{S/R}(x) \ran \ot q\big)  \cong \HH(M)$ is hyperbolic, and by Lemma~\ref{lem_metabolic}, $\HH(L \ot M) \cong \HH(M)$ because $L$ is free.
Thus \[  q \perp \HH(M) \cong \HH(M) \perp \big(\lan \rmN_{S/R}(x) \ran \ot q\big)\] and the claim \eqref{T1} follows by another application of Witt cancellation.
\end{proof}


\begin{prop}\label{prop_oddd} Let $R$ be a semilocal ring, let
$q$, $q'$ be regular quadratic forms and let $S$ be a finite \'etale  $R$--algebra of\/ {\em odd\/} rank. Then the analogue of \eqref{T1} holds for arbitrary $u\in S\ti$:
\begin{equation}\label{prop_oddd1}
 q_S \cong \lan u \ran_b \ot q'_S \quad \implies \quad q \cong \lan \rmN_{S/R}(u) \ran_b \ot q'.
\end{equation}
\end{prop}

\begin{proof}
Let $S/R$ have odd rank $d$ and let $e\in \NN$ be odd and satisfying $e \ge d^2 + d + 1$. By \ref{prop_tiroir} there exists a finite \'etale $S$--algebra $T$ of constant rank $e$ such that the $R$--algebra $T$ is unit-generated by $x=b^2u$ for some $b\in T\ti$ (note that we can view $S \subset T$ so that $b^2 u$ makes sense). Observe $\lan x \ran_b \cong \lan u \ran_b$ as bilinear forms over $T$. Hence $q_T \cong \lan x \ran_b  \ot q'_T$. By \ref{T}\ref{Ta} for $S$ replaced by $T$, we get $q \cong \lan \rmN_{T/R}(x)\ran_b \ot q'$. By standard properties of norms, see \ref{trno}, we have $\rmN_{T/R}(x) = \rmN_{T/R}(b^2) \rmN_{T/R}(u) = \rmN_{T/R}(b)^2 \rmN_{S/R}(u)^e \equiv \rmN_{S/R}(u) \mod R\ti{}^2$ and therefore $\lan \rmN_{T/R}(x)\ran_b \cong \lan \rmN_{S/R}(u)\ran_b$, proving \eqref{prop_oddd1}.    \end{proof}
\lv{
\begin{proof} We use again the abbreviation $\lan \,\cdot \, \ran_b = \lan \,\cdot \,\ran$ and first prove a special case. \sm

{\em Assume  all residue fields of $R$ are infinite.} By Corollary~\ref{corgen} there exists $v\in S\ti$ such that $x=uv^2$ is a primitive element of $S$. Observe $\lan u \ran \ot q'_S \cong \lan x \ran \ot q'_S$ and $\rmN_{S/R}(u) \equiv \rmN_{S/R}(x) \mod R\ti{}^2$. It is therefore enough to prove \eqref{prop_oddd1} for $u=x$. But this is Lemma~\ref{T}\ref{Ta}.
\sm

{\em General case.} Using Noetherian reduction \ref{snp-noe}, we can assume that $R$ is the semilocalization of a finitely generated $\ZZ$-algebra. Proposition~\ref{lem_gonflement} provides an ind-\'etale $R$--algebra $R_\infty=\limind R_j$, which is a semilocal ring with infinite  residue fields and which is a tower of cubic \'etale unit-generated extensions $R_j$:
\[\xymatrix{
 S \ar[r] & S_j = S \ot_R R_j \ar[r] & S_\infty = S \ot_R R_\infty  \\
 R \ar[r] \ar[u] & R_j  \ar[r]\ar[u] & R_\infty \ar[u]}
\]
Observe $\rmN_{S_\infty/R_\infty}(u \ot 1_{S_\infty}) = \rmN_{S/R}(u) \ot 1_{R_\infty}$. The special case above applies to the isometry $q_{S_\infty} \cong \lan u \ot 1_{S_\infty}\ran \ot q'_{S_\infty}$ and yields an isometry $q_{R_\infty} \cong \lan \rmN_{S/R} (u)\ot 1_{R_\infty}\ran \ot q_{R_\infty}$. Since all data involved in the latter isometry behave well with respect to direct limits, there exists $j\in \NN$ such that $q_{R_j} \cong \lan \rmN_{S/ R}(u)\ot 1_{R_j}\ran \ot q_{R_j} = ( \lan \rmN_{S/R}(u)\ran \ot q')_{R_j}$.
Since $R_j$ is an \'etale $R$--algebra of odd degree $3^j$, Corollary~\ref{prop_odd}\eqref{prop_odd-a} implies that  $q$ is  isometric to $\langle  \rmN_{S/R}(u) \rangle \otimes q'$. \end{proof}}

\comments{I re-instated the following example. It was in a previous version but was cut to save space. That was before we wrote the Springer paper.
}

%
%
\ms

The following is an immediate consequence of Proposition~\ref{prop_oddd}.

\begin{cor} \label{cor_oddd} Let $R$ be a semilocal ring.
Let $S$ be a finite \'etale  $R$--algebra of {\em odd} degree $d$. Let $q$ and $q'$ be {\em regular} quadratic forms. Then the following are equivalent:

   \begin{enumerate}[label={\rm (\roman*)}]
 \item  \label{cor_oddd_i} $q$ is similar to $q'$;
 \item  \label{cor_oddd_ii} $q_S$ is  similar to $q'_S$.
 \end{enumerate}
\end{cor}
\sm

Corollary~\ref{cor_oddd} is known in case $R$ is a field, $S/R$ is a field extension of odd degree and $q$ and $q'$ are arbitrary quadratic forms, see for example \cite[p.478]{BQ} or \cite[Prop.~1.1]{Sivatski} for fields of characteristic $\ne 2$ and \cite[Thm.~6.5]{Ho2} for fields of characteristic $2$ and arbitrary quadratic forms. Compare \ref{cor_oddd} with \ref{prop_odd}\eqref{prop_odd-a} where the equivalence is stated for isometry instead of similarity. \sm

\subsection{Scharlau's Norm Principle}\label{snp}
We put
\[ \rmG(q) = \rmG(q,q) = \{ u \in R\ti: q\cong u q  \} \]
and recall that $\rmG(q)$ is a subgroup of $R\ti$ satisfying
\begin{equation}
  \label{sf-1} R\ti{}^2 \subset \rmG(q) \subset R\ti.
\end{equation}
\lv{
Since $\vphi \co m \mapsto v\me m$, $v\in R\ti$, is an isometry between $v^2 q$ and $q$, because $(v^2q)(\vphi m) = q(v v\me m) = q(m)$.
}
We have seen in \ref{simi-defi}\eqref{simi-difi-y} that $R\ti{}^2=\rmG(q)$ if $(M,q)$ has constant odd rank and that $\rmG(q) = R\ti$ if $q$ is hyperbolic.
\sm

Let $S\in \Ralg$ be locally free of finite type as $R$--module and let $\rmN_{S/R}$ be the norm of $S$, as in \ref{trno}. We will say that the {\em Scharlau Norm Principle} holds for $(R,S)$
if for all quadratic $R$--spaces $(M,q)$  we have
    \begin{equation} \label{sf-snp1}
        \rmN_{S/R}\big( \rmG(q_S)\big) \subset \rmG(q).
    \end{equation}
In other words, for $x\in S\ti$,
\begin{equation}
  \label{sf-snp-2} q_S \cong x q_S \implies q \cong \rmN_{S/R}(x)q.
\end{equation}

Scharlau's Norm Principle holds in case $S/R$ is a finite extension of fields: It was proven by Scharlau for fields of characteristic $\ne 2$ in \cite{Scha}, see
\cite[2.8.6]{Sc} or \cite[VII, 4,3]{Lam-qf}. A proof for arbitrary fields, but even-dimensional forms is given in \cite[20.14]{EKM} (the restriction to even-dimensional forms is natural in view of \ref{sfod}). It was recently proven for fields of characteristic $2$ and arbitrary quadratic forms in \cite[7.3]{Ho2}. We will give a proof for nonsingular forms and finite field extensions in
Lemma~\ref{snp-fi} and for finite \'etale extensions of semilocal rings in Theorem~\ref{thm_scharlau}. 


The following Lemma~\ref{sfod} implies that Scharlau's norm principle holds for quadratic spaces of odd rank over semilocal rings.


\begin{lem}[Odd rank] \label{sfod} Scharlau's norm principle holds for any finite projective $S\in \Ralg$ and $(M,q)$ a quadratic $R$--space for which $M$ has constant {\em odd\/} rank.
\end{lem}

\begin{proof} Since $M_S$ has odd rank, we know $\rmG(q_S) = S\ti{}^2$ by \eqref{simi-difi-y3}. Hence $\rmN_{S/R}\big( \rmG(q_S)\big) = \rmN_{S/R}(S\ti{}^2) = \rmN_{S/R}(S\ti)^2 \subset \rmG(q)$.  \end{proof}

\comments{(2023-01) New lemma \ref{snp-fi}; it would be nice to have \ref{snp-fi} for Artinian $F$--algebras. By \ref{sf}\eqref{sf-dii}, it suffices to prove this for a local Artinian $F$--algebra.}

\pcomments{(2026-06-08) Tu peux utiliser le m\^eme argument que dans la lemme \ref{lem_norm_field}: pour une $F$--alg\`ebre artinienne local $R$ de corps r\'esiduel $L/F$, alors la norme d'un \'el\'ement $r\in R$ est une puissance de $\rmN_{L/k}(F)$ o\`u $a$ est l'image de $r$ dans $L$. Ainsi le r\'esultat sur les corps entraine celui sur les alg\`ebres artiniennes.}

\comments{(2023-11) Is Scharlau's Norm Principle true for absolutely flat (= von Neumann regular) rings? Note  over such a ring a finitely generated module is finite projective. Hence part of the proof of \ref{snp-fi} will work.}
\pcomments{(2026-06-08) Selon Wikipedia, un tel anneau absolument plat est un produit de corps, don la r\'eponse \`a cette question est positive, mais cela ne me semble pas tr\`es int\'eressant.}

\begin{lem}[Scharlau's norm principle for fields]\label{snp-fi} Let $F$ be a field and let $E\in \Falg$ be a finite-dimensional reduced $F$--algebra. Then Scharlau's norm principle holds for $E/F$. \end{lem}

\begin{proof} The $F$--algebra $E$ is a finite direct product of extension fields of $F$. 
By \ref{sf}\eqref{sf-dii} it thus suffices to prove the lemma for $E=K$ a finite field extension of $F$.

Let $(M,q)$ be a nonsingular quadratic space over $F$. By \ref{sfod} we can assume that $M$ has even dimension and hence that $q$ is regular. Let $x\in \rmG(q_K)$.

If $K=F(x)$, we can apply \ref{T}\ref{Tb} to conclude $\rmN_{K/F}(x) \in \rmG(q)$. Otherwise, $[K:F(x)]= n >0$. We then have $\rmN_{K/F}(x) = \rmN_{F(x)/F}\big( \rmN_{K/F(x)}(x)\big) = \rmN_{F(x)/F}(x^n) = \rmN_{F(x)/F}(x)^n$. If $n$ is even, $\rmN_{K/F}(x) \in F\ti{}^2 \subset \rmG(q)$ by \eqref{sf-1}, and if $n$ is odd, then $\rmN_{K/F}(x) \equiv \rmN_{F(x)/F}(z) \mod F\ti{}^2$ follows. It therefore suffices to show $\rmN_{F(x)/F}(x)\in \rmG(q)$. Since $(q_{F(x)})_K = q_K \cong x q_K \cong (x q_{F(x)})_K$, we can apply \ref{prop_odd} and obtain $q_{F(x)} \cong x q_{F(x)}$, so that $q\cong \rmN_{F(x)/F}(x)q$ follows from the case $K=F(x)$ considered before. Because $N_{F(x)}(x) \equiv N_{K/F}(x) \mod F\ti{}^2$, we are done. \end{proof}

The following lemma says that one can study the implication \eqref{sf3} using
noetherian reduction. We will use it in the proof of our main result \ref{thm_scharlau}.

\begin{lem} \label{snp-noe} Let $R=\limind_{\la \in \La} R_\la$, where $(R_\la)_{\la \in \La}$ is a directed system of $R_0$--rings for some $R_0\in \ZZalg$.
We suppose that we are given data $(M,q,q', S,u)$ consisting of two quadratic $R$--modules $(M,q)$ and $(M,q')$, and a finite projective $S\in \Ralg$ with $u\in S\ti$. Then there exists $\la \in \La$ and data $(M_\la, q_\la, q'_\la, S_\la, u_\la)$ over $R_\la$ of the same type such that
\[ (M_\la, q_\la, q'_\la, S_\la, u_\la)\ot_{R_\la} R  \cong (M,q,q',S,u)
\]
Moreover,
\begin{enumerate}[label={\rm (\roman*)}]
 \item  \label{snp-noei} if $q$ and $q'$ are both regular or both nonsingular, the same holds for $q_\la$ and $q'_\la$,

\item  \label{snp-noeii} if $S$ is finite \'etale or finite \'etale of degree $d$, the same holds for the $R_\la$--algebra $S_\la$.

\item If the $R_\la$--data satisfy
\begin{equation} \label{snp-no1}
 q_{S_\la} \cong \lan u_\la \ran_b \ot q'_{S_\la} \quad \implies \quad
     q_\la \cong \lan \rmN_{S_\la/R_\la}(u_\la) \ran_b \ot q_\la,
\end{equation}
i.e., \eqref{sf3} holds for $u_\la$, then the analogous implication holds for $(M,q,q',S,u)$.
 \end{enumerate}\end{lem}

\begin{proof}
  The existence of $(M_\la, q_\la, q'_\la, S_\la)$ with properties \ref{snp-noei} and \ref{snp-noeii} follows from \ref{dirquad}.
   The existence of $u_\la$ is obvious, as is the last statement.
 \end{proof}

\begin{thm}[Scharlau's Norm Principle for finite \'etale extensions] \label{thm_scharlau} Let $R$ be a semilocal ring and let $S$ be a finite \'etale extension of $R$.
Let $q$ be a nonsingular quadratic form over $R$ and let $\rmG(q) \subset R^\times$ be the subgroup of similarity factors. Then $\rmN_{S/R}\bigl( \rmG(q_S) \bigr) \subset \rmG(q)$.
\end{thm}

\begin{proof}
We denote by $\kappa_1,\dots, \kappa_c$ the residue fields of the maximal
ideals of $R$. We are given $u \in \rmG(q_S)$ and want to show that $\rmN_{S/R}(u) \in \rmG(q)$. \sm

{\em Reduction to $S$ and $M$ of constant rank.} We decompose $R=R_1 \times \cdots \times R_n$ such that $S=S_1 \times \cdots \times S_n$ with each $S_i$ being an \'etale $R_i$--algebra of finite constant rank. Then each $R_i$ is semilocal and the principle \eqref{sf-2} applies to the corresponding decomposition $q = q_1 \times \cdots \times q_n$. Since $\rmN_{S/R}(s_1, \ldots, s_n) = (\rmN_{S_1/R_1}(s_1), \ldots ,\rmN_{S_n/R_n}(s_n))$, it is no harm to assume that $S$ has constant finite degree $d$ and that $M$ has constant rank.  \sm

{\em Reduction to $q$ regular.} If $M$ has constant odd rank, the claim has been proved in Lemma~\ref{sfod}. We can therefore assume that $M$ has constant even rank, in which case $q$ is regular by \ref{qf}\eqref{quadfoe}. \sm

{\em Special case where all $\ka_i$ are infinite.} By Lemma~\ref{hunco} there exists $v\in S\ti$ such that $x=uv^{-2}$ is a primitive element of $S$. Then Lemma~\ref{T}\ref{Tb} proves $\rmN_{S/R}(u) \in \rmG(q)$. \sm

{\em General case.}  Using Noetherian reduction \ref{snp-noe}, we can assume that $R$ is the semilocalization of a finitely generated $\ZZ$-algebra \ref{noethred}\eqref{noethred-b}. Proposition~\ref{lem_gonflement} provides an ind-\'etale $R$--algebra $R_\infty=\limind R_j$, which is a semilocal ring with infinite  residue fields and which is a tower of cubic \'etale unit-generated extensions $R_j$:
\[\xymatrix{
 S \ar[r] & S_j = S \ot_R R_j \ar[r] & S_\infty = S \ot_R R_\infty  \\
 R \ar[r] \ar[u] & R_j  \ar[r]\ar[u] & R_\infty \ar[u]}
\]
Observe $\rmN_{S_\infty/R_\infty}(u \ot 1_{S_\infty}) = \rmN_{S/R}(u) \ot 1_{R_\infty}$. The special case above applies to the isometry $q_{S_\infty} \cong \lan u \ot 1_{S_\infty}\ran \ot q_{S_\infty}$ and yields an isometry $q_{R_\infty} \cong \lan \rmN_{S/R} (u)\ot 1_{R_\infty}\ran \ot q_{R_\infty}$. Since all data involved in the latter isometry behave well with respect to direct limits, there exists $j\in \NN$ such that $q_{R_j} \cong \lan \rmN_{S/ R}(u)\ot 1_{R_j}\ran \ot q_{R_j} = ( \lan \rmN_{S/R}(u)\ran \ot q)_{R_j}$.
Because $R_j$ is an \'etale $R$--algebra of odd degree $3^j$, Corollary~\ref{prop_odd}\eqref{prop_odd-a} implies that  $q$ is  isometric to $\langle  \rmN_{S/R}(u) \rangle \otimes q'$. \end{proof}



\subsection{Summary} \label{snp-sum} Let $R$ be a semilocal ring, let $S\in \Ralg$ be finite projective, and let $q\co M \to R$ be a quadratic form. We have established Scharlau's Norm Principle in the following settings:
\[
\begin{tabular}{c||c |c|c}
  & $S$ & $M$ & $q$ \\
  \hline \\
  \ref{T} & \text{unit-generated by $x$, for $x$} & \text{arbitrary} & \text{regular} \\
 \ref{prop_oddd} & \text{$S/R$ finite \'etale of odd rank} &  \text{arbitrary} & \text{regular}
 \\
 \ref{sfod} & \text{arbitrary} & \text{odd rank} & \text{nonsingular}
 \\
 \ref{snp-fi} & \text{$R=F$ field, $S$ reduced $F$--algebra} & \text{arbitrary} & \text{nonsingular}
 \\
 \ref{thm_scharlau} & \text{finite \'etale} & \text{arbitrary} & \text{nonsingular}
\end{tabular}
\]
\comments{(2022-05) Analysis:
\begin{enumerate}
 \item \ref{T} is used in the proof of \ref{prop_oddd} and \ref{sflem}, while \ref{sflem} is used in the proof of \ref{thm_scharlau};
 \sm

 \item Proposition~\ref{prop_oddd} is not used anywhere, except thast we use the proof principle in \ref{thm_scharlau}; \sm

 \item \ref{sfod} is used in the proof of \ref{thm_scharlau};
\sm

  \item \ref{prop_oddd} is not a special case of \ref{thm_scharlau} because it deals with two forms $q$ and $q'$;
\sm

 \item in the spirit of Springer's Theorem, we are {\color{red} missing the case
 \[ \text{$S$ unit-generated, but not \'etale, $M$ arbitrary, $q$ nonsingular}\]}
\end{enumerate}
}

\newpage

\section{Knebusch's Norm Principle}\label{sec:kneb}


\comments{(2021-03) A version of \ref{thm-kneb} for $S$ one-generated, not necessarily \'etale is missing. Such a version is proven in Kirill's paper \cite{Kirill} for $R$ a semilocal domain with all residue fields of characteristic $0$ and $S\in \Ralg$ a simple extension of degree $n$. }

\comments{(2022-06) A ``norm theorem" for semisingular quadratic forms  is established in the 2021-paper by Laghribi; need to investigate the relevance for our work}

The principal goal of this section is to prove Knebusch's norm principle \ref{thm-kneb} for semilocal rings $R$, relating the values of a quadratic space $(M,q)$ over $R$, the norm of a finite \'etale $S\in \Ralg$, and the values of the extended quadratic form $(M,q)_S$: 
\[ \xymatrix@C=50pt{
   M \ar[d]_q \ar[r] & M_S \ar[d]^{q_S} \\
   R\ar@<-0.5ex>[r] & S \ar@<-0.5ex> [l]_{\rmN_{S/R}} 
}\]
See \ref{bate}\eqref{bate-a} for some background. We recall that $\wS_q = \{m\in M: q(x) \in R\ti\}$ is the extended sphere of a faithful quadratic $R$--module $(M,q)$, defined in \eqref{quadco-aa2}. For the scheme $\uS_{q,u}$, $u\in R\ti$ and its smooth part $\uS\rmsm_{q,u}$ see \ref{ussl}.
Furthermore, we put
\begin{equation*}\label{thm-kneb1} \begin{split}
  \rmD(q) &= \textstyle\bigsqcup_{u\in R\ti}\, q\big(\uS_{q, u}(R)\big)
      = R\ti \cap q(M) = q(\wS_q), \\
  \rmD(q)^{[n]}  &= \rmD(q) \cdots \rmD(q) \quad (\text{$n$ factors}), \\
  \rmD(q)^{[\rm ev]} &= \textstyle \bigcup_{n >0 \text{ even}} \Dqn,  \quad \qquad
  \rmD(q)^{[\rm od]} = \textstyle \bigcup_{n \text{ odd}} \Dqn,  \\
  \rmD\rmsm(q) &= \textstyle\bigsqcup_{u\in R\ti}\, q\big(\uS_{q, u}\rmsm(R)\big),
 \end{split}\end{equation*}
and define $\rmD\rmsm(q)^{[n]}$ by replacing $\rmD(q)$ by $\rmD\rmsm(q)$ in the formula above. The main result of this section is 

\begin{thm}[Knebusch's norm principle for \'etale $S$]\label{thm-kneb} Let $R$ be a semilocal ring,  let $(M,q)$ be a quadratic space over $R$, and let 
$S\in \Ralg$ be a finite \'etale $R$--algebra of constant degree $d$. \sm

\begin{enumerate}[label={\rm (\alph*)}]
\item \label{thm-kneb-b} Then 
\begin{equation}      \label{thm-kneb-b1}
    \rmN_{S/R}(\rmD(q_S)\big)) \subset
       \begin{cases} \rmD(q)^{[\rm ev]}, & \text{if 
       $d$ is even,} \\
        \rmD(q)^{[\rm od]}, & \text{if $d$ is  odd}, \end{cases}
    \end{equation}
in particular
\begin{equation}\label{thm-kneb-b2}
\rmN_{S/R}\big( \rmD(q_S)^{[\rm ev]} \big) \subset \rmD(q)^{[\rm ev]}.
\end{equation}
Moreover, if $\rank M \ge 4$, we can replace $\rmD(q)$ by $\rmD\rmsm(q)$ on the right-hand sides of \eqref{thm-kneb-b1} and \eqref{thm-kneb-b2}. \sm

\item\label{thm-kneb-a} If all residue fields of $R$ are infinite,  then
\begin{equation} \label{thm-kneb-a1}
\rmN_{S/R}\big(\rmD(q_S)\big) \subset
 \begin{cases}  \rmD(q)^{[d]} & \text{if $q$ is regular;} \\
        \rmD(q)^{[2d]} & \text{if 
            $d$ is even;}\\
         \rmD(q)^{[2d+1]} & \text{if $d$ is odd, }
 \end{cases} \end{equation}
and
\begin{equation}  \label{thm-kneb-a11}
 \rmN_{S/R}\big( \rmD\rmsm(q_S)\big) \subset \rmD(q)^{[d]}.
\end{equation}
Moreover, if $\rank M \ge 4$, we can replace $\rmD(q)$ by $\rmD\rmsm(q)$ on the right hand sides of \eqref{thm-kneb-a1} and \eqref{thm-kneb-a11}.
\end{enumerate}
\end{thm}
\ms


We note that by Example~\ref{lem_smooth_locus_exam}\eqref{lem_smooth_locus_exam-c} the
condition $\uS_{q,1}\rmsm(R) \ne \emptyset$ can always be achieved by passing from $q$ to $uq$ for a suitable $u\in R\ti$. \sm

We will prove \eqref{thm-kneb-a11} in Lemma~\ref{pkn}, and the rest of 
Theorem~\ref{thm-kneb} in \ref{thm-kneb-proof}. A large part of the development leading up to the proof of Theorem~\ref{thm-kneb} holds for arbitrary $R$ and for a unit-generated finite free $S\in \Ralg$ of rank $d$. This is setting of \ref{dqd-ele}--\ref{lem_representability}.

\subsection{Some elementary facts about $\rmD(q)$ and $\rmD\rmsm(q)$}\label{dqd-ele} 
Unless stated otherwise, $R$ is arbitrary, $(M,q)$ is a faithful quadratic module over $R$, and $S\in \Ralg$ is faithfully projective. To avoid duplication of formulas we will use the abbreviation
\[ \rmD' \in \{ \rmD, \rmD\rmsm\} \]
and employ the rule that whenever two symbols $\rmD'$ appear in a formula then $\rmD'$ has the same value every time.

We suppose $\rmD(q) \ne \emptyset$, which is always the case if $R$ is an LG ring and $q$ is nonsingular, \ref{LGqdi}. The assumption $\rmD(q) \ne \emptyset$ implies that $(M,q)$ is primitive, see \ref{quadco}. 
Moreover, unless stated otherwise, in all formulas involving $\rmD\rmsm(q)$ it is also assumed that $\rmD\rmsm(q) \ne \emptyset$. Recall from Example~\ref{lem_smooth_locus_exam}\eqref{lem_smooth_locus_exam-c} that  $\rmD\rmsm(q)\ne \emptyset$ whenever $R$ is an LG ring and $(M,q)$ is a quadratic $R$--space with $\rank M \ge 2$, but see \ref{rkn}\eqref{dqd-ele-d} below. We also remind the reader that $\rmD(q) = \rmD\rmsm(q)$ if $q$ is regular, 
\ref{lem_smooth_locus-LG}\ref{lem_smooth_locus-LGd}.\sm

\begin{inparaenum}[(a)] 
 \item \label{dqd-ele-a} If $q_1$ is isometric to $q$, then $\rmD'(q) = \rmD'(q_1)$. Putting $\rmG(q) = \{ u \in R\ti: uq \cong q\}$, we therefore have $\rmG(q) \cdot  \rmD'(q) = \rmD'(q)$, 
      in particular
      \begin{equation}\label{dqd-ele-a1}
        R\ti{}^2 \cdot \rmD'(q) =  \rmD'(q).
     \end{equation}
     In general,
     \begin{equation}\label{dqd-ele-a2}
            1\in R\ti{}^2 \subset \rmD'(q)^{[2]}
     \end{equation}
      since for any $x\in \rmD'(q)$ and $u\in R\ti$ we have $(u/x)^2 x \in \rmD'(q)$ by \eqref{dqd-ele-a1} and then  
      $u^2 = x\big( (u/x)^2 x\big)\in \rmD'(q)^{[2]}$. 
      The inclusion \eqref{dqd-ele-a2} implies
      \begin{equation}\label{dqd-ele-moti}
       \rmD'(q)^{[n]} \subset \rmD'(q)^{[n+2m]}
      \end{equation}
       for any $n,m\in \NN_+$. Therefore
  \begin{equation} \label{dqd-ele-a4} \begin{split}
    R\ti{}^2 \subset \rmD'(q)^{[2]} &\subset \rmD'(q)^{[4]} \subset \cdots
         \rmD'(q)^{[2n]}\subset \cdots \subset \rmD'(q)^{[\rm ev]}, \quad \text{and} \\
 \rmD'(q) &\subset \rmD'(q)^{[3]} \subset \cdots
         \rmD'(q)^{[2n+1]}\subset \cdots \subset \rmD'(q)^{[\rm od]}.
  \end{split} \end{equation}
Another consequence of \eqref{dqd-ele-a1} is
     \begin{equation} \label{dqd-ele-a3}
      x\in \rmD'(q) \quad \implies \quad x\me \in \rmD'(q)
      \end{equation}
   since $x\me = x^{-2} x$. It follows that 
\begin{equation}\label{dqd-ele-a5}
 \rmD(q)^{[\rm ev]} = \lan \rmD(q)^{[2]}\ran \quad \text{and}\quad 
     \lan \rmD(q)\ran = \textstyle \bigcup_{n\in \NN_+} \rmD(q)^{[n]}
\end{equation} 
   are subgroups of $R\ti$, generated by $\rmD(q)^{[2]}$ and $\rmD(q)$ respectively. 
      \sm

\item \label{qdq-ele-aa}    {\em Suppose $R\ti{}^2 \cap \rmD'(q) \ne \emptyset$, equivalently $1\in \rmD'(q)$}. Then $R\ti{}^2 \subset \rmD'(q)$ by \eqref{dqd-ele-a1}. Hence $\rmD'(q)^{[n]} \subset \rmD(q)^{[n+m]}$ for all $m,n\in \NN_+$, so
   \[     \rmD'(q) \subset \rmD'(q)^{[2]}  \subset \cdots
         \rmD'(q)^{[n]}\subset \cdots \subset \rmD'(q)^{[\rm ev]} =\rmD'(q)^{[\rm od]} =\lan \rmD'(q) \ran
   \]
   is the subgroup generated by $\rmD'(q)$. \sm

\item\label{dqf-ele-scal} ({\em Scaling}) For $u \in R\ti$ and $d\in \NN_+$ we have $\rmD'(uq) = u \rmD'(q)$ and therefore, by \eqref{dqd-ele-a1}, 
\[ \rmD'(uq)^{[d]} = \begin{cases}
   u \rmD'(q)^{[d]}, & \text{if $d$ is odd,} \\
   \rmD'(q)^{[d]}, & \text{if $d$ is even}.
\end{cases}
\]
In particular $\rmD'(uq)^{[\rm ev]} =  \rmD'(q)^{[\rm ev]}$.

\inparcom{(2026-03-09) 
({\em Reduction to $\rad(q) =0 $}) By definition in \ref{radqf}, $\rmD(\bar q) = \rmD(q)$. Hence, since $\rad(\bar q) = 0$, in order  to prove results about $\rmD(q)$ like the strong/weak Knebusch principle, we can always assume that $\rad(q) = 0$. } 
\end{inparaenum}

\subsection{Background, terminology}\label{bate} \begin{inparaenum}[(a)] \item \label{bate-a} Let $F$ be a field, let $K/F$ be a finite field extension and let $q\co M \to F$ be an arbitrary  quadratic form.
In this setting, the formula
\begin{equation}\label{bate1}
  \rmN_{K/F}\big(\rmD(q_K)\big) \subset \rmD(q)^{[d]}
\end{equation}
has been proven by Knebusch in \cite{knebusch-norm}, using the consequence \ref{equi}\ref{equi-C} of Springer's Odd Degree Theorem, as explained in \ref{equisa}\eqref{equisa-c}. The result had been conjectured by Witt; it is presented in \cite[VII, 5.1]{Lam-qf}. An exposition of the weaker version  
\begin{equation}\label{bate1}
  \rmN_{K/F}\big(\rmD(q_K)\big) \subset \lan \rmD(q) \ran 
\end{equation}
given in \cite[18.10]{EKM}, using ideas from \cite{Knebusch73}. \sm


The formula \eqref{thm-kneb-b2} was proven by Zainoulline in 
\cite{Kirill} 
assuming that $S/R$ is a finite extension of semi-local regular rings containing a field  of characteristic $0$ and that $q$ is regular. It was proven by Ojanguren-Panin-Zainoulline \cite{OPZ} 
for $S/R$ a finite \'etale extension of a semilocal noetherian domain $R$ with infinite residue fields of characteristic $\ne 2$ and a regular quadratic form.%
\sm

\item \label{bate-b} Let $(M,q)$ be a quadratic module over an arbitrary base ring $R$ and let 
$S\in \Ralg$ is faithfully projective. 

The set $\rmD(q)$ is called the {\em set of values of $q$} or {\em value set}. We will refer to $\rmD\rmsm(q)$ as the {\em set of smooth values of $q$}, and to the relation
\[     \rmN_{S/R}\big( \rmD(q_S)^{[\rm ev]} \big)
        \subset \rmD(q)^{[\rm ev]}
\]
as the {\em weak Knebusch norm principle}. Since $\rmN_{S/R}$ is multiplicative and since by \eqref{dqd-ele-a5} the group $\rmD(q_S)^{[\rm ev]}$ is generated by $\rmD(q_S)^{[2]}$, the weak principle is equivalent to 
\[     \rmN_{S/R}\big( \rmD(q_S)^{[2]} \big)
        \subset \rmD(q)^{[\rm ev]}.
\]
We will see in \ref{rkn}\eqref{dqd-ele-e} that the weak principle holds as soon as it holds for $S$ of constant rank.
\sm 

\item\label{bate-c} Suppose that $S\in \Ralg$ has constant degree $d\in\NN_+$. We will say that the {\em strong Knebusch norm principle holds for $a\in \rmD(q_S)$}, if 
\[  \rmN_{S/R}(a) \in \rmD(q)^{[d]},  
\]
and we will say that the {\em strong Knebusch norm principle holds\/}, if it holds for all $a\in \rmD(q_S)$, i.e.
\[
   \rmN_{S/R}\big(\rmD(q_S)\big) \subset \rmD(q)^{[d]}.
\]
Clearly, by multiplicativity of $\rmN_{S/R}$, the strong version implies the weak version.

\inparcom{(2026-03) With this terminology, Theorem~\ref{thm-kneb} then says that for a finite \'etale extension of  a semilocal ring $R$ the strong Knebusch norm principle holds for regular quadratic forms if  all residue fields of $R$ are infinite, while the weak principle is valid for nonsingular $q$ and arbitrary semilocal $R$. }
\end{inparaenum}

\subsection{Reductions for Knebusch norm principles}\label{rkn} Let $R$ be arbitrary, let $S\in \Ralg$ be faithfully projective and let $(M,q)$ be a quadratic $R$--module. We study some reductions for proving the various Knebusch norm principles, mainly related to small values of the ranks of $S$ and $M$. 
As in \ref{dqd-ele} we will use $\rmD' \in \{ \rmD, \rmD\rmsm\}$. \sm  

\begin{inparaenum}[(a)] 
\item  \label{dqd-ele-e} ({\em Direct products}) Let $R=R_1 \times \cdots \times R_n$ be a direct product. By \ref{qf}\eqref{qf-redc}, the quadratic module $(M,q)$ uniquely decomposes as a direct product $(M,q) = (M_1, q_1) \times \cdots \times (M_n, q_n)$, where each $(M_i, q_i)$, $1\le i \le n$, is a quadratic $R_i$--module, which is nonsingular if $(M,q)$ is so. We have  $\rmD'(q) = \rmD'(q_1) \times \cdots \times \rmD'(q_n)$.   An $S\in \Ralg$ which is finitely generated projective (or of rank $d\in \NN_+$ respectively) uniquely decomposes as $S=S_1 \times \cdots \times S_n$ where each $S_i$ is a finitely projective $R_i$--module (of rank $d$ respectively).
    It is then easily seen that the smooth/general strong or weak Knebusch norm principle for $(R,q,S)$ is equivalent with the same principle for $(R_i, q_i, S_i)$. For the weak principle, one may use \eqref{dqd-ele-a5}.
    
    In particular, one can find a decomposition of $R$ such that both $M$ and $S$ have constant rank.      
\lv{
    \begin{align*}
    R\ti &= R_1\ti \times \cdots \times R_n\ti, \\
    q(M)  &= q(M_1) \times \cdots q(M_n), \\
     R\ti \cap q(M) = \rmD(q) &= \rmD(q_1) \times \cdots \times D(q_n), \\
     \rmD(q)^{[m]} &= \rmD(q_1)^{[m]} \times \cdots \times \rmD(q_n)^{[m]}, \quad (m \in \NN_+), \\
      S &= S_1 \times \cdots \times S_n, \\
     (M,q)_S &= (M_1, q_1)_{S_1} \times \cdots \times (M_n, q_n)_{S_n}, \\
     \rmD(q_S) &= \rmD(q_{1,S_1}) \times \cdots \times D(q_{n, S_n}), \\
     \rmN_{S/R} &= \rmN_{S_1/R_1} \times \cdots \times \rmN_{S_n/R_n} , \\
     \rmN_{S/R}\big(\rmD(q_S)\big) & = \rmN_{S/R}\big( \rmD(q_{1,S_1}) \times \cdots \times D(q_{n, S_n}) )\\
     & = \rmN_{S_1/R_1}\big( \rmD(q_{1,S_1})\big)  \times \cdots \times \rmN_{S_n/R_n}\big(\rmD(q_{n, S_n})\big)
    \\ 
    \text{hence}\qquad \qquad & 
    \\  
    \rmN_{S/R}\big(\rmD(q_S)\big) &\subset \rmD(q)^{[d]} 
      = \rmD(q_1)^{[d]} \times \cdots \times \rmD(q_n)^{[d]} 
   \\&\iff \quad \rmN_{S_i/R_1}(\rmD(q_{i, S_i}))\subset \rmD(q_i)^{[d]}
      , \quad 1\le i \le n
    \end{align*}
For the weak principle we have
\begin{align*}
   \rmN_{S/R}\big(\rmD(q_S)^{[2]}\big) & = \rmN_{S/R}\big( \rmD(q_{1,S_1})^{[2]} \times \cdots \times D(q_{n, S_n})^{[2]} )\\
     & = \rmN_{S_1/R_1}\big( \rmD(q_{1,S_1})\big)^{[2]}  \times \cdots \times \rmN_{S_n/R_n}\big(\rmD(q_{n, S_n})\big)^{[2]}
 \\ 
    \text{hence}\qquad \qquad & 
 \\  
\rmN_{S/R}\big(\rmD(q_S)^{[2]}\big) &\subset \rmD(q)^{\rm ev} 
      = \rmD(q_1)^{\rm ev} \times \cdots \times \rmD(q_n)^{\rm ev} 
   \\&\iff \quad \rmN_{S_i/R_1}(\rmD(q_{i, S_i})^{[2]})\subset \rmD(q_i)^{\rm ev}
     , \quad 1\le i \le n
 \end{align*}
OLD VERSION:
\begin{align*}
   \rmD(q)^{[\rm ev]} &= \textstyle \bigcup_{m\in \NN_+} \rmD(q)^{[2m]} =
        \textstyle \bigcup_{m\in \NN_+}
          \Big( \rmD(q_1)^{[2m]} \times \cdots \times \rmD(q_n)^{[2m]}\Big) \\
      &\subset   \Big(\textstyle \bigcup_{m\in \NN_+} \rmD(q_1)^{[2m]}\Big)
           \times \cdots \times \Big(\textstyle \bigcup_{m\in \NN_+} \rmD(q_n)^{[2m]}\Big)\\
     &= \rmD(q_1)^{[\rm ev]} \times \cdots \times \rmD(q_n)^{[\rm ev]}
   \end{align*}
However, assume $1\in \rmD(q)$, equivalently $1\in \rmD(q_i)$ for $1\le i \le n$. Then
\begin{align*} \rmD(q_1)^{[m_1]} \times \dots \times \rmD(q_n)^{[m_n]}
 &= \rmD(q_1)^{[m]} \times \dots \times \rmD(q_n)^{[m]} \\ &= \rmD(q)^{[m]}
\end{align*}
for $m = \max(m_1, \ldots, m_n)$.
Hence
\begin{align*}
  \rmD(q_1)^{[\rm ev]} \times \cdots \times \rmD(q_n)^{[\rm ev]}
 & = \textstyle\bigcup_{m_1, \ldots, m_n \in \NN_+}
      \rmD(q_1)^{[2m_1]} \times \ldots \times \rmD(q_n)^{[2m_n]}\\
& \subset  \textstyle \bigcup_{m\in \NN} \rmD(q)^{[2m]} = \rmD(q)^{[\rm ev]}.
 \end{align*}
Therefore \[
\rmD(q)^{[\rm ev]} = \rmD(q_1)^{[\rm ev]} \times \cdots \times \rmD(q_n)^{[\rm ev]}
\]
and then for any $m\in \NN_+$,
\begin{align*}
  \rmN_{S/R}\big(\rmD(q_S)^{[m]}\big) &=  \rmN_{S_1/R_1}\big( \rmD(q_{1,S_1})^{[m]}\big)  \times \cdots \times \rmN_{S_n/R_n}\big(\rmD(q_{n, S_n})^{[m]}\big) \\
  &\subset \rmD(q)^{[\rm ev]} =
     \rmD(q_1)^{[\rm ev]} \times \cdots \times \rmD(q_n)^{[\rm ev]}
 \\ \iff &  \rmN_{S_i/R_i}\big(\rmD(q_{i, S_i})^{[m]}\big)  \subset
 \rmD(q_i)^{[\rm ev]} \text{ for } 1\le i \le n.
   \end{align*}
In particular ($m=2$), the weak principle holds for $(R,S,q)$ iff it holds for all $(R_i, S_i, q_i)$.
} 
\lv{
\item\label{dqd-weak} ({\em Reduction to $S$ of constant degree}) Any faithfully projective $S\in \Ralg$ induces a decomposition of $R$ as a direct product $R=R_1 \times \cdots \times R_n$ and correspondingly a decomposition $S=S_1 \times \cdots \times S_n$, $S_i = S\ot_R R_i$, such that $S_i$  is an $R_i$--algebra of constant degree. The $R$--algebra $S$ is \'etale if and only if every $S_i$ is an \'etale $R_i$--algebra, \ref{fea}\eqref{fea-d}.

    To prove the weak Knebusch norm principle we may by \eqref{dqd-ele-c} assume that $1\in \rmD(q)$. Then, using the notation of \eqref{dqd-ele-e} we have
 \[
\rmD(q)^{[\rm ev]} = \rmD(q_1)^{[\rm ev]} \times \cdots \times \rmD(q_n)^{[\rm ev]}
\]
and then for any $m\in \NN_+$,
\begin{align*}
  \rmN_{S/R}\big(\rmD(q_S)^{[m]}\big) &=  \rmN_{S_1/R_1}\big( \rmD(q_{1,S_1})^{[m]}\big)  \times \cdots \times \rmN_{S_n/R_n}\big(\rmD(q_{n, S_n})^{[m]}\big) \\
  &\subset \rmD(q)^{[\rm ev]} =
     \rmD(q_1)^{[\rm ev]} \times \cdots \times \rmD(q_n)^{[\rm ev]}
 \\ \iff &  \rmN_{S_i/R_i}\big(\rmD(q_{i, S_i})^{[m]}\big)  \subset
 \rmD(q_i)^{[\rm ev]} \text{ for } 1\le i \le n.
   \end{align*}
In particular ($m=2$), the weak principle holds for $(R,S,q)$ iff it holds for all $(R_i, S_i, q_i)$.}
\sm 

\item \label{dqd-ele-c} ({\em Reduction to $1\in \rmD(q)$}) {\em If the strong or weak Knebusch norm principle holds for $(q,S)$, then it also holds for $(uq, S)$, $u\in R\ti$.}  We can therefore assume $1\in \rmD(q)$ for proving the Knebusch norm principle.

\lv{
    Proof for the strong principle: $\rmD'(uq) = u \rmD'(q)$, hence $\rmN_{S/R}\big( \rmD'(u_Sq_S)\big) = \rmN_{S/R}\big( u_S \rmD'(q_S)\big) = u^d \rmN_{S/R}\big( \rmD'(q_S)\big) \subset u^d \rmD'(q)^{[d]} = (u \rmD'(q))^{[d]} = \rmD'(uq)^{[d]}$. %

Proof for the weak principle:
$\rmN_{S/R}\big( \rmD'(u_Sq_S)^{[2]}\big) = \rmN_{S/R}\big( u_S^2 \rmD'(q_S)^{[2]}\big) = \rmN_{S/R}(u_S)^2 \rmN_{S/R}\big( \rmD'(q_S)^{[2]}\big)
\subset \rmN_{S/R}(u_s)^2 \rmD'(q)^{[\rm ev]} = \rmD'(q)^{[\rm ev]} $ by \eqref{dqd-ele-a1}. 
}
\sm 

\item \label{dqd-ele-b} (Universal forms) {\em If $\rmD(q) = R\ti$, the strong and hence also the weak Knebusch principle holds. The same is true in case $\rmD\rmsm(q) = R\ti$. }

The assumption $\rmD(q) = R\ti$ implies $\rmD(q)^{[d]}= R\ti$ and then $\rmN_{S/R}\big(\rmD(q_S)\big) \subset \rmN_{S/R}(S\ti) \subset R\ti = \rmD(q)^{[d]}$. The same argument can be used for the smooth case. 

We note that $q$ is universal if it contains a hyperbolic plane. In this case, even $\rmD\rmsm(q) = R\ti$ by Example~\ref{lem_smooth_locus_exam}\eqref{lem_smooth_locus_hypn}. We also recall from \ref{isotrop}\eqref{isotrop-d} that $q$ contains a hyperbolic plane whenever $(M,q)$ contains an isotropic vector.
\sm 

\item\label{qdq-ele-aaa} ({\em Example $S\in \Ralg$ has constant rank $1$}) This case is equivalent to $S=R$.
Therefore $\rmN_{S/R} = \Id_R$, hence $\rmN_{S/R}\big(\rmD'(q)\big) = \rmD'(q)$ and the strong Knebusch norm principle holds in the general as well as in the smooth case. Observe that all the formulas in \eqref{thm-kneb-a1} are true, e.g., $\rmD'(q)\subset \rmD'(q)^{[2]}$ if $\uS\rmsm_{q,1}\ne \emptyset$ by \eqref{qdq-ele-aa}. \sm

\item\label{dqd-ele-d}  ({\em One-dimensional forms}) Let $(M,q) = (R, \lan u \ran_q)$ for some $u\in R\ti$, as defined in \ref{qf}\eqref{qfba-one}. 
Since $\rmD(q) = u R\ti{}^2$, it follows from \eqref{dqd-ele-c} and \eqref{dqd-ele-b} that the strong and weak Knebusch principles hold. 
\lv{
 Indeed, $\rmD(q_S) = u S\ti{}^2$ and so $\rmN_{S/R}(\rmD(q_S)) = \rmN_{S/R}(u S\ti{}^2) = u^d \rmN_{S/R}(S\ti)^2 \subset u^d R\ti{}^2 = u^{d-1}(u R\ti{}^2) \subset \Dqd$ because $\rmD(q) = u R\ti{}^2$ and therefore also $u \in \rmD(q)$.
 }
We claim that {\em this also holds for the smooth version
of the strong Kne\-busch principle.} This follows from
\begin{equation} \label{dqd-ele-d1}
    \rmD\rmsm(\lan u \ran_S) \ne \emptyset \iff 2\in R\ti.
\end{equation}
Proof of \eqref{dqd-ele-d1}: If $\rmD\rmsm(\lan u \ran_S) \ne \emptyset$, there exists $s\in \uS_{q,a}\rmsm(S)$ for some $a \in S\ti$. For such an $s$ the equation $2uss'=1$ has a solution in $S$, which says that $2\in S\ti$. In other words, the left multiplication $L_{2}$ is invertible on $S$. But $L_{2} = 1_S \ot L_{2}|_R$ so that $L_2|_R$ is invertible on $R$ by faithfully flat descent, i.e., $2\in R\ti$. Conversely, if $2\in R\ti$, then $(M,\lan u \ran)$ is regular, and so $\rmD\rmsm(q_S) = \rmD(q_S) = u S\ti{}^2$.

We recall that $(M,q) = (R, \lan u \ran)$ if $(M, q)$ is a quadratic $R$--space of $\rank_R M = 1$ and $\Pic(R) = \{0\}$, e.g., $R$ is unimodular or even an LG ring, \ref{qf}\eqref{qfba-one}. \sm

\item \label{dqd-ele-dd} ({\em Two-dimensional forms}) We claim: {\em any quadratic space $(M,q)$ with $M$ of constant rank $2$ satisfies the
    strong Knebusch norm principle for any finite projective $S\in \Ralg$}.

 Since $q$ is primitive, we can identify $M$ with the odd part $\Cli_1(q)$ of the Clifford algebra associated with $(M,q)$. 
 Its even part $\Cli_0(q)=A$ is a quadratic \'etale $R$--algebra; let $n_A$ be its norm. By \cite[V, (2.2.1)]{K} we then get $(M,q) \cong (A, un_A)$ for some $u\in \rmD(q)$. Hence, by \eqref{dqd-ele-c}, it suffices to establish the strong Knebusch norm principle for $(A, n_A)$. After identifying $(M,q)$ and $(A, n_A)$ we have $\rmD(q) = \rmN_{A/R}(A\ti)$, a subgroup of $R\ti$. The base change $A_S= A\ot_R S$ is a quadratic \'etale $S$--algebra whose norm is the base change of $n_A$.
 Hence $\rmD(q_S) = \rmN_{A\ot S/ S}(A_S\ti)$ and, by transitivity and multiplicativity of norms \ref{trno},
 \begin{align*}
   \rmN_{S/R}\big( \rmD(q_S) \big)
   & = \rmN_{S/R}\big(\rmN_{A\ot S/S}((A_S\ti)\big) \\
   & = \rmN_{A\ot S/R} (A_S\ti)
   = \rmN_{A/R}\big( \rmN_{A\ot S/A}(A\ti_S)\big) \\
  &\subset \rmN_{A/R}(A\ti) = \rmD(q) = \Dqd.
 \end{align*}%

\item \label{dqd-ele-f} ({\em Rank reductions for proving the Knebusch norm principles, $R$ LG})
Let $R$ be an LG ring and let $(M,q)$ be a quadratic $R$--space. We then know (\ref{LGqdi}) that $\rmD(q) \ne \emptyset$, so that the observations above apply.

In view of the reduction  
\eqref{dqd-ele-e}, the verification of the Knebusch norm principles is reduced to that for $(M,q)$ and $S$ with $M$ and $S$ projective of constant rank $r\in \NN_+$ and $d\in \NN_+$ respectively. 
We know that any finite projective $R$--module of constant rank is free. Hence, by 
\eqref{qdq-ele-aaa}, \eqref{dqd-ele-d} and  \eqref{dqd-ele-dd} we can assume that 
$r>2$ and $d>1$ for proving the strong Knebusch norm principle. \sm 

\item\label{dqd-ele-h} ({\em Reduction to cosets $aS\ti$}) Let $R$, $(M,q)$ and $S$ be arbitrary, as in the beginning of this subsection, let $a\in S\ti$ such that $\rmN_{S/R}(a) \in \rmD'(q)^{[j]}$ for some $j\in \NN_+$, and let $\rmD'' \in \{\rmD, \rmD\rmsm \}$. Then 
\begin{equation}\label{dqd-ele-h1} 
 \rmN_{S/R}(\rmD''\big(q_S) \cap a S\ti{}^2 \big) \subset \rmD'(q)^{[j]}. 
\end{equation}
Indeed, we can assume $\rmD''\big(q_S) \cap a S\ti{}^2 \ne \emptyset$. Let $b\in  S\ti$ and $c= a b^2 \in \rmD''(q_S)$. Then $a = b^{-2} c \in S\ti{}^2 \, \rmD''(q_S) = \rmD''(q_S)$ by \eqref{dqd-ele-a1}. Hence $\rmN_{S/R}(c) = \rmN_{S/R}(a) \rmN_{S/R}(b)^2 \in \rmD'(q)^{[j]}\, R\ti{}^2 =  \rmD'(q)^{[j]}$, again by \eqref{dqd-ele-a1}. 
\end{inparaenum}
\sm 
We describe one further reduction. 

\comments{(2026-05-05) Previously, Lemma~\ref{redsm} was part of the old version of Proposition~\ref{prop_knebusch}. I put it here because it does not need the assumptions of that proposition, and also so that is clear what it is dependance on transitivity of $\SO(q)$ or $\orth(q)$, namely the assumption~\ref{redsm}\ref{redsm-i}. In the previous version this was  Proposition~\ref{prop_quadric}\ref{lem_sphere2_iii}, which is no longer available. It stated: \sm 

{\tt $R$ semilocal, $\rank M \ge 2$, $v\in \uS_{q,a}\rmsm(R)$
 and $m\in M$ unimodular. Then there exists $g\in \uSO(q)(R)$ such that $Rv \oplus Rg(m)m$ is a direct summand of $M$
} \sm 

We applied this to $u\in \uSO_{q_S, 1}\rmsm(S)$ which was non-empty by assumption
}

\begin{lem}[Reduction to smooth spheres]\label{redsm} Let $(M,q)$ be a quadratic $R$--space, and let $S\in \Ralg$ free of rank $d\in \NN_+$ and unit-generated by $a\in \rmD(q_S)$. Suppose \sm 

\begin{enumerate}[label={\rm (\roman*)}]
 \item\label{redsm-i} there exist $u,v \in M_S$ such that $Su + Sv$ is a complemented submodule of $M_S$ which is free of rank $2$ with basis $(u,v)$ satisfying $q_S(u) =  a$ and $q_S(v) = 1$. 
\end{enumerate}
Let $S' = S[X]/(X^2 - b_{q_S}(u,v) X + a)\in \Salg$ and let $a' \in S'$ be the image of $X \in S[X]$. Then, as $R$--algebra,  $S'$ is unit-generated by $a'$, free of rank $2d$ and 
\[ a'\in \rmD\rmsm(q_{S'}) \quad \text{with} \quad 
  \rmN_{S'/R}(a') = \rmN_{S/R}(a). 
\]
\end{lem}

\begin{proof} Since $q_S(u)=a$ and $q_S(v)=1$, the restriction of the quadratic form $q_S$ to this submodule has the form $q(xu+yv)= ax^2+ c xy+ y^2$ for $c = b_{q_S}(u,v) \in S$. The free quadratic $S$--algebra $S'=S[X]/(X^2- c X+ a)$ is generated by $a'$ and, by \eqref{trno-00}, its norm is  $\rmN_{S'/S}(a')=a \in S\ti$. Hence, $a'$ is invertible and 
\begin{equation}\label{prop_knebusch-1}
  \rmN_{S'/R}(a') = \rmN_{S/R}\big( \rmN_{S'/S}(a')\big) = \rmN_{S/R}(a)
\end{equation}
by transitivity of norms \ref{trno}\eqref{trno-b}. The $R$--algebra $S'$ is free of rank $2d\ge 2$. It is unit-generated by $a'$ since $a=a' (c-a')$ 
and $a$ is a generator of the $R$--algebra $S$. The vector $u + a'v \in M_{S'}$ is unimodular.
Because $q( u+ a'v)= a  - c + (a')^2 =0$, it follows that $q$ is $S'$-isotropic. By Example~\ref{lem_smooth_locus_exam}\eqref{lem_smooth_locus_hypn}, $\uS_{q_S,  a'}\rmsm(S') \ne \emptyset$.    
\end{proof}

\ms

In some favourable situations we have $S\ti = S\ti{}^2 \Prim_R(S)\ti$, see for example Lemma~\ref{hunco} or the related Proposition~\ref{prop_tiroir}. In this case, the formula \eqref{dqd-ele-h1} says we can concentrate on $\rmN_{S/R}\big( \Prim_R(S)\ti)$ when investigating versions of the Knebusch norm principle. This justifies the emphasis on unit-generators in the following. Proceeding by induction on the rank $d$ of $S\in \Ralg$, the following Lemma~\ref{knebp} describes how to pass from a unit-generated $S$ of rank $d$ to one of rank $d-1$.

\begin{lem}\label{knebp} Let $(S,a,q)$ be data as in the preliminary setting {\rm \ref{aps}}: $S\in \Ralg$ is unit-generated by $a\in S$ and free of rank $d$ as $R$--module, and $(M,q)$ is a quadratic $R$--module. Furthermore, we fix $T\in \Ralg$, put $a_T = a \ot_R 1_T\in S\ot_R T$ and use the Weil restriction of $S$--schemes \eqref{aps-0}:
\begin{align*} 
 \frR_{S/R}\big( \uW_S(M\ot_R S)\big)\,(T) 
   &= \uW_R\big(\frR_{S/R}(M\ot_R S)\big)\,(T) 
 \\ = \uW_R(M\ot_R S) \, (T) & = M \ot_R S \ot_R T = M_T \ot_T (S\ot_R T), 
\end{align*} 
see {\rm \ref{aps}} for notation, and analogously for $S$--spheres
\begin{align*}
  \frR_{S/R}(\uS_{q_S, a\me}) \,(T) 
 = \{ m \in M_T \ot_T (T\ot_R S) : q_{S\ot T} (m) = a_T \me \}.
\end{align*} 
Since $\uW_S(M_S)$ and $\uS_{q_S, a\me}$ are finitely presented affine schemes, {\rm \ref{q-Faser}\eqref{q-Faserb}}, both Weil restrictions exist as $R$--schemes, {\rm \ref{weilres}\eqref{weilres-a}}. We suppose that we have a vector 
\begin{equation}\label{knebp-00}  
w=\textstyle \sum_{i=0}^{d-1} m_i a_T^i\in \frR_{S/R}( \uS_{q_S, a\me}) \,(T), 
 \quad m_i \in M_T.
\end{equation}
We associate to w the vector 
\[ w(X)=\textstyle  \sum_{i=0}^{d-1} m_i X^i\in M_T \ot_T T[X] \]
so that $w(a_T) = w$. Then the conditions \ref{knebp-i} and \ref{knebp-ii} below are equivalent:

\begin{enumerate}[label={\rm (\roman*)}]

\item \label{knebp-i} $q_T(m_{d-1}) \in T\ti$,  

\item\label{knebp-ii} there exist a monic $Q_w\in T[X]$ of degree $d-1$ and $c_w\in T\ti$ such that
\begin{equation}\label{knebp-0} 
    X \, q_{T[X]}\big(w(X)\big)= 1+ c_w \,  P_T(X) \,  Q_w(X)
\end{equation} 
    where $P_T$ is the $T$--extension of the characteristic polynomial $P$ of $a$.%
    \sm 
\end{enumerate}
If \ref{knebp-i} and \ref{knebp-ii} hold, then: 
 \begin{enumerate}[label={\rm (\alph*)}]
 \item\label{knebp-a} $Q_w$ is the unique polynomial satisfying \eqref{knebp-0}. \sm

 \item\label{knebp-b} $c_w=q_T(m_{d-1})$ for $c_w$ as in \ref{knebp-i}. \sm 

 \item \label{knebp-d} $a_w= X +(Q_w) \in S_w=T[X]/Q_w(X)$ is an invertible primitive element of the one-generated free $T$--algebra $S_w$ of rank $d-1$, and  $w(a_w) \in M_T \ot_T S_w = M \ot_R S_w$ satisfies $q_{S_w}\big( w(a_w)\big) = a_w\me$. i.e.,  
\begin{equation}\label{knebp-d1}     
     w(a_w) \in \frR_{S_w/R}(\uS_{q_{S_w}, a_w\me})\, (T). 
 \end{equation}   

 \item \label{knebp-c} $\rmN_{S\ot T/T}(a_T) = \big(c_w\, \rmN_{S_w/T}(a_w)\big)\me$, in particular, \sm 

 \begin{enumerate}[label={\rm (\Roman*)}]
    \item\label{knebp-c1}    
     if $\rmN_{S_w/T}(a_w) \in \rmD(q_T)^{[d-1]}$, then $\rmN_{S\ot T/T}(a_T) \in \rmD(q_T)^{[d]}$, and 
     
    \item \label{knebp-c11} if $\rmN_{S_w/T}(a_w) \in \rmD\rmsm(q_T)^{[d-1]}$ and  $m_{d-1} \in \uS\rmsm_{q_T, c_w}(T)$,  then \allowbreak $\rmN_{S\ot T/T}(a_T) \in \rmD\rmsm(q_T)^{[d]}$. 
\end{enumerate} \sm 

\item \label{knebp-e} Suppose there exists $v(X)\in M_T \ot_T T[X]$ such that the polar form $b_{T[X]}$ of the quadratic form $q_{T[X]} \co M_T\ot_T T[X] \to T[X]$ satisfies
\[ b_{T[X]} \big(w(X), \, v(X)\big) \in T[X]\ti \]
  then $w\in \frR_{S/R}( \uS_{q_S, a\me}\rmsm)\,(T)$ and $w(a_w) \in \frR_{S_w/R}(\uS\rmsm_{q_{S_w}, a_w\me})\, (T)$.
\end{enumerate}
\end{lem}

\begin{proof} For the given $w(X)$ we have
\begin{equation}\label{knebp1}
 Xq_{T[X]}\big(w(X)\big) = q_T(m_{d-1})X^{2d-1} + \cdots (\text{lower terms in $X$}).
\end{equation}
If \ref{knebp-i} holds, then
\[ Xq_{T[X]}\big( w(X)\big) = c_w X^{2d-1} + \cdots (\text{lower terms in $X$})\]
with $c_w\in T\ti$. A comparison then shows $c_w=q_T(m_{d-1}) \in T\ti$. Conversely, assuming \ref{knebp-ii}, the polynomial $P_T$ divides $Xq_{T[X]}\big(w(X)\big) - 1$ because $a_Tq(w(a_T))= 1$, say $X q_{T[X]}\big(w(X)\big) - 1 = P_T(X) \wtl Q(X)$ for a unique $\wtl Q\in T[X]$ of degree $d-1$, the latter condition being obtained by comparison with \eqref{knebp1}. Moreover, the leading term $\wtl c$ of $\wtl Q$ satisfies $\wtl c=q_T(m_{d-1}) \in T\ti$, so that \eqref{knebp-0} 
holds  with $Q_w=c_w\me \wtl Q$. 

We have seen \ref{knebp-a} and \ref{knebp-b} in the proof above. Regarding \ref{knebp-d}, the $T$--algebra $S_w$ is clearly one-generated by $a_w$ and free of rank $d-1$. Specializing $X$ to $a_w$ in \eqref{knebp-0} yields 
\[ a_w \, q_{S_w}\big( w(a_w)\big) = 1 \]
since $Q_w(a_w) = 0$. This proves the remaining parts of \ref{knebp-d}. 

For the proof of \ref{knebp-c} note that $\rmN_{S\ot T/T}(a_T) = (-1)^d P_T(0)$ and $ \rmN_{S_w/T}(a_w)=(-1)^{d-1} Q_w(0)$ by \ref{trno}\eqref{trno-bou}. Evaluating \eqref{knebp-0} at $0$ yields $0 = 1 + c_w Q_w(0) P_T(0)$, and therefore
\begin{align*}
\rmN_{S\ot T/T}(a_T) &= (-1)^d P_T(0) = (-1)^d (-c_w\me Q_w(0)\me)
\\ &= c_w\me (-1)^{d-1} Q_w(0)\me = c_w\me \rmN_{S_w/T}(t)\me.
\end{align*} 
The final claim is then immediate using \ref{knebp-b} and $\rmD(q)=\rmD(q)\me$ by \eqref{dqd-ele-a3}. 

\ref{knebp-e} The evaluation map $T[X] \to T[X]/(P_T) = S\ot_R T$, $X \mapsto a_T$, is a $T$--algebra homomorphism and sends  $b_{T[X]} \big(w(X), \, v(X)\big)$ to $b_{S\ot T}\big(w, v(a)\big)\in (S\ot_R T)\me$ so that  $w\in \frR_{S/R}( \uS_{q_S, a\me}\rmsm)\,(T)$. The proof for $w(a_w)$ is analogous. 
\end{proof}
\sm

\textbf{Remarks.} Using the terminology of \ref{bate}\eqref{bate-c}, the formula \eqref{knebp-c1} says that if the strong Knebusch norm principle holds for $a_w$, then it also holds for $a$. This indicates that primitive elements are amenable to prove the strong Knebusch norm principle by induction. We will use this for example in the proof of Proposition~\ref{prop_knebusch}.  Of course, the point here is that it is a priori not clear that a $w$ satisfying the conditions of \ref{knebp} actually exists.

\subsection{Some preliminary $R$--functors and schemes.} \label{motkn}
The constructions of Lemma~\ref{knebp} are functorial and define various $R$--functors and schemes representing them. For $(S,a,q)$ as in Lemma~\ref{knebp} and $T\in \Ralg$ we consider elements
\begin{equation}\label{motkn-1}
 w(X) = \textstyle \sum_{i=0}^{d-1} m_i X^i \in M_T \ot_T T[X]
\end{equation}  
and the conditions 
\begin{enumerate}[label={\rm (Kn\roman*)}]
\item \label{motkni} \qquad $q_T(m_{d-1}) \in T\ti$, \sm 

\item \label{motknii} \qquad $w(a_T) \in \frR_{S/R}\big(\uS_{q_S, a\me}\big)\, (T)$, \sm 

\item \label{motkniii} \qquad $w(a_w) \in \frR_{S_w/R}\big(\uS_{q_{S_w}, a_w \me } \big)\, (S_w)$, \sm 
\end{enumerate}
as well as the smooth versions of \ref{motknii} and \ref{motkniii}: 
\begin{enumerate}[label={\rm (Kn\roman*)$_{\rm sm}$}] \setcounter{enumi}{1}
\item \label{motkniism} \qquad $w(a_T) \in \frR_{S/R}\big(\uS_{q_S, a\me}\rmsm \big)\, (T)$, \sm 

\item \label{motkniiism} \qquad $w(a_w) \in \frR_{S_w/R}\big(\uS\rmsm_{q_{S_w}, a_w \me } \big)\, (S_w)$. \sm 
\end{enumerate}
The notation used above is that of Lemma~\ref{knebp}: for any $w(X)$ as in \eqref{motkn-1} satisfying \ref{motkni} and \ref{motknii} we have a monic $Q_w \in T[X]$ and $a_w = X + (Q_w) \in S_w = T[X]/(Q_w)$. We then know that \ref{motkniii} holds, but not necessarily the smooth versions of these. Regarding \ref{motkniism}, one can proceed as follows.   
\sm 

As both $\uS_{q_S, a\me}$ and $\uW(M_S)$ are affine $S$--schemes of finite type, their Weil restrictions $\frR_{S/R}(\cdot)$ are representable by $R$--schemes, \ref{weilres}. By \cite[7.7/2(ii)]{BLR}, we have a closed immersion 
\[ \wtl X = \frR_{S/R}\big( \uS_{q_S, a\me} \big) \to
   \frR_{S/R} \big( \uW_S(M_S)\big)
\]
of $R$--schemes. Since the $S$--scheme $\uS_{q_S, a\me}$ is affine by \ref{lem_smooth_locus-LG}\ref{lem_smooth_locus-LGb}, the $R$--scheme 
$\wtl X$ is affine by \ref{weilres}\eqref{weilres-b}, but in general not smooth. 

Let $f$ be the restriction of the polynomial 
\[ \frR_{S/R} \big( \uW_S(M_S)\big) \xrightarrow{\; \bfq \circ \pi_{d-1}\;}
   \GG_a
\]
of \eqref{aps-2} to $\wtl X$. Then 
\[ \wtl X_f = \big( \frR_{S/R}\big( \uS_{q_S, a\me} \big)\big)_f
\]  
represents the $R$--functor whose $T$--points satisfy \ref{motkni}, \ref{motknii} and hence also \ref{motkniii}. 

By \ref{weilres}\eqref{weilres-b} and \ref{weilres}\eqref{weilres-d} we have an open immersion of $R$--schemes 
\[ 
    X = \frR_{S/R} \big( \uS_{q_S,a\me}\rmsm \big) \longto
     \wtl X = \frR_{S/R}\big( \uS_{q_S, a\me} \big)
\] 
with $X$ being a smooth, but not necessarily affine $R$--scheme. We note that $\frR_{S/R} \big( \uS_{q_S,a\me}\rmsm \big)$ is indeed representable by an $R$--scheme since 
$\uS_{q_S,a\me}\rmsm$ is a quasi-affine $R$--scheme of finite presentation, \ref{lem_smooth_locus-LG}\ref{lem_smooth_locus-LGe}. We define the $R$--scheme $X_f$ by the cartesian diagram 
\begin{equation} \label{motkn-1} \vcenter{ 
\xymatrix@C=50pt{\ar @{} [dr] |{\small\qed} %
X_f \ar[r] \ar[d] & \wtl X_f \ar[d] 
 \\ X \ar[r] & \wtl X } } . 
\end{equation}
We view $X_f$ as an open subscheme of both $X$ and $\wtl X _f$. It represents the $R$--functor whose $T$--points satisfy \ref{motkni}, \ref{motkniism} and \ref{motkniii}, and is a smooth scheme since $X$ is smooth. We will deal with \ref{motkniiism} in the following Lemma~\ref{knsch}. 

\begin{lem}\label{knsch} There exists an $R$--scheme $Z$ and a quasi-compact open immersion $Z \to \wtl X_f$ such that $Z$ represents the $R$--functor, whose  $T$--points, $T\in \Ralg$, are defined by the conditions \ref{motkni}, \ref{motknii} and \ref{motkniiism}. 
\end{lem}

\begin{proof} The problem with condition \ref{motkniiism} is that it depends $w(X)$. As in the proof of the representability of the Springer functor in \ref{respr} we will retreat to universal points. 

We know that the $R$--scheme $\wtl X_f$ of \ref{motkn} is affine, say 
$\wtl X_f = \Spec(C)$ for some $C\in \Ralg$. We apply the constructions of Lemma~\ref{knebp} to the universal point $w_\sharp \in \wtl X_f(C)$. We get $c_\sharp \in C\ti$ and a monic polynomial $Q_\sharp \in C[X]$ of degree $d-1$ such that 
$ q_{C[X]}\big( w_\sharp(X)\big) = c_\sharp \, P_C(X) \, Q_\sharp(x) \, \in C[X]$.
We put $a_w = X + (Q_\sharp) \in C_\sharp = C[X]/Q_\sharp$. The $C_\sharp$--scheme 
\[ S_\sharp : = \uS_{q_{C_\sharp}, a_\sharp\me}
\] 
is affine. Since $C_\sharp$ is a finite projective $C$--algebra, the Weil restriction $\frR_{C_\sharp/C}(S_\sharp)$ is represented by an affine $C$--scheme. Moreover, the open immersion 
\[ S_\sharp\rmsm :=  \uS_{q_{C_\sharp}, a_\sharp\me}\rmsm \to S_\sharp \]
gives rise to an open immersion $\frR_{C_\sharp/C}(S_\sharp\rmsm) \to \frR_{C_\sharp/C}(S_\sharp)$ of $C$--schemes. By Yoneda, 
\[ \Mor\big(\Spec(C), \frR_{C_\sharp/C}(S_\sharp) = \frR_{C_\sharp/C}(S_\sharp)\, (C) 
  = S_\sharp(C_\sharp).  
\]
Hence, we have a morphism of $C$--schemes 
$ 
s\co \Spec(C) \to \frR_{C_\sharp/C}(S_\sharp)
$
corresponding to $w_\sharp(X + Q_\sharp) \in S_\sharp(C_\sharp)$, which is a section of the structure morphism $\frR_{C_\sharp/C}(S_\sharp) \to \Spec(C)$. Finally, we define a $C$--scheme $Z$ by the fibre product 
\begin{equation*}  \vcenter{ 
\xymatrix@C=60pt{\ar @{} [dr] |{\small\qed} %
Z \ar[r]^{i'} \ar[d] & \Spec(C)= \wtl X_f \ar[d]    \\ 
 \frR_{C_\sharp/C}(S_\sharp\rmsm) \ar[r]^i & \frR_{C_{\sharp}/C}(S_\sharp)} } \quad . 
\end{equation*}
We will consider $Z$ as an $R$--scheme by composing $Z\to \Spec(C)$ with the structure morphism $\Spec(C) \to \Spec(R)$ of the $R$--algebra $C$. 

We claim that the $R$--scheme $Z$ is a scheme whose existence is claimed in the lemma. This can be proven in the same way as the proof of Proposition~\ref{respr}, replacing \eqref{respr01a} there by 
\begin{align*}
 \frR_{C_\sharp/C}(S_\sharp) \times _C T &\cong 
 \frR_{(C_\sharp \ot_C ^{w^\flat} T)/T}(S_\sharp \times_C T)
   \qquad(\text{by \ref{weilres-bc1}})
\\
&\cong \frR_{S_w/T}(\uS_{q_{S_w}, a_w\me})
\end{align*}
by $C_\sharp \ot_C ^{w^\flat} T \cong S_w$ and \eqref{q-Faserb1}. Also 
\eqref{respr01b} has to be replaced by 
\begin{align*}
 \frR_{C_\sharp/C}(S_\sharp\rmsm) \times _C T &\cong 
 \frR_{(C_\sharp \ot_C ^{w^\flat} T)/T}(S_\sharp\rmsm \times_C T)
   \qquad(\text{by \ref{weilres-bc1}})
\\
& \cong \frR_{S_w/T}(S_\sharp\rmsm \times_C T) 
 \cong  \frR_{S_w/T}(\uS\rmsm_{q_{S_w}, a_w\me})
\end{align*}
where $S_\sharp\rmsm \times_C T \cong \uS\rmsm_{q_{S_w}, a_w\me}$ follows from \eqref{q-Faserb1} again. We leave the other details to the reader. 

The morphism $S_\sharp\rmsm \to S_\sharp$ is a quasi-compact open immersion by \ref{lem_smooth_locus-LG}\ref{lem_smooth_locus-LGe}, hence, by \ref{weilres}\eqref{weilres-d}, so is the morphism $i$ in the diagram above. But then also $i'$ is a quasi-compact open immersion. \end{proof}

\sm 

In view of \ref{motkn} and \ref{knsch}, the following definition of the Knebusch functor $\ul{\Kneb}$ should not be surprising to the reader.


\subsection{Knebusch functors, Knebusch data}\label{knf}
Let $R$ be an arbitrary base ring and let $(S,a,q)$ be data as in Lemma~\ref{knebp}. 
By definition, the {\em Knebusch functor $\ul{\Kneb}_{(S,a,q)}$ associated with $(S,a,q)$\/} is the $R$--functor whose $T$--points, $T\in\Ralg$, is the set $\ul{\Kneb}_{(S,a,q)}\, (T)$ consisting of 
\begin{equation}\label{knf-1}
 w(X) = \textstyle \sum_{i=0}^{d-1} m_i X^i \in M_T \ot_T T[X]
\end{equation}  
satisfying the conditions \ref{motkni}, \ref{motkniism} and \ref{motkniiism} of \ref{motkn}.  

By \ref{knebp}\ref{knebp-b} the vector $m_{d-1} \in M_T$ lies in the sphere $\rmS_{q_T, c_w} = \uS_{q_T, c_w}(T)$, but not necessarily in the smooth part of this sphere. We therefore define the {\em special Knebusch functor $\ul{\Kneb}^+_{(S,a,q)}$ associated with $(S,a,q)$} as the subfunctor of $\ul{\Kneb}_{(S,a,q)}$ whose $T$--points are those $w\in \ul{\Kneb}_{(S,a,q)}(T)$
satisfying 
\begin{enumerate}[label={\rm (Kn\roman*)}]\setcounter{enumi}{3}
\item \label{motkniv} \qquad $m_{d-1} \in \uS_{q_T, c_w}\rmsm (T)$. 
\end{enumerate}

\comments{(2026-04-15) It seems that the old concept of Knebusch is still useful}

We will refer to an element of $\ul{\Kneb}_{(S,a,q)}(T)$ as a {\em Knebusch datum (over $T$)} and to those in  $\ul{\Kneb}_{(S,a,q)}^+(T)$ as a {\em special Knebusch datum (over $T$)}. We immediately make a few elementary observations.
\sm 

\begin{inparaenum}[(a)] 

\item\label{knf-b} ({\em $q$ regular}) Suppose $q$ is regular. Then spheres are smooth by \ref{lem_smooth_locus-LG}\ref{lem_smooth_locus-LGd}. Hence, \ref{motkniism} = \ref{motknii} and \ref{motkniiism} = \ref{motkniii}. Therefore, as pointed out in \ref{motkn},
\begin{equation}\label{knf-b1} 
\begin{split}&\text{\em if $q$ is regular, then the $R$--scheme 
   $\wtl X_f = X_f = \big(  \frR_{S/R} \big( \uS_{q_S,a\me}\big)\big)_f$}
\\ &\text{\em represents the $R$--functor $\ul{\Kneb}_{(S,a,q)}$. } 
\end{split}\end{equation}    
Moreover, all Knebusch data are special, so $X_f$ also represents the $\ul{\Kneb}^+_{(S,a,q)}$. 
We will address representability of $\ul{\Kneb}_{(S,a,q)}$  in Proposition~\ref{lem_representability}. \sm
 
\item\label{knf-c} ({\em Direct sums}) Suppose that $M=M'\oplus M''$ is a direct sum of $R$--modules, not necessarily orthogonal with respect to $q$. We put $q' = q|_{M'}$, and recall that for any $T\in \Ralg$ the $T$--modules $M'_T$ and $M''_T$ can be identified with submodules of $M_T$ such that $M_T = M'_T \oplus M''_T$. In particular, we can identify $M'_T\ot_T T[X]$ with a $T$--submodule of $M_T \oplus_T T[X]$.  Thus, we can view any $w'(X) \in M'_T \ot_T T[X]$ in the form \eqref{knf-1} also as an element of $M_T \ot_T T[X]$, denoted $\veps_T\big( w'(X)\big)$.  
    
    By construction, restriction of quadratic forms commutes with base change: $(q')_T = (q_T)'$, and this also holds for the associated polar forms. Thus, $w'(X)$ satisfies \ref{motkni} for $q'_T$ if and only if $\veps_T\big(w'(X)\big)$ satisfies \ref{motkni} for $q_T$. Moreover, the description of the smooth locus of spheres in \ref{lem_smooth_locus-LG}\ref{lem_smooth_locus-LGc}  shows that if $w'(X)$ satisfies \ref{motkniism} and \ref{motkniiism}, then $\veps_T\big( w'(X)\big)$ satisfies the analogous conditions. Altogether, we have an injective set map
    \[ \veps_T \co \ul{\Kneb}_{(S,a,q')} (T) \to \ul{\Kneb}_{(S,a,q)}(T). \]     
    Analogously, we also have an injective set map
    \[ \veps_T^+ \co \ul{\Kneb}^+_{(S,a,q')} (T) \to \ul{\Kneb}^+_{(S,a,q)}(T). \]     
    We leave it to the reader to verify that the family $(\veps_T)_{T\in \alg}$ constitutes a natural transformation between $R$--functors
    $ \veps \co \ul{\Kneb}_{(S,a,q')} \to \ul{\Kneb}_{(S,a,q)}$    
and that the analogous statement holds for special Knebusch functors. \end{inparaenum}

\begin{lem} \label{knevid} Let $S\in \Ralg$ be unit-generated by $a\in S\ti$ and free of rank $d\ge 2$, and let $(M,q)$ be a quadratic $R$--module. 
Assume that 
\begin{enumerate}[label={\rm  ($\star$)}]
\item \label{knevidstar}  the $R$--module $M$ is a direct, but not necessarily orthogonal sum $M=M'\oplus M''$ where $M'$ is free of rank $3$ with $R$--basis $\{e,f,g\}$ for which $q|_{M'}$ is given by \[ q(xe + yf + z g) = xy + cz^2\] for some $c\in R\ti$ and $x,y,z\in R$.
\end{enumerate}
Then $\Kneb(S,a,q) := \ul{\Kneb}_{(S,a,q)}(R) \ne \emptyset$. 
\sm

In particular, assume that $R$ is an LG ring and $(M,q)$ is an isotropic $R$--{\em space}. Then

\begin{enumerate}[label={\rm (\alph*)}]

\item\label{knevid-a} $\Kneb(S,a,q) \ne \emptyset$ if $M$ has rank $\ge 3$, and
\sm

\item\label{knevid-b} $\Kneb^+(S,a,q) \ne \emptyset$ if $M$ has rank $\ge 4$.
\end{enumerate} \end{lem}

\begin{proof} By \ref{knf}\eqref{knf-c} 
we can assume $M=M'$. Let $w\in M_S$ be the vector whose coefficients with respect to the $S$--basis $e,f,g$ of $M_S$ are $(1,w_2, a^{d-1})$ where 
$w_2\in S$ is determined by the equation $a\me = q_S(w) = w_2 + c a^{2d-2}$. 
Writing $w$ in the form $w = \sum_{i=0}^{d-1} m_i a^i$ with unique $m_i \in M$ and comparing the $a^{d-1}$ coefficient of $w$, we find $m_{d-1} = (0,b,1)$ for some $b\in R$. This implies $q(m_{d-1}) = c \in R\ti$.

We now lift $w$ to $w(X) = \sum_{i=1}^{d-1} m_i X^i\in M\ot_R R[X]$ and thus have $w(a) = w$. In the $R[X]$--basis $(e,f,g)$ of $M\ot_R R[X]$ we find $w(X) = (1, g(X), X^{d-1})$ for a unique $g\in R[X]$ of degree $\le d-1$. For $v(X) = (0,1, 0) \in M\ot_R R[X]$ we have 
$b_{R[X]} \big(w(X), v(X) \big) = 1$, and therefore $w \in \frR_{S/R}\big(\uS_{q,a\me}\rmsm\big)(R)$ and $w(a) = (1, g(a_w), a_w^{d-1}) \in  \frR_{S_w/R}(\uS\rmsm_{q_{S_w}, a_w\me})\, (R)$ by Lemma~\ref{knebp}\ref{knebp-e}.  
\sm

\ref{knevid-a} We will show that $(M,q)$ satisfies the assumptions of the general case \ref{knevidstar} treated above. Since by \ref{revLG}\eqref{revLG-a} the ring $R$ is unimodular, \ref{isotrop}\eqref{isotrop-unifap} shows that $(M,q)$ contains an isotropic vector. Hence, by \ref{isotrop}\eqref{isotrop-d}, the quadratic module $(M,q)$ contains a hyperbolic plane $\HH$. Since $q|_{\HH}$ is regular, $M$ has an orthogonal decomposition $M = \HH \oplus \HH^\perp$. By \eqref{qf-perp0}, the quadratic form $q|_{\HH^\perp}$ is nonsingular, and by \ref{LGqdi}
the quadratic space $(\HH^\perp, q|_{\HH^\perp})$ contains $m$ with $q(m) \in R\ti$.
Such an $m$ is unimodular by \ref{mx_lem}\eqref{mx_lema}, hence $Rm$ is free of rank $1$ with basis $\{m\}$ and a direct summand of $\HH^\perp$. Thus, $(M,q)$ satisfies the assumptions of the first part of the lemma.
\sm

\ref{knevid-b} As in the proof of \ref{knevid-a} we can split off a hyperbolic plane $\HH$ and then know that $q|_{\HH^\perp}$ is nonsingular of rank $\ge 2$. Applying the rank decomposition to the quadratic space $(\HH^\perp, q_{\HH^\perp})$ and then the classification of nonsingular spaces over LG rings in \ref{nqf-LG}, we see that $\HH^\perp$ contains a submodule $N$ which is free of rank $2$ and such that $q|_N$ is regular. Thus, we can decompose $\HH^\perp = N\perp M''$. 
Putting $M'=\HH \oplus N$, we get an orthogonal decomposition $M = M'\oplus M''$ with $q'=q|_{M'}$ being regular. Hence $\Kneb^+(S,a,q') = \Kneb(S,a,q') \ne \emptyset$ by \ref{knevid-a}. The claim then follows from $\Kneb^+(S,a,q') = \Kneb^+(S,a,q)$ by \ref{knf}\eqref{knf-c}. 
\end{proof}


\begin{prop}[Representability of $\ul{\Kneb}$] \label{lem_representability} As in {\rm \ref{knebp}--\ref{knevid}}, let $S$ be a finite free $R$-algebra of rank $d \geq 2$ which is unit-generated by $a \in \Prim_R(S)\ti$ and let $(M,q)$ be a quadratic $R$--module. \sm 

\begin{inparaenum}[\rm (a)] \item \label{lem_representability-r} Then the Knebusch functor $\ul{\Kneb}_{(S,a,q)}$ is representable by a quasi-compact open subscheme 
\[ \uKneb_{(S,a,q)} \] 
of the smooth $R$--scheme $\frR_{S/R}(\uS_{q,a\me}\rmsm)$, while the special Knebusch functor 
$\ul{\Kneb}^+_{(S,a,q)}$ is representable by the open subscheme
\[ \uKneb^+_{(S,a,q)} = \uKneb_{(S,a,q)} \cap \Sur_{b_q}\]
of $\uKneb_{(S,a,q)}$. \sm

\item \label{lem_representability-a} Let $(M,q)$ be a quadratic space
of rank $\ge 3$. Then $\uKneb(S,a,q)$ has geometrically integral fibres and is universally schematically dense in $\frR_{S/R}\bigl(\uS_{q,a^{-1}}\rmsm\bigr)$. \sm 
\sm

\item \label{lem_representability-b} Let $(M,q)$ be a quadratic space of rank $\ge 4$. Then  $\uKneb^+(S,a,q)$ is a smooth $R$--scheme, has geometrically connected fibres and is universally schematically dense in $\frR_{S/R}\bigl(\uS_{q,a^{-1}}^{\rm sm}\bigr)$.
\end{inparaenum}
\end{prop}

\begin{proof} \eqref{lem_representability-r} We define $\uKneb = \uKneb_{(S,a,q)}$ as the intersection of the open  subschemes  $X_f\to \wtl X_f$ and $Z\to \wtl X_f$ of \eqref{motkn-1} and Lemma~\ref{knsch} respectively:  
\[ 
\vcenter{ 
\xymatrix@C=50pt{\ar @{} [dr] |{\small\qed} %
\uKneb \ar[r] \ar[d] & Z\ar[d] \\ X_f \ar[r] & \wtl X_f 
    } } . 
\] 
We have seen in \ref{motkn} that $X_f$ is an open subscheme of the smooth scheme $X=\frR_{S/R}(\uS\rmsm_{q_S, a\me})$. 

Regarding the special Knebusch functor $\ul{\Kneb}^+$, we note that Lemma~\ref{knebp} shows that $m_{d-1} \in  \uS_{q_T, c_w}(T)$ for $T\in \Ralg$ and $w(X)\in \ul{\Kneb}$. Since $\uS\rmsm_{q_T, c_w}(T) = \uS_{q_T, c_w}(T) \cap \Sur_b(T)$ for $b=b_q$ by Lemma~\ref{q-Faser}\eqref{q-Faserb}, it follows that $\ul{\Kneb}^+$ is represented by the intersection of the $R$--scheme $\uKneb$ and the open subscheme $\Sur_b$ of $\uW(M)$. 
\sm 

\eqref{lem_representability-a} Without loss of generality, we can assume that $R=k$ is a field. 
By Lemma~\ref{lem_smooth_locus-LG}\ref{lem_smooth_locus-LGf}, the $S$--scheme $\uS\rmsm_{q_S, a\me}$ is  quasi-affine, hence also quasi-projective, and
has geometrically integral fibers. Moreover, the structure map $X\to \Spec(S)$ is surjective, since it is so over any algebraically closed field by Example~\ref{lem_smooth_locus_exam}\eqref{lem_smooth_locus_hypn}.
We can therefore apply \ref{weilres}\eqref{cgp-a} and obtain that the Weil restriction $X$ is geometrically integral. For the proof that $\uKneb$ has geometrically integral fibres, it is then enough to  show that the open subscheme $\uKneb$ of $X$ is not empty. To do so, we can extend scalars and assume that $q$ is isotropic. But then Lemma~\ref{knevid}\ref{knevid-a} finishes the proof. Finally, we use Lemma~\ref{todle}\eqref{todle-b} to conclude that $\uKneb$ is universally schematically dense in $X$. \sm
\sm 

\eqref{lem_representability-b} In view of \eqref{lem_representability-r} and \eqref{lem_representability-a}, it suffices to show that $\uKneb^+$ is not empty. Again, this can be checked over algebraically closed $R$--fields, in which case it follows from Lemma~\ref{knevid}\ref{knevid-b}. \end{proof}
\sm 



The following Lemma~\ref{lem_PR} is crucial for the further development. 

\begin{lem}\label{lem_PR} Let $R$ be a semilocal ring, let $(S,a)$ be a finite projective $R$--algebra of rank $d\ge 2$ and unit-generator $a$, and let $(M,q)$ be a quadratic $R$--space with $\rank_R M \ge 3$ for which $\uS_{q,a^{-1}}\rmsm(S) \not = \emptyset$.

Then $(S,a,q)$ admits  a Knebusch datum, i.e., $\Kneb(S,a,q) \ne \emptyset$.
Furthermore, if $\rank_R M \ge 4$, 
then $(S,a,q)$ admits  a special Knebusch datum.
\end{lem}

\comments{
(2023-05-03) I eliminated the assumption "or if all residue fields are infinite" in the last sentence of \ref{lem_PR} since for now it is not clear to me that this implies $\uKneb^+(S,a,q) \ne \emptyset$.
} 
\pcomments{(2026-06-08) Il me semble que c'est ce que l'on a fait le premier jour \`a Ottawa.}
\comments{(2026-06-10) Je ne rappel plus.}

\begin{proof} By assumption we can pick an element $u \in \uS_{q,a^{-1}}\rmsm(S)$.
\sm 

\begin{inparaenum}[(I)]
\item\label{lem_PR-I} We first consider the case that $R=k$ is an infinite field. According to 
    Proposition~\ref{smoLG}, the smooth $S$--sphere $\uS\rmsm := \uS_{q_S, a\me}\rmsm$ admits an open non-empty subscheme $\uU:=\uU_{q_S,u}$, which then is also an open subscheme of the $S$--sphere $\uS := \uS_{q_S, a\me}$. Furthermore, $\uU$ is isomorphic to an open subscheme of $\PP(M_S\ch)$, hence is smooth and quasi-projective, it is universally schematically dense, hence also $S$--dense  in $\uS\rmsm$, it has geometrically integral fibres, and $\uU(S) \ne \emptyset$. 
 
We can therefore apply \ref{weilres}\eqref{cgp-a} 
and get that $\uY = \frR_{S/k}(\uU)$ is a non-empty $k$--scheme. Also, by \ref{weilres}\eqref{weilres-b} and \ref{weilres}\eqref{weilres-d}, the $k$--scheme $\uY$ is an open subscheme of the affine $k$--scheme $\frR_{S/k}(\uS)$. It has a $k$--point since $\uU$ has an $S$--point. 
It follows that $Y(k)$ is Zariski dense in $\frR_{S/R}(\uS\rmsm)$.
On the other hand, we know from Proposition~\ref{lem_representability} that $\uKneb=\uKneb(S,a,q)$ is dense open in $\frR_{S/k}(\uS\rmsm)$ and from Lemma~\ref{knevid} that is is non-empty. But then
$\uY(k) \cap \uKneb(k) \ne \emptyset$.

If instead of $\rank M \ge 3$ we assume that $\rank M \ge 4$, we know
from Proposition~\ref{lem_representability} that $\uKneb^+=\uKneb^+(S,a,q)$ is dense open in $\frR_{S/k}(\uS\rmsm)$ and from Lemma~\ref{knevid} that it is non-empty. But then again $\uY(k) \cap \uKneb^+(k) \allowbreak \ne \emptyset$. We have now proven the claim for infinite fields. \sm

\item \label{lem_PR-II} Let us now consider the case that $R=k$ is an arbitrary field. If $q$ is isotropic, then, by Lemma~\ref{knevid}, we know that $(S,a,q)$ admits  a Knebusch datum, and even admits a special Knebusch datum, if $\rank M \ge 4$.  We may therefore assume that $q$ is anisotropic. Since $M$ is at least $3$--dimensional, this implies that $k$ is infinite (\cite[(12.3)]{Kneser}), and we are done by \eqref{lem_PR-I}. \sm

\item We proceed now to the general case, i.e., $R$ is semilocal. Let $J = \Jac(R)$ be the Jacobson radical of $R$ and let $\Max(R)$ be the set of maximal ideals of $R$. By \ref{slr},
    \[ \textstyle \prod_{\m \in \Max(R)} (R/\m) \cong R/J \quad \text{and} \quad
       S\ot_R (R/J) \cong S/SJ.
       \]  According to Proposition~\ref{prop_quadric1}, applied to the semilocal ring $S$ and the ideal $\fra=SJ\subset \Jac(S)$ of $S$, the canonical map $\uS\rmsm(S) \to  \uS\rmsm(S/\fra)$ is surjective. Hence so is
    \begin{align*}
    \big(\frR_{S/R}(\uS\rmsm)\big)(R) & = \uS\rmsm(S) \longto
     \uS\rmsm(S/\fra)\cong \big(\frR_{S/R}(\uS\rmsm)\big)(R/J)
      \\ & \cong \textstyle \prod_{\m \in \Max(R)} \frR_{S/R}(\uS\rmsm)(R/\m).
     \end{align*}
By \ref{lem_representability}, $\uKneb$ is an open subscheme of $\frR_{S/R}(\uS\rmsm)$ and by \eqref{lem_PR-II} we know that $(\uKneb)(R/\m) \ne \emptyset$ for all $\m \in \Max(R)$. Hence $\uKneb(R) \ne \emptyset$ by \eqref{opema-rem1}. The same argument works for $\uKneb^+(R)$ in case $\rank M \ge 4$.
\end{inparaenum} \end{proof}

\lv{
{\tt Old version of infinite field case:}

According to Lemma~\ref{lem_sphere} and Proposition \ref{prop_quadric}\ref{lem_sphere2_ii}, the sphere
$\uS\rmsm_{q_S,a^{-1}}$ admits an open $S$-subscheme $\uU_{q_S,u}$ which is isomorphic to
an open $S$-subscheme of $\PP(M_S^\vee)$ and for which $\uU_{q,u}(S)$ is not empty.

It follows from \ref{weilres}\eqref{weilres-c} that $Y=\frR_{S/k}(\uU_{q_S,u})$ is an open $k$-subscheme of $\frR_{S/k}(\uS_{q,a^{-1}}\rmsm)$ and is isomorphic to an open $k$--subscheme of $\frR_{S/k}\bigl(\PP(M_S^\vee)\bigr)$. We also know that $Y$
has a $k$--point. It follows that $Y(k)$ is dense in $\frR_{S/k}(\PP(M_S^\vee))$.
\sm

From Proposition~\ref{lem_representability} we know that  $\uKneb(S,a,q)$
is smooth, geometrically connected and is dense open in
$\frR_{S/k}(\uS_{q,a^{-1}}\rmsm)$, in particular, it is a non-empty scheme.
It follows that $Y(k) \cap \Kneb(S,a,q)(k)$ is non-empty.
\sm

We proceed now to the general case, i.e., $R$ is semilocal. Let $J = \Jac(R)$ be the Jacobson radical of $R$. Then $R/J = k_1 \times \dots \dots \times k_c$ is a product of fields. By the previously proven field case there exists $\overline{v}_j \in \uS_{q,a^{-1}}\rmsm(S \otimes_R k_j)$
defining a Knebusch datum.
This gives rise to an element \break
$\overline{v} \in \uS_{q,a^{-1}}\rmsm(S \otimes_R (R/J))$, $\overline{v}=\sum_{i=0}^{d-1} \overline{m_i} a^i$ with $\overline m_i \in M \otimes_R R/J$, such that $q( \overline{m}_{d-1}) \in (R/J)^\times$.
According to Proposition  \ref{prop_quadric}.\ref{prop_quadric1}, applied to the semilocal ring $S$ and the ideal $I = JS \subset \Jac(S)$, the element $\overline{v}$ lifts to an element $v= \sum_{i=0}^{d-1} m_i a^i \in \uS_{q,a^{-1}}\rmsm(S)$.  Since $\mathrm{Kneb}(S,a,q)$ is an open
subscheme of  $\frR_{S/R}(\uS_{q,a^{-1}}\rmsm)$, Lemma \ref{opema} enables us to conclude that $v$ defines a Knebusch datum over $R$.

If all residue fields of $R$ are infinite, then
$\mathrm{Kneb}^+(S,a,q)(k_i) \not = \emptyset$ for all $i$
(since each $\mathrm{Kneb}^+(S,a,q)$ is a $k_i$--rational variety).
The lifting  method above provides an element of $\mathrm{Kneb}^+(S,a,q)(R)$.

In the case  $r \geq 4$, Lemma \ref{lem_representability2}.\ref{lem_representability2+}
shows also that  $\mathrm{Kneb}^+(S,a,q)(k_i) \not = \emptyset$ for all $i$.
This case works then as well.

}
\ms

We can now prove \eqref{thm-kneb-a11}. It turns out to be an easy consequence of the strong Knebusch norm principle for unit-generators that are smooth values of $q_S$. 

\begin{lem}[Unit-generators that are smooth values]\label{pkn} Let $R$ be a semilocal ring, let $(M,q)$ be a  faithful quadratic $R$--space, and let $S$ be a finite projective $R$-algebra of degree $d\in \NN_+$ which is unit-generated. We denote the set of unit-generators of $S$ by $\Prim_R(S)\ti$.%
\sm

\begin{inparaenum}[\rm (a)] \item \label{prop_knebusch1} Let  $\uS_{q,a}\rmsm(S) \not = \emptyset$, i.e., $a\in \rmD\rmsm(q_S)$, and suppose that $a$ is a unit-generator. Then 
\begin{equation}
  \label{prop_knebusch1n}  \rmN_{S/R}(a) \in \rmD(q)^{[d]}.
\end{equation}
More generally, 
\begin{equation}\label{cor-kneb1}
 \rmN_{S/R}\big( \rmD\rmsm(q_S) \cap S\ti{}^2\Prim_R(S)\ti\big)
   \subset \rmD(q)^{[d]}.
\end{equation} 

\item \label{prop_knebusch1b} If all residue fields of $R$ are infinite,  then \eqref{thm-kneb-a11} holds:
\begin{equation}\label{prop_knebusch1b1}
\rmN_{S/R}\big( \rmD\rmsm(q_S)\big) \subset \rmD(q)^{[d]}.
\end{equation}  
  
\item \label{prop_knebusch-resa} If $d\ge 2$ and $\rank M \ge 4$, we can replace $\rmD(q)^{[d]}$ on the right-hand side of \eqref{prop_knebusch1n}, \eqref{cor-kneb1} and \eqref{prop_knebusch1b1} by $\rmD\rmsm(q)^{[d]}$.
\end{inparaenum} 
\end{lem}

\begin{proof} \eqref{prop_knebusch1} The equation~\eqref{prop_knebusch1n} says that the strong norm principle holds for a unit-generator $a\in \rmD\rmsm(q)$. To prove this, we can assume that $M$ has constant rank $r\ge 3$ and that $d\ge 2$, see \ref{rkn}\eqref{dqd-ele-f}. In this situation we will prove \eqref{prop_knebusch1n} and \eqref{prop_knebusch-resa} at the same time, using induction on $d$. Our proof works whenever $R$ satisfies the assumption 
\begin{enumerate}[label={\rm (\greek*)}]
\item \label{prop_knebusch11} {\em $\uKneb_{(S',a',q)}(R) \ne \emptyset$ for all finite free,  unit-generated $R$--algebras $(S',a')$, }     
\end{enumerate}
which by Lemma~\ref{lem_PR} is satisfied for a semilocal $R$. Assuming \ref{prop_knebusch11}, we know $\uS_{q,a\me}\rmsm(S) \ne \emptyset$  because $\uS_{q,a}\rmsm(S) \simlgr \uS_{q,a\me}\rmsm(S)$ via $m \mapsto a\me m$, $m\in M_S$. Hence, there exists a Knebusch datum $\uKneb_{(S,a,q)}(R)$ over $R$. Then Lemma~\ref{knebp}\ref{knebp-c} and induction finishes the proof of \eqref{prop_knebusch1n} and \eqref{prop_knebusch-resa} for \eqref{prop_knebusch1n}.  

Because 
\begin{align*}
  &\rmN_{S/R}\big( \rm\rmD\rmsm(q_S) \cap S\ti{}^2\Prim_R(S)\ti\big) 
 \\ &\quad   = \textstyle \bigcup_{a\in \Prim_R(S)\ti}
        \rmN_{S/R}\big(\rmD\rmsm(q_S) \cap S\ti{}^2 a \big),          
\end{align*}
the formula \eqref{dqd-ele-h1} says that \eqref{cor-kneb1} follows from \eqref{prop_knebusch1n}. 
\sm 
  
\eqref{prop_knebusch1b} If all residue fields of $R$ are infinite, Lemma~\ref{hunco} says that $S\ti = S\ti{}^2 \Prim_R(S)\ti$. Hence \eqref{prop_knebusch1b1} follows from \eqref{cor-kneb1}.  


\eqref{prop_knebusch-resa} It remains to deal with \eqref{cor-kneb1} and \eqref{prop_knebusch1b1}. But these follow immediately from the proof above.   
\end{proof}

\comments{(2023-12-07) ``Semilocal'' is only used to be able to apply Lemma~\ref{lem_PR}, see the assumption \ref{prop_knebusch11} in the proof of \ref{prop_knebusch}. The crucial point of \ref{prop_knebusch}\ref{prop_knebusch2} is the case $\uS_{q,1}\rmsm (R) \ne \emptyset$. This case implies all the remaining ones.}

An important step in the proof of Theorem~\ref{thm-kneb} is the following Lemma~\ref{prop_knebusch}, dealing with $\rmN_{S/R}(a)$ for unit-generators $a\in S$. 

\begin{lem} \label{prop_knebusch}
Let $R$ be a semilocal ring, let $(M,q)$ be a  faithful quadratic $R$--space,
and let $S$ be a finite projective $R$-algebra of degree $d\in \NN_+$ which is unit-generated by $a \in S\ti$. \sm

\begin{enumerate}[label={\rm (\alph*)}]
\item \label{prop_knebusch2}  Assume $\uS_{q,a}(S) \ne \emptyset$.
 \begin{enumerate}[label={\rm (\roman*)}]
    \item \label{prop_knebusch2i} If $\uS_{q,1}\rmsm (R) \ne \emptyset$ or if
       $d$ is even, then $\rmN_{S/R}(a) \in \rmD(q)^{[2d]}$.

    \item \label{prop_knebusch2ii} If $d$ is odd, then $\rmN_{S/R}(a) \in \rmD(q)^{[2d+1]}$.
 \end{enumerate} \sm

\item \label{prop_knebusch-res} Suppose $d\ge 2$  and $\rank M \ge 4$.
Then we can replace $\rmD(q)$ in \ref{prop_knebusch2} above by $\rmD\rmsm(q)$.
\end{enumerate}
\end{lem}

\begin{proof} 

\ref{prop_knebusch2} We first show that \ref{prop_knebusch2} holds whenever the strong norm principle is satisfied. Indeed, in this case we have $\rmN_{S/R}(a) \in \rmD(q)^{[d]}$. If $\uS_{q,1}\rmsm(R) \ne \emptyset$, then $\rmD(q)^{[d]} \subset \rmD(q)^{[2d]}$ by \ref{dqd-ele}\eqref{qdq-ele-aa}, and if $d$ is even or odd, then  $\rmD(q)^{[d]} \subset \rmD(q)^{[2d]}$ or 
$\rmD(q)^{[d]} \subset \rmD(q)^{[2d+1]}$ respectively by \eqref{dqd-ele-a4}.

In general, we can assume for the proof of \ref{prop_knebusch2} that $M$ has constant rank $r$. If $r=1$ or $r=2$, then the strong norm principle holds and we are done by the above. For the remainder of the proof we can therefore suppose that $r \ge 3$.  

Let us next deal with the case  $\uS_{q,1}\rmsm(R) \ne \emptyset \ne \uS_{q,a}(S)$. 
The remaining proof of \ref{prop_knebusch2} and \ref{prop_knebusch-res} works for any $R$  under the assumption~\ref{prop_knebusch11} of the proof of \ref{pkn} and the assumption \ref{prop_knebusch12}, which is 
\begin{enumerate}[label={\rm (\greek*)}]\setcounter{enumi}{1}
\item \label{prop_knebusch12} {\em  Suppose $\rank M \ge 3$. If $\uS\rmsm_{q_S, 1}(S) \ne \emptyset$ and $\uS_{q_S, a} \ne \emptyset$, there exist $u,v\in M_S$ such that $Su + Sv$ is a complemented submodule of $M_S$ which is free of rank $2$ with basis $(u,v)$ satisfying $q_S(u) = a$ and $q_S(v) = 1$. } 
\end{enumerate}
By Lemma~\ref{acsoLG}, the condition \ref{prop_knebusch12} holds whenever $R$ 
is semilocal. 

Thus, assuming \ref{prop_knebusch12}, we are in the situation of Lemma~\ref{redsm}. Using the notation there,  
and applying \eqref{prop_knebusch1n}  for $a'$, we have $\rmN_{S/R}(a)=\rmN_{S'/R}(a') \in \rmD(q)^{[2d]}$, and even $\rmN_{S/R}(a) \in \rmD\rmsm(q)^{[2d]}$ if \ref{prop_knebusch-res} holds.%
\sm

Having established the first alternative in \ref{prop_knebusch2i}, we will prove the second alternative in \ref{prop_knebusch2i}, i.e., $d$ even, part \ref{prop_knebusch2ii}, i.e., $d$ is odd, and the remaining part of \ref{prop_knebusch-res} all at the same time. 
By  \ref{lem_smooth_locus_exam}\eqref{lem_smooth_locus_exam-c}, there exists $x\in R\ti$ such that $\uS_{q,x}\rmsm(R) \ne \emptyset$. Then $xa$ is a unit-generator of $S$. We have $\uS_{xq, xa}(S) = \uS_{q,a}(S)\ne \emptyset$, and $\uS_{xq, 1}\rmsm (R) \cong \uS_{q,x}\rmsm (R) \ne \emptyset$.
Hence, we can apply the first alternative of \ref{prop_knebusch2i} to $(M,xq)$ and $xa$. Thus, $x^d \rmN_{S/R}(a) = \rmN_{S/R}(xa) \in \rmD'(xq)^{[2d]}= x^{2d}\rmD'(q)^{[2d]}$ where $\rmD'=\rmD$ in general or $\rmD'=\rmD\rmsm$ if \ref{prop_knebusch-res} holds. Therefore $\rmN_{S/R}(a) \in x^d\rmD'(q)^{[2d]}$. If $d$ is even, then $x^d\rmD'(q)^{[2d]}= \rmD'(q)^{[2d]}$ by \eqref{dqd-ele-a1}, and if $d$ is odd, then, by the same formula, $x^d\rmD'(q)^{[2d]}= x\rmD'(q)^{[2d]}\subset \rmD'(q)^{[2d+1]}$ because $x\in \rmD\rmsm(q)\subset \rmD(q)$. \end{proof}
\sm

The following corollary of Lemma~\ref{prop_knebusch} extends the validity of the formulas established there.


\begin{cor}\label{cor-kneb} We use the assumptions of Lemma~{\rm \ref{prop_knebusch}:} $R$ is a semilocal ring, $(M,q)$ is a faithful
quadratic $R$--space, 
and $S$ is a finite projective unit-generated $R$-algebra of degree $d\in \NN_+$. We further assume that the assumption~\ref{prop_knebusch12} of the proof of Lemma~{\rm \ref{prop_knebusch}} is satisfied. 
Then
\begin{equation} \label{cor-kneb2}
  \rmN_{S/R}\big( \rmD(q_S) \cap S\ti{}^2\Prim_R(S)\ti\big) \subset
  \begin{cases}
   \rmD(q)^{[d]}, & \text{if $q$ is regular,} \\
   \rmD(q)^{[2d]}, & \text{if $\uS_{q,1}\rmsm (R) \ne \emptyset$} \\
         & \text{ or if $d$ is even,} \\
    \rmD(q)^{[2d+1]}, & \text{if $d$ is odd}.
  \end{cases}
\end{equation}
Furthermore, if $d\ge 2$ and $\rank M \ge 4$,
we can replace $\rmD(q)$ by $\rmD\rmsm(q)$ on the right-hand side of \eqref{cor-kneb2}. 
\end{cor}

\begin{proof} Because 
\[ \rmN_{S/R}\big( \rmD(q_S) \cap S\ti{}^2\Prim_R(S)\ti\big) 
= \textstyle \bigcup_{a\in \Prim_R(S)\ti}
        \rmN_{S/R}\big(\rmD(q_S) \cap S\ti{}^2 a \big),          
\]  
the formula \eqref{dqd-ele-h1} says that \eqref{cor-kneb2} follow from the cases discussed in Lemma~\ref{prop_knebusch}. Analogously,  if $d\ge 2$ and $\rank M \ge 4$,  we can replace $\rmD(q)$ by $\rmD\rmsm(q)$ everywhere in view of  the refinement \ref{prop_knebusch-res} of Proposition \ref{prop_knebusch}. If $q$ is regular, then $\rmD\rmsm(q_S) = \rmD(q_S)$ by  Lemma~\ref{lem_smooth_locus-LG}\ref{lem_smooth_locus-LGd}. Therefore the first formula of \eqref{cor-kneb2} is a special case of \eqref{cor-kneb1}. \end{proof}

\subsection{Proof of Theorem~\ref{thm-kneb}} \label{thm-kneb-proof} Recall the setting: $R$ is semilocal, $S\in \Ralg$ is finite \'etale of degree $d\in \NN_+$, and $(M,q)$ is a faithful quadratic space over $R$.
\sm 

\ref{thm-kneb-b} We apply Proposition~\ref{prop_tiroir}. Taking an odd $n$ in loc.\ cit. and using transitivity of ``\'etale'', we obtain a tower $R\to S \to S'$ of \'etale $R$--algebras with $S'/R$ of constant degree $dn$. Moreover, for every $a\in \rmD(q_S)$ there exists $b\in S'$ and $c\in \Prim_R(S')\ti$ such that $a= b^2 c$, so that we can apply \eqref{cor-kneb2} to $a\in S'$. By transitivity and multiplicativity of the norm, $\rmN_{S'/R}(a) = \rmN_{S/R}\big( \rmN_{S'/S}(a \ot 1_{S'})\big) = \rmN_{S/R}(a^n) = \rmN_{S/R}(a)^n \equiv \rmN_{S/R}(a) \mod R\ti{}^2$. Hence, by 
\eqref{cor-kneb2} we get $\rmN_{S/R}(a) \in \rmD(q)^{[d']}$ where
\[ d' = \begin{cases}
     dn, & \text{$q$ regular}, \\ 2dn, & \text{$\uS_{q,1}\rmsm (R) \ne \emptyset$
             or $d$ even} \\ 2dn+1 , & \text{$d$ odd}.
\end{cases} \]
This implies \eqref{thm-kneb-b1}. Then \eqref{thm-kneb-b2} follows by multiplicativity of $\rmN_{S/R}$. \sm

\ref{thm-kneb-a} By assumption, all residue fields of $R$ are infinite. Hence
$S\ti = S\ti{}^2 \Prim_R(S)\ti$ by Lemma~\ref{hunco}. 

\comments{(2025-05) Lemma~\ref{hunco} requires only that $R$ satisfies the primitive criterion}

By \eqref{cor-kneb1} and  \eqref{cor-kneb2} we then get
\[ \rmN_{S/R}\big( \rmD\rmsm(q_S)\big) \subset \rmD(q)^{[d]}\]
and
\begin{equation*} \rmN_{S/R}\big(\rmD(q_S)\big) \subset
 \begin{cases}  \rmD(q)^{[d]} & \text{if $q$ is regular;} \\
        \rmD(q)^{[2d]} & \text{if $\uS_{q,1}\rmsm(R) \ne \emptyset$ or if $d$ is even;}\\
         \rmD(q)^{[2d+1]} & \text{if $d$ is odd}
 \end{cases}
\end{equation*}
which are  \eqref{thm-kneb-a1} and \eqref{thm-kneb-a11}. We also obtain that we can replace $\rmD(q)$ by $\rmD\rmsm(q)$  if $d\ge 2$ and $\rank M \ge 4$.
\qed

\lv{
\comments{(2026-05-12) We want to use Lemma~\ref{regp} in the proof of Lemma~\ref{prop_knebusch}, but we do not have the assumptions \ref{regp-i} and \ref{regp-ii} there. There are results similar to \ref{regp} using the action of $\orth(q)$ or $\SO(q)$. The formulation of Lemma~\ref{regp} is such that it can easily be changed, e.g. by elimination assumption~\ref{regp-ii}} 
\newpage
\begin{lem}\label{regp} Let $R$ be a semilocal ring, let $(M,q)$ be a quadratic space with $\rank_R M \ge 2$, and let $a,z\in R\ti$ such that $\uS_{q, a}(R) \ne \emptyset \ne \uS_{q,z}\rmsm(R)$. Furthermore, assume
\begin{enumerate}[label={\rm (\roman*)}]
  \item \label{regp-i} $|R/\gm| \ge 4$ for all maximal ideals $\gm \ideal R$, and 
      
  \item\label{regp-ii} $\uS_{q,a}(R) \to \uS_{q,a}\big(R/\Jac(R)\big)$ is surjective.    
 
\end{enumerate}
Then there exist $u,v\in M$ such that $Ru + Rv$ is a regular, hence complemented plane with basis $(u,v)$ satisfying $u\in \uS_{q,a}(R)$ and $v \in \uS_{q,z}\rmsm(R)$.   
\end{lem}
 
With regard to the assumption~\ref{regp-ii}, we note that $\uS\rmsm_{q,a}(R) \to \uS\rmsm_{q,a}\big(R/\Jac(R)\big)$ is surjective by Proposition~\ref{prop_quadric1}. Hence, this assumption concerns the singular points of $\uS_{q,a}\big(R/\Jac(R)\big)$. In particular, the assumption~\ref{regp-ii} is unnecessary, if $q$ is regular (since then $\uS_{q,a}(R) = \uS_{q,a}\rmsm(R)$) or if $\Jac(R) = 0$.    
\sm

\begin{proof} \begin{inparaenum}[(I)]
\item \label{regpI} We prove the lemma under the assumption that $(M,q)$ contains a hyperbolic pair $(e,f)$ and that $R$ is an LG ring. For $c\in R\ti$ we define $u_c = c\me a e + c f$ and let $v= e + zf$. Then $q(u_c) = a$ and $q(v) = z$ (note $v\in \uS_{q,z}\rmsm(R)$ by \ref{lem_smooth_locus_exam}\eqref{lem_smooth_locus_hypn}). Moreover, $(u_c, v)$ is a basis of $Re \oplus Rf$ if and only if $az - c^2 \in R\ti$. We are thus left with verifying that 
    the polynomial $ac - X^2$ represents a unit over $R$. By the defining property of an LG ring, it suffices to show that the polynomial represents a unit over every field $R/\gm$. But for them, the assumption \ref{regp-i}  guarantees the existence of an element as  required.%
    \sm 
    
\item\label{regpII} In this part of the proof we show that the lemma holds if $R=k$ is a field with $|k| \ge 4$. By \eqref{regpI} and \ref{isotrop}\eqref{isotrop-d}, we can assume that $(M,q)$ is anisotropic. 
    
    We can pick $v\in \uS\rmsm_{q,z}(R)$ and $u'\in \uS_{q,a}(R)$. If $u'\not\in kv$, we are done. So, let us assume that $u'\in kv$, say $u'= sv$ for $s\in k\ti$ satisfying $s^2 z = a$. As $v$ is a smooth point, we know $W = \{w\in M : b_q(w,v) = 1\} \ne \emptyset$. Also, by Lemma~\ref{mx_lem}\eqref{mx_lemb}, we have $M = kw \oplus M_v$ for any $w\in W$ and $M_v = (kv)^\perp$, which has positive dimension. Because $W= W + M_v$, it follows that there exists $y\in W \setminus kv$.
\lv{
Indeed, assume otherwise, $W \subset kv$ and let $w=tv\in W$. Then $1=b_q(w,v) = 2t q(v)$ forces $\Char(k) \ne 2$, whence $tv + w'\in W$ for any $w'\in M_v$, contradiction.
}
 Since $(M,q)$ is anisotropic, we can apply the reflection $\rho_y$ and get $u:= \rho_y(u') = u' - b_q(u', y) q(y) \me y \not\in kv$ because $b_q(u',y) = s b_q(v,y) = s \ne 0$. The plane $Ru \oplus Rv = Ry \oplus Rv$ is regular if and only if 
\[ \det \begin{pmatrix} 2q(y) & 1 \\ 1 & 2z \end{pmatrix} = 4 z q(y) - 1 \ne 0.
 \]
  
\sm 

\item\label{regpIII} Finally, we deal with the case of a semilocal $R$. Let $m \mapsto \ol m$ be the reduction $M \to \ol M = M / \Jac(R) M$, and use the analogous notation for $R \to \ol R = R/\Jac(R)$. We know $\uS_{q, a}(\ol R) \ne \emptyset \ne \uS_{q, z}\rmsm(\ol R)$. By a straightforward extension of \eqref{regpII}, there exist $m_1, m_2\in M$ such that $\ol m_1 \in \uS_{q,a}(\ol R)$, $\ol m_2 \in \uS_{q,z}\rmsm(\ol R)$ and $\ol R \ol m_1 + \ol R\ol m_2$ is a regular plane in $(\ol M, \ol q)$ with basis $(\ol m_1, \ol m_2)$. Since $\uS_{q,z}\rmsm \to \uS\rmsm_{\ol q, \ol z}$ is surjective by Proposition~\ref{prop_quadric1}, it follows that we can lift $\ol m_1$ to $v\in \uS_{q,z}\rmsm$. We can also lift $\ol m_2$ to some $u\in \uS_{q,a}(R)$ by assumption~\ref{regp-ii}. One easily checks that $Ru + Rv$ is a regular plane in $(M,q)$ with basis $(u,v)$, see the proof of \cite[I, Cor.~(0.4)]{Ba}. 
\end{inparaenum}
\end{proof}
}

\comments{(2026-05) There is more stuff here, 3 pages: a proof in case $R$ is a field, $q$ arbitrary,  and Knebusch elements: 
\begin{align*}
  &\rmKE(A/R) = \{ a\in A\ti: R[a]/R \text{ is finite projective of constant rank}
   \\ & \qquad \text{and $A$ is finite projective of constant rank as $R[a]$--module}\} .
\end{align*}
(not to be published) }

\lv{

\begin{lem}\label{moti} Let $R$ be arbitrary, let $(M,q)$ be a faithful quadratic $R$--module satisfying for all $m\in M$ the condition
\begin{equation} \label{moti0}
0 \ne m\in M \implies q(m)\in R\ti,
\end{equation}
and let $S\in \Ralg$ be finite projective of constant rank $d$ and unit-generated by $a\in \rmD(q_S)$. Then $\rmN_{S/R}(a) \in \Dqd$.
\end{lem}

\comments{(2026-05) For $R \ne 0$, $M$ faithful, condition \eqref{moti0} implies that $r^2 \cdot q(M) \subset R\ti$ whenever $rPM\ne 0$, hence $R^2 \subset R\ti$, linearize to get $R\subset R\ti$.}

\begin{proof} By Lemma~\ref{genolem}, $S=R[X]/(P)$ for some monic $P\in R[X]$ of degree $d\in \NN_+$. We will prove the lemma by induction on $d=\rank S = \deg(P)$. The case $d=1$ is trivial, since then $S=R$ and so $\rmN_{S/R} = \Id$.

Let therefore $d>1$. Since $\rmD(q_S) = \rmD(q_S) \me$, we have $a\me \in \rmD(q_S)$, i.e., there exists $v_S \in M_S$ such that $a\me = q_S(v_S)$, equivalently $1_S = a q_S(v_S)$. In view of $S= \bigoplus_{i=0}^{d-1} R a^i$, we have $M_S = \bigoplus_{i=0}^{d-1} M a^i$ and therefore can write $v_S= \sum_{i=0}^{d-1} v_i a^i$ with unique $v_i \in M$. Let us define $v(X) \in M_{R[X]}$ by $ v(X) = \sum_{i=0}^{d-1} v_i X^i$. Applying Euclidean division (\cite[IV, \S1.6, Cor]{BA5}), there exist unique polynomials $\wtl Q(X)$ and $U(X)\in R[X]$ such that
\begin{equation}\label{moti11}
 Xq\big(m(X)\big) - 1 = \wtl Q(X) P(X) + U(X), \quad \deg U(X) <d.
\end{equation}
Specializing $X$ to $a$ in the equation above, we get $v(X)|_{X\rightsquigarrow a} = v_S$, hence $Xq_{R[X]}\big(v(X)\big)|_{X\rightsquigarrow a} = a q_S(v_S) = 1$. Specializing $X\rightsquigarrow a$ on the right-hand side of \eqref{moti11}, we find $U(a) = 0$, and therefore $U(X) = 0$.
Define
\[ e:= \deg(Q) \qquad
    k := \max\{ i : v_i \ne 0 \}\le d-1.  \]
Also, let $c$ be the leading coefficient of $\wtl Q(X)$. We compare the leading coefficients on both sides of \eqref{moti11}:
\begin{align*}
  X q\big( m(X)\big) &= q(m_k) X^{2k+1} + (\text{lower terms}) \\
    &= c X^{d+e} + (\text{lower terms}),
\end{align*}
thus
\begin{equation}\label{moti4}  c=q(m_k)\in R\ti
\end{equation}
by \eqref{moti0}. We can therefore replace $\wtl Q(X)$ by $cQ(X)$ for a monic $Q\in R[X]$ and rewrite \eqref{moti11} in the form
\begin{equation}
  \label{moti1} X q(v(X)) = 1 + c P(X) Q(X) .
\end{equation}
Putting $X=0$ in \eqref{moti1} yields
\begin{equation}
  \label{moti3} cP(0) Q(0) = -1.
\end{equation}
The comparison of the leading terms also shows $d+ e = 2k + 1 \le 2d-2 + 1 = 2d-1$, therefore
\begin{equation}  \label{moti5}
 e \le d-1.
\end{equation}
Moreover, $2k + 1 = d+e$ implies $d+e \equiv 1 \mod (2)$, or
\begin{equation}  \label{moti6}
 e \equiv d +1 \mod (2).
\end{equation}
Define
\[ T=R[X]/(Q)\quad \text{and}\quad t= X + (Q)\in T,\] so that $T$ is one-generated by $t\in T$ and is free of rank $e<d$ by \eqref{moti5}.
Putting $X=t$ in \eqref{moti1} we get $t q(m(t)) = 1$, thus $t\me \in \rmD(q_T) = \rmD(q_T)\me$ and therefore
\begin{equation}\label{moti2}
  t \in \rmD(q_T),
\end{equation}
proving that $T$ is unit-generated by $t$.
We can now calculate the norms of $a$ and $t$: $\rmN_{S/R}(a) = (-1)^d P(0)$ and $\rmN_{T/R}(t) = (-1)^e Q(0)$ by \ref{trno}\eqref{trno-bou}. Hence, by \eqref{moti3} and \eqref{moti6},
\begin{align*}
  \rmN_{S/R}(a) &= (-1)^d P(0) = (-1)^d ( - c\me Q(0)\me) = c\me\big( (-1)^{d+1} Q(0)\big) \\ &= c\me \big( (-1)^e Q(0)\big)\me = c\me \rmN_{T/R}(t)\me.
\end{align*}
Because of \eqref{moti2} and \eqref{moti5}, we can apply induction to $(T,t)$ and get $\rmN_{T/R}(t) \in \rmD(q)^{[e]}$. Hence, by \eqref{moti4}, we have
\[ \rmN_{S/R}(a) = c\me N_{T/R}(t)\me \in (\rmD(q)^{e+1}){}\me = \rmD(q)^{[e+1]} \subset \rmD(q)^{[d]}, \]
where we used \eqref{dqd-ele-a3} and \eqref{dqd-ele-moti} because of \eqref{moti6}.
\end{proof}

\comments{(2023-11-16) We use \ref{moti-an}\eqref{moti-an-field} in \ref{knefi}}

\subsection{Remarks regarding \ref{moti}.} \label{moti-an}
\begin{inparaenum}[(a)] \item \label{moti-an-field} {\em Let $F$ be a field, let $(M,q)$ be a quadratic form over $F$ and let $S\in \Falg$ be $d$--dimensional and unit-generated by $a\in \rmD(q_S)$. Then $\rmN_{S/F}(a) \in \rmD(q)^{[d]}$.}
\sm

In order to prove this, we can assume that $\rad(q) = 0$ since $\rmD(q) = \rmD(\bar q)$ where $\bar q \co \ol M \to F$ is the quadratic form of \ref{radqf}. If $(M,q)$ is isotropic, $(M,q)$ contains a hyperbolic plane \cite[Prop.~7.13]{EKM}

\inparcom{(2023-12-16) Need to include \cite[Prop.~7.13]{EKM} in the discussion on "isotropy $\implies$ existence of hyperbolic plane", see \ref{isotrop}\eqref{isotrop-c}.}

and hence $\rmD(q) = F\ti$, so that the claim obviously holds. We can therefore assume that $q$ is anisotropic, i.e., the condition \eqref{moti0} holds. The claim then follows from Lemma~\ref{moti}.
\sm

\item\label{mori-an-arb} ({\em $R$ arbitrary }) On the other hand,
the condition \eqref{moti0} is rather strong for a general $R$. An analysis of the proof shows that \eqref{moti0} can be replaced by the following assumption: Given $(M,q)$, $S$ and $a$ as in {\rm \ref{moti}}, there exists $v(X) \in M\ot_R R[X]$ such that \end{inparaenum}
\begin{enumerate} [label={\rm (\roman*)}]
  \item \label{moti-ani} $av(a) = 1$, i.e., $v(a) \in \uS_{q,a^{-1}}(S)$ and

  \item \label{moti-anii} $X q\big(v(X)\big)= 1+ c P(X) Q(X)$ for some $c \in R\ti$ and a monic polynomial $Q(X) \in R[X]$.
\end{enumerate}
The conditions \ref{moti-ani} and \ref{moti-anii} capture the essence of Knebusch data, defined in \ref{defn_data}. However, in order to prove the existence of such a $v(X)$, equivalently the existence of a Knebusch datum, we need to impose stronger conditions on $v(X)$. We will prove the existence of a Knebusch datum under suitable assumption on $R$ and $q$ in Lemma~\ref{lem_PR}.

Proposition~\ref{prop_knebusch} is crucial for the proof of Theorem~\ref{thm-kneb}, Knebusch's norm principle for finite \'etale extensions. It is also the main step in the proof of the (strong) Knebusch norm principle for fields, given in \ref{knefi}.

\begin{prop}[Knebusch's norm principle over fields]\label{knefi}
Let $k$ be a field, let $(M,q)$ be a regular quadratic $k$--module of positive dimension and let $K/k$ be a finite field extension of degree $d$. Then
\[ \rmN_{K/k}\big(\rmD(q_K)\big) \subset \rmD(q)^{[d]}.\]
\end{prop}

\comments{(2023-11) We should try to get rid of the assumption ``regular'' in \ref{knefi}. `Regular'' is only needed for \eqref{knefi-I}.
\sm

Note that one can replace the use of Proposition~\ref{prop_knebusch}\ref{prop_knebusch1} by the more general
Remark~\ref{moti-an}\eqref{moti-an-field}: it holds for arbitrary $q$ and arbitrary $d$--dimensional extensions of $k$!}

\begin{proof} We will proceed by induction on the degree of $K/k$. Let $a\in \rmD(q_K)$. We first consider two special cases.

\begin{inparaenum}[(I)] \item\label{knefi-I} Suppose $a\in k$. Then $\rmN_{K/k}(a) = a^d$. If $d$ is even, then $a^d\in R\ti{}^2 \subset \rmD(q)^{[2]} \subset \cdots \subset \rmD(q)^{[d]}$ by \eqref{dqd-ele-a4}. On the other hand, if $d=2n+1$ is odd, then $a^d = a a^{2n} \in aR\ti{}^2\subset a \rmD(a)^{[2n]}$. Again by \eqref{dqd-ele-a4}, it then suffices to show $a\in \rmD(q_K) \implies a\in \rmD(q)$. Or, this is property \ref{equi-C} of Proposition~\ref{equi}, which by \ref{equisa}\eqref{equisa-c} holds for finite field extensions.

\item\label{knefi-II}  We can easily dispose of another extreme case: $K=k[a]$, the subfield generated by $a$. In this case $\rmN_{K/k}(a) \in \rmD(q)^{[d]}$ by Remark~\ref{moti-an}\eqref{moti-an-field} (or by Proposition~\ref{prop_knebusch}\ref{prop_knebusch1}).

\item  In view of \eqref{knefi-I} and \eqref{knefi-II} we are left with the case $k\subsetneq E:=k[a] \subsetneq K$. Then $d_1 = [E:k]< d$ and $d_2 = [K:E]< d$. By \eqref{knefi-I} for the extension $K/E$ we know $\rmN_{K/E}(a) \in \rmD(q_E)^{[d_2]}$, and by induction we also know $\rmN_{E/k}\big(\rmD(q_E)\big) \subset \rmD(q)^{[d_1]}$, which by multiplicativity of the norm implies $\rmN_{E/k}(\rmD(q_E)^{[d_2]}) \subset \big(\rmD(q)^{[d_1]}\big){}^{[d_2]} \subset \rmD(q)^{[d_1d_2]} = \rmD(q)^{[d]}$. Hence altogether $\rmN_{K/k}(a) = \rmN_{E/k}\big( \rmN_{K/E}(a)\big) \subset \rmD(q)^{[d]}$.
\end{inparaenum}
\end{proof}
\sm

\textbf{Remarks.} Our proof above follows Knebusch's original proof in \cite{knebusch-norm} which is reproduced in \cite[VII, (5.1)]{Lam-qf}. In the field case, the proof of Proposition~\ref{prop_knebusch}\ref{prop_knebusch1} is much simpler. On the other hand, Knebusch's norm principle for finite field extensions is true for arbitrary quadratic forms, not only regular ones. The proof in this generality, given in \cite[18.10]{EKM}, uses a completely different approach.

\subsection{Knebusch elements (not to be published)} \label{kel}  While Proposition~\ref{prop_knebusch} and Corollary~\ref{cor-kneb} were general enough to prove Theorem~\ref{thm-kneb}, they hint at a different approach to the Knebusch norm principle: rather than demanding the principle to hold for all elements of $\rmD(q_A)$, $A\in \Ralg$ finite projective, we can without difficulty prove the principle for some. Namely, let $A\in \Ralg$ be finite projective as $R$--module and put
\begin{align*}
  &\rmKE(A/R) = \{ a\in A\ti: R[a]/R \text{ is finite projective of constant rank}
   \\ & \qquad \text{and $A$ is finite projective of constant rank as $R[a]$--module}\} .
\end{align*}
Let $a\in \rmKE(A/R)$. Then the $R$--module $A$ has constant rank:
\[ \rank_R A = (\rank _R R[a]) \cdot (\rank_{R[a]} A), \]
e.g.\ by \cite[3.2.13]{Ford}. We then have: \sm

{\em Let $R$ be semilocal, let $(M,q)$ be a faithful regular quadratic $R$--module, and $A\in \Ralg$ is finite projective of constant rank $d$. Then }
\begin{equation}\label{kel1}
\rmN_{A/R}(a) \in \rmD(q)^{[d]} \quad\text{for every $a\in \rmD(q_{R[a]}) \cap \rmKE(A/R)$.}
\end{equation}
Indeed, let $d_1 = \rank_R R[a]$ and $d_2 = \rank_{R[a]} A$. Then $d=d_1 d_2$ by the formula above, $\rmN_{A/R}(a) = \rmN_{R[a]/R}\big(\rmN_{A/R[a]}(a)\big) = \rmN_{R[a]/R}(a^{d_2}) = \rmN_{R[a]/R}(a)^{d_2}$ by transitivity,  \ref{trno}\eqref{trno-b}, and $\rmN_{R[a]/R}(a) \in \rmD(a)^{[d_1]}$ by \ref{prop_knebusch}\ref{prop_knebusch1}. Thus, altogether, $\rmN_{A/R}(a) \in (\rmD(q)^{[d_1]})^{d_2} \subset \rmD(q)^{[d_1d_2]} = \rmD(q)^{[d]}$.
\sm

{\em If we assume that $A/R[a]$ is finite \'etale of constant rank for $a\in \rmKE(A/R)$, then}
\begin{equation}\label{kel2} \rmN_{A/R}\big(\rmD(q_A) \cap \rmKE(A/R)\big) \subset \rmD(q)^{[d]},
\end{equation}
Note the difference between $\rmD(q_{R[a]})$ and $\rmD(q_A)$.

Indeed, this can be proven along the lines of the proof of \eqref{kel1}, but with an additional argument: first, as explained in \ref{dqd-ele}\eqref{dqd-ele-c} we can assume $1\in \rmD(q)$. We know $a\in \rmD(q_A)\cap R[a]$. If  $d+2$ is even, then $\rmN_{A/R}(a) = \rmN_{R[a]/R}(a)^d_2 \in R\ti{}^2 \subset \rmD(q) \subset \rmD(q)^{[d]}$ by \ref{qdq-ele}\eqref{qdq-ele-aa}, and we are done. Otherwise, if $d_2$ is odd, then ...

The question then becomes: When is $\rmK\rmE(S/R) = S$?
}

\newpage

\section{Cohomological consequences}\label{sec:consequences}

In this section we apply the Knebusch norm principle to prove that in the semilocal case $\rmD(q)^{[\rm ev]}$ is the image of the spinor norm, \ref{spinim}\eqref{spinim-b}, which in turn is used to establish injectivity of
$H^1\fppf(R, \uSpin(q)) \to H^1\fppf(S, \uSpin(q_S))$ in case $S\in \Ralg$ is finite \'etale of odd degree, \ref{prop_sn}\eqref{prop_sn-b}. The second part of this section deals with the closed immersion $\uSO(q) \to \uAut\big( \Cli_0(q)\big)$, which will be used in the following section \ref{sec:normgroups} on norm groups.
\ms

The following Lemma~\ref{cohCT} is a preparation for Lemma~\ref{abcco} and
Proposition~\ref{prop_sn}, where we will consider $\uSpin(q)$ instead of $\uO(q)$ and $\uSO(q)$.

\comments{(2026-05-13) Added the assumption that $R$ is LG in \ref{cohCT}. The previous version of \ref{cohCT} did not have any assumption on $R$. But all the results in this section assume that $R$ is semilocal. So I guess, we just forgot to assume that $R$ is semilocal.}

\begin{lem} \label{cohCT} Let $(M,q)$ be a faithful quadratic space over an \new LG ring $R$, \enew and let $0 \ne S\in \Ralg$. We consider the following maps obtained by base change,
\begin{equation}\label{cohCT-0} \begin{split}
 \frb_{\uO, S} \co H^1\fppf(R, \uO(q)) &\longto H^1\fppf(S, \uO(q_S)), \\
\frb_{\uSO, S} \co H^1\fppf(R, \uSO(q)) &\longto H^1\fppf(S, \uSO(q_S)).
\end{split}\end{equation}

\begin{inparaenum}[\rm (a)] \item\label{cohCT-a}
Then
\[ \frb_{\uO, S} \text{ injective } \quad \implies \quad
   \frb_{\uSO, S} \text{ injective. }
\]

\item \label{cohCT-baa} Suppose $R$ is a normal domain, and let $K$ be its fraction field. Then
\begin{equation} \label{cohCT-aa1}
 \frb_{\uO, K} \text{ injective } \quad \iff \quad
   \frb_{\uSO, K} \text{ injective. }
\end{equation}
\end{inparaenum}\end{lem}

\begin{proof} We will prove \eqref{cohCT-a} and \eqref{cohCT-baa} at the same time, where, obviously, it suffices to show the implication from right to left in \eqref{cohCT-baa}. For simpler notation we abbreviate $H^1 = H^1\fppf$. After decomposing $R$ in its product of connected components, we can in view of \ref{revLG}\eqref{revLG-aa} assume that $R$ is connected, so  that $(M,q)$ has constant rank $r \geq 1$. We distinguish  $r$ odd and $r$ even. \sm

\noindent{\it $r$ odd.} In case $r=1$ we have $\uSO(q) = \{*\}$, so we can suppose $r\ge 3$, although this does not simplify the proof. We know $\uO(q)= \uSO(q) \times \bmu_2$, see \eqref{orthsc-2}. Since $H^1(R, \cdot)$ respects direct products, i.e., $ H^1 (R, \uSO(q) \times \bmu_2) = H^1(R, \uSO(q)) \times H^1(R, \bmu_2)$ (\cite[III, \S4, 4.2]{DG}) and since this decomposition is respected by base change, $\frb_{\uSO,S}$ is injective whenever this holds for $\frb_{\uO,S}$.

To prove the converse, i.e., \eqref{cohCT-baa} for odd $r$, it suffices to show that the base change map $H^1(R, \bmu_{2,R}) \to H^1(S, \bmu_{2,S})$ is injective. After identifying $H^1(R, \bmu_{2,R})$ with the group $\Disc(R)$ of isomorphisms classes of discriminant modules \cite[III, (3.2.2)]{K}, the needed injectivity is \cite[III, (3.3.1)]{K}.
\sm

\noindent{\it $r$ even.} The usual twisting argument reduces the proof to establishing that the maps have trivial kernel as maps between pointed sets, see \ref{poiset}. Since $q$ is regular in the even rank case, the exact sequence
\eqref{sogsc2a} induces a commutative diagram of pointed sets
\begin{equation*}\label{diag_sn2}
\vcenter{
 \xymatrix{ 
\uO(q)(R) \ar[d] \ar[r]^{\Di} & \ZZ/2\ZZ(R) \ar[r] \ar[d] & H^1(R, \uSO(q)) \ar[d]^{\frb_{\uSO}} \ar[r] & H^1( R, \uO(q)) \ar[d]^{\frb_\uO} \\
\uO(q)(S) \ar[r]^{\Di} & \ZZ/2\ZZ(S) \ar[r] & H^1(S, \uSO(q_S)) \ar[r] & H^1( S, \uO(q_S))
}}
\end{equation*}
with exact rows, where $\frb_{\uO}$ is injective by assumption.
Let $x\in M$ with $q(x) \in R^\times$ (by Lemma~\ref{LGqdi} such an element exists because $M$ is faithfully projective over an LG ring $R$). Since  $\Di(\sigma_x)=1 \in \ZZ/2\ZZ=(\ZZ/2\ZZ)_R(R)$ by \eqref{refso2}, the map $\Di$ is surjective,
which by exactness of the top row implies $\Ima\big( \ZZ/2\ZZ(R) \to H^1(R, \uSO(q))\big) = \{*\}$ and then that
\begin{equation}  \label{cohCT-1}
\text{the map $H^1(R, \uSO(q)) \to H^1( R, \uO(q))$ has trivial kernel.}
\end{equation}
By commutativity of the right square, the map $\frb_{\uSO}$ then has
trivial kernel too.

It remains to prove \eqref{cohCT-baa} for even $r$. We use again the exact sequence \eqref{sogsc2a} and the corresponding cohomology sequence, but this  time extended one step to the right. We obtain the following commutative diagram with exact rows:
\begin{equation*} \vcenter{
 \xymatrix{ \ZZ/2\ZZ(R) \ar[d]_\al \ar[r] &  H^1(R, \uSO(q)) \ar[d]^{\frb_{\uSO,K}} \ar[r] & H^1( R, \uO(q)) \ar[d]^{\frb_{\uO,K}} \ar[r] & H^1(R, \ZZ/2\ZZ) \ar[d]^\de\\
\ZZ/2\ZZ(K) \ar[r]  & H^1(K, \uSO(q_K)) \ar[r] & H^1( K, \uO(q_K)) \ar[r] & H^1(K, \ZZ/2\ZZ)
}}.
\end{equation*}
Because $R$ is a domain, $\ZZ/2\ZZ(R) = \{0_R, 1_R\}$, which is sent by $\al$ to $\{0_K, 1_K\}$, so that $\al$ is a bijection. Also, $H^1(R, \ZZ/2\ZZ)$ is isomorphic to the group $\De(R)$ of isomorphism classes of separable quadratic $R$--algebras \cite[III, (4.1.4)]{K}. Under the isomorphism of the quoted result, $\de$ becomes the canonical map  $\De(R) \to \De(K)$, which is injective by \cite[III, (4.4.3)]{K}. We are now in the situation of the Four-Lemma \ref{poinlem}\eqref{poinlem-d}, and can conclude that injectivity of $\frb_{\uSO(q), K}$ implies that $\frb_{\uO(q), K}$ has trivial kernel. 
\end{proof} 

\subsection{Remarks.}\label{cohCTrm}  \begin{inparaenum}[(a)]
\item \label{cohCT-aa} Our proof of Lemma~\ref{cohCT} follows the proof of \cite[Prop.~1.2]{Collio},
where \ref{cohCT}\eqref{cohCT-baa} is proved for a semilocal normal domain with $2\in R\ti$ (note that then $q$ is regular). The same technique is used in the proof of \cite[Thm.~4]{Fedorov}, where $R$ is assumed to be regular local.
The argument has been generalized to the setting of unitary spaces over a semilocal regular domain $R$ with $2\in R\ti$ in \cite[Prop.~8.5]{First}.
\sm

\item \label{cohCTrm-a} The proof above shows that the equivalence
\begin{equation}
  \label{cohCTrm-a1} \frb_{\uO, S} \text{ injective } \quad \iff \quad
   \frb_{\uSO, S} \text{ injective }
\end{equation}
holds as soon as the following somewhat technical conditions are satisfied:

\begin{inparaenum}[(i)]
  \item \label{cohCTrm-ai} $H^1(R, \bmu_{2,R}) \to H^1(S, \bmu_{2,S})$ is injective,

  \item \label{cohCTrm-aii} $H^1(R, (\ZZ/2\ZZ)_R) \to H^1(S, (\ZZ/2\ZZ)_S)$ is injective, and

  \item \label{cohCTrm-aiii} $(\ZZ/2\ZZ)(R) \to (\ZZ/2\ZZ)(S)$ is surjective.
  \end{inparaenum}

\noindent
Due to Deligne's trace criterion  (\cite[XVII, 6.3.13--6.3.15]{SGA4}, see \cite[B.3, B.4]{GN-Sp}), both \eqref{cohCTrm-ai} and \eqref{cohCTrm-aii} hold if $S/R$ is finite locally free of odd rank. Moreover, \eqref{cohCTrm-aiii} holds whenever $S$ is connected. If $2\in R\ti$, then $\ZZ/2\ZZ \simlgr \bmu_2$ via the map $\chi$ of
\eqref{dickhomc}, so that \eqref{cohCTrm-ai} and \eqref{cohCTrm-aii} coincide. \sm

\item One is of course interested in knowing when $\frb_{\uO, S}$ and hence $\frb_{\uSO, S}$ are injective. Injectivity of $\frb_{\uO,S}$ is condition~\ref{equi-D} of \ref{equi}; conditions implying \ref{equi-D} are stated there. A classical result where \ref{equi-D} holds is \cite[Lem.~2.1]{Knebusch73}, where $R$ is a valuation ring, $S$ its field of fractions and $q$ is regular \cite[Lem.~2.1]{Knebusch73}. Other instances are mentioned in \cite[1.4]{Cesna}. Two more are presented in Corollaries~\ref{cor-cesna} and \ref{prop_odd}, see also Lemma~\ref{hypdes}.
\end{inparaenum}

\begin{cor}[\v{C}esnavi\v{c}ius]\label{cor-cesna} Let $R$ be an unramified regular local ring, let $K$ be its field of fractions and let $(M,q)$ be a faithful quadratic $R$--space. Then the maps $\frb_{\uO, K}$ and $\frb_{\uSO, K}$ of \eqref{cohCT-0} are injective. In particular, for a faithful quadratic $R$--space $(M', q')$ we have $(M',q')_K \cong (M,q)_K \iff (M',q') \cong (M,q)$.
\end{cor}

\begin{proof}
   This is proven in \cite[Cor.~9.6]{Cesna} under the additional assumption that $2\in R\ti$. The proof in loc.\ cit. has two steps: \begin{inparaenum}[(I)] \item \label{cor-cesnaI} A reduction that it suffices to show injectivity of $\frb_{\uSO_n, K}$, where $\uSO_n= \uSO(q_n)$, $q_n$ the split quadratic form of rank $n$. This is where the author uses $2\in R\ti$. \item \label{cor-cesnaII} Applying the general result \cite[Thm.~9.1]{Cesna}, valid for quasi-split reductive $R$--groups, to the split $R$--group $\uSO(q_n)$. This step does not need $2\in R\ti$.

   We follow this approach, but replace step \eqref{cor-cesnaI} by the equivalence \eqref{cohCT-aa1} and then use twisting to replace $\uSO(q)$ by $\uSO(q_n)$, as already done in \cite{Cesna}. \end{inparaenum}
\end{proof} 

Without going into details, we note that Corollary~\ref{cor-cesna} holds for more general rings $R$, namely for those for which \cite[Thm.~9.1]{Cesna} is true.

\begin{cor}[$S/R$ odd rank] \label{prop_odd} Let $(M,q)$ be a faithful quadratic space over a semilocal ring $R$, and let $S\in \Ralg$ be a finite $R$--algebra of {\em odd\/} degree which satisfies one of the following conditions,
\begin{enumerate}[label={\rm (\roman*)}]
 \item \label{prop_odd-u} $S$ is a one-generated $R$--algebra, or

\item \label{prop_odd-pri} $R$ satisfies the primitive criterion and  $\uPrim_R(S) \ne \emptyset$, or  

\item \label{prop_odd-aa} $R$ is semilocal and $S$ is an \'etale $R$--algebra, or

\item \label{prop_odd-b}  $S/R$ is an extension of fields. \end{enumerate}
Then the following hold. \sm

\begin{inparaenum}[\rm (a)] \item \label{prop_odd-a} If $(M',q')$ is a quadratic space over $R$, then $q_S \cong q'_S \iff q \cong q'$, i.e., condition \ref{equi-D} of\/ {\rm \ref{equi}} holds. \sm

\item \label{cohCT-b} Both maps $\frb_{\uO,S}$ and $\frb_{\uSO,S}$ of \eqref{cohCT-0} are injective.
\end{inparaenum}
\end{cor}

\begin{proof} The cases \ref{prop_odd-u}, \ref{prop_odd-pri} and \ref{prop_odd-aa} are restatements of Springer's Odd Degree Theorem~\ref{spo} and its Corollary~\ref{spo-LGG}. In case \ref{prop_odd-b} it suffices to write the finite field extension $S/R$ as a finite tower of odd one-generated field extensions. Thus, condition \ref{equi-Dreg} of \ref{equi} holds, see \ref{equisa}\eqref{equisa-b}. Moreover, by \ref{equisa}\eqref{equisa-c} also condition \ref{equi-D1} is fulfilled. Therefore we also have condition \ref{equi-D}.  Since $H^1(R, \uO(q))$ classifies $R$--forms of $q$, the injectivity claim of $\frb_{\uO,S}$ in \eqref{cohCT-b} is a consequence of \eqref{prop_odd-a}. Then injectivity of $\frb_{\uSO, S}$ follows from \ref{cohCT}\eqref{cohCT-a}.
\end{proof}

\begin{lem}  \label{hypdes} Let $R$ be a semilocal domain and let $K$ be its fraction field. Suppose the following condition \ref{hypdes-i} on quadratic $R$-spaces is satisfied:

\begin{enumerate}[label={\rm (\Alph*)$_\hyp$}]\setcounter{enumi}{3}

\item  \label{hypdes-i} If $(N,q)$ is a quadratic $R$--space and if $U$ is a finite projective $R$--module such that $(N,q)_K \cong \HH(U)_K$, then $(N,q) \cong \HH(U)$.
\end{enumerate}
Then \ref{equi-D} holds for\/ {\em regular} quadratic modules, and hence \ref{equi-D} holds in full generality if $2\in R\ti$.

\end{lem}

\begin{proof} We need to show: if $(M,q)$ is a regular quadratic module and $(M',q')$ is a quadratic space such that $(M',q')_K \cong (M,q)_K$, then $(M',q') \cong (M,q)$. Thus let $(M',q')$ be such a quadratic space. By Corollary~\ref{carhyp-cor}\eqref{carhyp-cor-a} we have $(M,q) \perp (M,-q) \cong \HH(M)$ because $(M,q)$ is regular. Hence
\begin{align*}
  \HH(M)_K &\cong \big( (M,q) \perp (M,-q)\big)_K \cong (M,q)_K \perp (M,-q)_K
  \\&  \cong (M,q)_K \perp (M', -q')_K \cong \big( (M,q) \perp (M', - q')\big)_K.
\end{align*}
Applying \ref{hypdes-i}, we get $(M,q) \perp (M',-q') \cong \HH(M)$ and then $(M,q) \cong (M', q')$ by Corollary~\ref{carhyp-cor}\eqref{carhyp-cor-b}. \end{proof}

\subsection{Remarks.} The proof of Lemma~\ref{hypdes} goes back to \cite[Prop.~1.2, (F) $\implies$ (D)]{Collio}, stated there for a semilocal normal $R$ with $2\in R\ti$. The corollary has been retaken in \cite[Thm.~3]{Fedorov}, where $R$ is assumed to be a regular local ring with $2\in R\ti$.
\sm

A similar local-global approach as in \ref{hypdes} is used in \cite{PS} for $R$ a noetherian domain of dimension $1$ with $2\in R\ti$ and a quadratic $R$--space $(M,q)$.
\begin{enumerate}[label={\rm (\roman*)}]
\item (\cite[Thm.~2.1]{PS}) If in addition $R$ is semilocal and if $(M,q)\ot_R R_\p$ is isotropic for all $\p \in \Spec(R)$, then $(M,q)$ is isotropic.

\item  (\cite[Thm.~3.1]{PS}) If the singular set $\mathrm{Sing}(R)$ is finite and non-empty and if $(M,q)\ot_R R_\p$ is isotropic for all $\p \in \mathrm{Sing}(R)$, then $(M,q)$ is isotropic.
\end{enumerate}

\comments{(2026-04-29) The review of $\uSpin(q)$ and $\mathbf{S\Ga}(q)$, the old part of (a) of Lemma~\ref{spinim}, and the preparation for the proof of Lemma~\ref{abcco} (spinor norm $SN$) are now in the section \S\ref{sec:orthgroup-LG} on orthogonal transformations. }

For the notation $\rmsn$ and $\spi$ used in the \ref{spinim}, see \ref{mgs}.

\begin{lem} \label{spinim} Let $(M,q)$ be a faithful quadratic $R$--space.

\begin{inparaenum}[\rm (a)]

 \item \label{spinim-b} We always have
  $\rmsn\big( \spi\me(\Refl^+(q))\big) = \rmD(q)^{[\rm ev]}$,
  hence $\rmD(q)^{\rm[ev]} \subset \Ima(\rmsn)$. Moreover, if\/ $\SO(q) = \Refl^+(q)$, then $\Ima(\rmsn) =  \rmD(q)^{[\rm ev]}$.
\sm

\item\label{spinim-c} We have  $\Ima(\rmsn) = \rmD(q)^{\rm [ev]}$ whenever

 \begin{inparaenum}[\rm (i)]
    \qquad \item\label{spinim-ci} $\SO(q) = \Refl^+(q)$, or

     \qquad \item \label{spinim-cii} $R$ is semilocal.
 \end{inparaenum}
\end{inparaenum}
\end{lem}

\begin{proof} 
\eqref{spinim-b} Let $x\in \SG(q)$ with $\spi(x) = \rho_{m_1} \cdots \rho_{m_n}$ and $n$ even. By exactness of \eqref{mgs1}, $x=u  m_1 \cdots m_n$ for some $u\in R\ti$. Then $\rmsn(x) = u^2 q(m_1) \cdots q(m_n) \in \rmD(q)^{\rm [ev]}$ because $R\ti{}^2 \subset \rmD(q)^{[2]}$ by \eqref{dqd-ele-a2}. Conversely, let $r=q(m_1) \cdots q(m_n) \allowbreak \in \rmD(q)^{\rm [ev]}$. Then $r\in \rmsn\big( \spi\me(\Refl^+(q))\big)$ by \eqref{spinim1}.
\sm

\eqref{spinim-c} If $\SO(q) = \Refl^+(q)$, then $\spi\me(\Refl^+(q)) = \SG(q)$, so $\rmD(q)^{\rm [ev]} = \Ima(\rmsn)$.

Suppose $R$ is semilocal. We apply Corollary~\ref{cor-knex} and \eqref{spinim-b} to get a finite \'etale extension $T\in \Ralg$ of constant odd degree $2d+1\ge 1$ such that $\Ima(\rmsn_T) = \rmD(q_T)^{\rm [ev]}$.  Let $a \in \Ima(\rmsn)$. Then $a=\rmN_{T/R}(a_T) a^{-2d}$.  Since $a_T \in \Ima(\rmsn_T) = \rmD(q_T)^{\rm [ev]}$, the Knebusch norm principle \eqref{thm-kneb-b2} shows that $\rmN_{S/R}(a_T) \in \rmD(q)^{\rm [ev]}$, hence $a \in \rmD(q)^{\rm [ev]}$. Thus
$\Ima(\rmsn) \subset \rmD(q)^{\rm [ev]}$. By \eqref{spinim-b} the other inclusion is always true. \end{proof}

\comments{(2020-03-03) I have isolated the cohomology part of Prop.~\ref{prop_sn}  in \ref{abcco} for future applications and improvements of \ref{prop_sn}.\sm

(2023-05-12) I have removed the previously required condition
\begin{enumerate}[label={\rm (\roman*)}]
 \item \label{accco-i}  $\Ima(\vphi) = \rmD(q)^{[\rm ev]}/R\ti{}^2$, cf.\ \eqref{abc3}, and $\Ima(\vphi') = \rmD(q_S)^{[\rm ev]}/S\ti{}^2$ where
    $\vphi' \co \uSO(q)(S) \to S\ti/S\ti{}^2$ is analogous to $\vphi$,
\end{enumerate}
in \ref{abcco}. This is now proven in \eqref{abc-c1}, which is possible because of the new part \ref{spinim}\eqref{spinim-cii}, which previously was part of Proposition~\ref{prop_sn}. }

\begin{lem}\label{abcco} Let $(M,q)$ be a faithful quadratic space over a semilocal ring $R$ and let $S\in \Ralg$ be semilocal too.
Assume furthermore that
\begin{enumerate}[label={\rm (\roman*)}]
 \item \label{acco-iin} the group homomorphism
     \[ \al \co R\ti / \rmD(q)^{[\rm ev]} \to S\ti/\rmD(q_S)^{[\rm ev]},
    \]
   induced by $R\ti \to S\ti$, $u \mapsto u \ot 1_S$, is injective, and that
\sm

 \item \label{acco-iii} the map  $  \frb_{\uSO} \co H^1\fppf(R, \uSO(q)) \longto H^1\fppf (S, \uSO(q_S))$, obtained by base change, is injective.
\end{enumerate}
Then the analogous map
\begin{equation} \label{abcco-1}
   \frb_{\uSpin} \co H^1\fppf(R, \uSpin(q)) \longto H^1\fppf (S, \uSpin(q_S))
\end{equation}
is injective too.
\end{lem}

\begin{proof} We abbreviate $H^1 = H^1\fppf$. We consider the commutative diagram of pointed sets
\begin{equation*} 
\vcenter{
 \xymatrix{ 
\uSO(q)(R) \ar[d] \ar[r]^{\varphi} & R^\times/ (R^\times)^2 \ar[r] \ar[d] & H^1(R, \uSpin(q)) \ar[d]^{\frb_{\uSpin}} \ar[r] &
H^1( R, \uSO(q)) \ar[d]^{\frb_{\uSO}}  \\
\uSO(q)(S)  \ar[r]^{\varphi'} & S^\times/ (S^\times)^2 \ar[r]  & H^1(S, \uSpin(q_S))  \ar[r] & H^1(S ,\uSO(q_S))
}}
\end{equation*}
in which each row is part of the long cohomology sequence associated with the exact sequence $1 \to \bmu_2 \to \uSpin(q) \to \uSO(q)\to 1$, see \eqref{mgs3}, and where we used \eqref{abs2}.
By 
\eqref{abc-c1}, we know that $\Ima(\vphi) =  \rmD(q)^{[\rm ev]}/R\ti{}^2$ is  a subgroup of $R\ti/R\ti{}^2$, allowing us to form
\[ (R\ti/ R\ti{}^2) \big/ \Ima(\vphi) \cong R\ti/ \rmD(q)^{[\rm ev]} \]
and to replace the diagram above by the commutative diagram  of pointed sets with exact rows
\begin{equation*} 
\vcenter{
 \xymatrix{
1 \ar[r] &R^\times/ \rmD(q)^{[\rm ev]} \ar[r]^\psi \ar[d]^\al & H^1(R, \uSpin(q)) \ar[d]^{\frb_{\uSpin}} \ar[r] & H^1( R, \uSO(q)) \ar[d]^{\frb_{\uSO}}  \\
1 \ar[r] &S^\times/\rmD(q_S)^{[\rm ev]} \ar[r]^{\psi'}  & H^1(S, \uSpin(q_S))  \ar[r] & H^1(S ,\uSO(q_S)) } }.
\end{equation*}
Since ${\frb_{\uSO}}$ is an injective set map by assumption~\ref{acco-iii}, it is straightforward to show that $\psi$ induces a surjective map $\Ker(\al) \to \Ker({\frb_{\uSpin}})$ of pointed sets. 
It is therefore enough to show that $\Ker(\al)$ is trivial. But this is condition~\ref{acco-iin}.\end{proof}
\sm

The condition~\ref{abcco}\ref{acco-iin} is not always fulfilled. For example, let $R=\RR$, let $S = \CC$, and let $q$ be a positive-definition quadratic form of positive dimension. Then $\rmD(q) = \{ r\in R : r>0\}= \rmD(q)^{[\rm ev]}$, but $\rmD(q_\CC) = \CC\ti= \rmD(q_\CC)^{[\rm ev]}$ since either $q_\CC \cong \lan 1\ran$ or $q_\CC$ contains a hyperbolic plane,  in which case we can apply \ref{lem_smooth_locus_exam}\eqref{lem_smooth_locus_hypn}. Hence $\Ker (\al) \ne 0$. Note that in this example $[S:R]$ is even. On the other side, we will establish injectivity of $\al$ in case $[S:R]$ is odd in the proof of Proposition~\ref{prop_sn}\eqref{prop_sn-b}. It generalizes a result of Black \cite[Thm.~3.7]{Black} from fields of characteristic $\ne 2$ to semilocal rings in the case of spin groups.

\comments{(2023-05-13) Is the Grothendieck-Serre conjecture known to hold for Spin groups in our setting, i.e., $R$ a local regular ring, $S$ its field of fractions, $(M,q)$ faithful quadratic $R$--space, then $ \co H^1\fppf(R, \uSpin(q)) \longto H^1\fppf (S, \uSpin(q_S))$ is injective? \sm

By \ref{abcco} and Corollary~\ref{cor-cesna} (essentially \v{C}esnavi\v{c}ius), we only need to show \ref{abcco}\ref{acco-iin}, i.e.,
\[
   u \in R\ti, u_S \in \rmD(q_S)^{[ev]} \implies u \in \rmD(q)^{[ev]} .
\]
which is in itself an interesting question. However, unless there is a quick solution or a reference, I do not to embark on this for now.
\sm

(2021-03-18) The proof of the Grothendieck-Serre conjecture in \cite{OPZ} uses the approach of \ref{abcco} and verifies the conditions \ref{acco-iin} and \ref{acco-iii} of \ref{abcco}. Apparently, condition \ref{acco-iii} is proven in Ojanguren's 1980 paper on "Quadratic forms over regular rings". \sm

The main part of \cite[Thm.~4.1]{OPZ} (= Grothendieck-Serre Conjecture) is the proof of condition \ref{abcco}\ref{acco-iin}. The proof is not given in detail. Rather the authors refer to the 2001-paper by Kirill on ``Purity of functors with transfer''. It is interesting that the proof of \ref{abcco}\ref{acco-iin} in \cite{OPZ} uses the weak version of Knebusch's norm principle.}

\pcomments{(2026-05-13)  Le noyau trivial est connu dans le cas o\`u q est isotrope car c'est une cons\'equence facile du cas orthogonal. La conjecture est connue dans le cas o\`u l'anneau local contient un corps
(Fedorov-Panin), mais pas en g\'en\'eral. 

Si q est anisotrope, c'est exactement le propos de [OPZ] d'\'etablir pois d'utiliser le principe de norme de Knebusch pour d\'emontrer Grothendieck-Serre pour les groupes de spineurs. 

La semaine prochaine, on va essayer de suivre la m\^eme d\'emarche et obtenir peut-\^etre un r\'esultat important.  }

\begin{prop} \label{prop_sn} Let $(M,q)$ be a faithful quadratic $R$--space, where $R$ is a semilocal ring. \sm

\begin{inparaenum}[\rm (a)]
\item\label{prop_sn-a} If $S\in \Ralg$ is finite projective of odd rank, the condition {\rm \ref{abcco}\ref{acco-iin}} follows from the condition that
\end{inparaenum}
\begin{enumerate}[label={\rm (\Roman*)}] 
 \item \label{acco-ii} the Knebusch norm principle holds in the form $\rmN_{S/R}\big( \rmD(q_S)\big) \subset \rmD(q)^{[\rm ev]}$.
\end{enumerate}
\sm

\begin{inparaenum}[\rm (a)]\setcounter{enumi}{1}

\item\label{prop_sn-b} Let $S\in \Ralg$ be a finite \'etale extension of constant odd degree. Then the map
      \[ \frb_{\uSpin} \co H^1\fppf(R, \uSpin(q)) \longto H^1\fppf (S, \uSpin(q_S)) \]
 of \eqref{abcco-1} is injective.
\end{inparaenum} \end{prop}

\begin{proof} \eqref{prop_sn-a}
To prove that \ref{abcco}\ref{acco-iin} is implied by \ref{acco-ii}, let $a \in R^\times$ such that $a_S \in \rmD(q_S)^{[\rm ev]}$ and let $S/R$ have degree $d$. Because of assumption~\ref{acco-ii}, we have $a^d= {\rmN}_{S/R}(a_S) \in \rmD(q)^{[\rm ev]}$. Since $d$ is odd by assumption, it follows that $a \in \rmD(q)^{[\rm ev]}$. Thus $\Ker(\al) = \{\star\}$, and we are done.
\sm

\eqref{prop_sn-b} We verify the conditions \ref{acco-ii} above as well as
\ref{acco-iii} of Lemma~\ref{abcco}.
The map  $\frb_{\uSO} \co H^1(R, \uSO(q)) \to H^1(S, \uSO(q_S))$ is injective by   Corollary~\ref{prop_odd}\eqref{cohCT-b}. Thus \ref{acco-iii} holds.

For \ref{acco-ii} we note that changing $q$ by a scalar $u\in R\ti$ does not change the groups involved. Hence, as pointed out in the remarks after Theorem~\ref{thm-kneb}, we can assume that $\uS_{q,1}\rmsm (R) \ne \emptyset$. Then the Knebusch norm principle in the form \ref{acco-ii} holds by \eqref{thm-kneb-b1}.
\end{proof}


\subsection{The closed immersion $\mathbf{S\Ga}(q) \to \uGL_1((\Cli(q)_0))$} \label{cim} We consider again the setting of \ref{mgs}, so $R$ is arbitrary, $(M,q)$ is a faithful quadratic $R$--space, $\Cli= \Cli(q)$ is the associated Clifford algebra and $\Cli = \Cli_0 \oplus \Cli_1$ is the canonical $\ZZ/2\ZZ$--grading of $\Cli$. We abbreviate
\[ A = \Cli(q)_0, \]
and note that $A$ is a separable and finite projective $R$--algebra. Indeed, if $M$ has constant even rank, $A$ is an Azumaya algebra over its centre $\Dis(q)$ and if $M$ has odd rank, then $A$ is an Azumaya $R$--algebra. We will use the group homomorphism
\[ \sint \co A \ti \to \Aut(A), \quad x \mapsto i_x|_A \]
where $i_x$ is the conjugation of $x$ on $\Cli(q)$.

Any $g\in \orth(q)$ induces an automorphism $\Cli(g)$ of the Clifford algebra $\Cli$, uniquely determined by $\Cli(g)(m) = g(m)$ for $m\in M$. It stabilizes the $\ZZ/2\ZZ$-grading of $\Cli$, and thus gives rise to a group homomorphism
\[ \rho \co \SO(q) \to \Aut(A), \quad x \mapsto \Cli(g)|_A \]
which makes the diagram below commutative:
\begin{equation} \label{cim0} \vcenter{
  \xymatrix@C=40pt{ 1  \ar[r] & R \ti\ar@{=}[d] \ar[r]^\inc &
      \SG(q) \ar[d]^\inc \ar[r]^\spi &\SO(q)\ar[d]^\rho \\
  1  \ar[r] & R \ti \ar[r]^\inc & A\ti \ar[r]^\sint &\Aut(A)
    }}
\end{equation}
where $\inc$ is inclusion. Indeed, let $x\in \SG(q)$. Then $\Cli\big(\spi(x)\big)$ is an automorphism of $\Cli$ satisfying $\Cli\big(\spi(x)\big)(m) = \spi(x)(m) = x mx\me$ for all $m\in M$ and therefore $\Cli\big(\spi(x)\big)(c) = xcx\me$ for all $c\in \Cli$, in particular $\rho\big(\spi(x)\big) = \sint(x)$ holds. The commutativity of the left square is clear.

The construction of $\rho$ commutes with base change and so canonically extends to a homomorphism $ \rho \co \underline{\SO}(q) \to \underline{\Aut}(A)$ of $R$--group functors, hence to a homomorphism $\rho \co \uSO(q) \to \uAut(A)$ of the associated $R$--group schemes. Let $\uGL_1(A)$ be the $R$--group scheme representing the $R$--group functor $R' \mapsto (A\ot_R R')\ti$; it is a reductive $R$--group scheme \cite[3.1.0.50]{CF}. The diagram \eqref{cim0} extends to a commutative diagram of $R$--group schemes
\begin{equation} \label{cim1} \vcenter{
  \xymatrix@C=35pt{ 1  \ar[r] & \GG_m \ar@{=}[d] \ar[r]^\inc &
      \uSG(q) \ar[d]^\inc \ar[r]^\spi &\uSO(q)\ar[d]^\rho \ar[r] & 1\\
  1  \ar[r] & \GG_m  \ar[r]^\inc & \uGL_1(A) \ar[r]^\sint &\uAut(A)
    }}
\end{equation}
where $\inc$ and $\sint$ are the obvious group homomorphisms.
The top row is exact by \ref{mgs}.
We know that $\uSG(q)$ and $\uGL_1(A)$ are reductive group schemes. It then follows from \cite[Thm.~5.3.5]{Co1} that
\begin{equation} \label{cim-3}
\text{\em the inclusion $\inc\co \uSG(q) \to \uGL_1(A)$ is a closed immersion.}
\end{equation}
The following Lemma~\ref{cimc} shows that this also holds for $\rho$ in the odd rank case.
%
%

\begin{lem}  \label{cimc} We continue with the setting of {\rm \ref{cim}}, but let $(M,q)$ be a quadratic $R$--space of constant {\em odd} rank. \sm

\begin{enumerate}[label={\rm (\alph*)}]
\item \label{cimc-a} Then the group scheme homomorphism $\rho \co \uSO(q) \to \uAut(A)$ of \eqref{cim1} is a closed immersion.
\sm

\item \label{cimc-b}  If $M$ has rank $3$, then $A$ is a quaternion algebra, and \sm
\begin{enumerate}[label={\rm (\roman*)}]
  \item \label{cimc-bi} $\rmsn = \Nrd$, where $\Nrd$ is the reduced norm of the Azumaya algebra $A$, cf.\ {\rm \ref{rcp}}, \sm

  \item  \label{cimc-bii} $\uSG(A) = \uGL_1(A)$ and $\uSpin(q) = \Ker(\Nrd) = \uSL_1(A)$, \sm

  \item \label{cimc-biii} $\rho$ is an isomorphism: $\uSO(q) \cong \uAut(A)$.
\end{enumerate}\end{enumerate}
\end{lem}

\begin{proof} \eqref{cimc-a}
The assumption on $M$ implies that $A$ is an Azumaya $R$--algebra. Hence $\uGL_1(A) \to \uAut(A)$ is surjective in the flat topology. The diagram \eqref{cim1} can therefore be augmented to the commutative diagram
\begin{equation} \label{cimc1} \vcenter{
  \xymatrix@C=35pt{ 1  \ar[r] & \GG_m \ar@{=}[d] \ar[r]^\inc &
      \uSG(q) \ar[d]^\inc \ar[r]^\spi &\uSO(q)\ar[d]^\rho \ar[r] & 1\\
  1  \ar[r] & \GG_m  \ar[r]^\inc & \uGL_1(A) \ar[r]^\sint &\uAut(A) \ar[r] & 1
    }}
\end{equation}
with exact rows. It follows that $\rho$ is a monomorphism. We are again in a setting where we can apply  \cite[Thm.~5.3.5]{Co1} and thus conclude that $\rho$ is a closed immersion.
\sm

\ref{cimc-b} The algebra $A$ is an Azumaya $R$--algebra of constant rank $4$, i.e., is a quaternion algebra. Part \ref{cimc-bi} is proven in \cite[V, (3.1.1)]{K}, and part \ref{cimc-bii} is formula (3.1.1) in \cite[V, (3.3)]{K}, and then \ref{cimc-biii} follows from \eqref{cimc1}. \end{proof}

\subsection{Remarks on \ref{cimc}\ref{cimc-b}.} \label{cimrem} \begin{inparaenum}[(a)]
  \item \label{cimc-a} Let again $(M,q)$ be a quadratic space with $M$ of constant rank $3$. Then
     $A$ is a quaternion algebra and $ \uPGL_1(A)\simlgr \uAut(A)$ is a semisimple adjoint $R$--group scheme of type $\rmA_1$ by \cite[3.5.0.82]{CF}.  Also, $\uSO(q)$ is a semisimple adjoint $R$--group scheme of type $\rmB_1$,  \cite[6.4.0.30]{CF}. Since $\rmA_1 \cong \rmB_1$ as semisimple root data, the Isomorphism Theorem of semisimple group schemes \cite[XXIII, 2.3]{SGA3} says that $\uSO(q)$ and $\uAut(A)$ are isomorphic group schemes. The point of Lemma~\ref{cimc}\ref{cimc-b} is that it gives a concrete isomorphism. \sm

\item In \ref{cimc}\ref{cimc-b} we start with a quadratic space $(M,q)$ of constant rank $3$ and define a quaternion algebra $A$ in terms of $q$. Conversely, one can start with a quaternion algebra $Q$ with norm $n_Q$ and trace $\tr_Q$. Put $Q_\pure = \{ a\in Q : \tr_Q(a) = 0 \}$. Then $(Q_\pure , n_Q|_{Q_\pure})$ is a quadratic space of constant rank $3$. The quaternion algebra associated with this quadratic space is isomorphic to $A$ \cite[V, (3.2.4)]{K}.
    \end{inparaenum}

\comments{Deleted the following since it does not fit into the style of the memoirs (the main body is over rings); made a remark that this could be added in the intro.

{\sf Let $S$ be a scheme. It is clear to the experts that \ref{cim} and \ref{cimc} hold for quadratic spaces over $S$ and the corresponding $S$--group schemes. Since we only deal with quadratic forms over $S=\Spec(R)$, even mostly with $R$ semilocal, we have not pursued this generality here.}}

\subsection{Example $(M,q)$ of rank $3$ again}\label{cra} We continue with the setting of \ref{cimc}\ref{cimc-b}: $R$ is arbitrary and $(M,q)$ is a quadratic $R$--space of constant rank $3$. Thus the commutative diagram \eqref{cimc1} becomes
\begin{equation} \label{cra1} \vcenter{
  \xymatrix@C=35pt{ 1  \ar[r] & \GG_m \ar@{=}[d] \ar[r]^\inc &
      \uSG(q) \ar@{=}[d] \ar[r]^\spi &\uSO(q)\ar[d]^\rho_{\cong} \ar[r] & 1\\
  1  \ar[r] & \GG_m  \ar[r]^\inc & \uGL_1(A) \ar[r]^\sint &\uAut(A) \ar[r] & 1
    } }.
\end{equation}
By \eqref{mgs3} the homomorphism $\rho \circ \spi$ is still an epimorphism (in the flat topology) when restricted to the subgroup scheme
\[ \uSL_1(A) = \Ker(\Nrd) = \uSpin(q)\]
of $\uGL_1(A)$. Hence we get the commutative diagram
\[ \vcenter{
\xymatrix@C=40pt{
& 1 \ar[d] & 1 \ar[d] \\
 1 \ar[r] & \bmu_2  \ar[r] \ar[d]  & \uSL_1(A) \ar[r]^{\rho\, \circ \, \spi} \ar[d]^{\inc}& \uAut(A) \ar[r] \ar[d]^{\cong}& 1 \\
 1 \ar[r] & \GG_m \ar[d]^{\times 2} \ar[r] & \uGL_1(A) \ar[d]^{\Nrd} \ar[r]^{\sint} & \uAut(A) \ar[r] & 1 \\
& \GG_m \ar[d]  \ar[r]^{\sim} & \GG_m \ar[d]  \\
& 1  & 1
}
}\quad .\]
As in \ref{abc} we let $\vphi$ and $\de$ be the characteristic maps associated with the exact sequences $1 \to \bmu_2 \to \uSL_1(A) \to \uAut(A) \to 1$ and $1 \to \bmu_2 \to \GG_m \to \GG_m \to 1$ respectively. Lemma~\ref{lem_snake}\ref{lem_snake4}
says that the square in the diagram below
\[ \xymatrix@C=50pt{
   &  A\ti \ar[r]^{\sint}  \ar[d]_{\Nrd=\rmsn} & \Aut(A) \ar[d]^\vphi \\
  R\ti \ar[r]^{\times 2} & R\ti \ar[r]^{\de} & H^1\fppf(R, \bmu_2)
 } \]
anticommutes. Observe $\de\big(\Nrd(A\ti)\big) = \Nrd(A\ti)/R\ti{}^2$. If we suppose in addition that $R$ is semilocal, then Proposition~\ref{prop_sn}\eqref{prop_sn-a} yields the corollary below.

\begin{cor} \label{craco} Let $R$ be a semilocal ring, let $(M,q)$ be a quadratic $R$--space of rank $3$ and let $A = \Cli_0(q)$ be the associated quaternion $R$--algebra. Then $\Nrd(A\ti) = \rmsn(A\ti) = \rmD(q)^{\rm [ev]}$.
\end{cor}

\newpage 

\section{Reduced norms and traces for Azumaya algebras}\label{sec:redNorm}

\comments{(2026-03-29) The section \ref{sec:redNorm} is needed for the section \ref{sec:normgroups}, but nor before. In the previous version of the Knebusch notes, section \ref{sec:redNorm} was appendix D. I put it here to diminish the number of appendices.}

The purpose of this section is to generalize some known results on reduced characteristic  polynomials, norms and traces for Azumaya algebras from the field case, see e.g. \cite[\S 17]{BA8}, to the ring setting. Throughout $A$ is an Azumaya $R$--algebra.

\subsection{Separable algebras}\label{sep}
Let $B$ be a unital associative $R$--algebra. We denote by $B\op$ its opposite algebra and observe that $B$ becomes  $B \ot_R B\op$--module with respect to the $B\ot_R B\op$--action given by $(b_1 \ot b_2) \cdot b = b_1bb_2$. One calls $B$ separable if this makes $B$ a projective $B\ot_R B\op$--module.
Standard references for separable algebras are \cite{Ford} and \cite[III]{KO}, where the reader can find more characterizations of separable $R$--algebras.
We list some facts that we will be using; throughout $B$ is a separable $R$--algebra. We point out that the $R$--module underlying $B$ is not assumed to be finite projective here, contrary to what is required in \cite{CF}. \sm

\begin{inparaenum}[(a)]
  \item \label{azu-sep} ({\em Separable versus Azumaya}) A unital associative $R$--algebra $A$ is a separable $R$--algebra if and only if $A$ is  an Azumaya algebra over its centre $\rmZ(A)$ and $\rmZ(A)$ is a separable $R$--algebra (\cite[II, 3.8]{DI}, or see \cite[7.1.11]{Ford}, \cite[III, 5.5]{KO}).  \sm

\item\label{seb-bc} ({\em Base change}) The $S$--algebra $B \ot_R S$ is separable for any $S\in \Ralg$ \cite[4.3.2]{Ford}. \sm

   \item \label{sep-fap} ({\em Centre $\rmZ(B)$}) If the underlying $R$--module of $B$ is faithfully projective, then $\rmZ(B)$ is a faithfully projective $R$--module and hence a finite \'etale $R$--algebra \cite[4.4.6]{Ford}. \sm

   \item \label{sep-weil} ({\em $B$--modules}) Any exact sequence of $B$--modules, which is split-exact as sequence of $R$--modules, is also split exact as sequence of $B$--modules. Any $B$--module, which is projective as $R$--module, is also projective as $B$--module \cite[4.4.1]{Ford}. \sm

   \item \label{sep-li} ({\em Left ideals}) Let $L \subset B$ be a left ideal such that $B/L$ is projective as $R$--module. Then there exists an idempotent $c\in B$ such that $L = Bc$ and hence $L$ and $B/L \cong B (1_B -c)$ are projective as $B$--module.

       Indeed,  the exact sequence $0 \to L \to B \to B/L \to 0$ of $B$--modules is split-exact as exact sequence of $R$--modules. Hence, by \ref{sep-weil}, 
       it is split-exact as sequence of $B$--modules. Thus, there exists a left ideal $L' \subset B$ such that $L \oplus L' = B$. Decomposing $1_B = c + c'$ with $c\in L$, $c'\in L'$, one finds that $c$ is an idempotent with $L = Bc$ and $L' = B(1_B - c)$ \cite[1.1.20]{Ford}.

       The idempotent $c$ is in general not uniquely determined by $L$.%
\sm

\item\label{sep-gr} ({\em Groups}) If $B$ is finite projective as $R$--module, the $R$--group scheme $\uGL_1(B)$, defined in \ref{trno}\eqref{trno-ds}, is reductive \cite[3.1.0.50]{CF}.
\end{inparaenum}


\subsection{Azumaya algebras} \label{azu}
Recall that an $R$--algebra $A$ is an Azumaya $R$--algebra if it is a separable $R$--algebra whose centre $\rmZ(A)$ is isomorphic to $R 1_A$ under the map $R \to \rmZ(A)$, $r \mapsto r 1_A$, where $1_A$ is the identity element of $A$.
Some references for Azumaya algebras over rings are  \cite{DI}, \cite{Ford}, \cite{GroI},  \cite[III, \S5]{K}, \cite[Ch.~III and IV]{KO} and \cite{Salt}. 
Azumaya algebras over fields are for example considered in \cite{BA8} and \cite{GS}.

We review some facts used later.  We denote by $\rank_R A$ the rank of the faithfully projective $R$--module underlying the Azumaya $R$--algebra $A$ and by $\deg_R A$ the {\em degree} of $A$, defined by $(\deg_R A)^2 = \rank_R A$. Both are locally constant functions on $\Spec(R)$ with values in $\NN_+$. \sm

\begin{inparaenum}[(a)]

\item\label{azu-bc} ({\em Base change}) Let $S\in \Ralg$. Then the tensor product algebra $A_S = A\ot_R S$ is an Azumaya $S$--algebra \cite[7.1.9]{Ford}. \sm

\item\label{azu-dir} ({\em Direct products}) Let $R=R_1 \times \cdots \times R_n$ be a direct product of rings and let $A_i$, $i=1,\ldots, n$, be $R_i$--algebras. Then $A$ is a separable (an Azumaya) $R$--algebra if and only if every $A_i$ is a separable (an Azumaya) $R_i$--algebra. Using the rank decomposition of finite projective $R$--module, this allows to reduced questions on general Azumaya algebras to those of constant rank. \sm

\item \label{azu-cdt} ({\em Double Centralizer Theorem} \cite[7.2.2]{Ford})
We denote by $\rmC_A(B) = \{ a\in A: ab = ba \; \text{ for all $b\in B$}\}$ the {\em centralizer subalgebra\/} of an $R$--subalgebra $B$ of $A$. If $B$ is separable, then so is $B' = \rmC_A(B)$ and $\rmC_A(B') =B$ holds. Moreover, if $B$ is an Azumaya $R$--algebra, then so is $B'$ and the multiplication of $A$ gives rise to an isomorphism $B \ot_R B' \simlgr A$ of $R$--algebras.
\sm

\item\label{azu-ex} ({\em Examples}) Let $P$ be a faithfully projective $R$--module. Then $\End_R(P)$ is an Azumaya $R$--algebra.

    More generally, let $Q$ be a faithfully projective left $A$--module (called a ``progenerator'' in \cite{Ford}). The restriction $P= \sfR_{A/R}(Q)$ is a faithfully projective $R$--module \cite[1.1.8]{Ford}, hence $\End_R(P)$ is an
Azumaya $R$--algebra. It contains the $R$--subalgebra $L_A=\{L_a: a \in A\}$, where $L_a \co A \to A$, $x \mapsto ax$, is the left multiplication of $A$ by $a$.  It is isomorphic to $A$. By \eqref{azu-cdt}, the centralizer algebra $\rmC_{\End_R(P)}(L_A) = \End_A(Q)$ {\em is an Azumaya $R$--algebra.}
\sm

\item\label{azu-tensor} ({\em Tensor products}) If $A$ and $B$ are Azumaya $R$--algebras, the tensor product algebra $A\ot_R B$ is an Azumaya algebra \cite[7.1.3]{Ford}. In particular, the matrix algebra $\Mat_n(A) = \Mat_n(R) \ot_R A$, $n \ge 1$, is an Azumaya $R$--algebra.%
    \sm

\item \label{azu-hom} ({\em Isomorphisms}) Let $B$ be a second Azumaya $R$--algebra. Any $R$--algebra homomorphism $f\co A \to B$ is injective; it is an isomorphisms if and only if $\rank_R A = \rank_R B$ \cite[7.6.1]{Ford}.
\sm

\item \label{azu-spli} ({\em Neutralizing and splitting rings})
\new
We say that $A$ is {\em neutral\/} if $A\cong \End_R(M)$ for some faithfully projective $R$--module $M$. 
\enew 
Following \cite[III, \S6]{KO} we call $S\in \Ralg$ a {\em neutralizing ring\/} if
if there exists a faithfully  projective $S$--module $M$ and an $S$--algebra isomorphism $\al \co A \ot_R S \simlgr \End_S(M)$. In this case, we will call $\al$ a {\em neutralization of $A$\/} and say $A\ot_R S$ is {\em neutral}. Thus $A$ is neutral if and only if its Brauer class is the neutral element of $\Br(R)$.
If $S$ is a neutralizing ring of $A$ that is faithfully flat as $R$--module, we will refer to $S$ and $\al$ as a {\em faithfully flat neutralizing ring} and {\em faithfully flat neutralization\/} respectively.
Every Azumaya algebra $A$ has a faithfully flat neutralizing ring, even a neutralizing ring that is finite \'etale as $R$--algebra  \cite[10.3.9]{Ford}, \cite[III, 6.6]{KO}.

We say that $A$ is {\em split\/} if $A\cong \Mat_n(R)$ for some $n\in \NN$. In this case, $A$ has constant rank $n^2$. Conversely, if $A$ has constant rank, necessarily of the form $n^2$ for some $n\in \NN_+$, there  exists a faithfully flat \'etale $S\in \Ralg$ and an $S$--algebra isomorphisms $\al \co A \ot_R S \simlgr \Mat_n(S)$ \cite[10.3.10]{Ford}, here referred to as a {\em splitting ring} and a {\em splitting}.

The reader should be warned that in \cite{Ford,Salt} and several other references a neutralization as defined above is called a ``splitting''. Our terminology follows the general practice that a split object in a category should be  unique up to isomorphism, which neutralizations need not be. \sm

\item\label{azu-com} ({\em Commutative subalgebras}) Let $S\subset A$ be a commutative $R$--subalgebra. Viewing $A$ as a right $S$--module, the sandwich map $A\ot_R A\op \simlgr \End_R(A)$, $a\ot b \mapsto L_a R_b = (x \mapsto axb)$ restricts to a unital $S$--algebra homomorphism
\begin{equation} \label{azu-com1}
   \al_S \co A \ot_R S \longto \End_S(A), \quad a \ot s \mapsto L_a R_s.
\end{equation}
If $A$ is faithfully projective as $S$--module, $\al_S$ is injective by \eqref{azu-ex} and \eqref{azu-hom}. \sm

\item\label{azu-mcs} ({\em Maximal commutative subalgebras}) An $R$--subalgebra $S$ of $A$ is a {\em maximal commutative\/} subalgebra if it is maximal with respect to inclusion among the commutative $R$--subalgebras of $A$. Maximal commutative subalgebras exist by Zorn. \sm

\begin{inparaenum}[(I)] \item A commutative $R$--subalgebra $S$ of $A$ is maximal commutative if and only if $S=\rmC_A(S)$. Another criterion:
\begin{equation}\label{azu-com2}
 \text{\em $S$ is maximal commutative}  \quad \iff\quad
   \text{\em $\al_S$ is an isomorphism.}
\end{equation}
If $S$ is maximal commutative, then $\al_S$ is an isomorphism, see for example  \cite[7.4.2]{Ford} or \cite[III, 6.1]{KO}. Conversely, if $\al_S$ is an isomorphism, then the centre of $\End_S(A)$ is $\{s\Id_A: s\in S\} \cong S$ since $A\ot_R S$ is an Azumaya $S$--algebra by \eqref{azu-bc}. If $S'\supset S$ is a commutative subalgebra of $A$, then $R_{s'}$, $s'\in S$, commutes with all $L_aR_s$ for $a\in A$ and $s\in S$. Therefore $R_{s'} = R_s$ for some $s\in S$, forcing $s'=s\in S$.

We note that even when $\al_S$ is an isomorphism, it need not be a neutralization since it is not clear that the $S$--module $A$ is projective.  Details of the following results are given in \cite[7.4.2]{Ford} or \cite[III, 6.1]{KO}.\sm

\item \label{azu-coIIi} {\em Assume that $S$ is maximal commutative and that the $S$--module $A$ is projective}. It is then faithfully projective, hence $\al_S$ is a neutralization. The neutralizing ring $S$ is then even faithfully projective as $R$--module, in particular it is a faithfully flat neutralizing ring. Moreover, $S$ is a direct summand of the $S$--module $A$ (\cite[III, 1.9]{KO}).

    Denoting by $\rank_R M \co\Spec(R) \to \ZZ$ the rank function of a finite projective $R$--module $M$ and by $\vphi \co \Spec(S) \to \Spec(R)$ the canonical continuous map, which is surjective by faithful flatness, the following equalities hold:
\begin{equation}  \label{azu-2} \begin{split}
 (\rank_R A)  \circ \vphi &= (\rank_S A)^2, \qquad (\rank_R S) \circ \vphi = \rank_S A, \\
 \rank_R A &= (\rank_R S)^2.
\end{split}\end{equation}
The first two equations are proved in the quoted references; they imply the third.
The equations \eqref{azu-2} imply: {\em $\rank_R A$, $\rank_R S$ and $\rank_S A$ are constant, as soon as one of these three functions is constant.}
\lv{
This follows from the criterion: if $f \co \Spec(R) \to \ZZ$ is a function, then $f\circ \vphi $ is constant iff $f$ is constant. }

\inparcom{The equations \eqref{azu-2} are stated without $\vphi$ in \cite[III, 6.1]{KO}, which does not make sense. They are retaken in \cite[7.4.2]{Ford} as functions of $\Spec(R)$ with a condition that disappears when considered as functions of $\Spec(S)$. }
\end{inparaenum} \sm

\item\label{azu-groups} ({\em Groups}) By \ref{sep}\eqref{sep-gr}, $\uGL_1(A)$ is a reductive $R$--group scheme. We abbreviate $\uGL_n(A) = \uGL_1\big( \Mat_n(A)\big)$, and denote by $\uAut(A)$ the automorphism group scheme of $A$, cf.\ \ref{ag}\eqref{ag-ex}. It is related to $\uGL_1(A)$ by the sequence
\begin{equation}
  \label{azu-groups1}
  1 \to \GG_m \xrightarrow{\;\ze\;} \uGL_1(A) \xrightarrow{\;\Int\;} \uAut(A) \to 1 \end{equation}
where $\ze$ identifies $\GG_m$ with the centre of $\uGL_1(A)$ and where $\Int$ sends $u\in A_S\ti$ to the inner automorphism $\Int(u)$ of $A_S$. Every automorphism of an Azumaya $R$--algebra with $\Pic(R) = 0$ is inner \cite[7.8.14]{Ford}, in particular this holds for a local ring $R$ \cite[7.8.15]{Ford}. Hence, the sequence \eqref{azu-groups1} is exact in the Zariski topology by \ref{surlem}, and then a fortiori also in the \'etale or flat topology \cite[3.2.0.58, 3.5.0.93(5)]{CF}.

If $A$ has constant degree $n$, then $\uGL_1(A)$ has type $\rmA_{n-1}$ \cite[3.5.0.87]{CF}, while $\uAut(A)$ is semisimple adjoint of type $\rmA_{n-1}$ \cite[3.2.0.64, 3.5.8.82]{CF}.

For $A=\Mat_n(R)$ we put $\uGL_n = \uGL_1(A)$ and $\uPGL_n = \uGL_n/\GG_m$, following \cite[5.11]{GroI}. We then have an isomorphism  $\uPGL_n \simlgr \uAut\big(A)$ of $R$--group schemes.
\sm

\item\label{azu-idemp} ({\em Idempotents} \cite[Lem.~3.4]{Salt}, see also \cite[7.6.4]{Ford}) Let $c=c^2 \in A$ with $\ann_R(c) = 0$. Then $cAc$ is an Azumaya algebra, and $Ac$ is  faithfully projective as left $A$-- and also as $R$--module. Hence $\End_R(Ac)$ and $\End_A(Ac)$ are also Azumaya $R$--algebras by \eqref{azu-ex}. Right multiplication gives an $R$--algebra isomorphism
    \begin{equation} \label{azu-idem1}
      (cAc)\op \simlgr \End_A(Ac), \quad x \mapsto R_x|_{Ac}.
      \end{equation}

\item\label{azu-Br} ({\em Brauer equivalence})
Among the various characterizations of Brauer equivalence (\cite[\S7.3]{Ford}, \cite[III, 5.6]{KO}, \cite[Lem.~3.4]{Salt}) we will use the following: an Azumaya $R$--algebra $B$ is Brauer equivalent to $A$, in symbols $B \sim A$, if one of the following equivalent conditions holds:
\end{inparaenum}

\begin{enumerate}[label={\rm (\roman*)}]
  \item \label{rng-iv} $A \ot_R \End_R(P) \cong B \ot_R \End_R(Q)$ for some faithfully projective $R$--modules $P$ and $Q$,

  \item\label{rng-i} $B \cong \End_A(P)$ for some faithfully projective right $A$--module $P$,

  \item\label{rng-iii} ${A \ot_R B\op} \cong \End_R(M)$ for some  faithfully projective $R$--module $M$,

  \item \label{rng-ii} there exists $r\in \NN_+$ and an idempotent $c\in \Mat_r(A)$ with $\ann_R(c) = \{0\}$ and $B\cong c \, \Mat_r(A) \, c$.
\end{enumerate}
The equivalence of \ref{rng-iv}, \ref{rng-i} and \ref{rng-iii} is \cite[III, 5.6]{KO} (and also \cite[7.3.4]{Ford}). Given \ref{rng-i}, we choose $r\in \NN_+$ such that $A^r = P \oplus Q$ for a right $A$--module $Q$. The projection of $A^r$ onto $P$ is an idempotent $c\in \End_A(A^r) = \Mat_r(A)$ for which \ref{rng-ii} holds. Retracing the steps (or using \cite[Lem.~3.4]{Salt}) shows \ref{rng-ii} $\implies$ \ref{rng-i}.
\ms

The following Lemma~\ref{lemsepsub} clarifies the role of maximal \'etale subalgebras of an Azumaya $R$--algebra.

\begin{lem}[\'Etale subalgebras] \label{lemsepsub} Let $S$ be a commutative $R$--subalgebra of the Azumaya $R$--algebra $A$. \sm

\begin{inparaenum}[\rm (a)] \item \label{lemsepsup-a} Then
\begin{align*}
  \text{$S$ is a separable $R$--algebra} \quad &\iff \quad
  \text{$S$ is a finite \'etale $R$--algebra.}
\end{align*}
In this case, \end{inparaenum}

\begin{enumerate}[label={\rm (\roman*)}]
   \item \label{lemsepsup-ai} $A$ is a faithfully projective $S$--module,  hence {\rm \ref{azu}\eqref{azu-coIIi}\/} applies,
    in particular $\al_S$ is a faithfully flat neutralization, and

  \item \label{lemsepsup-aii} $S$ is a direct summand of the $S$--module $A$ and the $R$--module $A$.
\end{enumerate}
\sm

\begin{inparaenum}[\rm (a)] \setcounter{enumi}{1}
\item \label{lemsepsup-b} Let $S$ be a separable, hence finite \'etale $R$--algebra. Then
\begin{equation}\label{lemsepsup-b1}
  \text{$S$ is maximal commutative} \quad \iff \quad \deg_R A = \rank_R S.
  \end{equation}
\end{inparaenum}
\end{lem}
\sm

Regarding \eqref{lemsepsup-b1} we recall that both $\rank_R$ and $\deg_R$ are considered as locally constant functions on $\Spec (R)$.

\begin{proof} \eqref{lemsepsup-a} If $S$ is a separable $R$--algebra, \cite[4.4.1(2)]{Ford} says that $A$ is a projective $S$--module. It is then necessarily faithfully projective. It follows (\cite[III, 1.9]{KO}) that  $S$ is a direct summand of the $S$--module $A$ and therefore also a direct summand of the $R$--module $A$. As such, $S$ is a finite projective $R$--module and therefore a finite projective $R$--algebra, hence a finite \'etale $R$--algebra. \sm

\eqref{lemsepsup-b} If $S$ is maximal commutative, then $\deg_R A = (\rank_R A)^2 = \rank_R S$ by \eqref{azu-2}. For the proof of the converse, we note that $S$ is a direct summand of the $R$--module $A$ by \ref{lemsepsup-aii}. Also, it is no harm to assume that $\rank_R S$ is constant. Then \cite[5.4]{GroI} says that the $S$--algebra homomorphism $\al_S$ of \eqref{azu-com1} is an isomorphism. By \eqref{azu-com2}, $S$ is maximal commutative.
\end{proof}

\begin{defn}[Maximal \'etale] \label{maxetdef} Following Grothendieck \cite[5.6]{GroI} we call a finite \'etale $R$--subalgebra of $A$ a {\em maximal \'etale $R$--subalgebra of $A$}, if it satisfies \eqref{lemsepsup-b1}.

Every Azumaya $R$--algebra $A$ over a local ring contains a maximal \'etale subalgebra; one can even assume that it is one-generated
\cite[10.3.2]{Ford}. 
We will extend this to the semilocal case in \ref{maxet}, using the Hamilton-Cayley Theorem \eqref{rcp-22} for Azumaya algebras.
We point out that in general $A$ need not contain a maximal \'etale subalgebra
\cite[III, 6.5]{KO}, although $A$ always has a finite \'etale neutralizing ring \cite[10.3.9]{Ford}.
\end{defn}

\subsection{Reduced trace, reduced norm, and  reduced characteristic polynomial} \label{rcp}  For every $a\in A$ there exist
unique elements $\Nrd_{A/R}(a) \in R$ and $\Trd_{A/R}(a) \in R$, called {\em reduced norm\/} and {\em reduced trace\/} respectively, and a unique polynomial $\Pcrd_{A/R}(a; X)\in R[X]$, called the {\em reduced characteristic polynomial of $a$}, such that for every  faithfully flat neutralization $\al \co A \ot_R S \simlgr \End_S(M)$ in the sense of \ref{azu}\eqref{azu-spli} we have
\begin{equation}  \label{rcp-1}  \begin{split}
  \Nrd_{A/R}(a) \ot 1_S &= \det\big(\al(a \ot 1_S)\big), \\
 \Trd_{A/R}(a) \ot 1_S &= \Tr\big( \al(a \ot 1_S)\big)\quad \text{and} \\
\Pcrd_{A/R}(a;X) &= \det \big(X  - \al(a \ot 1_S)\big)
\end{split} \end{equation}
where $\det$ and $\Tr$ are the determinant and trace of an endomorphism \ref{trno} and where $\Pcrd(a;X)$ is the characteristic polynomial of the endomorphism $\al(a \ot 1_S)$. All three data are defined by faithfully flat descent, \cite[5.9]{GroI}, \cite[IV; \S2]{KO} or \cite[4.3]{Salt}.  We have the relation
\begin{equation} \label{rcp-2}
 \Pcrd_{A/R}(a;X) = \Nrd_{A[X]/R[X]}(X-a),
\end{equation}
analogous to \eqref{trno-0}. It is immediate that $\Nrd$, $\Trd$ and $\Pcrd$ are stable under base change and respect direct products. Each $a\in A$ satisfies its reduced characteristic polynomial (Hamilton-Cayley Theorem),
\begin{equation}\label{rcp-22}
 \begin{split}  \Pcrd(a; a) &= 0, \qquad \text{and} \\
 \Nrd_{A/R}(ab) &= \Nrd_{A/R}(a) \cdot \Nrd_{A/R}(b)
\end{split} \end{equation}
holds for $a,b\in A$ (\cite[IV, 2.3]{KO}). 
The reduced trace is a linear form on $A$ satisfying $\Trd(ab) = \Trd(b,a)$; the associated symmetric bilinear form
$A\times A \to R$, $(a,b) \mapsto \Trd(ab)$, is regular \cite[11.1.6]{Ford}.
\sm

If $A$ has constant rank $n^2$, $n\in \NN_+$, then
\begin{equation}\label{rcp-3} \begin{split}
 \Nrd_{A/R}(a) &= (-1)^n \Pcrd(a; 0) \qquad \text{and} \\
\Pcrd_{A/R}(a;X) &= X^n - \Trd_{A/R}(a) X^{n-1} + \cdots + (-1)^n \Nrd_{A/R}(a),
\end{split}\end{equation}
as in \eqref{trno-00}, and $\Nrd(r) = r^n$ for $r\in R$.

\subsection{Some examples.}\label{soex} Let $\al \co A \ot_R S \to \End_S(M)$ be a a faithfully flat neutralization of $A$.
\sm

\begin{inparaenum}[(a)] \item\label{soex-op}  ({\em Opposite algebra $A\op$}) The opposite algebra $A\op$ is an Azumaya $R$--algebra. The dual $S$--module $M^*$ is a faithfully projective $S$--module and, denoting by $^t \al(a \ot 1_S)\in \End_S(M^*)$ the transpose endomorphism, the map
\[ ^t \al \co A\op \ot_R S \simlgr \End_S(M^*); \quad a \ot 1_S \mapsto {^t \al}(a \ot 1_S) \]
is a faithfully flat neutralization of $A\op$. Since locally the matrix of ${^t\al}(a \ot 1_S)$ is the transpose of the matrix of $\al(a \ot 1_S)$, the reduced data of $A$ and of $A\op$ coincide, for example $\Nrd_{A\op/R}(a) = \Nrd_{A/R}(a)$.

\item \label{soex-sp} ({\em Split algebras}) Let $A = \End_R (P)$,  where $P$ is a faithfully projective $R$--module. Then $\Id$ is a faithfully flat neutralization of $A$ and the reduced characteristic polynomial, the reduced trace and the reduced norm are, respectively, the characteristic polynomial, the trace and the determinant.   \sm

\item \label{soex-mat} ({\em Matrix algebras}) The matrix algebra $\Mat_r(A)$, $r\in \NN_+$, is an Azumaya $R$--algebra, \ref{azu}\eqref{azu-tensor}. The map
    \begin{align*} \Mat_r(\al) \co  \Mat_r(A) \ot_R S &\simlgr \Mat_r \big( \End_S(M) \big) \cong
  \End_S(M^r) \\
  (a_{ij} ) \ot s &\;  \mapsto \quad \big( \al(a_{ij} \ot s)\big)
\end{align*}
is a faithfully flat neutralization of $\Mat_r(A)$. Let $d=\diag(a_1, \ldots, a_r)$ be the diagonal matrix in $\Mat_R(A)$ with entries $a_1, \ldots , a_r$. Then $\Mat_r(\al)$ sends $d$ to the diagonal matrix
$\diag\big( \al(a_1\ot 1_S), \ldots , \al(a_r \ot 1_S)\big)$,
whose characteristic polynomial and determinant are the products of the corresponding gadgets on the diagonal, while its trace is the sum of traces of its diagonal entries. Hence
\begin{equation} \label{rdp-mat2} \begin{split}
  \Pcrd_{\Mat_r(A)/R} (d ; X\big) &= \textstyle \prod_i \Pcrd_{A/R}(a_i),\\
  \Nrd_{\Mat_r(A)/R} (d ) &= \textstyle \prod_i \Nrd_{A/R}(a_i),\\
 \Trd_{\Mat_r(A)/R} (d ) &= \textstyle \sum_i \Trd_{A/R}(a_i).
\end{split}\end{equation}
In particular, for $E_r =\diag(1_A, \ldots, 1_A)$ the identity matrix in $\Mat_r(A)$ we obtain
\begin{equation}  \label{rdp-mat}
\begin{split}
  \Pcrd_{\Mat_r(A)/R}(aE_r; X) & = \Pcrd_{A/R}(a)^r , \\
   \Nrd_{\Mat_rA)/R}(aE_r) &= \Nrd_{A/R}(a)^r, \\
   \Trd_{\Mat_r(A)/R}(aE_r) &= r \Trd_{A/R}(a).
\end{split}
\end{equation}
Let $A_1$ and $A_2$ be two Azumaya $R$--algebras and let $a_i\in A_i$. Then 
\begin{equation}
  \label{soex-mat3} \Trd_{A_1 \ot_R A_2} (a_1 \ot a_2) = \Trd_{A_1}(a_1) \, \Trd_{A_2}(a_2).
\end{equation}
Indeed, we can assume that $A_1 = \Mat_n(R)$, hence $A_1 \ot_R A_2 \simlgr \Mat_n(A)$ under $(x_{ij} \ot a \mapsto (x_{ij}a)$. The formula then follows from \eqref{rdp-mat2}. 
\sm

\item \label{soex-mor} ({\em Morita contexts}) Let $e\in A$ be an idempotent with $\mathrm{Ann}_R(e) = \{0\}$, and put $f= 1_A - e$. By \ref{azu}\eqref{azu-Br} the subalgebras $B=eAe$ and $C=fAf$ are Azumaya $R$--algebras, both Brauer equivalent to $A$. The endomorphisms $\al(e\ot 1_S)$ and $\al(f\ot 1_S)$ are complementary projections, thus inducing a decomposition $M= P \oplus Q$ into faithfully projective $S$--submodules of $M$. The restrictions of $\al$ to $B \ot S$ and $C\ot S$ are faithfully flat neutralizations of $B$ and $C$ with values in $\End_S(P)$ and $\End_S(Q)$ respectively. Thus, given $b\in B$ and $c\in C$, the endomorphisms $\al( (b + c)\ot 1_S$ respects the decomposition $M=P\oplus Q$. Hence
    \begin{equation} \label{soex-mor1}
     \Nrd_{A/R}(b + c) = \Nrd_{B/R}(b) \, \Nrd_{C/R}(c), \end{equation}
    $\Pcrd_{A/R}(b + c;X) = \Pcrd_{B/R}(b;X) \Pcrd_{C/R}(c;X)$ and
    $\Trd_{A/R}(b + c) = \Trd_{B/R}(b) + \Trd_{C/R}(c)$.
\end{inparaenum}


\begin{lem}\label{maxet} Let $A$ be an Azumaya algebra over a semilocal ring $R$. Then $A$ contains a maximal \'etale subalgebra. If $A$ has constant rank, then $A$ contains a maximal \'etale subalgebra which is one-generated.
\end{lem}

\begin{proof} Azumaya algebras and \'etale algebras respect direct products, \ref{azu}\eqref{azu-dir} and \ref{fea}\eqref{fea-d}. It is therefore sufficient to prove the second assertion.

Let $\deg_R A = r$, let $\ka_1, \ldots, \ka_n$ be the residue fields of $R$, thus $R/\Jac(R) = \ka_1 \times \cdots \times \ka_n$. Then $ A/\Jac(R) A = A \ot_R (R/\Jac(R)) = A_1 \times \cdots \times A_n$ is a direct product of Azumaya $\ka_i$--algebras of degree $r$. By \cite[7.5.7]{Ford}, each $A_i$ contains a maximal \'etale $\ka_i$--algebra $S_i$, which is one-generated, say $S_i = \ka_i[s_i]$ for some $s_i \in S_i$. We know that each $S_i$ is $r$--dimensional, for example from \eqref{lemsepsup-b1}.

Let $s\in A$ be a lift of $(s_1, \ldots, s_n ) \in A/\Jac(R) A$ and let $S=R[s]= \Span_R \{ s^n : n\in \NN\}$. By \eqref{rcp-3}, the reduced characteristic polynomial is monic of degree $r$. Since $\Pcrd(s;s)= 0$ by \eqref{rcp-22}, the $R$--subalgebra $S$ is $S= \Span_R \{1_A, s, \ldots, s^{r-1}\}$; it is free of rank $r$ by Nakayama \ref{nak}, 
it is finite \'etale since every reduction $S\ot \ka_i = \ka_i[s_i]$ is finite \'etale, and it is then maximal \'etale by \eqref{lemsepsup-b1}. \end{proof}

The reduced and ordinary concepts are compared in the following lemma.

\begin{lem}\label{lem_nred} Let  $a\in A$.

 \begin{enumerate}[label={\rm (\alph*)}]
\item \label{lem_nreda} If $A$ has constant rank $n^2$, then

 \begin{enumerate}[label={\rm (\roman*)}]
 \item  \label{lem_nred_a} $\mathrm{Pc}_{A/R}(a;X)= \mathrm{Pcrd}_{A/R}(a;X)^n$;

\item \label{lem_nred_b} $\Tr_{A/R}(a)= n \, \Trd_{A/R}(a)$;

\item \label{lem_nred_c} $\rmN_{A/R}(a)=  \Nrd_{A/R}(a)^n$.

\end{enumerate}
\sm

\item \label{lem_nredb} We have
\begin{equation}\label{lem_nredb1} \begin{split}
     & a\in A\ti \iff \Nrd_{A/R}(a) \in R\ti \\
     & \iff a \ot 1_{R/\m} \in \big(A \ot_R (R/\m)\big)\ti \text{ for all maximal ideals $\m\ideal R$.}
\end{split}\end{equation}
\end{enumerate}
\end{lem}

\begin{proof} \ref{lem_nreda} We first prove \ref{lem_nred_c}. By \ref{azu} we can assume that $A=\Mat_n(R)$. In this case $\Nrd_A= \det$. The formula then follow from the fact that for
$a \in \Mat_n(R)$ the matrix of the left multiplication $L_a$ with respect to
the standard basis of $\Mat_n(R)$ is the block diagonal matrix $\mathrm{diag}(a, \ldots , a )$, cf.\ \cite[III, \S9.3, Ex.~3]{BA}.

Taking into account \eqref{rcp-2} and \eqref{trno-0}, the formula \ref{lem_nred_a} is a consequence of \ref{lem_nred_c} and in turn implies the formula \ref{lem_nred_b} by reducing to $A$ having  constant rank and then applying \eqref{rcp-3}.

For the proof of the first equivalence in \ref{lem_nredb} we can again assume that $A$ has constant rank, cf.\ \ref{azu}\eqref{azu-dir}. The claims then follows from $a\in A\ti \iff \rmN_{A/R}(a) = \Nrd_{A/R}(a)^n \in R\ti$.
Since $\Nrd$ commutes with base change, we get
$\Nrd_{A/R}(a) \in R\ti \iff \Nrd_{(A \ot \ka(\m))/ \ka(\m)} (a\ot 1_{\ka(\m)} )\in \ka(\m)\ti$
for every maximal ideal $\m \ideal R$, which is equivalent to the second condition of \eqref{lem_nredb1}. \end{proof}

\begin{lem} \label{prop_nred} Let $S \subset A$ be a maximal
commutative $R$--subalgebra of $A$ such that $A$ is a projective  $S$--module. Recall\/   {\rm \ref{azu}\eqref{azu-coIIi}\/} that then the $R$--module $S$ is faithfully projective. For every $s \in S$ we have the following formulae.
 \begin{enumerate} [label={\rm (\alph*)}]

\item \label{prop_nred_a} $\Pc_{S/R}(s;X)= \Pcrd_{A/R}(s;X)$;

\item \label{prop_nred_b} $\Tr_{S/R}(s)= \Trd_{A/R}(s)$;

\item \label{prop_nred_c} $\rmN_{S/R}(s)= \Nrd_{A/R}(s)$.
 \end{enumerate}
\end{lem}

\begin{proof}
As in the proof of Lemma~\ref{lem_nred} it suffices to prove \ref{prop_nred_c}.
Let $V$ be $A$ viewed as $S$--module. We know (\ref{azu}\eqref{azu-coIIi}\/) that $S$ is a  faithfully flat neutralizing ring of $A$ with neutralization $\al \co A \ot_R S \simlgr \End_S(V)$ induced by the sandwich map. Thus $\al(s \ot 1_S) = L^{(V)}_s$, the scalar multiplication of $V$ by $s\in S$. Hence by \eqref{rcp-1} we have $\Nrd_{A/R}(s) = \det_S\big(\al(s\ot 1_S)\big) = \det_S(L_s^{(V)})$ where $\det_S$ should indicate that the determinant is for an $S$--module homomorphism. We will also use $\det_R$ with the analogous meaning. For the left-hand side of \ref{prop_nred_c}, we have  $\rmN_{S/R}(s) = \det_R(L_s^{(S)})$, where $L_s^{(S)}$ is the left multiplication of the $R$--module $S$. Since $\det_R(L_s^{(S)})\ot_R 1_S= \det_S(L_s^{(S)} \ot_R \Id_S)$, our claim is $\det_S(L_s^{(V)})= \det_S(L_s^{(S)}\ot_R \Id_S)$.

For proving this claim we can in view of \eqref{azu-2} assume that $A$ is free of rank $n$ as $S$--module and that $S$ is free of rank $n$ as $R$--module. We then have an $S$--module isomorphism
\[ \psi \co S \ot_R S = S\ot_R R^n \simlgr S^n \cong V, \quad
t \ot (r_1, \ldots r_n) \mapsto (tr_1, \ldots , tr_n)
\]
satisfying $\psi \circ (L_s^{(S)} \ot_R \Id_S) = L_s^{(V)} \circ \psi$ and thus  proving our claim.  \end{proof}
\sm

\textbf{Remark.} Lemma~\ref{prop_nred} is known over fields. In this setting the reader can find two different proofs in \cite[\S17.5, Prop.~7]{BA8} and \cite[2.6.3]{GS}. The latter proof can be adapted to the setting here.
\sm

The following lemma serves as preparation for Proposition~\ref{prop_pc_cn}.

\comments{(2020-05-15) Is it sufficient to assume in \ref{propAD} that $Z/R$ is finite projective? I assume $Z/R$ \'etale to get that $A$ is faithfully projective $T$--module }

\begin{lem}\label{propAD}
Let $D$ a separable $R$--subalgebra of $A$ whose centre $\rmZ(D) = Z$ is finite \'etale as $R$--algebra, and let $T$ be a maximal \'etale subalgebra of $D$. To summarize:
\[\xymatrix@C=40pt{
 R \ar[r]^{\text{fin. \'et.}} \ar@/_1pc/[rrr]_{\text{separable}}
  & Z \ar[r]^{\text{max. \'et.}} &T \ar[r] &D \ar[r] &A
} \]
Then the following hold. \sm

\begin{inparaenum}[\rm (a)]
 \item\label{prop-AD-a} $T$ is a finite \'etale and hence faithfully flat $R$--algebra. With respect to the obvious right $T$--action, $A$ is a faithfully projective $T$--module, and the homomorphism
    \begin{equation}\label{prop-AD-a1}
      \al =  \al_T  \co A \ot_R T \longto \End_T(A), a\ot t \mapsto L_a R_t
     \end{equation}
     of \eqref{azu-com1} is an injective $T$--algebra homomorphism. \sm

\item\label{propAD-b} The algebra $S= Z\ot_R T$ is finite \'etale as $Z$--algebra and as $R$--algebra. The underlying additive group of $A$ together with the $S$--action
  \begin{equation} \label{propAD-b1}  (z \ot t) . a = zat \qquad(z\in Z, t\in T, a \in A),   \end{equation}
  makes $A$ a faithfully projective $S$--module, denoted $M$. The Weil restriction
  $\sfR_{S/T}(M)$ of $M$ is the standard right $T$--module $A$.
\sm

\item \label{propAD-c} The map
  \begin{equation} \label{propAD-c1}  \beta \co D\ot_Z S \longto \End_S(M), \quad d\ot s \mapsto (a \mapsto dsa)
  \end{equation}
  is an injective $S$--algebra homomorphism  satisfying
\begin{equation}\label{propAD-c2}
      \al (d \ot 1_T) = \sfR_{S/T}\big( \be(d\ot 1_S)\big)
\end{equation}
 for $d\in D$, where $\sfR_{S/T}\big( \be(d\ot 1_S)\big)$ is  the endomorphism $\be(d\ot 1_S)\in \End_S(M)$ considered as $T$--linear map.
\sm

\item \label{propAD-d} Let $\vphi\co \Spec(Z) \to \Spec(R)$ be the canonical map associated with the structure map $R \to Z$ and suppose
 \begin{equation}
   \label{propAD-d1} (\rank_R A) \circ \vphi = (\rank_R D)  \cdot
      \big( (\rank_R Z) \circ \vphi \big)^2.
 \end{equation}
Then $\al$ and $\beta$ are isomorphisms, and hence in particular faithfully flat neutralizations of $A$ and $D$ respectively.
\end{inparaenum}
\end{lem}

\begin{proof}
  \eqref{prop-AD-a} Since $T$ is a maximal \'etale subalgebra of $D$, it is finite \'etale over $Z$ by \ref{lemsepsub}. It is then also finite \'etale as $R$--algebra by transitivity of ``finite \'etale''. Because $A$ is faithfully projective as $R$--module, the $T$--module $A$ is faithfully projective by  \ref{sep}\eqref{sep-weil}. It then follows from \ref{azu}\eqref{azu-com} that $\al_T$ is an injective $T$--algebra homomorphism. \sm

\eqref{propAD-b} Since $T$ is a finite \'etale $R$--algebra, $S$ is finite \'etale as $T$--algebra because the property ``finite \'etale" respects base change. Also, as the tensor product of two finite \'etale $R$--algebras, $S$ is a finite \'etale $R$--algebra. 

It is straightforward to check that $A$ becomes an $S$--module with the $S$--action \eqref{propAD-b1}. The restriction of this $S$--action to $R$ is the given $R$--action on $A$. Hence, by \ref{sep}\eqref{sep-weil}, the $S$--module $A$ is projective, and then obviously faithfully projective.
\sm

\eqref{propAD-c} It is immediate that $\beta$ is a well-defined $S$--algebra homomorphism. Since $D$ is an Azumaya $Z$--algebra, $D\ot_Z S$ is an Azumaya $S$--algebra, \ref{azu}\eqref{azu-bc}. Also $\End_S(M)$ is an Azumaya $S$--algebra, \ref{azu}\eqref{azu-ex}. Hence $\beta$ is an injective $S$--algebra homomorphism by \ref{azu}\eqref{azu-com}. The equation \eqref{propAD-c2} holds by construction.
\sm

\eqref{propAD-d} Since both maps are $R$--linear, it suffices to prove bijectivity after localization in  a maximal ideal of $R$, in other words, we can assume that $R$ is local.

(I) {\em Bijectivity of $\al$}: Since the $R$--module $A$ is finite projective, it is free of constant square rank, say $\rank_R A = n^2$. By the same argument, the $R$--module $Z$ is free, say of $\rank_R Z = r$. Now \eqref{propAD-d1} becomes $n^2 = (\rank_Z D) \cdot  r^2$, implying that the $Z$--module $D$ has constant rank, say $\rank_Z D = d^2$, and then that $n = rd$. In fact, $Z$ being semilocal, $D$ is free of rank $d^2$ as $Z$--module. Applying \eqref{lemsepsup-b1} to the maximal \'etale subalgebra $T$ of the Azumaya $Z$--algebra $D$, shows that the $Z$--module $T$ has constant rank $d$. Hence $T$ is free of rank $d$ as $Z$--module and consequently free of rank $rd$ as $R$--module. Now \eqref{lemsepsup-b1} again, this time for $A$ and $T$, shows that $T$ is a maximal commutative $R$--subalgebra of $A$. Finally, $\al$ is an isomorphism by \eqref{azu-com2}.

(II) {\em Bijectivity of $\be$}: By \ref{azu}\eqref{azu-hom} it suffices to show $\deg_S (D\ot_Z S) = \rank_S M$.
By (I), $\rank _Z D = d^2$, so $\deg_S  D =d$. By \eqref{propAD-b}, the Weil restriction $\sfR_{S/T}(M)$ is the $T$--module $A$, which, as we have seen in (I), is free of rank $n = rd$. Moreover, since $Z$ is free of rank $r$ as $R$--module, the $T$--module $S = Z \ot_R T$ is free of rank $r$, By transitivity of ranks, $M$ has constant rank $d$, proving $\deg_S D = d \rank_S M$.
\end{proof}
\sm

We can now compare the reduced concepts \ref{rcp} of Azumaya algebras with that of separable subalgebras. We recall that $\rmN$ and $\Tr$ denotes the norm and trace  defined in \ref{trno}.

%
%
%

\begin{prop}\label{prop_pc_cn}
Let $D$ a separable $R$--subalgebra of $A$ whose centre $\rmZ(D) = Z$ is finite \'etale as $R$--algebra.  We further suppose \eqref{propAD-d1}, i.e.,
\begin{equation}   \label{prop_pc_cn1}
 (\rank_R A) \circ \vphi = (\rank_Z D) \cdot \big((\rank_R Z) \circ \vphi\big)^2
\end{equation}
Then the following hold for $d\in D$.
\begin{enumerate}[label={\rm (\alph*)}]
 \item  \label{prop_pc_c1} $\Pcrd_{A/R}(d;X)= \rmN_{Z[X]/R[X]}\bigl(
 \Pcrd_{D/Z}(d;X)\bigr)$;

 \item \label{prop_pc_c2} $\Trd_{A/R}(d)=  \Tr_{Z/R}\bigl( \Trd_{D/Z}(d) \bigr)$;

 \item \label{prop_pc_c3} $\Nrd_{A/R}(d)= \rmN_{Z/R}\bigl( \Nrd_{D/Z}(d) \bigr)$.
\end{enumerate}
\end{prop}

\begin{proof} It suffices to prove the equations after localization in a maximal ideal of $R$. Thus, we can assume that $R$ is a local ring. The $R$--algebra $Z$ is finite projective, hence a semilocal ring. We have seen in the proof of \ref{propAD}\eqref{propAD-d} that the condition  \eqref{prop_pc_cn1} implies that $D$ has constant rank as $Z$--module. It now follows from Lemma~\ref{maxet} that the Azumaya $Z$--algebra $D$ contains a maximal \'etale subalgebra $T$. Thus, we are in the setting of Lemma~\ref{propAD} and have
\begin{enumerate} [label={\rm (\roman*)}]
\item\label{prop_pc_cn-b} a faithfully flat neutralization $\be \co D \ot_Z S \simlgr \End_S(M)$, and

  \item\label{prop_pc_cn-c} a faithfully flat neutralization $\al \co A \ot_R T \simlgr \End_T\big( \sfR_{S/T}(M)\big)$ satisfying
      $   \al(d \ot 1_T) = \sfR_{S/T}\big( \be(d\ot 1_S)\big)$.
  \end{enumerate}
As in the proof of Lemma~\ref{lem_nred}, it suffices to prove \ref{prop_pc_c3}. This is a consequence of the following chain of equations, where the crucial step in the second line below is \cite[III, \S9.4, Prop.~6]{BA}; the last equality uses base change of the norm $\rmN$ and $S= Z\ot_R T$.
\begin{align*}
  \Nrd_{A/R}(d) \ot 1_T &= \det\big(\al(d\ot 1_T)\big)
        = \det\big(\sfR_{S/T}\big(\be(d\ot 1_S)\big)\big)
   \\ &= \rmN_{S/T}\big( \det\big(\be(d\ot 1_S)\big)\big)
   \\ & = \rmN_{S/T}\big( \Nrd_{D/Z}(d) \ot 1_S\big)
   \\ & =  \rmN_{Z/R}\big( \Nrd_{D/Z}(d)\big) \ot 1_T.
\end{align*}
Now \ref{prop_pc_c3} follows by faithfully flat descent.   \end{proof}

\begin{cor} \label{prop_pc} Let $B \subset A$ be a separable $R$--subalgebra of $A$. Thus, by {\rm \ref{sep}\eqref{azu-sep}}, $B$ is an Azumaya algebra over its centre $Z = \rmZ(B)$ and $Z$ is a separable $R$--algebra. \sm

\begin{inparaenum}[\rm (a)]
\item   \label{prop_pc_a}
  Let $D = \rmC_A(Z)$, the centralizer of $Z$ in $A$. Then $D$ is an Azumaya algebra over $Z$.\sm

 \item   \label{prop_pc_b} Let $B'=\rmC_A(B)$. Then $B'$ is an Azumaya $Z$--algebra and the multiplication of $A$ induces an isomorphism $B \ot_{Z} B' \simlgr D$ of Azumaya $Z$--algebras.%
    \sm

\item  \label{prop_pc_d} We assume that $B'$ has constant rank $r$ over $Z$, that $Z$ is finite projective as $R$--module, i.e. $Z$ is a finite \'etale $R$--subalgebra,  and that, using the notation of\/ {\rm \ref{prop_pc_cn}},
 \begin{equation} \label{prop_pc_dii}
  (\rank_R A) \circ \vphi = r^2 \cdot (\rank_Z B) \cdot \big((\rank_R Z) \circ \vphi\big)^2
  \end{equation}
holds. Then for each $b \in B$ we have the relations:
 \begin{align*}
   \Pcrd_{A/R}(b;X) &= \rmN_{Z[X]/R[X]}\bigl(  \Pcrd_{B/Z}(b;X)\bigr)^r, \\
  \Trd_{A/R}(b) &= r \, \Tr_{Z/R}\bigl( \Trd_{B/Z}(b) \bigr), \\
  \Nrd_{A/R}(b) &= \rmN_{Z/R}\bigl( \Nrd_{B/Z}(b) \bigr)^r.
 \end{align*} \end{inparaenum}
\end{cor}

\comments{(2021-04-26) According to \cite[\S24.5, page 256, Danger]{BA8}, even over fields an assumption like \eqref{prop_pc_dii} in \ref{prop_pc}\eqref{prop_pc_d} is necessary.}

\begin{proof} \eqref{prop_pc_a} Since $A$ is a faithfully projective $R$--module and $Z$ is a separable $R$--algebra, it follows from \ref{sep}\eqref{sep-weil} 
that $A$ is projective as $Z$--module. The claim then follows from \cite[4.1]{KOS} (= \cite[7.4.4]{Ford}). \sm

\eqref{prop_pc_b} By \ref{azu}\eqref{azu-cdt},  the $R$-algebra $B'$ is separable and we have $B=\rmC_A(B')$. It follows that the centre of $B'$ is $\rmZ(B') = B' \cap \rmC_A(B')= B' \cap B =\rmC_A(B) \cap B=Z$, so that $B'$  is an Azumaya algebra over $Z$.

Observe that both $B$ and $B'$ are Azumaya $Z$--subalgebras  of the Azumaya $Z$--algebra $D$. The isomorphism  $B \otimes_{Z} B' \simlgr D$ is therefore another consequence of the Double Centralizer Theorem.
\sm

\eqref{prop_pc_d} Taking into account the formulas in Lemma~\ref{prop_pc_cn}, it suffices to prove $\Pcrd_{D/Z}(b; X) = \Pcrd_{B/Z}(b;X)^r$, $\Trd_{D/Z}(b)=  r \Trd_{B/Z}(b)$ and $\Nrd_{D/Z}(b) = \Nrd_{B/Z}(b)^r$. After localizing in the \'etale topology for $Z$, we can assume that $B' = \Mat_r(Z)$, so that $D \cong B \ot_Z \Mat_r(Z) \cong \Mat_r(B)$. The claims then follow from \eqref{rdp-mat}.
\end{proof}

\subsection{Example for \ref{prop_pc}.} \label{prop_pcex} We give an example, in which the assumption \eqref{prop_pc_dii} is fulfilled. As in \ref{prop_pc} we assume that $A$ is an Azumaya $R$--algebra, $B$ is a separable $R$--subalgebra whose centre $Z$ is finite projective. In addition, we suppose that
\begin{enumerate}[label={(S)}]  
  \item[(S)] for every $\m \in \Specmax(R)$ the $\ka(\m)$--algebra $B_{\ka(\m)}$ is simple.
\end{enumerate}
Since centralizers commute with flat base change, it follows from \cite[\S15.5, Thm.~5]{BA8} that the centralizer subalgebra $B'$ also satisfies (S) and the dimension formula $\dim_{\ka(\m)} A\ot_R \ka(\m) = (\dim_{\ka(\m)} B_{\ka(\m)}) \cdot (\dim_{\ka(\m)} B'_{\ka(\m)})$ holds. Thus
 \begin{equation} \label{prop_pcex1}
  (\rank_R A) \circ \vphi = (\rank_Z B) \cdot (\rank_Z B') \cdot \big((\rank_R Z) \circ \vphi\big)^2,
  \end{equation}
which implies \eqref{prop_pc_dii} in case $B'$ has constant rank.

In particular, the above holds if $R$ is a field. In this case, the formulas in \ref{prop_pc}\eqref{prop_pc_d} are proven in \cite[\S17.5, Prop.~8]{BA8}.

\comments{I am not sure that the following remark regarding a possible generalization of \ref{prop_pc} is correct, because of the missing assumption \eqref{prop_pc_dii}. We can add this later: \\

{\sf If $Z=Z_1 \times \dots \times Z_u$, $B=B_1 \times \dots \times B_u$,
and $B'=B'_1 \times \dots \times B'_u$ with $B'_i$
of constant degree $r_i$ over $Z_i$ for $i=1,..,u$, we have more generally
for $b=(b_1, \dots, b_u)$
 \begin{enumerate}[label={\rm (\roman*)}]
 \item  \label{prop_pc_e1} $\mathrm{Pcrd}_{A/R}(b;t)= \prod\limits_{i=1}^u
 \rmN_{Z[t]/R[t]}\Bigl(  \mathrm{Pcrd}_{B_i/Z_i}(b_i;t)\Bigr)^{r_i}$;

  \item \label{prop_pc_e2} $\Trd_{A/R}(b)= \sum\limits_{i=1}^u r_i \,
  \tr_{Z_i/R}\bigl( \Trd_{B_i/Z_i}(b) \bigr)$;

  \item \label{prop_pc_e3} $\Nrd_{A/R}(b)= \prod\limits_{i=1}^u \rmN_{Z_i/R}\bigl(
  \Nrd_{B_i/Z_i}(b_i) \bigr)^{r_i}$.
 \end{enumerate} 
} }

Another application of Proposition~\ref{prop_pc_cn} is the following corollary.

\begin{cor}\label{prop_red_stable} Let $Z\in \Ralg$ be a faithfully projective $R$--algebra. We view $A\ot_R Z= A_Z$ canonically as a right $A$--module as well as an Azumaya $Z$--algebra, and put $B= \End_A(A_Z)$.   \sm

\begin{enumerate} [label={\rm (\alph*)}]
 \item \label{prop_red_stable1} $B$ is an Azumaya $R$--algebra, which is Brauer equivalent to $A$. \sm

 \item \label{prop_red_stable2}
Let $f \in A_Z$. The right multiplication
$R_f: A_Z \to A_Z$, $x \mapsto x \, f $, is an element of $B$ and satisfies
\begin{align*}
 \Pcrd_{B/R}(R_f; X) &= \rmN_{Z[X]/R[X]}\big( \Pcrd_{A_Z/Z}(f;X)\big), \\
  \Trd_{B/R}(R_f) &= \Tr_{Z/R}\big( \Trd_{A_Z/Z}(f) \big), \\
  \Nrd_{B/R}(R_f) &= \rmN_{Z/R}\big( \Nrd_{A_Z/Z}(f)\big).
\end{align*}

\item \label{prop_red_stable3} We have $\rmN_{Z/R}\big( \Nrd_{A_Z/Z}( A_Z{}\ti) \big) \, \subset \, \Nrd_{B/R}(B\ti)$.
\end{enumerate}
\end{cor}

\begin{proof} \ref{prop_red_stable1} We use here only that $Z$ is a faithfully projective $R$--module. There exists an $R$--module $Z'$ and $d\in \NN_+$ such that $Z \oplus Z' = R^d$ as $R$--modules. Then $A^d = A\ot_R R^d = A_Z \oplus (A\ot_R Z')$ implies that $A_Z$ is a faithfully projective $A$--module. By \ref{azu}\eqref{azu-ex}, $B$ is an Azumaya $R$--algebra, and it is Brauer equivalent to $A$ by \ref{azu}\eqref{azu-Br}.
\sm

\ref{prop_red_stable2} We find it instructive to give two proofs.

\comments{(date unclear) Since the technique of the second proof is the same as the technique used to prove \ref{stab}\ref{stab-ii}, we could eliminate the second proof here and make a remark that \ref{prop_red_stable2} could also be prove with the techniques used in the proof of \ref{stab}\ref{stab-ii}. Your opinion?}
\sm

{\em First proof.} The map $A_Z \to D := \{R_f: f \in A_Z\}$, $f\mapsto R_f$,  is an isomorphism of $R$--algebras. It follows that $D$ is an Azumaya $Z$--algebra with centre isomorphic to $Z$. Therefore the claim is a special case of \ref{prop_pc_cn} with $A$ there replaced by $B$ here, as soon as we have verified the condition \eqref{prop_pc_cn1}. Since $B$ is the centralizer algebra of a subalgebra of $\End_R(A_Z)$ isomorphic to $A$, the first equality below is a consequence of the Double Centralizer Theorem \ref{azu}\eqref{azu-cdt}:
\begin{align*}
 \rank_R A \cdot  \rank_R\big(\End_A(A_Z)\big)&= \rank_R\big( \End_R (A_Z)\big)
  \\ &= (\rank_R A)^2 \cdot (\rank_R Z)^2.
\end{align*}
 Therefore $\rank_R(B) = ( \rank_R A) \cdot (\rank_R Z)^2$, and \eqref{prop_pc_cn1} follows using surjectivity of $\vphi$ and $\rank_R A = \rank_Z(A_Z) \circ \vphi$. \sm

{\em Second proof.} As in the proof of Lemma~\ref{lem_nred} it sufficed to prove the formula for the norm. The statement is insensitive to \'etale localization. So we can assume that $A=\Mat_n(R)$. In that case $\Nrd_{A/R}= \det_n$.

\noindent{\it First case: $n=1$.} Then $A=R$, $B=\End_R(Z)$, and the formula is the definition of $\rmN_{Z/R}$.

\noindent{\it General case.} Without loss of generality we can suppose that
$R$ is local, so that $Z$ is semilocal. We have to prove that the map
\[ h\co A\otimes_R Z =\Mat_n(Z) \xrightarrow{f \mapsto R_f}  B \xrightarrow{\Nrd_B} R,\]
is $\rmN_{Z/R} \circ \Nrd_{A_Z/Z}= \rmN_{Z/R} \circ \det_{n,Z}$.
By density it is enough to prove the formula for $f \in A_Z^\times= \GL_n(Z)$.
We denote by $\uB$ (resp.\ $\uB^{-}$) the Borel  $R$--subgroup scheme of upper (resp.\ lower)
triangular matrices of $\uGL_n$ and by $\uU$ (resp.\, $\uU^{-}$) their unipotent radicals. We now use Demazure's decomposition theorem \cite[XXVI.5.2]{SGA3}
\[
\GL_n(Z) = \uU(Z) \cdot \uU^{-}(Z)  \cdot \uB(Z).
\]
and obtain
\[
\GL_n(Z)= \uU(Z) \cdot \uU^{-}(Z) \cdot \uU(Z) \cdot
\begin{pmatrix}
Z^\times & 0 & \cdots & 0 \\
0 & Z^\times & \cdots & 0\\
\vdots &  & \ddots & \vdots  \\
0 & 0& \cdots &Z^\times
\end{pmatrix}.
 \]
\noindent Since $\det_{n,Z}$ and $h$ are trivial morphisms on $\uU(Z)$ and
${\uU}^{-}(Z)$,  we can assume that $f= \mathrm{diag}(z_1,\dots, z_n)$ with $z_1, \dots, z_n \in Z^\times$. In other words, we are reduced to the case of $n=1$, which has been checked at the beginning.
\sm

\ref{prop_red_stable3} is a straightforward consequence of
 \ref{prop_red_stable2}. \end{proof}

\comments{In \ref{rng} and below I do not assume that $A$ has constant rank.

I changed the approach to the definition of $\Nrd^{\rm st}(A)$ since it is more in the spirit of algebraic K-theory. }

\subsection{The reduced (stable) norm group of an Azumaya algebra}\label{rng}  Unless some confusion is possible we write $\Nrd$ for $\Nrd_{A/R}$. The image of the reduced norm homomorphism $\Nrd \co A\ti \to R$ is called the {\em reduced norm group\/} of $A$ and denoted $\Nrd(A\ti)$.

A natural approach to the stable version of $\Nrd(A\it)$ is the following. Recall \ref{azu}\eqref{azu-Br} that $\Mat_n(A)$, $n\in \NN_+$, is an Azumaya $R$--algebra. We abbreviate $\GL_n(A) = \Mat_n(A)\ti$, identify $\GL_\ell(A)$, $\ell < n$, with a subgroup of $\GL_n(A)$ via  the embedding
\begin{equation}  \label{rng0}
\GL_\ell(A) \hookrightarrow \GL_n(A), \quad a \mapsto \diag(a, 1_A, \ldots, 1_A)
\end{equation}
and put
$ \GL(A) =  \textstyle \bigcup_{n\in \NN_+} \GL_n(A \big)$,
which is canonically a group. We observe
\begin{equation} \label{rng01}
\GL(A) =  \lan B\ti : B \sim A  \big \ran,
\end{equation}
where $\lan \cdots \ran$ denotes group generation. Indeed, it is well-known (and follows from \ref{azu}\eqref{azu-Br}) that $A \sim \Mat_n(A)$, proving the inclusion from left to right. Conversely, if $B \sim A$, we know from \ref{azu}\eqref{azu-Br} that there exists $n\in \NN_+$ and an idempotent $e\in \Mat_n(A)$ such that $b = e \Mat_n(A) e$. We embed $B\ti \hookrightarrow \GL_n(A)$ by $b \mapsto b + (1_{\Mat_n(A)} - e)$. The description \eqref{rng01} implies
that $\GL(A)$ only depends on the Brauer equivalence class of $A$:
\begin{equation} \label{rng00}
A \sim B \quad \implies \quad \GL(A) = \GL(B).
\end{equation}

The equation \eqref{soex-mor1} with $A$ replaced by $\Mat_n(A)$  shows that $\Nrd$ extends to a well-defined group homomorphism $\Nrd^{\rm st} \co \GL(A) \to R\ti$. We call its image
\begin{equation} \label{rng1}
  \Nrd^{\rm st}(A\ti) = \Nrd^{\rm st}\big( \GL(A) \big) = \textstyle \bigcup_{n\in \NN_+} \Nrd_{\Mat_R(A)/R} \big(\GL_n(A \big))
\end{equation}
the {\em stable reduced norm group} of $A$. It follows from \eqref{rng01} that
\begin{equation} \label{rng-2}
 \Nrd^{\rm st}(A\ti) = \big \langle \Nrd_{B/R}(B\ti) : B \sim A  \big \ran
\end{equation}
and from \eqref{rng00} that $\Nrd^{\rm st}(A\ti) = \Nrd^{\rm st}(B\ti)$ whenever $A \sim B$.
\sm

We denote by $\rmE_n(A)$ the subgroup of elementary matrices in $\GL_n(A)$. The embedding \eqref{rng0} respects elementary matrices. We thus get the subgroup $\rmE(A) = \bigcup_{n\in \NN_+} \rmE_n(A) \subset \GL(A)$, well-known to be the commutator subgroup of $\GL(A)$. Since $\Nrd$ is trivial on any commutator, $\Nrd^{\rm st}$ descends to a well-defined surjective group homomorphism
\begin{equation} \label{rng-K}
 \rmK_1(A) = \GL(A)/\rmE(A) \twoheadrightarrow \Nrd^{\rm st}(A).
\end{equation}

Clearly, $\Nrd(A\ti) \subset \Nrd^{\rm st}(A\ti)$. In Lemma~\ref{stab} we will establish equality in two cases. We will use the subgroup
\begin{equation} \label{rng-sl}
 \SL_n(A) = \{ a\in \Mat_n(A): \Nrd_{\Mat_n(A)/R} (a) = 1\},
 \end{equation}
a normal subgroup of $\GL_n(A)$ by \eqref{rcp-22} and the group of $R$--points of a reductive group scheme of type $\rmA$, \cite[3.5.0.92]{CF}.

\comments{(?) Previously a question. Do you know something about the map \eqref{rng-K}? Is it a split surjection? This is true for $A=R$.
Are there examples with $\Nrd{A\ti } \subsetneq \Nrd^{\rm st}(A\ti)$? }
\pcomments{(2026-05-13) Dans le cas d'un corps de  base, $\rmK_1(A) = A^\times/[A^\times, A^\times]$, voir le K-book de Weibel, 1.2.4. La norme r\'eduite induit bien une surjection $\rmK_1(A) \to \Nrd(A)^\times$, dont le noyau est not\'e $SK_1(A)$. Celui-ci est trivial si l'indice de $A$ est sans facteur carr\'es, mais est non trivial en g\'en\'eral sinon (examples de Platonov). Je ne sais pas si  $\rmK_1(A)$ est le produit direct de $SK_1(A)$ par 
$\Nrd(A)^\times$. Une r\'ef\'erence \`a regarder peut-\^etre est 
l'article de A.~Rapinchuk-Segev-Seitz, JAMS 2002.}

\begin{lem} \label{stab} We have
\begin{equation}\label{rng-3} \begin{split}
 &\Nrd(A\ti) = \Nrd^{\rm st}(A\ti) \quad \iff
\\ & \quad \GL_n(A) = \GL_1(A) \cdot \SL_n(A) \quad \text{holds for all $n\in \NN_+$}.
\end{split} \end{equation}
The equality \/ $\Nrd(A\ti) = \Nrd^{\rm st}(A\ti)$ holds in each of the following cases:
 \begin{enumerate} [label={\rm (\roman*)}]
    \item \label{stab-i} $A = \End_R(M)$ where $M$ is a faithfully projective
     $R$--module containing a unimodular vector;  \sm

\item \label{stab-ii} $R$ is semilocal.
\end{enumerate}
We have $\Nrd(A\ti) = R\ti = \Nrd^{\rm st}(A\ti)$ in case \ref{stab-i}, and $\Nrd(A\ti)= \Nrd(B\ti)$ for $A\sim B$ in case \ref{stab-ii}.
\end{lem}

\begin{proof} \eqref{rng-3} Assume $\Nrd(A\ti) = \Nrd^{\rm st}(A\ti)$ and let $g\in \GL_n(A)$. By assumption there exists $a \in \GL_1(A)$ such that $\Nrd(g) = \Nrd(a)$. For $d =\diag (a, 1_A, \ldots, 1_A)\in \GL_n(A)$ we have $\Nrd(d) = \Nrd(a)$ by \eqref{rdp-mat2}, whence $d\me g \in \SL_n(A)$. Conversely,
by \eqref{rng-2}, we have to show $\Nrd\big(  \GL_n(A) \big) \subset \Nrd(A\ti)$ for all $n\in \NN_+$. This follows from multiplicativity of $\Nrd$.
\sm

\ref{stab-i} Let $x\in M$ be unimodular. Then $M = Rx \oplus N$ as $R$--modules by \ref{unimod}. For $u\in  R\ti$ we define $a_u \in \GL(M)= A\ti$ by
$a_u(x) = ux$ and $a_u |_N = \Id_N$. Then $\Nrd(a_u) = \det(a_u) =  u$, implying
\begin{equation} \label{stab-1}
 \Nrd(A\ti) =  \det \big(\GL(M) \big) = R\ti.
 \end{equation}
Since $\Nrd(A\ti) \subset \Nrd^{\rm st}(A\ti) \subset R\ti$, we are done.
\sm

\ref{stab-ii} 
We consider the reductive group scheme $\uGL_n(A)$; its group of $R$--points is $\GL_n(A)$. We denote by $\uB$ (respectively $\uB^{-}$) the parabolic $R$--subgroup scheme of upper (respectively lower) triangular matrices of $\uGL_n(A)$ and by $\rad_u(\uB)$ and $\rad_u(\uB^{-})$ the corresponding unipotent radical. By Demazure's decomposition theorem \cite[XXVI.5.2]{SGA3},
\[
\GL_n(A) = \rad_u(\uB)(R) \cdot \rad_u(\uB^{-})(R)  \cdot \uB(R).
\]
Since
\[ \uB(R) = \rad_u(\uB)(R) \cdot \begin{pmatrix}
A^\times & 0 & \cdots & 0 \\
0 & A^\times & \cdots & 0\\
\vdots &  & \ddots & \vdots  \\
0 & 0& \cdots &A^\times
\end{pmatrix}
\]
and since $\Nrd$ is trivial on $\rad_u(\uB)(R)$ and
$\rad_u(\uB^{-})(R)$,  cf.\ \ref{soex}\eqref{soex-mat}, we conclude that $\GL_n(A) = \GL_1(A) \cdot \SL_n(A)$ and are done by \ref{rng-3}. That $\Nrd(A\ti) = \Nrd(B\ti)$ for $A\sim B$, then follows from \eqref{rng-2}. \end{proof}

\comments{(2023-02-25) Need to reformulate: The norm principle follows from \ref{prop_red_stable}\ref{prop_red_stable3} and \ref{stab}.
}

\subsection{Epilogue: another approach to reduced gadgets.} \label{epi-red} As noted in \ref{rcp}, our approach to the reduced characteristic polynomial, reduced trace and reduced norm goes back to Grothendieck \cite{GroI}. Another approach is based on the Hamilton-Cayley Theorem: every matrix satisfies its characteristic polynomial.
This approach, classical for associative algebras over fields, was pursued by Jacobson \cite{Jac1, Jac2} for finite-dimensional power-associative algebras over fields, which includes Jordan algebras over fields of characteristic $\ne 2$. That this second approach can also be followed in the setting of Azumaya algebras over schemes, has been pointed out by Grothendieck in \cite[5.13]{GroI}. Jacobson's work has been taken up by Loos in \cite{Lo-genalg}, where he considers Jordan algebras over rings $R$ whose underlying $R$--modules are finite projective and which are generically algebraic, thus in particular Azumaya algebras. In Loos' setting, a version of Lemma~\ref{prop_nred} is proven in \cite[2.9]{Lo-genalg}.


\section{Norm groups and \'etale norm  groups} \label{sec:normgroups}

\comments{(2026-03-30) We have decided yesterday to include this section on norm groups, since Philippe needs Proposition~\ref{prop_norm_SB} for local rings. 

This file needs results from the previous sections \S\ref{sec:semilocal}, \S\ref{sec:trans-semi}, \S\ref{sec:redNorm}, the results on GIT quotients \ref{gitqu}--\ref{thm_lee}, and results from \S\ref{sec:kneb} on Knebusch's norm principle (still under revision)}

\comments{(2026-03-30) This file was written in our ``pre-LG perod". It could well be that some results stated for semilocal rings are in fact true for LG rings. 

I changed \ref{SBs-lem}\eqref{SBs-lem-c} from semilocal to unimodular. 

We still use semilocal in Lemma~\ref{lem_regular2}, Proposition~\ref{prop_norm_SB}, Lemma~\ref{lem_norm_rk3}, and  Theorem~\ref{prop_norm_quad},}

Unless specified otherwise, in this section $A$ denotes an Azumaya $R$--algebra, not necessarily of constant rank.

The goal of this section is to study norm groups and \'etale norm groups of Severi-Brauer schemes and quadrics. We start with an introduction to Severi-Brauer schemes in \ref{SBs} and to norm groups for arbitrary schemes in \ref{nogr}.

\comments{(2026-03-30) I took the subsection \ref{SBs} from an old file: 'NormGroups-2021-06-03.tex', so that this section can be read independently of later generalizations of $\uSB(A)$. }

\subsection{Severi-Brauer schemes}\label{SBs} Throughout, $A$ is an Azumaya $R$--algebra of degree $\deg_R \co \Spec(R) \to \ZZ$, thus $(\deg_R A)^2 = \rank_R A$. We mainly follow van den Bergh's approach \cite{VB}. We first recall three descriptions of the $R$--points of the Severi-Brauer scheme of $A$. \sm

\begin{inparaenum}[(a)] \item We consider pairs $(P, \vphi)$ consisting of  a left $A$--module $P$, which is finite projective as $R$--module and satisfies $\rank_R P = \deg_R A$, and a surjective $A$--module homomorphism $\vphi \co A \to P$.
 Two such pairs $(P, \vphi)$ and $(P', \vphi')$ are called equivalent, if there exists an $A$--module isomorphism $u \co P \simlgr P'$ satisfying $\vphi' = u \circ \vphi$. We denote by $[P, \vphi]$ the equivalence class of a pair $(P, \vphi)$ and by $\SB(A)$ the set of equivalence classes of pairs $(P, \vphi$ as defined above. \sm

\item \label{SBs-c} We mention another description of $\SB(A)$.
We consider triples $(P, \theta, x)$ consisting of a finite projective $R$--module $P$ with $\rank_R P = \deg_R A$, a unital $R$--algebra homomorphism $\theta \co A \to \End_R(P)$, and an element $x\in P$ satisfying $\theta(A)\, Rx = P$. Two such triples $(P, \theta, x)$ and $(P', \theta', x')$ are called equivalent if there exists an $R$--module isomorphism $u \co P \to P'$ satisfying $u(x) = x')$ and $u \circ \theta(a) \circ u\me = \theta'$ for all $a\in A$. We let $\SB'(A)$ be the set of equivalence classes of triples $(P, \theta, x)$ as defined above.

It is shown in \cite[Lem.~3]{VB} (for more general algebras than Azumaya algebras) that
\[ (P, \vphi) \mapsto \big( \sfR_{A/R}(P), \theta, \vphi(1) \big), \quad  \theta(a) (p) = a \cdot p,  \quad (a\in A, p\in P) \]
induces a bijection $\SB(A) \simlgr \SB'(A)$. For the construction of the inverse map observe that the map $\theta$ of any triple $(P, \theta, x)$ makes $P$ an $A$--module via $a\cdot p = \theta(a) (p)$ such that $\vphi \co A \to P$, $a \mapsto \theta(a)(x)$ is a surjective $A$--module homomorphism.
\sm

\item \label{SBsn} We give a third description of $\SB(A)$ which is more in the spirit of the classical approach to Severi-Brauer varieties. We denote by $\SB''(A)$ the set of left ideals $I\subset A$, for which the $R$--module $A/I$ is finite projective of $\rank_R (A/I) = \deg_R A$. Given such an $I$, the canonical map $\can_I \co A \to A/I$ is a surjective $A$--module homomorphism, giving rise to a set map
\begin{equation} \label{SBs1}
 \SB''(A) \simlgr \SB(A), \quad I \mapsto [A/I, \can_I]
\end{equation}
which is easily seen to be a bijection with inverse map $[P, \vphi] \mapsto \Ker(\vphi)$.
\lv{
Details: The inverse map is well-defined since $\Ker(\vphi) = \Ker(\vphi')$ if $(P, \vphi) \cong (P', \vphi')$; also $\Ker(\vphi)$ is a left ideal of $A$ and it has the correct rank. We have $I \mapsto [A/I, \can_I] \mapsto \Ker(\can_I) = I$; and $[P, \vphi] \mapsto \Ker(\vphi) \mapsto [A/\Ker (\vphi), \can_I] = [P, \vphi]$ because of commutativity of the diagram
\[ \xymatrix@C=40pt{A \ar[r]^\vphi \ar[dr]_{\can_I} & P \\ & A/I \ar[u]_\cong } \]
}
Any $I \in \SB''(A)$ gives rise to a splitting
\[ \al \co A \simlgr \End_R(A/I), \quad \al(a) = (a'+ I \mapsto aa' + I)\]
since $\al$ is an $R$--algebra homomorphism between Azumaya algebras of the same rank, cf.\ \ref{azu}\eqref{azu-hom}.
\sm

\item \label{SBs-LGc} Finally, we define an $R$--functor $\underline{\SB}(A)$ by $\underline{\SB}(A)(S) = \SB(A\ot_R S)$ for any $S\in \Ralg$. By \cite{VB} this functor is representable by an $R$--scheme $\uSB(A)$, the {\em Severi-Brauer scheme of $A$}.  
\end{inparaenum}

\subsection{Norm groups -- definitions and basic properties.} \label{nogr}
\begin{inparaenum}[(a)] \item\label{def_norm_field}  ({\em The Kato-Saito definition}) The definition of  a norm group goes back to \cite[\S7]{KaSa}. It is therefore appropriate to first recall their definition. We recall that $\rmN_{S/R}$ denotes the norm of a finite projective extension $S$ of $R$, \ref{trno}.

Let $k$ be a field and let $Y$ be a $k$--scheme, locally of finite type.
The {\em norm group\/} of $Y$ is the subgroup $\rmN^\sharp_Y(k)$ of $k\ti $ generated by the subgroups $\rmN_{L/k}(L^\times)$ where $L$ varies over the finite field extensions $L$ of $k$ satisfying $Y(L) = \Mor(\Spec(L), Y)\ne \emptyset$. Provisionally, we define the {\em \'etale norm group\/} $(\rmN_Y^{\sharp}){}^{\rmet}(k)$ as the subgroup of $k\ti$ generated by all $\rmN_{E/k}(E\ti)$ where $E/k$ is a finite separable field extension with $Y(E) \ne \emptyset$. We will show in Lemma~\ref{lem_norm_field} that the two groups coincide.
\sm

\inparcom{(2021-0519) I deleted the assumption in \eqref{nogr-def} that $X$ be locally of finite presentation, since I did not see where it is used or useful}

\item \label{nogr-def} ({\em Norm groups for $R$--schemes}) As usual, let $R$ be an arbitrary base ring. By definition, the {\em norm group\/} of an $R$--scheme $X$ is the subgroup $\rmN_X(R)$ of $R^\times$  generated by the subgroups
$\rmN_{R'/R}\bigl( (R')^\times \bigr)$, where $R'$
varies over the  faithfully projective  extensions $R'$ of  $R$ for which $X(R') \not=\emptyset$.

Similarly, we define the {\em \'etale norm group} of $X$ as the
subgroup $\rmN_X^{\rmet}(R)$ of $R^\times$  generated by the subgroups
$\rmN_{R'/R}\bigl( (R')^\times \bigr)$ where $R'$
varies over the finite \'etale extensions $R'$ of $R$ of positive rank which satisfy $X(R') \not=\emptyset$.

In the remainder of this subsection we collect some basic properties of norm groups.
\sm

\item \label{nogr-c}  Of course $\rmN^{\et}_X(R) \subset \rmN_X(R)$; one of the themes of this section is that these groups coincide for certain schemes, see  \ref{lem_norm0},  \ref{prop_norm_SB} and  \ref{lem_norm_rk2}--\ref{prop_norm_quad}.
It is also clear that whenever $X(R) \ne \emptyset$, then $\rmN^{\et}_X(R)= \rmN_X(R)= R^\times$.
\sm

\item\label{nogr-nor} ({\em Direct products}) Let $R'= R'_1\times R'_2$ be a direct product of $R$--algebras. Since $X(R')= X(R'_1) \times X(R'_2)$ we have $X(R') \ne \emptyset \iff X(R'_1) \ne \emptyset \ne X(R'_2)$. Also,
    the product formula \ref{trno}\eqref{trno-c} for norms says
    \[ \rmN_{(R'_1 \times R'_2)/R}\big( (R')^\times \big) = \rmN_{R'_1/R}\big( (R_1')^\times \big) \cdot \rmN_{R'_2/R}\big( (R_2')^\times \big).\]
    Since $R'$ is faithfully projective as $R$--module if and only if $R_1'$ and $R'_2$ are finite projective and one of $R'_1$ or $R_2')$ is faithful, it follows that every element of $\rmN_X(R)$ has the form $\rmN_{R'/R}\big( y\big)$ for a faithfully projective $R'\in \Ralg$. The analogous considerations hold for finite \'etale extensions of positive rank. \sm

\item \label{nogr-bc} ({\em Base change, functoriality}) The assignments
\[ S \mapsto \rmN_X(S):= \rmN_{X_S}(S) \quad \text{and} \quad
 S \mapsto \rmN^{\et}_X(S) := \rmN_{X_S}^{\et}(S) \]
 extend to $R$--functors with values in the category of groups. \sm

\item \label{nogr-f} Let $Y  \hookrightarrow X$ be a closed or open imbedding of $R$--schemes. It is immediate from the definition that $\rmN_Y^{\et}(R) \subset \rmN_X^{\et}(R)$. This obvious relation will turn out to be very helpful later on. 
\end{inparaenum}
\ms

Let us verify that our definition \eqref{nogr-def} for $R$--schemes is compatible with the definition \eqref{def_norm_field} in the case of $R$ being a field.

\begin{lem} \label{lem_norm_field}
Let $k$ be a field and let $Y$ be a $k$-scheme locally of finite type.
Then $\rmN^\sharp_Y(k) = \rmN_Y(k)$ and $(\rmN_Y^\sharp){}^{\rmet}(k) = \rmN_Y^{\rmet}(k)$.
\end{lem}

\begin{proof} Clearly $\rmN^\sharp_Y(k) \subset \rmN_Y(k)$ and $(\rmN_Y^\sharp){}^{\rmet}(k) \subset \rmN_Y^{\rmet}(k)$. 

To prove $\rmN_Y(k) \subset  \rmN^\sharp_Y(k)$, recall that the group $\rmN_Y(k)$ is generated by the $\rmN_{R/k}(R\ti)$ for $R$ running over the finite $k$-algebras $R$ satisfying $Y(R) \not = \emptyset$.
Such an $R$ is artinian and is a product of local artinian $k$--algebras. 
Without loss of generality, we can therefore deal with a finite local artinian $k$-algebra $R$ satisfying $Y(R) \not = \emptyset$.  We denote by $\gm$ the maximal ideal of $R$ and by $L=R/\gm$ its residue field. Proposition \ref{prop_norm} provides a positive integer $d$ such that $\rmN_{R/k}(R^\times) \subseteq \bigl( \rmN_{L/k}(L^\times) \bigr){}^d$.  A fortiori we have  $\rmN_{R/k}(R^\times) \subset  \rmN_{L/k}(L^\times)$.
By using  the specialization map $Y(R) \to Y(L)$, it follows that 
 $ Y(L) \not = \emptyset$, so that  $\rmN_{L/k}(L^\times) \subseteq \rmN^\sharp_Y(k)$.
Thus  $\rmN_{R/k}(R^\times) \subset \rmN^\sharp_Y(k)$ as desired. 

The proof of $\rmN_Y^{\rmet}(k) \subset (\rmN_Y^\sharp){}^{\rmet}(k)$ is a simple modification of the proof above. As before, it suffices to consider an \'etale local artinian $k$--algebra $R$ with residue field $L$. Since $L$ is finite separable over $k$, we can conclude as above.  \end{proof}

\lv {OLD PROOF for the \'etale norm groups: To prove the other inclusion for the \'etale norm groups, we are given a finite \'etale $k$--algebra $B$ of positive dimension such that $Y(B) \not = \emptyset$ and we have to show that
 $\rmN_{B/k}(B^\times) \subseteq (\rmN^\sharp_Y){}^\rmet(k)$.
Since $B$ is an artinian $k$--algebra, it is a direct product of local \'etale artinian $k$--algebras, each of positive dimension. Hence, cf.\ \ref{nogr}\eqref{nogr-nor}, without loss of generality we can assume that $B$ is local. We denote  by $L$ the residue field of  $B$; it is finite \'etale over $k$, i.e., finite separable, and satisfies $Y(L) \not = \emptyset$.
According to \cite[\S13.8, Prop.~7(b)]{BA8}, the map $B \to L$ admits a splitting which makes $B$ an $L$-algebra. Since $\rmN_{B/k}(B^\times)=\rmN_{L/k}\big( \rmN_{B/L}(B^\times) \big)$, it  follows that $\rmN_{B/k}(B^\times) \subset \rmN_{L/k}(L^\times) \subset (\rmN^\sharp_Y){}^\rmet(k)$.}

\sm
We note that in the \'etale case one can replace the reference \ref{prop_norm} by the reference \cite[\S13.8, Prop.~7(b)]{BA8}, saying that $R \to L$ admits a splitting. 

We give another instructive example of norm groups below.

\begin{lem} \label{lem_norm0} Let $T\in \Ralg$ be a faithfully projective $R$--algebra, and let $X=\Spec(T)$. Then
\[ \rmN^{\et}_X(R) \subset \rmN_X(R) \subset \rmN_{T/R}(T\ti),\]
and all inclusions are equalities if $T$ is finite \'etale of positive rank.
\end{lem}

\begin{proof} Regarding the inclusions, only  $\rmN_X(R) \subset \rmN_{T/R}(T\ti)$ needs to be shown. To this end, let $R'$ be a faithfully projective $R$--algebra
satisfying  $X(R') \not =\emptyset$. Thus, denoting by $\si_T \co R \to T$ and $\si_{R'} \co R \to R'$ the structure maps of the $R$--algebras, there exists an $R$--algebra homomorphism $f\co T \to R'$ satisfying $f \circ \si_T = \si_{R'}$
because $X(R') \cong \Hom_{\Ralg}(T,R')$.
Hence the $R$--algebra  homomorphism $R' \to T'= T \ot_R R'$, $r' \mapsto 1_T \ot r'$ admits an $R'$--linear retraction, given by $a\ot r' \mapsto f(a)r'$. Equivalently, there exists an isomorphism $T' \cong R' \times B$, where  $B$ is an  $R'$--algebra which is locally free of finite rank. Since $T'{}\ti \cong R'{}\ti \times B\ti$ we get $R'{}\ti \subset \rmN_{T'/R'}(T'{}\ti) $ and therefore have the  inclusions
\begin{align*}
\rmN_{R'/R}(R'{}\ti) & \subset \rmN_{R'/R}\bigl( \rmN_{T'/R'}(T'{}\ti)\bigr)
= \rmN_{T'/R}(T'{}\ti) \\& = \rmN_{T/R}\big( \rmN_{T'/T}(T'{}\ti)\bigr) \subset  \rmN_{T/R}(T\ti).
\end{align*}
This implies $\rmN_X(R) \subseteq \rmN_{T/R}(T^\times)$.

If $T$ is finite \'etale of positive rank, $\Id_T$ defines a point of $X(T)$, so that $\rmN_{T/R}(T^\times) \subset \rmN_X^{\et}(R)$ follows. \end{proof}
\sm

We present two more properties of norm groups for schemes.

\begin{lem} \label{lemn_norm1} Let $X$ be an $R$--scheme. Recall that $\rmN_X(S) = \rmN_{X\times_R S}(S)$ for any $S\in \Ralg$. \sm

\begin{inparaenum}[\rm (a)]
 \item \label{lemn_norm1-a} Let $S$ be a faithfully projective $R$--algebra.
Then \[ \rmN_{S/R} \bigl( \rmN_{X}(S) \bigr) \, \subset \, \rmN_X(R),\]
and if $S$ is finite \'etale of positive degree, then also
\[ \rmN_{S/R} \bigl( \rmN^{\et}_{X}(S) \bigr) \, \subset \, \rmN^{\et}_X(R).  \]

\item \label{lemn_norm1-b} Let $X$ be locally of finite presentation and let $R_\infty = \varinjlim  R_\al$ be a direct limit (= colimit) of a directed system $(R_\al, \vphi_{\be\al})_{\al \in I}$ of algebras in $\Ralg$. 
    Then   $\big(\rmN_X(R_\al), \rmN_X(\vphi_{\be \al})\big)_{\al \in I} $ is a directed system of groups whose direct limit satisfies
    \[ \textstyle \varinjlim _{\al \in I} \rmN_X(R_\al) \cong \rmN_X(\varinjlim_\al R_\al).\]
     Analogously, $\varinjlim _{\al \in I} \rmN_X^{\et}(R_\al) \cong \rmN^{\et}_X(R_\infty)$.
\end{inparaenum}
\end{lem}

\begin{proof}
  \eqref{lemn_norm1-a} Let $s\in \rmN_X(S)$. Thus, by \ref{nogr}\eqref{nogr-nor}, $s= \rmN_{S'/S}(y)$ for some $y\in (S'){}\ti$ and $S'$ a faithfully projective $S$--algebra with $(X\times_R S)(S') \ne \emptyset$. Since we have a morphism $X\times_R S \to X$, we get $X\big( \sfr_{S/R}(S')\big) \ne \emptyset$.
  Because the $R$--algebra $\frR_{S/R}(S')$ is faithfully projective as $R$--module, transitivity of norms,  \ref{trno}\eqref{trno-b},  yields $\rmN_{S/R}(s) = \rmN_{S/R}\big( \rmN_{S'/S}(y)\big) = \rmN_{S'/R}(y) \in \rmN_X(R)$. The same proof works for the \'etale norm groups, using transitivity for finite \'etale algebras,  \ref{fea}\eqref{fea-trans}. \sm

\eqref{lemn_norm1-b} By functoriality \ref{nogr}\eqref{nogr-bc}, we get a directed system  $\big(\rmN_X(R_\al), \rmN_X(\vphi_{\be \al})\big)$ and a unique group homomorphism $f$ making the diagram below commutative:
\[
\vcenter{\xymatrix@C=10pt{
     R_\al \ar[dd]_{\vphi_{\be\al}} \ar[dr]^{\vphi_\al} \\
      &\limind R_\al  \\
     R_\be  \ar[ur]_{\vphi_\be} }}
\quad \rightsquigarrow \quad
\vcenter{\xymatrix@C=10pt{
     \rmN_X(R_\al) \ar[dd]_{\rmN_X(\vphi_{\be\al})} \ar[dr]
     \ar@/^/[drrr]^{\rmN_X(\vphi_\al)}\\
      &\limind \rmN_X(R_\al) \ar@{-->}[rr]^{\exists ! \, f} &&
        \rmN_X(R_\infty) \\
     \rmN_X(R_\be)  \ar[ur] \ar@/_/[urrr]_{\rmN_X(\vphi_\be)}
}}.
\]
Observe that the maps $\rmN_X(\vphi_{\be \al})$ and $\rmN_X(\vphi_\al)$ are the restrictions of $\vphi_{\be \al}$ and $\vphi_\al$ to the respective norm groups. This implies that $f$ is injective. To establish surjectivity of $f$, let $x \in \rmN_X(R_\infty)$. Thus, there exists a faithfully projective extension $S_\infty/R_\infty$ such that $x= \rmN_{S_\infty/R_\infty}(y)$ with $y \in S_\infty^\times$ and $X(S_\infty) \not = \emptyset$. According to \cite[IV$_3$, 8.5.5]{EGA},  
there exists $\alpha \in I$ and a faithfully projective $R_\alpha$--algebra  $S_\alpha$ such that $S_\infty= S_\alpha \otimes_{R_\alpha} R_\infty$.
For $\beta \geq \alpha$, we put $S_\beta= S_\alpha  \otimes_{R_\alpha} R_\beta$,
so that $S_\infty= \varinjlim S_\beta$.
Since $X$ is locally of finite presentation over $R$,  we have
$X(S_\infty)= \varinjlim  X(S_\beta)$ by \ref{ag}\eqref{ag-lp}.
Hence, after changing $\alpha$, we can assume that
$X(S_\alpha) \not = \emptyset$ and that $y$ arises from
some $y_\alpha \in S_\alpha$. It follows that $x\in  \rmN_{S_\infty/R_\infty}(y)$
arises from  $\rmN_{S_\alpha/R_\alpha}(y_\alpha)\in \rmN_X (R_\alpha)$,
hence  belongs to the image of $f$, thus finishing the proof for the norm groups.
The analogous proof works for the \'etale norm group using
\cite[IV$_4$, 17.7.8(ii)]{EGA}.
\end{proof}

In \ref{prop_norm_SB} we will determine the norm groups and \'etale norm groups of the Severi-Brauer scheme $\uSB(A)$ of an Azumaya $R$--algebra $A$ for $R$ semilocal, see \ref{SBs}\eqref{SBs-LGc} for the definition of $\uSB(A)$.  The following Lemma~\ref{SBs-lem} relates neutrality of an Azumaya algebra, cf.\ \ref{azu}\eqref{azu-spli}, to points of Severi-Brauer schemes.

\comments{(2026-03-30) Previously, \ref{SBs-lem}\eqref{SBs-lem-c} was stated for semilocal rings; I changed to unimodular. Also, the old proof of \ref{SBs-lem}\eqref{SBs-lem-c} did not sense. } 

\begin{lem}\label{SBs-lem}  Let $X = \uSB(A)$ be the Severi-Brauer scheme of the Azumaya $R$--algebra $A$. \sm

\begin{inparaenum}[\rm (a)]

\item \label{SBS-2} If $X(R) \ne \emptyset$, then $A$ is neutral, i.e., $A\sim  R$. \sm

\item \label{SBs-lem-b} If $A\sim R$, say $A = \End_R(M)$ for a faithfully projective $R$--module, then $X(R) \ne \emptyset \iff M$ contains a complemented rank-$1$-submodule. \sm

\item \label{SBs-lem-c} 
\new 
Let $R$ be a unimodular ring, {\rm \ref{unifap}}. Then $X(R) \ne \emptyset \iff A$ is neutral.
\enew
\sm

\item\label{Sbs-e} Let $S\subset A$ be a maximal commutative $R$--subalgebra, for which the canonical $S$--module $A$ is faithfully projective, e.g., assume that $S$ is maximal \'etale. Then $X(S) \ne \emptyset$.
\end{inparaenum} \end{lem}

\begin{proof} \eqref{SBS-2} If $X(R) \ne \emptyset$, then \ref{SBs}\eqref{SBsn} shows that $A$ contains  a left ideal $L$ such that $A/L$ is a projective $R$--module of rank $\deg_R(A)$. The left multiplication then gives rise to an $R$--algebra homomorphism $A \to \End_R(A/L)$, which is an isomorphism by \ref{azu}\eqref{azu-hom}.
\lv{OLD PROOF: 
Let $[(P, \theta, x)]\in \SB'(A)$. Since $\deg_R A >0$, the $R$--module $P$ is faithfully projective and $\theta \co A \to \End_R(P)$ is a neutralization by \ref{azu}\eqref{azu-hom}. }
\sm

\eqref{SBs-lem-b} We know from \ref{grap}\eqref{grape-dual} 
that $X(R) \ne \emptyset \iff \PP(M) \ne \emptyset$. \sm

\eqref{SBs-lem-c} The implication $\Rightarrow$ follows from \eqref{SBS-2}. 
\new
Conversely, let $A$ be neutral, i.e., $A\cong \End_R(M)$ for some faithfully projective $R$--module $M$. By definition of unimodularity, $M$ contains a unimodular vector, implying $X(R) \ne \emptyset$ by \eqref{SBs-lem-b}. 
\enew
\sm

\eqref{Sbs-e} Let $P$ be the $A_S$--module defined on $A$ by $(a\ot s)(p) = aps$. Then $A_S \to P$, $a\ot s \mapsto as$, is an epimorphism of $A_S$--modules.
Let $\vphi \co \Spec(S) \to \Spec(R)$ be the canonical map. Then $\rank_S P =  (\rank_R S) \circ \vphi = (\deg_R A)\circ \vphi = \deg_S (A\ot_R S)$ follows from \eqref{azu-2}, so that $(P, \vphi)$ gives rise to an element of $\SB'(A_S) \cong \uSB(S)$. If $S\subset A$ is maximal \'etale, then the $S$--module $A$ is faithfully projective by \ref{lemsepsub}\ref{lemsepsup-ai}. %
\end{proof}
\sm

\textbf{Remarks.} The implication in \ref{SBs-lem}\eqref{SBS-2} cannot be reversed. Indeed, if $A=\End_R(M)$ and $X(R) \ne \emptyset$, then $M$ is decomposable by \eqref{SBs-lem-b}. But there exist indecomposable faithfully projective $R$--modules, for example the tangent bundle of the $2$--dimensional real sphere is such an example (\cite[Remark after I, Prop.~4.15]{La}).
\ms

Lemma~\ref{lem_norm_stab} collects some preliminary results on norm groups and \'etale norm groups. The stable reduced norm group $ \Nrd^{\rm st}(A\ti)$ is defined in \eqref{rng1}:
\[
   \Nrd^{\rm st}(A\ti) = \textstyle \bigcup_{n\in \NN_+} \Nrd_{\Mat_R(A)/R} \big(\GL_n(A \big)) = \big \langle \Nrd_{B/R}(B\ti) : B \sim A  \big \ran.
\]

%
%

\begin{lem}\label{lem_norm_stab}  Let $X = \uSB(A)$ be the Severi-Brauer scheme of the Azumaya $R$--algebra $A$. \sm

\begin{inparaenum}[\rm (a)]
\item \label{lem_norm_stab-a} 
We always have $\rmN_X(R) \subset \Nrd^{\rm st}(A^\times)$. \sm

\item \label{SBs-lem-bb} If $A = \End_R(M)$ where $M$ is a faithfully projective $R$--module containing a unimodular vector, thus $X(R)\ne \emptyset$ by {\rm \ref{SBs-lem}\eqref{SBs-lem-b}}, then
  \begin{equation} \label{SBs-lem-b1}
     \Nrd^{\rm st}(A\ti)= \Nrd(A\ti) = R\ti = \rmN^{\et}_X(R)= \rmN_X(R). \end{equation}
\end{inparaenum}
 \end{lem}

\begin{proof} \eqref{lem_norm_stab-a} Let $S\in \Ralg$ be a faithfully  projective $R$--algebra such that $X(S) \ne  \emptyset$. We want to show that $\rmN_{S/R}(S^\times) \subset  \Nrd^{\rm st}(A^\times)$. By \ref{SBs-lem}\eqref{SBS-2} the Azumaya $S$--algebra $A_S$ is neutral. Hence, by \cite[7.4.3]{Ford}, there exists an Azumaya $R$--algebra $B$
which is Brauer-equivalent to $A$ and contains $S$ as a maximal commutative subalgebra such that $B$ is a projective (left) $S$--module. Hence $\rmN_{S/R}(S\ti) \subset \Nrd_{B/R}(B\ti) \subset \Nrd^{\rm st}(A\ti)$ by \ref{prop_nred} and \eqref{rng-2}.%
\sm

\eqref{SBs-lem-b} The equations $ \Nrd^{\rm st}(A\ti)= \Nrd(A\ti) = R\ti$ follow from \ref{stab}, and the equations  $R\ti = \rmN^{\et}_X(R)= \rmN_X(R)$ from \ref{nogr}\eqref{nogr-c}.
\end{proof}
\sm

If $R$ is semilocal, we will show in Proposition~\ref{prop_norm_SB} that the inclusion in Lemma~\ref{lem_norm_stab}\eqref{lem_norm_stab-a} is an equality.
To do so, we will use the concept of semisimple regular elements in $\uGL_1(A)$, see \ref{prop_regular-ex} and \ref{lem_regular_azumaya}, as well as Lemma~\ref{lem_regular2}.

\begin{lem} \label{lem_regular2} Let $A$ be an Azumaya algebra over the\/ {\em semilocal} base ring $R$.  \begin{enumerate} [label={\rm (\alph*)}]

 \item \label{lem_regular2_1} If all residue fields of $R$ are infinite, each element of $A\ti$ is a product of two semisimple regular elements. \sm

\item \label{lem_regular2_2}  If all residue fields $\ka(\m)=R/\m$ of $R$, $\m$ a maximal ideal of $R$, satisfy $|\ka(\m)| \geq 1 + \deg \big(A\ot_R \ka(\m)\big)$, then every element of $A^\times$ is a product of an even number of semisimple regular elements, in particular $A\ti$ contains semisimple regular elements.
 \end{enumerate}
\end{lem}

\begin{proof} Let $\m_1,\dots, \m_s$ be the maximal ideals of $R$ and let $\ka_i=R/\m_i$ for $i=1,...,s$ be the residue fields of $R$. For $a \in A\ti$ we put $a_i = a \ot 1_{\ka_i} \in A_{\ka_i}\ti$. We abbreviate $U = \uGL_1(A)\ssr$.
\sm

\ref{lem_regular2_1} We know from \ref{srgo}\eqref{srgo-c} that $U$ is an open subscheme of $\uGL_1(A)$ and from \ref{trno}\eqref{trno-ds} that $\uGL_1(A)$ is open in the affine space $\uW(A)$. Hence $U$ is an open subscheme of $\uW(A)$.
The canonical map $A \to A/\Jac(R) A = \prod_i A \ot_R \ka_i$ is surjective. Since semisimple regular elements are stable under base change (\ref{srgo}\ref{srgo-cii}), we conclude from \ref{opema-rem} that
  \begin{equation}\label{rse-1}
     U(R) \twoheadrightarrow \textstyle\prod_i U(\ka_i)
     \end{equation}
is surjective.

If $R$ is an infinite field, then \ref{adc} applies: $U(R)$ is Zariski dense in $U$. It follows that there exists $b\in U(R)$ such that $a b\me \in U(R)$. This settles the claim in the field case.
For a general semilocal $R$ we then get that there exist $b_i$, $c_i \in U(\ka_i) = U_{\ka_i}(\ka_i)$ such that $a_i = b_i c_i$. By \eqref{rse-1} there exists $b\in U(R)$ mapping onto $(b_i)_{i=1,..,s}$. The element $c=b^{-1} a\in A\ti$ specializes to $(c_i)_{i=1,\ldots,s}$, hence $c \in U(R)$ by \ref{opemaLG}.
\sm

\ref{lem_regular2_2} We first consider the case of a field $R=k$ and show that
every element in $A\ti$ is a product of an even number of semisimple regular elements. By \ref{lem_regular2_1} we can assume that $k$ is finite. Hence $A= \Mat_n(k)$ by a theorem of Wedderburn (\cite[7.5.4]{Ford}). It is known 
that  $\GL_n(k)$ is generated by the elementary matrices
and $\diag(v,1, \ldots, 1)$ where $v$ is a generator of $k\ti$.
Since $|k\ti| \geq n$ by assumption, the element  $d=\mathrm{diag}( 1,v, \dots, v^{n-1})\in A\ti$ is semisimple regular by \ref{seo}\eqref{ex_linear}.
If $e$ is an elementary matrix, then $e= d \, (d^{-1}\, e)$ is  a product of two semisimple regular elements. On the other hand, we have  $\diag(v,1 \ldots, 1)=  \diag(v, v^2,  \ldots, v^n) \cdot \diag( 1, v^2, \dots, v^{n})^{-1}$, so this element is also a product of two semisimple regular elements. Summarizing, any element of $A\ti$ is a product of an even number of semisimple regular elements in $A\ti$.

Let now $R$ be semilocal. By what we have already shown, every $U(\ka_i) \ne \emptyset$. Hence, by \eqref{rse-1}, there exists $g \in U(R)$.
Therefore, up to using products of $g g^{-1}$, there exists an even integer $N \geq 2$ such that for all $i$, $1\le i \le s$,
\[ a_i =  b_{i,1} \, b_{i,2} \cdots b_{i,N}, \]
where $b_{i,j} \in U(\ka_i)$ for $j=1,...,N$. For $j=1, \dots, N-1$ let $b_j \in U(R)$ be a lift of $(b_{1,j}, b_{2,j}, \dots,  b_{s,j})$
and define $b_N$ by $a=b_1 b_2 \dots b_{N-1} b_N$.
Since $b_N$ maps onto $(b_{1,N}, b_{2,N}, \dots,  b_{s,N})$,
it is semisimple regular, proving our claim.
\end{proof}
\sm

\textbf{Example.} The bounds in \ref{lem_regular2}\ref{lem_regular2_2} are sharp. For example, let $\FF_2$ be the field of two elements and let $A= \Mat_2(\FF_2)$. By the Hamilton-Cayley Theorem $a^2 - \Tr(a) + \det(a)E_2= 0$ holds for every $a\in A$, \ref{rcp}. Hence every $a\ne E_2$ generates a $2$-dimensional subalgebra $\FF_2[a] = \FF_2[X]/\big(X^2 - \Tr(a)X + \det(a)\big)$. By \ref{srel}\ref{srel-iii} and \ref{seo}\eqref{ex_linear}, $a$ is semisimple regular if and only if $\FF_2[a]$ is \'etale, which is the case if and only if $\Tr(a)^2 - 4 \det(a) \ne 0$, equivalently $\Tr(a) = 1$. Thus, the only semisimple regular elements in $\Mat_2(\FF_2)$ are $a_1 = \ppmatrix 1 1 1 0 $ and $a_2 = \ppmatrix 0 1 1 1$. They generate a subgroup $H$ of order $3$ in $\GL_2(\FF_2)$. But $|\GL_2(\FF_2)| = 6$, for example $\ppmatrix 0 1 1 0  \notin H$.
\ms

The following Proposition~\ref{prop_norm_SB} generalizes \cite[2.6.4, 2.6.6]{GS} and \cite[4.1.8]{G2019}, which treat the case of $R$ being a field. Again for $R$ a field, a proof of the equality $\Nrd(A\ti) = \rmN_X(R)$ in \ref{prop_norm_SB}\ref{prop_norm_SB3} is sketched in \cite[Rem.~p.~434]{Collio-Sansuc} and in \cite[Lem.~10]{KaSa}.

\begin{prop}[$R$ semilocal] \label{prop_norm_SB}
Assume that $R$ is semilocal and let $X=\uSB(A)$ be the Severi-Brauer scheme of the Azumaya $R$--algebra $A$. \sm
\begin{enumerate} [label={\rm (\alph*)}]

\item \label{prop_norm_SB1} $\Nrd(A^\times)=  \Nrd^{\rm st}(A^\times) = \Nrd(B\ti)$ for $B\sim A$.\sm

\item \label{prop_norm_SB2} $\rmN_{S/R}\big( \Nrd(A_S^\times) \big)  \subset  \Nrd(A^\times)$ for $S$ a faithfully projective $R$--algebra.%
     \sm

\item \label{prop_norm_SB3} $\Nrd(A^\times)= \rmN^{\et}_X(R)  = \rmN_X(R)$.
\end{enumerate}
\end{prop}

\begin{proof} \ref {prop_norm_SB1} is a restatement of Lemma~\ref{stab}. \sm

\noindent \ref{prop_norm_SB2} We have $\rmN_{S/R}\big( \Nrd(A_S^\times) \big)  \subset  \Nrd^{\rm st}(A^\times)$ according to \ref{prop_red_stable}\ref{prop_red_stable3} applied with $Z=S$. Hence  $\rmN_{S/R}\big( \Nrd(A_S^\times) \big)  \subset  \Nrd(A^\times)$ by
\ref{prop_norm_SB1}.
\sm

\noindent \ref{prop_norm_SB3}
Lemma~\ref{lem_norm_stab}\eqref{lem_norm_stab-a} shows that
$\rmN_X(R) \subset \Nrd^{\rm st}(A^\times)$, so we have
the inclusions $\rmN^{\et}_X(R) \subset \rmN_X(R)  \subset \Nrd(A^\times)$
by \ref{prop_norm_SB1}. Hence it remains to show that $ \Nrd(A^\times) \subseteq  \rmN^{\et}_X(R)$. \sm

Let $U$ be the open $R$--subscheme of $\uGL_1(A)$  consisting
of the regular semisimple elements of $A$, and let $a \in A^\times$. We want to show that $\Nrd(a) \in \rmN_X^{\et}(R)$.\sm

{\em First step: all residue fields of $R$ have cardinality $\geq 1 + N$, for $N = \max\{\deg A_\ka: \ka \text{ residue field of $R$}\}$.} If $a\in U(R)$, then $S= R[a]$  is a maximal \'etale subalgebra of $A$ by \ref{lem_regular_azumaya}. Hence $X(S) \ne \emptyset$ by \ref{SBs-lem}\eqref{Sbs-e}. Therefore $\Nrd_A(a)= \rmN_{S/R}(a) \in \rmN_X^{\et}(R)$ by \ref{prop_nred}. Consequently, if $\lan U(R) \ran$ denotes the subgroup of $A\ti$ generated by $U(R)$, then $\Nrd(\lan U(R)\ran) \subset \rmN^{\et}_X(R)$. But $\lan U(R) \ran = A\ti$ by Lemma~\ref{lem_regular2} under the assumption of this  first step.
\sm

{\em Second step: $R$ is the semilocalization of a finitely generated $\ZZ$--algebra}.
We apply Proposition~\ref{lem_gonflement}\eqref{semilocali-b} with $\ell = 2$: 
there exists a semilocal ring  $S$, which is an \'etale extension of $R$ of degree $2^u$ and whose residue fields are all of cardinality $\geq N+1$.
By the first step for the Azumaya $S$--algebra $A_S = A \ot_R S$ we know $\Nrd_{A\ot S/S}(A\ti_S) \subset \rmN_{X\times S}^{\et}(S)$. Lemma~\ref{lemn_norm1}\eqref{lemn_norm1-a} then yields $\rmN_{S/R}\big( \Nrd_{A\ot S/S}(A\ti_S)\big)\subset \rmN_X^{\et}(R)$. In particular for $a\in A\ti$ we get
$a^{2^u} = \rmN_{S/R}(a \ot 1_S) \in \rmN_X^{\et}(R)$. The analogous argument for $\ell = 3$ proves that $a^{3^v} \in \rmN_X^{\et}(R)$ for a suitable $v\in \NN_+$.
Now let $s$, $t\in \ZZ$ such that $1= s \,  2^u+ t \, 3^v$. Then
$ a= (a^{2^u})^s \,   (a^{3^v})^t \in \rmN^{\et}_X(R)$.
\sm

{\em Third step: $R$ general.} We apply noetherian reduction (\cite[7.7.2]{Ford}, \cite[I, 2.9]{KO}): there exists a finitely generated $\ZZ$--algebra $R_0$ and an Azumaya $R_0$--algebra $A_0$ such that $A \cong A_0 \ot_{R_0} R$. Let $\m_1, \ldots , \m_s$ be the maximal ideals of $R$. Replacing $R_0$ by the semilocalization \ref{semilocali} in the prime ideals $\m_i \cap R_0$ of $R_0$, we can assume that already $R_0$ is the semilocalization of a finitely generated $\ZZ$--algebra.

Let us fix $a_0 \in A\ti$. Enlarging $R_0$ if necessary, we can further suppose that $a_0 \in A_0\ti$. By the second step,  there exists a finite \'etale $R_0$--algebra $S_0$ and $s_0\in S_0$ such that $\uSB(A_0)(S_0) \ne \emptyset$ and $\Nrd_{A_0/R_0}(a_0) = \rmN_{S_0/R_0}(s_0)$. The $R$--algebra $S = S_0 \ot_{R_0} R$ is finite \'etale. By \ref{SBs-lem}\eqref{SBS-2} the Azumaya $S_0$--algebra $A_0\ot_{R_0}S_0$ is neutral, consequently so is $(A_0 \ot_{R_0} S_0) \ot_{S_0} S \cong A_0 \ot_{R_0} S \cong  A \ot_R S$. Since $S$ is in particular a finite $R$--algebra, it is semilocal and  therefore $\uSB(A)(S) \ne \emptyset$ by \ref{SBs-lem}\eqref{SBs-lem-c}. Moreover, by base change of $\Nrd$ and $\rmN$, we have $\Nrd_{A/R}(a_0 \ot 1_R) = \Nrd_{A_0/R_0}(a_0) \ot 1_R = \rmN_{S_0/R_0}(s_0) \ot 1_R = \rmN_{S/R}(s_0 \ot 1_R)$. This implies $\Nrd_{A/R}(A\ti) \subset \rmN^{\et}_X(R)$. \end{proof}

\ms

In the remainder of this section we will describe the norm groups and \'etale norm groups of quadrics $\uQ(q)$ for nonsingular $q$, recalled in \ref{hrq}. We will deal with quadratic forms of small ranks in \ref{lem_norm_rk2} and \ref{lem_norm_rk3}. In the rank-$2$-case we will use the following Lemma~\ref{PR} from \cite{PRbook}, characterizing split quadratic \'etale $R$--algebras. 
We include an ad-hoc proof of \ref{PR} for the sake of completeness. The equivalence \ref{PR}\ref{PR-i} $\iff$ \ref{PR}\ref{PR-ii} is \cite[V, (2.2.4)]{K}.

\comments{Is \ref{PR} or \ref{spq} known? The main point is \ref{PR-iii} $\implies$ \ref{PR-i}. Hyperbolic in \ref{PR-iii} does not mean split hyperbolic. Is there another reference than \cite{PRbook}? Observe that $R$ is not assumed to be semilocal.\sm

Well-known, but we do not know a published proof}

\comments{(2023-11-04) A similar characterization as \ref{PR} and more results are given in Scharlau's book, II, 11.8--11.15}

\begin{lem}[{\cite{PRbook}}]\label{PR} Let $C$ be a quadratic \'etale $R$--algebra. We denote by $n_C = \rmN_{C/R}$ its norm. Then the following are equivalent:
 \begin{enumerate}[label={\rm (\roman*)}]
  \item \label{PR-i} $C$ is split,

  \item\label{PR-ii} $n_C$ contains an isotropic vector,

  \item\label{PR-iii} $n_C$ is hyperbolic.
\end{enumerate}
In this case $n_C(C\ti) = R\ti$.
\end{lem}

\begin{proof} The implication \ref{PR-i} $\implies$ \ref{PR-ii} is clear and \ref{PR-ii}  $\implies$ \ref{PR-iii} is \eqref{isotrop-d1}. Let us assume that $n_C$ is hyperbolic. Thus, by \ref{carhyp}, 
there exist line bundles $L_1 $ and $L_2$ such that $C=L_1 \oplus L_2$ and $n_C(L_1) = 0 = n_C(L_2)$. We write $1_C = e_1 + e_2$ with $e_i\in L_i$, $i=1,2$. Then $n_C(e_i) = 0$ and $1 = n_C(1_C) = b_{n_C}(e_1, e_2)$ holds. Since any $c\in C$ satisfies the quadratic identity $c^2 - b_{n_C}(1_C, c) c + n_C(c)1_C = 0$, the elements $e_i$ are unimodular idempotents satisfying $e_1 e_2 = e_1 (1_C - e_1) = 0$, i.e. they are elementary idempotents as defined in \ref{eica}, see also \ref{elid}. 
Hence $C = Ce_1 \times Ce_2$ and $C_i = Re_i$ by rank considerations, proving \ref{PR-i}. In this case $n_C(C\ti) = R\ti$ is clear.
\end{proof}

We remind the reader of the notation $\rmD(q)$ and $\rmD(q)^{\rm [ev]}$ introduced in \ref{thm-kneb}.

\comments{(2021-06-07) Deleted the assumption that $R$ is semilocal. I changed the statement -- before it was claimed that $\rmN_X^{\et}(R)=  \rmN_X(R) = \rmN_{\calD/R}(\calD^\times)= \rmD^{\rm [ev]}(q)$ without the assumption $1\in \rmD(q)$. }

\begin{lem}[Rank $2$] \label{lem_norm_rk2}
Let  $(M,q)$ be a quadratic $R$--space of constant rank $2$, let $X= \uQ(q)$ be the associated quadric, and let $\calD= \Dis(q)$ be the discriminant algebra of $(M,q)$, a quadratic \'etale $R$--algebra by {\rm \ref{qfdi}\eqref{discralg-even}}; we let $\rmN_{\calD/R}$ be its norm, and assume $\rmD(q) \ne \emptyset$. \sm

\begin{enumerate} [label={\rm (\alph*)}]
 \item \label{lem_norm_rk2.1} Then $\rmN_X^{\et}(R)=  \rmN_X(R) = \rmN_{\calD/R}(\calD^\times)$.\sm

 \item \label{lem_norm_rk2.2} If\/ $1\in \rmD(q)$, then\/
 $\rmN_X^{\et}(R)=  \rmN_X(R) = \rmN_{\calD/R}(\calD^\times)= \rmD^{\rm [ev]}(q)$.
\end{enumerate}
\end{lem}

\begin{proof}
In the setting of the lemma, the discriminant algebra equals the even Clifford algebra $\Cli_0(q)$. The assumption $\rmD(q) \ne \emptyset$ implies that $q \cong u\rmN_{\calD/R} $ for some $u\in R\ti$, \cite[V, (2.2)]{K}.
\sm

\ref{lem_norm_rk2.1} We will prove $ \rmN_X(R) \subset \rmN_{\calD/R}(\calD\ti)\subset \rmN^{\et}_X(R)$. To show the first inclusion, let $S$ be a faithfully projective $R$--algebra such that  $X(S) \ne \emptyset$. Then $q_S$ is hyperbolic by \ref{carhyp}. Since this then also holds for $\rmN_{\calD\ot S/S}$,  Lemma~\ref{PR} shows that $\Dis(q_S) = \calD\ot_R S$ is split as $S$--algebra and $\rmN_{\calD\ot S/S}((\calD \ot S)\ti) = S\ti$. By transitivity of norms \ref{trno}\eqref{trno-b}, we then obtain
\begin{align*}
 \rmN_{S/R}(S\ti) &= \rmN_{S/R}\big( \rmN_{\calD \ot S/S}( (\calD \ot_R S)\ti)\big) \\ &= \rmN_{\calD/R}\big( \rmN_{\calD \ot S/\calD}( (\calD \ot_R S)\ti)\big)
  \subset \rmN_{\calD/R}(\calD\ti)
 \end{align*}
which implies $\rmN_X(R) \subset \rmN_{\calD/R}(\calD\ti)$.

Since $\calD \ot_R \calD$ is a split quadratic \'etale $\calD$--algebra \cite[III, (4.1)]{K}, its norm $\rmN_{\calD\ot_R \calD/\calD}$ is split hyperbolic, which forces $q_{\calD} = u \rmN_{\calD\ot \calD/\calD}$ to be hyperbolic too. Hence $X(\calD) \ne \emptyset$, and therefore $\rmN_{\calD/R}(\calD\ti) \subset \rmN_X^{\rm \et}(R)$. \sm

\ref{lem_norm_rk2.2} The assumption $1\in \rmD(q)$ implies $q \cong \rmN_{\calD/R}$. Therefore $\rmD(q) = \rmD(\rmN_{\calD/R}) = \rmN_{\calD/R}(\calD\ti)$. Since the latter is a subgroup of $R\ti$, it equals the subgroup generated by $\rmD(q)$, which is $\rmD(q)^{\rm[ev]}$ by \ref{dqd-ele}\eqref{qdq-ele-aa}. \end{proof}
\ms

In Lemma~\ref{lem_norm_rk3} we will prove an analogous result for quadratic spaces of rank $3$, replacing the quadratic \'etale $R$--algebra $\calD$ by a quaternion algebra, i.e., an Azumaya $R$--algebra of rank $4$.  The next two results investigate quaternion algebras; the first of them characterizes split quaternion algebras in analogy to Lemma~\ref{PR}.

\begin{lem}[Rank $3$]\label{lem_norm_rk3}
Let $R$ be a semilocal ring, and let $(M,q)$ be a quadratic space of rank $3$. We put $X= \uQ(q)$ and $Q = \Cli_0(M,q)$, a quaternion $R$--algebra. Then the norm groups satisfy $\rmN^{\et}_X(R) = \rmN_X(R) = \Nrd(Q\ti) = \rmD(q)^{\rm [ev]}$.
\end{lem}

\begin{proof} By \eqref{exqu1} (= \eqref{sqpc0}) we can replace $X$ by the Severi-Brauer scheme of $A$. By Proposition~\ref{prop_norm_SB} we then get $\rmN^{\et}_X(R) = \rmN_X(R) = \Nrd(Q\ti)$. The final equation
 $\Nrd(Q\ti) = \rmD(q)^{\rm [ev]}$ is Corollary~\ref{craco}.\end{proof}


%
%
\sm

\begin{thm}\label{prop_norm_quad} Let $R$ be a semilocal ring, let $(M,q)$ be a quadratic $R$--space of rank $\geq 2$ and let $X=\uQ(q)$ be the quadric associated with $q$. Then $\rmN_X^{\et}(R) =  \rmD(q)^{[\rm ev]}$.
\end{thm}

\begin{proof} We start by performing several reductions: $R$ is connected, $M$ has constant rank $\ge 3$, $\uS_{q,1}^{\rm sm}(R) \ne \emptyset$ and $ \rmD(q)^{[\rm ev]}= \lan \rmD(q) \ran$.

Indeed, both $\rmN_X^{\et}(R)$ and $\rmD(q)^{[\rm ev]}$ respect direct products, \ref{nogr}\eqref{nogr-nor} and \ref{dqd-ele}\eqref{dqd-ele-e}. We can therefore assume that $R$ is connected. Then $M$ has constant rank. The rank $2$ case has been handled in Lemma~\ref{lem_norm_rk2}, so that we can assume that $M$ has rank $\geq 3$. (By Lemma~\ref{lem_norm_rk3}  we could even assume $\rank M \ge 4$, but this would not lead to a simpler proof.) Recall from \ref{dqd-ele}\eqref{dqf-ele-scal}  that $\rmD(uq)^{[\rm ev]} = \rmD(q)^{[\rm ev]}$ for any $u\in R\ti$. Since the quadrics associated with $q$ and $uq$ coincide, we can scale $q$ and thus assume $\uS_{q,1}^{\rm sm}(R) \ne \emptyset$, cf.\ \ref{lem_smooth_locus_exam}\eqref{lem_smooth_locus_exam-c}. Then $1_R\in \rmD(q)$ and $\rmD(q)^{[\rm ev]}= \lan \rmD(q) \ran$ is the subgroup of $R\ti$ generated by $\rmD(q)$, see \ref{dqd-ele}\eqref{qdq-ele-aa}. \sm

Proof of $\rmN_X^{\et}(R) \subset \rmD(q)^{[\rm ev]}$: Let $a\in \rmN_X^{\et}(R)$. By \ref{nogr}\eqref{nogr-nor} we can assume that there exists a finite \'etale  $R$--algebra $S$ of positive rank such that $q_S$  is isotropic and $a\in \rmN_{S/R}(S\ti)$. By \ref{isotrop}\eqref{isotrop-d}, $q_S$ contains a hyperbolic plane, implying  $\rmD(q_S) = S\ti$. Since $R$ is connected, $S$ has constant rank. Therefore the Knebusch Norm Principle~\eqref{thm-kneb-b1} applies  and shows
$a  \in \rmN_{S/R}(S^\times) = \rmN_{S/R}\big( \rmD(q_S)\big) \subset \rmD(q)^{[\rm ev]}$.%
\sm

Proof of $ \rmD(q)^{[\rm ev]}\subset\rmN_X^{\et}(R)$: Since  $ \rmD(q)^{[\rm ev]}= \lan \rmD(q) \ran$ it suffices to prove $\rmD(q) \subset \rmN_X^{\et}(R)$.
We fix $a \in \rmD(q)$, write $a=q(m)$ and let $v\in \uS_{q,1}\rmsm(R)$. We will proceed in $4$ steps in increasing generality for  $m\in \uS_{q,a}(R)$. \sm


\begin{inparaenum}[(I)]\item\label{prop_norm_quadI} 
{\em  $m \in \uS^{sm}_{q,a}(R)$ and the residue fields of $R$ are all infinite.}
In this case, Proposition \ref{plar} provides $g \in \SO(q)(R)$ such that
$P = Ru \oplus R(g\cdot  v)$ is a regular plane; furthermore, it provides $u' \in P^\perp$ such that $W= P \perp Ru'$ is a direct summand of $M$ and $(W,  q|_W)$ is a quadratic space. We denote by $Y$ the projective $R$-quadric attached to $q_{\mid W}$.
Then $W$ represents $1$ and the rank 3 case of the statement (i.e.\ Lemma \ref{lem_norm_rk3})  yields the inclusion $\rmD^{[ev]}(q_{\mid W})  \subseteq  \rmN_Y^{\et} (R)$. Since $a =q(m) \, q(v) \in \rmD^{[ev]}(q_{\mid W}) $, we obtain that $q(m)  \in \rmN_Y^{\et} (R)$. \sm 


\item \label{prop_norm_quadII} {\em $m \in \uS^{sm}_{q,a}(R)$
and $R$ is the semilocalization of a finitely generated $\ZZ$-–algebra.}
According to Proposition \ref{lem_gonflement}, applied with an odd prime $\ell$,
there exists a tower
\[
R=R_0 \subset R_1 \subset R_2 \subset \cdots
\]
of semilocal rings such that for every $j \ge 0$ the following holds:
\end{inparaenum}

 \begin{enumerate}[label={\rm (\roman*)}]
  \item  $R_{j+1}$ is finite \'etale
  of degree $\ell$ over $R_j$;
  
  \item $R_\infty = \limind R_i$ is a semilocal ring all of whose residue fields are infinite.

\end{enumerate}
By step~\eqref{prop_norm_quadI}, we have $a \in \rmN_X^{\et}(R_\infty)$.
It follows that there exists $i \geq 0$ such that  $a \in \rmN_X^{\et}(R_i)$
(Lemma \ref{lemn_norm1}\eqref{lemn_norm1-b}). Since $R_i$ is finite \'etale of degree $\ell^i$, we have $\rmN_{R_i/R}\bigl( \rmN_X^{\et}(R_i)\bigr) \subset \rmN_X^{\et}(R)$.
It follows that $a^{\ell^i} \in  \rmN_X^{\et}(R)$.
If $\rmN_X^{\et}(R)$ contains $(R^\times)^2$, we can conclude that  $a \in \rmN_X^{\et}(R)$. But, indeed, $(R^\times)^2 \subset\rmN_X^{\et}(R)$: by Lemma~\ref{nqf-LG} we can write $(M,q)=(M_2,q_2) \perp (M',q')$ where $(M_2,q_2)$ is regular plane. Let $Z$ be the quadric associated with $(M_2, q_2)$. Then Lemma~\ref{lem_norm_rk2}\ref{lem_norm_rk2.1} shows $1\in \rmN_Z^{\et}(R)$ and therefore also $(R\ti)^2 \in \rmN_Z^{\et}(R)$. Finally, $\rmN_Z^{\et}(R)\subset \rmN_X^{\et}(R)$ in view of \ref{nogr}\eqref{nogr-f} for the closed embedding $Z \hookrightarrow X$.%
\sm

\begin{inparaenum}[(I)] \setcounter{enumi}{2}
\item\label{prop_norm_quadIII} {\em $m \in \uS^{sm}_{q,a}(R)$.} By  Example~\ref{noethred}\eqref{noethred-b}, $R$ is the directed limit of semilocalizations of finitely generated $\ZZ$--algebras, say $R = \limind_{\la \in \La} R_\la$, and by Corollary~\ref{dirquad}\eqref{dirquad-b} we can assume that there exists a directed system $(M_\la, q_\la)_{\la \in \La}$ of $R_\la$--quadratic spaces whose limit is $(M,q)$. Moreover, by  \ref{dirquad}\eqref{dirquad-c}, for some $\mu \in \La$ we know that there exist $a_\mu \in R_\mu \ti$ and $m_\mu \in \uS_{q_\mu, a_\mu}\rmsm(R_\mu)$ such that $m$ is the canonical image of $m_\mu$. Then, by step~\eqref{prop_norm_quadII},   
 $a_\mu  \in \rmN_X^{\et}(R_{\mu})$ and hence a fortiori $a \in \rmN_X^{\et}(R)$.
\sm 

\item\label{prop_norm_quadIV} $m\in \uS_{q,a}(R)$.  According to Proposition~\ref{embp}, there exists a complemented nonsingular $R$-submodule  $W \subset M$ free of rank $3$ which contains $m$ and an element $w$ such that $q(w) \in R^\times$ and
$w \in \uS^{sm}_{q_{\mid W},q(w)}$. We denote by $Y$ the projective $R$-quadric attached to $q_{\mid W}$. Using the rank $3$ case of the statement (Lemma \ref{lem_norm_rk3}), we have
$\rmN^{\et}_Y(R) = \rmD(q_{\mid W}^{\rm [ev]}$. Since  $q(m) \, q(w) \in \rmD(q_{\mid W}) ^{[\rm ev]}$, we obtain $q(m) \, q(w) \in \rmN_Y^{\et} (R)$.
On the other hand, we have the  closed embedding $Y  \hookrightarrow X$,  so that 
$\rmN_Y^{\et} (R) \subset \rmN^{\et}_X(R)$ by \ref{nogr}\eqref{nogr-f} again.
Thus $q(m) \, q(w) \in \rmN_X^{\et} (R)$. Finally, step \eqref{prop_norm_quadIII} tells us that $q(w) \in \rmN_X^{\et} (R)$. By then also $q(m)  \in \rmN_X^{\et} (R)$.
\end{inparaenum}\end{proof}

\subsection{Remarks on Theorem~\ref{prop_norm_quad}}\label{prop_norm_quadrem}
\begin{inparaenum}[(a)] \item The result is not true for $(M,q)$ of rank $1$. For instance, for $(M,q) = (R, \lan 1 \ran_q)$ we have $X= \emptyset$, so $\rmN_X^{\et} = \{1\}= \rmN_X(R)$, while $\rmD(q)^{[\rm ev]} = \rmD(q) = R\ti$.
\sm

\item Our proof 
is inspired by the proof of \cite[Lem.~2.2]{CSk}, which establishes the equality $\rmN_X(R) = \rmD(q)^{[\rm ev]}$ in case $R$ is a field of characteristic $\ne 2$
(the groups $\rmN_X(R)$ considered here and in \cite{CSk} coincide by Lemma~\ref{lem_norm_field}). As mentioned in \cite{CSk}, the equality $\rmN_X(R) = \rmD(q)^{[\rm ev]}$ is due to Rost (unpublished); it is discussed in \eqref{prop_norm_quadrem-c} below.%
\sm

\item \label{prop_norm_quadrem-c} If $(M,q)$ has rank $2$ or $3$,  we have shown in Lemmata~\ref{lem_norm_rk2} and \ref{lem_norm_rk3} that  $ \rmN^{\et}_X(R)=\rmD(q)^{[\rm ev]}= \rmN_X(R)$. In general, \ref{prop_norm_quad} says that $\rmN^{\et}_X(R)=\rmD(q)^{[\rm ev]} \subset \rmN_X(R)$. The proof of \ref{prop_norm_quad} shows that $\rmD(q)^{[\rm ev]} = \rmN_X(R)$ holds as soon as  the Knebusch Norm Principle is true for an arbitrary faithfully projective $R$--algebra  $S$ such that $q_S$ is isotropic. This is the case if $R$ is
    a field \cite[18.10]{EKM}.
\end{inparaenum}

\pcomments{(2021-04-17) The counterpart of $\mathrm{Ref}(q)^+ = \SO(q)$ on the quaternionic side is that $Q^\times$ is generated by $R^\times$ and even products of
invertible pure elements.\ms

(E, 2021-06-25) Should we include this somewhere? }

\newpage
\appendix

\section{Direct limits, Noetherian reduction} \label{sec:dirlim}

\comments{(2022-01) This is an appendix written for the Octonion book. I put it here with minor deletions, i.e., everything related to composition algebras.}

With the exception of the Examples~\ref{noethred},
throughout this appendix we assume that $R=\limind R_\la$ for a given directed system $(R_\la)_{\la \in \La}$ of $R_0$--rings, as reviewed in \ref{dirlim-rev}.

\subsection{Directed systems and direct limits.} \label{dirlim-rev} In this subsection we review directed (= inductive) systems and their limits,  here called direct limits but elsewhere sometimes referred to as inductive limits or colimits. 
We refrain from expounding a general theory of directed systems in categories, but  rather restrict ourselves to the situation needed here and elsewhere in the book. 
\sm

Let $(\La, \le)$ be a directed poset, i.e., $\La$ is a non-empty set with a binary relation $\le$, which is reflexive, transitive, and has the property that for $\la_i \in \La$, $i=1,2$,  there exists $\mu \in \La$ satisfying $\la_i \le \mu$.

Throughout we fix $R_0\in \ZZalg$, and a directed system $(R_\la, \vphi_{\mu\la})_{\la \in \La}$ be a directed system in $\Rnalg$.  Thus, by definition, every $R_\la$ is an $R_0$--algebra and the maps $\vphi_{\mu\la} \co R_\la \to R_\mu$, which are assumed to exist for $\la \le \mu$, are unital $R_0$--algebra homomorphisms satisfying $\vphi_{\la\la} = \Id_{R_\la}$ and $\vphi_{\nu\mu} \circ \vphi_{\mu\la} = \vphi_{\nu\la}$ whenever $\la \le \mu \le \nu$. Associated with such a system are  an $R_0$--algebra $R= \limind R_\la$, the {\em direct limit of $(R_\la)$}, and $R_0$--algebra homomorphisms $\vphi_\la \co R_\la \to R$ satisfying $\vphi_\mu \circ \vphi_{\mu \la} = \vphi_\la$ for $\la \le \mu$ and $R = \bigcup_{\la \in \La} \vphi_\la(R_\la)$. The $R_0$--algebra $R$ enjoys a universal property, analogous to the one recalled below for modules.
Examples of direct limits are given in \ref{noethred}.

To facilitate quoting \cite{EGA}, we note that the data $R_0$, $(R_\la, \vphi_{\mu \la})_{\la \in \La}$ and $R=\limind R_\la$ give rise to  an affine base scheme $S_0 = \Spec(R_0)$, an inverse (= projective) system $(S_\la, u_{\mu \la})_{\la\in \La}= (\Spec(R_\la), \Spec(\vphi_{\mu \la}))_{\la\in \La}$ of $S_0$--schemes, and an $S_0$--scheme $S= \Spec(R) = \limproj S_\la$, the limit of $(S_\la, u_{\mu \la})$. By Yoneda, the converse is also true, i.e., an affine base scheme $S_0$, an inverse system $(S_\la, u_{\mu\la})$ of affine $S_0$--schemes and its limit $S=\limproj S_\la$ arise from unique data $R_0$, $(R_\la, \vphi_{\mu \la})$ and $R$ as above.
\sm

A {\em directed system $(\scP_\la, f_{\mu\la})$ of modules\/}  associated with the directed system $(R_\la)$ consists of a family of right $R_\la$--modules $\scP_\la$ and dimorphisms $f_{\mu\la}$, $\mu\ge \la$, i.e., additive maps $f_{\mu \la} \co \scP_\la \to \scP_\mu$ satisfying $f_{\mu\la} (p_\la r_\la) = f_{\mu\la}(p_\la) \vphi_{\mu\la}(r_\la)$ for $p_\la \in \scP_\la$ and $r_\la \in R_\la$. Associated with such a system are  an $R$--module $\limind \scP_\la$, the {\em direct limit of $(\scP_\la)$\/}, and dimorphisms $f_\mu \co \scP_\mu \to \limind \scP_\la$ for which $f_\mu \circ f_{\mu \la} = f_\la$ whenever $\la \le \mu$ and $\limind \scP_\la = \bigcup_{\la \in \La} f_\la(\scP_\la)$ holds. The $R$--module $\limind \scP_\la$ has the following universal property: given an $R$--module $\scQ$ and dimorphisms $u_\la \co \scP_\la \to \scQ$ satisfying $u_\mu \circ f_{\mu\la} = u_\la$ for $\mu \ge \la$, there exists a unique $R$--module map $u \co \limind \scP_\la \to \scQ$ such that $u_\la = u \circ f_\la$ holds for all $\la \in \La$.
\begin{equation} \label{dirlim-rev0} \vcenter{
\xymatrix{
     \scP_\la \ar[dd]_{f_{\mu\la}} \ar[dr]^{f_\la} \ar@/^/[drrr]^{u_\la}\\
      &\limind \scP_\la \ar@{-->}[rr]^{\exists ! \, u} && \scQ \\
     \scP_\mu  \ar[ur]_{f_\mu} \ar@/_/[urrr]_{u_\mu}
}}\end{equation}
\lv{
In this case $u$ is injective if for all $\la$ and all $x,y\in \scP_\la$ the equality $u_\la(x) = u_\la (y)$ implies $f_{\mu\la}(x) = f_{\mu\la}(y)$ for some $\mu \ge \la$, and $u$ is surjective if $\scQ = \bigcup_\la \, u_\la(\scP_\la)$ \cite[8.2.1]{TY}. }

\subsection{Some standard facts} \label{dirlimfa} We keep the notation of \ref{dirlim-rev}: $(\La, \le)$ is a directed poset and $(R_\la, \vphi_{\mu\la})$ is a directed system in $\Rnalg$ with $R= \limind R_\la$.%
\sm

\begin{inparaenum}[(a)]
 \item \label{dirlimfa-a} ({\em Reduction of modules}) An $R$--module $M$ is finitely presented if and only if there exist $\la \in \La$ and a finitely presented $R_\la$--module $M_\la$ such that $M \cong M_\la \ot_{R_\la} R$ (\cite[IV$_2$, 5.13.7.1]{EGA}, \cite[Tag 05N7(1)]{St}. 
\sm

\item \label{dirlimfa-hom} ({\em Reduction of homomorphisms}) Let $\vphi \co M \to N$ be an $R$--linear map between finitely presented $R$--modules. By \eqref{dirlimfa-a} we can assume that there exist $\la \in \La$ and finitely presented $R_\la$--modules $M_\la$ and $N_\la$ such that $M=M_\la \ot_{R_\la} R$ and $N=N_\la \ot_{R_\la} R$. Then there exist $\mu \in \La$, $\mu \ge \la$, and an $R_\mu$--module map $\vphi_\mu \co M_\la \ot_{R_\la} R_\mu \to N_\la \ot_{R_\la} R_u$ such that $\vphi = \vphi_\mu \ot 1_R$,  
\cite[Tag 05N7(2)]{St}. 
\sm 

\item \label{dirlimfa-sur} ({\em Reduction of epi- and isomorphisms}) In the setting of \eqref{dirlimfa-hom}, suppose that $\vphi \co M \to N$ is a surjective $R$--linear map or an isomorphism. For $\mu' \in \La$, $\mu' \ge \mu$, we put
    $ M_{\mu'} = (M_\la \ot_{R_\la} R_\mu) \ot_{R_\mu} R_{\mu'} = M_\la \ot_{R_\la} R_{\mu'}$ and employ the analogous notation for $N$. Then there exists $\nu \in \La$, $\nu \ge \mu$ such that 
    \[ \vphi_\nu = \vphi\ot 1_{R_\nu} \co M_\nu \to N_\nu \]   
    is surjective or, respectively, an isomorphism. 
    
    Indeed, let $x_i$, $i=1, \ldots, m$ be generators of $M_\mu$. We can pick $\nu \ge \mu$ such that $N_\nu$ contains all $\vphi_\mu(x_i \ot 1_{R_\nu})$. Then $\vphi_\mu$ is surjective. In case $\vphi$ is an isomorphism, one applies \eqref{dirlimfa-hom} for $\vphi$ and $\vphi\me$ and uses \cite[Tag 05N7(3)]{St} to conclude. 
\end{inparaenum}

\lv{
(2026-06-02) OLD version: 
\begin{inparaenum}[(a)] \item 
Let $\scP_0$ be an $R_0$--module. Then
\begin{equation} \label{dirlim-rev1}
(\scP_\la, f_{\mu\la}) = (\scP_0 \ot_{R_0}  R_\la, \, \Id_{\scP_0} \ot \vphi_{\mu\la}).
\end{equation}
is a directed system with $\limind \scP_\la \cong \scP_0 \ot_R R$ \cite[II, \S6.3, Prop.~7]{BA}.\sm

\item 
Let $(\scP_\la, f_{\mu\la})$ be a directed system of modules for $(R_\la, \vphi_{\mu\la})$ and let $\La'\subset \La$ such that $(\La', \le)$ is directed with respect to the induced partial order.  Then $(\scP_{\la'}, f_{\mu' \la'})_{\la' \in \La'}$ is a directed system. By the universal property of $\limind$ there exists a unique $R$--linear map
    \begin{equation}\label{dirlim-rev-b1}  u \co \limind {}_{\la' \in \La'} \,  \scP_{\la'} \, \longto \,
    \limind {}_{\la \in \La}\,  \scP_\la.
    \end{equation}
  The map $u$ is an isomorphism if $\La' \subset \La$ is {\em cofinal\/}, i.e., if for every $\la \in \La$ there exists $\la' \in \La'$ such that $\la'\ge \la$ \cite[04E7]{Stacks}. 
  For example, for a fixed $\la \in \La$, the subset $\La_{\ge \la} = \{\mu \in \La: \mu \ge \la\}$ is cofinal.
%
\sm

\item
({\em Reduction of homomorphisms}) Let $M_0$ and $N_0$ be $R_0$--modules. Put $M_\la = M_0 \ot_{R_0} R_\la$, $M = M_0 \ot_{R_0} R$, $N_\la = N_0 \ot_{R_0} R_\la$ and $N = N_0 \ot_{R_0} R$.
Hence $M \cong \limind M_\la$ and $N \cong \limind N_\la$ by \eqref{dirlim-rev-a}.
Since $M_\mu = M_\la \ot_{R_\la} R_\mu$ for $\mu \ge \la$ and similarly for $N_\mu$, we get dimorphisms
\[ h_{\mu\la} \co \Hom_{R_\la} (M_\la, N_\la) \to \Hom_{R_\mu}(M_\mu , N_\mu), \quad u \mapsto u \ot \Id_{R_\mu},
\]
which define a directed system $(\Hom_{R_\la}(M_\la, N_\la), h_{\mu\la})$. Moreover, we have canonical maps
\[ h_\la \co \Hom_{R_\la}(M_\la, N_\la) \longto \Hom_R(M, N), \]
which form a directed system of maps and, by the universal property \eqref{dirlim-rev0}, give rise to a canonical $R$--linear map
\begin{equation}  \label{dirlim-rev-c1}
 \limind \Hom_{R_\la}(M_\la, N_\la) \longto \Hom_R(M,N).
\end{equation}
In particular: \sm

\begin{inparaenum}[(i)]
\item
 If $M_0$ is a finitely presented $R_0$--module, then \eqref{dirlim-rev-c1}  is an isomorphism \cite[IV$_3$, (8.5.2.1)]{EGA}.
    Hence, given an $R$--linear map $f\co M \to N$ there exist $\la\in \La$ and an $R_\la$--linear map $f_\la \co M_\la \to N_\la$ such that $f = f_\la \ot \Id_R$.
\sm

\item 
If both $M_0$ and $N_0$ are finitely presented $R_0$--modules and $f\co M \to N$ is an isomorphism of $R$--modules, there exist $\la \in \La$ and an isomorphism $f_\la \co M_\la \to N_\la$ of $R_\la$--modules such that $f = f_\la \ot \Id_R$ \cite[IV$_3$, (8.5.2.5)]{EGA} (or \cite[05LI]{Stacks}).
\end{inparaenum} \end{inparaenum}
}

\begin{lem}[Reduction of module properties]\label{dirlimfa-b} We continue with the setting of {\rm \ref{dirlim-rev}} and {\rm \ref{dirlimfa}}:  $(R_\la, \vphi_{\mu\la})_{\la \in \La}$ is a directed system in $\Rnalg$ with $R= \limind R_\la$. Let $\PP$ be one of the following properties of a module:

    \begin{inparaenum}[\rm (i)]
     \quad \item \label{dirlimfa-biv} flat,

      \quad \item \label{dirlimfa-bi} finite projective,

      \quad \item \label{dirlimfa-bii} projective of rank $r$,

      \quad \item \label{dirlimfa-biii} faithfully projective.
\end{inparaenum}

\noindent Let $M$ be a finitely presented $R$--module. Then $M$ has property $\PP$ if and only if there exist $\la \in \La$ and an $R_\la$--module $M_\la$ which has property $\PP$ and satisfies $M \cong M_\la \ot_{R_\la} R$.
\end{lem}

\begin{proof} Since the properties $\PP$ are stable under base change, we only have to prove the existence  of $M_\la$ with property $\PP$. By \ref{dirlimfa}\eqref{dirlimfa-a}, there exist $\la \in \La$ and a finitely presented $R_\la$--module $M_\la$ such that $M = M_\la \ot_{R_\la} R$. We thus get a directed system $(M_\mu = M_\la \ot_{R_\la} R_\mu)_{\mu \ge \la}$ of finitely presented modules. Then  \cite[Tag 02JO]{St} with $S=R$ says that there exists $\mu \ge \la$ such that $M_\mu$ is flat = finite projective, proving \eqref{dirlimfa-biv} and \eqref{dirlimfa-bi}. A direct proof of \eqref{dirlimfa-bi} and \eqref{dirlimfa-bii} is given in \cite[IV$_3$, 8.5.5]{EGA}. 
Finally, \eqref{dirlimfa-biii} follows by applying the rank decomposition of modules. 
\end{proof}
\sm

We also need to understand the behaviour of bilinear maps. 

\comments{(2020-10) \ref{dirlem-bi} is proven in \cite[Prop.~7.7.2]{Ford} for $M=N=P$, i.e. algebras, assuming finite projective. The proof in \cite{Ford} has a mistake in the diagrams (7.9) and (7.10). A proof of the general case along the lines of \cite{Ford} is in the file 'dirlim-2020-10-16'. The proof below is much shorter and more general. }

\begin{lem}[Bilinear maps] \label{dirlem-bi} 
As before, $R= \limind R_\la$.   Let $M$, $N$ and $P$ be finitely presented $R$--modules and let $b \co M \times N \to P$ be an $R$--bilinear map.
    Then there exist $\la \in \La$, finitely presented $R_\la$--modules $M_\la$, $N_\la$ and $P_\la$ and an $R_\la$--bilinear map $b_\la \co M_\la \times N_\la \to P_\la$ such that $M = M_\la \ot_{R_\la} R$, $N= N_\la \ot_{R_\la} R$, $P = P_\la \ot_{R_\la} R$ and $b =  (b_\la)_R$ are obtained by base change. Special cases of interest are the following:

    \begin{enumerate}[label={\rm (\roman*)}]
  \item \label{dirlem-ci} If $M=N$, we can assume $M_\la= N_\la$;

  \item \label{dirlem-cii} if $M=N=P$, we can assume $M_\la = N_\la = P_\la$, in other words, an $R$--algebra $(M,b)$ is obtained by base change from an $R_\la$--algebra for suitable $\la\in \La$;

  \item \label{dirlem-ciii} if $P=R$, we can assume $P_\la = R_\la$;

  \item \label{dirlem-civ} if $M$ is finite projective or projective of rank $r$, we can assume that so is $M_\la$; analogously for $N$ and $P$;

  \item \label{dirlem-cv} if the $R$--linear map $M \ot_R N \to P$, $m\ot n \mapsto b(m,n)$ is bijective, then we can choose $\la$, $M_\la$, $N_\la$, $P_\la$ and $b_\la$ such that 
      \[ M_\la \ot_{R_\la} N_\la \to P_\la, \quad m_\la \ot n_\la \mapsto b_\la(m_\la, n_\la),\] 
      is bijective.
  \end{enumerate}
\end{lem}

\begin{proof} By \ref{dirlimfa}\eqref{dirlimfa-a} there exist $\la \in \La$ and finitely presented $R_\la$--modules $M_\la$, $N_\la$ and $P_\la$ such that $M= M_\la\ot_{R_\la} R$, $N = N_\la \ot_{R_\la}R $ and $P = P_\la \ot_{R_\la} R$. 
We can then view $b$ as an $R$--linear map
\[ b \co M\ot_R N = (M_\la  \ot_{R_\la } N_\la )\ot_{R_\la } R \longto P_\la \ot_{R_\la} R.\]
Since $M_\la  \ot_{R_\la } N_\la $ is a finitely presented $R_\la $--module, \ref{dirlimfa}\eqref{dirlimfa-hom} says that there exists $\mu \in \La$, $\mu \ge \la$ and an $R_\mu$--linear map
\[ b_\mu  \co M_\la  \ot_{R_\la} N_\la \ot_{R_\la } R_\mu \longto P_\la \ot_{R_\la } R_\mu \]
such that $b= b_\mu \ot_{R_\mu } \Id_R$. Taking into account 
\[ 
M_\la \ot_{R_\la } N_\la \ot_{R_\la } R_\mu \cong (M_\la \ot_{R_\la } R_\mu  ) \ot_{R_\mu } (N_\la  \ot_{R_\la } R_\mu),
\]
this proves the general claim. The special cases \ref{dirlem-ci}--\ref{dirlem-ciii} hold by the construction above, \ref{dirlem-civ} is a consequence of  \ref{dirlimfa-b}\eqref{dirlimfa-bii}, and \ref{dirlem-cv} follows from
\ref{dirlimfa}\eqref{dirlimfa-sur}. \end{proof}

\begin{prop}[Flat cover] \label{flatcov} Let $S\in \Ralg$ be a flat cover of $R=\limind R_\la$. Then there exists $\la \in \La$ and a flat cover $S_\la$ of $R_\la$ such that $S= S_\la \ot_{R_\la} R$.
\end{prop}

\begin{proof}
Recall that $S\in \Ralg$ is a flat cover of $R$ if and only if $S$ is finitely presented as $R$--algebra and faithfully flat as $R$--module.   By \cite[IV$_1$, 1.8.4.2]{EGA} there exists $\la\in \La$ and a finitely presented $R_\la$--algebra such that ${S\cong S_\la \ot_{R_\la} R}$. 
From  \cite[IV$_3$, 11.2.6.1.(ii)]{EGA} (or \cite[Tag 02JO(iii)]{St} with $S=M$) we then get that there exists $\mu \in \La$, $\mu \ge \la$, such that $S_\mu=S_\la \ot_{R_\la} R_\mu$ is a flat $R_\mu$--module (and a finitely presented $R_\mu$--algebra). Moreover, by \cite[IV$_3$, 8.10.5(vi)]{EGA}, there exists $\nu \ge \mu$ such that the morphism $\Spec(S_{\nu} ) \to \Spec(R_{\nu} )$ is surjective. Thus  $S_{\nu}$ is faithfully flat as $R_\nu$--module. Since it is also finitely presented by the choice of $S_\la$, the algebra $S_\nu$ is a flat cover.\end{proof}


\begin{cor}[Bilinear and quadratic modules] \label{dirquad} 
As before, $R=\limind R_\la$. \sm 

\begin{inparaenum}[\rm (a)]
\item\label{dirquad-a} Let $(M,b)$ be a regular bilinear $R$--module in the sense of {\rm \ref{bfLG}}. Then there exist $\la \in \La$ and a regular bilinear module $(M_\la, b_\la)$ such that $(M,b) = (M_\la, b_\la)_R$. \sm

\item\label{dirquad-b}
Let $(M,q)$ be an $R$--quadratic module. Then there exist $\la \in \La$ and an $R_\la$--quadratic module $(M_\la, q_\la)$ such that $(M,q) = (M_\la, q_\la)_R$. If $(M,q)$ is regular (nonsingular respectively), one can choose $\la$ such that $(M_\la, q_\la)$ is regular (nonsingular respectively). \sm 

\item\label{dirquad-c}
Let $(M,q)$ be a quadratic $R$--space, let $m\in \uS_{q,a}\rmsm(R)$ for some $a\in R\ti$, {\rm \ref{lem_smooth_locus-LG}} and let $(M_\la, q_\la)$ be the directed system of quadratic spaces of {\rm \eqref{dirquad-b}}. Then there exist $\mu \in \La$ such that $m$ and $a$ are the canonical images of $m_\mu \in M_\mu$ and $a_\mu \in R_\mu\ti$ and that $m_\mu \in \uS_{q_\mu, a_\mu}\rmsm (R_\mu)$.   
\end{inparaenum}
\end{cor}

\begin{proof} \eqref{dirquad-a} The existence of a bilinear $R_\la$--module $(M_\la, b_\la)$ follows from Lemma~\ref{dirlem-bi}. Regularity means that the adjoint map $M \to M^*$, $m \mapsto b(m, \cdot)$, is invertible. Therefore, the claim follows from \ref{dirlimfa}\eqref{dirlimfa-sur}. 
\sm 

\eqref{dirquad-b} By \cite[I, (1.7)]{Ba} there exists a bilinear form $b \co M \times M \to R$ such that $q(m) = b(m,m)$ holds for all $m\in M$. By \ref{dirlem-bi},  there exist $\la\in \La$, a finite projective $R_\la$--module $M_\la$ and a bilinear form $b_\la \co M_\la \times M_\la \to R_\la$ such that $b=(b_\la)_R$. Define $q_\la \co M_\la \to R_\la$ by $q_\la(m_\la) = b_\la(m_\la, m_\la)$. Then $(M_\la, q_\la)$ is an $R_\la$--quadratic module such that $(M_\la, q_\la)_R = (M,q)$. For $\mu \ge \la$, we let $(M_\mu, q_\mu) = (M_\la, q_\la)_{R_\mu}$; its polar is $R_\mu$--extension of the polar of $q_\la$. 

Let $(M,q)$ be regular, i.e., the polar $b_q$ of $q$ is regular. 
Then  \eqref{dirquad-a} 
says that there exists $\mu\in \La$ for which $b_{q_\mu}$ is regular, i.e., $q_\mu$ is regular.

Finally, let $(M,q)$ be nonsingular. We have seen that there exist $\la\in \La$ and a quadratic $R_\la$--module $(M_\la, q_\la)$ such that $(M,q)= (M_\la, q_\la)_R$. We can decompose  $(M_\la, q_\la)$ and $(M,q)$ with respect to the rank decomposition of $M$ and $M_\la$. Since nonsingularity respects this decomposition, \ref{qf}\eqref{qf-redc}, it is no harm to assume that both $M$ and $M_\la$ have constant rank. This then also holds for the $R_\mu$--module $M_\mu$. By \ref{qf}\eqref{qfnsp}, there exist a unique split quadratic space $(N_\ZZ, h_\ZZ)$ over $\ZZ$ and a flat $R$--cover $S$ such that $(M,q)_S$ is isometric to $(N_\ZZ,h_\ZZ)_S = (N_S, h_S)$. 

By \ref{flatcov} we can assume that there exists a flat cover $S_\la$ of $R_\la$ such that $S\cong (S_\la)_R$. Putting $S_\mu = S_\la \ot_{R_\la} R_\mu$ and $(N_\mu, h_\mu) = (N_\ZZ, h_\ZZ)_{S_\mu}$ and replacing $\La$ by the cofinal subset $\{\mu \in \La: \mu \ge \la \}$ we then get a directed system $(S_\la)$ of flat covers $S_\la$ of $R_\la$ and a directed system $(N_\la, h_\la)$ of split $S_\la$--quadratic spaces whose limits are $S$ and $(N_S, h_S)$ respectively. Since $(M,q)_S \cong (N_S, h_S)$, it follows from the obvious isometry version of \ref{dirlimfa}\eqref{dirlimfa-sur} that for some $\mu$ already $(M_\mu, q_\mu)_{S_\mu}$ and $(N_\mu, h_\mu)$ are isometric. In particular, $(M_\mu, q_\mu)$ is nonsingular. \sm 

\eqref{dirquad-c} Our assumption means that $q(m) = a \in R\ti$ and $b_q(m, \cdot) \co M \to R$ is surjective. It is clear that these data descend to some $\mu \in \La$. 
\end{proof}


\comments{(2026-06-03) The following not needed, proof in  the file 'dirlim-LG-2026-06-02' and earlier: \sm 

We consider the following properties $\PP$ of an algebra, being a
\begin{enumerate}[label={\rm (\roman*)}]
  \item \label{redap-i} unital algebra,
   \item \label{redap-ii} associative algebra,
   \item\label{redap-comm} commutative algebra,
   \item \label{redap-iii}  Azumaya algebra,
   \item \label{redap-iv} Azumaya algebra of degree $n$,
   \item \label{redap-v} finite \'etale cover,
   \item \label{redap-vi} \'etale cover of degree $n$,
   \item \label{redap-vii} quasi-composition algebra,
   \item \label{redap-viii} composition algebra,
   \item \label{redap-ix} quasi-composition or composition algebra of rank $r\in \{1,2,4,8\}$.
  \end{enumerate}
Let $C$ be an $R$--algebra whose underlying $R$--module is finitely presented. Then $C$ has property $\PP$ if and only if there exist $\la \in \La$ and an $R_\la$--algebra $C_\la$ which has property $\PP$, satisfies $C_\la \ot_{R_\la} R \cong C$ as $R$--algebras and whose underlying $R_\la$--module is finitely presented.
}

\subsection{Remark.} 
Let $q_{0,r}$ be the split quadratic form of rank $r$ over $R$, defined in \ref{qf}\eqref{quadfoc}, and let $G = \uO(q_{0,r})$ be its orthogonal $R$--group scheme. By part \ref{qfnsp-iv}  of \ref{qf}\eqref{qfnsp}, the cohomology set $H^1\fppf(R, G)$ classifies nonsingular quadratic spaces over $R$ of constant rank $r$. Hence, the nonsingular part of \ref{dirquad}\eqref{dirquad-b} is a special case of \cite[Prop.~2.8.1]{G2}.

\subsection{Examples of direct limits.} \label{noethred} \begin{inparaenum}[(a)]
  \item\label{noethred-a} Let $R\in \ZZalg$ be arbitrary, and let $\Phi(R)$ be the set of finitely generated unital $\ZZ$--subalgebras of $R$, which is a directed poset with respect to inclusion. It is well-known that
      \begin{equation} \label{noethred-a1}
         R=\textstyle \limind_{\, F\in \Phi(R)} F.
      \end{equation}

      Since every $F\in \Phi(R)$ is a noetherian ring, the reductions in  \ref{dirlimfa-b} and \ref{dirquad} with respect to \eqref{noethred-a1} are usually referred to as {\em noetherian reduction}. \sm 
       
\item \label{noethred-b} We claim: {\em Every semilocal ring $R$ is a direct limit of semilocalizations of finitely generated unital $\ZZ$--subalgebras of $R$.}

  Indeed, let $\m_1, \ldots, \m_c$ be the maximal ideals of $R$. For $F\in \Phi(R)$ the ideals $\p_{Fi} = F \cap \m_i$ are prime ideals of $F$. Let $F'$ be the semilocalization of $F$  in $\p_{F1}, \ldots, \p_{Fc}$, see \ref{semilocali}. If $F_1\in \Phi(R)$ with $F\subset F_1$, the inclusion canonically extends to a ring homomorphism $F' \to F'_1$, giving rise to a directed system $(F')_{F\in \Phi(R)}$ whose direct limit is isomorphic to the semilocalization of $R$ in $\m_1, \ldots \m_c$, thus to $R$.
\end{inparaenum}

\comments{(2026-06-03) No longer used: \sm 

For easier reference we state two special cases of noetherian reduction
which follow from \ref{dirquad} and \ref{redap}.

\textbf{Corollary} \begin{inparaenum}[\rm (a)]  \item \label{noeth-rd-cor-a} Let $\scC$ be a composition $R$--algebra of fixed rank $r$. Then there exists a finitely generated unital $\ZZ$--subalgebra $R_0$ of $R$, thus a noetherian ring, and a composition $R_0$--algebra $\scC_0$ of rank $r$ such that $\scC_0 \ot_{R_0} R \cong \scC$ as $R$--algebras.
\sm

\item \label{noeth-rd-cor-b} Let $(M,q)$ be a quadratic space over a semilocal ring $R$. Then there exist  $R' \in \ZZalg$, which is the semilocalization of a finitely generated unital $\ZZ$--subalgebra of $R$, and a quadratic space $(M',q')$ over $R'$ such that $(M,q) \cong (M',q')_R$.
\end{inparaenum} 
} 
\comments{(2023-07-26) Noetherian reduction also holds for Galois extensions and crossed product algebras, see \cite[Prop.~6.4]{Salt}. }

\newpage


\section{Some concepts and results from algebraic geometry}\label{sec:ag}

\comments{This section has to be re-arranged in the end to make it coherent. For now, I put here everything in no particular logical order.}

We use concepts and results from algebraic geometry which can be found in \cite{EGA,GW,St}. In this appendix we review some of them, explain our notation and present some results that we could not find in the quoted references. 
\sm

\subsection{Some notation, terminology and facts}\label{ag} \begin{inparaenum}[(a)]
\item\label{ag-bas} As in the main body of the paper, $R$ always denotes a commutative ring, $R_r$ is the localization of $R$ at $r\in R$, $\Max(R)$ is the subset of $\Spec(R)$ consisting of maximal ideals of $R$, $\Ralg$ is the category whose objects are commutative associative unital $R$--algebras and whose morphisms are unital $R$--algebra homomorphisms. 

An $R$--functor is a functor from the category $\Ralg$ to the category of sets. We will usually denote $R$--functors by underlined letters and use the same letter in bold to denote the scheme representing a given $R$--functor, should it exist.
An $R$--scheme is a scheme over $\Spec(R)$. For $T\in \Ralg$ and an $R$--scheme $X$ we abbreviate $X_T = X\times_R T = X \times_{\Spec(R)} \Spec(T)$. For an affine scheme $X = \Spec(A)$, $A\in \Ralg$,  and $f\in A$ we denote by $X_f = D(f)\cong \Spec(A_f)$ the principal open subscheme determined by $f$. We have a canonical isomorphism of $T$--schemes
\begin{equation}
  \label{agbas1} X_f \times_R T \cong (X\times_R T)_{f\ot 1_T}. 
\end{equation} 
\lv{
Details: 
\begin{align*}
  X_f \times_R T &= \Spec(A)_f \times_{\Spec(R)}\Spec(T) \qquad (\text{Definition of $X_f$})
  \\ &= \Spec(A_f) \times_{\Spec(R)} \Spec(T)
     \qquad  (\text{\cite[00E4]{St}}) 
  \\ &= \Spec(A_f \ot_R T)\qquad \qquad \qquad (\text{\cite[01JQ]{St}})
  \\ & = \Spec( (A \ot_R T)_{f\ot 1_T})  
            \qquad  \qquad (\text{\cite[II, \S4.7, Prop.~18]{BAC}})
  \\&= \big(\Spec(A\ot_R T)\big)_{f\ot 1_T} 
     = (X \times_R T )_{f\ot 1_T}
\end{align*} 
}

\item\label{ag-ex} ({\em Some $R$--schemes}) An $R$--module $M$ defines an $R$--functor $\ulW(M)$ by $A(\in \Ralg) \mapsto M \ot_R A=M_A$ and given in the obvious way on morphisms in $\Ralg$. If $M$ is finitely presented, the $R$--functor $\ulW(M)$ is representable if and only if $M$ is finite projective \cite[5.4.5]{Romagny}. In this case we write $\uW(M)$ for the scheme representing $\ulW(M)$; it is given by $\Spec( \Sym(M\ch))$, where $M \ch$ is the $R$--dual of $M$ and $\Sym(\cdot)$ is its symmetric $R$--algebra. In particular, $\uW(M)$ is an affine smooth finitely presented $R$--scheme with geometrically connected fibres. For clarity, we sometimes write $\uW(M) = \uW_R(M)$. 
    
    For an $R$--module $M$ and $R'\in \Ralg$ we have 
    \begin{equation} \label{ag-ex1}
       \uW_R(M) \times_R R' \cong \uW_{R'}(M\ot_R R'). 
    \end{equation}  

    Let $B$ be an $R$--algebra, whose underlying $R$--module is finite projective. The $R$--functor of automorphisms of $B$, defined by $A \mapsto \Aut(B_A)$, is represented by an $R$--scheme, denoted $\uAut(B)$, \cite[II, \S1, 2.6]{DG}. 
    The $R$--functor $A \mapsto (B \ot_R A)\ti$ is representable by an $R$--scheme, denoted $\uGL_1(B)$ or $\uGL_{1,R}(B)$ if it is helpful to indicate $R$, \cite[2.4.2.2]{CF}.   
      \sm

\item In this appendix, we mostly consider schemes over a base scheme $S$, also called $S$--schemes. For an $S$--scheme $X$ and $x\in X$, we denote $\ka(x) = \calO_{X,x}/\m_x$ where $\calO_{X,x}$ is the local ring of $X$ at $x$ and $\m_x$ its maximal ideal. For $S$-schemes $X$ and $Y$ we put $X(Y) = \Mor_S(Y,X)$, the $S$--morphisms from $Y$ to $X$, and abbreviate $X(\Spec(A)) = X(A)$ for any affine scheme $Y=\Spec(A)$. 
\sm 

\item \label{ag-c} A scheme $X$ is {\em connected} ({\em irreducible} respectively) if its underlying topological space is connected (irreducible respectively). 
    It is {\em reduced\/} if all local rings are reduced (= no non-zero nilpotent elements); it is {\em integral\/} if it is reduced and irreducible.
    An affine scheme $X=\Spec(A)$ is reduced (integral) if and only if $A$ is a reduced ring (an integral domain), \cite[Tag 01J2]{St} for ``reduced''.

    Let $\sfP$ be one of the following properties of a scheme over a field: ``irreducible'', ``connected'', ``reduced'', ``integral''. One says that a scheme $X$ over a field $k$ has property $\sfP$ {\em geometrically\/} if for all fields $F\in \kalg$ the $F$--scheme $X_F = X\times_k F$ has property $\sfP$.     
     Let $K$ be an algebraically closed extension of $k$. Then the $k$--scheme $X$ has property $\sfP$ geometrically if and only if the $K$--scheme $X_K$ has property $\sfP$  (\cite[IV$_2$, (4.5.2) and (4.6.5)]{EGA} or \cite[Cor.~5.54]{GW}).
\sm

\item\label{ag-essen} ({\em Essentially free morphisms}) A morphism $f \co X \to S$ is {\em essentially free\/} if there exists an open affine covering $(S_i)_{i\in I}$ of $S$, for every $S_i$ an affine and faithfully flat scheme $S_i$--scheme $S_i'= \Spec(R_i')$ and an open affine covering $\big( \Spec(A_{ij})\big)$ of $X\times_S S_i'$ such that every $A_{ij}$ is a free $R_i'$--module. 
    See \cite[VI$_B$, D\'ef.~6.2.1]{SGA3}, or \cite[Exc.~14.20]{GW} where this type of morphism is called ``locally free''. 

\sm

\item \label{ag-smock} ({\em Smooth schemes over a field}) Let $k$ be a field and let $X$ be a smooth $k$--scheme. Then $X$ is locally of finite type, geo\-metrically reduced and geo\-metrically regular in the sense of \eqref{ag-c}, see \cite[Tag 056T]{St} for the last two properties. Any smooth scheme is locally of finite presentation and so in particular locally of finite type if the base is a field. \sm

\item\label{ag-lp} ({\em Schemes locally of finite presentation}) An $S$--scheme $X$ is locally of finite presentation if and only if for every directed set $I$ and every projective system $(T_i)_{i\in I}$ of affine schemes the canonical map
    \[ \textstyle \limind_{\, i\in I}  \Mor_S(T_i, X) \simlgr
         \Mor_S(\limproj T_i, X) \]
    is a bijection (\cite[IV$_3$, 8.14.2]{EGA}, see also \cite[Tag 01ZC]{St}). 
\sm

\item({\em Reducedness criterion}) \label{ag-red}
Let $R$ be an integral domain with fraction field $K$. If $X$ is a flat $R$--scheme and $X_K$ is a smooth $k$--scheme, then $X$ is reduced.

Indeed, reducedness being a local property, we can assume that $X=\Spec(A)$ for a flat $R$--algebra $A$ \cite[Tag 01U5]{St}. By \eqref{ag-smock}, all local rings of the $K$--scheme $X_K=\Spec(A\ot_R K)$ are reduced. Hence, by \eqref{ag-c}, the ring $A\ot_R K$ is reduced. Flatness of the $R$--module $A$ implies that $A=A\ot_R R$ embeds into $A\ot_R K$. Therefore $A$ is reduced and then $X=\Spec(A)$ is reduced by \eqref{ag-c} again.

A particularly interesting case where \eqref{ag-red} applies is that of a smooth scheme $X$ over the integral domain $R$. Then $X$ is flat and $X_K$ is smooth by  \cite[Tags 01VF and 01VB]{St}. 
Hence $X$ is reduced. \sm

\item\label{ag-des} ({\em Descent of smoothness}) Let $X$ be an $R$--scheme, locally of finite presentation and let $S\in \Ralg$ be faithfully flat. Then the $S$--scheme $X_S$ is smooth if and only if $X$ is smooth. This is a special case of \cite[IV$_4$, 17.7.3(ii)]{EGA}.%
\sm

\item \label{ag-du} ({\em Disjoint union}) Given a family $(X_\al)$ of schemes, there exists a unique scheme $X= \bigsqcup_\al X_\al$ whose 
underlying topological space is the disjoint union of the topological spaces underlying the $X_\al$ and whose structure sheaf canonically extends those of the $X_\al$, \cite[I. \S3.1]{EGA}. We refer to this scheme as the {\em disjoint union}. 

Assume $S= \bigsqcup S_\al$ for some family $(S_\al)$ of schemes. If every $X_\al$ is an $S_\al$--scheme, we view $\bigsqcup_\al X_\al$ canonically as a scheme over $S$. 

Suppose $X_1 = \Spec(A_i), \ldots, X_n = \Spec(A_n)$ is a finite family of affine schemes. Then $X_1 \sqcup \cdots \sqcup X_n$ is an affine scheme: 
\[ 
    \Spec(A_1) \sqcup \cdots \sqcup \Spec(A_n) = \Spec(A_1 \times \cdots \times A_n)
\] 
\cite[Tag 00ED]{St}. \sm 

\item \label{ag-d} ({\em Fibrewise Isomorphism (Open Immersion) Criterion} \cite[IV$_4$, 17.9.5]{EGA} or \cite[Exc.~3.4.3]{Co1}) Let $S$ be a scheme and let $f\co X \to Y$ be a morphism between $S$--schemes $X$ and $Y$ which are locally of finite presentation. If, moreover,  $X$ is a flat $S$--scheme, the following are equivalent:
\end{inparaenum}

\begin{enumerate}[label={\rm (\greek*)}]

\item \label{grofibi} $f$ is an isomorphism (an open immersion respectively);
    \sm

\item\label{grofibii} for every $s\in S $ the induced map $f_s \co X \times_S
    \ka(s) \to Y \times_S \ka(s)$ is an  isomorphism) of schemes  (an open immersion respectively); \sm

\item \label{grofibiii} for every $s\in S $ the induced map $f_{\bar s} \co X \times_S \overline{\ka(s)} \to S \times_S \overline{\ka(s)}$ is an isomorphism of $\overline{\ka(s)}$ (an open immersion respectively), where $\overline{\ka(s)}$ is an algebraic closure of $\ka(s)$.
\end{enumerate}
The equivalence of \ref{grofibi} and \ref{grofibii}
is \cite[loc.\ cit.]{EGA}. The equivalence of \ref{grofibii} and
\ref{grofibiii} follows by faithfully flat descent (\cite[IV$_2$, 2.7.1(viii)]{EGA}).

The following lemma is folklore.

\comments{(2025-12-05) Lemma~\ref{smok} is taken from the old appendix C; it justifies irreducible = connected in \ref{respr}\ref{respriii}. We need irreducible for the denseness criterion in \ref{zdens}\eqref{S-dense}.  \sm 


(2026-03-11) Lemma~\ref{smok} is not quoted in the text. There is a related result in \cite[Cor.~16.52]{GW}: 

{\tt Let $G$ be a group scheme locally of finite type over a field $k$. Then $G$ is geometrically irreducible if and only if $G$ is connected.} \sm 

(2026-03-12) New version of \ref{smok}, replaced "smooth" by "regular" and added first part. Note that ``regular'' is a property of a scheme, while ``smooth' is a property of a morphism. Therefore there is no "smooth" in the first part of \ref{smok}.}

\new
\begin{lem}\label{smok} Let $X$ be a regular scheme.  
Then the following are equivalent:
\begin{enumerate}[label={\rm (\roman*)}]
 \item  \label{smok-i} $X$ is integral;
 \item  \label{smok-ii} $X$ is irreducible;
 \item  \label{smok-iii} $X$ is connected.
\end{enumerate}
If $X$ is a regular scheme over a field, e.g., a smooth scheme, then the geometric version of \ref{smok-i}--\ref{smok-iii} are equivalent:  
\begin{enumerate}[label={\rm (\roman*)$^\prime$}]
 \item  \label{smok-i'} $X$ is geometrically integral;
 \item  \label{smok-ii'} $X$ is geometrically irreducible;
 \item  \label{smok-iii'} $X$ is geometrically connected.
\end{enumerate}
\end{lem}

\begin{proof} 
The implications \ref{smok-i} $\implies$ \ref{smok-ii} $\implies$ \ref{smok-iii} hold for any scheme. A connected regular scheme is integral by \cite[3.5.6]{Po}. The equivalence of the geometric properties in the field case then follows. Any smooth scheme over a field is geometrically regular by \cite[Tag 056T]{St}.  
%
\end{proof}
\enew


\begin{lem}\label{opemaLG} Let $U$ be an open subscheme of an $S$--scheme $X$, and let $\Spec(A)$ be an affine scheme over $S$. For $T\in \Aalg$ and $x\in X(A)$ we let $x_T = x\circ \can_T$ be the composition 
\[ \Spec(T) \xrightarrow{\can_T} \Spec(A) \xrightarrow{\;x\;} X 
\]
where $\can_T$ is the canonical structure morphism. Then 
\begin{align}  
\label{opemaLG1}  
U(A) &= \big\{ x\in X(A): \;  x_{A/\m}\in U(A/\m)
      \hbox{ for all $\m\in \Max(A)$}\big\}.
\end{align}
\end{lem}

Lemma~\ref{opemaLG} is likely well-known to the experts. But since we could not find a suitable reference, we include a proof.

\begin{proof} 
We can identify
\begin{align*}
U(A) & = \Mor(\Spec(A), U) 
\\ &\equiv \{ x\in \Mor(\Spec(A), X) : 
    x(\p) \in U \hbox{ for all $\p \in \Spec(A)$}\}
\end{align*} 
via the inclusion morphism $U \hookrightarrow X$. In particular, we have the inclusion from left to right in \eqref{opemaLG1}. Regarding the other inclusion, note that for $\m \in \Max(A)$ the morphism $\can_{A/\m}$ maps the unique point of $\Spec(A/\m)$ to $\m \in \Spec(A)$. Hence $x_{A/\m} = x(\m)$. An equivalent version of \eqref{opemaLG1} is therefore that for any $x\in X(A)$ we have
\begin{equation}  \label{opemaLG1a}
 x(\m) \in U \hbox{ for all $\m\in \Max(A)$} 
 \implies
  x(\p) \in U \hbox{ for all $\p \in \Spec(A)$}.   
\end{equation}
Suppose there exists $\p \in \Spec(A)$ such that $x(\p) \in X\setminus U$. There exists $\m \in \Max(A)$ such that $\p \subset \m$, i.e., $\m \in V(\p) = \overline{\{\p\}}$. 
By continuity of $x$ we then get $x(\m) \in x(\overline{\{\p\}}) \subset \overline{x(\p)}
\subset X\setminus U$, contradiction. We have thus established  \eqref{opemaLG1}.  
\end{proof}

\subsection{Application} \label{opema-rem} A typical way in which we will apply Lemma~\ref{opemaLG} is the following. In the setting of the lemma suppose that the canonical map
$\vphi \co X(A) \to \textstyle \prod_{\m \in \Max(A)} \, X(A/\m)$
is surjective. Then its restriction
\[ \vphi_U \co  U(A) \to \textstyle \prod_{\m \in \Max(A)} \, U(A/\m) \]
is surjective too. In particular,
\begin{equation}  \label{opema-rem1}
 U(A/\m) \ne \emptyset \;\; \forall \m \in \Max(A) \implies U(A) \ne \emptyset.
\end{equation}


\subsection{Content of a polynomial}\label{cop} Although we will here only use the example in \eqref{cop-c}, it is appropriate to first review the proper background. \sm 

\begin{inparaenum}[(a)] \item\label{cop-a} ({\em Kernel of two morphisms})
Let $f,g\co Y \to Z$ be two morphisms of $S$--schemes. Given an $S$--scheme $T$ and $\si \in S(T)$ we get the following commutative diagram
\[\xymatrix@C=60pt{
    Y\times_\si T \ar@<0.5ex>[dd]^{g\times_\si T}\ar@<-0.5ex>[dd]_{f\times_\si T}  \ar[rrr]\ar[dr] &&& Y 
    \ar@<0.5ex>[dd]^g\ar@<-0.5ex>[dd]_f \ar[dl] \\
        & T \ar[r]^\si & S \\
 Z\times_\si T \ar[ur]\ar[rrr] &&&  Z\ar[ul]
 } \]
where $Y\times_\si T$ and $Z\times_\si T$ are the fibre products with respect to $\si$. The assignment $\ul \Ker(f,g)(T) = \{ \si \in S(T): f\times_\si T  = g\times_\si T\}$ (equality of scheme morphisms) extends to an $S$--functor $\ul \Ker(f,g)$. 

Assuming that $Y \to S$ is \new an essentially free \enew morphism, \ref{ag}\eqref{ag-essen}, 
and that $Z$ is separated, $\ul \Ker(f,g)$ is represented by a closed subscheme $\Ker(f,g)$ of $S$ \cite[VI$_B$, 6.2.4(b)]{SGA3}, called the {\em kernel of $f$ and $g$}. Thus, $\Ker(f,g)$ is the largest closed subscheme $S_0$ of $S$ such that the restrictions $f|_{S_0}$ and $g|_{S_0}$ coincide. Note that there is a priori no reason for the existence of such a largest closed subscheme. 
  
Furthermore, by \cite[loc.~cit.]{SGA3}, if $Z$ is locally of finite type, then $\Ker(f,g)$ is a finitely presented  $S$--scheme. \sm 

\item\label{cop-b} ({\em Affine base $S=\Spec(R)$}) The assumption on $Y$ in \eqref{cop-a} is for example fulfilled for $Y=\uW(M)$ and $M$ a finite locally free $R$--module (recall that $\uW(M) = \Spec\big(\Sym(M\ch)\big)$ and that the symmetric algebra $\Sym(M\ch)$ is a projective $R$--module by \cite[\S6.6, Cor. of Thm.~1]{BA}). Given an $S$--morphism $f\co \uW(M) \to Z$ there exists a unique finitely generated ideal $I(f)\ideal R$ such that $\Ker(f, 0)$ is the affine scheme $\Spec\big( R/I(f)\big)$. Thus, $I(f)$ is the smallest ideal $J$ of $R$ such that $f|_{\Spec(R/J)} = 0$.  
    
    In particular, for $Z=\uW(R)=\Spec(R)$, a morphism $f \co \uW(M) \to \Spec(R)$ can be identified with an element $f\in \Sym_R(M\ch)$. In this case, we call $I(f) = \Cont(f)$ the {\em content ideal of $f$\/} because of the example \eqref{cop-c}. The defining property of $\Spec\big( R/\Cont(f)\big)$ translates to the following: $\Cont(f)$ is the smallest ideal $J$ of $R$ such that $\ze_J(f \ot 1_{R/J}) = 0$, where $\ze_J$ is the canonical isomorphism of $(R/J)$--algebras, 
\[ \ze_J \co \Sym_{R}(M\ch) \ot_R (R/J) \simlgr \Sym_{R/J}\big( (M\ot_R R/J)\ch\big) .
 \] %

\item\label{cop-c} ({\em Example $S=\Spec(R)$, $Y=\uW(R^n)$ and $Z=\uW(R)$}) In this case, 
\[ f \co \bbA_R^n = \uW(M) \to \uW(N) = \bbA_R^1 \] 
is nothing but a polynomial $P \in R[X_1, \ldots, X_n]$, and the content ideal $\Cont(f)$ is the content ideal in the traditional sense, i.e., the ideal of $R$ generated by the coefficients of $P$. Indeed, both ideals are the smallest ideals $J\ideal R$ such that $P=0\in (R/J)[X_1, \ldots, X_n]$. 
\end{inparaenum}

\subsection{Zariski density and $S$--density}\label{zdens} Let $X$ be a scheme over a base scheme $S$. 
 
\begin{inparaenum}[(a)] \item\label{std} We say that a subset $A$ of the topological space underlying $X$ is {\em Zariski dense\/} or simply is {\em dense\/} if it is so as topological space. 

{\em Standard example}: A open non-empty subscheme of an irreducible scheme is dense. 
\sm 

\new 
\item \label{adc} ({\em Zariski density of rational points over fields}) Let $k$ be a field. 
    Following \cite[2.1.1]{Po}, a {\em variety over $k$\/} or a {\em $k$--variety\/} is a separated scheme $X$ of finite type over $S=\Spec(k)$. We will add adjectives such as ``irreducible'' for more restrictive classes of varieties. For example, a ``variety'' in the sense of \cite[Tag 020D]{St} is an integral $k$--variety in our sense.      
    Also recall that a geometrically integral $k$--variety $X$ is {\em rational\/} ({\em unirational\/}), if there exists a rational map from an affine $k$-space $\bbA^n_k \to X$, which is birational (dominant respectively). In particular, any open subscheme of an affine $k$--space is rational and a fortiori unirational. 
We will use the following (of the many) denseness criteria for rational points. 

Given a field extension $F/k$ and a non-empty $k$--scheme $X$ which is locally of finite type,  the set $X(F)$ of $F$--rational points is Zariski dense in $X$ in any one of the following three cases: 
\end{inparaenum}
\begin{enumerate}[label={\rm (\roman*)}]
 \new
  \item\label{adci} $F$ is algebraically closed;
  
  \item\label{adcii} $F$ is separably closed and $X$ is a geometrically reduced $k$--variety;
      
  \item\label{adciii} $k=F$ is infinite and $X$ is a unirational $k$--variety.     
\end{enumerate}
\new
In particular, in any of these cases, if $X\ne \emptyset$, then $X(k)$ is infinite.

\inparcom{(2026-03-14) The assumption $k=F$ in \ref{adciii} is due to the fact that it is not clear to me that being unirational is stable under base change (did not have time to think about this). It likely is, see \cite[Def.~9.2.26]{Po}, so that we can discard $k=F$. }

For \ref{adci} and \ref{adcii} we can assume that $k=F$ after base change. Then \ref{adci} and \ref{adcii} are proven in \cite[Cor.~3.36]{GW} and \cite[Cor.~3.5.71]{Po} respectively or see \cite[11.2.5]{Springer} for $X$ irreducible. A proof of \ref{adciii} is given in \cite[AG 13.7]{Bo} and again in \cite[13.2.6]{Springer}.
\enew
\lv{
Background: Let $k$ be an infinite field, $X=\AA^n_k$, and $f\in k[X_1, \ldots, X_n]$ (in our case $k=\ka_i$ and $f = \dis \circ h^{u,m}$). Then $X_f$ contains a $k$--rations point, i.e., there exists $a\in k^n$ such that $f(a) \ne 0$. The easy proof is given on MathOverflow, question "Rational points on open subsets of affine space".
}
\sm 

\begin{inparaenum}[(a)]\setcounter{enumi}{2}
\item \label{S-dense} ({\em $S$-density}) Slightly generalizing \cite[\S 2.5]{BLR}, we say that an open subscheme $U$ of $X$ is {\em $S$-dense in $X$\/} if for every $s\in S$ the fiber $U_s = U\times_S \kappa(s)$ is Zariski dense in $X_s = X \times_S \kappa(s)$. If $S=\Spec(R)$, we use the term {\em $R$--dense\/} instead of $S$--dense. It is in this sense that we have used the term in our paper \cite{GN-LG}. 

\inparcom{(2026-01) Unfortunately, we have not defined ``$R$--dense'' in our paper \cite{GN-LG}. }    
    \sm 
     
    Of course, for schemes over a field $R$, Zariski dense equals $R$--dense.  But this is not so for arbitrary $R$: an $R$--dense open subscheme of a smooth scheme is dense,
but a dense open subscheme is not necessarily $R$--dense. For example, let $R$ be a DVR with fraction field $K$ and residue field $\ka$.  The open subscheme $U = \Spec(K)$ of $X = \Spec(R)$ is
dense, since $\Spec(K)$ is the generic point of $X$, but $U$ is not $R$--dense because $K \ot_R \ka = 0$
\lv{
Since $R$ is a DVR, then $\m=(p)$ for some $p\in R$ (the {\em uniformizer}). But $p$ is invertible in $K$ since $p\ne 0$. Thus for $x\in K$, $y\in R/\m$ we have
$x\ot y = x \cdot \frac{p}{p} \ot y = (x/p) \ot (py) = (x/p) \ot 0 = 0$.
We have $U_\ka = \Spec(K) \times_{\Spec(R)} \Spec(\ka) = \Spec(K \ot_R \ka) = \Spec(0) = \emptyset$. 
}
and so $U_\ka = \emptyset$. \sm 

%
\end{inparaenum}

\subsection{Schematic denseness} \label{tod}
We follow the definition of a {\em schematically dominant morphism} $Z \to X$ of schemes and a {\em schematically dense subscheme $Z$ of a scheme $X$} given in \cite[I, (5.4.2)]{EGA-neu}, which is a special case of that of \cite[IV$_3$, (11.10.2)]{EGA} and \cite[IX, D\'ef.~4.1]{SGA3} where a morphism is replaced by a family of morphisms. Thus, $f\co Z \to X$ is a schematically dominant morphism if and only if 
\begin{enumerate}[label={\rm (\roman*)}]  
\item \label{tod-0i} $f^\flat \co \scO_X \to f_*(\scO_Z)$ is injective, equivalently, 
\item for every open subset $V\subset X$ the homomorphism 
 $f^\flat(V) \co \Ga(V, \scO_X) \to \Ga(f\me(V), \scO_Z)$ is injective.
\end{enumerate}
A subscheme $j \co Z\to X$ is called schematically dense if $j$ is schematically dominant. 

The definition above is the same as that of \cite[9.18, 9.19]{GW}. The reader should be warned that Remark 9.20(1) of the published version of the book \cite{GW} is not correct, see the book's errata at {\tt www.algebraic-geometry.de}. The correct statement is given in \eqref{tod-a} below.

In \cite[Tag 01RB]{St}, an open subscheme $U$ of a scheme $X$ is called {\em scheme theoretically dense\/} if 
\begin{enumerate}[label={\rm (\roman*)}]\setcounter{enumi}{2}  
 \item \label{tod-0iii} for every open subscheme $V\subset X$ the scheme theoretic closure of $U\cap V$ is $V$.   
\end{enumerate}
That \ref{tod-0iii} is equivalent to \ref{tod-0i}, follows from \cite[Tag 01RE]{St}. 
\sm 

\begin{inparaenum}[(a)] \item\label{tod-a} ({\em Schematically dense versus Zariski dense}) Let $f \co Z \to X$ be a morphism of schemes. If $f$ is schematically dominant, it is dominant in the sense that $|f(Z)|$ is a dense subset of the topological space $|X|$. Conversely, if $X$ is reduced and $f$ is dominant, then $f$ is schematically dominant. Thus, {\em for a reduced $X$, dominant = schematically dominant}, see \cite[I, (5.4.3)]{EGA-neu}, or \cite[IV$_3$, (11.10.4)]{EGA}.  \sm 

\item\label{tod-b} Let $f\co Z\to X$ be a morphism of $S$--schemes. Then $f$ is called {\em universally schematically dominant relative to $S$} if for every base change $S'\to S$ the corresponding morphism $f'\co Z' = Z\times_S S' \to X' = X \times_S S'$ is schematically dominant.  
One calls an immersion $Z\to X$ or even just $Z$ {\em universally schematically dense relative to $S$\/} if it is universally schematically dominant relative to $S$, \cite[IV$_3$, (11.10.8)]{EGA}. We will abbreviate universally schematically dense relative to $S$ by {\em universally schematically dense\/} if $S$ is clear from the context or unimportant.  

{\em Example}: If $S=\Spec(k)$, $k$ a field, then schematically dominant = universally schematically dominant, \cite[IV$_3$, (11.10.6)]{EGA}. 
\sm 

\item\label{tod-c} (\cite[IV$_3$, (11.10.10)]{EGA}) Let $X \to S$ be flat and locally of finite presentation, and let $U$ be an open subscheme of $X$. Then $U$ is universally schematically dense relative to $S$ if and only if for every $s\in S$ the open subscheme $U_s = U \cap X_s$ of the $\Spec(\ka(s))$--scheme $X_s$ is schematically dense, equivalently, by \eqref{tod-a}, every $U_s$ is universally schematically dense in $X_s$. 
Hence, since schematically dense implies Zariski dense by \eqref{tod-a}, we have the implications  
\begin{equation}\label{tod-c1} \begin{split}
 &\text{\em $U$ is universally schematically dense relative to $S$} 
 \\ &\quad \implies 
     \text{\em $U$ is $S$--dense} 
\end{split} \end{equation}
for $S$--dense in the sense of \ref{zdens}\eqref{S-dense}. \sm


\item\label{tod-d} Suppose $X\to S$ is flat, locally of finite presentation and all fibres $X_s$, $s\in S$, are reduced. Then, by combining \eqref{tod-a} and \eqref{tod-c}, for any open subscheme $U$ of $X$ we have
\begin{equation}\label{tod-d1} \begin{split} 
&\text{\em $U$ is universally schematically dense relative to $S$} 
\\ &\qquad \iff \text{\em $U$ is $S$--dense.} 
\end{split} \end{equation} 

\item\label{tod-e} ({\em $X$ smooth})
Assume $f\co X\to S$ is smooth. Then $f$ is flat, locally of finite presentation and all fibres are smooth by \cite[Tags 01VE, 01VF, 01VB]{St}, hence they are in particular reduced by \cite[056T]{St}. It follows from \eqref{tod-d} that {\em $U$ is $S$--dense as defined in \cite[2.5]{BLR} if and only if $U$ is universally schematically dense relative to $S$.} 
Moreover, since universally schematically dense implies schematically dense, \cite[Lem.~2.5(b)]{BLR} follows from \eqref{tod-d1}. 
\sm 

\new 
\item\label{tod-f} (\cite[XVIII, Prop.~1.7]{SGA3}) Let $X\to S$ be faithfully flat and locally of finite presentation, and let $U$ be a universally schematically dense open subscheme of $X$. Then there exists an fppf covering $S'\to S$ for which $U(S') \ne \emptyset$. For example, one can take $S'= U$. \sm 

In the following Lemma~\ref{todle} we will consider other situations in which $U(Y)\ne \emptyset$, namely $Y=\Spec(F)$ for certain $S$--fields $F$. Recall that a field $F$ is an $S$--field if there exists a morphism of schemes $\Spec(F) \to S$. \end{inparaenum}
\enew
\ms

\begin{lem}\label{todle} Let $X$ be an $S$--scheme which is flat, separated, finitely  presented and has non-empty geometrically reduced fibres. Moreover, let $U$ be an open subscheme of $X$.
\sm 

\begin{inparaenum}[\rm (a)] \item \label{todle-b} Then the conditions \ref{todle-bi}--\ref{todle-bii} below  are equivalent. \end{inparaenum}

\begin{enumerate}[label={\rm (\roman*)}]
  \item \label{todle-bi} $U$ is universally schematically dense in $X$; 
  
  \item\label{todle-bis} $U(F) \ne \emptyset$ for all separably closed $S$--fields $F$; 
  
  \item \label{todle-bii} $U(K) \ne \emptyset$ for all algebraically closed $S$--fields $K$. 
\end{enumerate}  
\sm 

\begin{inparaenum}[\rm (a)]\setcounter{enumi}{1}
  \item \label{todle-c} Suppose in addition that $X$ has geometrically  integral fibres and that $X\to S$ is faithfully flat. Then the conditions \ref{todle-bi}--\ref{todle-bii} of {\rm \eqref{todle-b}} are equivalent to \end{inparaenum}

\begin{enumerate}[label={\rm (\roman*)}]\setcounter{enumi}{3}
 \item\label{todle-ciii} there exists an fppf covering $S'\to S$ for which $U(S') \ne \emptyset$.  
\end{enumerate}
\end{lem}

\begin{proof} \eqref{todle-b} \ref{todle-bi} $\implies$ \ref{todle-bis}: 
Since $F$ is an $S$--field, there exists a unique $s\in S$ such that the structure morphism $\Spec(F) \to S$ factors via the canonical morphism $\Spec(k) \to S$ where $k$ is the residue field at $s$. By \ref{tod}\eqref{tod-a}, we know that the fibre $U_k$ is Zariski dense in $X_k$. Since the latter $k$--scheme is not empty by assumption, so is $U_k$. By base change, $X_k$ and then also $U_k$ are separated, geometrically reduced and  finite type $k$--schemes. The claim then follows from \ref{adcii} of \ref{zdens}\eqref{adc}.

\sm 

The implication \ref{todle-bis} $\implies$ \ref{todle-bii} being trivial, it remains to prove that  \ref{todle-bii} $\implies$ \ref{todle-bi}.  Applying \ref{tod}\eqref{tod-c}, it suffices to show that all fibres $U_s \subset X_s$ are schematically dense. Thus, without loss of generality, we can suppose that $S=\Spec(k)$ for a field $k$. Let $K$ be an algebraically closed extension field of $k$. Then $U_K$ is an open non-empty  subscheme of the irreducible $K$--scheme $X_K$, and is therefore Zariski dense. But  $X_K$ is also reduced. Hence, $U_K$ is schematically dense  by \ref{tod}\eqref{tod-a}. By \cite[IV$_3$, Thm.~11.10.5(i)]{EGA}, schematically dense allows faithfully flat descent. Thus, $U$ is schematically dense, as required. 
\sm 

\eqref{todle-c} The implication \ref{todle-bi} $\implies$ \ref{todle-ciii} is \cite[XVIII, Prop.~1.7]{SGA3}, restated as \ref{tod}\eqref{tod-f}.  
For the converse,  let $S'\to S$ be an fppf covering as in \ref{todle-ciii}. 
As mentioned above, the property ``universally schematically dense'' allows flat descent. It is therefore sufficient to prove that the open subscheme $U'= U \times_S S'$ of $X'=X\times_S S'$ is universally schematically dense in $X'$. The assumption $U(S') \ne \emptyset$ implies $U'(K')\ne \emptyset$ for some algebraically closed $S'$-field $K'$, so that $U'$ is universally schematically dense by \eqref{todle-b}. 
\end{proof}
\sm

\textbf{Remarks.}
Lemma~\ref{todle} is inspired by \cite[Lem.~1.5]{Lo-ag}, retaken in \cite[0.14]{Lo-genalg},  where both times the result is stated without proof for $X$ a smooth separated finitely presented scheme with non-empty connected fibres over an affine base. 
\comments{(2026-01-06) I believe that Loos meant geometric fibres instead of fibres. Otherwise, I don't see how his claim makes sense. It is not clear to me why he assumed ``separated''. 
It is too bad that he did not give indications of a proof, besides the references to \cite{SGA3} and \cite{EGA}. Unfortunately, I never talked to him about this. 
}
\sm 

\begin{cor}[Example: Principal open subschemes of $\uW(M)$]\label{usd-d} Let $M$ be a finitely generated projective $R$--module, let $f\in \Sym(M\ch)$ be a polynomial and let $\Cont(f)$ be its content ideal defined in {\rm \ref{cop}\eqref{cop-b}}. Then the following are equivalent: \sm 

\begin{enumerate}[label={\rm (\roman*)}]
 \item\label{usd-di} The principal open subscheme $\uW(M)_f$ of $\uW(M) = \Spec\big(\Sym(M\ch)\big)$ is universally schematically dense. \sm

 \item \label{usd-dis} $\uW(M)_f(F) 
\ne \emptyset$ for all separably closed extension fields $F$ of $R$.%
\sm

\item \label{usd-dii} $\uW(M)_f(K) 
\ne \emptyset$ for all algebraically closed extension fields $K$ of $R$.%
\sm

\item\label{usd-diii} $\Cont(f) =R$. \sm

\item \label{usd-div} $f_\ka \ne 0$ for every residue field $\ka= R/\gm$, $\gm$ a maximal ideal of $R$.
\end{enumerate}
\end{cor}

\begin{proof} Since $X=\uW(M)$ is smooth with geometrically integral fibres and $U=\uW(M)_f$ is a principal open subscheme, the equivalence of \ref{usd-di}, \ref{usd-dis}  and \ref{usd-dii} is  a special case of \ref{todle}\eqref{todle-b}. 

If \ref{usd-di} holds, then $(\uW(M)_f)_\ka = \uW(M_\ka)_{f_\ka}$ is dense in $\uW(M_\ka)$ for every residue field $\ka = R/\gm$, and therefore necessarily \ref{usd-div} is fulfilled. In turn, \ref{usd-div} implies $\Cont(f)_\ka = \Cont(f_\ka) = \ka$, consequently \ref{usd-diii} follows. Finally, if \ref{usd-diii} holds,  then  $\Cont(f)_T = T$ for every $T\in \Ralg$. In particular, $f_K \ne 0$ for $K$ as in \ref{usd-dii}. But then $\uW(M_K)_{f_K}$ is well-known to be dense in $\uW(M_K)$, in particular we have \ref{usd-dii}. \end{proof}

\subsection{Zariski covers, Zariski sites} \label{Zarev}
We use the concept of a standard Zariski cover $\{ U_i \to Y = \Spec(A)\}_{i=1, \ldots, n}$ as defined in \cite[Tag 020R]{St}: every $U_i \to Y$ is an open immersion inducing an isomorphism with a standard affine open subscheme of $Y$, say $U_i \cong D(a_i) \subset Y$ for some $a_i \in A$, such that $A= Aa_1 + \cdots + Aa_n$. In this case, we will also refer to $(a_1, \ldots, a_n)$ as a standard Zariski cover of $\Spec(A)$. See \cite[Tag 020T]{St} for the definition of the big affine Zariski site of a scheme $S$ and \cite[Tags 00WM, 00WN]{St} for the definition of an epimorphism between the sheaves on a site, applied in \ref{surlem} to the big affine site of $S$ and the sheaves $X\equiv h_X$ and $Y\equiv h_Y$.

\comments{(2025-08-06) Generalized \ref{transi} to allow arbitrary schemes $G$ and $X$, needed for appendix, $X$ = projective scheme of parabolic subgroups of given type}

\comments{(2025-05-13) Adjusted \ref{surlem} to the terminology and definitions of \cite{St} and added some details in the proof. }

\begin{lem}[Fiberwise Surjectivity Criterion]\label{surlem} Let $S$ be a scheme, and let $f\co X \to Y$ be a morphism of $S$--schemes which are locally of finite presentation. Assume that $Y=\Spec(A)$ is an affine scheme and that for all\/ $\p \in \Spec(A)$ the induced map $f_\p \co X(A_\p) \to Y(A_\p)$, $g \mapsto f \circ g$, is surjective.

Then $f$ admits a section locally for the Zariski topology, i.e., there exist a  standard Zariski cover $(a_1, \ldots, a_n)$ of $A$ and morphisms $s_i \co U_i=\Spec(A_{a_i}) \to X$ such that $f \circ s_i = j_i$, the open immersion associated with the canonical homomorphism $A \mapsto A_{a_i}$:
\[ \xymatrix{X \ar[rr]^f && Y=\Spec(A) \\
&U_i = \Spec(A_{a_i})\ar@{-->}[ul]^{s_i} \ar[ur]_{j_i} }
\]
In particular, if also $X$ is affine, then $f$ is an epimorphism on the big affine Zariski site of\/ $S$.  
\end{lem}

\begin{proof} For $\p \in \Spec(A)$ we let $j_\p \co \Spec(A_\p) \to \Spec(A) = Y$ be the open embedding corresponding to the canonical homomorphism $A\to A_\p$. For $a\in A$ we use the analogous notation for the morphism $j_a \co \Spec(A_a) \to A$. 

By assumption there exists a local section $s_\p \co \Spec(A_\p) \to X$ such that $f \circ s_\p = j_\p$, i.e., the diagram below is commutative:
\[
\xymatrix{X \ar[rr]^f && \Spec(A) \\
&\Spec(A_\p)\ar[ul]^{s_\p} \ar[ur]_{j_\p} }
\]
Since $A_\p = \varinjlim_{a\notin \p} A_a$ and therefore 
$\Spec(A_\p) = \varprojlim_{a\notin \p} \Spec(A_a)$, and since $X$ and $Y$ are locally finitely presented schemes, it follows from \ref{ag}\eqref{ag-lp}
that 
\[\xymatrix@C=50pt{
 \varinjlim_{a\notin \p} \Mor_S (\Spec(A_a), X) \ar[r]^{\varinjlim_{a\notin \p} f_a} 
    \ar[d]_\cong &   \varinjlim_{a\notin \p} \Mor_S (\Spec(A_a), Y)
             \ar[d]^\cong
\\   
 \Mor_S(\varprojlim_{a\notin \p} \Spec(A_a), X) \ar[r]^{f_\p}
  & \Mor_S(\varprojlim_{a\notin \p} \Spec(A_a), Y)
}\]
is a commutative diagram. Since $f_\p$ is surjective, there exist $a[\p]\in A \setminus \p$ and a morphism $s_{a[\p]} \co \Spec(A_{a[\p]} \to X$ such $f\circ s_{a[\p]} = j_{a[\p]}$:  
\begin{equation}\label{surlem1} \vcenter{
\xymatrix{X \ar[rr]^f && \Spec(A) \\
&\Spec(A_{a[\p]})\ar[ul]^{s_{a[\p]}} \ar[ur]_{j_{a[\p]}} }
}\end{equation}
The ideal generated by all these $a[\p] $ is $A$. Thus, it is already
generated by finitely many $a_1, \ldots, a_n$ among the $a[\p]$. Thus, $(a_1, \ldots, a_n)$ is a standard Zariski cover of $A$ for which the first part of the lemma holds. 

The second claim is then an easy consequence: given any affine scheme $\Spec(T)\to S$ and a $y\in Y(T) = \Hom_S(\Spec T, Y)\cong \Hom_{\Ralg}(A, T)$ we need to find a standard Zariski cover $(t_1, \ldots, t_n)$ of $T$ and $x_i\in X(T_i)$ such that $f\circ x_i = y \circ \big(\Spec(T_i)\to \Spec(T)\big)$. Let $X = \Spec(B)$ and put $R=\calO_S(S)$.   After applying the bijection 
\[ 
\flat \co \Mor_S(\Spec(B), \Spec(A)\big) \simlgr \Hom_{\Ralg}(A,B), \quad f \mapsto f^\flat,
\] 
and analogously for the other scheme morphisms, our goal then is to find $\xi_i \in \Hom_{\Ralg}(B, T_{t_i})$ such that $\xi_i  \circ f^\flat = (T\to T_{t_i})\circ y^\flat$, i.e.,  
\begin{equation}
\vcenter{ \label{surlem2} 
 \xymatrix@C=70pt{
    A \ar[r]^{f^\flat} \ar[d]_{y^\flat} & B \ar@{-->}[d]^{\xi_i}\\
   T \ar[r]^{\can} & T_{t_i}  
}}\end{equation} 
is a commutative diagram. The commutative diagram \eqref{surlem1} for $a[\p]$ replaced by $a_i$ becomes
\[ \vcenter{\xymatrix@C=60pt{
B \ar[dr]_{s_{a_i}^\flat}  && A \ar[ll]_{f^\flat} \ar[dl]^{j_{a_i}^\flat} \\ &A_{a_i}
}}\quad . \]
We put $t_i = y^\flat(a_i)$. Then $(t_1, \ldots, t_n)$ is a standard Zariski cover of $T$, and since $T_{t_i} \cong T \ot_A A_{a_i}$ we obtain a well-defined homomorphism $\xi_i \co B \to T \ot_A A_{a_i}$, $b \mapsto 1_T \ot s_{a_i}^\flat (b)$ which makes \eqref{surlem2} commute.  
\end{proof}
\sm 

\comments{(2026-02-03) New Lemma~\ref{fibn}, used in the section on one-generated algebras. Should we change \ref{surlem} to write it for maximal ideals? }

\textbf{Remark.} Lemma~\ref{surlem} also holds under the weaker assumption that $f_\gm$ be surjective for every maximal $\gm \in \Spec(R)$. We leave the straightforward modification of the proof above to the reader. 
\sm 

We remind the reader that for any $R$--scheme $U$ and $T\in \Ralg$ we abbreviate $U_R (T) = \Mor_{\Spec(R)} (\Spec(T), U)$. 

\begin{cor}\label{fibn} 
Let $U$ be an $R$--scheme which is locally of finite presentation. We consider the following conditions:
\begin{enumerate}[label={\rm (\roman*)}]

\item\label{fibni} $U_R(R) \ne \emptyset$. 

\item\label{fibnii} There exists a standard Zariski cover $(a_1, \ldots, a_n)$ of $R$ such that $U_R(R_{a_i}) \ne \emptyset$ for every $i=1, \ldots, n$. 

\item \label{fibniii} $U_R(R_\gp) \ne \emptyset$ for every $\gp \in \Spec(R)$.

\item \label{fibniv} $U_R(R_\gm) \ne \emptyset$ for every maximal $\gm \in \Spec(R)$. 

\item \label{fibnv} $U_R(R/\gm) \ne \emptyset$ for every maximal ideal $\gm \in \Spec(R)$. 
    
\end{enumerate}
Then 
\[ \ref{fibni} \implies \ref{fibnii} \iff \ref{fibniii} \iff 
     \ref{fibniv} \implies \ref{fibnv}. \]
Moreover, suppose that $M$ is a finite locally free $R$--module and $U \to \uW(M)$ is a non-empty quasi-compact open subscheme of $\uW(M)$, then the following hold.

\begin{enumerate}[label={\rm (\alph*)}] 
\item \label{fibn-a} If $R$ is an LG ring, then all $5$ conditions are equivalent (but are not necessarily fulfilled). \sm 
 
 \item\label{fibn-b}  If $R$ is an LG ring and $|R/\gm |=\infty$ for all maximal ideals $\gm \ideal R$, then \ref{fibnv} and hence also the  conditions \ref{fibni}--\ref{fibniv} are fulfilled. 
\end{enumerate}
\end{cor}

\begin{proof} For every $a_i$ as in \ref{fibnii} there exists $\gp\in \Spec(R)$ such that $a_i \notin \gp$, and every $\gp\in \Spec(R)$ lies in a maximal $\gm \in \Spec(R)$. Hence, we have $R$--algebra homomorphisms $R \to R_{a_i} \to R_\gp \to R_\gm \to R/\gm$, therefore $R$--scheme morphisms 
\[ 
\Spec(R/\gm) \to \Spec(R_\gm) \to \Spec(R_\gp) \to \Spec(R_{a_i}) \to \Spec(R), \]
from which the implications   
$\ref{fibni} \implies \ref{fibnii} \implies \ref{fibniii} 
\implies\ref{fibniv} \implies \ref{fibnv}$ follow. 
     
For the proof of $\ref{fibniii} \implies \ref{fibnii}$ we specialize Lemma~\ref{surlem} with $X=U$, $A=R$ and $f$ the $\Spec(R)$--structure map of $U$. The assumption that $f_\gp\co U(R_\gp) \to \big(\Spec(R)\big)(R_\gp)$ be surjective, is fulfilled because there exists a unique $R$--algebra homomorphism $R \to R_\gp$. The first part of \ref{surlem} then proves \ref{fibnii}. As mentioned in the Remark~\ref{surlem}, the lemma also holds for maximal ideals, so we also have \ref{fibniii} $\implies$ \ref{fibnii}. \sm 

Finally, let $R$ be an LG ring and $U$ a non-empty quasi-compact (= finitely presented) open subscheme of some $\uW(M)$. The characterization \ref{prop_baire}\ref{prop_baire-a} of LG rings shows the implication \ref{fibnv} $\implies$ \ref{fibni}, proving \ref{fibn-a}. In the situation of \ref{fibn-b}, we apply the criterion~\ref{zdens}\eqref{adc} to the unirational $R/\gm$--variety $U_{R/\gm}$. Thus, $|U(R/\gm)| = \infty$ for every maximal $\gm \ideal R$, in particular \ref{fibnv} is fulfilled. \end{proof}

\comments{(2026-02-03) A proof that \ref{fibn}\ref{fibnii} $\iff$ \ref{fibn}\ref{fibniv} is here in 'lv'.}
\lv{
\begin{lem}\label{fip} Let $U$ be an $R$--scheme which is locally of finite presentation. Then the following are equivalent for $T\in \Ralg$: 
\begin{enumerate}[label={\rm (\roman*)}]
 \item \label{fibi} $U(T_\gm) \ne \emptyset$ for all maximal $\gm \in \Spec(T)$. 
     
\item \label{fibii} There exists a standard Zariski cover $(t_1, \ldots, t_n)$ of $T$ for which $U(T_{t_i}) \ne \emptyset$ for all $i$, $1\le i \le n$. 
\end{enumerate}
\end{lem}
\begin{proof} \ref{fibi} $\implies$ \ref{fibii} Fix a maximal $\gm\in \Spec(T)$. Then $T_\m = \limind_{\, t \not \in \m} T_t$, hence $\Spec(T_\gm) \allowbreak  = \limproj_{\, t\not\in \m} \Spec(T_t)$ and so 
$U(T_\m)= \limind_{\, t \not \in \m} U(T_t)$,
according to \cite[IV$_3$, 8.14.2]{EGA} (or \cite[Tag 01ZC]{St}). 
It follows that there exists $t_\m \in T \setminus \m $ such that $U(T_{t_\m}) \not=\emptyset$. The $t_\m$'s for $\gm$ running over the maximal ideals of $T$ generate $T$ as ideal. Hence, there exists finitely many maximal ideals $\m_1, \dots, \m_n$ of $T$ such that $T=T t_{\m_1}+ \cdots +  T t_{\m_n}$. 

\ref{fibii} $\implies$ \ref{fibi} For every maximal $\gm \in \Spec(T)$ there exists $i$ and an $R$--algebra homomorphisms $T_{t_i} \to T_\gm$, hence a morphism $\Spec(T_\gm) \to \Spec(T_{t_i})$, which implies \ref{fibi}. \end{proof}
}

\comments{(2026-02-04) Seems to be overkill that we use an infinite field in \ref{fibn-b}. For an infinite field we know that $|U(k)|= \infty$ while we are only interested in $|U(k)| \ge 1$. But I did not find a better criterion. 
For example, if we assume that $U$ is universally schematically dense and $M$ is faithful, is it possible that $U(k) = \emptyset$, while we know that $\uW(M)(k)\ne \emptyset$? }

\pcomments{(2026-03-01) We cannot improve anything. Take the affine line $X=\bbA^1_k$ over a finite field $k$, take $f$ to be the product of the $(x-i)$ for $i$ running over $k$, and consider the principal open subset $U=X_f$. Then $U(k)$ is empty. But $U$ is universally schematically dense, for example by \ref{usd-d}.}

\comments{(2026-04-03) Adjusted \ref{weilres} to arbitrary schemes since this is the setting of this appendix; added \ref{weilres}\eqref{weilres-bc}. }

\subsection{Weil restriction}\label{weilres} 
We recall the Weil restriction of schemes to the extend used here. 
More details can for example be found in \cite[7.6]{BLR}, \cite[A.5]{CGP}, \cite[I, \S1, 6.6]{DG}. \sm 

Let $h \co S'\to S$ be a morphism of schemes, let $X'$ be a contravariant functor from the category of $\Sch_S$ of $S$--schemes to the category of sets, like the functor of points of a scheme $X'$, and let $T$ be an $S$--scheme. The Weil restriction functor $\frR_{S'/S}$ is the right adjoint of the base change functor: there exists a canonical bijection
\[ \Mor_{\Sch_S} \big(T,\, \frR_{S'/S}(X')\big) \simlgr \Mor_{\Sch_{S'}}\, (T\times_S S',\, X')
\]
that is functorial in $T$ and $X'$. If $S=\Spec(R)$ and $S'=\Spec(R')$ are affine schemes, we abbreviate $\frR_{R'/R} = \frR_{S'/S}$. In this case, one has 
\begin{equation}\label{weilres0} 
 \frR_{R'/R}(X')(A)= X'(R'\ot_R A) 
 \end{equation}  
for $A\in \Ralg$. 

Even if $X'$ is representable by an $S'$--scheme, the $S$--functor $\frR_{S'/S}(X')$ need not be representable, see \eqref{weilres-ex} below. However, if it is representable, the representing scheme will also be denoted $\frR_{S'/S}(X')$.

We will use the following results, {\em assuming that $h \co S'\to S$ is finite locally free}, which in the affine case just means that $R'$ is a finite projective $R$--algebra. \sm 

\begin{inparaenum}[(a)]
  \item\label{weilres-a} ({\em Existence}) If $X'$ is a quasi-projective $S'$--scheme, then $\frR_{S'/S}(X')$ is representable by an $S$--scheme \cite[7.6/4]{BLR}. In particular, $\frR_{S'/S}(X')$ is representable if $X'$ is a quasi-affine scheme of finite type. \sm

  \item\label{weilres-b} ({\em Inheritance of scheme properties}) Let $X'$ be an $S'$--scheme for which $\frR_{S'/S}(X')$ exists as an $S$--scheme. If $X'$ has one of the properties below, then the $S$--scheme $\frR_{S'/S}(X')$ has the same property:

  \begin{inparaenum}[(i)]
  
   \quad \item\label{weilres-b0} finitely presented, 
  
    \quad  \item \label{weilres-bi} affine, 

   \quad \item \label{weilres-bii} locally of finite type (locally of finite presentation),

   \quad  \item \label{weilres-biii} smooth.

   \quad \item \label{weilres-biv} quasi-projective, assuming that $R$ and $R'$ are noetherian.
  \end{inparaenum}

\noindent Indeed, \eqref{weilres-bi} is established in the proof of \cite[7.6/4]{BLR}, 
for \eqref{weilres-b0}, \eqref{weilres-bii} and \eqref{weilres-biii} see \cite[7.6/5]{BLR}, 
and \eqref{weilres-biv} is established in \cite[A.5.8]{CGP}.
\sm

\item\label{cgp-a} ({\em Inheritance of geometrically integral fibres}) Let $k$ be a field, let $R'$ be a finite non-zero $k$--algebra, and let $X'$ be a smooth, quasi-projective $R'$--scheme whose structure map $X' \to \Spec(R')$ is surjective. If $X'$ has geometrically connected fibres, then $\frR_{R'/k}(X')$ is a smooth, non-empty, and geometrically integral $k$--scheme. 
    
    Indeed, by (\ref{weilres-biii}) and (\ref{weilres-biv}), $X=\frR_{R'/k}(X')$ is a smooth quasi-projective $k$--scheme, and by \cite[A.5.9]{CGP} it is non-empty and geometrically connected. That $X$ is then geometrically integral, is a consequence of Lemma~\ref{smok}. \sm 

\item \label{weilres-d} ({\em Quasi-compact open immersions}) Let $i' \co X' \to Y'$ be an open immersion of finitely presented $S'$--schemes. Suppose that the Weil restrictions $\frR_{S'/S}(X')$ and $\frR_{S'/S}(Y')$ exist. Then $\frR_{S'/S}(i') \co \frR_{S'/S}(X') \to \frR_{S'/S}(Y')$ is an open immersion, which is quasi-compact, equivalently, is of finite presentation. 
    
    Indeed, $\frR_{S'/S}(i) \co \frR_{S'/S}(X') \to \frR_{S'/S}(Y')$ is an open immersion by \cite[7.6/2]{BLR}. 
    The $R$--schemes $\frR_{R'/R}(X')$ and $\frR_{R'/R}(Y')$ are finitely presented by \eqref{weilres-b0}. Hence $\frR_{S'/S}(i)$ is a finitely presented morphism by \cite[Tag 02FV(2)]{St}. An open immersion is of finite presentation if and only if it is quasi-compact \cite[Tag 01TU]{St}. \sm

\item \label{weilres-bc} ({\em Base change}) Suppose $T\to S$ is another $S$--scheme. Then $T'=S'\times_S T \to T$ is a finite locally free morphism of $S$--schemes, and for every $S'$--scheme $X'$, there is a canonical isomorphism
    \begin{equation}\label{weilres-bc1} 
    \frR_{S'/S}(X') \times_S T \simlgr \frR_{T'/T} (X'\times_{S'} T')          
    \end{equation}%
(\cite[p.~192]{BLR}). \sm

\item \label{weilres-ex} ({\em Example $\uW_{R'}(M')$}) Let $S=\Spec(R)$ and $S'=\Spec(R')$ be affine schemes and let $M'$ be an $R'$--module. By \eqref{weilres0}, the Weil restriction of the $R'$--functor $\ulW_{R'}(M')$ of \ref{ag}\eqref{ag-ex} is given by 
\begin{equation}\label{weilres-3}  
           \frR_{R'/R}\big( \ulW_{R'}(M')\big) = \ulW_R\big( \frR_{R'/R}(M')\big)
\end{equation} 
where $\frR_{R'/R}(M')=: M$ is the scalar restriction of $M'$, i.e., $M'$ considered as $R$--module via the structure map $R\to R'$. If $M'$ is a finite projective $R'$--module, $\uW_{R'}(M')$ is representable. But the $R$--functor $\ulW_R(M)$ is in general not representable \cite[5.4.5]{Romagny}. However, it is representable, if $M'$ is finite projective as $R'$--module and $R'$ is finite projective as $R$--module, since then $M$ is a finite projective $R$--module too, \cite[1.1.8]{Ford}.

Let $R'\in \Ralg$ be finite projective, let $M'$ be a finite projective $R'$--module and let $A\in \Ralg$. Put $A'=R'\ot A \in \Rpalg$. We then have the well-known isomorphism of $A$--modules
\[ \frR_{R'/R}(M') \ot_R A \cong \frR_{A'/A}( M'\ot_{R'} A')
\]
and hence, by \eqref{weilres-bc1}, isomorphisms of $A$--schemes
\begin{equation}\label{weilres-ex2} \begin{split}
 & \frR_{R'/R}\big( \uW_{R'}(M') \big) \times_R A 
   \cong \frR_{A'/A} \big( \uW_{A'}(M'\ot_{R'} A') \big)
 \\ &\quad \cong  \uW_A\big(\frR_{A'/A}(M'\ot_{R'}A')\big)
    \cong  \uW_A ( \frR_{R'/R}(M') \ot_R A)  .
\end{split}\end{equation}
\end{inparaenum}

\subsection{Constant schemes ({\cite[I, \S1, 2.11, 6.9]{DG}})}\label{lcf} Let $S$ be a scheme, let $E$ be a set, and let $(S_e)_{e\in E}$ be a family of schemes with $S_e \cong S$ for all $e\in E$. The {\em constant $S$-scheme\/} based on $E$ is 
\[ E_S = \textstyle \bigsqcup_{e\in E}\,  S_e
\]
with structure morphism given by the isomorphisms $S_e\cong S$. Given an $S$--scheme $X$ and a morphism $f \co X \to E_S$ of $S$--schemes, we associate with $f$ the map $\wtl f \co X \to E$, given by $\wtl f (x) = f(e)$ if $x\in S_e$. Then, denoting by $\Top(X, E)$ the set of continuous (= locally constant) functions between the topological spaces $X$ and $E$, where $E$ is equipped with the discrete topology, we get a bijection 
\begin{equation} \label{lcf-1}  \Mor_{\Sch/S} (X, E_S) \simlgr \Top (X, E) , \quad f \mapsto \wtl f .
\end{equation}
Thus, $E_S$ represents the $S$--functor $\Top( \cdot , E)$. \sm 

We will mostly use the example $E= \NN$. If $\scM$ is a quasi-coherent $\scO_S$--module locally free of finite rank, its rank function $S \to \NN$, $s \mapsto \rank \scM_s$, is a locally constant function.     

\comments{(2025-08-25) My idea is to delete the appendix on ``Parabolic subgroups and Severi-Brauer schemes", and to put here in \ref{grap} what we actually use/need in the appendix on quadrics. The subsection \ref{grap} has to be adjusted to the ongoing project on Severi-Brauer schemes and parabolic subgroups.} 

\subsection{Grassmannians and projective spaces}\label{grap} Let $S$ be a scheme, and let $\scM$ be a quasi-coherent $\scO_S$--module locally free of finite type. If $S=\Spec(R)$ we will identify $\scM = \wtl M$, where $M$ is a finite projective $R$--module. \sm 

\begin{inparaenum}[(a)] \item\label{grap-a} Fix $n\in \NN$. We denote by $\uGr_n(\scM)$ the smooth projective $S$--scheme constructed in \cite[I, 9.7, 9.8]{EGA-neu}, see also \cite[Ex. 13.69]{GW} and \cite[Exc.~18.22]{GWII}. It represents the $S$--functor $\ulGr_n(\scM)$ which assigns to an $S$--scheme $h \co T \to S$ the $\scO_T$--submodules $\scN$ for which the quotient $h^* \scM / \scN$ is locally free of rank $n$.   

In case $S=\Spec(R)$, we write $\uGr_n(M)$ for $\uGr_n(\scM)$, and identify $\ulGr_n(\scM)$ with the $R$--functor $\ulGr_n(M)$ that assigns to $R' \in \Ralg$ the submodule $N$ of the $R'$--module $M \ot_R R'= M_{R'}$ for which $M_{R'}/N$ is locally free of rank $n$. \sm 

\item\label{grap-b}  Let $\ulGr\tot(\scM)$ be the $S$--functor which assign to a scheme $h \co X\to S$ the set of $\scO_X$--modules $\scN$ for which the quotient $h^*\scM/\scN$ is locally free of finite rank. It is represented by the $S$--scheme 
    \[ \uGr\tot(\scM) = \textstyle \bigsqcup_{n\in \NN} \, \uGr_n(\scM).\]

\inparcom{Reference [?] above is to the new paper by Philippe and Cameron, dito below in (c)}

For each $\scN \in \ulGr\tot(\scM)(X)$ we have the locally constant rank function, cf.~\eqref{lcf-1}, 
\[\rank_\scN \co X \to \NN_S(X), \quad x \mapsto \rank ( h^* \scM/\scN)_x  .
\] 
These rank functions fit together to give rise to a natural transformation $\ulGr\tot(\scM) \to \Mor(\cdot, \NN_S)$, hence to a morphism of $S$--schemes, 
\begin{equation}\label{grap-b1} 
\rank \co \uGr\tot(\scM) \to \NN_S, \quad \scN \mapsto \rank _\scN .
\end{equation}

\item \label{grape-nu} Let $\nu \in \NN_S$ be a locally constant function. We define 
\[ \uGr_\nu(\scM) \]
as the fibre of the rank function \eqref{grap-b1} over $\nu$. It is an smooth projective scheme. 

Over an affine scheme $\Spec(R)$, we identify the $R$--schemes $\uGr_\nu(\scM) = \uGr_\nu(M)$ and the $S$--functor $\ulGr_\nu(\scM)$ with the $R$--functor $\ulGr_\nu(M)$. The locally constant function $\nu = (\nu_{R'})_{R'\in \Ralg} \in \NN_{\Spec(R)}$ satisfies $\nu_{R'} = \nu \circ \Spec(R \to R')$, cf. ~\ref{lcf}. The $R'$--points, $R'\in \Ralg$, of the $R$--functor $\ulGr_\nu(M)$ are the $R'$--submodules $N\subset M_{R'}$ such that $M_{R'}/N$ is finite projective of rank $\nu_{R'} $.  One knows that $\ulGr_\nu(M)$ is represented by a smooth projective $R$--scheme. 
\sm 

\item \label{grape-dual} 
Taking annihilators of submodules of the dual $\scM\ch_T$, $T$ an $S$--scheme, of $\scM_T$, we can identify $\uGr_\nu (\scM\ch)(T)$ with the set of complemented submodules $\scN\subset \scM_T$ of rank $\nu$. 

In the affine case, we identify $\uGr_\nu (M\ch)(R')$, $R'\in \Ralg$, with the set of complemented submodules $N\subset M_{R'}$ of rank $\nu$. For example, we identify the projective space $\uP(M\ch)$ with the $R$--scheme whose $R'$--points are the complemented submodules $N\subset M_{R'}$ of rank $1$. 
\end{inparaenum} 

\newpage

\section{Some results on group schemes}

Unless specified otherwise, we consider group schemes over an arbitrary scheme $S$. 
\sm 

\subsection{Some known results}\label{skr} 
\begin{inparaenum}[(a)] \item \label{ag-ifp} ({\em Finite presentation of stabilizers})
Let 
$G$ be an $S$--group scheme acting on an $S$-scheme $X$. As in \cite[V, 10.2]{SGA3} let $H/X$ be the stabilizer of the diagonal $\De\co X \to X \times_S X$ for the action of $G \times_S X$ on $X \times_S X$, i.e., defined by the right cartesian square in \eqref{ag-i1} where $d_0$ is the action of $G$ on $X$ and $d_1$ is the projection onto the second factor. For any $x\in X(S)$ we define the
stabilizer $G_x$ of $x$ by the left cartesian square in \eqref{ag-i1}:
\begin{equation} \vcenter{\label{ag-i1}
 \xymatrix@C=40pt{ \ar @{} [dr] |{\Box} 
  G_x \ar[r] \ar[d] & H \ar[r] \ar[d]
  \ar@{} [dr] |{\Box}
  & G \times_S X \ar[d]^{(d_0, d_1)} \\
  S \ar[r]^x & X \ar[r]^\De & X \times_S X
}}\end{equation}
{\em If $G\to S$ and $X\to S$ are locally of finite presentation (of finite presentation respectively), the same holds for the morphisms $H\to X$ and $G_x\to S$ in the diagram \eqref{ag-i1}.}

Indeed, both properties allow base change and composition \cite[Tags 01TR, 01TS]{St}. Hence, the morphisms $G\times_S X \to S$ and $X\times_S S \to S$ are locally of finite presentation (of finite presentation respectively). But then so is $G\times_S X \to X \times_S X$ by cancellation \cite[Tag 02FV]{St}, and then $H \to X$ and $G_x \to S$ have these properties, again by base change. \sm 

\item \label{ag-f} ({\em Smoothness of group schemes over a reduced base})
Let $X$ be a scheme and let $H$ be an $X$--group scheme. Suppose
\end{inparaenum}
\begin{enumerate}[label={\rm (\Roman*)}]
\item \label{ag-fI} $X$ is reduced, e.g.\ $X=\Spec(R)$ for a reduced ring $R$ \cite[01J2]{St}, and 

\item \label{ag-fII} $H$ is locally of finite presentation over $X$, and
 the geometric fibers of $H$ are  smooth connected and have the same dimension.
\end{enumerate}
Then $H$ is smooth, separated and finitely presented (\cite[VI$_B$, 4.4 and 5.5]{SGA3}). See \ref{ag}\eqref{ag-red} 
for a criterion establishing the condition \ref{ag-fI}.%
\sm

\begin{inparaenum}[(a)] \setcounter{enumi}{2}
\item \label{ag-i} ({\em Smoothness of stabilizers}) In the setting of \eqref{ag-ifp}, i.e., $S$ is a scheme and $G$ is an $S$--group scheme acting on an $S$--scheme $X$,  assume that $X$ and the stabilizer $H$ satisfy the conditions \ref{ag-fI} and \ref{ag-fII} of \eqref{ag-f}.
Then $H$ is smooth, separated, and of finite presentation over $X$, and the analogous properties hold for the $S$--group scheme $G_x$ for any $x\in X(S)$.

Indeed, the properties of $H$ follow from the criterion \eqref{ag-f}, and they  imply the properties of $G_x$ mentioned above, since they are stable under base change.
\end{inparaenum}
\sm 

In \ref{gitqu}--\ref{thm_lee} we will consider GIT quotients.

\subsection{GIT quotients}\label{gitqu} 
Let $S$ be a scheme, let $G=\Spec(\calO_G)$ be an affine $S$--group scheme 
acting on the right on an affine $S$--scheme $X = \Spec(\calO_X)$.
We denote by $c\co \calO_X \to \calO_G \otimes_{\calO_S} \calO_X$ the coaction
and consider the quasi-coherent $\calO_S$--algebra
\[
 \calO_X^G= \bigl\{ f \in \calO_X : c(f)= 1 \ot f \}.
\]
The {\em GIT quotient\/} (GIT for ``Geometric Invariant Theory'')  is the
affine $S$-scheme $X\gq G= \Spec(\calO_X^G)$.
\sm

Let $S' \to S$ be a morphism of schemes and put $X'= X\times_S S'$ and $G'=G \times_S S'$. We then have a natural map $X'\gq G' \to X \gq G$, induced by the canonical maps $\calO_X^G \to \calO_X^G \ot_{\calO_S} \calO_{S'} \to \calO_{X'}^{G'}$, and hence a natural morphism $\al_{S'/S} \co X'\gq G' \to (X \gq G)\times_S S'$ making the diagram below commutative,
\begin{equation} \label{gitqu1} 
\vcenter{\xymatrix@C=40pt{
 X'\gq G' \ar[dr]^{\al_{S'/S}} \ar@/^1.2pc/[drr] \ar@/_1pc/[ddr]\\
   & (X\gq G)\times_S S' \ar[r] \ar[d] & S' \ar[d]\\
      & X\gq G \ar[r] & S
}}\quad . \end{equation}
The morphism $\al_{S'/S}$ is an isomorphism if $S' \to S$ is flat \cite[Lem.~2]{Seshadri}, in particular if $S'\to S$ is an open immersion, but not in general,  see for example the case of adjoint quotients of Lie algebras investigated in \cite{Bout-Cesna} and \cite{CR}.

We say that the GIT quotient $X\gq G$ is {\em universal\/} if its construction commutes with arbitrary base change, that is, $\al_{S'/S}$ is an isomorphism for any morphism $S' \to S$. Below we collect some results regarding GIT quotients and the morphisms $\al_{S'/S}$.  \sm

\begin{inparaenum}[(a)]
  \item \label{gitqu-a} ({\em Equivariant morphisms}) The construction of the GIT quotient respects equivariant group actions: if $\vphi_X \co X \to Y$ and $\vphi_G \co G \to H$ are morphisms of $S$--schemes and $S$--group schemes respectively satisfying $\vphi_X( x \cdot g) = \vphi_X(x) \cdot \vphi_G(g)$, then we obtain a natural induced morphism $X\gq G \to Y \gq H$, which is an isomorphism in both $\vphi_X$ and $\vphi_G$ are isomorphisms.    \sm

\item \label{gitqu-b} ({\em Transitivity}) Suppose $S'' \to S' \to S$ are scheme morphisms. For an $S$--scheme $Y$ we abbreviate $Y \times_S S' =: Y'$ and let
    \[
       a_Y \co  Y \times_S S'' \to (Y \times_S S') \times_{S'} S'' = Y' \times_{S'} S'' =: Y''
    \]
    be the canonical isomorphism. The isomorphism $a_X$ is equivariant with respect to the group scheme isomorphism $a_G$ and thus, by \eqref{gitqu-a}, gives rise to an isomorphism $\bar a$ of the corresponding GIT quotients. The diagram below is commutative:
    \begin{equation} \vcenter{
    \xymatrix@C=40pt{
    (X\times_S S'') \gq (G \times_S S'') \ar[d]_{\bar a}^{\cong}
        \ar[r]^{\qquad \al_{S''/S}} & (X\gq G) \times_S {S''} \ar[d]^{a_{X\gq G}}_\cong
    \\
    X'' \gq G'' \ar[d]^{\al_{S''/S'}}
    & (X\gq G)'\times_S S' \ar@{=}[d]
    \\
    (X' \gq G') \times_{S'} S'' \ar[r]_{\al_{S'/S} \times \Id_{S''}}
    & (X\gq G)'\times_S S'}
    }\quad . \end{equation}
  To see this, we can assume that $S$, $S'$ and $S''$ are affine, in which case commutativity is immediate.
\lv{
  In other words
    \begin{equation} \label{gitqu-b1}
        (\al_{S'/S} \times \Id_{S''}) \, \circ \, \al_{S''/S'}  \, \circ \, \bar a =
           a_{X\gq G}  \, \circ \, \al_{S''/S}
    \end{equation}
}
\sm

\item \label{gitqu-c} Let
\[  \xymatrix@C=40pt{ T' \ar[r] \ar[d] & T  \ar[d] \\ S' \ar[r] & S}  \]
be a commutative square, and let
\begin{align*}
 &\bar a_{S'} \co
   (X\times_S T') \gq (G \times_S T') \simlgr \big((X\times_S S') \times_{S'} T'\big)  \gq  \big((G\times_S S') \times_{S'} T')\big), 
   \\
   & \bar a_T \co (X\times_S T') \gq (G \times_S T') \simlgr \big((X\times_S T) \times_{T} T')\big)  \gq  \big((G\times_S T) \times_{T} T'\big)
 \end{align*}
be the canonical isomorphisms. Then
  \begin{equation} \label{gitqu-c1}
   (\al_{S'/S} \times \Id_{T'}) \, \circ \, \al_{T'/S'}  \, \circ \, \bar a_{S'} =
    (\al_{T/S} \times \Id_{T'}) \, \circ \, \al_{T'/T}  \, \circ \, \bar a_{T}.
 \end{equation}
Indeed, by \eqref{gitqu-b}, both sides equal $a_{X\gq G}  \, \circ \, \al_{T'/S}$.
\end{inparaenum}

\begin{lem}[Properties of universal GIT quotients] \label{lem_universal} In the setting of\/ {\rm \ref{gitqu}} the following hold.
\begin{enumerate}[label={\rm (\alph*)}]
 \item \label{lem_universal1} The universality property for $X\gq G$ is stable under arbitrary base change $S' \to S$. \sm

\item \label{lem_universal2} The universality property of $X\gq G$ is local for the fpqc topology, as for example defined in {\rm \cite[Tag 022B]{St}}. \sm

 \item \label{lem_universal3} The following are equivalent:

  \begin{enumerate}[label= {\rm (\roman*)}]

    \item\label{lem_universal3i} $X\gq G$ is universal;

    \item \label{lem_universal3ii}  For each morphism $\Spec(R) \to S$ and each $R$--algebra $R'$, the morphism $\al_{R'/R} \co X_{R'}\gq G_{R'}  \to (X_R \gq G_R) \times_{\Spec(R)} \Spec(R')$ is an isomorphism.
\end{enumerate}\end{enumerate}
\end{lem}

\begin{proof} \ref{lem_universal1} Assume that $X\gq G$ is universal and let $S'\to S$ be a morphism of schemes. To show that $(X\times_S S') \gq (G \times_S S')$ is universal, let $S'' \to S'$ be a second scheme morphism. By \ref{gitqu}\eqref{gitqu-b},
\begin{equation} \label{gitqu-b1}
       (\al_{S'/S} \times \Id_{S''}) \, \circ \, \al_{S''/S'}  \, \circ \, \bar a =
           a_{X\gq G}  \, \circ \, \al_{S''/S}
    \end{equation}
where $\bar a$ and $a_{X\gq G}$ are isomorphisms by construction, and $\al_{S'/S}$ and $\al_{S''/S}$ are isomorphisms by universality of $X\gq G$. Since then $ \al_{S'/S} \times \Id_{S''}$ is an isomorphism, $\al_{S''/S'}$ is an isomorphism too, establishing our claim. \sm

\ref{lem_universal2}  Let $(S_i)_{i \in I}$ be an fpqc cover of $S$,  put $X_i =X \times_S S_i$ and $G_i =G \times_S S_i$ for each $i \in I$, and assume that the GIT quotients $X_i\gq G_i$ are universal. To show that then $X\gq G$ is universal too, let $T$ be an $S$--scheme. We know that $(T_i = T\times_S S_i)_{i\in I}$ is an fpqc cover of $T$ \cite[Tag 022D]{St}. We specialize the formula \eqref{gitqu-c1} with $S'$ and $T'$ replaced by $S_i$ and $T_i$ respectively:
\[
   (\al_{S_i/S} \times \Id_{T_i}) \, \circ \, \al_{T_i/S_i}  \, \circ \, \bar a_{S_i} =
    (\al_{T/S} \times \Id_{T_i}) \, \circ \, \al_{T_i/T}  \, \circ \, \bar a_{T}.
\]
The morphism $T_i \to T$ is flat by definition of an fpqc cover, hence $\al_{T_i/T}$ is an isomorphism. The same argument applies to $S_i \to S$ and then implies that $\al_{S_i/S} \times_S \Id_{T_i}$ is an isomorphism too. Also, $\al_{T_i/S_i}$ is an isomorphism by assumption on the quotient $X_i \gq G_i$. Because the maps $\bar a_{S_i}$ and $\bar a_T$ are isomorphisms by construction, \[ \al_{T/S} \times_S \Id_{T_i} \co (X_T \gq G_T) \times_T T_i \simlgr
\big( (X\gq G) _T \big) \times_T T_i \]
is an isomorphism for each $i\in I$. Since the property of being an isomorphism is fpqc local \cite[Tag 02L4]{St},  it follows that $\al_{T/S}$ is an isomorphism, proving that $X\gq G$ is universal. \sm

\ref{lem_universal3} The implication \ref{lem_universal3i} $\implies$ \ref{lem_universal3ii} follows from \ref{lem_universal1}. For the other direction, we choose a Zariski cover of $S$ by affine schemes, which by \cite[Tag 022C]{St} is a fortiori an fpqc cover of $S$. Hence, to show that $X \gq G$ is universal, we can by \ref{lem_universal2} assume that $S=\Spec(A)$ is an affine scheme.

Let $T$ be an $S$--scheme and let $(T_i)_{i \in I}$ be a Zariski cover of $T$
by affine schemes $T_i=\Spec(B_i)$. The assumption \ref{lem_universal3ii} implies that each morphism   $\al_{B_i/A} \co X_{B_i} \gq G_{B_i}  \to (X \gq G) \times_{\Spec(A)} \Spec(B_i)$ is an isomorphism. In other words, the $T$--morphism $X_{T}\gq G_{T}  \to {(X\gq  G) \times_{\Spec(A)} T}$  becomes an isomorphism after  base change to the Zariski cover $(T_i)_{i \in I}$ of $T$. Applying again that a Zariski cover is an fpqc cover and  \cite[Tag 02L4]{St}, it follows that $X_{T} \gq G_{T}  \to (X \gq G) \times_{\Spec(A)} T$ is an isomorphism. Thus $X \gq G$ is universal.
\end{proof}

\
We can now slightly extend Lee's theorem \cite{Lee} on (Chevalley) adjoint quotients of reductive group schemes from an affine to an arbitrary base. See \ref{pare} for a review of reductive group schemes.

\begin{thm} \label{thm_lee}
Let $G$ be a reductive $S$--group scheme and consider its action on itself by inner automorphisms.

\begin{enumerate}[label={\rm (\alph*)}]
 \item \label{thm_lee1} The adjoint quotient $G\gq G$ is universal.

 \item \label{thm_lee2} Let $T$ be a maximal $S$--torus of $G$ and
 denote by $W=\rmN_G(T)/T$ its Weyl group. Then the map
 $T\gq W \to G\gq G$ is an isomorphism.
\end{enumerate}
\end{thm}

\begin{proof} \ref{thm_lee1} In the framework of affine schemes, this is has been proved in \cite[Thm.~4.1(1)]{Lee}, so that Lemma~\ref{lem_universal}\ref{lem_universal3}, \ref{lem_universal3ii} $\implies$ \ref{lem_universal3i} shows that $G\gq G$ is universal in general.
 \sm

\ref{thm_lee2} According to \cite[Tag 02L4]{St}, the statement
``$T\gq W \to G\gq G$ is an isomorphism'' is local for the fpqc topology.
We can therefore assume that $S$ is affine and that $G$ (resp.\ $T$) is a split reductive $S$--group (resp.\ a split $S$--torus). In this case, \ref{thm_lee2} has been established in \cite[Thm.~4.1(2)]{Lee}. \end{proof}

\comments{(2025-08-26) A short version of \ref{neutrevLG} is still here, commented out. The long version is needed for the section \ref{sec:quadrics}.  }
%

\subsection{Identity component of a group scheme} \label{neutrevLG} For the notion of the identity component $G^0$ of a group scheme $G$ over a scheme $S$, the reader is referred to \cite[VI$_B$, \S3]{SGA3}, see also \cite[II, \S5.1]{DG} for $S=\Spec(k)$, $k$ a field. Below we summarize what we will use. \sm 

\begin{inparaenum}[(a)] \item\label{neutrevLG-a}
Let $k$ be a field, and let $G$ be a $k$--group scheme locally of finite presentation, i.e., a locally algebraic $k$--group scheme in the terminology of \cite{DG}. We denote by $G^0$ the open subscheme whose underlying topological space is the connected component of the topological space $|G|$ containing $e\in G(k)$. One knows 
that $G^0$ is a characteristic $k$--subgroup scheme of $G$, and a closed and open $k$--subscheme of $G$. \sm

\item \label{neutrevLG-b} In the rest of this subsection, let $S$ be a general scheme, let $G$ be an $S$--group scheme locally of finite presentation, and let $\underline{G}^0$ be the $S$--functor defined in \cite[VI$_B$, \S3]{SGA3}. The following is known. 
\end{inparaenum} 

\begin{enumerate}[label={\rm (\roman*)}]
  \item \label{sec:neutral-bii} The functor $\underline{G}^0$ commutes with base change, \cite[VI$_B$, 3.3]{SGA3}.

  \item If $S$ is the spectrum of a field, then, by definition,  $\underline{G}^0$ is representable by the identity component  of the locally algebraic group $G$.

  \item \label{sec:neutral-biv} If $G$ has geometrically connected fibers, then, by definition,  
      $G$ represents $\underline{G}^0$.

  \item \label{sec:neutral-bv}  If $\underline{G}^0$ is representable by a scheme $G^0$, the canonical morphism $G^0 \to G$ is an open immersion, in particular $G^0$ is an open subgroup scheme of $G$ \cite[VI$_B$, 3.9]{SGA3}. Hence, by \eqref{opemaLG1}, for every affine scheme $\Spec(A)$ over $S$ we have 
\begin{equation}\label{somax1g} \begin{split}
    G^0(A) = \{ g\in G(A): \; & g_{A/\m}\in G^0(A/\m)\\
      &  \hbox{ for all maximal
                    $\m\in \Spec(A)$}\}. 
\end{split}\end{equation}

\item \label{sec:neutral-bvi}  If $G$ is smooth, then $\underline{G}^0$ is representable by an open smooth subgroup scheme of $G$, denoted $G^0$ \cite[VI$_B$, 3.10]{SGA3}.


  \item \label{sec:neutral-bw} If $G$ is a finitely presented smooth affine $S$--group scheme for which $G_{\bar s}^0$ is reductive for all $s\in \Spec(R)$, then $G^0$ is a reductive $S$--group scheme that is both open and closed in $G$, and in particular affine (\cite[Prop.~3.1.3]{Co1}). 
   \end{enumerate}

{\em Warning}: Even if $G$ is assumed to be smooth affine, it is in general not true that $G^0$ is a closed subscheme, nor that it is affine (\cite[XIX, 5.13]{SGA3}, \cite[VII, \S3, (iii)]{Raynaud}). \sm

{\em Example $G=\bmu_{2,R}$}: If $2\in R\ti$, then $G$ is smooth, hence $\underline{G}^0$ is representable. If $R$ is a field of characteristic $2$, then  $\underline{G}^0$ is not representable by a scheme.
\lv{
  Let $G = \bmu_2$, considered as an $R$--group scheme and let $G^0$ be the functor defined above. Over any field $\ka(q)$ we have $G_{\ka(q)} = \mu_2 (\ka(q))$, which therefore has 1 or 2 points. In the two cases, the 1-component is $ 1_{\ka(q)}$. Hence $G^0(S)$ consists of those $s\in S$ satisfying $s^2 = 1$ and which are $s=1$ after reduction to a field. Let $x$ be a non-zero element of the nil radical =  intersection of all prime ideals. Then $s \in G^0(S)$ but $s\ne = 1$, hence $G^0(S)$ is not the trivial functor. Hence $G^0$ is not representable by the trivial scheme.

On the other side, if $R$ is a ring with $2 \in R^\times$, then $G$ is smooth, hence representable. If $R$ is an algebraically closed field of characteristic $\ne 2$,  we can identify $G$ with its $G$--points, hence with $\{1, -1\}$ and $G^0 = {1}$ exists as scheme. }

\comments{(2022-01-10) Only \ref{lem_neutral0LG}\eqref{lem_neutral1} is used in the Knebusch paper as of today. }

\comments{(2025-08-26) I generalized \ref{lem_neutral0LG} from $S=\Spec(R)$ to arbitrary $S$. }

\begin{lem}\label{lem_neutral0LG} Let $G$ be an $S$--group scheme locally of finite presentation and assume that the $S$--functor $\underline{G}^0$ is representable
by an $S$--subgroup scheme $G^0$, necessarily open by {\rm \ref{neutrevLG}\ref{sec:neutral-bv}}, for example, by {\rm \ref{neutrevLG}\ref{sec:neutral-bvi}}, assume that $G$ is smooth.%
\sm

\begin{inparaenum}[\rm (a)]
\item \label{lem_neutral1} Let $f: H \to G$ be a morphism
of $S$--group schemes where $H$ is locally of finite presentation and has geometrically connected fibers. Then $f$ factors 
through $G^0$. \sm

\item \label{lem_neutral2}
The $S$--group scheme $G^0$ is the unique open subgroup scheme of $G$ such that for each  $s\in S$, 
$G^0 \times_S \kappa(s)$ is the identity component of the locally algebraic $\kappa(s)$-group $G_{\kappa(s)}= G \times_S \ka(s)$.
\end{inparaenum}
\end{lem}

\begin{proof}  \eqref{lem_neutral1} We start with the case that $S=\Spec(k)$, $k$  a field. Then $H$ is a connected locally algebraic group and $f: H \to G$
is a morphism of locally algebraic groups. In this case, $f^{-1}(G^0)$ is a closed and open $k$--subgroup scheme of $H$ by \ref{neutrevLG}\eqref{neutrevLG-a}.
By \ref{neutrevLG}\ref{sec:neutral-biv}, $\underline{H}^0$ is representable by $H$. It follows that $H = H^0\subseteq f^{-1}(G^0) \subseteq H$. Hence $f$ factorizes through $G^0$. The general case then follows readily from the definition of $\underline{G}^0$. \sm

\eqref{lem_neutral2} By  \ref{neutrevLG}\ref{sec:neutral-biv},
$(G^0)_s = (G_s)^0$ for each $s\in S$. Let $G'$ be an open $S$--subgroup scheme of $G$  such that $(G')_s= (G_s)^0$. Then  \eqref{lem_neutral1} 
shows that  the  open immersion  $i: G' \to G$ factors through $G^0$.  We thus get an open immersion $i' \co  G' \to G^0$ of $S$--schemes.  Since  $i'_s$ is an isomorphism for each $s\in S$, $i'$ is a surjective open immersion and therefore an isomorphism. 
\end{proof}

\subsection{Parabolic subgroups of reductive group schemes}\label{pare} 
For the convenience of the reader, we summarize here the results on parabolic subgroups of reductive groups over a scheme $S$ that we will use. We follow the terminology and notation of \cite{SGA3}, which coincides with the one used in \cite{Co1}. Although in the main part of the text we will only be  interested in reductive groups over an affine base, specializing here to the affine case does not lead to a great simplification. 
\sm

\begin{inparaenum}[(a)] \item \label{pare-aa} An $S$--group scheme $G$ is called {\em reductive} ({\em semisimple} respectively) if it is affine, smooth and all its geometric fibres $G_{\bar s}$, $s \in S$, are connected reductive groups (semisimple groups respectively) \cite[3.1.1]{Co1}, \cite[XIX, 2.7]{SGA3}. In the same vein, a {\em parabolic subgroup\/} of a reductive $R$--group $G$ is a  smooth subgroup $P$ of $G$ such that all its geometric fibres $P_{\bar s}$, $s \in S$, are parabolic subgroups of the reductive group $G_{\bar s}$, \cite[5.2.1]{Co1}, \cite[XXVI, 1.1 and 1.2]{SGA3}.
In the following, $P$ denotes a parabolic subgroup of a reductive group $G$. \sm

\item\label{pare-a} The canonical inclusion $P \to G$ is a closed immersion. The quotient sheaf $G/P$ in the \'etale topology on $S$--schemes is represented by a smooth projective $S$--scheme  \cite[5.2.3]{Co1}, \cite[XXII, 5.8.5]{SGA3}. \sm

\item\label{pare-urad} ({\em Unipotent radical}) The parabolic subgroup $P$ contains a unique smooth closed normal subgroup $\rad^u(P)$, the {\em unipotent radical of $P$}, whose geometric fibres $(\rad^u(P))_{\bar s}$, $s\in S$, coincide with the unipotent radical $\rad^u(P_{\bar s})$ of the parabolic subgroup $P_{\bar s}$ of the reductive group $G_{\bar s}$ \cite[5.2.5]{Co1}, \cite[XXII, 5.11.4]{SGA3}.  If $S=\Spec(R)$ is affine,  then, as $S$--scheme, $\rad^u(P)$ is isomorphic to $\uW(E)$, where $E$ is a finite projective $R$--module 
    \cite[XXVI, 2.5]{SGA3}. 
    \sm  

\item\label{pare-Levi} ({\em Levi subgroups}) A {\em Levi subgroup\/} is a smooth closed subgroup $L \subset P$ such that $L \ltimes \rad^u(P) \simlgr P$ \cite[5.4.2]{Co1}, \cite[XXVI, 1.7]{SGA3}. If $S$ is affine, every parabolic subgroup admits a Levi subgroup \cite[5.4.2]{Co1}, \cite[XXVI, 2.3]{SGA3}.  
\sm 

\item\label{pareopp} ({\em Opposite parabolic subgroup}) A parabolic subgroup $P'$ with the property that $P \cap P'$ is a Levi subgroup of $P$ and $P'$ is called {\em opposite to $P$}, \cite[XXVI, 4.3.2]{SGA3}. Let $P'$ be such an opposite parabolic subgroup, and let $U$ and $U'$ be the unipotent radicals of $P$ and $P'$ respectively. Then the canonical morphism $U' \to G/P$ is an open immersion, \cite[XXVI; Thm.~4.3.2(vi$^\prime$)]{SGA3}.   
    
    If $S$ is affine, every parabolic subgroup admits an opposite parabolic subgroup \cite[XXVI, 4.3.5(i)]{SGA3}. Moreover, if $S=\Spec(R)$, $R$ an LG ring, then $G(R)$ admits a decomposition 
 \begin{equation}\label{pareopp1}    
    G(R) = U(R)\, U'(R)\, P(R)
 \end{equation}   
    and if $P''$ is a second parabolic subgroup opposite to $P$,  there exist $u\in U(R)$ and $u'\in U'(R)$ such that $u'u$ conjugates $P'$ to $P''$, \cite[Thm.~4.1]{GN-LG}.  
    \sm 

\item \label{pare-b} ({\em The dynamic method}) 
Let $\la \co \GG_m \to G$ be a group morphism, sometimes called a one-parameter group or a cocharacter. It induces a $\GG_m$--action on $G$ by 
conjugation, and gives rise to a parabolic subgroup $\rmP_G(\la)$, see \cite[4.1.7, 5.2.2]{Co1} (or \cite[2.1]{CGP} where the condition in \eqref{pare-1} is abbreviated by ``$\lim_{t\to 0} \big(\la(t) h \la(t)\me\big)$ exists''), \cite[7.1.1]{G2}. Assume for simplicity that $S=\Spec(R)$ is affine. For $R'\in \Ralg$ and $t\in R'{}\ti$, we have $ \la(t) h \la(t)\me \in H(R'[t\me])$. Then $\rmP_G(\la$ is defined by requiring that all these elements  lie in $H(R')[t]$:  
\begin{equation}\label{pare-1}
\rmP_G(\la)\, (R') = \{ g\in G(R'): 
    \la(t) g \la(t)\me  \in G(R'[t]) \; \forall t\in R'{}\ti\},
\end{equation}
The centralizer $\Cent_G(\la)$ of $\la$ is a Levi subgroup of $\rmP_G(\la)$, \cite[7.1.1]{G2}. Moreover, the cocharacter $\la\me$ defines a parabolic subgroup opposite to $\rmP_G(\la)$ since the intersection $\rmP_G(\la) \cap \rmP_G(\la\me)$ is the centralizer of $\la$.   

\inparcom{(2025-08-11) I do not have a reference for the last claim (opposite parabolic).}

By \cite[7.3.2(1)]{G2}, if $S$ is affine, every parabolic subgroup of $G$ is of the form $\rmP_G(\la)$ for some group homomorphisms $\la \co \GG_m \to G$\sm

\item \label{pare-d} ({\em The scheme of parabolic subgroups}) The $S$--functor associating with $S'\to S$ the set of parabolic subgroups of the reductive $S'$--group $G\times_S S'$ is representable by a smooth projective $S$--scheme $\Par(G)$  \cite[XXVI, 3.5]{SGA3}, \cite[5.2.9]{Co1}. \sm 

\item\label{pare-Dyn} ({\em The Dynkin scheme and the type morphism}) The Dynkin scheme $\Dyn(G)$ of $G$ is defined in \cite[XXIV, 3.3]{SGA3}. We mainly use it in the case that $S=\Spec(R)$ is affine and $G$ has constant type.
    In this case, $G$ is a twisted form of a Chevalley $R$--group scheme $G_0$, say with associated Dynkin diagram $\De_0$, and $\Dyn(G)$ is an \'etale  twisted form of the constant scheme $\De_{0,R}$. 
\sm 

Following \cite[XXVI, 3.1]{SGA3}, $\Of\big(\Dyn(G)\big)$ denotes the scheme representing the $R$--functor that assigns to $R'\in \Ralg$ the set of open and closed subsets of $\Dyn(G)\times_R R' \cong \Dyn(G\times_R R')$. For example, if $\Dyn(G) = \De_R$, then $\Of\big(\Dyn(G)\big) = (\calP(\De))_R$, the constant scheme associated with the set $\calP(\De)$ of all subsets of the Dynkin diagram $\De$.

Finally, there exists a type morphism $\mathbf{t} \co \Par(G) \to \Of\big(\Dyn(G)\big)$, which in case of a split $G$ of constant type $\De$ associates to $P$ the root data needed to describe $P$ in terms of a pinning adapted to $P$ \cite[XXVI, 3.3]{SGA3}. Two parabolic subgroups $P$ and $P'$ have the same type in the sense that $\mathbf{t}(P) = \mathbf{t}(P')$ if and only if $P$ and $P'$ are conjugate (fpqc)--locally. For a given type $t\in \Of\big(\Dyn(G)\big)$ we use $\Par(G)_t$ to denote the fibre over $t$. If $P$ has type $t$, then
\begin{equation}  \label{pare-d2}
G/P \simlgr \Par(G)_t
\end{equation}
by \cite[XXVI, 3.6]{SGA3}. 
\sm 

\item\label{pare-cqu} ({\em Central quotients} \cite[Lem.~3.2.1]{G2}) Let $f \co G'\to G$ be a surjective homomorphism of reductive $S$--group schemes whose kernel is a central $S$--group of multiplicative type. Then the inverse image $P \mapsto f\me(P)$ defines an isomorphism $\Par(G) \simlgr \Par(G')$ of $S$--schemes, which preserves the types of parabolic subgroups. 
\end{inparaenum}

\subsection{Transitivity on the big affine Zariski site.}\label{transi} Let $S$ be a scheme and let $G$ be an $S$--group scheme acting on an $S$--scheme $X$. We say that {\em the action of $G$ on $X$ is transitive on the big affine Zariski site of\/ $S$\/} if the canonical map 
\[ G \times_S X \to X \times_S X, \quad (g, x) \mapsto (g\cdot x, x) 
\]
is an epimorphism of sheaves on the big affine Zariski site of $S$. Recall that this means that for every affine scheme $U=\Spec(A)$ over $S$ and for every pair $x,y \in X(U)$ there exists a standard Zariski cover $(a_1, \ldots, a_n)$ of $A$ such that for every $i$, $1\le i \le n$, there exists $g_i \in G(\Spec(A_{a_i}))$ satisfying $g_i \cdot x_i = y_i$ where $x_i = x \circ \big(\Spec(A_{a_i}) \to \Spec(A)\big)$ and where $y_i$ is defined analogously.

\subsection{Transitivity on the small affine Zariski site.}\label{transm} 
Let $S$ be a scheme and let $G$ be an $S$--group scheme acting on an $S$--scheme $X$. Replacing the big affine Zariski site of $S$ by the small affine Zariski site in \ref{transi}, we get the concept of a {\em transitive action of $G$ on $X$ on the small affine Zariski site of\/ $S$\/}: for every open immersion  $U=\Spec(A) \to S$ of an affine scheme $U$ and for every pair $x,y \in X(U)$ there exists a standard Zariski cover $(a_1, \ldots, a_n)$ of $A$ such that for every $i$, $1\le i \le n$, there exists $g_i \in G(\Spec(A_{a_i}))$ satisfying $g_i \cdot x_i = y_i$ with $x_i$ and $y_i$ as in \ref{transi}. 

Example: In the setting above suppose that for all $A$ and all $x,y \in X(A)$, the $A$--functor 
\[ T_{x,y}(B) = \{ g\in G(B): g \cdot x_B = y_B\}, \quad (B\in \Aalg) \]
is represented by an $A$--scheme locally of finite presentation, and that $T(A_\m) \ne \emptyset$ for all maximal ideals $\m \ideal A$. Then Corollary~\ref{fibn} 
implies that the action of $G$ on $X$ is transitive on the small affine Zariski site of $\Spec(R)$. 

\comments{(2025-08-07)  We use this example in our version of 
``Demazure's Conjugacy Theorem'' \ref{thm_conj_demazure} }

\ms

We will use Demazure's Conjugacy Theorem in the form stated in \cite[A.1, A.2]{GN-Sp} and partially generalized in \cite{GN-LG}. 


\begin{thm}[\bf Demazure's Conjugacy Theorem  {\cite[XXVI, \S3]{SGA3}}] \label{thm_conj_demazure} Let $G$ be a reductive $R$--group scheme, and let $X=\Par(G)_t$ be the $R$--scheme of parabolic subgroups of type $t$. \sm

\begin{inparaenum}[\rm (a)]
\item \label{thm_conj_demazure-a} Given parabolic subgroups $x$, $y$ in $X(R)$, there exist $f_1,\ldots, f_n \in R$ satisfying $f_1+\cdots+ f_n= 1$ and such that $y_{R_{f_i}} \in G(R_{f_i})\,.\, x_{R_{f_i}}$ for $i=1 \ldots n$. In other words, the action of $G$ on $X$ is transitive on the small affine Zariski site of $\Spec(R)$. \sm 
    
\item \label{thm_conj_demazure-b} If $R$ is an LG ring, then $G(R)$ acts transitively on $X(R)$. \sm 

\item \label{thm_conj_demazure-bii} Assume that $R$ is semilocal and that $R'$ is a finite $R$-algebra such that $X(R') \not = \emptyset$. Let $\Jac(R)$ be the Jacobson radical of $R$. Then the map
\[ X(R') \longto X(R'/\Jac(R)R') \cong \textstyle \prod_{\m \in \Spec(R) \text{ maximal}}\; X(R'/ \m R') \]
is onto.
\end{inparaenum}
\end{thm}

\begin{proof} We first prove (\ref{thm_conj_demazure-b}).
If $X(R)=\emptyset$, the statement is obvious. We can thus assume that $X(R) \ne \emptyset$ and pick a point $x \in X(R)$, i.e., a parabolic subgroup $P$ of $G$ of type $t$, and a parabolic subgroup opposite to $P$, \ref{pare}\eqref{pareopp}. Then the claim follows from \cite[Thm.~4.1(a)]{GN-LG} and \eqref{pare-d2}. \sm 

\eqref{thm_conj_demazure-bii} Recall that  $R'$ is semilocal by \ref{slr}\eqref{slr-b}. Our assumption is that $G_{R'}$ admits a parabolic subgroup $Q$ of type $t$. As noted in the proof of \eqref{thm_conj_demazure-b}, it then also admits an opposite  parabolic subgroup $Q'$. According to  \cite[XXVI, 5.2]{SGA3}, the product map $\rad_u(Q)(R')  \times \rad_u(Q')(R') \to X(R')$ is surjective (here $\rad_u(\cdot)$ denotes the unipotent radical). Applying this for the semilocal ring $R'$ as well as
for the semilocal ring $R'/\m R'$, $\m \in \Spec(R)$ maximal, shows that the horizontal maps in the commutative diagram below are surjective:
\[ \xymatrix{
    \rad_u(Q)(R')  \times \rad_u(Q')(R')  \ar[r] \ar[d] &  X(R') \ar[d] \\
     \prod\limits_{\m} \rad_u(Q)(R'/\m R')  \times \rad_u(Q')(R'/\m R')
     \ar[r] &  \prod\limits_{\m} X(R'/\m R')
} \]
Since $\prod_\m R'/\m R' \cong \prod_\m R' \ot_R (R/\m) \cong R' \ot_R (R/\Jac(R) \cong R'/\Jac(R)R'$, 
the map $R' \to \prod_{\m} R'/\m R'$ is onto. On the other hand, the $R'$--scheme $\rad_u(Q)$ (respectively $\rad_u(Q')$) is isomorphic to a vector group $R'$--scheme by \ref{pare}\eqref{pare-urad}, 
so that  the left vertical map is onto. Hence, by a simple diagram chase the right vertical map  is onto too.
\sm

\eqref{thm_conj_demazure-a} The group $G$ is an $R$--group of type (RR) by \cite[XXII, 5.1.3]{SGA3} and the points $x$, $y$ of $X(R)$ are parabolic subgroups, hence subgroups of type (R) by 5.2.3 of {\em loc.\ cit.}. It then follows from Theorem~5.3.9 of {\em loc.\ cit.} that the strict transporter $T$ of $x$ and $y$, defined by
\[ T(B) = \{ g\in G(B) : g \cdot x_{B} = y_{B}\} \quad (B\in \Ralg), \]
is a finitely presented $R$--scheme (among other properties). Since $T(R_\m) \ne \emptyset$ for any maximal $\m \in \Spec(R)$ by \eqref{thm_conj_demazure-b}, the claim follows from the Example in \ref{transm}. 
\end{proof}

%



\begin{lem} \label{lem_onto} Let $R$ be an LG ring, and let  $G$ be a reductive $R$-group scheme.%
\sm 

\begin{inparaenum}[\rm (a)] \item \label{lem_onto-a} Let $P$ be a parabolic $R$-group scheme equipped with a Levi $R$-subgroup scheme $L$, and let $\fra$ be an ideal of $R$. Then the following are equivalent:
\end{inparaenum}

\begin{enumerate}[label={\rm (\roman*)}]

\item \label{lem_ontoi} $L(R) \to L(R/\fra)$ is onto; \sm

\item\label{lem_ontoii} $P(R) \to P(R/\fra)$ is onto; \sm

\item \label{lem_ontoiii}  $G(R) \to G(R/\fra)$ is onto.
\end{enumerate} 
\sm 

\begin{inparaenum}[\rm (a)] \setcounter{enumi}{1}
\item \label{lem_onto-b}  Let $J=\Jac(R)$ be the Jacobson radical of $R$, let $k_1, \ldots, k_n$ be fields such that $R/J = k_1 \times \cdots \times k_n$, thus $R$ is semilocal, and let $(P_i)_{i \in I}$ be a family of parabolic subgroups of $G$ such that
\begin{equation}  \label{lem_onto-b1}
 G(k_j) = \big\lan \rad^u(P_i)(k_j): i \in I \big\ran 
\end{equation} 
for $j=1, \ldots, n$.  Then $G(R) \to G(R/J)$ is onto. 
\end{inparaenum}
\end{lem}

\begin{proof} \eqref{lem_onto-a} 
We will prove \ref{lem_ontoi} $\iff$ \ref{lem_ontoii} for any $R$ and any ideal $\fra \ideal R$. We put $\ol R = R /\fra$. 
\sm  

\ref{lem_ontoi} $\iff$ \ref{lem_ontoii}: Let $U=\rad^u(P)$ be the $R$-unipotent radical of $P$. 
We have $P=U \rtimes L$, by definition of a Levi subgroup, \ref{pare}\eqref{pare-Levi}. As $R$-scheme, $U$ is isomorphic to $\uW(E)$ for a finite projective $R$--module $E$, 
\ref{pare}\eqref{pare-urad},  
so that $U(R) \to U(\ol R)$ is obviously onto. It follows that  \ref{lem_ontoi} $\iff$ \ref{lem_ontoii}. 
\sm

\ref{lem_ontoiii} $\implies$ \ref{lem_ontoii}:  
Let $P'$ be a parabolic $R$-subgroup scheme opposite to $P$, and let $U'$ be its unipotent radical; thus $P'=U' \rtimes L$ by \ref{pare}\eqref{pare-Levi} again. We are given $p_0 \in P(\ol R)$. By assumption we have  $p_0= \underline{g}$ for some $g \in G(R)$, where $\underline{g}$ denotes the canonical image of $g$ in $P(\ol R)$. According to \eqref{pareopp1}, 
there exist $u\in U(R)$ and $u' \in U'(R)$ such that $g= u \, u' \, p$. Replacing $p_0$ by $\underline{p}^{-1} p_0$, we are reduced to the case $g =u \, u'$. It follows that   $\underline{u} \, \underline{u'}= p_0$, so that  $\underline{u}=p_0$ by \ref{pare}\eqref{pareopp}. 
Since $U(R) \to U(\ol R)$ is onto, we conclude that $p_0$ belongs to the image of $P(R) \to P(\ol R)$. \sm 

\ref{lem_ontoii} $\implies$ \ref{lem_ontoiii}: By \ref{revLG}\eqref{revLG-ideals}, the quotient $\ol R$ is also an LG ring. We thus have the decomposition \eqref{pareopp1} for $R$ and $\ol R$. We have seen in the proof of \ref{lem_ontoiii} $\implies$ \ref{lem_ontoii} that $U(R) \to U(\ol R$ and $U'(R) \to U'(\ol R)$ are onto, and we know that $P(R) \to P(\ol R)$ is onto. Thus $G(R) = U(R) U'(R) P(R) \to G(\ol R) = U(\ol R) U'(\ol R) P(\ol R)$ is onto too. 
\sm 

\eqref{lem_onto-b} We observe that $G(\ol R) = G(k_1) \times \cdots \times G(k_n)$. For $i\in I$ let $U_i = \rad^u(P_i)$. It follows that the assumption \eqref{lem_onto-b1} is equivalent to $G(\ol R) =  \big\lan U_i(\ol R) : i\in I \big\ran$. We have seen in the proof of \eqref{lem_onto-a} that the canonical homomorphism $U_i(R) \to U_i(\ol R)$ is onto. Hence $\big\lan U_i(R): i\in I\big\ran \to G(\ol R)$ is onto, in particular \eqref{lem_onto-b} holds. \end{proof}

\subsection{Reductive group schemes: isotropic versus reducible}\label{rgs}
Let $G$ be a reductive group scheme over a scheme $S$. Following \cite{G2}, see also \cite[5.1]{GN-LG}, we call $G$ {\em isotropic\/}, if $G$ admits a subgroup isomorphic to $\GG_{m,S}$. We say that $G$ is {\em reducible\/} if $G$ admits an everywhere proper parabolic subgroup $P$, i.e., $P_{\bar s} \subsetneq G_{\bar s}$ for all $s\in S$ and $P$ admits a Levi subgroup. We exhibit some known results: \sm 

\begin{inparaenum}[(a)]
  \item \label{rgs-a} (\cite[Thm.~7.3.1(2)]{G2}) If $S$ is connected \footnote{ The assumption that $S$ be connected is missing in \cite[5.1(a)]{GN-LG}}, a reductive $S$--group scheme $G$ is isotropic if and only if $G$ is reducible or the central torus $\rad(G)$ is isotropic. In particular, a semisimple $S$--group scheme is isotropic if and only if it is reducible. \sm 
      
 \item \label{rgs-aa} (\cite[Cor.~7.3.2(2)]{G2}) If $S$ is an affine scheme, a reductive $S$--group scheme $G$ is isotropic if and only if $G$ admits an  
     everywhere proper parabolic subgroup or the central torus $\rad(G)$ is isotropic.            
\end{inparaenum}

\subsection{Characters and Weil restriction}\label{cwr} 
\begin{inparaenum}[(a)] 
  \item\label{cwr-a} ({\em Definition}, see e.g. \cite[\S4]{Oesterle}) Let $S=\Spec(R)$ and let $G$ be an $R$--group scheme. We denote by $\GG_{m,R}$ the multiplicative group over $R$. A {\em character\/} is a homomorphism $\chi \co G \to \GG_{m,R}$ of $R$-group schemes, sometimes also called a {\em rational character\/}. The set of characters is in an obvious way a commutative group, denoted $\wdh G (R)$. 
%
    \sm 
      
 \item\label{cwr-b}  ({\em The homomorphism $\euN_{R'/R}\co \wdh G'(R') \to \wdh G(R)$}, \cite[II, 2.2]{Oes}) Let $R' \in \Ralg$ be faithfully projective as $R$--module, and let $G'$ be an affine $R'$--group scheme, hence the Weil restriction $G:= \frR_{R'/R}(G')$ exists as an $R$--scheme and is in fact an $R$--group scheme, \ref{weilres}. Using functoriality of the Weil restriction, a character $\chi' \co G' \to \GG_{m,R'}$ gives rise to an $R$--group homomorphism 
     $\frR_{R'/R}(\chi') \co \frR_{R'/R}(G') \to \frR_{R'/R}(\GG_{m,R'})$, which we can compose with the $R$--group homomorphism norm $\rmN_{R'/R} \co \frR_{R'/R}(\GG_{m,R'}) \to \GG_{m,R}$ of \eqref{trno-ds2} to get a character 
     \[ \euN_{R'/R}(\chi') := (\rmN_{R'/R} \circ \frR_{R'/R})(\chi') \co G \to \GG_{m,R} \]       
     of $G$. By multiplicativity of $\rmN_{R'/R}$ and $\frR_{R'/R}$, the map 
     \begin{equation}\label{cwr-b1}
       \euN_{R'/R} \co \wdh {{G'}}(R') \to \cE{\frR_{R'/R}(G')}(R), \quad \chi' \to \euN_{R'/R}(\chi')
     \end{equation} 
     is a group homomorphism. \sm

\new
\item \label{cwr-ex} ({\em The example $G'=\GG_{m,R'}$}) We specialize \eqref{cwr-b} for $G' = \GG_{m,R'}$. In this case, the $R'$--functor $A'\mapsto \wdh{G'}(A')$ is representable by the constant $R'$--group scheme associated with $\ZZ$. In particular, 
    $\wdh{G'}(R')$ can be identified with the abelian group of locally constant functions $\Spec(R') \to \ZZ$. 
    
    {\em Suppose now that $R'$ is connected.} Hence $\wdh{G'}(R') = \ZZ \cdot \Id_{R'}$ where we view $\Id_{R'}$ as the trivial character. Since $\frR_{R'/R}(\Id_{R'}) = \Id_R$, it follows that  $\euN_{R'/R}(\Id_{R'}) = \rmN_{R'/R}$, viewed as group scheme homomorphism as in \eqref{trno-ds2}, and thus
\begin{equation}\label{cwr-ex1} \begin{split}
  \euN_{R'/R} \co \wdh{\GG_{m,R'}} = \ZZ \cdot \Id_{R'} &\to \cE{\frR_{R'/R}(\GG_{m,R'})},  \\ n \cdot \Id_{R'} &\mapsto n \cdot \rmN_{R'/R}. 
\end{split}\end{equation} 
\enew
     
\item\label{cwr-c} ({\em The separable field case} \cite[II, Thm.~2.4]{Oes}) Let $k'/k$ be a finite field extension, and let $G'$ be a smooth affine $k'$--group scheme. Then the group homomorphism $\euN_{k'/k}$ of \eqref{cwr-b1} is a bijection in the following two cases  \sm 
    
    \begin{inparaenum}[(i)] \item \label{cwr-ci} $k'/k$ is separable, or
    
    \item\label{cwr-cii} $k'/k$ is purely inseparable of height $1$. 
     \end{inparaenum}      

\noindent In Theorem~\ref{thm_Lourdeaux} we will describe what happens for an arbitrary finite extension $k'/k$. \sm

\item Let again $k'/k$ be a finite field extension and let $G$ be an affine $k$--group. The canonical homomorphism
  \[ G \to \frR_{k'/k}(G_{k'}) =: G_0 \]  
   is a closed immersion \cite[A.5.7]{CGP}, and thus, by restriction, gives rise to a group homomorphism
   \[ \res_{|G} \co \wdh {G_0}(k) \to \wdh G(k), \quad \chi \mapsto \chi|_G. \] 
   
We will construct a group homomorphism in the other direction. First, given $\chi \in \wdh G(k)$, we have the group homomorphism $\chi_{k'} \co G_{k'} \to \GG_{m,k'}$, obtained by base change, and then, by applying the Weil restriction functor, the group homomorphism
\[ \frR_{k'/k}(\chi_{k'}) \co G_0 \to \frR_{k'/k}(\GG_{m,k'}).\]

\item \label{cwr-d} Suppose now that $k$ has characteristic $p>0$ and that $k'/k$ is a purely inseparable extension of height $h$, that is, $h\in \NN$ is the smallest integer such that $(k'){}^{p^h} \subset k$. In this case, we have a character 
\[ v_{k'/k} \co \frR_{k'/k}(\GG_{m,k'}) \to \GG_{m, k} \]  
given on $T\in \kalg$ by $(k'\ot_k T)\ti \to T\ti$, $x \mapsto x{}^{p^h}$.   
As before, let $G$ be a $k$--group scheme. Composing the two homomorphisms above yields a character of $G_0$, 
\[ v_{k'/k}  \circ \frR_{k'/k}(\chi_{k'}) \co G_0 \to \GG_{m,k} \]
Clearly,
\[ \tau \co \wdh G(k) \to \wdh{ G_0}(k), \quad \chi \mapsto 
v_{k'/k} \circ \frR_{k'/k}(\chi_{k'}) \]
is a group homomorphism. The composition 
$\res_{|G} \circ \tau \co \wdh G(k) \to \wdh G(k)$ 
is an endomorphism of the abelian group $\wdh G(k)$ whose image is contained in  
\begin{equation}\label{crw-d1} p^h \cdot \wdh G(k) := \{ \chi{}^{p^h} : \chi \in \wdh G(k) \}
\end{equation}  
It is shown in \cite[Thm.~A.4]{Lourdeaux} that {\em the restriction map $\res_{|G} \co \wdh{G_0}(k) \to \wdh G(k)$  is injective and identifies $\wdh{G_0}(k)$ with the subgroup $p^h \cdot \wdh G(k)$ if} 

\begin{inparaenum}[\rm (I)]
 \item\label{crw-dI} {\em $G$ is diagonalizable, or if
 
 \item $G$ is smooth and of multiplicative type.  }
\end{inparaenum}
\sm 

\item\label{cwr-g} ({\em The example $G=\GG_{m,k}$ in the inseparable case}) Let again $k'/k$ be a purely inseparable field extension of height $h$ and let $G=\GG_{m,k}$. Thus, as in \eqref{cwr-ex}, $\wdh G = \ZZ \cdot \Id_k$, and \eqref{crw-dI} applies. It follows that {\em  
    \[  
      \wdh {G_0}(k) = \Hom_{k-gp}(\frR_{k'/k}(\GG_{m.k'}), \GG_{m,k})
      \]  
      is free of rank $1$, generated by the homogeneous character $v_{k'/k}$ of degree $p^h$.} 
\end{inparaenum}

\lv{
\begin{cor}\label{corLourd} Let $k$ be a field of characteristic $p>0$, let $L/k$ be a finite field extension, let $K$ be the separable closure of $k$ in $L$ and let $h$ be the height of the purely inseparable extension $L/K$. Finally, let $H$ be a smooth $K$--group of multiplicative type. Then the composition of the group homomorphisms \eqref{crw-d1} and \eqref{cwr-b1}
\begin{equation} \label{corLourd1}  
 \cE{\frR_{L/K}(H_L)} (K) \xrightarrow{\; \res_{|H}\;}  \wdh H (K) \xrightarrow{\euN_{K/k}} \cE{\frR_{K/k}(H)}(k) 
\end{equation}
is injective and has image $\euN_{K/k} \big( p^h \cdot \wdh H(K)\big)$. In particular, for $H=\GG_{m,K}$ we get the following. \sm 

\begin{enumerate}[label={\rm (\roman*)}] 
  \item\label{corLourdi}
  $\cE{\frR_{L/K}(\GG_{m,L})}$ is free of rank $1$ with basis the character $v_{L/K}$. \sm
  
  \item\label{corLourdii} The closed immersion $\GG_{m,K} \to \frR_{L/K}(\GG_{m,L})$ induces an injective group homomorphism by restriction, 
      \[ \res_{|\GG_{m,K}} \co \cE{\frR_{L/K}(\GG_{m,L})}(K) \longto \cE{\GG_{m,K}}(K)  
      \]
      whose image is $\ZZ \cdot v_{L/K}|_{\GG_{m,K}}$.\sm
      
   \item \label{corLourdiii} The map 
     \[ \wdh{\GG_{m,K}}(K) \to \cE{\frR_{K/k}(\GG_{m,K})}(k), \quad 
        n \cdot \Id_K \mapsto n \cdot \rmN_{K/k} 
     \]
    is an isomorphism of abelian groups. 
\end{enumerate}
\end{cor}

\begin{proof}
  By \ref{cwr}\eqref{cwr-d}, the restriction map $\res_{|H}$ is an injective group homomorphism with image $p^h \cdot \wdh H(K)$, and by \ref{cwr}\eqref{cwr-c}  the group homomorphism $\euN_{K/k}$ is an isomorphism. This proves \eqref{corLourd1}. The examples \ref{cwr}\eqref{cwr-ex} and \ref{cwr}\eqref{cwr-g} yield  \ref{corLourdi}--\ref{corLourdiii}.  
\end{proof} 
}

\comments{(2025-06-25) It would be of interest to generalize \ref{thm_Lourdeaux} to the setting where $\GG_{m,K}$ is replaced by a smooth $K$--group of multiplicative type.  }

\begin{thm}\label{thm_Lourdeaux} 
 Let $k$ be a field of characteristic $p>0$, let $L/k$ be a finite field extension, let $K$ be the separable closure of $k$ in $L$ and let $h$ be the height of the purely inseparable extension $L/K$.
Then the composition of $k$--group homomorphisms
\[ 
\frR_{L/k}(\GG_{m,L}) 
\xrightarrow{\frR_{K/k}(v_{L/K})} \frR_{K/k}(\GG_{m,K})
\xrightarrow{\rmN_{K/k}} \GG_{m,k}
\]
induces an injective homomorphisms of character groups
\[ \rho \co  \cE{\frR_{L/k}(\GG_{m,L})}(k) \longto \cE{\GG_{m,k}}(k) = \ZZ \cdot \Id_k 
\]
whose image is 
\begin{equation} \label{thm_Lourdeaux2}
\Ima \rho = 
\ZZ \cdot \big( \rmN_{K/k} \circ \frR_{K/k}(v_{L/K})\big).
\end{equation}
\end{thm}

%
%

\begin{proof}  
We first set up some notation. Let $\widetilde K$ be  a Galois closure of $K/k$.
Hence $K\subset \wtl K$ and both $\wtl K/k$ and $\wtl K /K$ are Galois extension. We put $G=\Gal(\widetilde K/k)$ and $H=\Gal(\widetilde K/K)$, so that $K=\wtl  K^H$.
Let $\al$ be a primitive element of the separable extension $K/k$ and 
let $P$ be its minimal polynomial. Then $K=k[X]/P(X)$ where $P(X)= \prod_{\sigma \in G/H}(X-\sigma(\alpha))$. 
We then have an isomorphism 
\begin{equation}\label{thm_Lourdeaux22} 
K \otimes_k \widetilde K = \widetilde K[X]/P(X)\simlgr (\widetilde K)^{(G/H)}= \textstyle \prod_{\sigma \in G/H} \, \widetilde K_\si
\end{equation} 
which maps each $Q(X) \in \widetilde K[X]$ to $\bigl( Q(\sigma(\alpha)) \bigr)_{\sigma \in G/H}$ (of course, $\wtl K_\si = \wtl K$ as fields).
The natural action of $G$ on $K \otimes_k \widetilde K$ becomes 
on the right-hand side of \eqref{thm_Lourdeaux22} 
\[ 
 \ga \cdot \bigl( x_{\sigma}\bigr)_{\sigma \in G/H} 
 = \bigl( \ga( x_{\ga^{-1}\sigma})\bigr)_{\sigma \in G/H} \, .
\]
In particular,  the action of $G$ permutes the factors of $(\widetilde K)^{(G/H)}$. 

Next, we have  a compatible decomposition  
\begin{equation}\label{thm_Lourdeaux3}
 L \otimes_k \widetilde K 
  \simlgr \textstyle   \prod_{\sigma \in G/H} \widetilde L_\sigma
\end{equation} 
 where $\wtl L_\si = L \ot_k \wtl K_\si$, so that $\widetilde L_\sigma/ \widetilde K_\si$ is a purely inseparable extension of
height $h$. 
\quest{(2026-06-25) Why is $\wtl L_\si = L \ot_k \wtl K_\si = L \ot_k K$ a field? It contains $K \ot_k K$ which in general has $0$--divisors. Rather, it should be $L\ot_K \wtl K$ which is indeed a field, and a purely inseparable extension of $\wtl K$. On the other hand, changing $\ot_k$ to $\ot_K$ leads to problems below with the definition of $\chi_{\wtl K}$ because there we take $\ot_k$. 

It is also no longer clear to me why ew have \eqref{thm_Lourdeaux3}.

I am tired, and hope that I am mistaken.} 
Again the action of $G$ on $L \otimes_k \widetilde K$
permutes the $L_\sigma$'s. To summarize, we have the following diagram of algebras: 
\[\xymatrix@R=20pt{ & L\ar[rr] && L \ot_k \wtl K  = \prod_{\si \in G/H} \, \wtl L_\si \\ 
   && K \ot_k \wtl K = \prod_{\si \in G/H} \, \wtl K_\si \ar[ur]\\
k \ar[r] & K\ar[uu] \ar[ur] \ar[rr] && \wtl K \ar[ul] \ar[uu]
} \]

We are given a $k$-homomorphism $\chi \co \frR_{L/k}(\GG_{m,L}) \to \GG_{m,k}$. The extension of $\chi$ to $\wtl K$ becomes  
\[ \chi_{\wtl K} \co 
\prod_{\sigma \in G/H} \frR_{\widetilde L_\sigma/\widetilde K}(\GG_{m, \wtl L_\si}) \to \GG_{m,\widetilde K}
\] 
since $\frR_{L/k}(\GG_{m,L}) \ot_k \wtl K \cong \frR_{L \ot_k \wtl K/\wtl K}(\GG_{m, L\ot_k \wtl K})$ by \ref{weilres}\eqref{weilres-bc}. Using the purely inseparable case of the theorem, \ref{cwr}\eqref{cwr-g}, 
the component  $\frR_{\widetilde L_\sigma/\widetilde K}(\GG_{m,\wtl L_\si}) \to \GG_{m,\widetilde K}$ is $n_\sigma \, v_{\wtl L_\si /\wtl K}$ for a unique $n_\sigma \in \ZZ$. The action of $G$ on $L \otimes_k \widetilde K $
shows that all $n_\sigma$'s are equal to some $n \in \ZZ$.
In other words, $$
\chi_{\widetilde K}\bigl( (x_\sigma)\bigr)=
\prod\limits_{\sigma \in G/H} (x_\sigma^n)^{p^h}
= \rmN_{L \otimes_k \widetilde K/\widetilde K} \bigl( (x_\sigma^n)^{p^h}\bigr).
$$
It follows that  $\chi$ and $n \, \rmN_{K/k}\circ \frR_{K/k}(v_{L/K}$ agree
after extension to $\widetilde K$. By Galois descent 
we conclude that $f= n \, \rmN_{K/k}\circ \frR_{K/k}(v_{L/K})$. 
\end{proof} 

\begin{remark}\label{rem_norm} We put $v'_{L/K} = \frR_{L/K}(v_{L/K})$, and observe that $\rmN_{K/k} \circ v'_{L/K}$ is a homogeneous character of degree $p^h [K;k]$, while the norm $\rmN_{L/k}$, restricted to $\frR_{L/k}(\GG_{m,L})$, is a homogeneous character of degree $p^h [K:k]$, it follows from \eqref{thm_Lourdeaux2} that $\frac{[L:k]}{p^h [K:k]}$ is an integer and 
\begin{equation}  \label{rem_norm1}
\rmN_{L/k}= \frac{[L:k]}{p^h [K:k]} \, \big(\rmN_{K/k} \circ v'_{L/K}\big).
\end{equation}

One can of course prove directly that $\frac{[L:k]}{p^h [K:k]}$ is an integer, as follows.
We have $[L:k] = {[L:K]\cdot [K:k]}$, so the claim is that $p^h \; | \;  [L:K]$, which is a purely inseparable extension of height (= exponent) $h$. By definition of $h$, it contains an element of degree $p^h$, hence $p^h \;  | \; [L:K]$.
\end{remark}

\subsection{Norm map  for a  local finite $k$-algebra} \label{nomaprep} 
Let $k$ be a field and let $R$ be a local finite $k$-algebra. We denote by $\gm$ the maximal ideal of $R$ and by $L=R/\gm$ its residue field. The commutative $k$-group $G=\frR_{R/k}(\GG_{m,R})$ is  affine smooth and (geometrically) connected
\cite[A.5.11.(1) and (3)]{CGP}, 
and we have a norm map $\rmN_{R/k} : \frR_{R/k}(\GG_{m,R}) \to \GG_m$,  cf.~\eqref{trno-ds2}. 
On the other hand, we have a natural $k$-homomorphism
$s: \frR_{R/k}(\GG_{m,R}) \to \frR_{L/k}(\GG_{m,L})$, arising from the structure morphism $R \to L$. 

\comments{(2026-06-25) The claim \eqref{prop_norm1} can likely be generalized to the setting of \cite[A.3.5]{Oes}, but with $L \ne k$.}

\begin{prop} \label{prop_norm} We use the setting of {\rm \ref{nomaprep}}. The norm  character $\rmN_{R/k} :
 \frR_{R/k}(\GG_{m,k}) \to \GG_m$ factorizes as
follows
\[
\xymatrix{
\frR_{R/k}(\GG_{m,R})  \ar[rr]^{s} \ar[rd]_{\rmN_{R/k}} & &  
   \frR_{L/k}(\GG_{m,L})  \ar[ld]^{(\rmN_{L/k})^d} \\
& \GG_{m,k} 
}
\] 
where $d=[R:k]/[L:k] \in \ZZ$.
\end{prop}

\begin{proof} We claim that we have an exact sequence of commutative $k$--algebraic groups
\begin{equation}\label{prop_norm1}
 1 \to U \to G=\frR_{R/k}(\GG_{m,R}) \to \frR_{L/k}(\GG_{m,L}) \to 1.
\end{equation}
with $U$ being a split (commutative) $k$--unipotent group. 
We postpone the proof of \eqref{prop_norm1} until the second part of the proof. 

Assuming \eqref{prop_norm1} and using that $\Hom_{k-gp}(U, \GG_m)=1$ by \cite[XVII, Prop~2.4(ii)]{SGA3}, the norm character $\rmN_{R/k}: G \to \GG_{m,k}$ factorizes through  a $k$-homomorphism $w: \frR_{L/k}(\GG_m) \to \GG_{m,k}$. Since $\rmN_{R/k}$ is $[R:k]$-homogeneous, so is $w$. 
Hence, by  Theorem~\ref{thm_Lourdeaux}, there exists an integer $n$ such that $w= n \cdot (\rmN_{K/k} \circ v'_{L/K})$ where $v'_{L/K} = \frR_{L/K}(v_{L/K})$. By comparing the degrees of homogeneity of $w$ and $\rmN_{K/k} \circ v'_{L/K}$, we have
$[R:k]= n \, [K:k] \, p^h$. It then follows that 
\begin{align*}  
w &= \frac{[R:k]}{[K:k] \, p^h} \cdot (\rmN_{L/k} \circ v'_{L/K}) =
\frac{[R:k]}{[L:k]} \cdot \big( \frac{[L:k]}{[K:k] \,  p^h} \, (\rmN_{L/k} \circ v'_{L/K})\big) \\ & = \frac{[R:k]}{[L:k]} \cdot \rmN_{L/k} = (\rmN_{L/k})^d
\end{align*}  
by \eqref{rem_norm1}.
\sm 

It remains to prove \eqref{prop_norm1}, which for $L=k$ is proven in
\cite[A.3.5]{Oes}, see also \cite[Prop.~A.5.12]{CGP}. The proof  needs some preparation. Let $n\geq 0$ be the smallest integer such that $\gm^{n+1}=0$ and put $R_i=R/ \gm^i$ for $i=1,\ldots, n$, so that $R_1 = L$ and $R_{n+1}= R$. We consider $G_i=\frR_{R_i/k}(\GG_{m,R_i})$, for example $G=G_{n+1}$, and observe that $\Lie(G_i) \cong  \Lie(\GG_{m,R_i}) \cong R_i$ by \cite[Cor.~A.7.6]{CGP}. 
Since $G_i(A) = \GG_{m,R_i}(R_i \ot_k A) = (R_i \ot_k A)\ti$ for $A\in \kalg$ we have natural $k$--homomorphisms
\[ 
u_i: G_{i+1} \to G_i \qquad (i=1,...,n-1), 
\]
induced by the canonical epimorphism $R_{i+1} \to R_i$. 
Because $\Lie(u_i): R_{i+1} \to R_i$ is onto, it follows, e.g.~\cite[II, \S5, 5.1 and 5.3]{DG}, that each $u_i$ is smooth and surjective. 
\lv{
Reason: $G_{i+1}$ is smooth, affine, finitely presented $k$--groups, hence algebraic $k$--groups in the sense of \cite{DG}. It then follows from \cite[II, \S5, 5.1(d)]{DG} that $u_i(G_{i+1})_{\red} \to G_u$ is a closed immersion and from \cite[II, \S, 5.3]{DG} that $u_i$ is smooth and $u_i(G_{i+1})_{\red} \to G_u$ is an open immersion. But $G_i$ is also connected, so that $u_i(G_{i+1})_{\red}= G_i$, hence $u_i(G_{i+1}) = G_i$ because $u_i(G_{i+1})_{\red}$ and $u_i(G_{i+1})$ have the same underlying topological space. }

We claim: For $i=1,\ldots, n$, the $k$--group  $\ker(u_i)$ is isomorphic  to the $k$-vector group $V_i$ associated to $\gm^i/\gm^{i+1}$. To see this, we define a $k$-homo\-mor\-phism $v_i: V_i \to G_{i+1}$ as follows. For each $k$-algebra $A$, we map $[x] \in V_i(A) =(\gm^i/\gm^{i+1}) \otimes_k A$ to $1 + x \in {G_{i+1}(A)}=  (A \otimes_k R/\gm^{i+1})^\times$. It induces an exact sequence
$1 \to V_i(A) \to G_{i+1}(A) \to G_i(A)$. It follows that we have an exact sequence 
\[ 
1 \to V_i \to G_{i+1} \xrightarrow{\; u_i\; }  G_i \to 1
\] 
for $i=1,...,n$. Since $G=G_{n+1}$ and $G_1=\frR_{L/k}(\GG_m)$
the composite map
$u: G=G_{n+1} \to \frR_{L/k}(\GG_m)=G_1$ is smooth surjective
and its kernel is a split (commutative) $k$-unipotent $k$-subgroup $U$.
Thus, we have shown \eqref{prop_norm1}. \end{proof}

\newpage


\section{Cohomology}\label{sec:cohomology}
\comments{(2022-06-05) Need to compare with Auel's 2011 MPI preprint, last section; he seems to have similar results. }

\comments{(2026-04-29) Very preliminary version. We use \ref{lem_snake}\ref{lem_snake4} in \ref{abc}}

\subsection{Pointed sets}\label{poiset} We review the basic terminology and some easy results regarding pointed sets, see for example 
\cite[\S28]{KMRT}. \sm

A {\em pointed set\/} is a set with a distinguished element, usually denoted $(A, *_A)$ where $A$ is a set and $*_A$ is the distinguished element. For simpler notation we will abbreviate $A=(A,*_A)$. We view $\{*\}$ as an obviously pointed set.

A {\em morphism\/} $f\co A \to B$ of pointed sets is a set map $f \co A \to B$ satisfying $f(*_A) = *_B$. The {\em kernel $\Ker(f)$\/} of a morphism $f\co A \to B$ of pointed sets is $\Ker(f) = \{a \in A : f(a) = *_B\}$.
A sequence of morphisms of pointed sets
\[ A \xrightarrow{\;f\;} B \xrightarrow{\; \be\;}C\]
is {\em exact\/} if $\Ker(\be) = \Ima(f)$ where $\Ima (f) = \{f(a) : a\in A\}$ is the set-theoretic image. Observe that $A\xrightarrow{f} B \to \{*\}$ is exact if and only if $f$ is a surjective set map, while $\{*\} \to A \xrightarrow{f} B$ is exact if $\Ker(f) = \{*_A\}$, which does not imply that $f$ is an injective set map. Obviously, a morphism $f\co A \to B$ of pointed sets, whose underlying set map is injective, has trivial kernel.

\subsection{Some commutative diagrams}\label{poinlem} \begin{inparaenum}[\rm (a)]
\item \label{poinlem-a} Let
\[ \xymatrix@C=40pt{
   A \ar[r]^\psi \ar[d]_\al & B \ar[d]^\be \\
   A' \ar[r]^{\psi'} & B'
}\]
be a commutative diagram of pointed sets. Then there exists a unique homomorphism $\psi_1 \co \Ker(\al) \to \Ker(\be)$ of pointed sets such that
\[ \xymatrix@C=40pt{
  \Ker(\al) \ar[r]^{\psi_1} \ar[d]_\inc & \Ker(\be) \ar[d]^\inc \\
   A \ar[r]^{\psi} & B
}\]
is a commutative diagram of pointed sets, namely $\psi_1 = \psi|_{\Ker(a)}$.
\lv{
{\em Proof.} \eqref{poinlem-a}
For $a\in \Ker(\al)$ we have $(\be \circ \psi)(a) = \psi'(\al(a)) = \psi'(*_{A'}) = *_{B'}$, thus $\psi\big(\Ker(\al)\big) \subset \Ker (\be)$. \qed
}
\sm

\item \label{poinlem-d} ({\em Four-Lemma}, see for example \cite[Proof of Thm.~8.5]{First})
Let
\[
\xymatrix@C=40pt{
      &&& \{*\} \ar[d] \\
     A \ar[r]^\vphi \ar[d]^\al & B \ar[r]^\psi \ar[d]^\be & C \ar[d]^\ga \ar[r]^\xi
        & D \ar[d]^\de\\
   A' \ar[d] \ar[r]^{\vphi'} & B' \ar[r]^{\psi'} & C' \ar[r]^{\xi'} & D' \\
   \{*\}
}\]
be a commutative diagram of pointed sets with exact rows and columns, thus $\al$ is surjective and $\Ker(\de) = \{*\}$.  Then
\[ \text{$\be$ injective as set map} \quad \implies \quad \Ker(\ga) = \{*\}.\]
\sm

{\em Proof.} Let $c\in \Ker(\ga)$. By \eqref{poinlem-a}, $\xi(c) \in \Ker(\de) = \{*\}$. Thus $c\in \Ker(\xi) = \Ima(\psi)$, say $c=\psi(b)$. From $* = \ga(c) = \ga\psi(b) = \psi'\be(b) $ we find $b'= \be(b) \in \Ker(\psi') = \Ima(\vphi')$, thus $b' = \vphi'(a')$ for some $a'\in A'$ and then $a' = \al(a)$ for some $a\in A$ by surjectivity of $\al$. Furthermore, $\be(b) = b' = \vphi' \al(a) = \be \vphi(a)$. By injectivity of $\be$, we get $b= \vphi(a)$ and then $c=\psi(b) = (\psi\circ \vphi)(a) = *$. \qed
\end{inparaenum}

\comments{(2021-02-12) The following lemma to be used in Lemma~\ref{lem_snake}.
should be known. }

\begin{lem}\label{omc}
  Let $S$ be a scheme and let
\begin{equation}
  \label{omc1} \begin{split}
    \xymatrix{ 1 \ar[r] & A \ar[r] \ar[d]_{f_1} &G \ar[r]^q \ar[d]^f & B \ar[d]^{f_2}\ar[r] &  1 \\      1 \ar[r] & A' \ar[r] &G' \ar[r]^{q'} & B' \ar[r] &  1
  } \end{split}
\end{equation}
   be a commutative diagram of $S$--sheaves in the topology $\euT$ = flat, \'etale or Zariski with exact rows. Then
 \begin{equation}
   \xymatrix@C=45pt{ B(S) \ar[r]^{f_2(S)} \ar[d]_\vphi& B'(S) \ar[d]^{\vphi'}  \\
       H^1_\euT (S, A ) \ar[r]^{H^1_\euT(S, f_1)} & H^1_\euT(S, A')
   }
 \end{equation}
is a commutative diagram of pointed sets, where $\vphi$ and $\vphi'$ are the characteristic maps associated with the top and bottom row of \eqref{omc1}.
\end{lem}

\begin{proof} We fix $b\in B(S)$ and put $b' = f_2(b) \in B'(S)$. By definition, $\vphi(b)$ is the class of the $A$-sheaf torsor $q\me(b) \subset G(S)$. Similarly, $\vphi(b')$ is the class of the $A'$-sheaf torsor $q'{}\me(b') \subset G'(S)$, while $H^1_\euT(S, f)([q\me(b)])$ is the class of the contracted product $q\me(b) \we^A A'$ with $A$ acting on the left of $A'$ via $f_1$. Since the category of $A'$-sheaf torsors is a groupoid, it suffices to construct a morphism of $A'$-sheaf torsors $q\me(b) \we^A A' \to q'{}\me (b')$.

Recall that $q\me(b) \we^A A'$ is the quotient sheaf $(q\me(b) \times_S A')/A$ where $A$ acts on $q\me (b) \times_S A'$ by $(x,a') \cdot a = (x\cdot a, f_1(a\me) a')$ with $x\in q\me(b)$, $a'\in A'(S)$ and $a\in A(S)$. We have a map $\wtl \Phi \co q\me(b) \times_S A' \to q'{}\me(b')$ of sheaves over $S$, defined on the $S$--points by $\wtl\Phi(x,a') = f(x)a'$, which is constant on the $A$--orbits
and therefore descends to a well-defined sheaf map $\Phi \co q\me(b) \we^A A' \to q'{}\me(b')$.  We let $A'$ act on $q\me(b) \times_S A'$ by $(x,a') \cdot a_1'$. Since $\wtl \Phi$ is equivariant under this $A'$--action, it follows that $\Phi$ is equivariant too, i.e., $\Phi$  is a morphism of $A'$--torsors.
\end{proof}

\comments{(2021-02-18) Replaced ``fppf flat sheaves in groups" by " group sheaves over $S$ with exact rows and columns in the flat topology" }

\subsection{}\label{lem_snake-prep} Let $S$ be a scheme. We consider the following commutative diagram of group sheaves over $S$ with exact rows and columns in the flat topology:
\[
\xymatrix@C=40pt{
& 1 \ar[d] & 1 \ar[d] \\
 1 \ar[r] & A  \ar[r] \ar[d]  & \widetilde G \ar[r]^{q} \ar[d]& G \ar[r] \ar[d]^{\wr}& 1 \\
 1 \ar[r] & A' \ar[d]^{p} \ar[r] & G' \ar[d]^{p'} \ar[r]^{q'} & G \ar[r] & 1  .\\
& B \ar[d]  \ar[r]^{\sim} & B \ar[d]  \\
& 1  & 1
}
\]

\begin{lem}\label{lem_snake} The commutative diagram~{\rm \ref{lem_snake-prep}} gives rise to the following diagram of pointed sets
\begin{equation}\label{diag_giraud} \vcenter{
 \xymatrix@C=40pt{
 \ar @{} [dr] |{\rm (I)}
  G'(S) \ar[d]^{p'} \ar[r]^{q'}  &  \ar @{} [dr] |{(\rm II)} G(S)   \ar[d]^{\varphi} \ar[r]^{\varphi'} & H^1\fppf(S,A') \ar[d]^\cong  \\
  \ar @{} [dr] |{\rm (III)}
B(S) \ar[r]^{ \delta} \ar[d]^{\psi}& H^1\fppf(S, A )  \ar[d] \ar[r] & H^1\fppf(S,A')  \\
H^1\fppf( S, \widetilde G) \ar[r]^\sim & H^1\fppf( S, \widetilde G)
}}
\end{equation}
where the middle horizontal line arises from the exact sequence
$1 \to A \to A' \to B \to 1$ and  where $\psi$, $\varphi$ and $\varphi'$ are the characteristic maps attached to the sequences $1 \to  \widetilde G \to G'  \to B   \to 1$, $1 \to A \to \widetilde G \to G \to 1$ and $1 \to A' \to  G' \to G \to 1$ respectively. Then the following hold.
\begin{enumerate}[label={\rm (\alph*)}]
 \item \label{lem_snake1} The squares {\rm (II)} and {\rm (III)} are commutative.

 \item \label{lem_snake2} Let  $g' \in G'(S)$ with  images $g=q'(g') \in G(S)$ and $b=p'(g') \in B(S)$.
 Then $\varphi(g)=  \delta(b^{-1}) \in H^1\fppf(S,A)$.

 \item \label{lem_snake3} If $A$ is commutative, we have  $\varphi(g)=  - \, ^b \! \delta(b)$,
 where the action of $B(S)$ on $H^1\fppf(S,A)$ is induced by the outer
 action of $B$ on $A$.

 \item \label{lem_snake4} If $A$ is central in $A'$, the diagram {\rm (I)} anticommutes.  
\end{enumerate}
\end{lem}

\begin{proof} \ref{lem_snake1} The commutativity of the square (II) follows from Lemma~\ref{omc} by replacing $(G,B,B')$ there by $(\widetilde G, G, G)$ here. Similarly, replacing $(G,A',B')$ in Lemma~\ref{omc} by $(A', \widetilde G, B)$ proves the commutativity of the square (III).
\sm

\noindent \ref{lem_snake2}
We are given $g' \in G'(S)$ with images $g=q(g') \in G(S)$ and $b=p'(g') \in B(S)$.
By definition $\varphi(g)$ is the class of the sheaf
$A$--torsor ${q}^{-1}(g) \subset  \widetilde G$.
Similarly $\delta(b)$ is the class of the sheaf $A$--torsor $p^{-1}(b) \subset A'$.
We define an $A$-equivariant map $f :{q}^{-1}(g) \to G'$, $x \mapsto  (g')^{-1} x$.
For a flat $S$--scheme $T$ of finite presentation and  $x \in q^{-1}(g)(T)$, we have
  $q'(  (g')^{-1}_T x)=1$, so that $f(x) \in A'(T)$; furthermore we have $p'( (g')^{-1}_T x )=p'(g'_T)^{-1} =b^{-1}_T$. We have thus defined  an $A$--map $f : q^{-1}(g) \to p^{-1}(b^{-1})$. Since the category of $A$--torsors is a groupoid,
this map is an isomorphism. Thus $\varphi(g)= \delta(b^{-1})$.
\sm

\noindent \ref{lem_snake3} If $A$ is commutative, we have
$0=\delta(1)= \delta( b \, b^{-1} )= \, {^{b}\!\delta(b)} + \delta(b^{-1})  $
according to \cite[III.3.4, formula (3)]{Gir}, hence $\delta(b^{-1}) = \, - \, ^{b}\!\delta(b)$.
\sm

\noindent \ref{lem_snake4}. If $A$ is central in $A'$, the action of $B$ on $A'$
is trivial so that $^b \! \delta(b)=\delta(b)=-\varphi(g)$.
Thus the diagram (III) anticommutes.
\end{proof}
\newpage


\begin{thebibliography}{ABHSII}



\bibitem[ABHS]{ABHS} S.~Arpin, S.~Bozlee, L.~Herr, and H.~Smith, {\em The scheme of monogenic generators I: representability}, Res. Number Theory \textbf{9} (2023), no. 1, Paper No. 14, 33 pp.

\bibitem[ABHS-II]{ABHSII} \bysame, {\em  The scheme of monogenic generators II: local monogenicity and twists}, 
Res.~Number Theory \textbf{9} (2023), no. 2, Paper No. 43, 39 pp.

\bibitem[ABB]{ABB} A.~Auel, M.~Bernardara, and M.~Bolognesi, {\em Fibrations in complete intersections of  quadrics, Clifford algebras, derived categories, and rationality problems}, J.~Math.~Pures Appl.~\textbf{102} (2014), 249--291.

\bibitem[Bae]{Ba} R.~Baeza, {\it  Quadratic Forms over Semilocal Rings},
    Lecture Notes in Mathematics {\bf 655} (1978).

\bibitem[BC]{BC} P.~Balmer and B.~Calm\`es, {\em Bases of total Witt groups and lax-similitude}, J.~Algebra Appl.~\textbf{11} (2012), no.~3, 1250045.

\bibitem[BW]{Balmer-Walter} P.~Balmer and C.~Walter, {\em  A Gersten-Witt spectral sequence for regular schemes}, Ann. Sci. École Norm. Sup. (4) \textbf{35} (2002), no. 1, 127–152.


\bibitem[Bas1]{BR} H.~Bass, {\it Lectures on topics in algebraic
    K-theory}, Notes by Amit Roy, Tata Institute of Fundamental Research Lectures on Mathematics, No. 41 (1967), Bombay.

\bibitem[Bas2]{Bas2} \bysame, {\em Algebraic K--theory}, Benjamin, New York 1968. 

\bibitem[Bas3]{Bass-69} \bysame, {\em Modules which support nonsingular forms}, J.~Algebra \textbf{13} (1969), 246--252.

\bibitem[Bas4]{Bass-74} \bysame, {\em Clifford algebras and spinor norms over a commutative ring}, Amer.~J.~Math. \textbf{96} (1974), 156--206.

\bibitem[BFP]{BFP} E.~Bayer-Fluckiger, U.~A.~First, and R.~Parimala, {\em  On the Grothendieck-Serre Conjecture for Classical Groups}, J.~Lond.~Math.~Soc.~(2) \textbf{106} (2022), 2884--2926. 

\bibitem[BFL]{BFL} E.~Bayer-Fluckiger and H.~W.~Lenstra, Jr., {\em Forms in odd degree extensions and self-dual normal bases}, Amer.~J.~Math, \textbf{112} (1990), 359--373.

\bibitem[Bl]{Black} J.~Black, {\em Zero cycles of degree one on principal homogeneous spaces}, J.~Algebra \textbf{334} (2011), 232--246.

\bibitem[BQ]{BQ} J.~Black  and A.~Qu\'eguiner-Mathieu, {\em Involutions, odd degree extensions and generic splitting}, Enseign.~Math.~\textbf{60} (2014),  377--395.

\bibitem[Bo]{Bo} A. Borel, {\it  Linear Algebraic Groups (Second enlarged  edition)}, Graduate Text in Mathematics {\bf 126} (1991), Springer.


\bibitem[BT]{BT} A.~Borel and J.~Tits, {\em Homomorphismes ``abstraits'' de groupes algébriques simples}, Ann.~of Math.~(2) \textbf{97} (1973), 499--571.

\bibitem[BLR]{BLR} S.~Bosch, W.~L\"utkebohmert, M.~Raynaud, {\em N\'eron models}, Ergebnisse der Mathematik und ihrer Grenzgebiete, Band \textbf{21}, Springer-Verlag 1990.

\bibitem[B:A1]{BA} N.~Bourbaki, {\it Alg\`ebre}, Ch.1 \`a 3, Springer.

\bibitem[B:A2]{BA5} \bysame, {\em Alg\`ebre}, Ch.~4 \`a 7, Masson, Paris, 1981.

\bibitem[B:A3]{BA8} \bysame, {\em Alg\`ebre}, Ch.~8, second revised edition of the 1958 edition, Springer, Berlin 2012.

\bibitem[B:A4]{BA3} \bysame, {\em Algèbre}, Ch.~9, (French), reprint of the 1959 original, Springer-Verlag, Berlin, 2007.

\bibitem[B:AC]{BAC} \bysame, {\it Alg\`ebre commutative}, Ch.~1 \`a 4, Masson,
Paris 1985.

\bibitem[B:AC2]{BAC2} \bysame, {\em Alg\`ebre commutative}, Ch. 5\`a 7, Masson, Paris 1985.


\bibitem[BC]{Bout-Cesna} A.~Bouthier and K.~\v{C}esnavi\v{c}ius, {\em Torsors on loop groups and the Hitchin fibration}, 
    Ann.~Sci.~\'Ec.~Norm.~Supér.~(4) \textbf{55} (2022), no. 3, 791-–864.

\bibitem[CF]{CF} B. Calm\`es and J. Fasel, {\it Groupes classiques},
    in {\em Autour des sch\'emas en groupes. A celebration of SGA3}, vol II, Panoramas et Synth\`eses \textbf{46}, 1-133, Soci\'et\'e Math\'ematique de France, Paris, 2015.   

\bibitem[Ca]{carter} R.~W.~Carter, {\em Finite Groups of Lie Type: Conjugacy Classes and Complex Characters}, Wiley-Interscience, New York, 1985.

\bibitem[\v{C}e]{Cesna} K.~\v{C}esnavi\v{c}ius, {\em Grothendieck-Serre in the quasi-split unramified case},  Forum Math.~Pi~\textbf{10} (2022), Paper No.~9, 30 pp. 


\bibitem[CR]{CR} P.-E.~Chaput and M.~Romagny, {\em On the adjoint quotient of Chevalley groups over arbitrary base schemes}, J.~Inst.~Math.~Jussieu \textbf{9} (2010), 673--704.

\bibitem[Coh]{Cohn} P.~M.~Cohn, {\em On the decomposition of a field as a tensor product},  Glasgow Math.~J.~ \textbf{20} (1979) 141--145.

\bibitem[Col]{Collio} J.-L. Colliot-Th\'el\`ene, {\em Formes quadratiques sur les anneaux semi-locaux r\'eguliers}. Bull.~Soc.~Math.~France \textbf{59}, 13--31, 1979.

\bibitem[CoSa]{Collio-Sansuc} J.-L.~Colliot-Th\'el\`ene and J.-J.~Sansuc, {\em On the Chow groups of certain rational surfaces: a sequel to a paper of S. Bloch},   Duke Math. J. \textbf{48} (1981), 421-–447.

\bibitem[CoSk]{CSk} J.-L.~Colliot-Th\'el\`ene and A.~N.~Skorobogatov, {\em Groupe de Chow des z\'ero-cycles sur les fibr\'es en quadriques},  K-Theory {\bf 7} (1993), 477--500.

\bibitem[Con1]{Co1} B.~Conrad, {\em Reductive group schemes\/}, in {\em Autour
    des  sch\'emas en groupes, vol. I}, Panoramas et Synth\`eses \textbf{42-43},
    Soc. Math. France 2014.

\bibitem[Con2]{Co3} \bysame,  
  {\em Standard parabolics subgroups. Theory and
    Examples}, \\  http://virtualmath1.stanford.edu/~conrad/249BW16Page/handouts/stdpar.pdf
    
\bibitem[CGP]{CGP} B.~Conrad, O.~Gabber, and G.~Prasad, {\em Pseudo-reductive  groups}, second edition. New Mathematical Monographs \textbf{26}, Cambridge University Press, Cambridge 2015.


\bibitem[CRW]{CRW} T.~Craven, A.~Rosenberg, and R.~Ware, {\em The map of the Witt ring of a domain into the Witt ring of its field of fractions}, Proc. Amer. Math. Soc. \textbf{51} (1975), 25-–30.

\bibitem[DG]{DG} M.~ Demazure and P.~Gabriel, {\em Groupes alg\'ebriques},   North-Holland (1970).

\bibitem[DI]{DI} F.~DeMeyer and E.~Ingraham, {\em Separable Algebras Over     Commutative rings}, Lecture Notes in Mathematics \textbf{181} (1971),  Springer.

\bibitem[Di]{Dieu}   J. Dieudonn\'e,  {\it La g\'eom\'etrie des groupes
    classiques}, Ergebnisse der Mathematik und ihrer Grenzgebiete \textbf{5}, 3-i\`eme \'edition, Springer-Verlag Berlin-New York, 1971.

\bibitem[EGA]{EGA} A.~Grothendieck (avec la collaboration de J. Dieudonn\'e),     {\it El\'ements de G\'eom\'etrie Alg\'ebrique}, Publications  math\'ematiques de l'I.H.\'E.S. no.~4, 8, 11, 17, 20, 24, 28, 32, 1960--1967.

\bibitem[EGA-I]{EGA-neu} A.~Grothendieck and J.~Dieudonn\'e, {\em \'El\'ements de g\'eom\'etrie alg\'ebrique: I. Le langage des schémas}, Grundlehren der Mathematischen Wissenschaften  \textbf{166} (2nd ed.). Berlin; New York: Springer-Verlag, 1971.

\bibitem[EKM]{EKM} R.~Elman, N.~Karpenko, and A.~Merkurjev, {\em The
    algebraic and geometric theory of quadratic forms},
    Amer.~Math.~Soc.~Colloq.~Publ. \textbf{56},  Amer. Math. Soc., Providence,
    RI, 2008.

\bibitem[EG]{EG} D.~R.~Estes and R.~M.~Guralnick, {\em Module equivalences: Local to global when primitive polynomials represent units}, J.~Algebra \textbf{77} (1982), 138-–157.

  
\bibitem[Fed]{Fedorov} R.~Fedorov, {\em On the Grothendieck-Serre conjecture on principal bundles in mixed characteristic}, Trans.~Amer.~Math.~Soc. https://doi.org/10.1090/tran/8490, article electronically published on Nov.~5, 2021.

\bibitem[Fer]{Fer} D.~Ferrand, {\em Mongeneous Algebras. Back to Kronecker}, \arxiv{math/0310260}.

\bibitem[Fi1]{Fir21} U.~First, {\em On the non-neutral component of outer forms of the orthogonal group}, J.~Pure Appl.~Algebra \textbf{225} (2021), no.~1, Paper  No.~106477, 6 pp.

\bibitem[Fi]{First} U.~First, {\em An 8-periodic exact sequence of Witt groups of Azumaya algebras with involution}, Manuscripta Math. \textbf{170} (2023),  313--407.


\bibitem[FRS]{FRS} U.~First, Z.~Reichstein and S.~Salazar, {\em On the number of generators of a separable algebra over a finite field}, \arxiv{1709.06982}

\bibitem[Fo]{Ford} T.~J. Ford, {\it Separable Algebras}, Graduate     Studies in Mathematics {\bf 183}, Amer.~Math.~Soc, Providence, RI, (2017).



\bibitem[GPR1]{GPR} S.~Garibaldi, H.~P.~Petersson and M.~L.~Racine, {\em Albert algebras over $\ZZ$ and other rings},  Forum Math.~Sigma \textbf{11} (2023), Paper No.~e18, 38 pp. 

\bibitem[GPR2]{PRbook} \bysame,  
    {\em Albert algebras over commutative rings}, New Mathematical Monographs vol.~\textbf{48}, Cambridge University Press 2024.

\bibitem[Gi1]{Gille-Bourbaki} P.~Gille, {\em Le probl\`eme de Kneser-Tits},  (French, French summary) [The Kneser-Tits problem], 
S\'eminaire Bourbaki, Vol.~2007/2008, Ast\'erisque No.~\textbf{326} (2009), Exp.~ No.~983, vii, 39--81 (2010).


\bibitem[Gi2]{G2}  \bysame, {\it Sur la classification des groupes   semi-simples}, Autour des sch\'emas en groupes (III), Panoramas et  Synth\`eses \textbf{47}, Soci\'et\'e Math\'ematique de France 2015.

\bibitem[Gi3]{G2019} \bysame, {\em Groupes alg\'ebriques semi-simples en dimension cohomologique $\leq 2$}, Lecture Notes in Mathematics \textbf{2238}, Springer Nature Switzerland AG 2019.

\bibitem[Gi4]{G2021} \bysame, {\em When is a reductive group scheme linear?}, \arxiv{2103.07305}, to appear in Michigan Mathematical Journal (2021).

\bibitem[GN1]{GN-Sp} P.~Gille and E.~Neher, {\em Springer's Odd Degree Extension Theorem for Quadratic Forms over semilocal Rings},  Indag. Math.~(N.S.) \textbf{32} (2021), no. 6, 1290--1310.   

\bibitem[GN2]{GN-LG} \bysame, {\em Group schemes over LG-rings and applications to cancellation theorems and Azumaya algebras}, Eur.~J.~Math. \textbf{11} (2025), no. 3, Paper No. 42, 67 pp. 
    
\bibitem[GS]{GS} P.~Gille and T.~Szamuely, {\em Central simple algebras and Galois cohomology}, Second edition, Cambridge Studies in Advanced Mathematics \textbf{165}, Cambridge University Press, Cambridge, 2017.


\bibitem[Gir]{Gir} J.  Giraud, {\em Cohomologie non-ab\'elienne}, 
  Grundlehren der mathematischen Wissenschaften \textbf{179},  Springer-Verlag, Berlin-New York, 1971. 

\bibitem[Go]{Goldman}O.~Goldman, {\em Determinants in projective modules}, Nagoya Math.~J. \textbf{18} (1961), 27--36.

\bibitem[GW]{GW} U.~G\"ortz and T.~Wedhorn, {\em Algebraic geometry I. Schemes -- with examples and exercises}, Second edition, Springer Studium Mathematik—Master, Springer Spektrum, Wiesbaden 2020.
    
\bibitem[GWII]{GWII} \bysame, {\em Algebraic geometry II: Cohomology of schemes -— with examples and exercises}, Springer Studium Mathematik—Master, Springer Spektrum, Wiesbaden 2023.  

\bibitem[Gr]{GroI} A.~Grothendieck, {\em Le groupe de Brauer I. Alg\`ebres d'Azumaya et interpr\'etations diverses}, S\'em.~Bourbaki 1964/65, no.~290, in {\em Dix Expos\'es sur la Cohomologie des Sch\'emas}, North-


\bibitem[Ho]{Ho2} D.~W.~Hoffmann, {\em Similarity of quadratic and symmetric bilinear forms in characteristic 2}, Indagationes Mathematicae (2020), https://doi.org/10.1016/j.indag.2020.08.008.

\bibitem[Hu]{Hum} J.~Humphreys, {\em Conjugacy Classes in Semisimple Algebraic Groups}, Mathematical Surveys and Monographs {\bf 43} (1995), American Mathematical Society, Providence, RI.

\bibitem[Ja1]{Jac1} N.~Jacobson, {\em Some groups of transformations defined by Jordan algebras. I.}, J.~Reine~Angew.~Math.~\textbf{201} (1959), 178--195.

\bibitem[Ja2]{Jac2} \bysame, {\em Generic norm of an algebra}, Osaka Math.~J.~\textbf{15} (1963), 25--50.


\bibitem[Ja]{Jan} J.~C.~Jantzen, {\em Representations of Algebraic Groups}, Mathematical Surveys and Monographs \textbf{107}, American Mathematical Society, Providence R.I., second edition 2003.

\bibitem[KS]{KaSa} K.~Kato and S.~Saito, {\it Unramified class field theory of arithmetical surfaces}, Annals of Math. {\bf 118} (1983), 241-275.

\bibitem[Knb1]{Knebusch} M.~Knebusch, {\em Isometrien \"uber semilokalen Ringen},  Math.~Z.~\textbf{108} (1969), 255--268.

\bibitem[Knb2]{Knebusch-habil} \bysame, {\em Grothendieck- und Wittringe von nichtausgearteten symmetrischen Bilinearformen}, (German) S.-B. Heidelberger Akad.~Wiss.~Math.-Natur.~Kl. 1969(70) (1969/70), 93--157. 

\bibitem[Knb3]{knebusch-norm} \bysame, {\em Ein Satz \"uber die Werte von quadratischen Formen \"uber K\"orpern}, Invent.~Math.~\textbf{12}, 300-–303 (1971)

\bibitem[Knb4]{Knebusch73} \bysame, {\em Specialization of quadratic and symmetric bilinear forms, and a norm theorem}, Acta Arith.~\textbf{24} (1973), 279--299.

\bibitem[Knb5]{Knebusch-Queens} \bysame, {\em Symmetric bilinear forms over algebraic varieties},  Conference on Quadratic Forms 1976 (Proc. Conf., Queen's Univ., Kingston, Ont., 1976), pp. 103–283. Queen's Papers in Pure and Appl. Math., No. 46, Queen's Univ., Kingston, Ont., 1977.
    
\bibitem[Kns1]{kneser-Sem}  M.~Kneser, {\em Witt's Satz \"uber quadratische Formen und die Erzeugung orthogonaler Gruppen durch Spiegelungen}, Math.-Phys. Semesterber. \textbf{17}, 33--45 (1970)

\bibitem[Kns2]{Kneser} \bysame,  {\it Quadratische Formen},
    Springer-Verlag, Berlin, 2002.
    
\bibitem[Knu]{K} M.-A.~Knus, {\it Quadratic and Hermitian Forms over
    Rings}, Grundlehren der mathematischen Wissenschaften {\bf 294}
    (1991), Springer.

\bibitem[KMRT]{KMRT} M.-A. Knus, A. Merkurjev, M. Rost and J.-P. Tignol, {\it   The  Book of Involutions}, Amer.~Math.~Soc.~Colloquium Publ.,  \textbf{44}, Amer.~Math.~Soc., Providence, RI (1998).

\bibitem[KO]{KO} M.-A. Knus, and M. Ojanguren, {\em Th\'eorie de la Descente et   Alg\`ebres d'Azumaya}, Lecture Notes in Mathematics \textbf{389} (1974), Springer

\bibitem[KO1]{KO77} \bysame, {\em Modules and Quadratic Forms over Polynomial Algebras}, Proc.~Amer.~Math.~Soc.~\textbf{66} (1977), 223--226.


\bibitem[KOS]{KOS} M.-A. Knus, M.~Ojanguren, and D.~J.~Saltman, {\em On Brauer groups in characteristic p}, Brauer groups (Proc.~Conf., Northwestern Univ., Evanston, Ill., 1975), Springer, Berlin, 1976, pp. 25--49. Lecture Notes in Math., Vol. \textbf{549}.




\bibitem[Lam]{Lam-qf} T.~Y.~Lam, {\em An Introduction to Quadratic Forms
    over Fields}, Graduate Studies in Mathematics \textbf{67}, Amer.~Math.~ Soc., Providence, RI, 2005.

\bibitem[Lam2]{La} \bysame, {\it Serre's Problem on Projective Modules}, Springer Monograph in Mathematics (2006), Berlin.

\bibitem[Lan]{Lang} S.~Lang, {\em Algebra\/}, 3rd edition, Addison-Wesley, 1993.

\bibitem[Lee]{Lee} T.-Y.~Lee, {\em Adjoint quotients of reductive groups}, Autour des sch\'emas en groupes, Vol.~III, 131--145, Panor.~Synth\`eses \textbf{47} (2015), Soc.~Math.~France, Paris

\bibitem[Len]{Lenstra} H.~W.~Lenstra, {\em Galois Theory for Schemes}, electronic 3rd edition, online resource.


\bibitem[Lo1]{Lo-ag} O.~Loos, {\em On algebraic groups defined by Jordan pairs}, Nagoya Math.~J.~\textbf{74} (1979), 23--66.

\bibitem[Lo2]{Lo0} \bysame, {\it Tensor products and discriminants of
    unital quadratic forms over commutative rings}, Mh. Math. \textbf{122}
       (1996), 45--98.

\bibitem[Lo3]{Lo-genalg} \bysame, {\em Generically algebraic Jordan algebras
    over commutative rings}. J.~Algebra \textbf{297} (2006), 474--529.

\bibitem[Lo4]{Lo2} \bysame, {\it Cubic and symmetric compositions over
    rings}, manuscripta math. \textbf{124} (2007), 195--236.
    

\bibitem[Lou]{Lourdeaux} A.~Lourdeaux, {\em On geometry of pseudo-reductive groups}, 
Comm. Algebra \textbf{50} (2022), no.~12, 5371-–5386.
    
\bibitem[MW]{MW} B.~R.~McDonald and W.~C.~Waterhouse, {\em Projective modules over rings with many units}, Proc.~Amer.~Math.~Soc.~\textbf{83}(3) (1981), 455-–458.

\bibitem[MR]{MicRev} A.~Micali and P.~Revoy, {\em Modules quadratiques}, Bull.~Soc.~Math. France, M\'emoire \textbf{63} (1979).

\bibitem[Mi]{milne:ag} J.~Milne, {\em Algebraic Groups -- The Theory of Group Schemes of Finite Type over a Field}, Cambridge Studies in Advanced Mathematics \textbf{170}, Cambridge University Press, Cambridge, 2017.

\bibitem[NSW]{NSW} J.~Neukirch, A.~Schmidt and J.~Wingberg, {\it  Cohomology of Number Fields}, second edition, Grundlehren der mathematischen Wissenschaften 323 (2008), Springer.

\bibitem[Oe1]{Oes} J.~Oesterl\'e, {\em Nombres de Tamagawa et groupes unipotentes en caract\'eristique $p$}, Invent.~Math.~\textbf{78}(1), 13-–88.

\bibitem[Oe2]{Oesterle} \bysame, {\em Sch\'emas en groupes de type multiplicatif}, (French) [Group schemes of multiplicative type], Autour des sch\'emas en groupes. Vol. I, 63--91, Panor. Synth\`eses, 42/43, Soc.~Math.~France, Paris, 2014.

\bibitem[Oj1]{O} M.~Ojanguren, {\em Quadratic forms over regular rings},
J.~Indian Math.~Soc.~(N.S.) \textbf{44} (1980), no. 1--4, 109–-116).

\bibitem[Oj2]{Ojan82} \bysame, {\em A splitting theorem for quadratic forms}, Comment.~Math.~Helv. \textbf{57} (1982), no. 1, 145-–157.

\bibitem[Oj3]{Oj82} \bysame, {\em Unit\'es repr\'esent\'ees par des formes quadratiques ou par des normes r\'eduites}, in {\em Algebraic K-theory, Part II} (Oberwolfach, 1980), pp. 291--299, Lecture Notes in Math. \textbf{967}, Springer, Berlin-New York, 1982.

\bibitem[OP]{OP} M.~Ojanguren and I.~Panin, {\em A purity theorem for the Witt group},  Ann.~Sci.~\'Ecole Norm.~Sup.~ \textbf{32} (1999), 71--86.
  
\bibitem[OPZ]{OPZ} M.~Ojanguren, I.~Panin, and K.~Zainoulline, {\em On the norm principle for quadratic forms}, J.~Ramanujan Math.~Soc.~\textbf{19} (2004) 1-–12.
  
\bibitem[OS]{OS} M.~Ojanguren and R.~Sridharan, {\em Cancellation of Azumaya algebras}, J.~Algebra \textbf{18} (1971), 501--505.

\bibitem[PP]{PP} I.~Panin and K.~Pimenov, {\em Rationally isotropic quadratic spaces are locally isotropic: II}, Doc.~Math.~2010, Extra vol.: Andrei A.~Suslin sixtieth birthday, 515–523.

\bibitem[PR]{panin-rehmann} I.~Panin and U.~Rehman, {\em A variant of a theorem by Springer}, St.~Petersburg Math.~J. \textbf{19} (2008), 953--959.

\bibitem[PS]{PS} R.~Parimala, and R.~Sridharan, {\em  Quadratic forms over rings of dimension $1$},  Comment.~Math.~Helv.~\textbf{55} (1980), 634--644.



\bibitem[Pe1]{Pet} H.~P.~Petersson, {\em Idempotent $2\times 2$ matrices}, Seminarberichte aus dem Fachbreich Mathematik der FernUniversit\"at Hagen \textbf{78} (2007), pp.~157--172 = paper 242 on the Jordan Theory Preprint Archives,  http://agt2.cie.uma.es/~loos/jordan/

\bibitem[Pe2]{P-Fields} \bysame, {\em A survey on Albert algebras}, Transform.~Groups \textbf{24} (2019), 219-–278.


\bibitem[Po]{Po} B.~Poonen, {\em Rational points on varieties}, Graduate Studies in Mathematics \textbf{186}, Amer.~Math.~Soc., Providence, RI, 2017.

\bibitem[Ra1]{Raynaud} M.~Raynaud, {\em Faisceaux amples sur les sch\'emas en groupes et les espaces homog\`enes}, Lect.~Notes in Math.~\textbf{119}, Springer--Verlag, 1970.

\bibitem[Ra2]{Ray-hensel} \bysame, {\em Anneaux Locaux Hens\'eliens}, Lect.~Notes in Math.~\textbf{169}, Springer--Verlag, 1970.


\bibitem[Rob]{Roby} N.~Roby, {\em Lois polynomes et lois formelles en th\'eorie des modules}, Ann. Sci. \'Ecole Norm. Sup. \textbf{80} (1963), 213-–348.


\bibitem[Rom]{Romagny} M.~Romagny, {\em G\'eom\'etrie Alg\'ebrique 2},  cours Universit\'e Pierre et Marie Curie, 2011-2012.

\bibitem[Sah]{Sah} C.~H.~Sah, {\em Symmetric bilinear forms and quadratic forms},
    J. Algebra {\bf 20} (1972), 144-160.


\bibitem[Sal]{Salt} D.~J.~Saltman, {\em Lectures on division algebras}, CBMS Regional Conference Series in Mathematics, vol. \textbf{94}, Amer.~Math.~Soc., Providence, RI, 1999.


\bibitem[Sch1]{Scha} W.~Scharlau, {\em Zur Pfisterschen Theorie der quadratischen Formen}, Invent.~Math.~\textbf{6} (1969), 327--328.

\bibitem[Sch2]{Sc} \bysame, {\it Quadratic and Hermitian Forms}, Grundlehren der Mathematischen Wissenschaften \textbf{270}, Springer-Verlag Berlin 1985.



\bibitem[SGA3]{SGA3} {\it S\'eminaire de G\'eom\'etrie alg\'ebrique de l'I.H.E.S., 1963-1964, sch\'emas en groupes, dirig\'e par M. Demazure et A. Grothendieck},  Lecture Notes in Math. 151-153. Springer (1970).

\bibitem[SGA4$_3$]{SGA4} {\em Th\'eorie des topos et cohomologie \'etale des sch\'emas}, Tome 3 (French),  S\'eminaire de G\'eom\'etrie Alg\'ebrique du Bois-Marie 1963--1964 (SGA 4). Dirig\'e par M.~Artin, A.~Grothendieck et J. L. Verdier. Avec la collaboration de P.~Deligne et B.~Saint-Donat. Lecture Notes in Mathematics  \textbf{305}, Springer-Verlag, Berlin-New York, 1973.

\bibitem[SGA7]{SGA7} P.~Deligne and N.~Katz, {\em Groupes de monodromie en g\'eom\'etrie alg\'ebrique II}, Lecture Notes in Mathematics \textbf{340}, Springer--Verlag, New York, 1973.

\bibitem[Ses]{Seshadri}  C.~S.~Seshadri, {\em Geometric reductivity over arbitrary base},  Advances in Math.~\textbf{26} (1977), 225--274.

\bibitem[Si]{Sivatski}  A.~S.~Sivatski, {\em Similarity of quadratic forms and related problems}, J.~Pure Appl.~Algebra \textbf{223} (2019), 4102-–4121.

\bibitem[Sp1]{springer-qf} T.~A.~Springer, {\em Sur les formes quadratiques d'indice z\'ero}, C.~R.~Acad.~Sci.~Paris \textbf{234} (1952), 1517--1519.

\bibitem[Sp2]{Springer} \bysame, {\it Linear algebraic groups}, Second edition (1998), Birkh\"auser.


\bibitem[St]{St} The Stacks Project Authors, {\em Stacks project}, \url{http://stacks.math.columbia.edu/}

\bibitem[Ste]{Ste} R.~Steinberg, {\em Regular elements of semisimple algebraic groups}, Inst.~Hautes \'Etudes Sci.~Publ.~Math.~\textbf{25} (1965), 49--80.

\bibitem[Ste2]{Ste-mem} \bysame, {\em Endomorphisms of linear algebraic groups}, Mem.~Amer.~Math.~Soc., vol \textbf{80}, 1968.


\bibitem[Sw]{Sw} R.~G.~Swan, {\em $K$-theory of quadric hypersurfaces},
    Ann. of Math. \textbf{122} (1985), 113--153.

\bibitem[Ti]{Tits} J.~Tits, {\em Groupes de Whitehead de groupes alg\'ebriques simples sur un corps (d’apr\`es V.~P.~Platonov et al.)}, (French)
      S\'eminaire Bourbaki, 29e ann\'ee (1976/77), Exp.~No.~505, pp. 218--236, Lecture Notes in Math., \textbf{677}, Springer, Berlin, 1978.

\bibitem[vdB]{VB} M.~Van~den~Bergh, {\em The Brauer-Severi scheme of the trace ring of generic matrices}, Perspectives in ring theory (Antwerp, 1987), 333--338, NATO Adv. Sci. Inst. Ser. C Math. Phys. Sci., \textbf{233}, Kluwer Acad. Publ., Dordrecht, 1988.

\bibitem[vdK]{vdK} W.~van der Kallen, {\em The $K_2$ of rings with many units}, Ann.~Sci.~\'Ecole Norm.~Sup.~\textbf{10} (1977), 473--515.


\bibitem[Za]{Kirill} K.~Zainoulline, {\em On Knebusch's norm principle for quadratic forms over semi-local rings}, Math.~Z.~\textbf{251} (2005), 415-–425.

\end{thebibliography}
\end{document}